\documentclass[10pt]{amsbook}

\usepackage{graphicx}        
\usepackage{multicol}        
\usepackage{hyperref}

\usepackage{tikz-cd}
\usepackage{calc,xassoccnt}
\NewTotalDocumentCounter{totalfigures}
\NewTotalDocumentCounter{totaltables}
\NewTotalDocumentCounter{appendixchapters}
\DeclareAssociatedCounters{figure}{totalfigures}
\DeclareAssociatedCounters{table}{totaltables}

\makeindex             

\definecolor{uclwhite}{rgb}{1, 1, 1}

\definecolor{ucl1dkpurple}{RGB}{82,66,91}
\definecolor{ucl2dkpurple}{RGB}{134,122,140}
\definecolor{ucl3dkpurple}{RGB}{168,160,173}
\definecolor{ucl4dkpurple}{RGB}{220,217,222}

\definecolor{ucl1dkred}{RGB}{90,27,49}
\definecolor{ucl2dkred}{RGB}{139,95,110}
\definecolor{ucl3dkred}{RGB}{172,141,152}
\definecolor{ucl4dkred}{RGB}{222,209,214}

\definecolor{ucl1dkblue}{RGB}{0,67,89}
\definecolor{ucl2dkblue}{RGB}{76,123,138}
\definecolor{ucl3dkblue}{RGB}{127,161,172}
\definecolor{ucl4dkblue}{RGB}{204,217,222}

\definecolor{ucl1dkgreen}{RGB}{75,70,32}
\definecolor{ucl2dkgreen}{RGB}{129,125,98}
\definecolor{ucl3dkgreen}{RGB}{165,162,143}
\definecolor{ucl4dkgreen}{RGB}{219,218,210}

\definecolor{ucl1black}{RGB}{0,0,0}
\definecolor{ucl2black}{RGB}{75,75,75}
\definecolor{ucl3black}{RGB}{128,128,128}
\definecolor{ucl4black}{RGB}{205,205,205}

\definecolor{ucl1pink}{RGB}{145,24,83}
\definecolor{ucl2pink}{RGB}{178,93,134}
\definecolor{ucl3pink}{RGB}{200,139,169}
\definecolor{ucl4pink}{RGB}{233,209,221}

\definecolor{ucl1mdred}{RGB}{195,58,45}
\definecolor{ucl2mdred}{RGB}{213,117,108}
\definecolor{ucl3mdred}{RGB}{225,156,150}
\definecolor{ucl4mdred}{RGB}{243,216,213}

\definecolor{ucl1mdblue}{RGB}{69,156,189}
\definecolor{ucl2mdblue}{RGB}{124,186,209}
\definecolor{ucl3mdblue}{RGB}{162,205,222}
\definecolor{ucl4mdblue}{RGB}{218,235,242}

\definecolor{ucl1mdgreen}{RGB}{130,141,55}
\definecolor{ucl2mdgreen}{RGB}{167,175,115}
\definecolor{ucl3mdgreen}{RGB}{192,198,155}
\definecolor{ucl4mdgreen}{RGB}{230,232,215}

\definecolor{ucl1orange}{RGB}{215,123,35}
\definecolor{ucl2orange}{RGB}{227,162,101}
\definecolor{ucl3orange}{RGB}{235,189,145}
\definecolor{ucl4orange}{RGB}{247,229,211}

\definecolor{ucl1ltpurple}{RGB}{191,175,188}
\definecolor{ucl2ltpurple}{RGB}{210,199,208}
\definecolor{ucl3ltpurple}{RGB}{223,215,221}
\definecolor{ucl4ltpurple}{RGB}{242,239,242}

\definecolor{ucl1yellow}{RGB}{229,175,0}
\definecolor{ucl2yellow}{RGB}{237,199,76}
\definecolor{ucl3yellow}{RGB}{242,215,127}
\definecolor{ucl4yellow}{RGB}{250,239,204}

\definecolor{ucl1ltblue}{RGB}{168,192,209}
\definecolor{ucl2ltblue}{RGB}{194,211,223}
\definecolor{ucl3ltblue}{RGB}{211,223,232}
\definecolor{ucl4ltblue}{RGB}{238,242,246}

\definecolor{ucl1brtgreen}{RGB}{204,209,88}
\definecolor{ucl2brtgreen}{RGB}{219,223,138}
\definecolor{ucl3brtgreen}{RGB}{229,232,171}
\definecolor{ucl4brtgreen}{RGB}{245,246,222}

\definecolor{ucl1stone}{RGB}{217,214,204}
\definecolor{ucl2stone}{RGB}{228,226,219}
\definecolor{ucl3stone}{RGB}{236,234,229}
\definecolor{ucl4stone}{RGB}{255,255,255}

\definecolor{ucl1ltgreen}{RGB}{185,193,147}
\definecolor{ucl2ltgreen}{RGB}{206,211,179}
\definecolor{ucl3ltgreen}{RGB}{220,224,201}
\definecolor{ucl4ltgreen}{RGB}{241,243,233}

\usepackage{amssymb,mathtools}
\usepackage{tikz}
\usepackage{tikz-qtree}
\usepackage{tkz-euclide}
\usepackage{pgfplots}
\usepackage{enumerate}
\makeatletter
\let\@@enum@org\@@enum@
\def\@@enum@[#1]{\@@enum@org[\normalfont #1]}
\makeatother

\DeclarePairedDelimiter{\ceil}{\lceil}{\rceil}
\DeclarePairedDelimiter{\floor}{\lfloor}{\rfloor}
\DeclarePairedDelimiter{\form}{\langle}{\rangle}

\newcommand\ba{\begin{align*}}
\newcommand\ea{\end{align*}}
\newcommand\be{\begin{enumerate}}
\newcommand\ee{\end{enumerate}}
\newcommand\bp{\begin{proof}}
\newcommand\ep{\end{proof}}
\newcommand\bpp{\begin{prop}}
\newcommand\epp{\end{prop}}
\newcommand\bd{\begin{defn}}
\newcommand\ed{\end{defn}}

\numberwithin{equation}{subsection}

\newcommand\fform[1]{\langle\!\langle #1\rangle\!\rangle}

\newcommand\bC{\mathbb{C}}
\newcommand\bP{\mathbb{P}}

\newcommand\bN{\mathbb{N}}
\newcommand\bR{\mathbb{R}}

\newcommand\bF{\mathbb{F}}
\newcommand\bQ{\mathbb{Q}}
\newcommand\bZ{\mathbb{Z}}
\newcommand\Z{\mathbb{Z}}
\newcommand\bH{\mathbb{H}}
\newcommand\AAA{\mathcal{A}}
\newcommand\BB{\mathcal{B}}

\newcommand\CC{\mathcal{C}}

\newcommand\FF{\mathcal{F}}

\newcommand\GG{\mathcal{G}}
\newcommand\DD{\mathcal{D}}

\newcommand\II{\mathcal{I}}
\newcommand\JJ{\mathcal{J}}
\newcommand\KK{\mathcal{K}}
\newcommand\OO{\mathcal{O}}

\newcommand\PP{\mathcal{P}}
\newcommand\CS{\mathcal{S}}

\newcommand\UU{\mathcal{U}}
\newcommand\VV{\mathcal{V}}
\newcommand\MM{\mathcal{M}}

\newcommand\Sym{\operatorname{Sym}}
\newcommand\onto{\twoheadrightarrow}
\newcommand\into{\hookrightarrow}

\DeclarePairedDelimiter{\abs}{\lvert}{\rvert}
\DeclarePairedDelimiter{\abss}{\lVert}{\rVert}

\newcommand\Aff{\operatorname{Aff}}
\newcommand\PL{\operatorname{PL}}
\newcommand\GL{\operatorname{GL}}
\newcommand\Hom{\operatorname{Hom}}

\DeclareMathOperator\Inv{Inv}
\DeclareMathOperator\Int{int}
\DeclareMathOperator\cl{cl}
\DeclareMathOperator\ext{ext}
\DeclareMathOperator\comp{comp}
\DeclareMathOperator\Conf{Conf}

\newcommand\rot{\operatorname{rot}}

\newcommand\supp{\operatorname{supp}}
\newcommand\suppo{\operatorname{supp}^\circ}

\newcommand\suppc{\operatorname{\overline{supp}}}
\newcommand\suppe{\operatorname{{supp}^{\mathrm{e}}}}
\newcommand\Ro{\operatorname{Ro}}

\newcommand\Id{\operatorname{Id}}

\newcommand\Var{\operatorname{Var}}
\newcommand\Diffb{\operatorname{Diff}_+^{1+\mathrm{bv}}}

\newcommand\diam{\operatorname{diam}}

\newcommand\gam{\Gamma}
\newcommand\Mod{\operatorname{Mod}}

\DeclareMathOperator\Homeo{Homeo}

\newcommand\eps{\epsilon}
\newcommand\yt{\widetilde}
\newcommand\sse{\subseteq}
\newcommand\co{\colon}
\newcommand\der[2]{#1^{(#2)}}

\DeclareMathOperator\tr{tr}
\DeclareMathOperator\Fix{Fix}
\DeclareMathOperator\BS{BS}
\DeclareMathOperator\Out{Out}
\DeclareMathOperator\Aut{Aut}

\DeclareMathOperator\PSL{PSL}

\DeclareMathOperator\SL{SL}
\DeclareMathOperator\SO{SO}
\DeclareMathOperator\Diff{Diff}

\DeclareMathOperator\Stab{Stab}

\DeclareMathOperator\CL{CovLen}
\DeclareMathOperator\CD{CovDist}
\DeclareMathOperator\CR{CritReg}
\DeclareMathOperator\CM{{CM}}
\DeclareMathOperator\Nbr{{Nbr}}
\DeclareMathOperator\TF{{TopFilt}}

\newtheorem{thm}{Theorem}[section]
\newtheorem{lem}[thm]{Lemma}
\newtheorem{cor}[thm]{Corollary}
\newtheorem{prop}[thm]{Proposition}
\newtheorem{con}[thm]{Conjecture}
\newtheorem{que}[thm]{Question}
\newtheorem{claim}{Claim}

\newtheorem{notation}[thm]{Notation}
\newtheorem{prin}[thm]{Principle}
\theoremstyle{definition}
\newtheorem{rem}[thm]{Remark}
\newtheorem{defn}[thm]{Definition}
\newtheorem{exmp}[thm]{Example}

\begin{document}

\author{Sang-hyun Kim, Thomas Koberda}
\title{Structure and regularity of group actions on one-manifolds\\
}
\maketitle

\frontmatter


%
%

\chapter*{Preface}

This book represents an account and contextualization of a research program 
that was executed by the authors and their collaborators over the period
of several years. Both authors began their careers in classical geometric group theory and began collaborating on
projects about right-angled Artin groups and their relationship with mapping class groups of surfaces, at a time
when the first author was a visiting assistant professor at Tufts University and while the second
author was a graduate student at Harvard University.

Around 2014, the authors became interested in a question attributed by M.~Kapovich to V.~Kharlamov, 
one that asks which right-angled Artin groups can act faithfully by smooth
diffeomorphisms on the circle. Because of the profusion of right-angled Artin subgroups of mapping class groups
 and because of
the robust persistence of right-angled Artin groups under passing to finite index subgroups,
an answer to Kharlamov's question had the potential
to shed light on a question of F.~Labourie that he posed in his 1998 ICM 
talk and that was reiterated by B.~Farb and J.~Franks in several places:
are there finite index subgroups of (sufficiently complicated) mapping class groups acting smoothly
on the circle?

At first, it seemed to the authors that an easy solution to Kharlamov's question was available, and that all right-angled Artin groups admitted
such faithful actions. It was soon pointed out by J.~Bowden that the actions constructed by the authors and H.~Baik in fact failed
to be differentiable at certain accumulations of fixed points, and so only furnished faithful smooth actions on the real line. 
After approximately two
years of work, inspired by Bowden's observation, Baik and the authors managed to prove that that most right-angled Artin groups
do not admit faithful $C^2$ actions on the compact interval nor on the circle. As a consequence, Labourie's question could be completely
answered in the case of $C^2$ actions of finite index subgroups of mapping class groups.

About a year and a half later, the authors were able to develop a tool that they dubbed the \emph{$abt$--Lemma}, which is a certain
combinatorial--algebraic obstruction for smoothness of a group action on a compact one--manifold. As a result, they were able
to answer Kharlamov's question by giving a concise characterization of right-angled Artin groups acting
faithfully by smooth diffeomorphisms on the circle.

Because of the technical details
involved in the analysis of smooth actions of right-angled Artin groups on one--manifolds, the authors became interested
in Thompson's group $F$, which occurs naturally in these contexts. Inspired by the dynamical theory of Thompson's group, 
the authors investigated the class of
chain groups in subsequent work with Y.~Lodha.
These objects form a highly diverse class of groups of homeomorphisms of the interval, while nevertheless exhibiting remarkable uniformity
properties.

Over the years that the preceding story unfolded, the authors had long wondered about the existence of finitely generated groups of
homeomorphisms which can act by $C^k$ diffeomorphisms, but which admit no faithful actions by $C^{k+1}$ diffeomorphisms, especially
in the case $k\geq 2$. 
The question of the existence of such groups was posed explicitly by A.~Navas in his 2018 ICM talk. Having built a
foundation consisting of the theory of right-angled Artin group actions on one--manifolds and the theory of chain groups, the authors
were able to construct explicit examples of finitely generated groups that act on the circle and on the interval with prescribed
\emph{critical regularities}\index{critical regularity}, in the following strong sense: 
every smoother action of a finite index subgroup factors through an abelian quotient.

More recent work of the authors with C.~Rivas has improved this state of knowledge, particularly concerning right-angled Artin group
and mapping class group actions on the circle. It is now known that 
for all $\epsilon>0$, a finite index subgroup of a mapping class group (that is not virtually free--by--cyclic) does not
admit a faithful action by $C^{1,\epsilon}$ diffeomorphisms on the interval nor on the circle.

The historical development of mathematical ideas is rarely systematic, which is why we recount the story of the authors' contribution to
the theory in this book in a chronological way, so as to complement the systematic exposition in the body of the monograph. Moreover, the
drive within academia to publish research articles that continually break new ground is in direct tension with the cultural and 
intellectual need of the mathematics community
to have access to complete, detailed, and polished accounts of mathematical subjects. It is thus that we approached the project of writing
this book: we construct a coherent narrative that gives a systematic overview of the theory of critical regularity of groups, and at the same
time provide sufficient background and context so that a reasonably experienced graduate student could traverse the gap from introductory
courses to the cutting edge in the span of these pages. The result is a book that is half lecture notes, half research monograph.

Herein, the reader will find all the results mentioned above, with complete proofs. In many cases, the proofs we have provided are
more efficient and more general than the ones given in the original publications. Notably, the precise statement of the result that furnishes
finitely generated groups of prescribed critical regularity is substantially stronger than the one in the authors' original paper.

In addition to the account of the authors' contributions,
there is an exposition on background and contextual results that are scattered through both
the published literature as well as various unpublished sources. 
One chapter is dedicated to what is known as Denjoy's theory of diffeomorphisms.
Many excellent books and expository articles give proofs of Denjoy's Theorem, though our presentation is tailored to our agenda.
In the commentary surrounding Denjoy's Theorem, we have included a self--contained discussion of invariant and stationary measures.
We avoid long digressions into amenability and measurable dynamics, proceeding merely with our narrative needs in mind. 

The section on
differentiable Denjoy counterexamples includes results of the authors which have hitherto appeared only in research journals.

A large chapter of the book
is devoted to results of F.~Takens, R.~Filipkiewicz, and M.~Rubin, which for us mostly serve to contextualize 
the theory of critical regularity. To the authors'
knowledge, this will be the first instance of a treatment of these results as parts of a coherent whole. Moreover, the authors hope that
this chapter will serve to clear up certain inaccuracies, misattributions, and inefficiencies present in the literature.

One chapter of the book, together with the appendices, introduces standard tools for analyzing group actions on one--manifolds. The exposition
therein concentrates on the $C^0$, the $C^1$, and the $C^2$ theories. We have recounted proofs of nearly all the results, sometimes
reproducing proofs that already exist in the literature, and other times streamlining known proofs.
The remaining chapters of the book were written mostly in service to proving the existence of groups with prescribed critical regularity,
and to a commentary on the applications of the construction.

The authors sincerely hope that the reader will find this book intellectually stimulating and satisfying, and humbly wish that they
may attract the attention
of fresh minds to this beautiful subject.
\vspace{\baselineskip}
\begin{flushright}\noindent
Seoul and Charlottesville,\hfill {\it Sang-hyun Kim}\\
June 2021\hfill {\it Thomas Koberda}\\
\end{flushright}

\chapter*{Acknowledgements}

The authors thank Adam Clay, Benson Farb, Andreas Holmsen, Nam-Gyu Kang, Yash Lodha, Curtis T.~McMullen, Crist\'obal Rivas and 
William Winston for helpful discussions and corrections.
The authors are especially grateful to Andrés Navas for his inspiring book and pioneering works in the direction of this work.
The authors are also grateful to Elena Griniari for guiding the authors throughout the creation of this book.
The first author is supported by Samsung Science and Technology Foundation (SSTF-BA1301-51).
The second author is partially supported by  an Alfred P. Sloan Foundation Research Fellowship, by NSF Grant DMS-1711488, and by 
NSF Grant DMS-2002596.

\tableofcontents

\mainmatter


%
%
%
\chapter{Introduction}\label{ch:intro} 

\begin{abstract}
In this monograph, we give an account of the relationship between the algebraic structure of finitely generated and countable
groups and the regularity with which they act on manifolds. We concentrate on the case of one--dimensional manifolds, culminating
with a uniform construction of finitely generated groups acting with prescribed regularity on the compact interval and on the circle. We develop
the theory of dynamical obstructions to smoothness, beginning with classical results of Denjoy, to more recent results of Kopell, to
modern results such as the $abt$--Lemma. We give a classification of the right-angled Artin groups that have finite critical
regularity, and discuss their exact critical regularities in many cases, and we compute the 
virtual critical regularity of most mapping class groups of orientable surfaces.\end{abstract}

This book is a discussion of the relationship between algebraic structure of groups, dynamics of group actions, and regularity of diffeomorphisms. Particularly,
we are interested in finitely generated, or more generally countable groups, acting on a fixed compact manifold, and how the level
of differentiability of the action controls the possible algebraic structure of the groups.

The following diagram illustrates the basic slogans of this book.
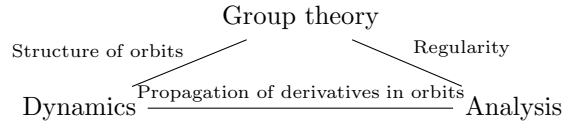
\[\begin{tikzcd}
 & \textrm{Group theory} \arrow[dr,dash,"\textrm{Regularity}"] \\
\textrm{Dynamics} \arrow[ur,dash,"\textrm{Structure of orbits}"] \arrow[rr,dash,"\textrm{Propagation of derivatives in orbits}"] 
&& \textrm{Analysis}
\end{tikzcd}
\]

At one of the vertices, we have group theory, which for us means the description of the algebraic structure of groups. We are
primarily concerned with infinite groups that are classically of interest in geometric group theory and related areas, such as mapping
class groups of hyperbolic surfaces, right-angled Artin groups, and the Higman--Thompson type groups such as Thompson's group $F$.
These groups are so interesting to us because they embody a rich combination of commutativity and non--commutativity,
that is, \emph{partial commutativity}\index{partial commutativity}.

Algebraically,
partial commutativity is interesting and informative, because it is tractable (in the sense that many defining relations are easy to write
down and understand intuitively) and because it is complicated enough to exhibit a wide range of phenomena; indeed, the theory
of right-angled Artin groups shows that in a precise sense, ``group with partial commutation" are at least as complicated as the class
of all finite graphs, which in turn are a model for all of finitary mathematics. Partial commutativity is also useful for studying groups
up to commensurability. Indeed, commutation of group elements is robust under replacing a group element with a nonzero power.
Thus, partial commutativities in an ambient group generally propagate to finite index subgroups. This is very helpful in understanding
groups whose finite index subgroups are too complicated for contemporary technology, such as mapping class groups of hyperbolic surfaces.

Much of this book is concerned with groups of diffeomorphisms. Commutation
of diffeomorphisms, as we will explore in great detail, is a rare phenomenon, and imposes a high 
degree of structure on pairs of diffeomorphisms
which exhibit it, especially in higher regularity. In fact, 
the following could be considered as a principle that informs a large portion of the writing of this
book:

\begin{prin}\label{prin:main}
The Mean Value Theorem makes it hard for diffeomorphisms of a compact manifold to commute with each other.
\end{prin}

Principle~\ref{prin:main} illustrates an interaction between the vertices labeled ``group theory" and ``analysis", though it
conceals the important role of the third vertex. Principle~\ref{prin:main}
highly constrains actions of groups that exhibit a high degree of partial commutativity, by forcing group
elements to have a lot of fixed points, and then higher derivatives at accumulation points of the fixed points
(which is where compactness comes in!) often fail to be continuous, due
to the orbit structure of the group action.

Dynamics' domain is that of the global structure of orbits, and interacts in an essential way with both the structure of a group that is acting,
and with the topology of the underlying manifold. All of our discussion would be greatly hobbled, if not completely invalidated, were the
manifolds under consideration not restricted to one dimension. The interval and circle are ordered structures, and their topologies are
completely determined by the orders. The consequence of this interaction between the order relation and the topology is roughly that,
away from fixed points, homeomorphisms of the interval and of the circle have a direction. This last statement means that the combinatorics
of orbits are much more tractable than they would be on an arbitrary manifold, or even worse, on an arbitrary topological space. Without
a systematic way to analyze and classify possible orbit structures, the edge connecting the ``group theory" and the ``analysis" vertices
would be inscrutable.

By investigating the three vertices and edges of the diagram, we are able to carry out a number of constructions and to illustrate several
completed mathematical programs. Among the highlights of the book are the following:

\begin{itemize}
\item
A complete classification of right-angled Artin groups which admit faithful actions on the interval and the circle with regularity at least two,
in terms of their defining graphs.
\item
A complete classification of mapping class groups that (virtually) admit faithful actions on the interval and the circle with regularity at least two,
in terms of the topology of the underlying surface.
\item
A uniform construction of finitely generated groups with prescribed \emph{critical regularity}\index{critical regularity}
when acting on the interval or on the circle. Here,
the critical regularity of a group $G$ is the supremum of the real numbers $r=k+\tau$ with $k\in\bZ_{>0}$ and $\tau\in [0,1)$, where $G$
admits a faithful action by $C^k$ diffeomorphisms whose $k^{th}$ derivatives are $\tau$--H\"older continuous.
\item
New, robust constructions of codimension one foliations on closed $3$--manifolds subject to mild topological hypotheses, which have
prescribed regularity properties.
\item
Several technical obstructions to $C^1$ smoothability of group actions that do not rely on Thurston's Stability Theorem, which start
with a relatively weak hypothesis (the orbit structure is subjected to a naturally occurring dynamical constraint) and end with a relatively strong
conclusion (the action was not differentiable).
\end{itemize}

We move on now beyond philosophical musings to discuss the content and context of this book in more precise terms.

\subsection{Some conventions and notation}

Throughout this book, we let $\bZ_{>0}$ denote the set of positive integers, 
and let $\bN$ denote the set of \emph{natural numbers}\index{natural numbers}. That is, we have
\[
\bN:=\{0,1,2,\ldots\}=\{0\}\cup \bZ_{>0}.\]
We denote by $I$ a compact interval; typically we mean $I=[0,1]$ unless specified otherwise.
A map or manifold is said to be \emph{smooth}\index{smooth manifold} if it is $C^\infty$.

\section{Groups of manifold diffeomorphisms}

In this section, we give an overview of the main actors in the sequel. A reader familiar with the basics of manifold theory and homeomorphism
groups can safely skim or skip this section.

One of the basic objects in this book is a \emph{manifold}\index{manifold} $M$. 
For us, a manifold is always connected and second countable. Since $M$
is connected and hence path--connected, the \emph{dimension}\index{dimension of a manifold}
of $M$ will be well-defined, as this quantity is locally constant.

The entire discussion of this book occurs in relation to, and usually entirely within, the group $\Homeo(M)$ of homeomorphisms
\[\phi\colon M\longrightarrow M.\] When $M$ is orientable, there is a natural subgroup $\Homeo_+(M)$ consisting of orientation
preserving homeomorphisms of $M$. When $\Homeo_+(M)$ and $\Homeo(M)$ do not coincide, we will usually restrict our attention
to $\Homeo_+(M)$, since this latter group is just an index two subgroup in the full group of homeomorphisms.

Manifolds, \emph{a priori}, are only topologically locally homeomorphic to Euclidean spaces, and there are many subtle questions about
differentiable manifolds, piecewise geometric manifolds, and the relationships between these categories. 

Since these are not the concern
of this book, we will always assume that manifolds are smooth, which is to say that they are implicitly equipped with an atlas of charts
whose transition functions are smooth  (that is, $C^\infty$) local diffeomorphisms.

With the sole exception of Chapter~\ref{sec:filip-tak}, we will be working in low dimensions. Within the context of low dimensions, we will
be working almost exclusively in dimension one. Dimension two will only occur in relation to mapping classes of surfaces, and dimension
three will occur only in the context of the existence of non-smoothable codimension one foliations. Thus, the distinction between the
topological category, the $C^1$ category, the smooth category, and the piecewise linear category of manifolds is nil, for our purposes
(cf.~\cite{Moise77book,Thurston1997book}).

Thus, for a manifold $M$, we will be justified in filtering $\Homeo(M)$ by subgroups that are analytically defined. The most coarse filtration
we will use is by integral regularities. We will write $\Diff^k(M)$ for the group of $C^k$ diffeomorphisms of $M$. Thus, $f\in\Diff^k(M)$
if and only if $f$ and $f^{-1}$ are continuous and have continuous derivatives up to order $k$. That $\Diff^k(M)$ is indeed a group
is a consequence of the chain rule. We thus obtain a descending chain of subgroups \[\Homeo(M)>\Diff^1(M)>\cdots>\Diff^k(M)>\cdots.\]
It is convenient to allow $\Diff^0(M)=\Homeo(M)$, and we will implicitly adopt this convention throughout.

When $M$ is orientable, then for all $k$ we can also decorate the group $\Diff^k(M)$
with the plus sign by restricting to the orientation preserving diffeomorphisms, identifying $\Diff_+^k(M)$ with a subgroup of $\Homeo_+(M)$.
The group $\Diff^{\infty}(M)$ is defined
by \[\Diff^{\infty}(M)=\bigcap_{k\geq 1} \Diff^k(M),\] and $\Diff_+^{\infty}(M)$ is defined in the obvious way. 

We remark that there are interesting groups of diffeomorphisms of manifolds that lie even deeper than $\Diff^{\infty}(M)$. For instance,
if $M$ has an analytic or algebraic or symplectic structure, then one can consider the group of diffeomorphisms that respect this structure. If
$M$ is equipped with a smooth Riemannian metric then one can consider the isometries of $M$, which is to say the diffeomorphisms
which respect this metric. One can consider other geometric structures that are not necessarily metric, strictly speaking, such as projective
structures, and thus get many more interesting groups of diffeomorphisms. These kinds of groups, while rich in examples and questions,
are beyond the scope of this book.

One group which will occur a number of times and which does not fit easily into any of the groups of diffeomorphisms we have defined so far
is the group $\PL_+(M)\le \Homeo_+(M)$
of \emph{piecewise linear}\index{piecewise linear} homeomorphisms of $M$, where $M\in\{I,\bR,S^1\}$. This name, though standard,
is a bit of a misnomer, since elements of $\PL(M)$ are technically \emph{piecewise affine}\index{piecewise affine}.
An element $f\in\PL_+(M)$ described as being locally defined
by affine maps of $M$ (where here the affine structure on $S^1$ is by the identification of $S^1=\bR/\bZ$) outside of a finite set of
points (called the \emph{breakpoints}\index{breakpoint} of $f$).
The derivative of $f$ is locally constant away from the breakpoints of $f$ and is usually
not continuous at the breakpoints, but we insist that $f$ be continuous. The fact that $\PL_+(M)$ is a group is an easy exercise.
The relationship between $\PL_+(M)$ and $\Diff_+^{\infty}(M)$ will be investigated in Chapter~\ref{ch:chain-groups}.

We can refine the filtration of $\Homeo(M)$ by diffeomorphism groups of integral regularity to non--integral regularities. Let
$f\colon X\longrightarrow Y$ be continuous, where $X$ and $Y$ are metric spaces.
Recall that $f$ is \emph{H\"older continuous}\index{H\"older continuity} 
with exponent $\tau>0$ if there exists a constant $C>0$ such that for all $a,b\in X$, we have
\[d_Y(f(a),f(b))\leq C\cdot \left(d_X(a,b)\right)^{\tau}.\] The function $x\mapsto x^{\tau}$ is called a
\emph{modulus of continuity}\index{modulus of continuity}. The most interesting values of $\tau$ to consider lie in $(0,1]$, where the case
of $\tau=1$ is called \emph{Lipschitz continuity}\index{Lipschitz continuity}.
We say $f$ is \emph{locally $\tau$--H\"older continuous}\index{local H\"older continuity}
if it is H\"older continuous with the exponent $\tau$ in some neighborhood at each point. 
A local modulus of continuity coincides with a global one if the space is compact
or if the function is compactly supported, i.e.~when the function is uniformly continuous. 

For $\tau\in(0,1)$, we write $\Diff^{k,\tau}(M)$ for the set of $C^k$ diffeomorphisms of $M$ whose $k^{th}$ derivatives 
are locally $\tau$--H\"older
continuous. 
It is true that for $k\geq 1$, the set $\Diff^{k,\tau}(M)$ does indeed form a group, though some argument is needed to establish
this fact; we relegate further details to Appendix~\ref{ch:append1}.
The elements of this group are called $C^{k,\tau}$ diffeomorphisms.
For $\tau=1$, we write $\Diff^{k,\mathrm{Lip}}(M)$ for the $C^k$ diffeomorphisms 
of $M$ whose $k^{th}$ derivatives are locally Lipschitz continuous.

We emphasize again that a regularity ($C^k$ or $C^{k,\tau}$) is a local property.
This convention makes the theory more consistent; for instance, $C^2$ diffeomorphisms are 
necessarily $C^{1,\mathrm{Lip}}$ although they may not have globally Lipschitz first derivatives.  

The set of homeomorphisms $\Diff^{k,\mathrm{Lip}}(M)$ of $M$ forms a group, even when $k=0$.
While it is true that Lipschitz maps are differentiable almost everywhere,
the distinction between almost everywhere differentiability and everywhere differentiability is significant. Whereas it is typical practice to write
$\Diff^{k+\tau}(M)$ for $\Diff^{k,\tau}(M)$, one must be careful when $\tau=1$.
Indeed, note that
\[\Diff^{k+1}(M)\neq\Diff^{k,1}(M)=
\Diff^{k,\mathrm{Lip}}(M).\]

These various groups of homeomorphisms and diffeomorphisms of manifolds have their own intrinsic topology, coming
inherited from the Whitney $C^k$
topology on all differentiable functions $M\longrightarrow M$. 
Oftentimes, these groups can be thought of as infinite dimensional
Lie groups that are locally Banach or locally Fr\'echet, though we will usually not need this aspect of their structures.

Sometimes
the point--set topology of these groups does play a role, however, as we have already suggested. Indeed, these groups are often
not connected, this being typified by the existence of an orientation--reversing homeomorphism of $M$. Even restricting to
orientation preserving homeomorphisms, the group $\Homeo_+(M)$ and its differentiable subgroups may not be connected. This arises
from the fact manifolds often have nontrivial \emph{mapping class groups}\index{mapping class group},
which is to say groups of isotopy classes of homeomorphisms.
The mapping class group $\Mod(M)$ is generally defined to be $\pi_0(\Homeo(M))$, which has the structure of a group since $\Homeo(M)$
itself has the structure of a group. The topology on $\Mod(M)$ is the quotient topology, which in many cases of interest is simply the
discrete topology. In some parts of this book, primarily in Chapter~\ref{sec:filip-tak}, we will need connectedness of groups of homeomorphisms
and diffeomorphisms; for this reason, when it is relevant, we will restrict to the connected component of the identity in these groups.

A further complication arises for noncompact manifolds, which we will address when necessary: in dealing with noncompact
manifolds, we will usually consider only groups of compactly supported homeomorphisms or diffeomorphisms when investigating the
noncompact manifolds for their own sakes. Occasionally, we will consider noncompactly supported homeomorphisms of noncompact
manifolds that are lifts of a homeomorphism of a compact quotient. Of course, for compact manifolds, these preceding remarks are
generally irrelevant.

The connected component of the identity oftentimes has nontrivial topology, such as a nontrivial fundamental group. For the interval,
the group of orientation preserving homeomorphisms and $C^k$ diffeomorphisms is contractible, where as for the circle the corresponding
group is homotopy equivalent to the circle itself. This nontriviality of the fundamental group of $\Homeo_+(S^1)$ will be used to
pass freely between homeomorphisms of the circle and periodic homeomorphisms of the real line.

\section[The Mather--Thurston Theorem, the Epstein--Ling Theorem]{The Mather--Thurston Theorem, the Epstein--Ling Theorem, and lattice-like rigidity}

The Mather--Thurston Theorem, while itself not a subject of detailed discussion in this book, 
has deep implications regarding Haefliger's classifying spaces for foliations and the homology of diffeomorphism groups. 
In its essence, the theorem says that diffeomorphism groups of manifolds are algebraically rigid, 
in the sense that the homomorphisms between them are rather limited.
We state a version of this fact that is particularly helpful for contextualizing the content
of this monograph. 

\begin{thm}[See~\cite{Mather1,Mather2,Mather3,Thurston1974BAMS}]\label{thm:mather-thurston}
Let $M^n$ be a smooth connected boundaryless manifold of dimension $n\ge1$. 
For $k\in\bZ_{>0}\cup\{\infty\}$, let $\Diff_c^k(M)_0 $ denote the group of compactly
supported diffeomorphisms of $M$ that are isotopic to the identity by a compactly supported isotopy. Then the group $\Diff_c^k(M)_0 $ is
simple, provided that $k\neq n+1$.
\end{thm}

We remark that in the literature, the Mather--Thurston Theorem often refers to the perfectness of the groups $\Diff_c^k(M)_0 $, i.e.~the fact
that
\[[\Diff_c^k(M)_0,\Diff_c^k(M)_0 ]=\Diff_c^k(M)_0 .\] The simplicity of the full group $\Diff_c^k(M)_0$ then follows from much more general
principles due to Higman, Epstein--Ling, and others; see Chapter~\ref{sec:filip-tak} for a detailed discussion of the simplicity of commutator
subgroups of diffeomorphism groups, and the discussion of Theorem~\ref{thm:epstein-ling-intro} below.

For the reader who is not familiar with Theorem~\ref{thm:mather-thurston}, the hypothesis on the dimension may appear bizarre,
and has to do with finding a fixed point of a certain operator. For finite $k$, when $k\neq n+1$ then this operator can be shown to admit
a fixed point, but the fixed point finding procedure fails otherwise. It is not clear whether this failure of the method of proof when $k=n+1$
is an artifact or a feature; indeed, Mather proves in~\cite{Mather3} that the group $\Diff_c^{1+\mathrm{bv}}(\bR)$ of compactly supported
diffeomorphisms of $\bR$ whose first derivatives have bounded variation is not perfect. By decomposing the distributional derivative 
$D(\log f')$ for $f\in \Diff_c^{1+\mathrm{bv}}(\bR)$ into its regular and singular parts and by integrating against the regular part of the
measure, Mather shows that one obtains a surjective homomorphism from $\Diff_c^{1+\mathrm{bv}}(\bR)$ to $\bR$.
We direct the reader to~\cite{CKK2019} for a detailed exposition of
Theorem~\ref{thm:mather-thurston} in the case $n=1$.

More classically, it was known from the work of Epstein and Ling that the commutator subgroup
of $\Diff_c^k(M)_0$ is simple.
\begin{thm}\label{thm:epstein-ling-intro}
Let $M^n$ be a smooth connected boundaryless manifold of dimension $n$.
For $k\in\bZ_{>0}\cup\{\infty\}$, let $\Diff_c^k(M)_0 $ be as in Theorem~\ref{thm:mather-thurston}. Then the group
 \[\Diff_c^k(M)_0 '=[\Diff_c^k(M)_0 ,\Diff_c^k(M)_0 ]\]
is simple.
\end{thm}

We shall prove Theorem~\ref{thm:epstein-ling-intro}
in Chapter~\ref{sec:filip-tak}. The argument that $\Diff_c^k(M)_0 '$ is simple is not sensitive to the dimension
of $M$, and is general enough to be applicable to a wide range of other groups. The step from Theorem~\ref{thm:epstein-ling-intro}
to Theorem~\ref{thm:mather-thurston} is in computing the homology group $H_1(\Diff_c^k(M)_0 ,\bZ)$, where $\Diff_c^k(M)_0 $ is viewed
as a discrete group.

As a consequence of Theorem~\ref{thm:epstein-ling-intro}
and Theorem~\ref{thm:mather-thurston}, we have the following observation:

\begin{cor}\label{cor:diffeo-homo}
Let $M^m$ and $N^n$ be smooth connected boundaryless manifolds of dimensions $m$ and $n$, and let $k,\ell\in\bZ_{>0}\cup\{\infty\}$.
If  \[\phi\colon\Diff_c^k(M)_0 \longrightarrow\Diff_c^{\ell}(N)_0\]
is a non-injective homomorphism, then its image is abelian; if in addition $k\neq \dim M+1$, the image is trivial.
\end{cor}

In practical applications, it is not always so easy to show that the map $\phi$ in Corollary~\ref{cor:diffeo-homo} fails to be injective. Sometimes
it is injective: indeed, if $M=N$ and if $\ell\leq k$, then there is a natural inclusion between 
the diffeomorphism groups. Were one to try to prove that
there is no injective homomorphism between these diffeomorphism groups when $\ell>k$ from first principles, it is not immediately apparent
how to proceed. Why, after all, should there be no way to find a simultaneous smoothing of all $C^k$ diffeomorphisms of $M$? The answer
is that one cannot, and we will return to this question in Section~\ref{sec:filip-tak-intro} below;
the reader may also consult~\cite{HurtadoGT15,Kramer11,Mann2015,Mann2016GT,Mann-Notices,RoseSol07,Rybicki95}
for some background results in this direction.

Before proceeding, we will take this opportunity to clarify some terminology that is relevant to the discussions in this book. There are
two distinct notions of smoothability of groups of homeomorphisms or diffeomorphisms, namely 
\emph{topological smoothability}\index{topological smoothability} and
\emph{algebraic smoothability}\index{algebraic smoothability}.
If $G$ is a group of homeomorphisms or diffeomorphisms of a manifold $M$, then a topological smoothing
of $G$ into $\Diff^{k,\tau}(M)$ is given by a topological conjugacy of $G$ into $\Diff^{k,\tau}(M)$.

That is, there is a homeomorphism
$h\co M\longrightarrow M$ such that \[hGh^{-1}\le \Diff^{k,\tau}(M).\] 
If the regularity of the elements of $G$ is lower than $(k,\tau)$ then clearly
$h$ cannot be an element of $\Diff^{k,\tau}(M)$, and it is best to simply record $h$ as a homeomorphism. We remark that 
\emph{in the literature
on group actions, a smoothing of a group is usually a topological smoothing}.

An algebraic smoothing of $G$ is simply an injective homomorphism into $\Diff^{k,\tau}(M)$. Clearly a topological smoothing is an
algebraic smoothing, though the converse many not hold. A topological smoothing of $G$ respects all the dynamical properties of the action
of $G$, whereas an algebraic smoothing is blind to everything but the algebraic structure of $G$. There will be several places in this book
where we will explicitly be interested in topological smoothings, and these will be clearly noted as such; see Corollary~\ref{c:slp-0}, for instance. Otherwise, the primary
interest herein will be algebraic smoothings of groups.

One of the primary goals of this book is to give examples of finitely generated groups of diffeomorphisms that have a prescribed level
of regularity. We will defer stating the full statement of the relevant results until Section~\ref{sec:critreg-intro} below, stating only a relevant
consequence.

\begin{prop}\label{prop:unsmooth}
For  $k\in\bZ_{>0}$ and $\tau\in [0,1)$
and for $M\in\{I,S^1\}$, there exists a finitely generated nonabelian group $G_{k,\tau}$ with a simple commutator subgroup
such that there is an injective homomorphism
\[\phi\colon G_{k,\tau}\longrightarrow \Diff_+^{k,\tau}(M),\] but such that
for every finite index subgroup $H_{k,\tau}\le G_{k,\tau}$, every homomorphism \[\psi\colon H_{k,\tau}\longrightarrow
\Diff_+^{\ell,\eps}(M)\] has abelian image, where $\ell\in\bZ_{>0},\epsilon\in[0,1)$ and $k+\tau<\ell+\eps$.
\end{prop}

The conclusion of Proposition~\ref{prop:unsmooth} cannot be improved to triviality of the image. This is a consequence of Thurston's
Stability Theorem, which will be discussed in Appendix~\ref{ch:append3}. Thurston proves that a finitely generated group of diffeomorphisms
of $\Diff^1[0,1)$ always surjects to $\bZ$ and hence cannot be simple.

Proposition~\ref{prop:unsmooth} gives a finitary (in the sense that it is determined by finitely many group elements) obstruction to the
algebraic smoothability of a group. We thus obtain the following consequence from Theorem~\ref{thm:mather-thurston} and
Theorem~\ref{thm:epstein-ling-intro}:

\begin{cor}\label{cor:diffeo-homo-fg}
Let $M\in\{\bR,S^1\}$, and suppose that $k+\tau<\ell+\eps$. Then every homomorphism \[\Diff_c^{k,\tau}(M)_0\longrightarrow\Diff_+^{\ell,\eps}(M)\]
has trivial image, with the possible exception when $k=2$ and $\tau=0$, in which case every such homomorphism has abelian image.
\end{cor}

The absence of an analogue of Corollary~\ref{cor:diffeo-homo-fg} in higher dimensions is one of the motivations for developing a theory
of critical regularity for finitely generated groups in dimensions two and higher. In dimension one, the question of whether $\Diff_c^2(\bR)$
is simple remains open, and our finitary methods do not shed any light on this question.

The groups furnished by Proposition~\ref{prop:unsmooth} enjoy rigidity properties that are reminiscent of various rigidity results for
lattices in higher rank Lie groups. The most general of these is Margulis' Superrigidity Theorem. We state a simplified version of this
result, since the full generality would not illuminate the principle any further.

\begin{thm}[Margulis Superrigidity, simplified statement~\cite{MargulisBook1991,Witte2015}]\label{thm:margulis-super}
Let $\Gamma$ be a non-cocompact lattice in $\mathrm{SL}_{k\ge3}(\bR)$, and let \[\rho\colon \gam\longrightarrow\mathrm{GL}_n(\bR)\]
be a homomorphism. Then there is a continuous homomorphism \[\widehat\rho\colon \mathrm{SL}_k(\bR)\longrightarrow\mathrm{GL}_n(\bR)\]
such that $\rho$ and $\widehat\rho$ coincide on a finite index subgroup of $\Gamma$.
\end{thm}

Here, a \emph{lattice}\index{lattice} is a discrete subgroup in $\mathrm{SL}_k(\bR)$ whose covolume (with respect to the Haar measure) is
finite.

Theorem~\ref{thm:margulis-super} relies on a great deal of structure to which we simply do not have access. The analogy between the
groups furnished by Proposition~\ref{prop:unsmooth} and nonuniform (i.e.~non-cocompact) lattices in $\mathrm{SL}_k(\bR)$ is 
simply that homomorphisms from them to a group of more regular diffeomorphisms are extremely limited. 

\begin{cor}\label{cor:super-analogue}
Let $G_{k,\tau}$ be as in Proposition~\ref{prop:unsmooth}.
Suppose $k+\tau<\ell+\eps$. If $H\le \Diff_+^{\ell,\eps}(M)$ is a subgroup, then every homomorphism \[G_{k,\tau}\longrightarrow H\] has
abelian image.
\end{cor}

\section{The Takens--Filipkiewicz Theorem and Rubin's Theorem}\label{sec:filip-tak-intro}

Another facet to the line of inquiry pursued in this book emerges from the following basic question: to what degree is a mathematical object
determined by its group of symmetries? If $M$ is a manifold, it is reasonable then to wonder the degree to which $M$ is determined by
its natural group of symmetries, which in the context of the foregoing discussion would be its group of homeomorphisms $\Homeo(M)$.
If $M$ has further structure, such as a smooth structure, a complex structure, a symplectic structure, etc.~then one can pose an analogous
question as to the degree to which $M$ with the further structure is preserved by the group of symmetries of $M$ that preserve that structure.

In 1963, Whittaker~\cite{Whittaker1963} proved that an abstract isomorphism \[\phi\colon\Homeo(M)\longrightarrow\Homeo(N)\] between
two arbitrary manifolds arises from a homeomorphism between $M$ and $N$. That is,
there is a homeomorphism $\alpha\colon M\longrightarrow N$ such that for all $f\in\Homeo(M)$, we have that
\[\phi(f)=\alpha\circ f\circ\alpha^{-1}.\] 
This reconstruction theorem completely determines
isomorphisms between homeomorphism groups of manifolds.
While this book is about smooth structures on manifolds, we remark that
Whittaker's theorem is also a trivial consequence of Rubin's Theorem, as stated below.

It turns out that the smooth structure
of a manifold is completely determined by the group $\Diff^r_c(M)_0$ up to $C^r$ diffeomorphism, at least when $M$ has no boundary.
The first result we will discuss in this direction is as below, which is the weaker of the two; a special case of this theorem is originally due to F.~Takens. 

\begin{thm}[Takens' Theorem]\label{thm:takens-intro}
If a set-theoretic bijection
 $w\colon M\longrightarrow N$
 between smooth connected boundaryless manifolds $M$ and $N$
  conjugates $\Diff_c^{p}(M)_0$ to $\Diff_c^{p}(N)_0$
  for some $p\in\bN\cup\{\infty\}$,
 then $w$ is a $C^p$ diffeomorphism.\end{thm}

Theorem~\ref{thm:takens-intro} is an example of a reconstruction theorem that takes as an input a relatively weak hypothesis (i.e.~a
set theoretic bijection inducing a certain group isomorphism) and promotes that bijection to a map with a high level of structure
(i.e.~a $C^r$ diffeomorphism). 

We will (nearly) replicate the complete 
original proof of Theorem~\ref{thm:takens-intro} in the case of one--dimensional
manifolds in Section~\ref{sec:takens} and $p=\infty$. The only part of the proof we will not spell out is a result of Sternberg on local linearization of
diffeomorphisms with hyperbolic fixed points. We will not reproduce Takens' argument in higher dimensions. 
One reason is to avoid introducing a large amount of additional background. 
Another is that Taken's theorem is subsumed by Filipkiewicz's Theorem, which itself can be further generalized as follows:

\begin{thm}[Takens--Filipkiewicz Theorem]\label{thm:filip-intro}
Let $M$ and $N$ be smooth connected boundaryless manifolds.
If there exists a group isomorphism
\[ \Phi\co \Diff_c^{p}(M)_0\longrightarrow \Diff_c^{q}(N)_0\] 
for some $p,q\in\bN\cup\{\infty\}$,
then 
we have that 
$p=q$
and that
 there exists a $C^p$ diffeomorphism
$w\colon M\longrightarrow N$ inducing $\Phi$ by conjugation.
\end{thm}

Whereas Takens' Theorem presumes the existence of a bijection between $M$ and $N$ that induces an isomorphism between the
corresponding diffeomorphism groups, Filipkiewicz starts with a very weak assumption (i.e.~the groups 
$\Diff_c^p(M)_0 $ and $\Diff_c^q(N)_0$ are
isomorphic) and ends with a strong conclusion (i.e.~$p=q$ and $M$ and $N$ are diffeomorphic via a $C^p$ diffeomorphism that
induces the isomorphism). As a consequence, Filipkiewicz's Theorem implies that there are no exotic automorphisms of $\Diff_c^p(M)_0 $
for a smooth manifold $M$; they are all induced by $C^p$ diffeomorphisms of $M$.

In our exposition of Filipkiewicz's result, we will retain a version of his argument for the sake of historical record, though in the main body
of the exposition, we will obtain Filipkiewicz's Theorem as a consequence of a simplified version of a much more general result due to Rubin.

\begin{thm}[Rubin's Theorem]\label{thm:rubin-intro}
If we have an isomorphism between two locally dense groups of homeomorphisms on perfect, locally compact, Hausdorff topological spaces,
then there exists a homeomorphism of those spaces that conjugates the isomorphism.
\end{thm}

In Theorem~\ref{thm:rubin-intro}, let $G$ be a group of homeomorphisms of a locally compact  Hausdorff  topological
space $X$ and let $U\sse  X$. We write $G[U]\le G$ for the subgroup consisting of homeomorphisms that are the identity outside of $U$.
We say $G$ is \emph{locally dense}\index{locally dense action} 
if for each point $x\in X$ and for each open neighborhood $U$ of $x$, 
the closure of the orbit $G[U].x$ has nonempty interior. We will give a self-contained proof of Rubin's Theorem, which the authors hope
the reader will find intellectually satisfying. 
Our presentation of Filipkiewicz's Theorem is based on Rubin's Theorem, which will imply that an isomorphism between
diffeomorphism groups of manifolds automatically induces a homeomorphism between them. To upgrade this homeomorphism to
diffeomorphism of the desired regularity, one uses a simultaneous continuity result for actions of $\bR^n$ on manifolds that is due to 
Bochner--Montgomery, and which itself is often
incorrectly attributed to Montgomery--Zippin. The reader will find a full account of the details in Section~\ref{sec:boch-mont}.

Returning to the main thread of discussion, let $M$ be a given manifold.
Even though for $p\neq q$ we have that $\Diff_c^p(M)_0 $ and 
$\Diff_c^q(M)_0$ are not isomorphic to each other as
algebraic groups, it is very difficult to imagine distinguishing between two such enormous infinite--dimensional continuous topological groups,
at least from an algebraic point of view. There are many different ways that one could hope to distinguish between two groups
$G$ and $H$, which is to say in an attempt
to find a checkable certificate that they are not isomorphic to each other. One such certificate might come from looking at the class of
subgroups of $G$ and $H$, and to find a subgroup of $G$ that does not occur as a subgroup of $H$, or vice versa. The simplest
subgroups to consider which have any hope of distinguishing $G$ and $H$ would be the finitely generated subgroups, and so we can
formulate a motivating question.

\begin{que}\label{que:subgroup-main}
Let $p< q$. Do the finitely generated subgroups of $\Diff^p_c(M)_0$ distinguish it from $\Diff^q_c(M)_0$? That is, is there a finitely generated group
of $C^p$ diffeomorphisms of $M$ that is not algebraically $C^q$ smoothable, i.e.~realized as a subgroup of $\Diff^q_c(M)_0$? More generally,
can we distinguish between $\Diff^p_c(M)_0$ and $\Diff^q_c(N)_0$ by finitely generated subgroups in order to conclude that either $p\neq q$ or
$p=q$ and $M$ is not $C^p$--diffeomorphic to $N$?
\end{que}

Returning to the results in this book that we have already announced, we restate Proposition~\ref{prop:unsmooth} in a way that
gives a satisfactory answer to Question~\ref{que:subgroup-main}, at least in the case of one--dimensional compact manifolds.
In the case of $M=I$ or $M=S^1$, every element of $\Homeo_+(M)$ is already isotopic to the identity, and compactness of isotopies
is automatic by the compactness of the ambient manifold, and so we can safely suppress the 
$0$--subscripts in the notation for diffeomorphism
groups. We will write $r=k+\eps$ and $s=\ell+\tau$, where $k,\ell\in\bZ_{>0}$ and where $\eps,\tau\in [0,1)$, and $M$ and $N$ will
both denote $I$ or $S^1$.

\begin{prop}\label{prop:unsmooth-refined}
For compact connected one--manifolds $M$ and $N$,
the following conclusions hold.
\begin{enumerate}[(1)]
\item
If $r<s$ then there is a finitely generated subgroup $G_r \le \Diff^r_+(M)$ such that no finite index subgroup $G_r$ is isomorphic to a subgroup
of $\Diff^s_+(M)$.
\item
For all $r$ and $s$ and $M\neq N$, there is a finitely generated subgroup $G_{r,M}\le \Diff_+^r(M)$ such that no finite index subgroup of
$G_{r,M}$ is isomorphic to a subgroup of $\Diff_+^s(N)$.
\end{enumerate}
\end{prop}

Indeed, the answer to Question~\ref{que:subgroup-main} is yes, in the case of compact one--manifolds. The reader may object to a claim
that the subgroups $G_r$ or $G_{r,M}$ are truly ``certificates", in the sense of some sort of easy checkability. Even though each group $G_r$
or $G_{r,M}$ is finitely generated, Proposition~\ref{prop:unsmooth-refined} implicitly furnishes an uncountable collection of pairwise
non--isomorphic subgroups, for example $\{G_r\}_{r\in\bR_{\geq 1}}$. There are only countably many isomorphism classes of finitely
presented groups, and in fact only countably many isomorphism classes of recursively presented groups (i.e.~finitely generated groups
whose relations are enumerable by a Turing machine) since there are only countably many different Turing machines. It follows that
the vast majority of the groups $\{G_r\}_{r\in\bR_{\geq 1}}$ are not even recursively presentable, 
which leads to a completely valid philosophical
question as to what it means to record such a group. The reader will find in the proof of Proposition~\ref{prop:unsmooth-refined}
(or, really, the proof of the most general results from which Proposition~\ref{prop:unsmooth-refined} follows), we will write down
generators for the groups $\{G_r\}_{r\in\bR_{\geq 1}}$, though again the claim that we are really ``writing down" diffeomorphisms is debatable,
since we will express them as certain limits which are not truly explicit. The proof that $G_r$ admits no homomorphism with
nonabelian image to $\Diff_+^s(M)$ for $s>r$ will inevitably have to take a detour through some infinitary methods.

Even the assertion that most of the groups $\{G_r\}_{r\in\bR_{\geq 1}}$ are not recursively 
presentable is non--constructive, since it says nothing
about $G_r$ for any specific value of $r$ and rather relies on the set theoretic properties of the indexing set. If one restricts to a countable
subset of the indexing set, we do not know if the corresponding groups are recursively presentable. We are left with the following question:

\begin{que}\label{que:fp}
For $M\in \{I,S^1\}$ and $k\in\bN$, is there a finitely presentable subgroup $G_k\le \Diff_+^k(M)$ that is not isomorphic to
a subgroup of $\Diff_+^{k+1}(M)$? What about a recursively presentable subgroup?
\end{que}

We direct the reader to the last paragraph of Appendix C for an answer when $k=0$.
Unfortunately, the methods given in this book do not give any insight into Question~\ref{que:fp}, at least as far as the authors can see.
The question of whether or not there exists an easily checkable certificate to distinguish two diffeomorphism groups remains, in this sense,
open.

\section{Critical regularity: history and overview}\label{sec:critreg-intro}

Much of the discussion in this book will be framed in terms of \emph{critical regularity}\index{critical regularity},
which we introduce and survey in this section.
The critical regularity of a group should be thought of as the sharp bound on the level of smoothness with which the group
can act on a smooth manifold, or more precisely, the upper limit on the algebraic smoothability of a group.

\subsection{Definitions and remarks}\label{ss:defn-rem}

To formulate a precise definition, let $G$ be a group, and let $M$ be a smooth manifold. The \emph{critical regularity}
of $G$ is defined to be \[\CR_M(G)=\sup\{r\in\bR_{\geq 0}\mid G\le \Diff_c^r(M)_0\}=\sup\{k+\eps\mid G\le \Diff_c^{k,\eps}(M)_0\}.\]
Here, $G$ is an abstract group, and so
$G\le \Diff^r(M)$ is shorthand for there existing an injective homomorphism $G\longrightarrow\Diff^r(M)$.
One could in principle extend the definition of critical regularity to the full group of $C^r$ diffeomorphisms of $M$, though in the cases
$M=I$ and $M=S^1$ that are the primary subject of this book, restricting the definition is the same as assuming that the diffeomorphisms
under consideration are orientation preserving.

The definition of critical regularity we have given here can be taken to be the more refined 
\emph{algebraic critical regularity}\index{algebraic critical regularity}, which
is then contrasted with the \emph{topological critical regularity}\index{topological critical regularity},
which we define now for later use. The topological critical regularity
refers not to a group but to a group with an action on $M$, and is then the supremum of regularities $r$ for which the action is topologically
conjugate into a group of $C^r$ diffeomorphisms.

From the definition of (algebraic) critical regularity, we have that:
\[\textrm{if}\quad G\le \Homeo_0(M)\quad \textrm{then}\quad \CR_M(G)\geq 0.\] If there is no injective homomorphism $G\longrightarrow
\Homeo_0(M)$ then we have \[\CR_M(G):=\inf\{\}=-\infty.\] When $G\le \Diff^{\infty}_c(M)_0$ we have $\CR_M(G)=\infty$. 
Thus, \emph{a priori}, the possible
range of $\CR_M$ is $\{-\infty\}\cup [0,\infty]$. 

For an arbitrary  given manifold $M$, it is not clear  that this range is fully populated; that is, it is often
very difficult to decide if for finite values of $r$ if there is a group $G$ such that $\CR_M(G)=r$, especially if
one insists that $G$ be a countable
or finitely generated group. 
Note from Filipkiewicz's Theorem that 
\[
\CR_M(\Diff^p_c(M)_0)=\CR_M(\Diff^p(M))=p\]
for all smooth connected boundaryless $n$--manifold $M$ and for all $p\in\bN$.
The set of actual values achieved by $\CR_M(G)$ for finitely generated groups $G$ is called the \emph{critical regularity spectrum}\index{critical
regularity spectrum} of $M$.

As the critical regularity of a group is a supremum, the statement that \[\CR_M(G)=r\] does not
necessarily mean that there is an injective homomorphism
$G\longrightarrow\Diff^r_c(M)_0$, but rather only that there is an injective homomorphism into $\Diff^s_c(M)_0$ for all $s<r$. Thus, the
finite part of the critical regularity
spectrum of $M$ bifurcates into the \emph{achieved spectrum}\index{achieved spectrum} and the 
\emph{unachieved spectrum}\index{unachieved spectrum}, where $r$ lies in the achieved
spectrum if there is in fact an injective homomorphism
$G\longrightarrow\Diff^r_c(M)_0$.

For the rest of this section, we will survey the state of knowledge about critical regularity for $M\in\{I,S^1\}$. The first examples of
groups of homeomorphisms of $M$ that one produces generally have infinite critical regularity. Indeed, a single homeomorphism
of $M$ generates a copy of $\bZ$ (or possibly $\bZ/n\bZ$ in the case $M=S^1$), and even if this homeomorphism is not continuously
differentiable, the group generated by it is algebraically smoothable. In the case of $M=I$, the group generated by the homeomorphism
is even topologically smoothable, though in the case of $M=S^1$, it may not be topologically smoothable.

From a flow on $M$, one can build free abelian groups of diffeomorphisms of arbitrarily large rank. Since one can build a smooth flow
from a smooth vector field on $M$, the critical regularity of a free abelian group is infinite. For free abelian groups of rank two or more,
one sees the difference between algebraic smoothability and topological smoothability of groups acting on $I$.
Whereas given a faithful action
by $\bZ^n$ for $n\geq 2$, one can always find another faithful action of this group that is as smooth as one likes, though one might not
be able to conjugate a given faithful action of $\bZ^n$ to a smoother one; see Kopell's Lemma (Theorem~\ref{thm:kopell}), for instance.

Nonabelian free groups also have infinite critical regularity, at least if they have at most countable rank. This fact is somewhat harder to
see than the infinite critical regularity of free abelian groups, and results essentially from a Baire Category type argument. We will
spell out some of the details in Corollary~\ref{cor:free-abundant}.

For groups that are not abelian or free, the computation of critical regularity quickly becomes difficult, as it can be hard to produce
faithful actions of a given group in the first place, and understanding the regularity properties of the group action usually involves very
subtle convergence questions and correspondingly subtle techniques. The definition of critical regularity implicitly quantifies over all
possible actions of a group on $M$, and so an additional layer of complication is added by the fact that it can be quite difficult to tell
if one has addressed all possible action of a group. To deal with this complication, one has to have good control over the orbit structure
for group actions on $M$, and this is where dynamical methods play an essential role.

\subsection{Cyclic groups, topological versus algebraic critical regularity}

We have already alluded to one of the earliest manifestations of the interplay between regularity, dynamics, and group actions, and
has to do with one of the simplest groups (i.e.~$\bZ$) acting freely on the circle $S^1$. Here, a \emph{free action}\index{free action}
is one where
no nontrivial group element fixes a point in $S^1$. The easiest examples of such actions come from the action of $S^1$ on itself, when
the manifold is viewed as an abelian group. Writing $S^1=\bR/\bZ$, an arbitrary element of $S^1$ acts additively on $S^1$, and the
subgroup generated by that element clearly acts freely. It is not difficult to see that an element $\theta\in S^1$ will have infinite order
as a homeomorphism of $S^1$ if and only if $\theta$ is not rational under the identification $S^1=\bR/\bZ$.
We write the corresponding element of $\Homeo_+(S^1)$ as $R_{\theta}$ and call it an \emph{irrational rotation}\index{irrational rotation}.

Clearly not every infinite order homeomorphism of $S^1$ acting freely is equal to an irrational rotation. Indeed, if $f\in\Homeo_+(S^1)$ is
arbitrary then $\form{ f^{-1}R_{\theta}f}\cong\bZ$ will be a free action of $\bZ$ on $S^1$ which will not be a rotation itself unless
$f$ commutes with $R_{\theta}$, which is quite rare among homeomorphisms that are not themselves rotations. Nevertheless, this is 
a somewhat ``trivial" perturbation of a rotation, since $f^{-1}R_{\theta}f$ differs from $R_{\theta}$ only via an identification of $S^1$ with
itself, which is to say $f^{-1}R_{\theta}f$ is \emph{topologically conjugate}\index{topological conjugacy} to $R_{\theta}$.

For an arbitrary element $g\in\Homeo_+(S^1)$, one can define an invariant called its \emph{rotation number}\index{rotation number},
and written $\rho(g)$.
This invariant is quite old, having been introduced by Poincar\'e, and has applications in geometric group theory and bounded cohomology
of groups as well as in dynamics~\cite{FrigerioBook,Ghys1987,Ghys2001,KKM2019}.
The rotation number takes values in $S^1=\bR/\bZ$, and it measures the tendency of $g$ to make points in $S^1$ wind about the circle.
We will avoid giving a precise definition of the rotation number here, as the definition is not so straightforward, and since
it will be discussed in great detail in Chapter~\ref{sec:denjoy}.

Several of the most important features of the rotation number are that it is invariant under topological conjugacy,
that it is nonzero precisely in the absence of a fixed point, that it is a homomorphism when restricted to the cyclic group generated
by a single homeomorphism (i.e.~it is \emph{homogeneous}\index{homogeneity}), and that it satisfies
$\rho(R_{\theta})=\theta$.
It is therefore natural to wonder, at least for irrational values of $\theta$, whether or not
the rotation number is a complete topological conjugacy invariant; that is, if $\rho(g)=\theta$, is it true that $g$ is topologically conjugate
to $R_{\theta}$?

The answer, perhaps surprisingly at first, is that it depends on the regularity of $g$. In broad strokes, we have the following:

\begin{thm}[Denjoy's theory of rotations, coarse summary]\label{thm:denjoy-intro}
If $\theta\in\bR/\bZ$ is irrational, then the following conclusions hold.
\begin{enumerate}[(1)]
\item
There exist continuous, and even $C^1$ diffeomorphisms $g$ of $S^1$ such that $\rho(g)=\theta$ but such that $g$ is not topologically
conjugate to $R_{\theta}$.
\item
If $g\in\Diff_+^2(S^1)$ and $\rho(g)=\theta$, then $g$ is topologically conjugate to $R_{\theta}$.
\end{enumerate}
\end{thm}

For a commonly satisified regularity assumption on $g$ that guarantees topological conjugacy to a rotation, one can relax
twice--differentiability to $\Diff_+^{1+\mathrm{bv}}(S^1)$, the group of $C^1$ diffeomorphisms of $S^1$ whose first derivatives have
bounded variation. Chapter~\ref{sec:denjoy} will give a self--contained and complete account of Theorem~\ref{thm:denjoy-intro}, and
will in the process will provide a full proof.

As we have mentioned before, there is a rather large gap between $\Diff_+^1(S^1)$ and
$\Diff_+^2(S^1)$, or even $\Diff_+^{1+\mathrm{bv}}(S^1)$, one that is manifested
not least by considering all possible diffeomorphisms whose derivatives satisfy
various H\"older continuity conditions. It is still an area of active research to determine the exact conditions on the derivatives of a
diffeomorphism of $S^1$ which guarantee topological conjugacy to a rotation, or conversely that allow for the construction of diffeomorphisms
with prescribed rotation numbers that are not conjugate to rotations (so--called \emph{Denjoy counterexamples}\index{Denjoy
counterexample}).

Before proceeding, we recall the notion of $\alpha$--continuity, which is a generalization of Lipschitz and H\"older continuity. Let
\[\alpha\colon\bR_{\geq0}\longrightarrow\bR_{\geq0}\] be a homeomorphism that is concave as a function. We say that a real valued function
$f$ is $\alpha$--continuous if there is a universal constant $C$ such that \[|f(x)-f(y)|\leq C\cdot \alpha(|x-y|)\] for all $x$ and $y$ in the domain
of $f$. Note that $\tau$--H\"older continuity is just $\alpha$--continuity for $\alpha(x)=x^{\tau}$. We call $\alpha$ a \emph{concave modulus of
continuity}\index{concave modulus of continuity}.

To give the reader an idea of some of what is known, we will give a full proof of the following result due to the authors in
Chapter~\ref{sec:denjoy}. In broad strokes, this theorem
shows that there are many Denjoy counterexamples that are very close to being $C^2$:

\begin{thm}\label{thm:int-denjoy}
Let $\alpha$ be a concave modulus of continuity and let $\theta\in\bR/\bZ$ be irrational.
If \[\int_{(0,1]}\frac{dx}{\alpha(x)}<\infty,\] then there exists a Denjoy counterexample $f\in\Diff_+^1(S^1)$ such that $f'$ is $\alpha$--continuous.
\end{thm}

\subsection{Foliations}

We now move beyond cyclic groups to describe the provenance and state of knowledge concerning critical regularity for finitely generated
groups, and more generally, countable groups.

One of the original motivations for studying regularity of group actions arises from foliation theory. A \emph{foliation}\index{foliation}
on a manifold $M$
of dimension $n$ is an equivalence relation on $M$, where the equivalence classes are given by $p$--dimensional immersed submanifolds,
for $p\leq n$. Thus, a foliated manifold locally looks like a direct product decomposition $\bR^n\cong\bR^p\times \bR^q$, where $q=n-p$,
and the immersed $p$--dimensional manifolds are called \emph{leaves}\index{leaf of a foliation} of the foliation.
The numbers $p$ and $q$ are the \emph{dimension}\index{dimension of a foliation} and 
\emph{codimension}\index{codimension of a foliation} of the foliation, respectively.
Imagine the French \emph{mille
feuille}\index{mille feuille} (thousand leaves, or thousand sheets) pastry.

The data describing a foliation is a coherent atlas, much like one defining a manifold, which keeps track of the local product decomposition,
and which satisfies suitable compatibility conditions. We will not recall the precise definition of a foliation here, postponing it
instead to Section~\ref{sec:foliation}.

Foliation theory originates in differential topology, with its foundational result perhaps being
Ehresmann's Submersion Theorem~\cite{Ehresmann51,DundasBook},
which was seminal not only in the theory of foliations but also in the study of fibrations and fiber bundles. Foliation theory developed into
a vast and deep subject, with connections to relativity, differential geometry, $C^*$--algebras, noncommutative geometry, and dynamical
systems. We will not attempt to summarize this story here, and we will only note the facts we need
when they are required in Chapter~\ref{ch:app}. We refer the reader to~\cite{CandelConlonI,CandelConlonII} for an encyclopedic reference on
foliation theory.

Even when $M$ is a smooth manifold, the data of the foliation may have a degree of regularity which is much lower than that of $M$.
Indeed, even when $M$ is a smooth manifold, one can find foliations on $M$ whose transition functions (i.e.~the ones encoding the local
product decomposition) are not differentiable, and such that there is no homeomorphism (of foliated manifolds) to another manifold
wherein the target foliation is differentiable.

The smoothness of a foliation is closely related to a certain pseudo-group associated to the foliation, called its 
\emph{holonomy}\index{holonomy}. In
Section~\ref{sec:foliation}, we will describe a particular type of construction of foliations, called 
\emph{suspension of a group action}\index{suspension of a group action},
which realizes many groups as holonomy groups of foliations. The suspension of a group action is essentially the same as the construction
of a flat bundle $E$ with fiber $F$ and base space $M$,
whose monodromy group lies in $\Diff^r(M)$ and some suitable regularity $r$, and results in a
foliated bundle with exactly that regularity. More precisely, the data of the foliation is specified by a representation \[\phi\colon\pi_1(F)
\longrightarrow \Diff^r(M),\] and the homeomorphism class of the foliated bundle is completely determined by the topological conjugacy class
of $\phi$. In particular, the construction of suspended group actions that are unsmoothable beyond some regularity $r$ is equivalent to
constructing representations $\phi$ into diffeomorphism groups that are not topologically smoothable. This is one specific reason that
regularity of group actions is of interest in foliation theory.

Since topological smoothability is a stronger condition than algebraic smoothability, it is perhaps not surprising that it is somewhat easier
to find examples of groups with a given topological critical regularity (i.e.~a particular action of the group cannot be topologically smoothed
to a higher degree of regularity) than it is to find examples of groups with a given algebraic critical regularity (i.e.~no action in higher degree
of regularity exists). In Section~\ref{sec:foliation}, we will discuss in detail a construction due to Tsuboi of such a group with given topological
critical regularity; in all degrees, the group is fixed, and so the algebraic critical regularity of the abstract group in question is infinite,
but the particular actions Tsuboi constructs are not topologically conjugate to smoother ones. Another construction in a very similar
vein was also found by Cantwell and Conlon. We will delay a discussion of the details until Chapter~\ref{ch:app}.

Obtaining control on the topological critical regularity of a group action already has consequences for foliation theory, as we see, though
one can ask if a certain abstract group occurs as the holonomy of a codimension one foliation, for example. When it comes to
analyzing the smoothness of a foliation with such a holonomy group, the critical regularity of the group plays an essential role.

\subsection{Nilpotent groups}

We have seen that free abelian groups may have finite topological critical regularity, even in rank one, but their (algebraic) critical
regularity as we are discussing it is infinite. The next level of algebraic complication for groups is nilpotence. It is not
difficult to show that finitely generated torsion--free nilpotent groups are \emph{orderable}\index{orderable group},
which is to say that they admit total orderings
that are compatible with the group structure. We include a quick proof of this fact for the convenience of the reader, and it will be clear
that the method of proof generalizes to a wide class of groups.

\begin{prop}
Let $N$ be a finitely generated torsion--free nilpotent group. Then $N$ admits a left invariant ordering, which is to say a total ordering $\leq$
such that for all $a,b,c\in N$, we have $a\leq b$ if and only if $ca\leq cb$.
\end{prop}
\begin{proof}
We proceed by induction on the \emph{Hirsch length}\index{Hirsch length}
of $N$, or the \emph{polycyclic length}\index{polycyclic length} of $N$. That is, we look at the longest subnormal
chain of subgroups \[N=N_0>N_1>\cdots>N_k=\{1\},\] such that $N_i$ is normal in $N_{i-1}$ for all suitable indices, and so that
$N_{i-1}/N_i\cong\bZ$. The length of this series is defined to be $k$.
It is a straightforward exercise to show that such a sequence exists and that it terminates.

Clearly, we have that the group $\bZ$ admits a left invariant ordering, as we may consider the inherited ordering from $\bR$. Next, if
\[1\longrightarrow K\longrightarrow G\longrightarrow Q\longrightarrow 1\] is a short exact sequence of groups such that $K$ and $Q$
are both left orderable, then so is $G$. Indeed, we let $g_1,g_2\in G$. We we say $g_1\leq g_2$ if the images of these elements in $Q$
satisfy $\overline{g_1}\leq \overline{g_2}$ in the ordering on $Q$, or if $g_1^{-1}g_2\in K$ and $1\leq g_1^{-1}g_2$ in the ordering on $K$.
The reader may check that this is a left invariant ordering on $G$. It is clear now that $N$ admits a left invariant ordering.
\end{proof}

Among other things, we see that if $N$ is finitely generated, torsion--free, and nilpotent,
then $N\le \Homeo_+(M)$ for $M\in\{I,S^1\}$. The above result extends to solvable groups as well.
We direct the reader to 
Appendix~\ref{ch:append2} for more detail on orderability, and in particular
a proof of the fact that a countable group admitting a left invariant ordering admits an injective homomorphism into $\Homeo_+[0,1]$.

While actions by homeomorphisms arising from orderability are usually not differentiable, it is not clear \emph{a priori} that nilpotent groups
do not have infinite critical regularity. It turns out that their critical regularity is finite, by a
1976 result of Plante--Thurston~\cite{PT1976}.

\begin{thm}\label{thm:pt-intro}
If $M\in\{I,S^1\}$ then every nilpotent subgroup of  $\Diff_+^{1+\mathrm{bv}}(M)$ is abelian.
\end{thm}

As an immediate consequence, using the fact that $\Diff^2_+(M)\le \Diff_+^{1+\mathrm{bv}}(M)$, we have:

\begin{cor}
If $N$ is a nonabelian, torsion--free nilpotent group, then $\CR(N)\leq 2$. If $\CR(N)=2$ then the critical regularity is not achieved.
\end{cor}

We will give a full proof of Theorem~\ref{thm:pt-intro} in Chapter~\ref{ch:c2-thry}. If $N$ is a finitely generated torsion--free nilpotent group
that is not abelian, we have now that its critical regularity is finite. The problem of determining the critical regularity precisely is very
challenging, and only partial results are known.

The first question to ask is whether or not a finitely generated, nonabelian, torsion--free nilpotent group $N$ can be realized as a group
of diffeomorphisms of $M$. The answer is yes, though the proof is far from obvious.

\begin{thm}\label{thm:jor-farb-franks}
Every finitely generated, torsion--free nilpotent group embeds into $\Diff_+^1(M)$.
\end{thm}

It is remarkable that Theorem~\ref{thm:jor-farb-franks} appeared much later than Theorem~\ref{thm:pt-intro}, appearing first in~\cite{FF2003}.
Related constructions were given by Jorquera~\cite{Jorquera}. It is not very difficult to produce explicit, faithful actions of a
finitely generated, torsion--free
nilpotent group via appeals to homeomorphisms that take certain intervals to certain other intervals, and it is in gluing these homeomorphisms
together in a way that is smooth which necessitates most of the work in proofs of Theorem~\ref{thm:jor-farb-franks}. To illustrate these ideas
to the reader, we will give the constructions
carried out in Farb--Franks' and Jorquera's papers in Section~\ref{sec:nilpotent}, and partially write down the diffeomorphisms
they use. We will explain the steps needed to complete the proof, but we will not reproduce the calculations.

Whereas Farb--Franks and Jorquera produce an action of a given nilpotent group (or in Jorquera's case, an action of a certain direct limit
of nilpotent groups), a remarkable result of Parkhe gives a general orbit structure theory for 
\emph{all} nilpotent group actions on $I$ and $S^1$.
We will not reproduce the statement of Parkhe's structure theorem here, which the reader may find as Theorem~\ref{thm:parkhe-main}.

As a consequence, Parkhe is able to deduce a lower bound on the critical regularity of nilpotent groups that is uniform in the growth
rate. Here, the \emph{growth rate}\index{growth of groups} 
of a finitely generated group measures the number of elements in the group that can be written as products
of at most a fixed number $n$ of generators. To define the growth rate precisely, let $S$ be a finite generating set for a group $G$ such
that $S=S^{-1}$ (i.e.~$S$ is closed under taking inverses). We write $b_S(n)$ for the number of distinct elements of $G$ that can be
written as a product of at most $n$ elements of $S$. Whereas $b_S$ depends on $S$, its asymptotic behavior does not. We say
that $G$ has \emph{polynomial growth}\index{polynomial growth} of degree at most $d$ if $b_S(n)=O(n^d)$.
Here, the notation $O(n^d)$ is the usual big-oh notation, meaning that the growth rate of the group
as a function of $n$ is bounded above by $C\cdot n^d$, for a suitable constant $C$.

Clearly, an infinite group has at least polynomial growth, and no group
can have growth rate faster than an exponential function of $n$. A group $G$ has \emph{exponential growth}\index{exponential growth}
if the growth rate is bounded below by an exponential function, and has \emph{intermediate growth}\index{intermediate growth}
if it has superpolynomial and
subexponential growth. We refer the reader to~\cite{dlHarpe2000} for more background.

A foundational result of Gromov characterizes finitely
generated groups that have a nilpotent subgroup of finite index as precisely those that have polynomial growth;
see Subsection~\ref{ss:parkhe} for more detail and references.
Parkhe proves the following,
which will occur also as Theorem~\ref{thm:parkhe-conj} below:

\begin{thm}\label{thm:parkhe-intro}
Let $N\le \Homeo_+(M)$ be a finitely generated group of polynomial growth $O(n^d)$. Then for all $\tau<1/d$, the group $N$ is topologically
conjugate into $\Diff_+^{1,\tau}(M)$.
\end{thm}

Theorem~\ref{thm:parkhe-intro} has several
interesting consequences. For one,
in conjunction with Theorem~\ref{thm:pt-intro}, it says that the critical regularity of a finitely generated, torsion--free, nonabelian nilpotent groups
always lies in the interval $(1,2]$. Moreover, it says that if $N$ acts on $M$ by homeomorphisms, then this action is secretly differentiable;
a failure to be $C^1$ is simply concealed in a poor choice of identification of $M$ with itself. Finally, it says that a foliated bundle given by
suspending a nilpotent group action on a compact one--manifold is always smoothable to a $C^{1,\tau}$ action for a suitable $\tau$ that
is bounded below by a coarse geometric invariant of the group. In the interest of space, we will not give a complete proof of Parkhe's results.

Whereas Parkhe's work gives uniform lower bounds on critical regularity for nilpotent groups that complement the Plante--Thurston Theorem,
the bounds produced by Theorem~\ref{thm:parkhe-intro} are generally not optimal; one need only consider the case of abelian groups.
Even for nilpotent groups, the precise determination of the critical regularity of nilpotent groups is very difficult and known only in a few
cases. To state them, let $N_m\le \mathrm{SL}_{m+1}(\bZ)$ denote the group of unipotent upper triangular matrices. That is, $N_m$ consists of
all integer $(m+1)\times (m+1)$ matrices which are upper triangular, and whose only diagonal entry value is $1$. It is a standard fact that $N_m$ is
nilpotent, and that every finitely generated torsion--free nilpotent group can be realized as a subgroup of $N_m$ for some sufficiently large $m$.
The group $N_1$ is isomorphic to $\bZ$, and the group $N_2$ is isomorphic to the integral Heisenberg group \[\mathrm{Heis}=\form{ x,y,z\mid [x,y]=z,\,
[z,x]=[z,y]=1}.\] The following result is due to Jorquera--Navas--Rivas (cf.~Theorem~\ref{thm:jnr} below):

\begin{thm}\label{jnr-intro}
We have $\CR_I(N_3)=1.5$.
\end{thm}

It is unknown whether the critical regularity of $N_3$ is achieved. The following result is due to Castro--Jorquera--Navas
(cf.~Theorem~\ref{thm:cjn} below):

\begin{thm}\label{cjn-intro}
For all $d\geq 2$, there exists a metabelian group of nilpotence degree $d$ whose critical regularity is equal to $2$. The critical regularity
for these groups is not achieved. In particular, the integral Heisenberg group satisfies $\CR_M(H)=2$ for $M\in\{I,S^1\}$.
\end{thm}

Here, a group is metabelian if its commutator subgroup is abelian. The nilpotence degree of a group is the length of its lower central series.
We will not prove Theorems~\ref{jnr-intro} nor~\ref{cjn-intro}, mostly for reasons of space.
The results mentioned here document the current state of knowledge concerning critical regularity for nilpotent groups.

\subsection{Right-angled Artin groups and mapping class groups}

There are very few classes of groups for which critical regularity is understood to the degree that it is for nilpotent groups. One class
for which a significant amount of progress has been made is \emph{right-angled Artin groups}\index{right-angled Artin group}.
A right-angled Artin group is determined by
a finite simplicial graph $\Gamma$ with vertex set $V$ and edge set $E$. The right-angled Artin group \[A(\gam)=\form{ V(\gam)\mid
[v,w]=1\textrm{
if and only if } \{v,w\}\in E(\gam)}.\] A reader unfamiliar with right-angled Artin groups may check that the class accommodates free
abelian groups, nonabelian free groups, and many groups in between. It turns out that right-angled Artin groups exhibit a combination
of diversity and uniformity of behavior. The reader may consult~\cite{Charney2007,KK2013,KK2013b,Koberda-Yale,Wise2012} for instance,
for background on the subject for which we will not have space to discuss
in detail. 

Since the class of right-angled Artin groups contains free abelian groups and free groups, clearly some right-angled Artin groups have infinite
critical regularity. It is not very difficult to show that free products of free abelian groups, and direct products of free products of free abelian
groups are all right-angled Artin groups with infinite critical regularity (see Section~\ref{sec:raag} below). The underlying graphs of these
right-angled Artin groups are disjoint unions of complete graphs (free products of free abelian groups) and joins of these kinds of graphs
(direct products of free products of free abelian groups). In general right-angled Artin groups are residually
torsion--free nilpotent~\cite{DK1992a},
whence it will follow from Theorem~\ref{thm:ff2003} below that they admit faithful $C^1$ actions on $M$:

\begin{prop}\label{prop:raag-c1-intro}
Every finitely generated right-angled Artin group embeds into  $\Diff_+^{1}(M)$.
\end{prop}

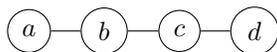
\begin{figure}[h!]
\tikzstyle {bv}=[black,draw,shape=circle,fill=black,inner sep=1pt]
  \begin{center}
\begin{tikzpicture}[main/.style = {draw, circle}]
\node[main] (1) {$a$};
\node[main] (2) [right of=1] {$b$};
\node[main] (3) [right of=2] {$c$}; 
\node[main] (4) [right of=3] {$d$};

\draw (1)--(2)--(3)--(4);
\end{tikzpicture}%
\caption{The graph $P_4$, one of the simplest defining graphs of a right-angled Artin group with finite critical regularity.}
\label{f:p4-intro}
\end{center}
\end{figure}

For complicated defining graphs, it is somewhat harder, or more directly, impossible, to construct smooth actions of right-angled Artin
groups. Historically, the first graph shown to be an obstruction to infinite critical regularity for right-angled Artin groups was the path $P_4$
of length three~\cite{BKK2019JEMS}. See Figure~\ref{f:p4-intro}.

Note that \[A(P_4)=\form{ a,b,c,d\mid [a,b]=[b,c]=[c,d]=1}.\] A reader familiar with three--manifold topology will observe that
$A(P_4)$ is the fundamental group of a link complement in $S^3$, precisely the link consisting of a chain of four unknots, successive
ones linked with linking number one. We have the following result, which will appear as Theorem~\ref{thm:a4-c2} below:

\begin{thm}\label{thm:bkk-intro}
There is no injective homomorphism \[A(P_4)\longrightarrow\Diff_+^{1+\mathrm{bv}}(M).\] In particular, $\CR_M(A(P_4))\leq 2$.
\end{thm}

Graphs which contain no copy of $P_4$ as a full subgraph are called \emph{cographs}\index{cograph}, and right-angled Artin groups on
cographs are characterized by not containing copies of $A(P_4)$ as subgroups, as we will show in Section~\ref{sec:raag}.
Cographs are organized in a hierarchy called the \emph{cograph hierarchy}\index{cograph hierarchy}.
Roughly, all cographs are generated recursively from a singleton
vertex via the graph operations of disjoint union and join; dually, right-angled Artin groups on cographs are generated recursively from $\bZ$
via the group operations of direct product and free product.

It is not true that a right-angled Artin group on a cograph has infinite critical regularity; in fact, it turns out that the group $(F_2\times\bZ)*\bZ$
is an obstruction to infinite critical regularity, and its presence characterizes right-angled Artin groups with finite critical regularity.

\begin{thm}\label{thm:raag-class-intro}
Let $\Gamma$ be a finite simplicial graph. Then the following are equivalent.
\begin{enumerate}[(1)]
\item
$\CR_M(A(\gam))<\infty$.
\item
$\CR_M(A(\gam))\leq 2$.
\item
We have $(F_2\times\bZ)*\bZ\le A(\gam)$.
\item
We have that $A(\gam)$ is not a direct product of free products of free abelian groups.
\end{enumerate}
\end{thm}

Theorem~\ref{thm:raag-class-intro} is strictly stronger than Theorem~\ref{thm:bkk-intro} and will be proved fully in this book.
It follows essentially from Theorem~\ref{thm:f2-int} and Theorem~\ref{thm:f2-circ}, together with the
remainder of the discussion in Section~\ref{sec:raag}. The reader may also consult Theorem~\ref{thm:kharlamov} for a more refined
statement.

The precise value of the critical regularity of right-angled Artin groups is unknown in general. It seems like as good a guess as any other
that the critical regularity should depend on the combinatorics of the defining graph, though there is strong evidence that when the
critical regularity of a right-angled Artin group is finite, then it should be exactly equal to $1$. The following will appear as
Theorem~\ref{thm:kkr-2020} below:

\begin{thm}\label{kkr-intro}
For a compact connected one--manifold $M$,  we have $\CR_M((F_2\times F_2)*\bZ)=1$.
\end{thm}

If $\Gamma$ is a graph that contains a square as a full subgraph and if its complement graph is connected,
then Theorem~\ref{kkr-intro} implies that
$A(\gam)$ also has critical regularity one. Graphs $\Gamma$ for which $(F_2\times F_2)*\bZ$ is not a subgroup of $A(\gam)$ resist
an easy characterization, unlike cographs and graphs for which $A(\gam)$ contains $(F_2\times\bZ)*\bZ$. The critical regularity
of $(F_2\times\bZ)*\bZ$ remains tantalizingly open, though it would be not be inappropriate to hope for the intellectually satisfying
answer that it should also be exactly one, in which case right-angled Artin group critical regularities would bifurcate neatly.

\begin{que}
Let $\Gamma$ be a finite simplicial graph. Is it true that \[\CR_M(A(\gam))\in\{1,\infty\}?\]
\end{que}

Critical regularity questions for mapping class groups of surfaces are also understood to some degree, mostly because of their
entanglement with right-angled Artin groups. Let $S$ be an orientable surface with genus $g$ and $n$ punctures,
marked points, or boundary components.
We let $\Mod(S)$ denote the \emph{mapping class group}\index{mapping class group}
of $S$, consisting of isotopy classes of orientation preserving homeomorphisms
of $S$. A comprehensive overview of mapping class groups and their properties can be found in~\cite{FM2012}.
Classical results of Nielsen, Thurston, and Handel show that if $S$ has a marked point, then $\Mod(S)$ acts faithfully on the circle,
and if $S$ has a boundary component then $\Mod(S)$ acts faithfully on the interval. These facts are discussed as Theorem~\ref{thm:nielsen}
and Theorem~\ref{thm:thurston-handel} below.

In~\cite{FF2001}, Farb and Franks proved that if the genus of $S$ is sufficiently large, then $\CR_M(\Mod(S))\leq 2$. For reasons that
we will explain shortly, it is natural to replace the full mapping class group with groups that are 
\emph{commensurable}\index{commensurability}
with mapping class groups
(i.e.~groups with finite index subgroups that are isomorphic to a finite index subgroup of the mapping class group).

Since their introduction by Dehn~\cite{Dehn-collected,Birman-MCG},
many authors have observed that mapping class groups of surfaces share many properties with
lattices in semisimple Lie groups. Here as before, a \emph{lattice}\index{lattice}
$\Gamma$ in a Lie group $G$ is a discrete subgroup such that $G/\gam$ has
finite volume with respect to the Haar measure on $G$.
It is known that the mapping class group itself is not a lattice in a semisimple Lie group, and it shares
some properties with rank one lattices, and some properties with higher rank lattices (see~\cite{FLM01}, for example).

Actions of lattices on compact manifolds are the subject of the \emph{Zimmer Program}\index{Zimmer Program},
which asserts that ``large" groups cannot act on
``small" manifolds in interesting ways. Here, a ``small" manifold is usually a compact manifold of dimension $n$, and a ``large" group
is an irreducible lattice in a semisimple Lie group whose rank is sufficiently large in comparison with $n$. Recent years have seen huge
progress on this program, and we refer the reader to~\cite{BM1999,Witte1994,Ghys1999,BFH20,BFH16} for further information.

For actions on the circle and interval, ``large" has usually meant an irreducible lattice of rank at least two. One of the first results regarding
lattices and actions on compact one--manifolds was obtain by Witte Morris:

\begin{thm}[See~\cite{Witte1994}]\label{thm:wm-intro}
Let $n\geq 3$ and let $\gam\le \mathrm{SL}_n(\bZ)$ be a subgroup of finite index. 
Then every homomorphism
\[\gam\longrightarrow\Homeo_+(M)\] has finite image, for $M\in\{I,S^1\}$,
\end{thm}

For more general lattices, we have the following result that was obtained by Ghys~\cite{Ghys1999} and Burger--Monod~\cite{BM1999}:

\begin{thm}\label{thm:bm-intro}
For an irreducible lattice $\Gamma$ in a semisimple Lie group of rank at least two, we have the following.
\begin{enumerate}[(1)]
\item
Every $C^0$ action of $\Gamma$ on $S^1$ has a finite orbit.
\item
Every $C^1$ action of $\Gamma$ on $S^1$ factors through a finite group.
\end{enumerate}
\end{thm}

Theorem~\ref{thm:bm-intro} implies that there are no interesting differentiable actions of higher rank lattices on the circle, and that
after possibly passing to a finite index sublattice, every continuous action of a higher rank lattice on the circle is actually just an action
on the interval. As we have suggested already, and as the reader will find from consulting Appendix~\ref{ch:append2}, the existence
of nontrivial actions of higher rank lattices on the circle thus reduces to finding left invariant orderings on lattices. This is generally
a very difficult problem; recent progress has been made by Deroin--Hurtado~\cite{DeroinHurt}.

Lattices in rank one Lie groups, however, often do admit highly regular faithful actions on compact one manifolds. A lattice in
$\mathrm{PSL}_2(\bR)$, for instance, admits a faithful analytic action on the circle, since the group $\mathrm{PSL}_2(\bR)$
is itself a group of analytic diffeomorphisms of the circle. All hyperbolic $3$--orbifold groups, which are precisely the lattices
in $\mathrm{PSL}_2(\bC)$, admit faithful $C^1$ actions on the circle and on the interval. This is a consequence of the work
of Agol, Wise, and Kahn--Markovic~\cite{Agol2008,Agol2013,Wise2011,Wise2012,KM2012}
that shows that these groups embed in right-angled Artin groups (possibly after
passing to a finite index subgroup), combined with Theorem~\ref{thm:ff2003}.

Since a finite index subgroup of a lattice is again a lattice, any analogy between mapping class groups and lattices in semisimple Lie groups
should only consider the \emph{commensurability class}\index{commensurability class} of the mapping class group, which is to say of the
mapping class group up to passing to a finite index
subgroup. The difficulty here is that finite index subgroups of the mapping class group are very opaque.
See Section~\ref{sec:mcg} for a more detailed discussion.

Right-angled Artin groups furnish a tool to study finite index subgroups of mapping class groups because of the following straightforward
observation: if $\Gamma$ is a finite simplicial graph, if $G$ is a group, and if $A(\gam)\le G$, then every finite index subgroup of $G$ contains
an isomorphic copy of $A(\gam)$. Moreover, most mapping class groups contain copies of $A(P_4)$ or $(F_2\times \bZ)*\bZ$. This
can be seen easily from a result of the second author that exhibits a profusion of right-angled Artin subgroups of mapping class groups.
See Theorem~\ref{thm:mcg-raag} below. As a consequence, even though mapping class groups often do admit faithful continuous
actions on compact one--manifolds, these actions are never $C^{1+\mathrm{bv}}$ or better. The following result will be discussed and
contextualized further in Section~\ref{sec:mcg}.

\begin{thm}\label{thm:bkk-mcg-intro}
Let $S$ be an orientable surface and let $M\in\{I,S^1\}$. Then one of the two mutually exclusive conclusions holds.
\begin{enumerate}[(1)]
\item
There is a finite index subgroup $G\le \Mod(S)$ such that $G$ is a product of a free group and an abelian group. In this case,
$\CR_M(G)=\infty$.
\item
For all finite index subgroups $G\le \Mod(S)$, we have $\CR_M(G)\leq 2$.
\end{enumerate}
\end{thm}

Thus, mapping class groups again exhibit behavior of a lattice of intermediate rank, that is something between rank one and higher
rank. They often admit faithful continuous actions on compact one--manifolds,
but not very smooth ones.

For most mapping class groups, one can improve the bound in Theorem~\ref{thm:bkk-mcg-intro} to $\CR_M(G)\leq 1$. We will postpone
further discussion until Section~\ref{sec:mcg}.

\subsection{Property (T) and intermediate growth}\label{sec:ggt}

Most of the obstructions to smooth actions on one--manifolds we have considered up to this point arise from
commutation, and how partially commutative phenomena interact with differentiability.
There are some deep properties of groups which do not have much to do with partial commutativity but which nevertheless seem
to be incompatible with smooth actions on one--manifolds. We mention relevant results here for completeness of the overview of the
theory, though we will not discuss them in further detail in the body of the book.

\subsubsection{Kazhdan's property (T)}
Property (T), introduced by Kazhdan~\cite{Kazhdan64},
is a somewhat abstract property on the unitary dual of a group. Qualitatively, groups with property (T)
are ``large" groups, and are hence expected not to admit interesting actions on ``small" manifolds, according to the principles guiding
the Zimmer program. Indeed, property (T) is a property enjoyed by irreducible lattices in higher rank semisimple Lie groups, and was
used to prove many facts about them.

To give a rough definition, the set of unitary representations of a group has a natural topology on it, called the Fell topology. A group has
property (T) if the trivial representation is isolated in this topology, which gives an etymology of the name of the property. There are many
other equivalent characterizations of property (T). A more intuitive one, in the authors' opinion, is that $G$ has property (T) if and only if
any isometric affine action of $G$ on a Hilbert space has a fixed point. There is yet another characterization of property (T) in terms of
cohomology; we avoid discussing this subject at length here, directing the reader instead to the books~\cite{bekka-valette,Zimmer84}.

Given the conjectural picture that lattices in higher rank should admit no interesting actions on compact one--manifolds, it is natural
to ask whether groups with the more abstract Kazhdan's property (T) can admit such actions. Progress in this direction is made
by Navas~\cite{Navas02}:

\begin{thm}\label{thm:navas-kazhdan}
Let $G$ be a countable group with property (T), and let $\tau>1/2$. If \[\phi\colon G\longrightarrow\Diff_+^{1,\tau}(S^1)\] is a homomorphism
then $\phi$ has finite image. In particular, if $G$ is an infinite countable group with property (T) then $\CR_{S^1}(G)\leq 1.5$.
\end{thm}

Bader--Furman--Gelander--Monod~\cite{BFGM2007AM} show that if $G$ is as in Theorem~\ref{thm:navas-kazhdan}, then
$G$ is not a subgroup of $\Diff_+^{1,1/2}(S^1)$. It remains an open question as to whether a group with Kazhdan's property (T) can
admit an infinite image action by homeomorphisms on the circle or on the interval.

\subsubsection{Intermediate growth}

Recall that if $G$ is finitely generated by a symmetric generating set $S$ then we can define the growth function $b_S(n)$, measuring
the number of elements of $G$ that can be written as a product of at most $n$ elements of $S$. All the explicit groups we have considered
up to this point are either of polynomial growth (i.e.~virtually nilpotent groups) or of exponential growth (everything else). Since both
these classes of groups contain groups with infinite critical regularity and with finite critical regularity, it is hard to suppose that the
growth rate of a group should be related to critical regularity in any way.

It turns out that groups of intermediate growth, however, always have finite critical regularity when they act on the interval or on the circle.
Groups of intermediate growth are not trivial to construct in the first place, the first of which being constructed by Grigorchuk
(see~\cite{dlHarpe2000}). Groups of intermediate growth admitting faithful actions on the interval were produced by
Grigorchuk--Machi~\cite{GM1993}. Navas~\cite{Navas2008GAFA} proved the following result:

\begin{thm}\label{thm:navas-intermediate}
Let $\tau>0$, let $M\in\{I,S^1\}$, and let $G$ be a finitely generated subgroup of $\Diff_+^{1,\tau}(M)$ with subexponential growth. Then $G$ is
virtually nilpotent. In particular, a finitely generated group $G$ of intermediate growth satisfies $\CR_M(G)\leq 1$.
\end{thm}

Navas also proved that $\Diff_+^1(M)$ contains subgroups of intermediate growth, and so the critical regularity bound
given by Theorem~\ref{thm:navas-intermediate} is achieved by some
groups of intermediate growth. The essential connection between growth of groups and one--dimensional dynamics is through
free semigroups. If $G$ is a group of subexponential growth then $G$ cannot contain a free semigroup on two generators.
Navas is then able to exploit the absence of the free semigroup in $G$ to exclude 
\emph{crossed homeomorphisms}\index{crossed homeomorphisms} from
any action of $G$ on $M$. Manifestations of crossed homeomorphisms will arise in Chapter~\ref{ch:c2-thry}, and their absence
is one characterization of Conradian group actions. We will not discuss this part of the story any further, since Navas already
gives an account in~\cite{Navas2011}.

\subsection{Lipschitz lower bounds}

In the preceding sections, we have glossed over groups of critical regularity between zero and one. All the explicit examples we have
discussed have had critical regularity at least one. There is a good reason for this, which arises from a result of
Deroin--Kleptsyn--Navas~\cite{DKN2007}:

\begin{thm}\label{thm:dkn-intro}
If $M\in\{I,S^1\}$, then 
every countable subgroup of $\Homeo_+(M)$ is topologically conjugate into a group of bi--Lipschitz
homeomorphisms of $M$.
\end{thm}

Since the Lipschitz modulus of continuity is stronger than every H\"older modulus of continuity, it is reasonable to say in light of
Theorem~\ref{thm:dkn-intro} that:
\begin{cor}\label{cor:lip-bound}
Let $G$ be a countable group. If \[\CR_M(G)\geq 0\quad \textrm{then}\quad \CR_M(G)\geq 1.\]
\end{cor}

Corollary~\ref{cor:lip-bound} is not an unreasonable statement, especially
since bi--Lipschitz homeomorphisms are $C^1$ almost everywhere. However, it leads to the antinomy that groups of critical regularity
one might still admit no differentiable actions on $M$. Indeed, such groups are discussed in Appendix~\ref{ch:append3}.

From the dynamical point of view that we develop in this book, the quantitative behavior of a bi--Lipschitz homeomorphism is no
different from that of a diffeomorphism. More generally, the quantitative behavior of a $C^k$ diffeomorphism whose $k^{th}$ derivative
is Lipschitz is no different from that of a $C^{k+1}$ diffeomorphism. In the authors' paper~\cite{KK2020crit}, it was unclear whether
the classes of finitely generated groups of integer critical regularity $k\geq 2$ were truly bifurcated into $(k,\mathrm{Lip})$ and $(k+1)$
diffeomorphisms; that is, the methods used could not produce a group of $C^k$ diffeomorphisms with Lipschitz $k^{th}$ derivatives that
were not $C^{k+1}$ diffeomorphisms already. We will comment on this further in Subsection~\ref{ss:prescribed} below.

Theorem~~\ref{thm:dkn-intro} implies that for one--manifolds, the spectrum of critical regularities of countable groups is contained in
$\{-\infty\}\cup [1,\infty]$. Another technical complication that Theorem~\ref{thm:dkn-intro} avoids is that the set of homeomorphisms
of $M$ that are $\tau$--H\"older continuous do not form a group;
for $\tau,\tau'\in[0,1)$, the composition of a $C^{\tau}$ and a $C^{\tau'}$ homeomorphism
need only be $C^{\tau\tau'}$. Of course, groups of such homeomorphisms exist, but a foundational
issue arises in analyzing compositions of such homeomorphisms. It is a nontrivial fact then that for $\tau\in [0,1)$
and for integers $k\geq 1$, the set $\Diff_+^{k,\tau}(M)$
does indeed form a group, as we will show in Appendix~\ref{ch:append1}.

We will provide a complete proof of Theorem~\ref{thm:dkn-intro} in this book, and it appears as Theorem~\ref{thm:lip-conj} in the body
of the monograph.

\subsection{Groups of prescribed critical regularity}\label{ss:prescribed}

The reader will note that in most of the discussion in the preceding sections that dealt with explicit, the natural classes of groups
that occur all have critical regularity in the set $\{1,2,\infty\}$. Before the appearance of obstructions to topological smoothing,
many authors even asserted that there seemed to be little qualitative difference between foliations (and hence, via suspensions,
group actions) in regularity $C^2$ and $C^{\infty}$ (see Cantwell--Conlon's preprint~\cite{CC-unpublished}, for instance).

Even after the appearance of topological obstructions, such as Tsuboi's and Cantwell--Conlon's  as we have mentioned already,
known groups of non-integral algebraic critical regularity were few, and there were no known examples of groups with integral critical
regularity three or more. The question of the existence of such groups was explicitly posed by Navas in his 2018 ICM
address~\cite{Navas-icm}, which
was written just before the paper~\cite{KK2020crit} appeared as a preprint.

We are now ready to state the existence theorem for groups of specified critical regularity in its
(nearly) fullest generality, as will occur again
as Theorem~\ref{t:optimal-all} below and the various variations on it.
The condition on the integral below may seem \emph{ad hoc}, though in fact there is some significant
leeway in making these choices, as we will see in Chapter~\ref{ch:optimal}

\begin{thm}\label{thm:crit-intro}
Let $k\in\bZ_{>0}$, 
and let $\alpha,\beta$ be concave moduli of continuity satisfying
\[
\int_{(0,1]} \frac1x \left( \frac{\beta(x)}{\alpha(x)}\right)^{1/k}dx<\infty.\]
If $k=1$, we further assume that $\beta$ is \emph{sub-tame}\index{sub-tame modulus} in the sense that
\[\lim_{t\to 0^+}\sup_{x>0}\frac{\beta(tx)}{\beta(x)}=0.\]
Then there is a finitely generated subgroup \[G=G_{\alpha,\beta,M}\le \Diff_+^{k,\alpha}(M)\] such that the following conclusions hold.
\begin{enumerate}[(1)]
\item
The commutator subgroup $[G,G]$ is nonabelian and simple.
\item
For all finite index subgroups $H\le G$, we have $[G,G]'\le H$.
\item
For all finite index subgroups $H\le G$ and for all \[\psi\colon H\longrightarrow\Diff_+^{k,\beta}(M),\] we have that the image of
$\psi$ is abelian.
\end{enumerate}
\end{thm}

A large portion of this book is devoted to giving a complete proof of Theorem~\ref{thm:crit-intro}. For H\"older moduli of continuity, we will
write $r=k+\tau$ and $\GG^r(M)$ for the set of isomorphism types of countable subgroups of $\Diff_+^{k,\tau}(M)$.
Applying the above theorem to the H\"older modulus $x^r$ and another modulus that is a slight perturbation or $x^r$  (as in Lemma~\ref{l:omega-st}), 
we obtain the result below.

\begin{cor}\label{cor:crit-intro}
Let $r\geq 1$. Then the following conclusions hold.
\begin{enumerate}[(1)]
\item
The set \[\GG^r\setminus\bigcup_{s>r}\GG^s\] contains continuum many isomorphism types of finitely generated groups and
countable simple groups.
\item
The set \[\left(\bigcap_{s<r}\GG^s\right)\setminus \GG^r\] contains continuum many isomorphism types of finitely generated groups and
countable simple groups.
\end{enumerate}
\end{cor}

For $r=1$, the second part of Corollary~\ref{cor:crit-intro} should be interpreted as groups admitting bi--Lipschitz but non--differentiable
actions. Note that the two parts of Corollary~\ref{cor:crit-intro} furnish groups for which the critical regularity is achieved and for which
it is not achieved, respectively.
We remark that the simple groups in Corollary~\ref{cor:crit-intro} are necessarily infinitely generated, as follows from Thurston's
Stability Theorem (see Appendix~\ref{ch:append3}).

In short, as to whether there exist groups acting with prescribed algebraic critical regularity on the interval and on the circle, the answer
is yes, and there are tons of examples. As might be expected, there are many more questions raised than there are answered by
Theorem~\ref{thm:crit-intro}. For instance:

\begin{itemize}
\item
What happens in higher dimension? Can finitely generated subgroups of diffeomorphism groups give insight into exotic smooth
structures on manifolds?
\item
For which moduli of continuity (even just restricting to H\"older moduli) are there finitely presented examples of groups with a prescribed
critical regularity? What about recursively presented examples?
\item
What are the coarse geometric properties of the groups furnished in Theorem~\ref{thm:crit-intro}? Do they form pairwise
distinct quasi--isometry classes?
\item
Is there a finitely generated group acting faithfully by $C^k$ diffeomorphisms on $M$ for all $k$, but that admits no faithful
$C^{\infty}$ action on $M$?
\end{itemize}

The last of these questions was explicitly posed in~\cite{KK2020crit}, and it seems rather difficult.

We close this section with a brief overview of another recent perspective on Theorem~\ref{thm:crit-intro} due to
Mann--Wolff~\cite{Mann:aa}, which uses different methods
to build groups with prescribed critical regularity. The basic idea is that one can combine groups that are differentially rigid with
diffeomorphisms with prescribed regularity properties in order to get groups with given critical regularity. 
\emph{Differential rigidity}\index{differential rigidity}
of a group $G$ means that $G\le \Diff^{\infty}(M)$ and that above some cutoff $r$ and $s\geq r$, all faithful $C^{s}$--actions of $G$ on $M$
are conjugate to the $C^{\infty}$ action by a $C^{s}$ diffeomorphism. Examples of differentially rigid group actions
on the interval (with cutoff $r=1$) appear in~\cite{BMNR2017MZ}, and on the circle (with cutoff $r=3$) appear in~\cite{Ghys93IHES}.
Suppose now we wish to construct a group with prescribed critical regularity $r=k+\tau$.
Roughly speaking, taking a differentially rigid group and adjoin a properly supported $C^r$ diffeomorphism
$f$, calling the resulting group $G_f$.
Then arbitrary realization of $G_f$ as a group of $C^s$ diffeomorphisms forces the differentially rigid group to be $C^s$--conjugate to
a $C^{\infty}$ action. If the orbit structure of the differentially rigid group is sufficiently rigid, this conjugacy has to send $f$ to a
conjugate of itself by a $C^s$ diffeomorphism, which will again be only $C^r$. It follows then that $G_f$ has critical regularity exactly
$r$.

While Mann--Wolff's construction does not achieve the algebraic control that Theorem~\ref{thm:crit-intro} does, it is remarkable for
other reasons. For one, it shows that the set of
isomorphism classes of finitely generated subgroups of $\Diff_+^{k,\mathrm{Lip}}(M)$ and $\Diff_+^{k+1}(M)$
do not coincide, something which comparisons of moduli of continuity in Theorem~\ref{thm:crit-intro} cannot achieve. Moreover,
their result shows that the set of isomorphism classes of subgroups of $\Diff_+^{k,\tau}(M)$ is as rich and diverse as the set
$\Diff_+^{k,\tau}(M)$ itself.

\section{What this book is about}

The purpose of this book is to give a self--contained account of the construction of groups with a prescribed critical regularity.
In the interest of a coherent and motivated discussion, we have included a large number of peripheral and related results, all of which
serve to contextualize and complete the theory we are describing. In the interest of space and readability, we have been compelled to make
choices to omit many subjects. Some discussion of what is omitted can be found in the next section.

Some of the material we have decided to include is already described in other books and surveys, and so there is some overlap
between this book and existing literature. Our primary aim in recounting such results is to illustrate how the theory of critical regularity fits in the
theory of dynamical systems and geometric group theory. Our secondary aim is to provide a reference that is useable as a standalone volume,
and sometimes repetition is necessary. We summarize the content of the book as follows.

\begin{enumerate}
\item
Chapter~\ref{sec:denjoy} gives a full account of dynamics and regularity for single homeomorphisms of the circle. One of the highlights
of the chapter is a complete proof of Denjoy's Theorem. There are many other secondary results and background that are developed along
the way. These include:
\begin{itemize}
\item
The existence of invariant measures for solvable group actions on the circle.
\item
The existence of stationary measures
for diffusion semigroups.
\item
Bi--Lipschitz conjugacy of countable groups of homeomorphisms.
\item
Poincar\'e's theory of rotation numbers.
\item
Unique ergodicity of homeomorphisms with irrational rotation number.
\item
Additivity of the rotation number on subgroups with an invariant measure.
\item
The construction of $C^1$ Denjoy counterexamples for a large class of moduli of continuity.
\end{itemize}

\item
Chapter~\ref{sec:filip-tak} gives an account of the results of Takens, Filipkiewicz, and Rubin,
which in particular show that the abstract isomorphism type of the
$C^p$ diffeomorphism group of a manifold determines the manifold up to $C^p$ diffeomorphism. Thus, while classical, this chapter
illustrates the continuous version of the existence of groups with prescribed critical regularity.
Takens' Theorem is proved in its original form in dimension one,
with the exception of one linearization result of Sternberg. Filipkiewicz's Theorem is proved in a generalized form; one proof of Filipkiewicz's Theorem
is obtained from Rubin's Theorem, which itself is also proved in its entirety. Among the additional results
established along the way are the Bochner--Montgomery result on simultaneous continuity, and the simplicity of the 
commutator subgroup of $\Diff_c^p(M)_0 $ for arbitrary $M$ and $p$; the latter of these figures prominently
in the discussion surrounding the Mather--Thurston Theorem, and also in the discussion of chain groups.

\item
Chapter~\ref{ch:c2-thry} develops many of the fundamental tools and methods used to investigate groups of 
$C^1$ and $C^2$ diffeomorphisms
of $I$ and $S^1$. Some of the ideas and methods contained therein are classical or semi--classical, and some were developed in recent
years by the authors. Among the highlights of the discussion are:
\begin{itemize}
\item
Kopell's Lemma on commuting diffeomorphisms (classical).
\item
The Plante--Thurston Theorem (classical).
\item
Farb--Franks, Jorquera, and Parkhe's theory of (residually) nilpotent groups acting on compact one manifolds (modern).
\item
The Two-jumps Lemma, which is developed by the authors and Baik (semi--classical)
and which generalizes an earlier result of Bonatti--Crovisier--Wilkinson.
\item
The $abt$--Lemma (modern).
\end{itemize}

\item
Chapter~\ref{ch:chain-groups} develops the theory of chain groups. These are a family of finitely generated groups of homeomorphisms
of the interval that exhibit a range of uniformity and diversity, and which are also an important technical tool in constructing groups with
specified critical regularity. Among the topics discussed here are the following:
\begin{itemize}
\item
Covering distances in groups of homeomorphisms.
\item
Algebraic and dynamical properties of chain groups.
\item
Ghys--Sergiescu's smooth realization of Thompson's group $F$.
\end{itemize}

\item
Chapter~\ref{ch:slp} is one of the two mostly technical chapters in the book. The purpose of this chapter is to establish the
so--called \emph{Slow Progress Lemma}, which is a highly technical result that makes precise the intuition that ``more regular"
diffeomorphisms are ``less expansive". The Slow Progress Lemma originally appeared in~\cite{KK2020crit}, in a weaker form than
that given in this book.

\item
Chapter~\ref{ch:optimal} is the other technical chapter in the book. Its function is to construct diffeomorphisms with optimal expansivity
properties, and to use them to construct finitely generated groups with prescribed critical regularity. In this chapter, the ideas from
Chapters~\ref{ch:c2-thry},~\ref{ch:chain-groups}, and~\ref{ch:slp} all come together to furnish the announced groups.

\item
Chapter~\ref{ch:app} reaps the benefits of the theory developed in the previous chapters. Among the results established there are:
\begin{itemize}
\item
The construction of foliations on closed $3$--manifolds (satisfying mild topological hypotheses) that have a prescribed critical regularity.
\item
Characterization of right-angled Artin groups and mapping class groups that have finite critical regularities.
\end{itemize}

\item
The appendices to the book cover relevant results which, in the opinion of the authors, do not fit as neatly into the narrative of the book
but which they believe to be central enough to merit inclusion. Some of these, such as H\"older's Theorem, are so central that one could
reasonably argue that they belong in the main body of the book, though because the proof known to the authors relies on orderability, we
have relegated it to the appendices. The topics in the appendices include:

\begin{itemize}
\item
Medvedev's Theorem on smooth moduli of continuity,
and the Muller--Tsuboi trick on flattening interval diffeomorphisms at the boundary.
\item
Background on orderable groups and H\"older's Theorem.
\item
The Thurston Stability Theorem.
\end{itemize}
\end{enumerate}

The topics included in the book undoubtedly reflect the personal biases and preferences of the authors. We have striven to provide
as complete and nuanced a picture of the theory of critical regularity as possible, and we have included many references in places
where we omit details.

\section{What this book is not about}
As we have already mentioned, we have excluded a large number potential topics of discussion.
This book is not a general treatise on dynamical systems, nor is it even a general reference on groups acting on the circle.
Here we note some major omissions from our monograph, all of which are discussed at length in other sources.

\begin{enumerate}
\item
Orderability of groups. Many excellent books on the subject already exist~\cite{BMR77,CR2016,DDRW08,DNR2014,Navas2011},
and we only need orderability in a few places. All the relevant
background for us is contained in Appendix~\ref{ch:append2}.
\item
The structure of full diffeomorphism groups, viewed as discrete groups.
This is the discussion surrounding the Mather--Thurston Theorem, and besides
the simplicity of commutator subgroups and the Takens--Filipkiewicz--Rubin theory, 
we will not be including an exposition on full diffeomorphism groups at all.
\item
$C^0$ theory. There is much to be said about actions of many types of groups, from closed surface groups to mapping class groups,
hyperbolic manifold groups, simple groups, and beyond on the circle and on the interval. There is also a vast theory of piecewise linear
groups acting on the interval and the circle. We will mostly exclude these topics from our account.
\item
Groups of analytic diffeomorphisms. The theory of subgroups of $\mathrm{PSL}_2(\bR)$ and its covers is well--developed (e.g.~\cite{KKM2019,KatokBook}), and we will
not discuss it. In general, groups of analytic diffeomorphisms do not furnish interesting examples in the circle of ideas on which
we are focussing.
For us, partial commutativity yields the most interesting examples and subsequent theory; among analytic diffeomorphisms, commutativity
is transitive, and so partially commutative groups of analytic diffeomorphisms are too simple to be interesting from our perspective.
\item
Solvable groups of diffeomorphisms. There is a structure theory for solvable subgroups of diffeomorphism groups of the
interval~\cite{Navas04solv},
which since it is peripheral to the discussion of critical regularity and since it is worked out in detail in~\cite{Navas2011}, we have decided
to omit it. Among the highlights of this theory are that polycyclic subgroups of $\Diff_+^{1+\mathrm{bv}}[0,1]$ are metabelian, but there
are solvable groups of diffeomorphisms of arbitrary long derived series length. We also direct the reader to~\cite{BursWilk04} for a discussion
of the rigidity of solvable group actions on the circle and a classification of solvable groups of analytic diffeomorphisms of the circle.
\item
Lattices in Lie groups. There is much to say about lattices in semisimple Lie groups in higher rank, and how these can act on manifolds;
indeed, this is a major concern of the Zimmer program. Other than what has already been said in the introduction of this book,
we have opted to avoid a discussion of this vast and rich subject.
\item
Geometric group theory. Many excellent books on geometric group theory and various classes of groups of interest in geometric group
theory already exist~\cite{dlHarpe2000,ClayMarg17,DM2018,Loeh2017,GGTIAS,Bowditch06,GH1990,GGT1990,BH1999}.
Right-angled Artin groups, mapping class groups, and nilpotent groups will be
sources of examples for us, though we will
only reference facts about them as needed.
\end{enumerate}

\section{What we will assume of the reader}

The authors have taken great pains to make this book as self-contained as possible. The writing contained herein should be accessible
to a beginning--to--intermediate level graduate student in mathematics. We will assume familiarity with standard notions from point--set
topology, and insofar as algebraic topology is concerned, we will only assume familiarity with the fundamental group, basic
covering space theory, and the classification of surfaces.
We will assume only the most basic notions about differentiable manifolds. Familiarity with hyperbolic geometry
is helpful but not necessary.

We will require some more sophisticated background in analysis,
to the level which is typically covered in a first graduate course in measure theory
and functional analysis. We will assume that the reader is comfortable with measure theory, the Riesz Representation Theorem for
functionals on the Banach space of continuous functions on a compact Hausdorff topological space, and the Banach--Alaoglu Theorem.
Sufficient background is found in the first few chapters of Rudin's standard textbook~\cite{Rudin87book}
and the first few chapters of~\cite{Zimmer-funct}.

We will assume that the reader is familiar with infinite groups and standard constructions with them. We will make free use of group
presentations and manipulations of them, and the reader should be familiar with the basics of groups acting on trees, to the level
of the early chapters of~\cite{Serre1977}. We will make free reference to nilpotent and solvable groups, as well as to the lower central
and derived series of a group. Familiarity with notions from coarse geometry and geometric group theory such as quasi--isometry
is useful but not necessary.


%
%
%
\chapter{Denjoy's Theorem and exceptional diffeomorphisms of the circle}\label{sec:denjoy}

\begin{abstract}
This chapter is a mostly self-contained account of the Denjoy--Herman--Yoccoz theory of circle diffeomorphisms, together with some
new contributions of the authors. Many beautiful expositions of Denjoy's theory are available in the literature -- the reader may consult
~\cite{athanassopoulos,Navas2011}, for instance. We do not pretend to improve on the work of those authors, and we reproduce many of the
ideas therein for completeness' sake and for the convenience of the reader.\end{abstract}

\section{The minimal set and the rotation number}
Let $f\in\Homeo_+(S^1)$. One of the most basic questions one can ask about $f$ concerns a characterization of the \emph{dynamics}
of $f$. Dynamics, broadly speaking, investigates the long term behavior of systems endowed with a prescribed transformation rule.
In the case of a homeomorphism $f\in\Homeo_+(S^1)$, or more generally a subgroup $\gam\le \Homeo_+(S^1)$, the system consists of a phase
space, namely the circle $S^1$, and the transformation rule is encoded by the function (or functions) $f$ (or $\Gamma$).

For a dynamical system like this, dynamical questions might be of the following form:
\begin{enumerate}
\item
What are the global fixed points of the group $\Gamma$? That is, characterize the points $x\in S^1$ such that $\gamma(x)=x$ for all $\gamma\in
\gam$.
\item
What are the finite orbits (i.e.~periodic points) of $\Gamma$?
\item
What are the long term itineraries of points (under the action of $\Gamma$) like?
\item
What are the $\Gamma$--invariant, closed subsets of $S^1$?
\item
What sorts of probability measures (or measure classes) on $S^1$ are $\Gamma$--invariant?
\end{enumerate}

This chapter will give a complete or at least extensive answer to all of these questions, at least in the case where $\gam\cong\bZ$ is
generated by a single homeomorphism of $S^1$. We will begin by introducing the minimal set of a group of homeomorphisms of $S^1$,
and then introduce the rotation number as a fundamental and powerful dynamical invariant of a homeomorphism. The interaction between
the minimal set, rotation number, and analytic concerns (i.e.~regularity) will be the subject of the remainder of the chapter.

\subsection{Minimal sets and exceptional diffeomorphisms}
We begin by giving an answer to the following question: let $\gam\le \Homeo_+(S^1)$. What are the possible closed, $\Gamma$--invariant subsets
of $S^1$? Are there canonical such $\Gamma$--invariant subsets?

\subsubsection{Periodic orbits}
To get a feeling for the foregoing question and possible answers to it, we note some examples.
First, consider the homeomorphism of the real line given by
$x\mapsto x+1$. This is clearly orientation preserving, and by adding a point at infinity, we obtain an element $f\in\Homeo_+(S^1)$ that has
exactly one fixed point. This point has ``parabolic" dynamical behavior, in the sense that points move away from infinity in the far negative
part of $\bR$, and move towards infinity in the far positive part of $\bR$. To get a sense of what the
proper and nonempty closed invariant subsets are, we note that infinity is clearly invariant. Moreover, infinity together with the orbit of a point
in $\bR$ is closed and invariant. More generally, one can consider a closed subset of $[0,1]$,
take the union of its translates by $x\mapsto x+1$
and infinity, and thus produce a closed invariant subset of $S^1$.

For a slightly more complicated example that retains some of
the same flavor, we may consider a rotation of $S^1$ by a rational number. Clearly, every point of such a homeomorphism is periodic. One can combine such a rotation with an arbitrary homeomorphism
of the real line by noting that the interval between two consecutive points in an orbit of the rotation
(in the direction of rotation) is homeomorphic
to $\bR$. 
Thus, one can take an arbitrary homeomorphism of $\bR$ and propagate it around the circle via the rational rotation,
arranging for example (if the homeomorphism of $\bR$ has no fixed points) a homeomorphism
of $S^1$ that has exactly one periodic orbit. The resulting homeomorphism's closed invariant subsets can be analyzed in a manner similar
to the previous example.

\subsubsection{Generalities on rotations}

The preceding discussion introduces the notion of a rotation of $S^1$. Certain concrete realizations of the circle allow for straightforward
descriptions of rotations. For instance, if $S^1$ is realized as the unit complex numbers in $\bC$, then rotation by an angle $\theta$
is realized by mutiplication by the complex number $e^{i\theta}$. If the circle is realized as $\bR/\bZ$, then rotation by an angle
$\theta$ is simply realized by addition of $\theta \pmod{\bZ}$. The centrality of rotations as
one of the most basic examples of a homeomorphism
of the circle arises from a more abstract characterization of the circle as a one--dimensional compact Lie group. The circle,
viewed as an abstract group, is circularly orderable (see the Appendix~\ref{ss:circular}), and there is a unique topology that is compatible
with this circular order (in that open intervals in the circular order are a basis of open sets in the topology). It turns out in fact that this topology
determines a smooth manifold structure on the circle that is unique up to diffeomorphism (see~\cite{Tao-Hilbert}, for instance).
With this setup, rotations of the circle are merely the left (or right) regular action of the circle on itself, viewed as a group. A
fundamental question that is investigated in this chapter is, how does one distinguish rotations of the circle (i.e.~homeomorphisms of
the circle arising from the action of the circle on itself) via their dynamics?

\subsubsection{Irrational rotations and a digression on ergodic theory}
To introduce a central actor in the subsequent discussion,
consider a rotation of $S^1$ through an irrational number.
Not only are there
no periodic points, but in fact every orbit is dense. A homeomorphism for which every orbit is dense is called 
\emph{minimal}\index{minimal homeomorphism}. Equivalently, the
closure of every orbit is the whole phase space, and so the only invariant closed subspaces are the empty set and the whole space.

\begin{prop}\label{prop:irr-rot-min}
Let $f\in\Homeo_+(S^1)$ be a rotation through an irrational number. Then $f$ is minimal.
\end{prop}
\begin{proof}
Note that the \emph{radial}\index{radial metric} (i.e.~arclength) metric on $S^1$ is 
invariant by $f$, so that $f$ acts by an \emph{isometry}\index{isometry} 
of the radial metric.
Now, since the angle $\theta$ through which
$f$ rotates is irrational, no point has a finite orbit under $f$. This is simply because if $f^n(x)=x$ then, writing
$S^1$ additively as $\bR/\bZ$, we have \[f^n(x)=x+n\theta\equiv x\pmod {\bZ},\] which means that $n\theta$ is an integer, a
contradiction. A straightforward compactness argument shows that, if $\OO$ is the orbit of $f$ on
$x$ and $\epsilon>0$, then there are two points
in $\OO$ at distance at most $\eps$ from each other.

Now, suppose that $\OO$ is not dense. Let $U\sse S^1$ be a component
of $S^1\setminus\overline{\OO}$, and assume that $U$
has length $\ell>0$. Choose points \[x_n=f^n(x),\quad x_m=f^m(x)\] such that the radial distance between $x_n$ and $x_m$ is less than
$\ell$. Then $f^{n-m}(x)$ is at a distance smaller than $\ell$ from $x$. Now, since both endpoints of $U$ are accumulation points of the closure
of $\OO$, there is an $N$ such that $f^N(x)$ is as close to either endpoint of $U$ as we like. It follows then that $f^{N+n-m}(x)$ will lie in
$U$, for a suitable $N$. Therefore, $\OO$ is dense.
\end{proof}

Irrational rotations have yet another important property, closely related to their minimality: they are \emph{ergodic}\index{ergodic}.
Recall that the radial
length function extends to a unique Borel measure on $S^1$, which is a multiple of Lebesgue measure on $S^1$.
Since rotations act by isometries on $S^1$ with respect to the radial metric and because Lebesgue measure is determined by its value
on open intervals,
we have that rotations also preserve Lebesgue measure. A measurable map $f\colon S^1\longrightarrow S^1$ is called \emph{ergodic}
if the only measurable
$f$--invariant subsets of $S^1$ have zero measure or full measure. To see the relationship between ergodicity and minimality, the reader
will be able to show (after absorbing Theorem~\ref{thm:minimal-set}) that a homeomorphism of $S^1$ that preserves Lebesgure measure
and that is ergodic with respect to Lebesgue measure is necessarily minimal.

\begin{prop}\label{prop:irr-rot-erg}
Let $f\in\Homeo_+(S^1)$ be an irrational rotation. Then $f$ is ergodic with respect to Lebesgue measure.
\end{prop}
\begin{proof}
Let $g\in L^2(S^1)$ be a square--integrable function that is invariant under rotation, i.e.~$g\circ f=g$. It suffices to show that $g$ is constant
almost everywhere with respect to Lebesgue measure. This will suffice to prove the ergodicity of $f$,
since one can set $g$ to be the indicator function
of a measurable set.

Now, we expand $g$ in $L^2(S^1)$ by Fourier theory, viewing $S^1=\bR/\bZ$.
For $n\in \bZ$, we have that the $n^{th}$ Fourier coefficient of $g$ is 
\[c_n=\int_{S^1}g(t)e^{-2\pi i n t}\,dt,\] so that \[\sum_{n\in \bZ} c_n e^{2\pi i n t}\longrightarrow g\] in $L^2(S^1)$. Now, if $g$ is invariant under
$f$, we must have that the Fourier coefficients of $g$ and $g\circ f$ are the same. Viewing the circle additively and computing, we have
\[c_n=\int_{S^1}g(t+\theta)e^{-2\pi i n t}\,dt=\int_{S^1} g(t)e^{- 2\pi i n (t-\theta)}\,dt=e^{2\pi i n\theta}c_n.\] Since $\theta$ is irrational, we have
that $e^{2\pi i n\theta}\neq 1$ for all $n\neq 0$, so that all Fourier coefficients vanish, other than the zeroth one. It follows that $g$ is constant
almost everywhere.
\end{proof}

In fact, even more is true about irrational rotations: they are \emph{uniquely ergodic}\index{uniquely ergodic}.
That is, they admit a unique invariant probability
measure on $S^1$, which therefore must be Lebesgue measure. If a probability measure $\lambda$ is the unique invariant measure for
a measurable
transformation
$f$, then $f$ is automatically ergodic. To see this, if $A$ and $B$ measurably partition the phase space and both have positive measure
with respect to $\lambda$, then
one can build a measure $\lambda_A$, which is just $\lambda$ restricted to $A$ and rescaled to have measure one, and such that
$\lambda_A$ assigns measure zero to $B$. A measure $\lambda_B$ can be built similarly. 
Clearly both $\lambda_A$ and $\lambda_B$
are invariant, as is every convex combination of them. So, a uniquely ergodic transformation is automatically ergodic.

The proof of Proposition
\ref{prop:irr-rot-erg} actually shows that if $\mu$ and $\lambda$ are
arbitrary Borel probability measures that are invariant under a rotation through an
irrational number, and if $f$ is a continuous function on $S^1$, then \[\int_{S^1}f\,d\mu=\int_{S^1}f\,d\lambda.\] Indeed,
these integrals are merely the $0^{th}$ Fourier coefficients of $f$. Since Borel measures are separated by integrating against continuous
functions (by Lusin's Theorem for example~\cite{Rudin87book}), it follows that $\mu=\lambda$.

\begin{cor}\label{cor:leb-ue}
Let $f\in\Homeo_+(S^1)$ be a rotation through an irrational number. Then $f$ is uniquely ergodic with respect to Lebesgue measure.
\end{cor}

In Subsection~\ref{ss:rot} below, we will introduce the rotation number, and we will show that in fact all
homeomorphisms with irrational rotation number
(of which irrational rotations themselves are examples) are uniquely ergodic (see Theorem~\ref{thm:irr-ue}). This fact is a crucial ingredient
in Denjoy's theory of circle diffeomorphisms as we will develop it here.

One of the general and useful properties of a uniquely ergodic transformation $T$ of a phase space $X$ is that every sequence of
probability measures that converges to a $T$--invariant measure must converge to the unique $T$--invariant probability measure on $X$.
One such sequence of probability measures can be described as follows, and is particularly intuitive when $T$ is
a homeomorphism of a compact Hausdorff topological space (as we shall
assume).

For $n\in\bN$ and $x\in X$, let $\mu_{n,x}$ be a probability measure on $X$ defined by
\[\mu_{n,x}(A)=\frac{1}{2n+1}\left|A\cap \left\{\bigcup_{i=-n}^n \{f^i(x)\}\right\}\right|,\] where $A$ is a Borel subset of $X$. Note that the
set $\{\mu_{n,x}\}_{n\ge1}$ is precompact by the Banach--Alaoglu Theorem (see~\cite{Zimmer-funct}, for instance),
and so we may extract a weak--$\ast$ limit $\mu$.
A straightforward manipulation shows that $\mu$ is $T$--invariant and is hence the unique $T$--invariant probability measure on $X$.
We will spell out these ideas in more detail in Lemma~\ref{lem:kaku-abel}.

\subsubsection{Continuous Denjoy counterexamples}
There are examples of homeomorphisms of $S^1$ which have neither periodic points, nor are they minimal. In fact, these are the generic
types of homeomorphisms, though imagining them at first is not always so easy. One of the most important examples comes from
a \emph{blow-up}\index{blow-up} of an orbit of a minimal homeomorphism. 
We now describe this construction for an irrational rotation of the circle, though
the method is broadly applicable. This type of example is called 
a~\emph{(continuous) Denjoy counterexample}\index{Denjoy counterexample}, for reasons which will become
more apparent later in this chapter.

Let $f$ be an irrational rotation, and let $\OO$ be the orbit of a point $x\in S^1$.
Let $\{\ell_n\}_{n\in\bZ}$ be a collection of positive real numbers
such that \[\sum_{n\in\bZ}\ell_n<\infty.\] We construct a sequence of intervals $J_n=[0,\ell_n]$ for $n\in\bZ$, and homeomorphisms
\[\phi_{n,m}\colon J_n\longrightarrow J_m\] given by \[\phi_{n,m}(x)=\ell_n^{-1}\cdot \ell_m\cdot x.\] Now, let $X$ be the space obtained
by inserting a copy of $J_n$ at $f^n(x)\in\OO$. More precisely, we cut $S^1$ open at $x_n$, and glue the two endpoints of $J_n$ to
the two preimages of $x_n$. The resulting topological space is easily seen to be homeomorphic to $S^1$. One now defines a map
$f_X\colon X\longrightarrow X$ by setting $f_X(y)=f(y)$ if $y\notin\OO$, and $f_X(y)=\phi_{n,n+1}(y)$ if $y\in J_n$. It is again
straightforward to check that $f_X$ defines a homeomorphism of $X\cong S^1$, and evidently $f_X$ has no periodic points. We claim that
\[C=X\setminus\bigcup_{n\in\bZ} J_n^0\] is homeomorphic to a Cantor set, where here $J_n^0$ denotes the interior of $J_n$. Because
this argument is on the long side and because the result is important, we isolate it formally.

\begin{prop}\label{prop:blowup-cantor}
The set $C$ is homeomorphic to a Cantor set.
\end{prop}

The reader may compare the statement and proof of Proposition~\ref{prop:blowup-cantor} to that of Theorem~\ref{thm:minimal-set} below.

\begin{proof}[Proof of Proposition~\ref{prop:blowup-cantor}]
Note first that $C$ is compact, since $X$ is compact Hausdorff and $C$ is closed. It is also clear that $C$ is nonempty.
The subspace $C$ is clearly metrizable since it inherits the radial metric from $X\cong S^1$. To complete
the verification that $C$ is homeomorphic to a Cantor set, it suffices to show that it is perfect and totally disconnected. To see that $C$ is
totally disconnected, note that there is a map $\kappa \colon X\longrightarrow S^1$ given by collapsing the intervals $\{J_n\}_{n\in\bZ}$ each to
a point. This map is one--to--one except on the endpoints of the collapsed intervals, where it is two--to--one. If $y\in\OO$ then it is clear
that the two points in the preimage $\kappa^{-1}(y)$ are separated by open sets. Now, if $y,z\in S^1$ are arbitrary and distinct
then $\OO$ meets
both components of $S^1\setminus\{y,z\}$, and so that arbitrary points
\[y_X\in\kappa^{-1}(y),\quad z_X\in\kappa^{-1}(z)\] are separated by open sets that meet the complement of $C$.
It follows that $C$ is totally disconnected.

To see that $C$ is perfect, let $y\in S^1$ be arbitrary, and choose $y_X\in\kappa^{-1}(y)$. There is a sequence of points $\{x_n\}\sse\OO$
such that $x_n\longrightarrow y$, by Proposition~\ref{prop:irr-rot-min}. 
Note that for each $M\in\bZ$, there is an $N$ such that
for all $n\geq N$, the only points in $\OO$ that meet the (shorter) interval $(x_n,y)$ are contained in \[\OO_M=\{f^m(x)\mid |m|\geq M\}.\]
Since the total length of the intervals $\{J_n\}_{n\in\bZ}$ is finite, we have that for each $\eps>0$ there is an $M\gg 0$ such that
\[\sum_{|m|\geq M}\ell_m<\eps.\] It follows that the radial distance between $\kappa^{-1}(x_n)$ and $y_X$ tends to zero as $n$ tends to
infinity. Therefore, $y_X$ is not an isolated point of $C$, and the proof is complete.
\end{proof}

It is clear that the set $C$ and the union $J=\cup J_n^0$
of the interiors of the intervals $\{J_n\}_{n\in\bZ}$ are both invariant under $f_X$. 
The set $J$ is called the \emph{wandering set}\index{wandering set} and is
characterized by the fact that if $I\sse J_n^0$ for some $n\in\bZ$ then $f^m(I)\cap I\neq\varnothing$ implies $m=0$. The maximal
connected components of $J$ are called \emph{wandering intervals}\index{wandering interval}.
It turns out
that $C$ is the unique nonempty closed invariant subset of $X$ on which $f_X$ acts minimally, but we will not justify this claim just yet.

\subsubsection{Minimal sets for general group actions}
The preceding examples give examples of the typical minimality phenomena observed among groups of homeomorphisms of $S^1$. This
is made precise by the following basic result.

\begin{thm}[cf.~\cite{Navas2011}, Theorem 2.1.1 or ~\cite{athanassopoulos}, Proposition 2.5]\label{thm:minimal-set}
Let \[\gam\le \Homeo_+(S^1).\] Exactly one of the following conclusions holds.
\begin{enumerate}[(1)]
\item
There is a finite $\Gamma$--orbit.
\item
The group $\Gamma$ acts minimally on $S^1$.
\item
There is a unique, nonempty, closed, invariant $C\sse S^1$ such that $\Gamma$ acts minimally on $C$ and such that $C$ is homeomorphic
to a Cantor set. Moreover, if $\OO$ is an arbitrary orbit of $\Gamma$ then the closure $\overline{\OO}$ contains $C$.
\end{enumerate}
\end{thm}
\begin{proof}
We use Zorn's Lemma, applied to the partially ordered set $\PP$ consisting of closed, nonempty, $\Gamma$--invariant subsets of $S^1$, ordered
by reverse inclusion. We have
that $\PP$ is obviously nonempty, since $S^1\in\PP$. Chains have upper bounds in $\PP$, since closed subsets of $S^1$ are compact.
In particular, if $\{K_i\}_{i\in I}$ is a collection of nested compact subsets of $S^1$ then the finite intersection property characterization of
compactness implies that \[\bigcap_{i\in I} K_i\neq\varnothing.\] Zorn's Lemma allows us to extract a minimal element $C\in\PP$.

Observe first that the $\Gamma$--action on $C$ is (topologically) minimal.
This is more or less immediate, since otherwise there would be a point $x\in C$
whose orbit $\OO(x)$ is not dense in $C$. Then, the closure $\overline{\OO(x)}$ is a closed $\Gamma$--invariant subset of $S^1$ that is
properly contained in $C$, violating the minimality (with respect to the partial order) of the choice of $C$.

We write $C'$ for the derived set of $C$. Recall that this is the set of accumulation points of $C$, or equivalently the subset of $C$ obtained
by discarding all isolated points. Observe that $C'$ is a closed and $\Gamma$--invariant subset of $S^1$, and so $C'$ is either empty or equal to
$C$, by the minimality of the choice of $C$. If $C'$ is empty then $C$ is discrete and hence finite, which means that $\Gamma$ has a finite orbit.

Otherwise, we have that $C'=C$, and so we consider the boundary $\partial C$. Again, $\partial C$ is closed and $\Gamma$--invariant, and so
either \[\partial C=C'\quad \textrm{or} \quad \partial C=\varnothing.\]
The only way the latter equality can hold is if $C'=C=S^1$, in which case the original
$\Gamma$--action on $S^1$ is minimal.

Thus, we are left with the case \[C'=C=\partial C.\] In this case, we have that $C$ is compact (since it is closed),
perfect (since it has no isolated points), and totally disconnected (since it is equal to its boundary and hence contains no
intervals). It follows that
$C$ is homeomorphic to a Cantor set.

It remains to show that $C$ is contained in the closure of the orbit of an arbitrary point $x\in S^1$. This will also prove that $C$ is unique,
since $C$ is then just the intersection of all orbit closures.
Let $y\in C$ and $x\in S^1$ be arbitrary. If $x\in C$ as well, then the fact that the $\Gamma$--action on $C$ is minimal implies that there
is a sequence of elements $\{\gamma_i\}_{i\ge1}\sse\gam$
such that \[\gamma_i(x)\longrightarrow y\] as $i$ tends to infinity. Otherwise, we may
assume that $x$ is contained in an open interval $J\sse S^1\setminus C$. Writing \[\{x_1,x_2\}=\partial J\sse C,\] we have that 
the $\Gamma$ orbits of both $x_1$ and $x_2$ are dense in $C$. So, there are group elements $\{\gamma_i\}_{i\ge1}\sse\gam$ such that
$\gamma_i(x_1)\neq\gamma_j(x_1)$ for $i\neq j$ and such that
\[\gamma_i(x_1)\longrightarrow y\] as $i$ tends to infinity. Since $\gamma_i(J)$ is an interval in $S^1\setminus C$ for all $i$, and since
$\Gamma$ acts by orientation preserving homeomorphisms, we have that \[\gamma_i(J)\cap\gamma_j(J)=\varnothing\] for $i\neq j$. It follows
then that the length $|\gamma_i(J)|$ tends to zero as $i$ tends to infinity, so that the radial distance between $\gamma_i(x)$ and
$\gamma_i(x_1)$ also tends to zero. The triangle inequality immediately implies that the radial distance from $\gamma_i(x)$ and $y$
also tends to zero, so that the $\Gamma$--orbit of $x$ accumulates on all of $C$. This completes the proof of the theorem.
\end{proof}

The set $C$ furnished by Theorem~\ref{thm:minimal-set} in the absence of a finite orbit is called the \emph{minimal set}\index{minimal set}. 
If $C\neq S^1$,
we call the set $C$ the \emph{exceptional minimal set}\index{exceptional minimal set}
of the $\Gamma$--action on $S^1$. If \[\gam\cong\bZ=\form{ f}\] admits an
exceptional minimal set, then we call $f$ \emph{exceptional}\index{exceptional homeomorphism}.
The continuous Denjoy counterexample we have introduced in this subsection
is an example of an exceptional homeomorphism of the circle.

\subsection{The rotation number}\label{ss:rot}

The rotation number is a basic dynamical invariant of a homeomorphism of $S^1$ that measures how the orbit of a point travels around
the circle. The definition is simple enough. Let $f\in\Homeo_+(S^1)$. Choose an arbitrary lift $\tilde f\in\Homeo_+^{\bZ}(\bR)$.
Choose an arbitrary $x\in\bR$ and let \[\rot(f)=\lim_{n\to\infty}\frac{F^n(x)}{n}\in\bR/\bZ.\] 
The point $\rot(f)\in S^1$ is defined to be the \emph{rotation number}\index{rotation number} of $f$.

Of course, there are a number of things to check to make sense of this definition.
We need to verify the existence of the limit,
and the independence of the definition from the choice of $\tilde f$ and $x$.
For this, we first record a basic fact as  below.

Let $G$ be a group.
A real valued function $\varphi$ on $G$ is called a \emph{quasimorphism}\index{quasimorphism}
if
there exists a $C\in\bR$ such that for all $g,h\in G$, we have
\[ |\varphi(gh)-\varphi(g)-\varphi(h)|\le C.\]
The infimum of such a constant $C$ is the \emph{defect}\index{defect of a quasimorphism} of the quasimorphism.

This definition is a mild variation on a
\emph{subadditive sequence}\index{subadditive}, which is easiest to define for 
\[
\varphi\co \bN\longrightarrow\bR_+,\]
in which case we would simply require
$\varphi(m+n)\leq \varphi(m)+\varphi(n)$. For subadditive functions, the following result is known as Fekete's Lemma.

\begin{lem}[cf.~\cite{Navas2011}, Lemma 2.2.1, for instance]\label{lem:subadditive}
If $f\co \bN\longrightarrow\bR$ is a quasimorphism,
then 
 the limit \[\rho=\lim_{n\to\infty} f(n)/n\] exists and is the unique
real number for which the sequence $\{f(n)-n\rho\}_{n\ge1}$ is bounded.
\end{lem}
\begin{proof}
The idea is to realize $\rho$ as the intersection of a nested sequence of closed intervals, whereby the finite intersection property
characterization of compactness says that the intersection of these intervals is nonempty. With this in mind, let $C$ denote the defect
of the quasimorphism $f$ and set \[J_n=\left[\frac{f(n)-C}{n},\frac{f(n)+C}{n}\right].\] 

If $k\geq 1$ is an integer then we apply the triangle inequality to conclude that
\[k f(n)-k C
\le f(k n)-C \le 
f(k n)+C\leq k f(n)+k C.\]
Dividing through by $k n$, we have $J_{kn}\sse J_n$.
It follows that every finite collection \[\{J_{n_1},\ldots,J_{n_i}\}\] has nonempty intersection,
and so the intersection \[J=\bigcap_{n\ge1} J_n\neq\varnothing,\] by the finite intersection property characterization of compactness.

If $\rho\in J$ then clearly $|f(n)-n\rho|\leq C$, so that the sequence $\{f(n)-n\rho\}_{n\ge1}$ is bounded. We claim that $J$ consists of
exactly one point. Indeed, suppose that $\delta\neq \rho$ are two distinct elements of $J$. Then $|\rho-\delta|>0$, which we can compute
\[|f(n)-n\delta|=|f(n)+n(\rho-\delta)-n\rho|\geq n|\rho-\delta|-C,\] since $|f(n)-n\rho|$ is bounded by $C$. It follows then that $|f(n)-n\delta|$
tends to infinity, a contradiction. Thus, $\rho$ is unique. That \[\rho=\lim_{n\to\infty}\frac{f(n)}{n}\] is now immediate.
\end{proof}

The following is immediate from Lemma~\ref{lem:subadditive},
and we leave the detail to the reader.
The map $\bar\varphi$ will be called the \emph{homogenization}\index{homogenization} of $\varphi$.

\begin{lem}[cf.~\cite{Calegari2009,Navas2011,KKM2019}]\label{lem:qm}
If $\varphi$ is a real--valued quasimorphism defined on a group $G$,
then for each $g\in G$ the limit \[\bar \varphi(g)=\lim_{n\to\infty} \varphi(g^n)/n\] exists. 
Furthermore, the following hold.
\be[(1)]
\item The value $\bar \varphi(g)$ is the unique
real number for which the sequence $\{\varphi(g^n)-n\bar \varphi(g)\}_{n\ge1}$ is bounded;
\item $\bar\varphi(g^m)=m\bar\varphi(g)$ for all $g\in G$ and $m\in\bZ$;
\item $\bar\varphi(hgh^{-1})=\bar\varphi(g)$ for all $g,h\in G$.
\ee
\end{lem}

Now let us now consider an arbitrary $F\in\Homeo_+^\bZ(\bR)$, where here $\Homeo_+^\bZ(\bR)$ denotes the group of $\bZ$--periodic
homeomorphisms of $\bR$.
For a given point $x\in \bR$, 
we define
\[\varphi_x(F):=F( x)- x\in\bR.\]
We note the following.
\be[(i)]
\item If $x-y\in\bZ$ then $\varphi_x(F)=\varphi_y(F)$.
\item For $x,y\in\bR$, we have $\varphi_x(F)-\varphi_y(F)< 1$.
\item The map $F\mapsto \varphi_x(F)$ is a quasi-morphism on $\Homeo_+^\bZ(\bR)$.
\ee
Part (i) is obvious by letting $y=x+m$ and noting $F(y)=F(x)+m$.
For part (ii), we may first assume by part (i) that $y\le x<y+1$. Then we have
\[
\varphi_x(F)-\varphi_y(F)=F(x)-F(y)-(x-y)< F(y+1)-F(y)-(y-y)=1.\]
It is now easy to deduce part (iii), since for $F,G\in\Homeo^\bZ_+(\bR)$ we have
\[
\varphi_x(F\circ G)
=F\circ G(x)-x
= \varphi_{G(x)}(F)+\varphi_x(G)\in(\varphi_x(F)+\varphi_x(G)-1,\varphi_x(F)+\varphi_x(G)+1).\]

By Lemma~\ref{lem:qm} we can homogenize $\varphi_x$ and obtain
a homogeneous class function $\bar\varphi\co \Homeo_+^\bZ(\bR)\longrightarrow\bR$
defined by
\[
\bar\varphi(F):=\lim_{n\to\infty}\varphi_x(F^n)/n=\lim_{n\to\infty} (F^n(x)-x)/n
=\lim_{n\to\infty} F^n(x)/n.\]
The condition (ii) above implies that $\bar\varphi$ is independent of the choice of $x\in\bR$, justifying the reference-point--free notation.

Let us now suppose $f\in \Homeo_+(S^1)$.
If $F$ and $G$ are two lifts of $f$ in $\Homeo^\bZ_+(\bR)$ then 
the value $F(x)-G(x)$ is a fixed integer, and hence
$\bar\varphi(F)$ and $\bar\varphi(G)$ differ by that integer.
Quotienting by $\bZ$, we obtain a well-defined map
\[
\rot(f):=\bar\varphi(F)=\lim_{n\to\infty} \frac{F^n(x)}n \in\bR/\bZ,\]
where $F$ is an arbitrary lift of $f$, and $x\in\bR$ is also arbitrary.

As a matter of terminology, we will say that the rotation number of a homeomorphism is \emph{rational}\index{rational rotation} 
if it belongs to $\bQ/\bZ$, and \emph{irrational}\index{irrational rotation} otherwise. 
Let us summarize important properties of the rotation number.

\begin{prop}\label{prop:rot-easy}
The rotation number \[\rot\colon\Homeo_+(S^1)\longrightarrow S^1\] enjoys the following properties.
\begin{enumerate}[(1)]
\item
Homogeneity: that is, for all $f\in\Homeo_+(S^1)$ and $n\in\bZ$, we have \[\rot(f^n)=n\rot(f).\]
\item
Class function: $f,g\in\Homeo_+(S^1)$ then $\rot(g^{-1}fg)=\rot(f)$.
\item
If $R_{\theta}$ is a rotation of $S^1$ by an angle $\theta\in S^1$ then $\rot(R_{\theta})=\theta$.
\item
The map $f$ has a fixed point if and only if $\rot(f)=0$.
\item The map $f$ has a periodic point if and only if $\rot(f)$ is rational;
moreover, all periodic points have the same period.
\end{enumerate}
\end{prop}
\begin{proof}
The first three parts are immediate consequences of Lemma~\ref{lem:qm}.
For the fourth part of the proposition, the well-definedness of the rotation number implies that we need only show that if $\rot(f)=0$ then
$f$ has a fixed point. To show the contrapositive, let $f\in\Homeo_+(S^1)$ be fixed-point-free. Let $F$ be the unique left $f$ such that
$F(0)\in (0,1)$. Because $f$ has no fixed points, the function $F(x)-x$ does not achieve values in the set $\bZ$.
Since $F$ is continuous and periodic, it
satisfies the conclusions of the extreme value theorem, whence we may assume that there is an $\epsilon$ such that
\[0<\eps\leq F(x)-x\leq 1-\epsilon\] for all $x\in\bR$.
Observe that for all $n$ we have the estimate \[0<n\cdot\eps\leq \sum_{i=0}^{n-1}F(F^i(0))-F^i(0)=
F^n(0)-0\leq n(1-\eps).\] Dividing through by $n$, we see that $F^n(0)/n$ is bounded away from $0$ and $1$, whence the rotation
number of $f$ is nonzero. 

The homogenity of $\rot$ now implies that 
$\rot(f)$ is rational if and only if $f^n$ fixes a point for some nonzero $n\in\bZ$.
To see that all periodic points have the same period, let $x,y\in S^1$ have periods $n$ and $m$ respectively,
with $n<m$. Then $f^n$ admits $x$ as a fixed point, and so $\rot(f^n)=0$. However, $f^n(y)\neq y$ and $f^{nm}(y)=y$, and since
$f$ is orientation preserving, we have that $f^n$ has no global fixed points, which is a contradiction.\end{proof}

The continuous Denjoy counterexample we have introduced above gives an example of a homeomorphism of $S^1$ with irrational rotation
number, as the reader may check, and which is obviously not minimal. It turns out that examples of this ilk, which is to say ones constructed
by blowing up orbits, are essentially the only examples of homeomorphisms of the circle which have irrational rotation number and which
are exceptional.

\begin{thm}\label{thm:irr-rot-semi}
Let $f\in\Homeo_+(S^1)$ be such that $\rot(f)=\theta$ is irrational. Then there exists a $\bZ$--periodic map \[H\colon\bR\longrightarrow\bR\]
such that:
\begin{enumerate}[(1)]
\item
We have $H(0)=0$.
\item
The map $H$ is continuous and nondecreasing.
\item
The map $H$ descends to a self--map \[h\colon S^1\longrightarrow S^1\] such that:
\begin{enumerate}[(i)]
\item
The map $h$ is a degree one, orientation preserving, continuous surjection.
\item
If $R_{\theta}$ denotes rotation by $\theta$, we have $h\circ f=R_{\theta}\circ h$.
\item
The map $h$ is a homeomorphism if and only if $f$ is minimal.
\end{enumerate}
\end{enumerate}
\end{thm}

The maps $H$ and $h$ are oftentimes called \emph{semi-conjugacies}\index{semi-conjugacy}.
We will give the proof of Theorem~\ref{thm:irr-rot-semi} below,
after Theorem~\ref{thm:kaku-mar}. 
The notion of a semi-conjugacies is intuitively straightforward, because topology in one dimension is so limited.
If $H$ is a semi-conjugacy, then it is not difficult to see that the preimage $H^{-1}(x)$ of a point is either a singleton point or a compact
interval. Thus, a semi-conjugacy behaves like a homeomorphism near some points (and therefore restricts to a bijection
on those points), and the remaining intervals are collapsed to
points.

Before closing this subsection, we give the following fact, which shows that the rotation number of a homeomorphism
is a semi-conjugacy invariant.

\begin{prop}\label{prop:semi-rot}
Let $f,g\in\Homeo_+(S^1)$ be semi-conjugate. Then $\rot(f)=\rot(g)$.
\end{prop}
\begin{proof}
Let $h$ be a semi-conjugacy from $f$ to $g$, and let $\{F,G\}$ be arbitrary lifts of $\{f,g\}$ to $\Homeo_+^{\bZ}(\bR)$ and similarly
let $H$ be a lift of $h$ (which of course may not be a homeomorphism). Note that
\[H\circ F^n\equiv G^n\circ H \pmod{\bZ}\] for all $n\in\bZ$. It is convenient to write $\Delta(x)=H(x)-x$. Then, we can compute:
\[F^n(x)+\Delta(F^n(x))\equiv G^n(H(x))-H(x)+x+\Delta(x)\pmod{\bZ}.\] It follows that \[\frac{F^n(x)-x)}{n}+\frac{\Delta(F^n(x))}{n}\equiv
\frac{G^n(H(x))-H(x)}{n}+\frac{\Delta(x)}{n}\pmod{\bZ}.\] Since both $\Delta\circ F^n$ and $\Delta$ are bounded, the limits of the two sides
as $n\to\infty$ are equal, as claimed.
\end{proof}

\subsection{Rotation numbers, invariant measures, and amenability}

Another important perspective on rotation numbers for circle homeomorphisms is provided from the angle of invariant probability measures.
Suppose $\mu$ is a probability measure on $S^1$, and suppose that $\mu$ is invariant under the action of a group $\gam\le \Homeo_+(S^1)$.
Let $f\in\gam$.
For $x\in S^1$, let \[\rho_{\mu}(f)=\mu[x,f(x))\pmod{\bZ}.\]
The very notation $\rho_{\mu}$ suggests that it should be independent of the choice of
$x$, which indeed it is thanks to the invariance of $\mu$ under the action of $\Gamma$. Indeed, choose a point $y$ such that $y<f(x)<f(y)$ if
it exists, where these points are viewed as lying in the interval $[0,1)$.
Then, we compute: \[\mu[y,f(y))=\mu[y,f(x))+\mu[f(x),f(y)).\] Since \[\mu[f(x),f(y))=\mu[x,y),\] we obtain \[\mu[y,f(y))=\mu[x,f(x)).\] We leave
the verification of the general case to the reader. The interest in introducing $\rho_{\mu}$ comes from the following basic result which will
be used several times in the sequel.

\begin{thm}[cf.~\cite{Navas2011}, Theorem 2.2.10]\label{thm:invt-homo}
Let $\mu$ be an invariant probability measure for a group $\gam\le \Homeo_+(S^1)$. The the following hold.
\begin{enumerate}[(1)]
\item
For all $f\in \gam$, we have that $\rho_{\mu}(f)=\rot(f)$.
\item
The rotation number furnishes a homomorphism \[\rot\colon\gam\longrightarrow S^1.\]
\end{enumerate}
\end{thm}
\begin{proof}
We write \[\exp\colon [0,1)\longrightarrow S^1\] for the exponentiation map $x\mapsto e^{2\pi ix}$. The map $\exp$ has a measurable inverse,
and we write such an inverse $\exp^{-1}$.

Let $\mu$ be as in the statement of the theorem and let $f\in\gam$. First, renormalize $\mu$ to have total measure $1$. We pull back $\mu$
along $\exp$ (that is, push forward along $\exp^{-1}$) to obtain a measure $\exp^{-1}_*\mu$ on $[0,1)$, and we extend this measure
periodically to obtain a measure $\nu$ on $\bR$ that is invariant under an arbitrary lift $F$ of $f$. The reader may check that the measure
$\nu$ has the property that for all $x\in \bR$ and $k\in\bN$, we have \[\nu[x,x+ k)= k.\] 
We may thus make some straightforward estimates as
follows. If $F$ is a lift of $f$ and if \[F^n(x)\in [x+ k,x+ k+1)\] then we have \[F^n(x)-x-1\leq  k\leq \nu[x,F^n(x))\leq  k+
1\leq F^n(x)-x+1.\]
Dividing by $n$ and passing to the limit, we see that \[\lim_{n\to\infty}\frac{F^n(x)-x}{n}=\lim_{n\to\infty}\frac{\nu[x,F^n(x))}{n}.\]
Breaking up $[x,F^n(x))$ along points in the $F$--orbit of $x$, we see that the right hand side of this last equality coincides with
\[\lim_{n\to\infty}\frac{1}{n}\sum_{k=0}^{n-1}\nu[F^k(x),F^{k+1}(x))=\lim_{n\to\infty}\frac{1}{n}\cdot n\cdot\nu[x,F(x))=\nu[x,F(x)),\]
where here we
are using the $F$--invariance of $\nu$. It follows immediately now that $\rho_{\mu}(f)=\rot(f)$, which establishes the first claim of the theorem.

To see that the rotation number is a homomorphism when restricted to $\Gamma$, it suffices to see that
\[\rho_{\mu}(f\circ g)=\rho_{\mu}(f)+
\rho_{\mu}(g).\] Note that
\begin{align*}\rho_{\mu}(f\circ g)=\mu[x,fg(x))\pmod{\bZ}\\ \equiv\mu[x,g(x))+\mu[g(x),fg(x))\pmod{\bZ}\\ \equiv\rho_{\mu}(f)+
\rho_{\mu}(g)\pmod{\bZ},\end{align*} as required.
\end{proof}

The rotation number is not a homomorphism when restricted to an arbitrary subgroup of $\Homeo_+(S^1)$, as the reader may have
suspected. For an explicit example of a such a subgroup, we consider \[\PSL_2(\bR)\le \Homeo_+(S^1)\] where $\PSL_2(\bR)$ acts
on $S^1\cong \bR\bP^1$ by fractional linear transformations \[\begin{pmatrix}a&b\\c&d \end{pmatrix}\colon x\mapsto\frac{ax+b}{cx+d}.\] The
subgroup $\Gamma\le \PSL_2(\bR)$ generated by the matrices
\[\left\{\begin{pmatrix}2&0\\0&\frac{1}{2} \end{pmatrix},\begin{pmatrix}t&0\\0&t^{-1} \end{pmatrix},\begin{pmatrix}1&1\\0&1 \end{pmatrix}\right\},\]
where here $t\in\bR$ is arbitrarily chosen so that $\form{ \log 2,\log t}$ is a dense additive subgroup of $\bR$,
has the property that the three generators of $\Gamma$ each have fixed points in $S^1$ and hence have rotation number zero. If the rotation
number were a homomorphism when restricted to $\Gamma$ then the rotation number would be trivial for all elements of $\Gamma$.

However,
general facts about Lie groups imply that $\Gamma$ is in fact dense in $\PSL_2(\bR)$. Indeed, since $\Gamma$ is not virtually solvable (as it 
contains a nonabelian free group), it is Zariski dense in $\PSL_2(\bR)$. Since $\Gamma$ is not discrete in $\PSL_2(\bR)$ by the choice
of $t$, its topological closure is a Lie subgroup of positive dimension which must therefore be all of $\PSL_2(\bR)$. It follows that
$\Gamma$ meets every nonempty open subset of $\PSL_2(\bR)$, and in particular there is an element $\gamma\in\gam$
with trace contained in the
interval $(-2,2)$. It is not difficult to show that such an element $\gamma$ is conjugate to the Euclidean rotation
\[\begin{pmatrix}\cos\theta &-\sin\theta\\ \sin\theta&\cos\theta \end{pmatrix}\] for some $\theta\notin \bZ$, as
\[\tr(\gamma)=2\cos\theta\in (-2,2).\] It follows that the rotation number of $\gamma$ is nonzero, a contradiction.

So far, the results concerning invariant measures and rotation numbers have been conditional, in the sense that they assume the existence
of an invariant measure before they can be applied. Most subgroups of $\Homeo_+(S^1)$ do not admit invariant measures, though if one is
willing to impose some algebraic conditions, they will. In general if $G$ is a countable discrete group,
one says that $G$ is \emph{amenable group}\index{amenable group} if every
action of $G$ on a compact Hausdorff topological space admits an invariant probability measure. Examples of amenable groups include
abelian groups and solvable groups, and the class of amenable groups is closed under extensions, subgroups, quotients, and direct limits.
Amenability is not a primary concern of this book, and we will content ourselves to direct the interested reader
to~\cite{Lub94,Paterson88,Runde02,Zimmer84} for more on amenability.

We will prove amenability for solvable groups since the amenability of these groups, especially with respect to applying
Theorem~\ref{thm:invt-homo}, will arise in later chapters. Specifically, we will prove a variation on the classical Kakutani--Markov
Fixed Point Theorem.

\begin{thm}\label{thm:kaku-mar}
Let $G$ be a countable discrete solvable group acting by homeomorphisms on a compact, Hausdorff  space $X$. Then $G$
preserves a probability measure on $X$.
\end{thm}

The proof of Theorem~\ref{thm:kaku-mar} requires some functional analysis which we will simply assume here, directing the reader
to~\cite{Rudin87book,Zimmer-funct} for the relevant background.
Consider the space $C(X)$ of complex--valued continuous functions on $X$. 
The Riesz Representation
Theorem implies that the dual space $V=C(X)^*$ of bounded linear functionals on $C(X)$ is the space of complex regular 
Borel measures on $X$ of finite total measure,
and the space of positive linear
functionals on $C(X)$ consists precisely of
positive regular Borel measures on $X$, with the pairing being given by \[(\mu,f)\mapsto\int_X f\,d\mu.\] For our purposes, the following
statement is sufficient, which the reader may find as Theorem 2.14 in~\cite{Rudin87book}.

\begin{thm}[Riesz Representation Theorem]\label{thm:riesz}
Let $X$ be a compact Hausdorff topological space, and let $\Lambda$ be positive linear functional on $C(X)$ (i.e.~$\Lambda(f)\ge 0$
whenever $f(x)\geq 0$ for all $x$). Then there is a regular Borel measure $\mu$ such that $\mu(X)<\infty$ and such that
\[\Lambda(f)=\int_X f\,d\mu\] for all $f\in C(X)$.
\end{thm}

General (bounded) functionals on $C(X)$ can be recovered as complex linear combinations of measures on $X$.

The natural topology to consider on $V$ is the weak--$\ast$ topology, so that $\{\mu_n\}_{n\ge1}$ converges to $\mu$ if for every
$f\in C(X)$, we have \[\int_X f\,d\mu_n\longrightarrow\int_X f\,d\mu.\] The fundamental property of the weak--$\ast$ topology that
we will use is the following.

\begin{thm}[Banach--Alaoglu Theorem, Theorem 1.1.28 of~\cite{Zimmer-funct}]\label{thm:banach-al}
Let $E$ be a normed linear space. Then the unit ball in $E^*$ is compact in the weak--$\ast$ topology.
\end{thm}

We set $\PP(X)\sse V=C(X)^*$ to be the subset consisting of regular Borel probability measures on $X$. The subset $\PP(X)\sse V$ is compact
in the weak--$\ast$ topology, by the Banach--Alaoglu Theorem, since it is a closed subset of the unit ball.
Moreover, it is clear from its definition that $\PP(X)$ is convex. We remark briefly
that $\PP(X)$ is obviously nonempty if $X$ has at least one point, since a Dirac mass at a point of $X$ lies in $\PP(X)$.

Now let $G$ be a group acting on a compact Hausdorff space $X$ by homeomorphisms. We will always assume that $G$ is equipped with the discrete topology,
so that there are no issues with continuity when it comes to the associated actions on $V$ and $\PP(X)$. The group $G$ acts on the
space of all measures on $X$ by pushforward. Since $G$ acts continuously on $X$, the group $G$ preserves the Borel measures and the
regular measures. It is clear that $G$ preserves the set of positive regular Borel measures on $X$, and among those it preserves the
probability measures. Thus, $\PP(X)$ is $G$--invariant.

The following fact is sometimes called the Kakutani--Markov Fixed Point Theorem, when one allows $G$ to be a general abelian topological
group leaving invariant a compact, convex subset of a topological vector space whose topology is defined by a sufficient family of seminorms
on which $G$ acts continuously (see~\cite{Zimmer-funct}, for example).

\begin{lem}[cf.~\cite{Zimmer-funct}, Theorem 2.1.5]\label{lem:kaku-abel}
Let $G$ be an abelian group and $X$ be a compact Hausdorff space.
If $Y\sse\PP(X)$ is a nonempty, compact, convex, $G$--invariant set,
then there is a measure $\mu\in Y$ that is $G$--invariant. Moreover, the set $Y_G$ of $G$--invariant measures in 
$Y$ is convex and compact.
\end{lem}
\begin{proof}
For notational convenience, we write \[\phi\colon G\longrightarrow \Aut(V)\] for the action of $G$ on $V$. For $g\in G$ and $n\in\bZ_{>0}$,
we obtain operators
\[A_{n,g}=\frac{1}{n}\sum_{i=0}^{n-1}\phi(g^i),\] which are clearly continuous
endomorphisms of $V$. Since $Y$ is convex and invariant under the
action of $G$, it is immediate to check that $A_{n,g}$ preserves $Y$.

Consider the semigroup $S$ generated by \[\{A_{n,g}\mid n\in\bZ_{>0},\, g\in G\}.\] An arbitrary element $T\in S$ is simply a finite composition
of operators of the form $A_{n,g}$. Since $Y$ is invariant under each of endomorphisms, we have that $T(Y)\sse Y$ for
arbitrary elements $T\in S$. Note that since $G$ is an abelian group, the semigroup $S$ is commutative.

{\bf Claim 1:} We have that \[Y_G=\bigcap_{T\in S} T(Y)\neq\varnothing.\] 

It is clear that $Y_G$ coincides with the intersection above.
Since each $T$ is continuous and since $Y$ is compact, it suffices to
verify the finite intersection property characterization of compactness.
Since $Y$ is $S$--invariant, for $T_1,T_2\in S$ we have that
\[
T_1\circ T_2(Y)\sse T_1(Y)\]
and that
\[T_1\circ T_2(Y)=T_2\circ T_1(Y)\sse T_2(Y).\]
In general, if \[\{T_1,\ldots,T_n\}\sse S\] are arbitrary, then we have that
\[T_1\circ\cdots\circ T_n(Y)\sse\bigcap_i T_i(Y).\] It follows that $Y_G$ is nonempty as claimed.

{\bf Claim 2:} If $\mu\in Y_G$ then $\mu$ is $G$--invariant. 

By the very definition of $Y_G$, we have $\mu=s(\lambda)$ for some $s\in S$ and
 $\lambda\in Y$.
It suffices to consider the case that $s= A_{n,g}$.
  Observe then that \[\phi(g)(\mu)-\mu=\frac{1}{n}\left(\phi(g^n)(\lambda)-\lambda\right).\] Let $f\in C(X)$
with $\sup_{x\in X} |f(x)|\leq 1$.
Evidently, by the triangle inequality, we have that \[\left| \int_X f\,d(\phi(g)(\lambda))-\int_X f\,d\lambda\right| \leq \frac{2}{n},\] since $\lambda$
is a probability measure. Since $n$ is arbitrary, we have that \[\left| \int_X f\,d(\phi(g)(\mu))-\int_X f\,d\mu\right| =0\] for all continuous
functions $f$, so $\mu$ is $G$--invariant. 

{\bf Claim 3:} The set $Y_G$ is convex and compact. These properties of $Y_G$ follow easily from the definitions.
\end{proof}

We immediately obtain:
\begin{cor}\label{cor:homeo-invt-meas}
Let $f\in\Homeo_+(S^1)$. Then $f$ preserves a Borel probability measure on $S^1$.
\end{cor}

This gives an alternative definition of a rotation number, by the formula given in Theorem~\ref{thm:invt-homo}.
We can prove Theorem~\ref{thm:kaku-mar} for solvable groups by bootstrapping Lemma~\ref{lem:kaku-abel}, and inducting on the length
of the derived series of $G$.

\begin{proof}[Proof of Theorem~\ref{thm:kaku-mar}]
Let $D_0=G$ and $D_k=[D_{k-1},D_{k-1}]$ denote terms of the derived series of $G$. By assumption, $G$ is solvable and so for some
$n$ we have $D_{n+1}=\{1\}$. We proceed by induction on $n$. The set $Y_n$ consisting of $D_n$--invariant Borel probability measures on
$X$ is nonempty, compact, and convex by Lemma~\ref{lem:kaku-abel}, and this is the base case of the induction.

We assume that the set $Y_1$ consisting of $D_1$--invariant Borel probability measures on $X$ is nonempty, compact, and convex.

{\bf Claim 1:} The set $Y_1$ is $D_0$--invariant. Indeed, let $\mu\in Y_1$, let $g\in G$, and consider $\phi(g)(\mu)$.
This is clearly a Borel probability
measure on $X$. If $h\in D_1$ then
the normality of $D_1$ in $D_0$ implies that
\[\phi(h)(\phi(g)(\mu))=\phi(hg)(\mu)=\phi(g)(\phi(g^{-1}hg)(\mu))=\phi(g)(\mu),\] so that $\phi(g)(\mu)$
is $D_1$--invariant. This proves the claim.

{\bf Claim 2:} The set $Y_0$ of $D_0$--invariant Borel probability measures is nonempty. Indeed, the action of $D_0$ on $Y_1$
factors through the quotient $D_0/D_1$, since $D_1$ acts trivially on $Y_1$. We have that $Y_1$ is nonempty, compact, and convex,
and $D_0/D_1$ is abelian. Lemma~\ref{lem:kaku-abel} shows that the set $Y_0$ is nonempty.
\end{proof}

We are now ready to give a proof of the fact that every orientation preserving homeomorphism of $S^1$ with irrational rotation number is
semi-conjugate to an irrational rotation.

\begin{proof}[Proof of Theorem~\ref{thm:irr-rot-semi}]
We begin with a setup identical to that in Theorem~\ref{thm:invt-homo}.
Let $f$ admit an invariant probability measure $\mu$ as furnished by Corollary~\ref{cor:homeo-invt-meas} which as before we renormalize
to have total mass $1$,
and let \[\exp\colon [0,1)\longrightarrow S^1\] be the usual exponential map. The map
$\exp^{-1}$ is measurable. So, we let $\nu$ be the pushforward of $\mu$ by $\exp^{-1}$, and we extend $\nu$ periodically to all of $\bR$.

Since $g$ has no periodic points by assumption, $\mu$ has no atoms and hence every point
has measure zero. Clearly
$\nu$ is also atomless. Then by construction, we have $\nu(F(X))=\nu(X)$ for every lift $F$ of $f$ and every Borel set $X\sse \bR$. We
set \[H(x)=\int_0^x\, d\nu,\] so that $H(0)=0$ and $H$ is periodic with respect to translation by $1$. Since $\nu$ has no atoms, the function
$H$ is continuous, and since $\nu$ is a positive measure, we have that $H$ is non-decreasing.

We now claim that $H$ semi-conjugates $f$ to a rotation, which in light of Proposition~\ref{prop:semi-rot} must be rotation by $\rot(f)=\theta$.
Computing, we note that \[H(F(x))-H(x)=\int_0^{F(0)}\,d\nu+\int_{F(0)}^{F(x)}\,d\nu-\int_0^x\, d\nu.\] Since $\nu$ is $F$--invariant, the last
two terms cancel out, so that \[H(F(x))-H(x)=H(F(0))-H(0)\] for all $x$. It follows that $H\circ F$ is a translation of $\bR$, so that pushing
$H$ down to a self-map $h$ of $S^1$, we have that $h\circ f$ is a rotation.

Note that if $h$ is a homeomorphism then $f$ must be minimal since minimality is preserved under conjugacy. Conversely, if $f$ is minimal, let
$x\in S^1$ be a point such that every neighborhood of $x$ has positive measure with respect to $\mu$. Such a point exists since otherwise
$S^1$ would admit a finite cover by open sets with zero measure. We have that the orbit of $x$ is dense, so that every
nonempty open set has positive
measure. It follows then that the map $H$ is strictly increasing, so that $h$ must be a homeomorphism.
\end{proof}

\subsection{Invariant measures and free subgroups}
To close this section,
we mention one more result about invariant measures due to Margulis, though we will not give a proof here. Recall that a group satisfies the
\emph{Tits alternative}\index{Tits alternative} if every subgroup is either virtually solvable or contains a 
nonabelian free group. The terminology arises from Tits'
famous result that linear groups enjoy the Tits alternative~\cite{tits-jalg}.

Homeomorphism groups of one--manifolds do not enjoy the Tits alternative, though we are not quite ready to show this yet. Thompson's
group $F$ is a group of homeomorphisms of the interval which is neither virtually solvable nor contains a nonabelian free subgroup; see
Theorem~\ref{thm:f-subgp} below. A version of the Tits alternative for homeomorphisms of the circle was established by Margulis:

\begin{thm}[See~\cite{marg-CR00,Ghys2001} and also Theorem 2.3.2 in~\cite{Navas2011}]\label{thm:marg-ghys}
Let \[G\le \Homeo_+(S^1)\] be a subgroup. Then either $G$ preserves a probability measure on $S^1$, or $G$ contains a nonabelian free
subgroup.
\end{thm}

Theorem~\ref{thm:marg-ghys} is only interesting for groups acting on the circle without global fixed points, since a Dirac mass concentrated
at a global fixed point is an invariant probability measure. Also, whereas the Tits alternative is a true dichotomy, Theorem~\ref{thm:marg-ghys}
is not. There are many interesting actions by groups containing nonabelian free subgroups that preserve nonatomic probability measures,
and this is closely related to the theory of \emph{Conradian orderings}\index{Conradian ordering},
which we will not discuss in detail in this book. The interested
reader may begin by consulting~\cite{navasrivas09}, for instance.

\section{Denjoy's theorem}
In this section, we prove the following result of Denjoy:

\begin{thm}[Denjoy's Theorem]\label{thm:denjoy}
Let $f\in\Diffb(S^1)$, and suppose that $\rot(f)=\theta$ is irrational. Then $f$ is conjugate to the rotation $R_{\theta}$.
\end{thm}

Thus, Denjoy's Theorem relates exceptional minimal sets, rotation numbers, and analysis in a decisive way. The moment $f$ is a
diffeomorphism and the derivative of $f$ has bounded variation (a conclusion which is satisfied, for example, when $f$ is twice differentiable),
then $f$ cannot be an exceptional diffeomorphism. We emphasize that Denjoy's Theorem is false if $f$ is merely assumed to be a
diffeomorphism, so that smooth Denjoy counterexamples exist. The exact regularities in which there exist Denjoy counterexamples is a subtle
and complicated question, and is the subject of Section~\ref{sec:denjoy-crit} below. The account we follow here uses an equidistribution
phenomenon called the Denjoy--Koksma inequality, which combines with certain derivative estimates to yield Denjoy's Theorem. Another
approach based on Sacksteder's Theorem and the existence of hyperbolic fixed points can be found in~\cite{Navas2011}.

\subsection{Rational approximations of real numbers}
If $f\in\Homeo_+(S^1)$ has rational rotation number, then clearly $f$ need not be conjugate to a rational rotation, regardless of the level of
regularity one imposes on $f$. Thus, one must use irrationality in some essential way. Indeed, rational approximation of arbitrary
real numbers by rational numbers plays a key role in the proof of Denjoy's Theorem. Specifically, we have the following classical result:

\begin{thm}[Dirichlet's Diophantine Approximation Theorem]\label{thm:dirichlet-1}
Let \[\theta\in\bR\setminus\bQ\] and $N\in\bZ_{>0}$. There exist relatively prime integers $p$ and $q$ with $1\leq q<N$ such that
\[|q\cdot\theta-p|<\frac{1}{N}.\]
\end{thm}

The rational number $p/q$ is called a \emph{rational approximation}\index{rational approximation} of $\theta$.
We restrict to irrational $\theta$, since the result for rational values of $\theta$ is not relevant to our purposes, and because the proof is
easier to the point of being almost trivial.

\begin{proof}[Proof of Theorem~\ref{thm:dirichlet-1}]
This is just a consequence of the pigeonhole principle. We consider the set of real numbers \[S=\{0,1,\{n\theta\}\mid 1\leq n\leq N-1\},\]
where here $\{x\}$ denotes the fractional part of $x$, the unique real number $\eps\in [0,1)$ such that $x-\eps$ is an integer.

Clearly $|S|=N+1$ and $S\sse [0,1]$. Cutting $[0,1]$ into disjoint intervals of length $1/N$, the pigeonhole principle implies that at least
one such interval meets $S$ in two points. Observe that since $\theta$ is not rational, none of the points \[\{k/N\mid 1\leq k\leq N-1\}\] meet
the set $S$, so that the distance between two such points is strictly less than $1/N$.
If $s_1,s_2\in S$ are two such points then \[s_1-s_2=q\cdot\theta-p\] for a suitable choice of integers
$p$ and $q$, where $0<q<N$. It follows that \[|q\cdot\theta-p|<\frac{1}{N}.\] If $p$ and $q$ fail to be relatively prime then we may simply divide
through by the greatest common divisor, thus obtaining the claimed result.
\end{proof}

We will use the following version of Dirichlet's theorem:

\begin{cor}\label{cor:dirichlet}
Let $\theta\in\bR\setminus\bQ$. There exists a sequence $\{(p_n,q_n)\}_{n\ge1}$ with $p_n\in\bZ$ and $q_n\in\bZ_{>0}$ for all $n$, such that:
\begin{enumerate}[(1)]
\item
For all $n$, we have $p_n$ and $q_n$ are relatively prime.
\item
We have \[\lim_{n\to\infty} q_n=\infty.\]
\item
We have \[\left| \theta-\frac{p_n}{q_n}\right|<\frac{1}{q_n^2}.\]
\end{enumerate}
\end{cor}
\begin{proof}
Let $N_1\in\bZ_{>0}$ be arbitrary, and let $(p_1,q_1)$ be a pair as furnished by Theorem~\ref{thm:dirichlet-1}, so that \[\left|\theta-\frac{p_1}{q_1}
\right|<\frac{1}{N_1\cdot q_1}<\frac{1}{q_1^2}.\] We then choose a natural number \[N_2>\frac{1}{|\theta-p_1/q_1|}.\]
Theorem~\ref{thm:dirichlet-1} furnishes a new pair $(p_2,q_2)$,
with $N_2$ as the input parameter. Note that we may estimate \[\left|\theta-\frac{p_2}{q_2}\right|<\frac{1}{N_2\cdot
q_2}<\left|\theta-\frac{p_1}{q_1}\right|\cdot\frac{1}{q_2}\leq \left|\theta-\frac{p_1}{q_1}\right|.\] It follows that $p_2/q_2$ is a strictly better
approximation of $\theta$ that $p_1/q_1$. Thus, we may recursively produce a sequence $\{(p_n,q_n)\}_{n\ge1}$ as claimed, and since
for all $i$ we have $p_{i+1}/q_{i+1}$ is a better approximation to $\theta$ than $p_i/q_i$, it follows that $q_n\to\infty$ as $n\to\infty$.
\end{proof}

The usefulness of rational approximations is that they give some control over the distribution of orbits of an irrational rotation for powers
that do not exceed the denominator of the approximating rational number. 
Again, let \[\exp\colon [0,1]\longrightarrow S^1\] be the exponentiation map
$\exp(t)=e^{2\pi i t}$.

\begin{cor}\label{cor:dirichlet-circle}
Let $\theta\in (0,1)$ be irrational, and let $p/q\in (0,1)$ be a rational approximation of $\theta$. Consider the sequence of points
\[S_q=\{\exp( k\theta)\mid 1\leq k\leq q\}.\] Then for all $0\leq j<q$ we have \[\left|S_q\cap \exp\left(\left[\frac{ j}{q},\frac{(j+1)}{q}\right]
\right)\right|=1.\]
\end{cor}
\begin{proof}
First, observe that \[e^{ k\theta}\in \exp\left(\left[\frac{ kp}{q},\frac{(k+1)p}{q}\right]\right)\quad \textrm{or}\quad e^{ k\theta}\in 
\exp\left(\left[\frac{ (k-1)p}{q},\frac{ kp}{q}\right]\right),\] according to whether $\theta$ is greater than or smaller than $p/q$,
since we have \[\left|k\theta-\frac{kp}{q}\right|<\frac{k}{q^2}<\frac{1}{q}.\] The intervals
\[\exp\left(\left[\frac{ k}{q},\frac{(k+1)}{q}\right]
\right)\] have disjoint interiors, as $k$ varies between $0$ and $q-1$, since $p$ and $q$ are relatively prime, and each contains exactly
one point in $S_q$. This establishes the corollary.
\end{proof}

\subsection{The Denjoy--Koksma inequality and unique ergodicity}\label{sss:denjoy-koksma}

A key observation about about
homeomorphisms with irrational rotation number is the Denjoy--Koksma inequality, which says that if $\phi$ is a reasonably nice function
on $S^1$, then averaging the values of $\phi$ on orbits of the homeomorphism converges to the average value of $\phi$. The reader
familiar with ergodic theory will notice a resemblance with Birkhoff's Ergodic Theorem~\cite{EW11book,petersen83}.
Here, ``reasonably
nice" means \emph{bounded variation}\index{bounded variation}.
Recall that a function $\phi\colon S^1\longrightarrow \bR$ has bounded variation if
\[
\Var(\phi;S^1):=\sup_{n\ge1}\sum_{\{x_0,\cdots,x_n\}\sse S^1} |\phi(x_i)-\phi(x_{i-1})|\] is finite, where here
\[x_0<x_1<\cdots <x_n<x_0,\] where the indices are considered
modulo $n$, and where the ordering is the cyclic ordering on $S^1$ determined by a choice of orientation.

\begin{thm}[Denjoy--Koksma inequality]\label{thm:denjoy-koksma}
Suppose $f\in\Homeo_+(S^1)$ with irrational rotation number $\theta$. Let $p/q$ be a rational approximation of \[{\theta}\in (0,1),\]
and let $\phi\colon S^1\longrightarrow\bR$ be a function of bounded variation. If $x\in S^1$ is arbitrary and if $\mu$
is an $f$--invariant Borel probability measure, then, we have
\[\left|\frac{1}{q}\sum_{k=0}^{q-1}\phi(f^k(x))-\int_{S^1}\phi\,d\mu\right|\leq \frac{\Var(\phi;S^1)}{q}.\]
\end{thm}

Even though the Denjoy--Koksma inequality superficially resembles the Birkhoff Ergodic Theorem, we note that the factor of $q$ is of
crucial importance in the application of the inequality. Thus, the Denjoy--Koksma inequality should be viewed as
having both number--theoretic
and dynamical content.

\begin{proof}[Proof of Theorem~\ref{thm:denjoy-koksma}]
We follow the argument in~\cite{athanassopoulos}.
Theorem~\ref{thm:irr-rot-semi} furnishes a semi-conjugacy $h$ between $f$ and the rotation $R_{\theta}$ by $\theta$. Observe that
the pushforward $\nu=h_*\mu$ of $\mu$ by $h$ is an invariant measure for $R_{\theta}$, and thus
by Corollary~\ref{cor:leb-ue} must be Lebesgue measure, which we normalize for the purpose of this proof to have total mass one.

Choose an arbitrary point $z_0\in S^1$, let \[\{z_0,z_1,\ldots,z_{q-1}\}\sse S^1\] be the orbit of $z_0$ under the rotation by
$1/q$, and let \[\{x_0,\ldots,x_{q-1}\}\sse S^1=h^{-1}(S^1)\]
be points such that $x_i\in h^{-1}(z_i)$ for all $i$. In the circular order on $S^1$, we may
clearly assume that \[x_0<x_1<\cdots<x_{q-1}<x_0,\] and so we may consider the oriented interval $[x_i,x_{i+1}]$, where the indices are
again regarded modulo $q$. By definition, we have \[\int_{[x_i,x_{i+1}]}\,d\mu=\int_{[z_i,z_{i+1}]}\,d\nu=\frac{1}{q}.\] By Corollary
\ref{cor:dirichlet-circle}, for all $1\leq k\leq q$, there is a unique $i=i(k)$ such that \[R_{\theta}^k(h(x_0))\in [z_i,z_{i+1}].\] We may thus
write \[\left|\frac{1}{q}\sum_{k=1}^{q}\phi(f^k(x_0))-\int_{S^1}\phi\,d\mu\right|=
\left|\sum_{k=1}^{q}\left(\frac{1}{q}\phi(f^k(x_0))-\int_{[x_{i(k)},x_{i(k)+1}]}\phi\,d\mu\right)\right|.\] Writing $c_k=\phi(f^k(x_0))$
and $J_k=[x_{i(k)},x_{i(k)+1}]$, this last
expression is bounded above by \[\sum_{k=1}^q\left|\int_{J_k}(c_k-\phi)\,d\mu\right|\leq
\frac{1}{q}\sum_{k=1}^q\sup_{x,y\in J_k}|\phi(x)-\phi(y)|.\] The right hand side of this
last expression is clearly bounded above by $\Var(\phi;S^1)/q$, whence the
inequality follows after setting $x=f^{-1}(x_0)$.
\end{proof}

The Denjoy--Koksma inequality has the following corollary, which was promised after Proposition~\ref{prop:irr-rot-erg}.
\begin{thm}\label{thm:irr-ue}
Let $f\in\Homeo_+(S^1)$ have irrational rotation number. Then $f$ is uniquely ergodic.
\end{thm}
\begin{proof}
Let $\theta$ be the rotation number of $f$ and let \[\frac{p_n}{q_n}\longrightarrow{\theta}\] be a sequence of rational approximations
such that $q_n\to\infty$ as $n\to\infty$. Theorem~\ref{thm:denjoy-koksma} implies that if $\phi$ has bounded variation and $\mu$
is an $f$--invariant Borel probability measure, then
\[\lim_{n\to\infty}\frac{1}{q_n}\sum_{k=0}^{q_n-1}\phi\circ f^k=\int_{S^1}\phi\,d\mu,\] where this convergence is uniform on $S^1$. If
$\lambda\neq\mu$ were another $f$--invariant Borel probability measure, then there would be a Borel subset $A\sse S^1$ such that
$\lambda(A)\neq\mu(A)$. The classical result Lusin's Theorem from measure theory~\cite{Rudin87book}
implies that there is a continuous function $\phi$
(which can easily be arranged to have bounded variation) approximating the characteristic function of $A$ and
such that \[\int_{S^1} \phi\,d\mu\neq\int_{S^1}\phi\,d\lambda.\] This is clearly a contradiction.
\end{proof}

\subsection{Completing the proof of Denjoy's Theorem}

We now complete the last steps in proving Denjoy's Theorem. In this subsection, we always assume that $f\in\Diff_+^1(S^1)$ is an
orientation preserving diffeomorphism. The next lemma says that the average value of the logarithm of the derivative of a diffeomorphism
is zero.

\begin{lem}[cf.~\cite{athanassopoulos}, Proposition 3.3]\label{lem:log-int}
Let $f$ have an irrational rotation number and let $\mu$ be the unique $f$--invariant Borel probability measure. Then
\[\int_{S^1}\log (f')\,d\mu=0.\]
\end{lem}
\begin{proof}
Note that \[\log((f^n)')=\sum_{k=0}^{n-1}(\log (f'))\circ f^k,\] as follows from the chain rule. It follows that \[\lim_{n\to\infty}\frac{1}{n}\log((f^n)')=
\int_{S^1}\log (f')\,d\mu,\] where the equation is valid when the left hand side is evaluated at an arbitrary point of $S^1$.
Indeed, since $f$ is uniquely ergodic with respect to $\mu$ by Theorem~\ref{thm:irr-ue}, the sequence of probability measures given by
measuring the (renormalized) number of points in the $f$--orbit of an arbitrary
point that meet a given subset of $S^1$ has a unique weak--$\ast$
limit that must be $\mu$ (see the remarks after Corollary~\ref{cor:leb-ue} above). The left hand side of the equation is
the limit of a sequence of
integrals of $\log(f')$ with respect to these measures.

Now, suppose that \[\int_{S^1}\log(f')\,d\mu>0.\] We must have that $(f^n)'$ tends to $+\infty$, uniformly on $S^1$.
If $F$ represents an arbitrary lift of $f$ to $\bR$ then \[\int_{S^1}(f^n)'(\theta)\,d\lambda=\int_0^{1}(F^n)'(x)\,dx\longrightarrow +\infty,\] where
$\lambda$ denotes Lebesgue measure on $S^1$ and where $dx$ denotes Lebesgue measure on $\bR$.

If, on the other hand \[\int_{S^1}\log(f')\,d\mu<0,\] we have that $(f^n)'$ tends to $0$, uniformly on $S^1$, so that
\[\int_{S^1}(f^n)'(\theta)\,d\lambda=\int_0^{1}(F^n)'(x)\,dx\longrightarrow 0.\] However, $f^n$ is a diffeomorphism, and so the integral
over $S^1$ of its derivative is simply $1$. In either of the previous two cases, we obtain a contradiction.
\end{proof}

The following is the last crucial observation before the proof of Denjoy's Theorem.

\begin{lem}[cf.~\cite{athanassopoulos}, Proposition 3.5]\label{lem:denjoy-var}
Suppose $\rot(f)=\theta$ is irrational, that $p/q$ is a rational approximation of $\theta$, and that $\log(f')$ has bounded variation $V$.
Then \[|\log((f^q)')|\leq V.\]
\end{lem}
\begin{proof}
The Denjoy--Koksma inequality shows that for $x\in S^1$ arbitrary, we have \[\left|\sum_{k=0}^{q-1}\log(f'(f^k(x)))-q\cdot\int_{S^1}\log(f')\,d\mu
\right|\leq V.\]
By Lemma~\ref{lem:log-int}, the second term on the left hand side is zero. The chain rule says that \[\sum_{k=0}^{q-1}\log(f'(f^k(x)))=
\log((f^q)').\] The lemma is now immediate.
\end{proof}

The conclusion of Lemma~\ref{lem:denjoy-var} is equivalent to  the statement that \[e^{-V}\leq (f^q)'\leq e^V.\]
Finally, we can give a proof of Denjoy's Theorem.

\begin{proof}[Proof of Theorem~\ref{thm:denjoy}]
If $f$ is minimal then $f$ is conjugate to the rotation $R_{\theta}$, by Theorem~\ref{thm:irr-rot-semi}. So, we may assume that $f$ is not
minimal. We have that $f$ has no periodic points by Proposition~\ref{prop:rot-easy}, and so Theorem~\ref{thm:minimal-set} implies
that there is a minimal invariant Cantor set $C\sse S^1$. If $J\sse S^1\setminus C$ is a component the $J$ is an open interval, and
we have $f^n(J)\cap J=\varnothing$ for $n\neq 0$.

Let \[\frac{p_n}{q_n}\longrightarrow{\theta}{}\]
be a sequence of rational approximations, where $q_n\to\infty$ as $n\to\infty$, as furnished by
Corollary~\ref{cor:dirichlet}. Clearly, we may assume that $q_n<q_m$ for $n<m$.
Lemma~\ref{lem:denjoy-var} bounds the derivative of $f^{q_n}$ away from zero, 
independently of anything other than the variation $V$ of $\log (f')$.
Therefore, the Mean Value Theorem implies that \[\frac{|f^{q_n}(J)|}{|J|}\geq e^{-V},\] where $|J|$ denotes the length of the interval $J$.
We finally obtain \[+\infty=\sum_{n\ge1}e^{-V}|J|\leq \sum_{n\ge1}|f^{q_n}(J)|\leq 1,\] which is a contradiction.
\end{proof}

\section{Exceptional diffeomorphisms and integrable moduli}\label{sec:denjoy-crit}

This section is devoted to investigating the failure of Denjoy's Theorem in the case of elements $f\in\Homeo_+(S^1)$ that do not lie in
$\Diffb(S^1)$. We have already explicitly constructed elements of $\Homeo_+(S^1)$ that have no periodic points and exceptional minimal
sets, which we called continuous Denjoy counterexamples. These homeomorphisms are defined by equivariant families of affine
homeomorphisms on the complement of the exceptional minimal set,
and so it might be rather difficult to force a continuous Denjoy counterexample to have desirable regularity properties.
Indeed, Denjoy's Theorem implies that a continuous Denjoy counterexample cannot have a first derivative of bounded variation.
If $f$ is a $C^2$ diffeomorphism then $f'$ and $(f^{-1})'$ are both differentiable, and so that the total variation of $\log (f')$
for example is just
\[\int_{S^1}\left|\frac{f''}{f'}\right|\,d\lambda,\] where $\lambda$ is Lebesgue measure normalized to have total mass $1$.
This integral is finite since $f'$ is continuous and therefore
is bounded away from zero on $S^1$.

We will show that many (but not all) regularities that are weaker than $C^2$ do admit Denjoy counterexamples, which is to say exceptional
diffeomorphisms with prescribed irrational rotation number.

\subsection{Stationary measures and Lipschitz homeomorphisms}
This subsection will largely follow Section 2.3.2 in~\cite{Navas2011}. We reproduce just the details we need in order to give a complete
account of critical regularity in the sequel.

Recall that $f\in\Homeo_+(S^1)$ (or more generally, a real valued function on a metric space)
is called \emph{Lipschitz}\index{Lipschitz map} if there is a constant $C$ such that for all $x,y\in S^1$, we have
\[|f(x)-f(y)|\leq C|x-y|,\] where here as usual the circle is considered as the additive reals modulo $\bZ$.
If an identical inequality also holds for $f^{-1}$, we say that $f$ is \emph{bi-Lipschitz}\index{bi-Lipschitz map}.
The infimum of such constants $C$ for a function $f$
is called the \emph{Lipschitz constant}\index{Lipschitz constant} of $f$.
It is not difficult to show that the set of bi-Lipschitz homeomorphisms of $S^1$ forms a subgroup of $\Homeo_+(S^1)$, as the
Lipschitz constants just multiply under composition. 

It is not difficult to produce a continuous Denjoy counterexample that is bi-Lipschitz, simply by choosing the lengths of the wandering
intervals $\{J_n\}_{n\in\bZ}$ with sufficient care. We will prove a fact originally observed by
Deroin--Kleptsyn--Navas~\cite{DKN2007} that establishes a much more
general fact: an arbitrary countable subgroup $\gam\le \Homeo_+(S^1)$ is topologically conjugate to group of bi-Lipschitz homeomorphisms
of $S^1$. This level of generality is not necessary for us to produce bi-Lipschitz exceptional diffeomorphisms,
but the existence of a such a conjugacy is an important part of the story of critical regularity of groups
that we will recount later in this book.

\begin{thm}[See~\cite{Navas2011}, Proposition 2.3.15; cf.~\cite{Calegari2007}, Theorem 2.118]\label{thm:lip-conj}
Suppose $\gam\le \Homeo_+(S^1)$ is countable.
Then there is an element $\phi\in\Homeo_+(S^1)$ such that $\phi^{-1}\gam \phi$ is a group of bi-Lipschitz
homeomorphisms of $S^1$.
\end{thm}

Before proving Theorem~\ref{thm:lip-conj}, we develop some basic facts about stationary measures. Let $\bP$ be a probability measure
on $\Gamma$, where $\Gamma$ is a countable group. We write $\supp\bP$ for the \emph{support}\index{support of a probability
measure on a group}
of $\bP$, which is the set
\[\supp\bP=\{\gamma\in\gam\mid\bP(\gam)>0\}.\] Recall that the support of a Borel measure on a topological space $X$
is typically defined to be the largest closed subset $C$ of $X$ such that every open subset of $C$ has positive measure. To avoid inconsistent
definitions, we will implicitly endow $\Gamma$ with the discrete topology.

To avoid the degeneracy of a probability measure that is essentially supported on
a proper subgroup of $\Gamma$, one usually requires $\bP$ to be \emph{nondegenerate}\index{nondegenerate probability measure}. 
That is, $\form{\supp\bP}=\gam$. The group
action of $\Gamma$ on $S^1$ and the measure $\bP$ give rise to a bounded operator (which we restrict to real valued continuous functions)
\[\Delta\colon C(S^1)\longrightarrow C(S^1)\] defined by
\[\Delta(f)(x)=\int_{\gam}f(\gamma(x))\,d\bP(\gamma).\] The intuitive meaning of $\Delta$ is clear: it averages the values of a function $f$
among its orbit under $\Gamma$, weighted by the measure $\bP$. Since the measure $\bP$ is a probability measure, one can think
of $\Delta$ as modeling a random evolution of functions under a $\Gamma$ action, and this
is why $\Delta$ is often called a \emph{diffusion operator}\index{diffusion operator}.

To proceed, we consider the dual operator to $\Delta$, which is defined as the formal adjoint
\[\Delta^*\colon C(S^1)^*\longrightarrow C(S^1)^*,\] acting on bounded linear functionals on $C(S^1)$. 
By the Riesz Representation Theorem (Theorem~\ref{thm:riesz}), we have that the space of positive functionals on $C(S^1)$
coincides with the space
of positive Borel measures on $S^1$ with finite total measure. So, if $\mu$ is such a measure on $S^1$ then $\Delta^*\mu$ is also
a measure on $S^1$, defined by \[\int_{S^1}f\,d(\Delta^*\mu)=\int_{S^1}\Delta(f)\,d\mu=\int_{S^1}\int_{\gam}f(\gamma(x))\,d\bP(\gamma)d\mu.\]
The measure $\Delta^*\mu$ is called the \emph{convolution}\index{convolution of measures} 
of $\bP$ and $\mu$, and is also written $\bP*\mu$. A measure $\mu$
on $S^1$ such that $\bP*\mu=\mu$ is called \emph{stationary}\index{stationary measure} for the data of $\Gamma$ and $\bP$.

One can make several immediate observations about convolution. First, suppose $f\in C(S^1)$ is a positive function, i.e.~$f(x)\geq 0$ for
all points $x\in S^1$. It is immediate that
\[\int_{S^1}\Delta(f)\,d\mu\geq 0,\] since at every point $x\in S^1$ and for every $\gamma\in\gam$, the value of $f(\gamma(x))$ is non-negative.
It follows that convolution with $\bP$ preserves the set of positive Borel measures of finite total mass on $S^1$. Integrating the
constant function $f(x)\equiv
1$ against $\Delta^*\mu$ shows that $\Delta^*$ preserves the set of Borel probability measures on $S^1$.

\begin{lem}\label{lem:stat-meas}
Let $\{\gam,\bP\}$ be as above. Then there exists a Borel probability measure $\mu$ on $S^1$ that is stationary for $\bP$.
\end{lem}
\begin{proof}
The convolution operator preserves the nonempty, convex, and compact set of Borel probability measures on $S^1$,
by the remarks preceding the lemma. 
To prove that there exists a fixed point of the convolution operator, one can simply follow
the proof of Lemma~\ref{lem:kaku-abel} above. Instead of using a group $G$, we simply use a semigroup $S$ generated by a single element
(namely $\Delta^*$).
The proof goes through without any difficulty, and we leave the details to
the reader.
\end{proof}

The support of a stationary measure\index{support of a measure}
$\mu$ for $\{\gam,\bP\}$ turns out to be closely related to the dynamical trichotomy given in
Theorem~\ref{thm:minimal-set}.

\begin{lem}[See~\cite{Navas2011}, Lemma~2.3.14]\label{lem:stat-support}
Let $\{\gam,\bP\}$ be as above, with $\bP$ nondegenerate,
and let $\mu$ be a stationary probability measure on $S^1$ for this data. Then the following conclusions hold.
\begin{enumerate}[(1)]
\item
If $\mu$ has an atom then $\Gamma$ has a finite orbit.
\item
If $\Gamma$ acts minimally then $\mu$ has full support; that is, $\mu$ gives positive measure to every open subset of $S^1$.
\item
If $\Gamma$ admits an exceptional minimal set $C$ then $\mu(C)=1$.
\end{enumerate}
\end{lem}
\begin{proof}
Suppose first that $x\in S^1$ is a point with $\mu(x)>0$, and where $x$ has maximal mass.
Since $\mu$ is stationary, we have that
\[\mu(x)=\int_{\gam}\mu(\gamma^{-1}(x))\,d\bP(\gamma).\] 
If $\mu(\gamma^{-1}(x))<\mu(x)$ for some 
$\gamma$
then we obtain \[\sum_{\gamma\in\gam}\mu(\gamma^{-1}(x))\bP(\gamma)<\mu(x),\] a contradiction.
Thus, $\mu(\gamma^{-1}(x))=\mu(x)$ for $\gamma\in\supp\bP$. Since $\supp\bP$ generates
all of $\Gamma$, we see that $\mu$ is constant on the $\Gamma$--orbit of $x$. Since $\mu$ is a probability measure, it follows that the $\Gamma$--orbit
of $x$ is finite.

The support of $\mu$ is a closed subset of $S^1$ by definition, and it is also $\Gamma$--invariant, as is a straightforward consequence of the
stationarity of $\mu$ and the non-degeneracy of $\bP$. It follows that if $\Gamma$ acts minimally on $S^1$ then $\supp\mu=S^1$.

Now, suppose that $\Gamma$ admits an exceptional minimal set $C$. Then $C\sse\supp\mu$, by 
the minimality of $C$. It suffices to show that
if $J\sse S^1\setminus C$ is an interval then $\mu(J)=0$.

Assume for a contradiction that $\mu(J)>0$. Without loss of generality, $J$ is a maximal connected component of the
complement of $C$, so that for
all $\gamma\in\gam$ we either have \[\gamma(J)\cap J=\varnothing\quad \textrm{ or }\quad \gamma(J)=J.\] As in the case of an atomic
measure, we have \[\mu(J)=\int_{\gam}\mu(\gamma^{-1}(J))\,d\bP(\gamma),\] by stationarity. It follows that the $\Gamma$--orbit of $J$ is
finite. However, if $x\in\partial J$ then the $\Gamma$--orbit of $x$ is dense in $C$, which is the desired contradiction.
\end{proof}

In the case of a minimal action, it follows that the stationary measure $\mu$ can be pushed forward to Lebesgue measure by a circle
homeomorphism:

\begin{lem}\label{lem:leb-push}
Let $\mu$ be a probability measure on $S^1$ with no atoms and whose support coincides with $S^1$.
Then there is a homeomorphism $\phi\in\Homeo_+(S^1)$
such that $\phi_*\mu$ is Lebesgue measure.
\end{lem}
\begin{proof}
As we have seen several times already, the measure $\mu$ induces a periodically defined measure $\nu$ on $\bR$. One can define
a homeomorphism of $\bR$ by
\[\Phi(x)=\int_0^x \,d\nu.\] 
Since $\mu$ has no atoms, neither does $\nu$, and $\Phi$ is continuous. Since $\mu$ is fully supported and
since an open interval separates every pair of distinct points in $\bR$, the map $\Phi$ is injective. Since $\Phi$ is periodic, it is in fact a
homeomorphism and descends to a homeomorphism $\phi$ of $S^1$. 
The homeomorphism $\phi$ pushes $\mu$ forward to Lebesgue
measure.
\end{proof}

We can now prove that an arbitrary countable group action on $S^1$ is conjugate to a bi-Lipschitz action.

\begin{proof}[Proof of Theorem~\ref{thm:lip-conj}]
First, note that one may assume $\Gamma$ acts on $S^1$ minimally. If not, then one can add an irrational rotation to $\Gamma$ to obtain a larger
(but still countable) group of homeomorphisms.

We set $\bP$ to be a non-degenerate probability measure on $\Gamma$ that is symmetric, which is to say $\bP(\gamma)=\bP(\gamma^{-1})$
for all $\gamma\in\gam$. Lemma~\ref{lem:stat-meas} says that there exists a stationary measure $\mu$ on $S^1$ associated to
$\Gamma$ and $\bP$. Writing $\supp\bP\sse\gam$ for the support of the measure $\bP$, the stationarity of $\mu$ implies that
\[\sum_{\gamma\in\supp\bP}\mu(\gamma^{-1}(J))\bP(\gamma)=\mu(J)\] for all intervals $J\sse S^1$. For an arbitrary $\gamma\in \gam$,
we obtain that \[\mu(\gamma(J))\bP(\gamma^{-1})\leq \mu(J),\] and so
since $\bP$ is assumed to be symmetric, we get \[\mu(\gamma(J))\leq\frac{\mu(J)}{\bP(\gamma)},\] provided $\gamma$
lies in the support of $\bP$.

By Lemma~\ref{lem:stat-support} and Lemma~\ref{lem:leb-push},
we have that there is an element $\phi\in\Homeo_+(S^1)$ such that $\phi_*\mu=\lambda$,
where $\lambda$ denotes Lebesgue measure. We claim that $\phi$ conjugates $\Gamma$ to a group of bi-Lipschitz homeomorphisms.
To this end, let $\gamma\in\supp\bP\sse \gam$, let $C=\bP(\gamma)^{-1}$, and let $J\sse S^1$ be an interval. We estimate:
\[|\phi\circ\gamma\circ\phi^{-1}(J)|=\mu(\gamma\circ\phi^{-1}(J))\leq C\mu(\phi^{-1}(J))=C|J|,\] where here as before the absolute value
denotes Euclidean length. The first and last equalities are justified since $\phi$ pushes $\mu$ forward to Lebesgue measure, and the 
inequality is justified in the previous paragraph. It follows that $\phi$ conjugates a generating set for $\Gamma$ to a set of bi-Lipschitz
homeomorphisms of $S^1$, which establishes the result.
\end{proof}

We remark that a verbatim analogue of Theorem~\ref{thm:lip-conj} holds for the interval. Since $\Homeo_+[0,1]$ can be embedded in
$\Homeo_+(S^1)$ by identifying the endpoints of the interval, it is immediate that one can conjugate countable subgroups of $\Homeo_+[0,1]$
to bi-Lipschitz homeomorphisms. Some care needs to be taken to show that the conjugacy can be realized within $\Homeo_+[0,1]$. We leave
the details to the reader.

\subsection{Smooth Denjoy counterexamples and the spectrum of their moduli of continuity}\label{ss:smooth-denjoy}

We now concentrate on exceptional diffeomorphisms $f\in\Diff_+^1(S^1)$. By definition, the derivatives of $f$ and $f^{-1}$ are
continuous functions. Bohl~\cite{Bohl1916} and Denjoy~\cite{Denjoy1932} have constructed exceptional diffeomorphisms of the circle,
and their results were extended significantly by Herman~\cite{Herman1979}.
In light of Denjoy's Theorem, there is a dividing line (or dividing lines) in the regularity properties of the derivatives $f'$ and $(f^{-1})'$
that distinguish between the possibility of exceptional diffeomorphisms and the proscription thereof.

To introduce a precise setup, let \[\alpha\colon [0,\infty)\longrightarrow [0,\infty)\] be a homeomorphism that is concave as a function.
Sometimes we will refer to $\alpha$ as a \emph{concave modulus of continuity}\index{concave modulus of continuity}. 
Let $X$ and $Y$ be metric spaces.
We say a continuous map 
\[
f\co X\longrightarrow Y\]
is $\alpha$--continuous, or $C^\alpha$, if there exists a $C>0$ 
such that for all $x,y\in X$ we have
\[d_Y(f(x),f(y))|\leq C\cdot\alpha(d_X(x,y)).\] 
The infimum of possible values of $C$ is called the $\alpha$--norm of $f$, and is written $[f]_{\alpha}$.
We let $C^\alpha(X,Y)$ denote the space of $\alpha$--continuous maps from $X$ to $Y$.

The reader will note that if $\alpha(x)=x$ then $\alpha$--continuity is the same as Lipschitz continuity.
More generally, if $\alpha(x)=x^{\tau}$ for $\tau\in (0,1]$ then $\alpha$--continuity is the same as $\tau$--H\"older continuity.

The reason for considering concave moduli of continuity is because nonconstant functions generally do not satisfy non-concave modulus
of continuity bounds. The reader may check, for instance, that if $\tau>1$ then only constant functions are $\tau$--H\"older continuous.
We say that $f\in\Diff_+^{1,\alpha}(S^1)$ if $f$ is a diffeomorphism and if \[[f']_{\alpha},[(f^{-1})']_{\alpha}<\infty.\]

If $f\in\Diff_+^1(S^1)$
is arbitrary then $f\in\Diff_+^{1,\alpha}(S^1)$ for some suitable \emph{smooth} concave modulus $\alpha$.
 This follows from the following somewhat more general fact.
A \emph{geodesic}\index{geodesic} in a metric space $X$ is a continuous path
\[ \gamma\co [0,D]\longrightarrow X\]
such that \[d_X(\gamma(s),\gamma(t))=|s-t|\] for all $s,t\in[0,D]$.
A \emph{geodesic space}\index{geodesic metric space} is a metric space every pair of points in which can be joined by a geodesic.

\begin{lem}[cf.~\cite{CKK2019}, Proposition 2.7]\label{lem:modulus-existence}
Let $X$ and $Y$ be metric spaces. If $X$ is a geodesic space and if $f\co X\longrightarrow Y$ is a uniformly continuous function,
then there exists a concave modulus\[\beta\co[0,\infty)\longrightarrow [0,\infty)\] which is smooth on $(0,\infty)$ such that the following hold:
\be[(i)]
\item $f$ is $\beta$--continuous;
\item if $\alpha$ is a concave modulus such that $f$ is $\alpha$--continuous,
then  there exist constants $\delta>0$ and $C\geq 1$ such that
 \[ \beta(x)\le C \alpha(x)\] for all $x\in[0,\delta]$.
\ee
\end{lem}

The concave modulus $\beta$ above is called an \emph{smooth optimal concave modulus of continuity for $f$}\index{optimal modulus of
continuity}.
We will postpone a detailed proof of the above lemma to Appendix~\ref{ch:append1}.

\begin{rem}
The regularity of a function is naturally treated as a local property.
More precisely, we say a map $f\co X\longrightarrow Y$ is \emph{locally $\alpha$--continous}\index{locally $\alpha$--continuous}
if for each $x_0\in X$ there exists a ball $B$ of finite diameter centered at $x_0$ 
such that the restriction $f\restriction_B$ is $\alpha$--continuous. We denote by
\[
C^\alpha(X,Y)\]
the space of locally  $\alpha$--continuous maps from $X$ to $Y$.
On the other hand,
it is an easy exercise~\cite[Lemma 1.1]{CKK2019} to see that if $X$ is compact then a locally $\alpha$--continuous map from $X$ to $Y$ is $\alpha$--continuous.
As in this book we focus on the compact spaces $S^1$ and $[0,1]$ or on compactly supported diffeomorphisms,
the notion of local $\alpha$--continuity coincides with the notion of $\alpha$--continuity most of the time.
We  note also that in the above lemma, with an extra hypothesis that $X$ is compact the optimal modulus $\beta$ satisfies
\[
C^\beta(X,Y)=\bigcap\left\{ C^\alpha(X,Y)\mid \alpha\text{ is a concave modulus such that }f\in C^\alpha(X,Y)\right\}.\]
\end{rem}

Lemma~\ref{lem:modulus-existence} shows that \[\Diff_+^1(S^1)=\bigcup_{\alpha}\Diff_+^{1,\alpha}(S^1),\] as $\alpha$ ranges over
all concave moduli of continuity. Therefore we may formulate the following question:

\begin{que}\label{que:denjoy-counter}
For $\theta$ irrational and $\alpha$ a concave modulus of continuity,
does there exist an exceptional diffeomorphism $f\in\Diff^{1,\alpha}_+(S^1)$ with rotation number $\theta$?
\end{que}

Some partial answers to Question~\ref{que:denjoy-counter} are immediate. For instance, if $\alpha(x)=x$ then there are no such
diffeomorphisms. Indeed, if $f'$ is Lipschitz then it is also of bounded variation, so that $\log(f')$ has bounded variation; therefore
Denjoy's Theorem says that if $\rot(f)=\theta$ then $f$ is conjugate to $R_{\theta}$.

Weakening the modulus of continuity to
$\tau$--H\"older continuity for $\tau\in (0,1)$ then the answer is completely different: Herman~\cite{Herman1979}
proved that for $\alpha=x^{\tau}$
there are exceptional diffeomorphisms in $\Diff^{1,\alpha}_+(S^1)$ for all irrational rotation numbers. Herman's methods applied to a wider
range of moduli. For instance, he produced exceptional diffeomorphisms with arbitrary irrational rotation number for the modulus
\[\alpha(x)=x(\log 1/x)^{1+\eps},\] where $\epsilon>0$ is arbitrary. 
This last modulus is stronger than every $\tau$--H\"older modulus of
continuity for $\tau<1$ since for small values of $x$ we hae \[x(\log 1/x)^{1+\epsilon}\leq C\cdot x^{\tau}\] for a suitably chosen constant
$C$, but is weaker than the Lipschitz modulus, as the reader may check.

A definitive answer to Question~\ref{que:denjoy-counter} is by no means the end of the discussion about the existence and nonexistence
of exceptional diffeomorphisms of the circle. Sullivan~\cite{Sullivan1992} and Hu--Sullivan~\cite{HuSu1997},
for example, studied other regularity conditions one can place
on the derivative of a diffeomorphism with irrational rotation number in order to guarantee conjugacy to an irrational rotation, including the
\emph{Zygmund condition}\index{Zygmund condition} and 
\emph{bounded quadratic variation}\index{quadratic variation}. We will not be discussing this aspect of the theory here, as it would
cause us to stray too far from the main narrative, and instead we direct the reader to the literature.

The primary result that we shall discuss for the remainder of this section is the following, which recovers Herman's results and furnishes
new moduli of continuity that admit exceptional diffeomorphisms.

\begin{thm}[See~\cite{KK2020-DCDS}, Theorem 1.2]\label{thm:kk-denjoy}
Let $\theta$ be an irrational rotation number, and suppose $\alpha$ is a concave modulus of continuity such that
\[\int_0^1\frac{dx}{\alpha(x)}<\infty.\] Then there exists an exceptional $f\in\Diff^{1,\alpha}(S^1)$ with $\rot(f)=\theta$.
For all $\eps>0$, we may arrange \[\sup_{x\in S^1} |f(x)-R_{\theta}(x)|+\sup_{x\in S^1} |f'(x)-1|<\eps,\]
so that $f$ is $\eps$--close to $R_{\theta}$ in
the $C^1$--topology.
\end{thm}

Notice that for \[\alpha(x)\in\{x(\log 1/x)^{1+\epsilon},\, x^{\tau}\mid \eps>0,\,\tau<1\},\] we have that $1/\alpha(x)$ is integrable near zero.
We can even get integrability of $1/\alpha(x)$ near zero for
\begin{align*}
\alpha(x)=x(\log 1/x)(\log\log 1/x)^{1+\epsilon}\\ \alpha(x)=x(\log 1/x)(\log\log 1/x)(\log\log\log 1/x)^{1+\epsilon}\\
\alpha(x)=x(\log 1/x)(\log\log 1/x)(\log\log\log 1/x)(\log\log\log\log 1/x)^{1+\epsilon}\ldots,
\end{align*}
and this furnishes many moduli of continuity, each stronger than the previous, that admit exceptional diffeomorphisms. Of course,
the function $1/x$ is not integrable near zero, consistent with the nonexistence of exceptional diffeomorphisms with Lipschitz moduli.
The function $\alpha(x)=x(\log 1/x)$ is not integrable near zero, and it is a well-known open question within the field whether or not there
exist exceptional diffeomorphisms with this modulus of continuity. It is even suspected that the answer to this last question depends on the
diophantine properties of the rotation number $\theta$ in question. The reader is directed to
Section 4.1.4 of~\cite{Navas2011} for a more detailed discussion.

\subsubsection{Herman's construction of exceptional diffeomorphisms}

The main idea behind the construction of exceptional diffeomorphisms is to reduce the problem to function theory. That is, in order to
construct the diffeomorphism $f$, one should construct the function $g=f'-1$ by hand, and then $f$ is obtained just by integration. The
function $g$ will consist of a sequence of bumps that are sufficiently smooth, and that are supported on a carefully chosen shrinking
family of intervals, which are then glued into the circle as in our original description of a continuous Denjoy counterexample. To see
how this might work, we have the following fact.

\begin{prop}[See~\cite{KK2020-DCDS}, Proposition 4]\label{prop:herman-1}
Let $\alpha$ be a concave modulus of continuity, and let \[\{J_i=(x_i,y_i)\}_{i\in\bZ}\] be a collection of disjoint subintervals of $S^1$. Suppose
that $g$ is a positive, $\alpha$--continuous function on $S^1$ that satisfies the following conditions:
\begin{enumerate}[(1)]
\item
We have
\[\int_{S^1} g\,d\theta=1;\]
\item
For all $i\in\bZ$, we have \[\int_{J_i}g\,d\theta=|J_{i+1}|;\]
\item
For all $i$, we have \[\int_{x_{i-1}}^{x_i}g\,d\theta=|x_{i+1}-x_i|.\]
\end{enumerate}
Then the function defined by \[f(x)=x_1+\int_{x_0}^x g\,d\theta\] is an element of $\Diff_+^{1,\alpha}(S^1)$ that satisfies $f(J_i)=J_{i+1}$ for
all $i\in\bZ$.
\end{prop}
\begin{proof}
This is a straightforward verification. First, check that $f'$ is $\alpha$--continuous and that $f$ is a bijection, and that $(f^{-1})'$ is
$\alpha$--continuous by an easy application of the inverse function theorem. One also has, by an easy
calculation, that \[x_i+\int_{x_{i-1}}^x g\,d\theta=x_{i+1}+\int_{x_i}^x g\,d\theta,\] and so $f$ has the desired properties.
\end{proof}

Let $\{J_i\}_{i\in\bZ}$ be a sequence of compact intervals in $S^1$. We say that $\{J_i\}_{i\in\bZ}$ is 
\emph{circular order preserving}\index{circular order} if
\[J_i<J_k<J_{\ell}\quad \textrm{if and only if}\quad J_{i+1}<J_{k+1}<J_{\ell+1}.\] Recall that we use the notation $\{x\}\in(0,1)$ for the fractional
part of $x\in\bR$.

The continuous Denjoy examples as we have
originally constructed them can be obtained via integration as in Proposition~\ref{prop:herman-1}, though not by integrating a function
but rather an atomic measure. Let $\{\ell_i\}_{i\in\bZ}$ be a collection of positive lengths, whose total sum adds up to at most $1$,
and let $\theta$ be an irrational number. Let $\mu$ be the measure on $S^1$ defined by
\[\mu=\left(1-\sum_{i\in\bZ}\ell_i\right)\lambda+\sum_{i\in\bZ}\ell_i\delta_{i\cdot\theta},\] where $\lambda$ denotes Lebesgue measure
(or total length $1$) and $\delta$ denotes a Dirac mass. Writing
$x_i=\mu([0,\{i\cdot\theta\})$ for $i\in\bZ$, we may write \[J_i=[x_i,x_i+\ell_i].\]
It is easy to see that the collection of intervals $\{J_i\}_{i\in\bZ}$ is a
disjoint, circular order preserving collection of compact intervals in $S^1$.

Let $f\in\Diff_+^{1,\alpha}(S^1)$ be an exceptional diffeomorphism, where as before the function $\alpha$ is a concave modulus of continuity.
Let $C$ denote the exceptional minimal set, and let $J\sse S^1\setminus C$ be a maximal open subinterval. We have that for
all $n\in\bZ\setminus \{0\}$, the equality $f^n(J)\cap J=\varnothing$. We let $J=J_0$ and write $J_i=f^i(J)$ for $i\in\bZ$. For compactness
of notation, we write $\ell_i=|J_i|$ for the length of $J_i$.

The following estimate is a crucial in analyzing smooth exceptional diffeomorphisms.

\begin{lem}[Fundamental Estimate;~\cite{KK2020-DCDS}, Lemma 3.1]\label{lem:fund-est}
Let $f$ be as above. Suppose that \[\lim_{n\to\infty}\frac{\ell_{n+1}}{\ell_{n}}=1.\] Then
\[\sup_{n\in\bZ}\frac{1}{\alpha(\ell_n)}\left|1-\frac{\ell_{n+1}}{\ell_n}\right|<\infty.\]
\end{lem}

To unpack the meaning of Lemma~\ref{lem:fund-est}, recall that $\ell_n\longrightarrow 0$ as $n\to\infty$, and similarly
\[\left|1-\frac{\ell_{n+1}}{\ell_n}\right|\longrightarrow 0\] as $n\to\infty$ as the successive ratios of lengths of the $\{J_i\}_{i\in\bZ}$ tends to one.
Thus, in order for $f$ to lie in $\Diff_+^{1,\alpha}(S^1)$ with this simple limiting behavior for lengths of wandering intervals, then
$\alpha(\ell_n)$ must control the size of $(1-\ell_{n+1}/\ell_n)$.

\begin{proof}[Proof of Lemma~\ref{lem:fund-est}]
Let $C$ denote the minimal set of $f$. We claim that the derivative satisfies $f'\equiv 1$ on $C$. Indeed, this follows from the Mean Value
Theorem. For all $n\in\bZ$ there exists $p_n\in J_n$ such that $f'(p_n)=\ell_{n+1}/\ell_n$. By Theorem~\ref{thm:minimal-set} and since
\[\lim_{n\to\pm\infty}\ell_n=0,\]
for an arbitrary $x\in C$, we may extract a subsequence $S_x=S\sse \{p_n\}_{n\in\bZ}$ where $f'\longrightarrow 1$ along $S$ and
$S$ converges to $x$. It follows from the continuity of the derivative that $f'$ is identically $1$ on $C$.

Now, note that if $x\in\partial J_n$ then \[\left|1-\frac{\ell_{n+1}}{\ell_n}\right|=|f'(x)-f'(p_n)|\leq [f']_{\alpha}\cdot\alpha(\ell_n).\] The lemma
now follows.
\end{proof}

Lemma~\ref{lem:fund-est} shows where much of the the difficulty lies in constructing exceptional diffeomorphisms
with a particular modulus of continuity.
In order to build such diffeomorphisms, one needs to find a sequence of finite lengths of finite total length, so that the inequality in
the lemma holds. This may not always be possible. For instance, if $\ell_{i+1}/\ell_i\to 1$ as $i\to\infty$, then one can show that if $f'$ is
$\alpha$--continuous then there is a positive constant $A$ such that $\alpha(\ell_i)\geq A/i$. This shows that one cannot find such an
exceptional diffeomorphism such that $f'$ is Lipshitz, or even $(x(\log 1/x))$--continuous.
See~\cite{KK2020-DCDS} for a more thorough discussion of this point.

The fundamental estimate in Lemma~\ref{lem:fund-est} above can be used in conjunction with Proposition~\ref{prop:herman-1} to produce
many smooth Denjoy counterexamples, once lengths of wandering intervals are carefully chosen. The reader will note that the following
proposition provides a kind of converse to the fundamental estimate.

\begin{prop}[\cite{KK2020-DCDS}, Proposition 5, cf.~\cite{Herman1979} and Chapter 12 of~\cite{KH1995}]\label{prop:herman-2}
Let $\alpha$ be a concave modulus of continuity, and let \[J_k=[x_k,y_k],\, k\in\bZ\] be a sequence of disjoint, circular order preserving intervals
in $S^1$. Suppose furthermore that for all $i\in\bZ$ we have \[\lambda\left( [x_{i-1},x_i]\setminus \bigcup_{k\in\bZ} J_k\right)=
\lambda\left( [x_{i},x_{i+1}]\setminus \bigcup_{k\in\bZ} J_k\right),\] where as before $\lambda$ denotes Lebesgue measure.
Writing $\ell_k=|J_k|$ as before, assume that \[\sup_k\frac{1}{\alpha(\ell_k)}
\left(1-\frac{\ell_{k+1}}{\ell_k}\right)<\infty,\] and that \[\inf_{k}\frac{\ell_{k+1}}{\ell_k}>0.\]
Then there exists an $f\in\Diff^{1,\alpha}_+(S^1)$ such that $f(J_k)=J_{k+1}$ for all $k$.
\end{prop}
We note that the condition that the infimum of successive lengths is bounded away from zero is automatically satisfied if the successive
ratios converge to $1$, for example.

\begin{proof}[Proof of Proposition~\ref{prop:herman-2}]
For $k\in\bZ$, let $\rho_k$ be an arbitrary smooth function on $[0,1]$ satisfying \[\int_{[0,1]}\rho_k\,dx=1,\] and which satisfies
\[1-\left(1-\frac{\ell_{k+1}}{\ell_k}\right)\rho_k(x)>0\] for all $x\in [0,1]$. We will not write down a formula for $\rho_k$, as a construction
of such functions is standard from differential topology. Since \[\inf_k\frac{\ell_{k+1}}{\ell_k}>0\] by assumption, we may assume that
$\rho_k$ is uniformly bounded, independently of $k$, and in fact we may choose a uniform Lipschitz constant.
We write \[g(x)=1-\sum_{k\in\bZ}\left(1-\frac{\ell_{k+1}}{\ell_k}
\rho_k\left(\frac{x-x_k}{\ell_k}\right)\right).\] Note that $g$ is positive by construction.
Integrating, $g$ over $S^1$ results in a telescoping sum that gives a total integral of $1$. The remaining hypotheses of
Proposition~\ref{prop:herman-1} are checked similarly. To get the result via Proposition~\ref{prop:herman-1}, it suffices to
check that $g$ is $\alpha$--continuous.

Note that \[\sup_{x,y\in J_k}\frac{|g(x)-g(y)|}{\alpha(|x-y|)}=1-\frac{\ell_{k+1}}{\ell_k}
\sup_{s,t\in [0,1]}\frac{|\rho_k(s)-\rho_k(t)|}{\alpha(\ell_k|s-t|)}.\] Note that
\[\frac{|\rho_k(s)-\rho_k(t)|}{\alpha(\ell_k|s-t|)}\leq C_k\frac{|s-t|}{\alpha(\ell_k|s-t|)},\] where here $C_k$ denotes the Lipschitz constant of
$\rho_k$. The concavity of $\alpha$ implies that the function $x/\alpha(x)$ is monotone increasing, so we get that
\[C_k\frac{|s-t|}{\alpha(\ell_k|s-t|)}\left(1-\frac{\ell_{k+1}}{\ell_k}\right)\leq \frac{C}{\alpha(\ell_k)}\left(1-\frac{\ell_{k+1}}{\ell_k}\right),\] where $C$
is the supremum of the Lipschitz constants of $\{\rho_k\}_{k\in\bZ}$. This bounds the $\alpha$--norm of $g$ on the intervals $\{J_k\}_{k\in\bZ}$,
and since $g$ is identically $1$ outside of the union of these intervals, we see that $[g]_{\alpha}$ is bounded.
\end{proof}

As a consequence, we have the following:

\begin{cor}\label{cor:herman-3}
Let $\alpha$ be a concave modulus of continuity. Suppose there exists a sequence $\{\ell_k\}_{k\in\bZ}$
of positive real numbers such that \[\sum_{k\in\bZ}\ell_k\leq
1\] and such that \[\sup_k\frac{1}{\alpha(\ell_k)}\left(1-\frac{\ell_{k+1}}{\ell_k}\right)<\infty,\] then there exists an exceptional diffeomorphism
$f\in\Diff_+^{1,\alpha}(S^1)$ with a wandering interval $J\sse S^1$ such that $|f^k(J)|=\ell_k$ for all $k$.
\end{cor}

\subsubsection{Integrability of moduli and exceptional diffeomorphisms}\label{ss:integrability}
We are now ready to prove Theorem~\ref{thm:kk-denjoy}. 

\begin{proof}[Proof of Theorem~\ref{thm:kk-denjoy}]
By Propostition~\ref{prop:herman-2} and Corollary~\ref{cor:herman-3},
in order to construct the exceptional diffeomorphism $f$, it suffices to use the integrability of $\alpha$ near zero to construct
suitable lengths. We will use the modulus $\alpha$ to construct the lengths of the wandering intervals as functions of their index.
By Lemma~\ref{lem:modulus-existence}, we may assume that $\alpha$ is a smooth function. Write
\[K=\max\{2,1/\alpha(1)\}\quad \textrm{and}\quad v(x)=x^2\alpha(1/x).\]
Observe that for $t\geq 1$ we have $v(x/t)\geq v(x)/t^2$. Indeed, we compute:
\[v(x/t)=(x^2/t^2)\alpha(t/x)\geq (x^2/t^2)\alpha(1/x)=v(x)/t^2.\] Since $\alpha$ is concave, we have $x/\alpha(x)$ is monotone increasing, and
therefore \[x\alpha(1/x)=\frac{\alpha(1/x)}{1/x}\] is also monotone increasing and hence has a nonnegative derivative.
Therefore, we can compute \[(v(x)/x)'=\left(\frac{\alpha(1/x)}{1/x}\right)'\geq 0.\] It follows that when $x\geq K$ then \[v(x)\geq x\cdot v(1)\geq
x/K.\]

We also have that $\alpha(x)/x$ has nonpositive derivative, so that \[\left(\frac{\alpha(x)}{x}\right)'=\frac{x\alpha(x)'-\alpha(x)}{x^2}\leq 0.\]
It follows that \[\alpha(x)\geq x\alpha'(x)>0\] for all $x$, and so that in particular \[v(x)=x^2\alpha(1/x)\geq x\alpha'(1/x)>0.\]

We combine these observations to see that
\[xv'(x)=x(2x\alpha(1/x)-\alpha'(1/x))=2v(x)-x\alpha'(1/x)\in [v(x),2v(x)].\]

Armed with these estimates, we can define the lengths of the wandering intervals. We write \[\ell_k=\frac{1}{v(|k|+K)}\] for $k\in\bZ$. The
integrability assumption on $\alpha$ is precisely what allows us to choose such lengths to be finite in total. Namely, we have 
\[\int_K^{\infty}v(x)\,dx=\int_0^{1/K}\frac{dx}{\alpha(x)}<\infty.\] So, we may set \[\sum_{k\in\bZ}\ell_k\leq 1,\] at the cost of possibly
increasing the value of $K$.

To complete the construction of $f$, we need only to estimate the size of \[\frac{1}{\alpha(\ell_k)}\left(1-\frac{\ell_{k+1}}{\ell_k}\right),\] 
and show that it is bounded independently of $k$.

For $k\in \bZ$, we will write $j=|k|$ for ease of notation. We have shown that \[v(j+K)\geq \frac{j+K}{K},\] and observe that
\[v(j+K\pm 1)\geq v\left(\frac{j+K}{2}\right)\geq\frac{v(j+K)}{4}.\] The second inequality follows from the inequality $v(x/t)\geq v(x)/t^2$ that
holds for $t\geq 1$, as we have already established.

We may now apply the Mean Value Theorem to a suitable point \[y\in [j+K,j+K+1]\quad \textrm{or} \quad y\in [j+K-1,j+K].\] We have
\begin{align*}
&\frac{1}{\alpha(\ell_j)}\left|1-\frac{\ell_{j+1}}{\ell_j}\right|=\frac{1}{\ell_j^2v(1/\ell_j)}\left|1-\frac{v(j+K)}{v(j+K\pm 1)}\right|\\
&=\frac{v(j+K)^2}{v\circ v(j+K)}\cdot \frac{v'(y)}{v(j+K\pm1)}=\frac{v(j+K)}{v\circ v(j+K)}\cdot\frac{v(j+K)}{v(j+K\pm 1)}\cdot v'(y)\\
&\leq \frac{j/K+1}{v(j/K+1)}\cdot 4\cdot\frac{2v(y)}{y}\leq \frac{j/K+1}{v(2j+2K)/(2K)^2}\cdot\frac{8v(y)}{y}\\
&=32K\cdot \frac{j+K}{y}\cdot\frac{v(y)}{v(2j+2K)}\leq 64K.
\end{align*}

This furnishes a uniform bound for \[\frac{1}{\alpha(\ell_k)}\left(1-\frac{\ell_{k+1}}{\ell_k}\right),\] independently of $k$ and implies the
existence of the claimed exceptional diffeomorphism.

To complete the proof of the result, we will show that $f$ can be chosen as close to the rotation $R_{\theta}$ as we like, in the $C^1$
topology. To see that $f$ can be chosen to be $C^0$--close to $R_{\theta}$, we may simply increase the cutoff $K$ so that
\[\sum_{k\in\bZ}\ell_k\approx 0.\] 

Achieving this, we have that for each $\eps>0$, we may choose a cutoff so that
\[\sup_{x\in S^1}|f(x)-R_{\theta}(x)|<\epsilon.\]

To get $f$ to be $C^1$--close to $R_{\theta}$, we need to show that for each $\eps>0$ we may arrange \[\sup_{x\in S^1}|f'(x)-1|<\eps.\]
We have that $f'(x)-1$ is identically zero outside of the wandering intervals, and using the construction of $f$ from
Proposition~\ref{prop:herman-2}, we see that on the component $J_k$ wandering set, the function $f'(x)-1$ is simply given by
\[\left(1-\frac{\ell_{k+1}}{\ell_k}\right)\rho_k\left(\frac{x-x_k}{\ell_k}\right).\]

The function $\rho_k$ is bounded independently of $k$. We have that $\ell_{k+1}/\ell_k\to 1$ as $k\to\infty$, and so by setting the cutoff
$K$ to be sufficiently large, we may arrange for $|f'(x)-1|<\eps$ on $J_k$, which completes the proof.
\end{proof}

\subsubsection{On the lengths of wandering intervals}
We end this chapter with a remark on an open problem in the theory of exceptional diffeomorphisms. If $f\in\Diff^1_+(S^1)$ is an exceptional
diffeomorphism with minimal set $C$, it might not generally be the case that $f'\equiv 1$ on $C$. That is, the estimate in
Lemma~\ref{lem:fund-est} might fail. In that case, it is not true that $\ell_{k+1}/\ell_k\to 1$ for successive wandering intervals, as $k\to\infty$.
 Without
this assumption, the construction of exceptional diffeomorphisms we have carried out here does not work at all.

There is another angle from which to probe lengths of wandering intervals for exceptional diffeomorphisms. Consider the set
$\{\ell_k\}_{k\in\bZ}$, and reorder it according to decreasing lengths. Thus, we get a sequence $\{\lambda_i\}_{i\ge1}$ of real numbers,
where $\lambda_i$ is the $i^{th}$ largest of the real numbers $\{\ell_k\}_{k\in\bZ}$. This sequence has very little to do with the original
order of the lengths $\{\ell_k\}_{k\in\bZ}$. It is a well-known open question in the field as to whether or not $\lambda_{i+1}/\lambda_i\to 1$
as $i\to\infty$. This question was originally posed by D.~McDuff, and she was able to prove that $1$ is an accumulation point of successive
quotients of lengths, and that the quotients are bounded independently of $i$. The reader is directed
to~\cite{McDuff1981,athanassopoulos,KK2020-DCDS} for a more detailed discussion.


%
%
%
\chapter{Full diffeomorphism groups determine the diffeomorphism class of a manifold}\label{sec:filip-tak}

\begin{abstract}In this chapter, we establish that the algebraic structure of the full $C^p$ diffeomorphism group of a smooth connected boundaryless manifold can recover the regularity $p$, and the manifold up to $C^p$ diffeomorphisms. 
This was proved by Filipkiewicz in 1980s, following a preceding, weaker result of Takens. Filipkiewicz's proof uses a strategy pursued by Whittaker,
and conclude
by the Bochner--Montgomery Theorem regarding Hilbert's Fifth Problem.
Instead of using Whittaker's method, we give a proof using a powerful reconstruction theorem due to Rubin. We also give a generalization of Filipkiewicz original approach as an alternative proof. This chapter is mostly self-contained, containing complete proofs of Bochner--Montgomery Theorem and Rubin's Theorem.\end{abstract}

In this chapter we leave the world of one-manifolds and discuss general diffeomorphism groups of manifolds, before returning
to one--manifolds in the next chapter. We include this discussion now because it gives a natural context to introduce several tools
which will be useful in later chapters.

One of the main goals of this monograph will be to give a self--contained proof that the isomorphism types of
 finitely generated subgroups
of $\Diff_+^p(M)$ determine both $M$ and $r$, at least in the case where $M$ is a $1$--manifold. It is therefore natural to wonder
the degree to which this result generalizes to higher dimensional manifolds. We will begin with some musings about critical regularity of
finitely generated groups acting on manifolds in dimension two or more. The bulk of this chapter will be devoted to results of F.~Takens and
R.~Filipkiewicz, which show that for general manifolds, the isomorphism type of the group $\Diff_c^p(M)_0$, where here $1\leq p\leq\infty$,
determines the manifold $M$ up to $C^p$--diffeomorphism. 
We make a standing assumption throughout this book that, unless explicitly
noted to the contrary,\emph{ all manifolds are Hausdorff and second countable, hence admit partitions of unity.}
We will also let $M$ denote a smooth connected boundaryless manifold of dimension $d$. 

We will recall the original proofs that were given by Takens and Filipkiewicz in their original papers. Filipkiewicz's paper is not self-contained,
requiring background results in the algebraic structure of diffeomorphism groups, together with a subtle regularity result due to
Bochner--Montgomery. We will include a complete account of the simplicity of commutator subgroups of diffeomorphism groups of manifolds,
and a complete proof of the Bochner--Montgomery Theorem. From the Bochner--Montgomery Theorem, we will be able to prove
Takens' Theorem, and we will give Takens' proof in one dimension. To prove Filipkiewicz's Theorem, we will give a self--contained account
of Rubin's Theorem, which is an extremely powerful reconstruction result for objects with sufficiently complicated groups of automorphisms.
We will retain a polished version of Filipkiewicz's result as a matter
of historical record and to provide his perspective on this flavor of reconstruction result.

\section{Diffeomorphism groups of general manifolds and critical regularity}
For a $C^p$ manifold $M$ with $p\in\bZ_{>0}\cup\{0,\infty\}$, we let 
$\Diff^p(M)$ denote the group of $C^p$--diffeomorphisms of $M$. 
The group of compactly supported diffeomorphisms in $\Diff^p(M)$ is denoted as
$\Diff^p_c(M)$. 
If $f\in\Diff^p_c(M)$ is $C^p$--diffeotopic to the identity via a diffeotopy supported in a compact set, then we write $f\in\Diff^p_c(M)_0$. 
All these groups are equipped with the Whitney $C^p$--topology~\cite{Hirsch1994}, although we will essentially care only about its $C^1$--topology.

\begin{rem}\label{rem:diffeotopy}
The group $\Diff_c^p(M)_0$ coincides with the identity component of $\Diff_c^p(M)$, although this fact does not play a role for our purpose of this chapter. See~\cite[Corollary 1.2.2]{Banyaga1997} for a proof.\end{rem}

The ultimate goal of this chapter is to prove the following result.

\begin{thm}[Generalized Takens--Filipkiewicz Theorem]\label{thm:filip-gen}
Let $M$ and $N$ be smooth connected boundaryless manifolds, and let $p,q\in\bZ_{>0}\cup\{0,\infty\}$.
If there exists a group isomorphism 
\[ \Phi\co G\longrightarrow H\] 
for some groups $G,H$ satisfying
\[
\Diff_c^{p}(M)_0\le G\le \Diff^{p}(M)\]
and
\[
\Diff_c^{q}(N)_0\le H\le \Diff^{q}(N)\]
then we have that 
$p=q$
and that
 there exists a $C^p$ diffeomorphism
$w\colon M\longrightarrow N$ satisfying
\[\Phi(g)=w\circ g\circ w^{-1}\]
for all $g\in G$.\end{thm}

The case $p=q=0$ of the above theorem is due to Whittaker~\cite{Whittaker1963}.
Takens~\cite{Takens79} originally proved the theorem when $p=q=\infty$
and when $G$ and $H$ are full $C^\infty$ diffeomorphism groups of respective manifolds, with an additional assumption that such a 
map $w$ is already given as a bijection, and not \emph{a priori} by a 
diffeomorphismi. Filipkiewicz extended this to the theorem below, which is also a special case of Theorem~\ref{thm:filip-gen}~\cite{Filip82}.

\begin{thm}[Filipkiewicz's Theorem]\label{thm:filip}
If $p,q\in\bZ_{>0}\cup\{\infty\}$ and if 
 \[\Phi\co\Diff^p(M)\longrightarrow\Diff^q(N)\] is a group isomorphism,
then we have that $p=q$ and that there exists a $C^p$ diffeomorphism
$w\colon M\longrightarrow N$ satisfying
\[\Phi(g)=w\circ g\circ w^{-1}\]
for all $g\in \Diff^p(M)$.
\end{thm}

We end this introduction by formulating a conjecture that is natural in light of the previous remarks (cf.~Question~\ref{que:subgroup-main}
from the introduction).

\begin{con}\label{con:critreg-hd}
Let $M$ and $N$ be compact connected $C^p$ manifolds. 
Then $M$ and $N$ are $C^p$--diffeomorphic if and only if the
set of isomorphism classes of finitely generated subgroups
of $\Diff^p_c(M)_0$ and $\Diff^p_c(N)_0$ coincide.
\end{con}

We remark that, as discussed in the introduction,
homeomorphisms in dimension one are controlled by the fact that, away from fixed points, they always have a well--defined
direction. This structure imposed by the ambient topology makes many of the arguments made in this book possible. In higher dimension,
such structure is lost, and behavior of diffeomorphism groups of even two--dimensional manifolds is much more complicated than the
one--dimensional case.

Conjecture~\ref{con:critreg-hd}, if correct, has the potential to yield a fresh perspective on classical questions in differential topology,
especially relating to exotic differentiable structures, and can give reformulations of these problems in terms of combinatorial group theory.

Conjecture~\ref{con:critreg-hd} in its stated form is wide open. As such, we will not comment on it any further here.

\section{Generalities from differential topology}
In this section, let us collect standard facts on manifolds and their diffeomorphism groups. Two main results to be
established are the existence of a compactly supported smooth group action of $\bR^d$ on $M$,
and the fragmentation of $C^p$ diffeomorphisms on $M$. Canonical references include~\cite{Hirsch1994},
~\cite{deRham1984} and~\cite{PS1970}.

\subsection{Open balls in manifolds}\label{s:balls}
Let $g$ be a homeomorphism of a topological space $X$.
To avoid confusion with the existing literature, we  use a notation
\begin{align*}
\suppo g &:=\supp g=X\setminus \Fix g,\\
\suppc g&:=\overline{\supp g}.
\end{align*}

\begin{rem}
The latter of the above is called the \emph{closed support}\index{closed support} of $g$,
which is often denoted simply as $\supp g$ in the literature of differential topology. Though, some authors reserve the
notation $\supp g$ for the \emph{open support}\index{open support} $\suppo g$, as we have done in other chapters of this book.
\end{rem}

Let $X$ be a topological space, and let  $G\le\Homeo(X)$ be a subgroup.
We  write $G_c$ the subgroup of $G$ consisting of compactly supported homeomorphism groups.
When $G$ is given with a topology a priori,
we let $G_0$ denote the identity component of $G$.
Consequently, $G_{c0}$ is interpreted as $(G_c)_0$.

Let $M$ be a smooth connected boundaryless manifold of dimension $d$. 
In the case when  $X=M$
and when $G=\Diff^p(M)$, the group of $C^p$ diffeomorphisms of $M$,
 we use the notation $\Diff_c^p(M)$, $\Diff^p(M)_0 $ and $\Diff_c^p(M)_0$ slightly changing the positions of the subscripts. 

We let $B^d(a;r)\sse\bR^d$ denote the open ball with radius $r$ centered at $a\in\bR^d$.
It will be useful for us to have an open basis of a manifold $M$ that is invariant under $\Diff^p(M)$, defined by the following open sets.

\bd\label{defn:ball}
For a manifold $M$ of dimension $d$, we say an open set $U\sse M$ is a \emph{$C^p$--open ball}\index{$C^p$--open ball}
if there exists a $C^p$--embedding 
\[ h\co \bR^d\longrightarrow M\] satisfying that $U=h(B^d(0;1))$.
\ed
We denote by $\BB^p(M)$ the collection of all the $C^p$ open balls in $M$.
We simply say $U\in\BB^p(M)$ is a \emph{$C^p$--open ball}, or simply an \emph{open ball} if the meaning is clear.

Let $G$ be a group acting on a set $X$. Then $G$ induces the \emph{pseudo-group topology}\index{pseudo-group topology}
on $X$, for which the basic open sets are
subsets $U\sse X$ for which $U=\supp g$ for some $g\in G$. 
The lemma below shows that the affine action of $\bR^d$ on the Euclidean space can be 
conjugated into the compactly supported smooth diffeomorphism group of $M$.
This lemma  implies that the original Euclidean topology of $M$, the pseudo-group topology coming from 
$\Diff^p_c(M)_0$, and the topology generated by $\BB^p(M)$ all coincide. 
In particular, the collection $\BB^p(M)$ is a $\Diff^p(M)$--invariant open basis of $M$.

\begin{lem}[\cite{deRham1984}]\label{lem:affine}
Let $d\ge1$.
Then there exists a smooth diffeomorphism 
\[
h\co  B^d(0;1)\longrightarrow\bR^d\]
such that for each fixed $v\in\bR^d$ the map
\[
a_v(x):=\begin{cases}
 h^{-1}\left( v+h(x)\right), &\text{ if }x\in B^d(0;1),\\
 x,& \text{ otherwise.}\end{cases}\]
is a smooth diffeomorphism of $\bR^d$. 
In particular, the map
\[v.x:=a_v(x)\]
defines a smooth group action of $v\in \bR^d$ on $x\in \bR^d$ supported in $B^d(0;1)$, which is conjugate to the translation action $x\mapsto x+v$. 
\end{lem}
\bp
We will provide only the recipe of the map $h$ as the computational verification is elementary and worked out in~\cite[Section 15]{deRham1984};
one can also give an alternative proof that uses Theorem~\ref{thm:bs-cpt}.
Pick an arbitrary smooth diffeomorphism
\[
\phi\co [0,1)\longrightarrow[0,\infty)\]
such that $\phi(x)=x$ for $x\in(0,1/3]$
and such that 
\[
\phi(x)=e^{1/(1-x)^2}\]
for $x\in[2/3,1)$.
Define $h\co B^d(0;1)\longrightarrow\bR^d$ by
\[
h(x):= \phi(|x|)\cdot \frac{x}{|x|}\]
and $h(0)=0$. If $X\co\bR^d\longrightarrow\bR^d$ is the constant vector field having the value $v$, then one can check that the pull back 
\[
h^*(X)(p):=(D_p h)^{-1} X\circ h(p)\]
extends to a smooth vector field on $\bR^d$ supported in $B^d(0;1)$.
\ep

Let us now consider a general context of topological spaces.
Let $X$ be a topological space.
If $G\le\Homeo(X)$ and if $U\sse X$, then we let $G_U$ be the elements $g\in G$ such that  \[\suppc g\sse U.\]

\bd\label{defn:lt-ld}
We say a group $G$ of homeomorphisms on a space $X$ is \emph{locally transitive (LT)}\index{locally transitive} if every point in
$X$ has a local basis consisting of open sets $U$ such that the group $G_U$ acts transitively on $U$.\ed

In other words, we require that for each pair $x,y\in U$ there exists $g\in G$ such that $\suppc g\sse U$ and such that $g(x)=y$.
It is trivial that if $G$ acts locally transitively on $X$ then the pseudo-group topology recovers the original topology on $X$.
Lemma~\ref{lem:affine} gives a rigorous proof of the following, which is somewhat obvious anyway from the intuition.
\begin{lem}\label{lem:ltld-diff}
An arbitrary subgroup of $\Homeo(M)$ containing $\Diff_c^\infty(M)_0$ is locally transitive.\end{lem}

\subsection{Fragmentation of diffeomorphisms}\label{ss:frag}
We say a group $G\le\Homeo(X)$ is \emph{fragmented}\index{fragmentation} if for every open cover $\UU$ of $X$, we have that 
\[G=\form{ G_U\mid U\in\UU}.\]
As the reader might have guessed, we wish to show that $\Diff^p_c(M)_0$ is fragmented. This fact can be traced at least back to M.~Hirsch.

\begin{lem}\label{lem:frag}
If $p\in\bZ_{>0}\cup\{\infty\}$, then the group  $\Diff^p_c(M)_0$ is fragmented.
\end{lem}

\begin{proof}[cf.~\cite{PS1970,Banyaga1997}]
Let $\UU$ be an arbitrary cover of $M$.
Pick an element $g\in\Diff_c^p(M)_0$, which is isotopic to the identity via diffeomorphisms $\{g_t\}_{t\in[0,1]}$
supported in a compact set $K\sse M$.
Since every topological group is generated by an arbitrary neighborhood of the identity, it suffices for us to find a
$C^1$--neighborhood $V\sse\Diff_K^p(M)_0$ of the identity that consists of fragmented diffeomorphisms;
it will then follow that $g$ is fragmented.

We have some finite cover \[\{U_1,\ldots,U_m\}\sse \UU\] of $K$.
Since $M$ is paracompact, we have a partition of unity $\{\phi_i\}_{1\leq i\leq m}$ subordinate to the cover $\{U_1,\ldots,U_m\}$ of $K$.
We set $\phi_0=0$ and
 \[\omega_j(x)=\sum_{i=0}^j \phi_i(x).\] 
 Note that $\omega_m(x)=1$ on $K$ and $\omega_i(x)=\omega_{i-1}(x)$ if $i\ge1$ and $x\in M\setminus U_i$.

Assuming we have some identity neighborhood $V\sse \Diff^p_K(M)_0$ at hand, pick $h\in V$ and a $C^p$--diffeotopy \[\{h_t\}_{t\in[0,1]}
\sse\Diff^p_K(M)_0;\] see Remark~\ref{rem:diffeotopy} or the paragraph below.
For $i=1,\ldots,m$ we have $C^p$ maps
 \[h_i(x):=h_{\omega_i(x)}(x)\]
 that are homotopic to the identity.
 
The set $\Diff^p_K(M)$ can be included in the set of all $C^p$
functions from $M$ to itself 
\[C^p_K(M,M)\]
that are the identity outside $K$; see~\cite{Hirsch1994}.
Moreover, $\Diff^p_K(M)$ is defined by $C^1$--open conditions
in $C^p_K(M,M)$.
If $V$ is small enough a priori then one can require that each $h_i$ is $C^1$--close enough to the identity after defining a priori
\[
h_t(x)=\exp_x (t \exp^{-1}_x(h(x)))\]
by the Riemannian exponential map. Hence, each $h_i$ is a $C^p$--diffeomorphism.

Observe that $h_{i-1}(x)=h_i(x)\in\Diff_K^p(M)_0$
whenever $x\notin U_i$. It follows that $h_{i-1}^{-1}h_i$ is the identity outside of $U_i$. Finally, \[h=(h_0^{-1}h_1)(h_1^{-1}h_2)\cdots
(h_{m-1}^{-1}h_m),\] which furnishes the desired fragmentation of $h\in V$.
\end{proof}

\section{Simplicity of commutator subgroups}
The following is a well-known result due essentially to D.~B.~A. Epstein, which will be required in the proof of Filipkiewicz's Theorem and
which we will prove here for the sake of completeness. 
We follow the general outline of the authors' expository work with J.~Chang in~\cite{CKK2019}, 
based on the simplification of Epstein's method by W.~Ling~\cite{Ling1984}.
 
\begin{thm}[Epstein~\cite{Epstein1970}]\label{thm:epstein}
If $M$ is a smooth connected boundaryless manifold
and if $p\in\bZ_{>0}\cup\{\infty\}$, 
then the group $[\Diff^p_c(M)_0,\Diff^p_c(M)_0]$ is simple.
\end{thm}

In fact, one can strengthen Theorem~\ref{thm:epstein} in various contexts to diffeomorphism groups with non-integral regularities.
See~\cite{CKK2019} for a self-contained account, for instance.

It is a crucial point in Theorem~\ref{thm:epstein} to restrict to the connected component of the identity of $\Diff^p(M)$. For instance,
if $\Sigma$ is a closed orientable surface of genus at least two then $[\Diff^p(\Sigma),\Diff^p(\Sigma)]$ admits an infinite discrete group
as a quotient; indeed the group $\Diff^p(\Sigma)/\Diff^p_0(\Sigma)$ is the group of isotopy classes of $C^p$ diffeomorphisms of
$\Sigma$, which is identified with the mapping class group $\Mod(\Sigma)$. This latter group has finite abelianization, which implies that
$[\Diff^p(\Sigma),\Diff^p(\Sigma)]$ surjects to a finite index subgroup of $\Mod(\Sigma)$.

There are two main ideas for the proof of Theorem~\ref{thm:epstein}.
The first is Higman's Theorem, which is a very useful result for studying the simplicity of
groups that act on a set in a highly transitive manner.
In order to apply techniques in the vein of Higman's Theorem
to groups of diffeomorphisms, one needs to be able to generate the full group by elements which are supported on sufficiently
small subsets of the manifold, which is where the technology of fragmentation comes into play; this will be the second
main idea of this section.

\subsection{Higman's Theorem}
Higman's Theorem has applications
beyond diffeomorphism groups, and is a standard tool in the study of groups of piecewise linear homeomorphisms. In the authors' opinion,
it is not nearly as well known as it should be, and so we will include a complete discussion of it.
Recall our convention that for an element $g$ in a permutation group $\Sym(X)$ of a set $X$, we denote its
\emph{support}\index{support of a permutation} by
\[\supp g=X\setminus\Fix g.\]
In a later subsection when $X$ is equipped with a topology, we will make a distinction between $\supp g$ and its closure.

\begin{thm}[Higman's Theorem]\label{thm:higman}
Let $G$ be a nontrivial group acting faithfully on a set $X$. Suppose that for all triples $\{r,s,t\}$ of nonidentity elements of $G$, there is an element
$g\in G$ such that \[g(\supp r\cup \supp s)\cap tg(\supp r\cup\supp s)=\varnothing.\] Then:
\begin{enumerate}[(1)]
\item
The commutator subgroup $G'=[G,G]$ is nonabelian and simple.
\item
Every proper quotient of $G$ is abelian.
\item
The center of $G$ is trivial.
\end{enumerate}
\end{thm}

Generally speaking, if $G$ satisfies the hypotheses of Higman's Theorem then $G'$ could be an infinitely generated simple group,
even if $G$ happens to be finitely generated.

Higman's Theorem is tailor-made to investigate certain groups of homeomorphisms of topological objects. The moment that for all
nontrivial $t\in G$ there
always exists a group element $g\in G$ that shrinks and moves the total supports of two nontrivial elements of $G$ so that $t$ moves the
shrunken and translated supports of themselves, then $G$ satisfies the conclusions of Higman's Theorem.

Higman's Theorem is sufficient to prove Theorem~\ref{thm:epstein} for
Euclidean spaces (see Remark~\ref{rem:euclidean})
although one could face
some issues which have to do with the global topology in the case of a general manifold $M$.
For two group elements $t$ and $g$ we write as usual
\[
t^g:=g^{-1}tg.\]

\begin{proof}[Proof of Theorem~\ref{thm:higman}]
We prove a sequence of claims which will allow us to then easily obtain the conclusion of the theorem.

\begin{claim}
The commutator subgroup $G'$ is nonabelian.\end{claim}
 First note that $G$ itself must be nonabelian. 
If $G$ were abelian, we set $r=s=t$ to be
nontrivial elements of $G$. For each $g\in G$ 
we have
\[
g(\supp r\cup \supp s)\cap tg(\supp r\cup\supp s)
= g\supp r\cap gt\supp r=g\supp r\ne\varnothing.\]
This violates the hypothesis on $G$. 
To see $G'$ is nonabelian, we now choose $r=s=t$ to be nontrivial elements of $G'$ and let $g\in G$ be a corresponding element from the hypothesis.
For each $x\in g\supp r=\supp grg^{-1}$, we have that 
\[x\not\in rg\supp r=\supp \left(rgr(rg)^{-1}\right).\]
It follows that
\[
x= rgr(rg)^{-1}(x)\ne rgrg^{-1}r^{-1} \left(gr^{-1}g^{-1}(x)\right).\]
This implies that $r$ and $grg^{-1}$ are non-commuting elements in $G'$,
proving the claim.

\begin{claim}
For all nontrivial $\{r,s,t\}\sse G$, there is an element $g\in G$ such that \[\left[r,t^gs(t^g)^{-1}\right]=1.\]\end{claim}
By applying $g^{-1}$ to the hypothesis of Higman's Theorem,
we see that there is an element $g\in G$ satisfying \[\supp r\cap t^g(\supp s)=\varnothing,\]
so that $r$ and $t^gs(t^g)^{-1}$ have disjoint supports. This proves the claim.

\begin{claim}If $T$ is a nontrivial normal subgroup of $G$ then $T$ contains $G'$.\end{claim}
Let $t$ be a nontrivial element of $T$, and let $r,s\in G$ be arbitrary and nontrivial. By the second claim, there is an element $g\in G$ such that
$(t^g)^{-1}r t^g$ and $s$ commute. It follows that 
\[[r,s]=[r\cdot (t^g)^{-1}r^{-1}t^g,s].\] 
Since the right hand side is a product of conjugates of $t^{\pm1}$ and $t^{-1}$,
we deduce that $[r,s]\in T$. It follows that
$G'\le T$, as claimed.

From the first and the third claims, we readily see that $G''=G'$ and that every proper quotient of $G$ is abelian.
One can also promote the given hypothesis to further require that $g$ belongs to the commutator subgroup. 

\begin{claim}For all nontrivial $\{r,s,t\}\sse G$, there exists a $u\in G'$ such that
\[u(\supp r\cup \supp s)\cap t(u(\supp r\cup\supp s))=\varnothing.\]
\end{claim}
Let $g\in G$ be as in the hypothesis of the theorem for the triple $\{r,s,t\}$. If $g\in G'$ we are done, so we may assume $g\notin G'$,
and so $g$ is nontrivial. Applying the second claim for the triple 
$(t,g,t)$, we have an element $h\in G$ such
that \[[t,{t^h}g(t^h)^{-1}]=1.\] Set $u:=[t^h,g^{-1}]\in G'$.
Since \[[t,ug^{-1}]=[t,t^hg^{-1}(t^h)^{-1}]=1,\] we have that
$t^u=t^g$. 
It follows that $u$ gives us a suitable element of $G'$ since
\[(\supp r\cup \supp s)\cap t^u(\supp r\cup\supp s)
=(\supp r\cup \supp s)\cap t^g(\supp r\cup\supp s)=
\varnothing.\]
This proves the claim.

We can now complete the proof of the theorem. Let $Z$ be the center of $G$. We have that $Z$ is an abelian and normal subgroup of $G$.
By the third claim, if $Z$ is nontrivial then $Z$ contains $G'$. This latter group is nonabelian by the first claim, so $Z$ is indeed trivial.

It only remains to show that $G'$ is simple. 
Suppose that $T$ is a nontrivial normal subgroup of $G'$, and let $t\in T$ be a nontrivial element. If $r,s\in G'$ are
arbitrary nontrivial elements, then the last claim shows that for some $u\in G'$ we have \[[r,{t^u}s(t^u)^{-1}]=1.\]
Applying the same trick as in the third claim, 
we see that \[[r,s]=[r\cdot (t^u)^{-1}rt^u,s]\in T.\] It follows then that 
\[G''=[G',G']\le T.\] As we have seen $G'=G''$, we conclude $T=G'$. This proves the simplicity of $G'$.
\end{proof}

Higman's Theorem is formulated in an essentially combinatorial manner. For us, it is useful to recast it in a topological manner,
and in the process making more precise the remarks following the statement of Theorem~\ref{thm:higman}. The
following definition appears in the authors' work with Lodha~\cite{KKL2019ASENS} (cf.~\cite{CKK2019}).

\begin{defn}
Let $X$ be a topological space and let $G\le \Homeo(X)$ be a subgroup. We say that $G$ acts
\emph{compact-open--transitively}\index{compact-open--transitive}, or
\emph{CO--transitively}\index{CO--transitive},
if for each proper compact subset $A\sse X$ and for each nonempty open $B\sse X$, there is an element
$g\in G$ such that $g(A)\sse B$.
\end{defn}
In the case where $X=S^1$, being CO--transitive is equivalent to being 
\emph{minimal and strongly proximal}\index{minimal and strongly proximal} as defined in~\cite{BF2014ICM}.

\begin{lem}[See~\cite{KKL2019ASENS}]\label{lem:co-trans}
Let $G$ be a group of compactly supported homeomorphisms acting on a noncompact Hausdorff topological space $X$. If $G$ acts
CO--transitively on $X$ then $G'$ is nonabelian and simple, every proper quotient of $G$ is abelian, and the center of $G$ is trivial.
\end{lem}
\begin{proof}
Since $X$ is noncompact and Hausdorff, it follows easily that $G$ is nontrivial. To apply Theorem~\ref{thm:higman}, let $\{r,s,t\}$ be
nontrivial elements of $G$. We may find a nonempty open set $U\sse X$ such that $t(U)\cap U=\varnothing$, since $X$ is
Hausdorff. Let $A$ be a compact set such that \[\supp r\cup\supp s\sse A.\] Since the action of $G$ is CO--transitive, we know that
there is an element $g\in G$ such that $g(A)\sse U$. As $g(A)\cap t(g(A))=\varnothing$, we obtain
\[g(\supp r\cup\supp s)\cap tg(\supp r\cup\supp s)=\varnothing.\] Thus, the hypotheses of Theorem~\ref{thm:higman} are satisfied and
the conclusion of the lemma follows.
\end{proof}

\begin{rem}\label{rem:euclidean}
The property of being CO--transitive transfers from a subgroup to a supergroup.
Since the action of $\Diff^\infty_c(\bR^d)_0$ on $\bR^d$ is CO--transitive,
and so are the actions of the groups in the family
\[
\left\{
\Diff^k_c(\bR^d)_0,\Diff^k_c(\bR^d)\right\}
\]
for all $k\in\bN\cup\{\infty\}$. It follows that
$[G,G]$ is simple  for such a group $G$ in the family. 
Actually, the group
\[\Diff^k_c(\bR^d)_0\]
is known to be simple for $k\ne d+1$, while the case $k=d+1$ is still  open at the time of this writing (cf.~\cite{Mather1,Mather2}). 
\end{rem}

\subsection{The Epstein--Ling Theorem}\label{sss:methods}
As remarked above, in order to apply a result like Higman's Theorem to general diffeomorphism groups of manifolds,
we need to be able to write
a diffeomorphism as a product of ones which have small support. When the manifold in question is a Euclidean space then the topology
of the manifold does not get in the way of CO--transitivity. For a general manifold, one needs to be more careful.

Recall that a Hausdorff topological space $X$ is \emph{paracompact}\index{paracompact}
if for every open cover $\VV$ of $X$, there is a refinement
$\UU$ such that for all $A,B\in\UU$, either $A\cap B=\emptyset$ or there is a third $C\in\VV$ such that $A\cup B\sse C$. We will
always assume that manifolds are paracompact. This kind of refinement is sometimes called an \emph{open star refinement}\index{open
star refinement} of $\VV$.
Regarding the simplicity of the commutator subgroups of certain homeomorphism groups, 
the following topological property played a key role in Epstein's work~\cite{Epstein1970} and in Ling's simplification~\cite{Ling1984}.

\bd
Let $X$ be a topological space.
We say a group $G\le \Homeo(X)$ is  \emph{transitive--inclusive}\index{transitive--inclusive}
if there exists a basis of open sets $\BB$ for the topology of $X$
such that for all basic open sets $U,V\in\BB$
some element $g\in G$ satisfies $g(U)\sse V$.\ed

\begin{thm}[Epstein--Ling~\cite{Epstein1970,Ling1984}]\label{thm:trans-i}
If $X$ is a paracompact space and if a subgroup $G$ of $\Homeo(X)$ is fragmented and transitive--inclusive, then
$[G,G]$ is simple.
\end{thm}

The reader will note the strong resemblance between Theorem~\ref{thm:trans-i} and Lemma \ref{lem:co-trans}, which is not accidental.
Let us give a proof of the above theorem in the remainder of this section.
For this, we let $X$ and $G$ be as in Theorem~\ref{thm:trans-i},
and let $\BB$ be an open basis of $X$ with respect to which $G$ is transitive--inclusive.
We will make an additional assumption that $G$ is nonabelian, without loss of generality.

In the sequel, for a group $G$, we will write $G^{(n)}$ for the $n^{th}$ term of the \emph{derived series}\index{derived series}
of $G$. Thus, $G=G^{0}$,
and $G^{(k+1)}=[G^{(k)},G^{(k)}]$.

\begin{lem}[\cite{CKK2019}, Lemma 3.13]\label{lem:derived}
If $x\in X$, if $n\geq 0$ is an integer, and if $U\sse X$ is an
open neighborhood of $x$ then there is an element $g\in G_U^{(n)}$ such that $g(x)\neq x$.
\end{lem}

Lemma~\ref{lem:derived} in particular implies that $G_U$ cannot be solvable.

\begin{proof}[Proof of Lemma~\ref{lem:derived}]
Since $G$ is nonabelian, $X$ has at least three points. Since $X$ is Hausdorff and transitive--inclusive, the action of $G$ on $X$ does not
have a global fixed point. It follows that there is an element $g\in G$ such that $g(x)\neq x$. The open sets $\{U,X\setminus \{x\}\}$ form an
open cover of $X$, and so we may fragment $g$ with respect to this cover. It follows immediately that there is an element $g_U\in G_U$
such that $g_U(x)\neq x$.

By induction, suppose that we have shown that there exists an element $g\in G_U^{(n-1)}$ that does not fix $x$. There is a neighborhood
$V\sse U$ of $x$ such that $g(V)\cap V=\varnothing$, since $X$ is Hausdorff. Observe that by induction there is an element 
$h\in G_V^{(n-1)}$ such that $h(x)\neq x$. Note that \[g(h(x))\neq g(x)=h(g(x)).\] The product $(ghg^{-1})h^{-1}$ is therefore nontrivial
and lies in $G_U^{(n)}$, which establishes the lemma.
\end{proof}

\begin{lem}[\cite{CKK2019}, Lemma 3.14]\label{lem:higman-epstein}
We have the following.
\begin{enumerate}[(1)]
\item
Every proper quotient of $G$ is abelian.
\item
We have $G'=G''$.
\item
The center of $G$ is trivial.
\end{enumerate}
\end{lem}
\begin{proof}
Let $K\le G$ be a normal subgroup and let $t\in K$ be a nontrivial element. We may find an element $V(t)\in\BB$ such that $t(V(t))\cap V(t)=
\varnothing$. If $V\in \VV$ is arbitrary, we let $f,g\in G_V$ be arbitrary. By transitive--inclusivity, we have that there exists an element
$h\in G$ such that $h(V)\sse V(t)$. In particular, we have \[h(\supp f)\cap t(h(\supp g))\sse h(V)\cap t(h(V))=\varnothing.\] The reader
will note the similarity between this expression and the hypotheses of Theorem~\ref{thm:higman}.

Write $k=(t^h)(g^{-1})(t^h)^{-1}$. Computing, we find that $[f,k]=1$. Moreover, $[f,g]$ lies in the normal closure of $t$, since
\[[f,g]=[f,gk]=[f,[g,t^h]],\] and this latter commutator is a product of conjugates of $t$ and $t^{-1}$. Since $K$ is normal, we have
that $G_V'=G_V^{(1)}\le K$.

If $f,g\in G$ are arbitrary, then we apply fragmentation with respect to the chosen open star refinement $\UU$ of $\VV$. We have that
$f$ is a (finite) product of elements $\{f_i\}_{i\in I}$ with $f_i\in G_{U_i^f}$ for suitable elements $U_i^f\in\UU$, and similarly $g$ is a product
of element $\{g_j\}_{j\in J}$ with $g_j\in G_{U_j^g}$ for suitable $U_j^g\in \UU$.

Standard manipulation of commutators shows that in a group $G$ and for arbitrary elements $f,g,h\in G$,
there exist elements $u,v\in\form{ f,g,h}$ such that \[[fg,h]=[f,h]\cdot [g,h]^u,\quad [f,gh]=[f,g]\cdot [f,h]^v.\] Applying this observation
to the commutator of $[f,g]$, we have that there exist suitable elements
\[\{u_{i,j}\}_{i\in I,j\in J}\sse \form{ \{f_i\}_{i\in I},\{g_j\}_{j\in J}}\] such that $[f,g]$ is a product of commutators of the form
$[f_i,g_j]^{u_{i,j}}$.

Note that the commutators $[f_i,g_j]$ is trivial if $U_i^f\cap U_j^g=\varnothing$. Thus, the commutators is nontrivial only if
$U_i^f\cap U_j^g\neq \varnothing$. Now, we use the paracompactness. The fact that $\UU$ is an open star refinement of $\VV$ implies
that if $U_i^f\cap U_j^g\neq\varnothing$ then there is an open set $V\in\VV$ such that $f_i,g_j\in G_{V_{i,j}}$. We have already shown that
$G_{V_{i,j}}'\le K$, whence it follows that $[f,g]\in K$. It follows that $G'\le K$.

The second claim of the lemma follows from the first. Indeed, $G''$ is nontrivial, since $G$ is nonsolvable by Lemma~\ref{lem:derived}. We
have that $G''$ is a normal subgroup of $G$ and hence $G''$ contains $G'$. The other inclusion is immediate.

The third claim follows similarly. The center $Z$ of $G$ is an abelian normal subgroup and hence must contain $G'$ if it is
nontrivial. Since $G$ is not
solvable, $G'$ is nonabelian. Thus, $Z$ is trivial.
\end{proof}

For each $x\in X$, there is an open $V_x\in\BB$ containing $x$ and an element
$g_x\in G$ such that $g_x(V_x)\cap V_x=\varnothing$. We let \[\VV=\{V_x\mid x\in X\}\sse \BB\] and set $\UU$ to be an open star refinement of $\VV$.

Before proving Theorem~\ref{thm:trans-i}, we will require one more technical lemma which shows that for the cover $\UU$ (i.e.~ the
chosen open star refinement of the cover $\VV$ above), the action of the group $G$ can be simulated by the action of the group $G'$.

\begin{lem}[\cite{CKK2019}, Lemma 3.15]\label{lem:simulate}
Let $G$ be as before, let $g\in G$, and let $U\in \UU$. There exists an element $h\in G'$ such that $g(U)=h(U)$.
\end{lem}
\begin{proof}
We may fragment $G$ for the covering $\UU$, and so there exist open sets \[\{U_1,\ldots,U_m\}\sse \UU\] such that $g$ is a finite
product of elements $g_i\in G_{U_i}$, say $g=g_m\cdots g_1$.

We now produce elements $h_i\in G$ as follows. If $U_i\cap U=\varnothing$ then let $h_i$ be the identity. If $U_i\cap U\neq\varnothing$
then there is an element $V\in\VV$ such that $U_i\cap U\sse V$, and we may choose $h_i$ such that $h_i(V)\cap V=\varnothing$,
so that in particular $U\cap h_i(U_i)=\varnothing$.

Set \[k_i=(g_{i-1}\cdots g_1)h_i\in G.\] Since $h_i(U_i)\cap U=\varnothing$, we obtain that \[k_i(U_i)\cap (g_{i-1}\cdots g_1(U))=\varnothing.\]

It follows from these calculations that \[\supp(k_ig_ik_i^{-1})\cap (g_{i-1}\cdots g_1(U))=\varnothing.\] So,
\[\prod_{i=m}^1[g_i,k_i](U)=\prod_{i=m}^2 [g_i,k_i]\cdot g_1\cdot (k_1g_1^{-1}k_1^{-1})(U)=\prod_{i=m}^2 [g_i,k_i]\cdot g_1(U).\]
By an easy induction on $m$, we have that \[\prod_{i=m}^1[g_i,k_i](U)=g(U),\] which establishes the lemma.
\end{proof}

We  now show that for a fragmented, transitive--inclusive group $G\le \Homeo(X)$, the commutator subgroup $G'$ is simple. 

\begin{proof}[Proof of Theorem~\ref{thm:trans-i}]
It suffices to prove that the normal closure of a nontrivial element of $G'$ coincides with all of $G$.
Retaining the notation from the previous lemmas, let $H$ be the subgroup generated
by $G_{g(U)}^{(2)}$, where $g$ varies over $G$ and $U$ varies over $\UU$. By Lemma~\ref{lem:simulate}, we have that $g(U)=h(U)$ for
some suitable $h\in G'$, so that $H$ coincides with the subgroup generated by $G_{h(U)}^{(2)}$, where $U\in\UU$ and $h\in G'$. Observe that
for $k\in G$, we have that \[kG_{g(U)}^{(2)}k^{-1}=G_{k(g(U))}^{(2)}\le H,\] so that $H$ is in fact normal in $G$. Lemma~\ref{lem:higman-epstein}
then implies that $H=G'$.

Let $s\in G'$ be an arbitrary nontrivial element, which we wish to show normally generates $G'$, where here we mean normally in $G'$
and not normally in $G$. Let $V(s)\in\BB$ be a basic open set such
that $V(s)\cap s(V(s))=\varnothing$. If $U\in\UU$ is arbitrary, it suffices to show that $G_U^{(2)}$ is contained in the normal closure of $s$.
Indeed, then $G_{g(U)}^{(2)}$ will be contained in the normal closure of $s$ for all $g\in G$ and all $U\in \UU$, whence $s$ normally generates
$G'$.

So, let $f,g\in G_U'$ be arbitrary. Transitive--inclusivity implies that there is an element $h\in G$ such that $h(U)\sse V(s)$, and by
Lemma~\ref{lem:simulate}, we may assume that $h\in G'$. We have that \[h(\supp(f))\cap s(h(\supp g))\sse h(U)\cap s(h(U))=\varnothing,\]
so that \[[f,g]=[f,[g,s^h]],\] and as have argued before, this last commutator is a product of conjugates of $s$ and $s^{-1}$. It follows that
$[f,g]$ lies in the normal closure of $s$, as required.
\end{proof}

Now we can finally prove the simiplicity of the commutator subgroup of $\Diff^k_c(M)_0$.

\begin{proof}[Proof of Theorem~\ref{thm:epstein}]
We have that $\Diff^k_c(M)_0$ is fragmented by Lemma \ref{lem:frag}. The manifold $M$ is paracompact by assumption, and
the faithfulness of the action of $\Diff^k_c(M)_0$ is by definition. 
We may choose a basis for the topology of $M$ consisting of the $C^k$--open balls (Definition~\ref{defn:ball}). 
If $U$ and $V$ are such open balls, we wish to find a compactly
supported diffeomorphism of $M$ that sends $U$ into $V$. Since $M$ is connected, there is an embedded path
$\gamma$ from $U$ to $V$, whereby there is a
tubular neighborhood $N$ of $\overline{U\cup V\cup \gamma}$ such that $\overline{N}$ is compact and a diffeomorphism supported in $N$
that sends $U$ to $V$. This diffeomorphism then extends to $M$ by the identity. Thus, the action of $\Diff^k_c(M)_0$ is transitive--inclusive.
The theorem now follows from Theorem~\ref{thm:trans-i}.
\end{proof}

\begin{cor}[{cf.~\cite[Theorem 2.5]{Filip82}}]\label{cor:minimal}
The subgroup \[[\Diff^k_c(M)_0,\Diff^k_c(M)_0]\] is the unique nontrivial minimal normal subgroup
of the four groups below:
\[
\Diff_c^k(M)_0,
\Diff_c^k(M),
\Diff^k(M)_0,
\Diff^k(M).
\]
\end{cor}
\begin{proof}
Let $G_{c0}, G_{c},G_0$ and $G$ denote the given four groups respectively.
Note that $G_{c0}$ is nonabelian, fragmented, and transitive--inclusive.
By
part (1) of Lemma~\ref{lem:higman-epstein}, 
the group $[G_{c0},G_{c0}]$ contains
every nontrivial normal subgroup of $G_{c0}$.
This implies the conclusion for the case of $G_{c0}$.

Assume $G$ is one of the remaining three groups, and let $1\ne N\unlhd G$.
It suffices for us to show that $N$ contains $  [G_{c0},G_{c0}]$.
Pick an $C^k$ open ball $U\sse M$ and $h\in N$ such that $U$ and $h(U)$ are disjoint.
We can find 
\[f,g\in \Diff_c^k(U)_0\le G_{c0}\unlhd G\] such that $[f,g]\ne1$. For instance, 
we can pick distinct points \[\{a,b,c,d\}\sse  U\] and use the path--transitivity (Lemma~\ref{lem:t3}) of $G_{c0}$ to find
$f,g\in\Diff_c^k(U)_0$ such that
\[
f(a)=b, f(b)=c, g(a)=a, g(b)=c.\]
Then we have 
\[ gf(a)=c=f(b)\ne fg(a).\]

Setting $\bar f:=[f,h]$ and $\bar g:=[g,h]$ we see that $\bar f, \bar g\in N\cap G_{c0}$.
Moreover, we have 
\[ [F,G](x)=[f,g](x)\]
for all $x\in U$. It follows that
\[
1\ne [F,G]\in N\cap [G_{c0},G_{c0}]\unlhd  [G_{c0},G_{c0}].\]
Since $ [G_{c0},G_{c0}]$ is simple
by
Theorem~\ref{thm:trans-i}, it must coincide with $N\cap  [G_{c0},G_{c0}] $.
\end{proof}
We will later see that all of the four groups in the above corollary are pairwise non-isomorphic (Corollary~\ref{cor:normal}).

\section{The Bochner--Montgomery Theorem on continuous group actions}\label{sec:boch-mont}
The Hilbert's fifth problem, in one interpretation, asks whether or not every locally Euclidean topological group is a Lie group.
This was confirmed affirmatively by Montgomery and Zippin~\cite{MZ1952} and Gleason~\cite{Gleason1952}. 
An analogous result for a transformation group (that is, an action of a Lie group $G$ on a manifold $X$ which defines a simultaneously
 continuous map $(g,x)\mapsto g.x$) was proved before by Bochner and Montgomery~\cite{BM1945}; we remark this result is often
 attributed to Montgomery and Zippin in the literature, although this is not quite right.

\begin{thm}[Bochner--Montgomery, Theorem 4 of~\cite{BM1945}]\label{thm:BM}
Let $k\in\bN\cup\{\infty\}$, and let $X$ be a $C^k$--manifold.
If a Lie group $G$ acts on $X$ by $C^k$--diffeomorphisms
such that the map $H\co G\times X\longrightarrow X$ given by
\[ H\co (g,x)\mapsto g.x\]
is continuous, then $H(g,x)$ is a $C^k$ map on $(g,x)$ with respect to an analytic parameterization of $g\in G$.\end{thm}

We will give a proof of the above theorem only for the case of Euclidean group actions, based on the original proof of Bochner and Montgomery~\cite{BM1945,MZ1955}. This special case is sufficient for our purpose of proving Filipkiewicz's theorem (Theorem~\ref{thm:filip}).
Moreover, it still captures all the essential ideas of the theorem in the full generality while maintaining relatively simple notation.
As usual, a \emph{manifold}\index{manifold} for us means a Hausdorff,  paracompact, locally Euclidean space.
\begin{thm}[Bochner--Montgomery theorem for Euclidean actions]\label{thm:BM-simple}
Let $d,k\in\bZ_{>0}$, and let $X$ be a smooth connected boundaryless $n$--manifold.
If
\[H\co \bR^d\times X\longrightarrow X\]
is a continuous map inducing an action of $\bR^d$ on $X$ by $C^k$--diffeomorphisms, then $H$ is a $C^k$ map.
\end{thm}
Note that the case $k=\infty$ is an immediate consequence.

\subsection{Analytic prerequisites}
Let us first collect some classical facts on real analysis.
The following result (this version being due to G.~Peano)
describes a condition that allows for switching the order of partial derivatives.
\begin{thm}[{Clairaut's Theorem, \cite[Chapter 9]{Rudin-PMA}}]\label{thm:clairaut}
If $f$ is a real--valued continuous function on $\bR^2$ such that 
the maps
\[
D_1f,\, D_2f,\, D_{12}f
\]
exist and are continuous everywhere,
then $D_{21}f$ exists everywhere and coincides with $D_{12}f$.\end{thm}

Here, the subscripts denote the indices of the variables that are being differentiated.
The interchangeability of integration and differentiation
is an easy consequence of the Clairaut's theorem.
\begin{thm}[{Leibniz Integral Rule, \cite[Chapter 9]{Rudin-PMA}}]\label{thm:leibniz}
If $f(x,y)$ is a real--valued continuous function on $\bR^2$ such that 
$D_1f$
exists and is continuous everywhere,
then for all $a<b$ we have that
\[
\int_a^b D_1f(x,y)dy = \frac{d}{dx} \int_a^b f(x,y)dy.\]
\end{thm}

Recall that a \emph{$G_\delta$ set}\index{$G_{\delta}$ set} is a countable intersection of open sets.
The following classical result is a consequence of the Baire Characterization Theorem of
Baire class one functions and the Baire Category Theorem.
We include a succinct proof, following the lines of~\cite{Kechris1995}, for the convenience of the reader.

\begin{thm}[{Baire's Simple Limit Theorem, cf.~\cite[Theorem 24.14]{Kechris1995}}]\label{thm:baire-simple}
If $\{f_i\}_{i\ge1}$ is a sequence of real--valued continuous functions on a locally compact Hausdorff topological space $Z$ such that the pointwise limit
\[
f(z):=\lim_{i\to\infty} f_i(z)\]
exists for all $z\in Z$, then the continuity set
\[
\operatorname{Cont}(f):=\{z\in Z\mid f\text{ is continuous at }z\}\]
is a dense $G_\delta$ subset of $Z$.\end{thm}

\bp
Pick a countable open basis $\{U_i\}_{i\ge1}$ of $\bR$. Observe that
\[f^{-1}(U_i)=\bigcup\left\{
\bigcap_{j\ge N} f_j^{-1}(\bar U_m)
\middle\vert
{m,N\ge1\text{ and }\bar U_m\sse U_i}
\right\}
.\]
It follows that the set
\[
Z\setminus\operatorname{Cont}(f)
=\bigcup_i \left( f^{-1}(U_i)\setminus \Int f^{-1}(U_i)\right)\]
is a countable union of closed sets with empty interior.
The proof is completed by an application of the Baire Category Theorem.
\ep

\subsection{Standing assumptions}
To prove Theorem~\ref{thm:BM-simple}, we will make the following assumptions. 
We fix positive integers $d$ and $k$, and set
\[
T:=\bR^d\]
We let $X$ be a smooth connected boundaryless $n$--manifold.
The coordinates of points in $T$ will be called \emph{time variables}\index{time variable},
while those in $X$ are called \emph{spatial variables}\index{spatial variable}, specified locally.
We fix a continuous map
\[
H\co T\times X\longrightarrow X\]
which defines a $C^k$--action of $T$ on $X$;
we also write $g(x):=H(g,x)$ for $g\in T$ and $x\in X$ when there is little danger of confusion. We  write
\[H=(H_1,\ldots,H_n)\]
in a suitable local coordinate of $X$.

For an arbitrary function 
\[
F\co T\times X\longrightarrow\bR\]
and for $t=(t_1,\ldots,t_d)$ and $x=(x_1,\ldots,x_n)$,
we write
\begin{align*}
D_jF(t,x)&:=\frac{\partial F}{\partial x_j}(t,x),\\
\bar D_qF(t,x)&:=\frac{\partial F}{\partial t_q}(t,x),
\end{align*}
whenever they are defined in suitable local coordinates.
We denote by $e_j$ the $j^{th}$ standard basis vector in $\bR^d$.

\subsection{On spatial derivatives}
The following  consequence of Baire's Simple Limit Theorem
plays a crucial role in transferring the regularity of spatial variables to that of time variables. 
It roughly says that
\emph{first spatial derivatives are simultaneously continuous at almost every time moment}.
Let us denote $B^n(a;r)\sse\bR^n$ (or $B(a;r)$ more simply) the open ball with radius $r$ centered at $a\in\bR^n$.

\begin{lem}\label{lem:MZ-baire}
If $F$ is a real--valued continuous map on $T\times \bR^n$ such that the map
\[
x\mapsto D_1F(g,x)\]
exists and is continuous for each $g\in T$,
then for each fixed $x_0\in \bR^n$, the set
\[
Z^*(x_0):=\{g_0\in T\mid D_1F(g,x)\text{ is continuous at the point }(g_0,x_0)\}\]
is a dense $G_\delta$ subset in $T$.
\end{lem}
\bp
Let us set $f:=D_1F$ and 
\[
H(g,x,t):=\frac1t\left( F(g,x+te_1)-F(g,x)\right).\]
The map $H$ is continuous on $(g,x,t)$ for $t\ne0$; moreover,
\[
H(g,x,t)=f(g,x+t\theta)\]
for some $\theta=\theta(g,x,t)\in[0,1]$.

Since
\[f(g,x)=\lim_{t\to0} H(g,x,t),\]
Baire's Simple Limit Theorem implies that
the set of continuity points 
\[
\operatorname{Cont}(f)\sse T\times\bR^n\] 
is a dense $G_\delta$ set.
If we denote by $p_1$ and $p_2$ the projections from $T\times\bR^n$ to the two factors, we see that
\[
Z^*(x_0)=p_1\left(p_2^{-1}(x_0)\cap\operatorname{Cont}(f)\right)\]
is a $G_\delta$ subset of $T$ as well.

Let us fix an arbitrary relatively compact open set $Z_0^*$ of $T$.
We claim that there exists some $g_0\in Z_0^*$,
 open neighborhoods $\{\OO_i\}$ of $x_0$,
open  neighborhoods $\{Z_i^*\}$ of $g_0$
 such that 
\[
\diam f(Z_i^*\times\OO_i)\le 1/i.\]
Once this claim is proved, we see that $(g_0,x_0)\in\operatorname{Cont}(f)$
and that $g_0\in Z^*(x_0)$.
Since the choice of $Z_0^*$ is arbitrary, we conclude that $Z^*(x_0)$ is dense.

It only remains for us to show the claim. We fix $r=1/2$.
Since $x\mapsto f(g,x)$ is continuous at each $g\in T$, we have
\[
Z_0^*=\bigcup_{\delta>0} \left\{ g\in Z_0^* \middle\vert
f\left(g,B(x_0,\delta)\right)\sse B(x_0,r/5)\right\}.\]
By the Baire Category Theorem, there exists some $\delta_1>0$
such that for
\[
\OO_1:=B(x_0,\delta_1),\]
the following set is not nowhere dense:
\[
Y_1:=\{g\in Z_0^*\mid f(g,\OO_1)\sse B(x_0,r/5)\}.\]
In other words, we can find a compact neighborhood $Z_1\sse Z_0^*$ such that
$Y_1\cap Z_1$ is dense in $Z_1$.

Applying Baire's Simple Limit Theorem to the sequence $\{H(g,x_0,1/i)\}_i$, we see that $\Int\, Z_1$ contains a continuity point $g_1$
of the map $g\mapsto f(g,x_0)$.
We have an open neighborhood $Z_1^*\sse Z_1$ of $g_1$ such that 
\[
f(Z_1^*,x_0)\sse B\left( f(g_1,x_0),r/5\right).\]

Let $(h,x)\in Z_1^*\times\OO_1$. 
We establish the following inequalities for some $t>0$ and $g\in Y_1\cap Z_1^*$.
\begin{align*}
|f(h,x)-H(h,x,t)|&\le r/5,\\
|H(h,x,t)-H(g,x,t)|&\le r/5,\\
|H(g,x,t)-f(g,x)|&\le r/5,\\
|f(g,x)-f(g,x_0)|&\le r/5,\\
|f(g,x_0)-f(g_1,x_0)|&\le r/5.\end{align*}

The first inequality holds for all sufficiently small $t>0$, depending on $h$ and $x$. We use the density of $Y_1\cap Z_1^*$ in $Z_1^*$
to find $g$ satisfying the second. The third is guaranteed by reducing $t>0$ further, depending on $g$.
The fourth follows from the choice of $\OO_1$. The choices of $g_1$ and $Z_1^*$, for which we used Baire's Simple Limit Theorem,
makes the last inequality hold.
Combining the five inequalities, we have that
\[
f(Z_1^*\times\OO_1)\sse B(f(g_1,x_0),r).\]

Inductively, we can find open sets $Z_i^*$ and elements
$g_i\in T$ such that
\[ g_i\in Z_i^*\sse \overline{Z_i^*}\sse Z_{i-1}^*,\]
and such that
\[
f(Z_i^*\times\OO_i)\sse B(f(g_i,x_0),r/i)\]
for some open neighborhood $\OO_i$ of $x_0$.
Then we have that
\[\diam f(Z_i^*\times\OO_i)\le 2r/i=1/i.\]
Choosing an element
\[
g_0\in \bigcap_i Z_i^*= \bigcap_i \overline{Z_i^*},\]
we obtain the claim.

\ep

We show that \emph{first spatial derivatives are simultaneously continuous}:
\begin{lem}\label{lem:spatial-c1}
The map $D_jH_i$ is continuous on $ T\times X$ for all $i,j$.
\end{lem}
\bp
Let us fix $x_0$. By Lemma~\ref{lem:MZ-baire} we can find densely many $g_0\in T$ such that
$D_jH_i(g,x)$ is continuous at $(g_0,x_0)$ for all $i,j$. 
For a small vector $h\in T$ and for all $x\in X$ we have
\[
H_i(g_0+h,x)=H_i(g_0,h(x)).\]
In a suitable coordinate neighborhood, we have
\[
(D_jH_i)(g_0+h,x)=\sum_k D_k H_i(g_0,h(x))\cdot (D_jH_k)(h,x).\]
Some parentheses above are not strictly necessary for the correct interpretation of the above equation, 
but we have added them nevertheless for readers' convenience.

As $(h,x)$ approaches $(0,x_0)$ the left hand side approaches $D_jH_i(g_0,x_0)$
since $(g_0,x_0)$ is a point of continuity for $D_j H_i$. 
The matrix \[\left( D_jH_i(g_0,h(x))\right)_{i,j=1,\ldots,n}\] is continuous at $(h,x)=(0,x_0)$
since we have made a standing assumption that $H(h,x)$ is continuous,
and that $D_kH_i$ is continuous in the spatial variables;
in particular,  this matrix is nonsingular for $(h,x)$ near $(0,x_0)$
as $x\mapsto H(g_0,x)$ is a $C^k$--diffeomorphism.
It follows that $(0,x_0)$ is a continuity point of $D_jH_k$ for all $x_0\in X$.
Replacing $g_0$ by an arbitrary $g\in T$,
a very similar argument shows that $(g,x_0)$ is a continuity point of $D_jH_k$.\ep

Bootstrapping the above argument, we have that \emph{$k$--th spatial derivatives are simultaneously continuous}:
\begin{lem}\label{lem:spatial-ck}
For $j_1,\ldots,j_k\in\{1,\ldots,n\}$,
the $k$--th partial derivative 
\[D_{j_1}\cdots D_{j_k}H_i\]
 is continuous on $T\times X$.
\end{lem}
\bp
The existence of such a spatial partial derivative is a part of the standing assumption.
We use the Baire's Simple Limit Theorem again to find 
a point $(g_0,x_0)\in T\times X$ at which all such spatial partial derivatives are simultaneously continuous. 
One can then proceeds inductively, in a similar manner as in the proof of Lemma~\ref{lem:spatial-c1};
namely, one establishes the continuity at $(0,x_0)$, and then at $(g,x_0)$ for an arbitrary $g\in T$. This completes the proof.
\ep

\subsection{On time derivatives}
We now turn our attention to the time derivatives, and prove that 
\emph{the first time derivatives  exist and spatially $C^{k-1}$ at time zero}:
\begin{lem}\label{lem:time-c1-0}
For all $i$ and $q$,
the time derivative at time zero
\[\bar D_qH_i(0,x)\]
exists and is $C^{k-1}$ with respect to the spatial variable $x\in X$. 
Moreover, if $h$ is sufficiently close to $0$ then 
for all $x\in X$ we have
\[
H_i(he_q,x)-H_i(0,x)=\sum_j \left(\int_0^h D_jH_i(\tau e_q,x)d\tau\right) \cdot \bar D_q H_j(0,x).\]
\end{lem}
As usual, the above expression and our computation below make sense  in suitable local coordinates.
\bp
We may assume $q=1$. 
Consider the auxiliary function
\[
T_i(h,x):=\int_0^h H_i(\tau e_1,x)d\tau.\]
By Lemma~\ref{lem:spatial-c1} the map $D_jH_i$ exists and is simultaneously continuous.
It follows from the Leibniz Integral Rule (Theorem~\ref{thm:leibniz}) that
\[
D_jT_i(h,x)=\int_0^h D_jH_i(\tau e_1,x)d\tau.\]

For each $t$ sufficiently close to $0$, we can find $z=z_i(t)$
on the line segment between $x$ and $H(te_1,x)$ such that
\[T_i(h,H(te_1,x))-T_i(h,x)
=\sum_j D_j T_i(h,z_i(t))\cdot (H_j(te_1,x)-x_j).\]
The left hand side can be rewritten as
\[\int_t^{h+t} H_i(\tau e_1,x)d\tau - \int_0^h H_i(\tau e_1,x)d\tau
=\int_0^t\left( H_i((h+\tau)e_1,x)-H_i(\tau e_1,x)\right)d\tau.\]

After defining the matrix 
\[
A(h,y)=\left(A(h,y)_{ij}\right):=\left(D_jT_i(h,y)\right),
\]
we conclude that
\begin{equation}\label{eq:aij}
\frac1t \int_0^t\left( H_i((h+\tau)e_1,x)-H_i(\tau e_1,x)\right)d\tau
=
\sum_j A(h,z_i(t))_{ij}\cdot \frac{H_j(te_1,x)-H_j(0,x)}{t}.\end{equation}

The matrix $A(h,z_i(t))$ is arbitrarily close to $h\cdot\Id$ if $h$ and $t$ are sufficiently close to $0$; in particular,
it is non-singular for such $h$ and $t$.
It follows that
\[
\frac{H_i(te_1,x)-H_i(0,x)}{t}
=\sum_j \left(A(h,z_i(t))^{-1}\right)_{ij}\cdot\frac1t \int_0^t\left( H_j((h+\tau)e_1,x)-H_j(\tau e_1,x)\right)d\tau.\]
Sending $t$ to $0$, we obtain that
\[
\bar D_1 H_i(0,x)
=\sum_j \left(A(h,x)^{-1}\right)_{ij} (H_j(he_1,x)-x_j).\]

Since the map
\[
A_{ij}(h,x)=D_jT_i(h,x)=\int_0^h D_jH_i(\tau e_1,x)d\tau\]
is $C^{k-1}$ on $x$, so is $\bar D_1H_i(0,x)$.
By sending $t\to 0$ in the equation~\eqref{eq:aij}, we  obtain 
the second conclusion of the lemma as well.
\ep

Recall our abbreviation
\[
g(x):=H(g,x)\]
for $g\in T$ and $x\in X$.
We now establish that \emph{the first time derivatives of $H$ exist, are simultaneously continuous, 
and are spatially $C^{k-1}$ everywhere}:
\begin{lem}\label{lem:time-c1}
For all $q$ and $i$, the map $\bar D_qH_i$ 
exists and is continuous on $T\times X$; moreover, it
satisfies
\[
\bar D_q H_i(g,x)=\bar D_q H_i(0,g(x))\]
for $(g,x)\in T\times X$;
furthermore, the map $x\mapsto \bar D_qH_i(g,x)$ is $C^{k-1}$ with respect to the spatial variable $x\in X$ for each $g\in T$.
\end{lem}
\bp
We may only consider the case $q=1$.
Using Lemma~\ref{lem:time-c1-0} we have that
\begin{align*}
\frac{H_i(g+he_1,x)-H_i(g,x)}h
&=
\frac{H_i(he_1,g(x))-H_i(0,g(x))}h\\
&=\sum_j \frac1h
\left(\int_0^h D_jH_i(\tau e_1,g(x))d\tau\right) \cdot \bar D_1 H_j(0,g(x)).
\end{align*}
Sending $h\to 0$, we see from 
Lemma~\ref{lem:spatial-c1} 
that
\[
\bar D_1 H_i(g,x)
=\sum_j D_jH_i(0,g(x))\cdot \bar D_1 H_j(0,g(x))
=\bar D_1H_i(0,g(x)).\]
Here, we used that $H(0,y)=y$ for all $y\in X$.
The second conclusion also follows 
since $y\mapsto \bar D_1H_i(0,y)$ is $C^{k-1}$ 
by Lemma~\ref{lem:time-c1-0},
and since $x\mapsto g(x)$ is $C^k$ by the standing assumption.
\ep

We are now ready to complete the proof of the Bochner--Montgomery Theorem for Euclidean actions.
\bp[Proof of Theorem~\ref{thm:BM-simple}]
By Lemmas~\ref{lem:spatial-c1} and~\ref{lem:time-c1}, we have that $H$ is $C^1$ on $T\times X$.
Assume inductively that $H$ is $C^{K-1}$ for some $1\le K\le k$.
We saw that
\[
\bar D_q H_i(g,x)=\bar D_q H_i(0,g(x)).\]
Since $(g,x)\mapsto g(x)$ is $C^{K-1}$
and since the map $y\mapsto \bar D_j H_i(0,y)$ is $C^{k-1}$,
the map $(g,x)\mapsto \bar D_j H_i(g,x)$ is $C^{K-1}$.

It remains to consider a $(K-1)^{th}$  partial derivative mixed with time and spatial variables, say written as $D_0$, of the map $D_jH_i$.
The map $D_jH_i(g,x)$ is $C^{k-1}$ with respect to the spatial variables. If $D_0$ involves a time variable, then the Clairaut's
Theorem (Theorem~\ref{thm:clairaut})
gives the desired conclusion when combined the preceding paragraph; that is, $D_0 D_j H_i$ exists and is continuous.
This proves that $H_i$ is $C^K$, and by induction, that  $H$ is $C^k$.
\ep

\section{Takens' Theorem}\label{sec:takens}
For the purposes of this section, let $M$ and $N$ be
smooth, connected, boundaryless manifolds.
The (generalized) Takens--Filipkiewicz theorem below itself falls slightly short of proving that the group $\Diff^{r}(M)$ determines
$M$ up to $C^{r}$ diffeomorphism,
but which illustrates an important step in proving that fact, namely that a bijection between two manifolds inducing a bijection between
diffeomorphism groups is itself a diffeomorphism; cf.~Section~\ref{sec:filip}. 
The result for full diffeomorphism groups was originally proved by Takens when $p=q=\infty$ and when $M$ and $N$ have the
same dimensions.
R.~P.~Filipkiewicz extended this to all $p,q\in\bZ_{>0}\cup\{\infty\}$ for 
$\Diff^p(M)$ and $\Diff^q(N)$.

\begin{thm}[Generalized Takens' Theorem]\label{thm:takens}
Let $M$ and $N$ be smooth connected boundaryless manifolds, and let $p,q\in\bZ_{>0}\cup\{0,\infty\}$.
Suppose we have groups $G$ and $H$ satisfying
\[
\Diff_c^{p}(M)_0\le G\le \Diff^{p}(M)\]
and
\[
\Diff_c^{q}(N)_0\le H\le \Diff^{q}(N).\]
If a set--theoretic bijection 
 $w\colon M\longrightarrow N$ satisfies
\[ w G w^{-1} =H,\]
then we have that
$p=q$
and that
$w$ is a $C^p$ diffeomorphism.
\end{thm}

As observed by Takens, Theorem~\ref{thm:takens} implies the following statement: in the group of set--theoretic bijections
$M\longrightarrow M$, the group $\Diff^p(M)$ is equal to its own normalizer. 
We will present two different proofs of the theorem, after establishing that $w$ is a homeomorphism by elementary considerations.
The first proof in this book uses the Bochner--Montgomery Theorem, based on~\cite{Filip82}. This readily gives a simple and full
account of Theorem~\ref{thm:takens2}. 

For the second proof, which comes from Takens' original article~\cite{Takens79}, we reduce the
case to maps between Euclidean spaces and then argue locally that the bijection $\phi$ is a diffeomorphism.
This gives a more computationally explicit argument. 
We present only a part of the proof for the sake of brevity, avoiding the introduction of many new concepts of smooth dynamics; namely,
after reducing to the Euclidean case we will consider only the dimension one case with full diffeomorphism groups.

\begin{rem}
It is very interesting to note that both of the proofs presented here are based on  non-constructive existence theorems from real analysis.
The first proof is deduced from  the Bochner--Montgomery theorem, which in turn is based on Baire's Simple Limit Theorem.
The second is based on Lebesgue's theorem on the almost everywhere differentiability of monotone functions, together with
Sternberg's Linearization Theorem, which in turn is a consequence of the Schauder--Tychonoff Fixed Point Theorem
(see Chapter V of~\cite{DSbook1}).
In both proofs, one finds a reasonably regular point by applying the suitable nonconstructive result, 
and then ``transfers'' this regularity to other points using the transitivity of group actions. 
\end{rem}

\subsection{Promoting a bijection to a homeomorphism}\label{ss:bij-homeo}
The following lemma in particular implies the ``easy'' part of Takens--Filipkiewicz Theorem. Namely, it asserts that the map
$w$ given in the hypothesis of the Takens' Theorem must be a homeomorphism.

\begin{lem}\label{lem:bij-homeo}
Let $X_1$ and $X_2$ be topological spaces, 
and let $G_i\le\Homeo(X_i)$ be locally transitive for $i=1,2$.
If there exists an isomorphism
\[
\Phi\co G_1\longrightarrow G_2\]
and a bijection
\[
w\co X_1\longrightarrow X_2\]
such that
\[wgw^{-1}=\Phi(g)\]
for all $g\in G_1$,
then $w$ is a homeomorphism.
\end{lem}
\bp
Let $\UU_i$ be a local basis as in the definition of local transitivity for $i=1,2$.
Consider an arbitrary $x\in X_1$ and its open neighborhood $U\in\UU_1$.
By the local transitivity one has some $g\in G_1$ such that
\[
x\in \suppo g\sse \suppc g\sse U.\]
This implies that the collection
\[
\{\suppo g\mid g\in G_1\}\]
is an open basis of $X_1$. 
Furthermore, we have that
\[
w(\suppo g)
=\suppo wgw^{-1}
=\suppo \Phi(g)\]
for all $g\in G_1$. This implies that $w$ is open. By symmetry, 
we conclude that $w$ is a homeomorphism.\ep

\subsection{Takens' theorem via Bochner--Montgomery}
It is relatively a simple task to deduce Theorem~\ref{thm:takens} 
from the Bochner--Montgomery Theorem (Theorem~\ref{thm:BM-simple}),
which we do in this subsection. We will prove the following stronger result,
and show how this implies Theorem~\ref{thm:takens}.

\begin{thm}\label{thm:takens2}
Let $p\in\bZ_{>0}\cup\{\infty\}$, and let $M$ and $N$ be smooth, connected, boundaryless  manifolds.
If $w\colon M\longrightarrow N$ is a homeomorphism such that 
\[ w^{-1} \Diff^\infty_c(N)_0 w  \sse \Diff^p(M),\]
then $w^{-1}$ is a $C^p$ map.
\end{thm}

\bp
We let $d=\dim M=\dim N$. 
As the statement is purely local, it suffices to prove the desired regularity $C^p$ at a given fixed point, say $x_0\in M$ and $y_0:=w(x_0)$.
Pick a sufficiently small smooth open ball $U\sse N$ containing $y_0$ so that $U$ and $w^{-1}(U)$ are 
contained in charts of $N$ and $M$. By Lemma~\ref{lem:affine} we can find 
a smooth embedding $h\co \bR^d\longrightarrow N$
such that $h(\bR^d)=U$ and such that $h(0)=y_0$.
This embedding induces a smooth action $\bR^d$ on $N$.
Namely, for $v\in \bR^d$ and $y\in N$ we define 
\[
B_v(y):=
\begin{cases} 
h\left(v+h^{-1}(y)\right),&\text{ if }y\in U,\\
y,& \text{otherwise}.\end{cases}\]

Now we can investigate the relationship between $w$ and $B$. 
For $v\in\bR^d$ and $x\in M$, we define
\[
A_v:=w^{-1}\circ B_v\circ w \in\Diff^p(M).\]
The continuity of the map
\[
(v,x)\mapsto A_v(x)= w^{-1}\left( 
h\left(v+h^{-1}(w(x))\right)
\right)
\]
 is obvious. 
Bochner--Montgomery Theorem implies that the map
\[ (v,x)\mapsto A_v(x)\] 
must be a $C^p$ map.

Note that $w^{-1}(B_v(y_0))=A_v(x_0)$ for all $v\in\bR^d$.
Since
the map
\[
B_v(y_0)\mapsto v\]
is a smooth chart about  $y_0$ 
and 
the map 
\[
v\mapsto A_v(x_0)\]
is $C^p$, we conclude that $w^{-1}$ is a $C^p$ map.
\ep

We can now complete the first proof of the Takens' Theorem.

\bp[Proof of Theorem~\ref{thm:takens}]
Since $\Diff^\infty_c(M)_0$ and $\Diff^\infty_c(N)_0$ are locally transitive, we have that $w$ is a homeomorphism by
Lemma~\ref{lem:bij-homeo}.
Applying Theorem~\ref{thm:takens2}, we also have that $w$ and $w^{-1}$ are $C^q$ and $C^p$ maps, respectively. 
By symmetry, we may assume $p\le q$, so that $w$ is a $C^p$ diffeomorphism. 

It only remains to show $p=q$. 
If $p<q$, then the conjugation by the $C^p$--diffeomorphism $w$ induces a group isomorphism between
$\Diff^p_c(M)_0$ and $\Diff^p_c(N)_0$.
This is a restriction of $\Phi$, whose image $H$ is contained in  $\Diff^q(N)$. This is contradiction,
since $\Diff^p_c(N)_0\not\sse \Diff^q(N)$ for $q>p$.
\[\begin{tikzcd}
G\arrow{r}{\Phi} \arrow[dash]{d}{} &
H=H\cap  \Diff^q(N)\\
\Diff^p_c(M)_0\arrow[swap]{r}{f\mapsto wfw^{-1}}& \Diff^p_c(N)_0
\end{tikzcd}
\]
 We thus conclude that $p=q$.
\ep

\subsection{Addendum: Takens' original proof}
In this section, we give an alternative, partial proof of Theorem~\ref{thm:takens} based on Takens' original argument~\cite{Takens79}.
Namely, we will give a proof of the following result for $p\ge2$ and $\dim M=\dim N=1$. The hypothesis on $p$ arises from the
hypotheses in Sternberg's Linearization Theorem.

\begin{thm}[Takens  Theorem~\cite{Takens79}]\label{thm:takens3}
Let $M$ and $N$ be smooth connected boundaryless manifolds.
If  $w\colon M\longrightarrow N$
 is a set--theoretic bijection satisfying that
 \[
w \Diff^\infty(M)w^{-1}= \Diff^\infty(N),\]
then $w$ is a smooth diffeomorphism.\end{thm}

Takens' proof of Theorem~\ref{thm:takens3} can be broken up into several smaller pieces.
The first is to consider the case $M=\bR^d$, as is a common first step in such results.
This step requires some care to avoid issues with exotic smooth structures on Euclidean spaces when $d=4$.

\begin{thm}[\cite{Takens79}, Theorem 2]\label{thm:takens-reduction}
If a bijection $\phi\colon\bR^d\longrightarrow\bR^d$ satisfies that
\[
\phi^{-1}\Diff_+^\infty(\bR^d)\phi=\Diff_+^\infty(\bR^d),\]
then $\phi$ is a smooth--diffeomorphism.
\end{thm}

For the purposes of Theorem~\ref{thm:takens3}, the global topology of $M$ and $N$ are irrelevant. This is made precise by the following
fact:

\begin{lem}\label{lem:takens-reduction}
Theorem~\ref{thm:takens3} follows from Theorem~\ref{thm:takens-reduction}.
\end{lem}

The second step is to prove that if $d=1$ then for $M=N=\bR$, the map $\phi$ in Theorem~\ref{thm:takens3} is
automatically a diffeomorphism.
This is an essentially one--dimensional argument, as it relies fundamentally on the almost everywhere differentiability of a monotone
real--valued function on $\bR$.

The third step is the general case, which uses the second step as a bootstrap. The proof is less self-contained, and relies on invariant
manifold theory. It is not possible for us to give a self--contained proof of the general case
here, as it would take us too far afield into smooth dynamics.
We will therefore content ourselves to giving a proof of Theorem~\ref{thm:takens3} in dimension one only.

\subsubsection{Reduction to Euclidean spaces}
The proof of Lemma~\ref{lem:takens-reduction} is mostly formal, though as we alluded to above, there is one issue to be avoided, namely
the existence of exotic differential structures on $\bR^d$. 

To begin, we note from Lemma~\ref{lem:bij-homeo} that $\phi$ is necessarily a homeomorphism, for $M$ and $N$ arbitrary. 
If $W\sse M$ is diffeomorphic to $\bR^d$, then clearly
$\phi(W)$ is homeomorphic to $\bR^d$. The possibility of exotic smooth structures precludes concluding that $\phi(W)$ is diffeomorphic to $\bR^d$ and thus that $\phi$ is a diffeomorphism.

\begin{proof}[Proof of Lemma~\ref{lem:takens-reduction}]
It suffices to show that for every point $x\in M$, there is an open neighborhood $W$ containing $x$ such that $\phi$ is a diffeomorphism
when restricted to $W$. Without loss of generality, we may assume that $W$ is diffeomorphic to $\bR^d$. We write $V=\phi(W)\sse N$.
Let $D\sse V$ be a closed ball, and let $E=\phi^{-1}(D)$. We let $\psi\in\Diff^\infty(M)$ be such that $E\sse\psi(E)^{\circ}$, the
interior of $\psi(E)$, and such that \[W=\bigcup_{i\geq 0}\psi^i(E).\] The conditions on $\phi$ imply that $\phi\psi\phi^{-1}=\eta$ is
a diffeomorphism of $N$ such that $D\sse\eta(D)^{\circ}$, and \[V=\bigcup_{i\geq 0}\eta^i(D).\] Because $\eta^i(D)$ is diffeomorphic to
a closed ball for all $i$ and is compactly contained in the interior of $\eta^{i+1}(D)$, we see from the
Isotopy Extension Theorem (\cite{Palais1960,Lima64}) that
\[\eta^{i+1}(D)\setminus\eta^i(D)^{\circ}\cong S^{d-1}\times [0,1].\] From this it follows easily
that $V$ is in fact diffeomorphic to $\bR^d$. This argument is symmetric with respect to replacing $\phi$ by its inverse,
and so we may conclude that $\phi$ is a diffeomorphism.
\end{proof}

\subsubsection{Theorem~\ref{thm:takens-reduction} for $n=1$}
Establishing Theorem~\ref{thm:takens} for $\bR$ (and hence for $1$--manifolds) again breaks into two steps. First, we show that
if $\phi$ is as in the hypotheses, then $\phi$ and $\phi^{-1}$ both have a first derivative everywhere in $\bR$.
Then, we bootstrap this fact to conclude that
$\phi$ and $\phi^{-1}$ are smooth. For the remainder of this subsection, we will assume that $\phi$ is as in
Theorem~\ref{thm:takens}, with $M=N=\bR$.

\begin{lem}[\cite{Takens79}, Lemma 3.1]\label{lem:phi-der}
For all $x\in \bR$, the functions $\phi$ and $\phi^{-1}$ both have nonzero first derivative at $x$.
\end{lem}
\begin{proof}
Since $\phi$ is a homeomorphism, we have that $\phi$ and $\phi^{-1}$ are monotone
real-valued functions on $\bR$.
A standard fact about such monotone functions is that they are differentiable almost everywhere.

If $f\colon\bR\longrightarrow\bR$ is a diffeomorphism and $x\in \bR$, then $\phi^{-1}f\phi$ is differentiable at $x$ and has a nonzero
derivative. Recording this fact, we have that
\[(\phi^{-1}f\phi)'(x)=\lim_{h\to 0}\frac{(\phi^{-1}f\phi)(x+h)-(\phi^{-1}f\phi)(x)}{h}.\]

The ratio in the limit can be formally written as a product of three ratios, all of which are defined for $h\neq 0$:
\[
q_1=\frac{(\phi^{-1}f\phi)(x+h)-(\phi^{-1}f\phi)(x)}{(f\phi)(x+h)-(f\phi)(x)},\] \[q_2=\frac{(f\phi)(x+h)-(f\phi)(x)}{\phi(x+h)-\phi(x)},\]
\[q_3=\frac{\phi(x+h)-\phi(x)}{h}.\]
We will argue that for all three ratios, the limits
as $h\to 0$ exist and are nonzero. The continuity of $\phi$ and the fact that $f$ is a diffeomorphism of $\bR$ implies that the limit
as $h\to 0$ of $q_2$ exists and is nonzero, since we are merely expressing a change of coordinates for $f$.

Since $\phi$ is differentiable almost everywhere, we may select a point $x$ so that the limit for $q_3$ exists. Suppose
first that $\phi'(x)=0$. We then compute the limit of $q_1$, and conclude that if $\phi'(x)=0$ then $\phi^{-1}$ does not have a derivative
at $(f\phi)(x)$. Since $f$ is arbitrary and diffeomorphisms of $\bR$ act transitively on $\bR$, we obtain that $\phi^{-1}$ is nowhere
differentiable. This contradicts the fact that $\phi^{-1}$ is almost everywhere differentiable. So, $\phi'(x)\neq 0$, whence $\phi^{-1}$
has a finite and nonzero first derivative at $(f\phi)(x)$. Again using the transitivity of the diffeomorphism group, we conclude that $\phi^{-1}$
has nonzero first derivative everywhere. Symmetrically, we find the same conclusion for $\phi$.
\end{proof}

We now upgrade Lemma~\ref{lem:phi-der} to prove that $\phi$ is actually $C^{\infty}$ with a $C^{\infty}$ inverse. We have the following
technical construction:

\begin{lem}[\cite{Takens79}, Lemma 3.2]\label{lem:phi-conj}
Let $2\le p\le\infty$, and let $\phi\colon \bR\longrightarrow\bR$ be a homeomorphism such that $\phi$ and $\phi^{-1}$ are differentiable for all points in $\bR$.
Suppose furthermore that $f_1$ and $f_2$ are $C^{p}$ orientation preserving diffeomorphisms of $\bR$ such that:
\begin{enumerate}[(1)]
\item
We have $f_1(0)=0$;
\item
We have $f_1'(0)<1$;
\item
For each $x\in\bR$, we have that $f_1^k(x)\longrightarrow 0$ as $k\longrightarrow\infty$;
\item
We have $\phi^{-1}f_1\phi=f_2$.
\end{enumerate}
Then $\phi\in\Diff^{p}(\bR)$.
\end{lem}

In the statement of Lemma~\ref{lem:phi-conj}, $\phi$ is not assumed to be a $C^1$ diffeomorphism, since we do not assume that the
derivatives of $\phi$ and $\phi^{-1}$ are continuous.

\begin{proof}[Proof of Theorem~\ref{thm:takens-reduction} for $n=1$, assuming Lemma~\ref{lem:phi-conj}]
Let $f_1(x)=x/2$ and let $\phi$ be as in the statement of the theorem. Then Lemma~\ref{lem:phi-der} implies that $\phi$ and $\phi^{-1}$
are differentiable at all points of $\bR$. By assumption, we have that \[f_2=\phi^{-1}f_1\phi\in\Diff^{\infty}(\bR).\]
Lemma~\ref{lem:phi-conj} immediately implies the desired conclusion.
\end{proof}

Before we can give the proof of Lemma~\ref{lem:phi-conj}, we will require a result due to Sternberg, which implies that the diffeomorphisms
$f_1$ and $f_2$ in the statement of the lemma are $C^{p}$--conjugate to linear diffeomorphisms. We state Sternberg's result here
in the very special case of one--manifolds, and we omit the proof since the statement itself is believable,
since accessible proofs are available in the literature, since a full proof would
take us far afield, and since we have already given a complete proof of Takens' result, and since the result will
also be subsumed by Filipkiewicz's result in the next section.

\begin{thm}[See~\cite{Sternberg57,Sternberg58}]\label{thm:sternberg}
Let $2\le p\le\infty$.
If $f\colon\bR\longrightarrow\bR$ is a $C^p$ diffeomorphism with $f(0)=0$ and $|f'(0)|\neq 1$, then = the
diffeomorphism $f$ is locally $C^{p}$--conjugate to
a linear diffeomorphism. That is, there exists a $C^p$ diffeomorphism $\psi$ such that $\psi f\psi^{-1}$ is given by a linear map on some
neighborhood $U$ that contains $0$.
\end{thm}

\begin{proof}[Proof of Lemma~\ref{lem:phi-conj}]
Recall we have assumed that $p\ge2$.
We first prove that $\phi$ is a $C^p$ diffeomorphism in a neighborhood of $0$.
By replacing $\phi$ by a translate if necessary, we may assume that $\phi(0)=0$. So, we have that $f_2(0)=0$.
Since dynamics are
preserved under conjugacy, we have that $f_2^k(x)\longrightarrow 0$ for all $x\in \bR$. 
Using the Chain Rule, we see that
\[
D_0:=f_1'(0)=f_2'(0)<1.\]

Theorem~\ref{thm:sternberg} implies the existence
of diffeomorphisms $\psi_i\in\Diff^{p}(\bR)$ for $i\in\{1,2\}$ that fix $0$ and that conjugate $f_i$ to a linear map in a neighborhood
$U$ of the origin, i.e.~
$\psi_i f_i\psi_i^{-1}$ is linear on $U$ for $i\in\{1,2\}$.

We define a map $\eta\colon\bR\longrightarrow\bR$ by $\psi_1\phi\psi_2^{-1}$. It is clear from the definitions that $\eta$ and
$\eta^{-1}$ are both differentiable at zero. We claim that $\eta$ is linear near the origin, whence it will follow that $\phi$ is a
$C^p$ diffeomorphism near the origin.

Note that $D_0$ is the slope of the linearization of $f_i$ for
$i\in\{1,2\}$.
Without loss of generality we may assume that $1\in U$, and we write $\eta(1)=\alpha$.
Choose $\beta\in U$, and write $\gamma=\eta(\beta)$. 
Considering the $k$--th power $D_0^k$ of $D_0$, we have
\begin{align*}
\eta(D_0^k)=\eta(\psi_2f_2^k\psi_2^{-1}(1))=\psi_1\phi\psi_2^{-1}\psi_2f_2^k\psi_2^{-1}(1)=\\
=\psi_1f_1^k\psi_1^{-1}\psi_1\phi\psi_2^{-1}(1)=D_0^k\eta(1)=D_0^k\alpha.
\end{align*}
Similarly, we obtain that  $\eta(D_0^k\cdot\beta)=D_0^k\cdot\gamma$.

Since $D_0<1$, we may compute as follows:
\begin{align*}
\eta'(0)&=\lim_{k\to\infty}\frac{\eta(D_0^k)-\eta(0)}{D_0^k}=\alpha,\\
\eta'(0)&=\lim_{k\to\infty}\frac{\eta(D_0^k\cdot \beta)-\eta(0)}{D_0^k\cdot\beta}=\frac{\gamma}{\beta}.
\end{align*}

It follows that $\alpha=\gamma/\beta$. Since $\beta$ is arbitrary
and $\eta$ is continuous, we have that $\eta(\beta)=\beta\eta(1)$ for all $\beta\in U$. This implies that $\eta$ is linear in a neighborhood of
$0$, as claimed.

To prove that $\phi$ is a $C^p$ diffeomorphism everywhere, we simply change coordinates. If $f_1$ is a diffeomorphism which satisfies
the hypotheses of the lemma with $y$ in place of $0$, we can translate $\phi$ in order to have $\phi$ fix $y$ as well.
Then, we can conjugate the entire picture by a translation to move the fixed point to $0$
and without changing the absolute value of the derivative of $f_1$.
\end{proof}

This completes the proof of Theorem~\ref{thm:takens} for $1$--manifolds and concludes our discussion of Takens' result.

\section{Rubin's Theorem}\label{sec:rubin}
In this section, we give a (simplified, for our purposes)
self--contained proof of Rubin's powerful reconstruction theorem~\cite{Rubin1989,Rubin1996}. This  will
then give an efficient proof of Filipkiewicz's Theorem.

Let $G$ be a group acting on a space $X$. 
Recall our notation below that is used throughout this book:
\begin{align*}
\supp g&:=X\setminus\Fix g,& \text{ if }g\in G,\\
G[U]&:=\{g\in G\mid \supp g\sse U\},& \text{ if }U\sse X.
\end{align*}
The group $G[U]$ is sometimes called the \emph{rigid stabilizer}\index{rigid stabilizer} of $U\sse X$.

\bd\label{d:locally-moving}
Let $X$ be a topological space, and let $G\le\Homeo(X)$.
\be[(1)]
\item
We say $G$ is \emph{locally moving}\index{locally moving action} 
if for each nonempty open set $U\sse X$ the group $G[U]$ is nontrivial.
\item
We say $G$ is \emph{locally dense}\index{locally dense action} 
if for each point $x\in X$ and for each open neighborhood $U$ of $x$, 
the closure of the orbit $G[U].x$ has nonempty interior.
\ee
\ed
Note that a locally dense action is locally moving if the space is
\emph{perfect}\index{perfect space} (i.e.~has no isolated points).

\begin{thm}\label{thm:rubin}
Let $X_1$ and $X_2$ be perfect, locally compact, Hausdorff topological spaces,
and let $G_i\le\Homeo(X_i)$ be  locally dense groups, for $i=1,2$. Suppose that there exists an isomorphism of groups
\[
\Phi\co G_1\longrightarrow G_2.\]
Then there exists a homeomorphism \[\phi\co X_1\longrightarrow X_2\]
such that for all $g\in G_1$ and for all $x\in X_1$ we have
\[
\phi(g.x)=\Phi(g).\phi(x).\]\end{thm}

An account of this result, following the original arguments given by Rubin~\cite{Rubin1989,Rubin1996}, will occupy the remainder
of this section.
\subsection{First order expressibility of rigid stabilizers}
Suppose $A$ is a subset of the space $X$.
We denote the interior and the closure of a set $A\sse X$ by 
$\Int_X A=\Int A$ and $\cl_X A=\cl A$, respectively. 
The notion of the \emph{extended support}\index{extended support} of a homeomorphism will be essential for our discussion, 
and we define it as
\[\suppe g:=\Int\cl \supp g,\]
for $g\in G$.
Note the inclusion
\[
\supp g\sse \suppe g\sse \cl{\supp g}.\]

For a subset $A$ of a group $G$, we let
\[
Z_G(A):=\{h\in G\mid [a,h]=1\text{ for all }a\in A\}, \text{ if }A\sse G.\]
We also write $Z_G(g):=Z_G(\{g\})$.

We remark that the reader may think of $Z_G(A)$ as a \emph{definable set}\index{definable set} in the sense of model theory 
(see~\cite{hinman-book,marker-book,tz-book} for background, for instance). This is not truly a definable set in the sense of classical
model theory in general since $A$ is allowed
to be infinite, and so there may not be a single formula with parameters in $A$ that defines $Z_G(A)$.
But this subtle point will
not cause any trouble for us since we will only consider the case where $A$ is defined by a first order formula.

A remarkable fact that is crucially used in the proof of Rubin's Theorem is that, 
for a locally moving group $G$, 
the rigid stabilizer of the extended support of a $g\in G$ can be expressed in 
terms of a first order formula in the language of group theory (i.e.~the language that admits a single binary operation and
a distinguished constant corresponding to the identity; sometimes it is useful to include the [definable] inversion function).
To state this sentence concretely, for each $f\in G$ we define 
\begin{align*}
\eta_G(f):=\{g\in G \mid & \text{ for all }h\in G\setminus Z_G(f)\text{ there exist }f_1,f_2\in Z_G(g)\\
&\text{ such that }1\ne [[h,f_1],f_2]\in Z(g)\}.\end{align*}

We remark that for each $f\in G$, the set $\eta_G(f)$ is indeed definable from the parameter $f$, and the reader is encouraged
to write the sentence down in formal syntax.

\begin{thm}\label{thm:rubin-supp}
If $G$ is a locally moving group of homeomorphisms on a Hausdorff topological space,
then for all $f\in G$ we have that
\[
G[\suppe f]=Z_G\left(\left\{ g^{12}\mid g\in \eta_G(f)\right\}\right).\]
\end{thm}


We will establish Theorem~\ref{thm:rubin-supp} in this subsection. 
Let us first make elementary observations from point set topology.
For a set $A$ in a space $X$, we let $\comp_X A=\comp A=X\setminus A$,
and 
\[
\ext_X A=\ext A=A^\perp=\comp \cl  A.\]
The latter of the above is called the 
\emph{exterior}\index{exterior of a set} of $A$.
We say $A$ is a \emph{regular open set}\index{regular open set} of $X$ if 
\[
A = A^{\perp\perp}.\]
We let $\Ro(X)$ denote the set of all regular open sets.

\begin{lem}\label{lem:ro-basic}
Let $X$ be a topological space. 
\be[(1)]
\item For all subset $A$ of $X$, we have
\begin{align*}
 A^{\perp}&= \Int \comp A,\\
 A&\sse A^{\perp\perp}=\Int \cl A.
 \end{align*}
\item\label{p:ppp} If $A\sse X$ is open, then $A^{\perp}=A^{\perp\perp\perp}$; in particular, we have $A^\perp\in\Ro(X)$.
\item If $A\sse X$ is open,  then $\cl A = \cl \left(A^{\perp\perp}\right)$.
\item\label{p:ab-disj} If $A$ and $B$ are disjoint open subsets of $X$,
then $A^{\perp\perp}$ and $B$ are also disjoint.
\item\label{p:uv-perp} 
Let $U\in\Ro(X)$.
If an open set $V\sse X$ satisfies
that $V\cap U^\perp=\varnothing$,
then $V\sse U$.
\item\label{p:uv-haus}
If $X$ is  Hausdorff  and if $x_0$ and $x_1$ are distinct points in $X$
then there exist a disjoint pair of regular open neighborhoods $U_i$ of $x_i$ for $i=0,1$.
\ee
\end{lem}

\bp
All are simple to check. For instance, to see (\ref{p:ab-disj}), it suffices to note
\[A^{\perp\perp}=\Int \cl A\sse \Int \cl \comp B=\Int \comp B\sse \comp B.\]
For (\ref{p:uv-perp}), we note
\[
V=\Int V \sse\Int \comp U^{\perp}= U^{\perp\perp}=U.\]
To see (\ref{p:uv-haus}), pick a disjoint pair of open neighborhoods $V_i$ of $x_i$ for $i=0,1$. By part (\ref{p:ab-disj}) we have that
\[
V_i\cap V_{1-i}^{\perp\perp}=\varnothing.\]
Parts (\ref{p:ppp}) and (\ref{p:uv-perp}) then imply that 
\[ V_i \sse V_{1-i}^\perp.\]
Setting $U_i:=V_i^{\perp\perp}\cap V_{1-i}^\perp\in\Ro(X)$, 
we see that
\[
U_0\cap U_1\sse V_0^{\perp\perp}\cap V_0^{\perp}=\varnothing,\]
and that
\[ x_i\in V_i\cap V_{1-i}^{\perp}\sse U_i.\]

The rest are simple mental exercises, the fun of which we will not ruin for the reader.
The book \cite{GH2009} is also a good reference. 
\ep

By the lemma above, we see that $\suppe g$ is a regular open set whenever $g\in\Homeo(X)$.
Let us now note some basic facts about locally moving group actions.

\begin{lem}\label{lem:LM-basic}
Let $G$ be a locally moving group faithfully acting on a Hausdorff topological space $X$.
\be[(1)]
\item\label{p:move-supp}
If $U$ is a nonempty open subset of $\supp f$  for some $f\in G$ then there exists some nontrivial element 
$g\in G[U]\setminus Z_G(f)$ such that $\supp g\cap f\supp g=\varnothing$.
\item\label{p:move-cap}
If $A$ is a nonempty open subset of $\bigcap_{i=1}^m \supp g_i$ for suitable elements \[\{g_1,\ldots,g_m\}\sse  G,\] then there exists a 
nonempty open subset $B\sse A$ such that 
\[
B\cap \left(\bigcup_{i=1}^m g_i (B)\right)=\varnothing.\]
\item\label{p:move-fn}
For each open set $U\sse X$ and for each nonzero integer $n$, there exists an element $g\in G[U]$ and a 
nonempty open set $V\sse U$ such that 
\[
\{V, g(V), g^2(V),\ldots, g^n(V)\}\]
is a disjoint collection of open sets.
\item\label{p:move-fg}
If $f,g\in G\setminus\{1\}$ satisfy $\supp f\sse\supp g$,
then there exists some $h\in G[\supp f]$ such that $[hfh^{-1},g]\ne1$.
\item\label{p:move-sse}
For two regular open sets $U,V$ of $X$, we have that
\[
U\sse V\Longleftrightarrow G[U]\le G[V].\]
\ee
\end{lem}
\bp
(\ref{p:move-supp}) Pick $x\in U$ such that $f(x)\ne x$. We can find some open neighborhood $V$ of $x$ such that $V\cap f(V)=\varnothing$.
Since $G$ is locally moving we have some $g\in G[V]\le G[U]$. Then 
\[
\supp g\cap f\supp g\sse V\cap f(V)=\varnothing.\]
Since $f\supp g\ne\supp g$, we have $g\not\in Z_G(f)$.

(\ref{p:move-cap})
Pick $a\in A$. Since $a\ne g_i(a)$ for each $i$, we can find some open neighborhood $B_i$ of $a$ such that $B_i\cap g_i(B_i)=\varnothing$.
Then $B:=\bigcap_i B_i$ satisfies the conclusion.

(\ref{p:move-fn}) 
This will be strengthened in Theorem~\ref{thm:lawless}.
Let us give a short proof here.
Since $U$ is Hausdorff and $G[U]$ is locally moving,  it suffices to consider the case that $U=X$.
The case $n=1$ is trivial since $X$ is Hausdorff. 
Let us assume the conclusion for $n$, and inductively prove the case $n+1$. 
We have some nonempty open set $V\sse X$ and some element $g\in G$ such that
$\{g^i(V)\}_{0\le i\le n}$ is a disjoint collection.
If we have some $v\in V$ such that $g^{n+1}(v)\ne v$, then some open neighborhood of $v$ satisfies the conclusion for $n+1$.
Hence, we may assume that $g^{n+1}\restriction_V=\Id_V$. Pick an arbitrary element $h\in G[V]$, and a suitable nonempty open 
set $W\sse V$ such that $W\cap h(W)=\varnothing$.

We note that 
\[ (gh)^i(W)=g^i\circ h(W)\sse g^i(V)\]
forms a disjoint collection
for $i\in\{0,\ldots,n\}$. Since
\[
(gh)^{n+1}(W)\cap W = g^{n+1}h(W)=h(W)\]
is disjoint from $W$, we have that
\[
\{W,gh(W),\ldots,(gh)^{n+1}(W)\}\]
is a disjoint collection. 
We just found an open set $W$ satisfying the conclusion with the element $gh\in G$.

(\ref{p:move-fg}) We may assume $[f,g]=1$, for otherwise we can set $h:=1$. 
By the above, we can choose some nonempty open set $A\sse\supp f$ such that
\[
A\cap(f(A)\cup g(A))=\varnothing.\]
Depending on whether or not $\supp (fg)$ intersects $A$, 
we can find a nonempty open set $B\sse A$ such that 
one of the following holds:
\begin{itemize}
\item
 $B\cap fg(B)=\varnothing$;
\item
 $\supp (fg)\cap B=\varnothing$.\end{itemize}
Pick a nontrivial element $h\in G[B]$ and nonempty open set $C\sse B$ such that 
\[
\{C,h(C),h^2(C)\}\]
is a disjoint collection; see Lemma~\ref{lem:LM-basic} (\ref{p:move-fn}).
By the choice of $A$, we see that
\begin{align*}
 x_1&:=g\cdot hfh^{-1}(hC)=ghf(C)=gf(C)=fg(C),\\
 x_2&:=hfh^{-1}\cdot g(hC)=hfg(hC).\end{align*}
In the case when $B\cap fg(B)=\varnothing$, we have that 
\[x_2=fg(hC)\ne fg(C)=x_1.\]
If $\supp (fg)\cap B=\varnothing$, then we have
\[x_1=fg(C)=C\ne h^2(C)= hfgh(C)=x_2.\]
This completes the proof.

(\ref{p:move-sse}) The forward direction is trivial. For the reverse direction, 
assume that $U\not\sse V$. By Lemma~\ref{lem:ro-basic} (\ref{p:uv-perp}), we have that 
$W:=U\cap V^{\perp}\ne\varnothing$. By the locally moving hypothesis, we have some $1\ne f\in G[W]\le G[U]$. 
Since \[W\sse V^{\perp}=\comp\cl V\sse \comp V,\]
we see that $\supp f\cap V=\varnothing$. In particular, we have $G[U]\not\sse G[V]$.
\ep

Note the following simple combinatorial fact.
\begin{lem}\label{lem:sets}
Let $m\ge1$ be an integer and let  $Z$ be a set.
Suppose a collection of subsets \[\{A_{ij}\sse Z\}_{0\le i\le m,1\le j\le m}\]
has the property that for each $j=1,\ldots,m$
the sets
\[
A_{0j}, A_{1j},\ldots,A_{mj}\]
are pairwise disjoint.
Then we have that
\[\bigcap_{i=0}^m \left(\bigcup_{j=1}^m A_{ij}\right)=\varnothing.\]
\end{lem}
\bp
Let $x\in Z$. We imagine to have $m$ pigeons and $m+1$ pigeonholes.
We let the $j$--th pigeon stay in the $i$--th hole if $x\in A_{ij}$; it is clear from the hypothesis that each pigeon can stay in at most one hole.
For some $0\le i\le m$, the $i$--th hole will be vacant.
This implies that $x$ does not belong to $\bigcup_j A_{ij}$.\ep


The following lemma asserts that the group theoretic condition $g\in\eta_G(f)$ ``almost'' detects the topological condition that $f$ and $g$ have disjoint supports.

\begin{lem}\label{lem:eta}
Let $G$ be a locally moving group faithfully acting on a Hausdorff topological space $X$.
For two elements $f,g\in G$ we have the following conclusions.
\be[(1)]
\item
If $\supp f\cap\supp g=\varnothing$, then $g\in \eta_G(f)$.
\item
If $g\in\eta_G(f)$ then 
\[\supp f\cap \supp g^3\cap \supp g^4=\varnothing.\]
\ee
\end{lem}
We will often omit the subscript $G$ in $Z_G$ or  $\eta_G$ when the meaning is clear. 
\bp[Proof of Lemma~\ref{lem:eta}]
(1) In order to show $g\in\eta(f)$, pick an arbitrary $h\in G\setminus Z(f)$.
In particular, we have $\supp f\cap \supp h\ne\varnothing$.
By Lemma~\ref{lem:LM-basic} (\ref{p:move-supp}) we can pick  elements $f_1,f_2\in G$
 satisfying the following conditions:
 \begin{itemize}
 \item $f_1\in G[\supp f\cap \supp h]\setminus Z(h)$ and $\supp f_1 \cap h\supp f_1=\varnothing$;
 \item $f_2\in G[\supp f_1]\setminus Z(f_1)$ and $\supp f_2 \cap f_1\supp f_2=\varnothing$.
 \end{itemize}
 Since $\supp f_2\sse \supp f_1\sse \supp f$, we have $f_1,f_2\in Z(g)$.
Note that 
\[\supp hf_1h^{-1}=h\supp f_1\] is disjoint from $\supp f_1$.
It follows that 
\[
[[h,f_1],f_2]=[hf_1h^{-1}\cdot f_1^{-1},f_2]=[f_1^{-1},f_2]
=\left(f_1^{-1}f_2f_1\right)\cdot f_2^{-1}\in Z(g).\]
The last term of the equality is nontrivial, since it is written as the product of two nontrivial homeomorphisms with 
disjoint supports. This shows that $g\in \eta_G(f)$.

(2) Let us assume $g\in \eta(f)$. Assume for contradiction that
\[ A:=\supp f\cap \supp g^3\cap \supp g^4\ne\varnothing.\]
In particular, we have $A\sse \bigcap_{i=1}^4 \supp g^i$. 
By Lemma~\ref{lem:LM-basic} (\ref{p:move-cap}) there exists a 
nonempty open set $B\sse A$ satisfying
\[
B\cap \left( \bigcup_{i=1}^4 g^i(B)\right)=\varnothing.\]
It follows that $\{g^i (B)\}_{0\le i\le 4}$ is a disjoint family of nonempty open sets.
By the same lemma, we can also find an element $h\in G[B]\setminus Z(f)$. 
Using the hypothesis that $g\in\eta(f)$, we know of the existence of elements $f_1, f_2\in Z(g)$ such that 
\[ u:=[[h,f_1],f_2]\in Z(g)\setminus\{1\}.\]
Setting $S:=\{1,f_1,f_2,f_2f_1\}\sse Z(g)$, we easily see that
\[
\bigcup_{i=0}^4 g^i(\supp u) = \supp u
=
\supp \left(h\cdot (f_1 h^{-1}f_1^{-1})\cdot (f_2f_1 h f_1^{-1}f_2^{-1})(f_2 h^{-1}f_2^{-1})\right)\sse  S(B),\]
and that
\[
\supp u\sse \bigcap_{i=0}^4 \bigcup_{s\in S} g^{-i}s(B).\]
This contradicts Lemma~\ref{lem:sets}, since for each $s\in S$ the collection
\[\{ g^{-i}s(B)=sg^{-i}(B)\}_{0\le i\le 4}\]
consists of disjoint sets.
\ep

We are now ready to express the set $G[\suppe f]$ in the language of group theory. It will be convenient for us to introduce the notation
\[\xi_G^m(f):=Z_G\left(\left\{ g^m\mid g\in \eta_G(f)\right\}\right),\]
which is defined for every group $G$ and every element $f\in G$.

\bp[Proof of Theorem~\ref{thm:rubin-supp}]
Let $m\in\bZ_{>0}$.
Let us first pick arbitrary $q\in G[\suppe f]$ and $g\in \eta(f)$.
Since \[\supp g^{12}\sse \supp g^3\cap\supp g^4,\] we see from Lemma~\ref{lem:eta} that
\[\supp f\cap\supp\left(g^{12}\right)=\varnothing.\]
Lemma~\ref{lem:ro-basic} implies that
\[\suppe f\cap\supp\left(g^{12}\right)=\varnothing.\]
Since $\supp q\sse \suppe f$, we see that $q$ commutes with $g^{12}$.
We have thus shown that \[G[\suppe f]\le \xi_G^{12}(f)\le\bigcup_{m\ge1}\xi_G^m(f).\]

To see the opposite inclusion, pick arbitrary $m\in\bZ_{>0}$ and $q\in G\setminus G[\suppe f]$.
Since $\supp q\not\sse\suppe f$, Lemma~\ref{lem:ro-basic} (\ref{p:uv-perp}) implies that the open set
\[
A:=\supp q\cap \left(\suppe f\right)^\perp\]
is nonempty. We can therefore pick some $g\in G[A]$ such that $g^m\ne1$, 
using Lemma~\ref{lem:LM-basic} (\ref{p:move-fn}).
Summarizing, we have
\[
\varnothing\ne \supp g^m\sse \supp g\sse A\sse \supp q.\]
Using part (\ref{p:move-fg}) of the same lemma, we can pick 
some $h\in G\left[\supp g^m\right]$ such that 
$q$ does not commute with $hg^mh^{-1}$.

We claim that the elements $f$ and $hgh^{-1}$
have disjoint supports. 
Indeed, 
applying Lemma~\ref{lem:ro-basic}, we see that
\begin{align*}
\supp hgh^{-1}&=h\supp g=\supp g
\sse A\sse (\suppe f)^\perp\\
&=(\supp f)^{\perp\perp\perp}=(\supp f)^{\perp}
\sse \comp\supp f.\end{align*}
The claim is thus proved.

Lemma~\ref{lem:eta} implies that $hgh^{-1}\in\eta(f)$,
and hence that $q\not\in \xi_G^m(f)$. 
This implies that $G\setminus G[\suppe f]\sse G\setminus \xi_G^m(f)$, for an arbitrary $m\in\bZ_{>0}$.

So far, we have proved that 
\[
\bigcup_{m\ge1} \xi_G^m(f)\le G[\suppe f]\le \xi_G^{12}(f)\le\bigcup_{m\ge1} \xi_G^m(f).\]
This completes the proof.

\ep

\begin{rem}\label{rem:rubin-supp}
\be[(1)]\item In the above   we have  established a purely group theoretic statement \[
Z_G\left(\left\{ g^{m}\mid g\in \eta_G(f)\right\}\right)\sse Z_G\left(\left\{ g^{12}\mid g\in \eta_G(f)\right\}\right)\]
for all $m\ge1$. The mysterious number 12 seems to be intrinsically related to the locally moving hypothesis of $G$ (and of course
comes out of the proof of Theorem~\ref{thm:rubin-supp} as the least common multiple of $3$ and $4$).
\item
We used in the proof the fact that if $\supp q\setminus \suppe f\ne\varnothing$
then $\supp q\setminus \supp f$ contains a nonempty open set. 
The same is not true under the weaker hypothesis that $\supp q\setminus\supp f\ne\varnothing$.
This is one of the places where an extended support turns out to be more useful for us than an open support. 
\ee
\end{rem}

From now on, we write
\[
\xi_G(f):=\xi_G^{12}(f)=Z_G\left(\left\{ g^{12}\mid g\in \eta_G(f)\right\}\right)\le G.\]
We can now express inclusions between the sets
\[\{\suppe f\mid f\in G\}\] purely in terms of group theory.
\begin{cor}\label{cor:rubin-supp}
Let $G$ be a locally moving group of homeomorphisms on a Hausdorff topological space $X$. Then for all $f,g\in G$,
we have $\suppe f\sse \suppe g$ if and only if $\xi_G(f)\le\xi_G(g)$.
\end{cor}
\bp
This is a trivial consequence of Theorem~\ref{thm:rubin-supp} and Lemma~\ref{lem:LM-basic} (\ref{p:move-sse}).
\ep

\subsection{From groups to Boolean algebras}\label{sec:gp-bool}
The next step in the proof of Rubin's Theorem is to extract a bijection between regular open sets of spaces that respects their
 \emph{Boolean structures}\index{Boolean structure}. 

Let $B$ be a set containing two distinguished elements $0$ and $1$. 
If there exist binary operations $\wedge,\vee$
and a unary operation $\perp$ satisfying the following natural set-theoretic axioms, then we say the structure 
$(B,\wedge,\vee,\perp,0,1)$ is a \emph{Boolean algebra}\index{Boolean algebra}. 
\begin{itemize}
\item both of $\wedge$ and $\vee$ are commutative and associative;
\item $\wedge$ and $\vee$ are distributive, in the sense that
\begin{align*}
u\wedge(v\vee w) &=(u\wedge v)\vee(u\wedge w),\\
u\vee(v\wedge w) &=(u\vee v)\wedge(u\vee w).\end{align*}
\item we have $u\vee(u\wedge v)=u=u\wedge(u\vee v)$.
\item $u\vee u^\perp=1$ and $u\wedge u^{\perp}=0$.
\end{itemize}
Note that
 $0,1,\wedge,\vee,\perp$ are axiomatized to reflect the behavior of the set theoretic objects $\varnothing,X,\cap,\cup,\comp$, respectively.
See~\cite{GH2009} or~\cite{Jech2003} for a concise introduction to Boolean algebras.

A Boolean algebra $B$ naturally comes with 
a partial order $\le$ defined
by
\[v\le u\Longleftrightarrow v\wedge u^{\perp}=0.\]

It is helpful to regard $v\wedge u^{\perp}$ as the ``subtraction'' of $u$ from $v$. 
If a Boolean algebra $B$ admits a least upper bound and a greatest lower bound for every nonempty 
subset $A\sse B$ then we say $B$  is \emph{complete}\index{complete Boolean algebra}.
We say $B$ is \emph{atomless}\index{atomless Boolean algebra} if it does not contain a nonzero minimal element.

For the purpose of this book, one may only consider the Boolean algebra of regular open sets defined as follows.
Let $X$ be a Hausdorff topological space. 
Recall that we denote the collection of regular open sets in $X$ by $\Ro(X)$, namely the sets $U\sse X$ satisfying
\[ U =  \Int\cl U=U^{\perp\perp}.\]
The set $\Ro(X)$ is equipped with the unary operation $\perp$, and a binary operation $\vee$ defined by
\[
U\vee V:=(U\cup V)^{\perp\perp}.\]
Setting \[\wedge:=\cap,\quad 0:=\varnothing,\quad 1:=X,\] 
we have that $\Ro(X)$ acquires the structure of a Boolean algebra.

A particularly useful fact for us is that the ``subtraction'' of two regular open sets can be defined by
\[
V\cap U^{\perp}\in\Ro(X).\]
We have seen in Lemma~\ref{lem:ro-basic} (\ref{p:uv-perp}) that if $V\not\sse U$ for regular open sets $U,V$ then $V\cap U^{\perp}\ne\varnothing$.
In other words, the Boolean partial order for $\Ro(X)$ coincides with the set theoretic inclusion $\sse$.

Let us note that a self-homeomorphism induces a Boolean automorphism, in the following strong sense.
\begin{lem}\label{lem:top-bool}
Let $X$ is a Hausdorff topological space, and let $G\le\Homeo(X)$.
If $U\in\Ro(X)$,
then we have that
\[G[U] = \{g\in G\mid g(V)=V\text{ for all regular open set }V\sse U^\perp\}.\]
\end{lem}
\bp
If $g\in G[U]$, then $g$ restricts to the identity on $U^\perp$ and so, $g$ belongs to the right handside.
For the reverse direction, suppose $g\not\in G[U]$. By Lemma~\ref{lem:ro-basic} we can find some  $z\in \supp g\setminus U^\perp$.
By the same Lemma~\ref{lem:ro-basic} we can pick a pair of disjoint regular open sets $V_0$ and $V_1$ of $z$ and $g(z)$, respectively. 
Since $g(z)\in g(V_0)$ we see that $g(V_0)\ne V_0$. This implies that $g$ does not belong to the right hand side, completing the proof.\ep

By setting $U=\varnothing$ we see that the group $\Homeo(X)$ acts naturally and faithfully on the Boolean algebra $\Ro(X)$ 
by Boolean automorphisms. For $G\le\Homeo(X)$ we continue to denote this new action on $\Ro(X)$ still as $G$. 
In particular, the group $G[U]\le G$ can be described purely in terms of the Boolean algebra $\Ro(X)$, 
the Boolean automorphic action of $G$ and its element $U$.

For a family $F\sse \Ro(X)$ we define the supremum and the infimum of this family by
\begin{align*}
\sup F&:=\left(\bigcup F\right)^{\perp\perp}\in \Ro(X),\\
\inf F&:=\left(\bigcap F\right)^{\perp\perp}\in \Ro(X).
\end{align*}
The Boolean algebra $\Ro(X)$ is {complete},
since these are indeed the least upper bounds and the greatest lower bounds
with respect to the inclusion. 
If we further assume that $X$ is a perfect, then $\Ro(X)$ is atomless.
Observe that if $G$ is a locally moving group of homeomorphisms on $X$ then $X$ must be perfect.

We now prove that an abstract group isomorphism between two locally moving groups of homeomorphisms of Hausdorff topological spaces
induces a Boolean isomorphism of regular open sets of those spaces, that is equivariant with respect to this group isomorphism.

\begin{thm}\label{thm:gp-bool2}
Let $X_1$ and $X_2$ be Hausdorff topological spaces.
If a group isomorphism \[\Phi\co G_1\longrightarrow G_2\] is given,
where  $G_i\le\Homeo(X_i)$ are locally moving subgroups for $i=1,2$,
then there uniquely exists a Boolean isomorphism  \[\Psi\co \Ro(X_1)\longrightarrow \Ro(X_2)\]
such that the following diagram of set-theoretic maps are commutative for all $g\in G_1$:
\[\begin{tikzcd}
G_1 \arrow{r}{\suppe}\arrow{d}[swap]{\Phi}  & \Ro(X_1)\arrow{r}{g}\arrow{d}{\Psi}& 
\Ro(X_1)\arrow{d}{\Psi}\arrow{r}{G_1[\cdot]} & \operatorname{Subgroups}(G_1)\arrow{d}{\Phi}\\
G_2 \arrow{r}{\suppe}  & \Ro(X_2)\arrow{r}{\Phi(g)}& \Ro(X_2)\arrow{r}{G_2[\cdot]} & \operatorname{Subgroups}(G_2).\end{tikzcd}
\]
\end{thm}

The first horizontal maps are sending $g$ to $\suppe g$, for $g$ in $G_1$ or $G_2$. 
The third horizontal maps send a regular open set $U\in\Ro(X_i)$ to the corresponding rigid stabilizer $G_i[U]$. 
The bijection $\Psi$ above respects the Boolean objects and operations
\[ \cap, \vee,\perp,\varnothing,X_i.\]
Note that suprema and  infima are respected as well, since the partial order $\sse$ is respected.
The horizontal maps in the middle square are actions of $G_i$ on the Boolean algebras $\Ro(X_i)$.

The deduction of Theorem~\ref{thm:gp-bool2} from Theorem~\ref{thm:rubin-supp}  
is relatively simple and may be regarded as a Boolean algebra ``language game''. We remark one logical subtlety:
Theorem~\ref{thm:rubin-supp} characterizes rigid stabilizers in terms of first order formulae with parameters, and
Theorem~\ref{thm:gp-bool2} follows in some sense because of ``preservation of logical implication". Whereas the language of
group theory is retained, the logic that is being used to establish the implication is different.
Indeed, since suprema and quantification over subgroups will be
required, passage to second order logic is warranted in the proof of
Theorem~\ref{thm:gp-bool2}.

Suppose $(Q,\le)$ is a partially ordered set. For each $q\in Q$ we write
\[Q\restriction_q:=\{x\in Q\mid x\le q\}.\]
The main example that are considered here is the case when $Q$ is the set of nonzero elements in a Boolean algebra.

We say that a subset $P$ of $(Q,\le)$ is  \emph{dense}\index{dense subset of a partial order} 
if every $q$ in $Q$ satisfies $Q\restriction_q\cap P\ne\varnothing$.
If $B$ is a complete Boolean algebra and if $P\sse B\setminus\{0\}$ is a dense subset,
then every element $b\in B\setminus\{0\}$ can be expressed as 
\[ b:=\sup (B\restriction_b\cap P);\] cf.~Chapter 25 of~\cite{GH2009}.

We note the following result on Boolean algebras.

\begin{lem}\label{lem:ba-complete}
Let $B$ be a complete Boolean algebra,
and let $P\sse B\setminus\{0\}$ be a dense subset.
If there exists an order--preserving dense embedding 
\[f\co P\longrightarrow B'\setminus\{0\}\]
for some complete Boolean algebra $B'$,
then $f$ uniquely extends to a Boolean isomorphism from $B$ to $B'$.
\end{lem}

Lemma~\ref{lem:ba-complete} can be succinctly expressed by the slogan that
every \emph{separative}\index{separative subset of a partially ordered set}
dense partially ordered set admits a unique completion to a Boolean algebra. Here, a \emph{separative} partially ordered 
set is one that is a dense subset in the set of nonzero elements of a Boolean algebra. We omit the proof, which is elementary to verify. 
The interested reader is directed to~{\cite[Theorem 14.10]{Jech2003}} for a further discussion.

\bp[Proof of Theorem~\ref{thm:gp-bool2}]
For each $i=1,2$ we let
\[
P_i:=\{\suppe g\mid g\in G_i\setminus\{1\}\}\sse \Ro(X_i)\setminus\{\varnothing\}.\]
Since $G_i$ is locally moving, we see that $P_i$ is dense in $\Ro(X_i)\setminus\{\varnothing\}$. Define a map
\[\Psi\co P_1\longrightarrow P_2\]
by the formula
\[
\Psi(\suppe g):=\suppe \Phi(g).\]
Recall the notation $\xi_G(g)$ from the previous subsection.
Using the fact that
\[\Phi(\xi_{G_1}(g))=\xi_{G_2}(\Phi(g)),\]
we can apply Corollary~\ref{cor:rubin-supp} 
and deduce that $\Psi$ is a well-defined, order--preserving bijection
between $P_1$ and $P_2$. 
Since $\Ro(X_i)$ is complete, we see from Lemma~\ref{lem:ba-complete}
that $\Psi$ uniquely extends to a Boolean isomorphism 
\[\Psi\co \Ro(X_1)\longrightarrow\Ro(X_2).\]
The first square in the diagram commutes by definition.
The third one commutes since for $g\in G_1$ and $U\in\Ro(X_1)$ we have
\[g\in G_1[U]
\Longleftrightarrow \suppe g\sse U
\Longleftrightarrow \suppe \Phi(g)\sse \Psi(U)
\Longleftrightarrow \Phi(g)\in G_2[\Psi(U)].\] 

It suffices now to prove that $\Psi$ intertwines the two actions.
For $\suppe f\in P_1$ and $g\in G_1$ we note 
\[
\Psi(g.\suppe f)=\Psi(\suppe gfg^{-1})=\suppe \Phi(gfg^{-1})=\Phi(g)\suppe\Phi(f).\]
Note also that $P_i$ is dense $\Ro(X_i)$, and that $\Psi$ respects the suprema.
Since $\Ro(X_i)$ is complete, so is the proof.

\ep

\subsection{From Boolean algebras to topologies}
So far, we have extracted a Boolean isomorphism between regular open sets of Hausdorff topological spaces,
under the hypothesis that they admit isomorphic locally moving groups. We now consider the stronger hypotheses of local compactness 
of the spaces, together with local density of groups, and deduce that the spaces are in fact homeomorphic. For this, 
we will need to identify a point in terms of regular open sets, thus crucially using the notion of an ultrafilter.

Let $B=B(\wedge,\vee, \perp,0,1)$ be a Boolean algebra. Recall that a \emph{filter}\index{filter} is a nonempty subset $F\sse B$ such that 
the following two conditions hold: 
\be[(i)]
\item If $a,b\in F$ then $a\wedge b\in F$;
\item If $a\in F$ and $a\le b$ for $b\in B$, then $b\in F$.
\ee
A proper filter $F\subsetneq B$ is an \emph{ultrafilter}\index{ultrafilter} if it satisfies the following additional condition:
\be[(i)]
\addtocounter{enumi}{2}
\item For all $b\in B$ either $b\in F$ or $b^\perp\in F$.
\ee
It is routine to check that a proper filter of $B$ is an {ultrafilter} if and only if it is maximal. 

A subset $F_0$ of $B$ is said to have the \emph{finite intersection property}\index{finite intersection property for a filter}, if we have
\[
b_1\wedge \cdots \wedge b_k\ne0\]
for all $k\ge1$ and $b_i\in F_0$. It is a standard fact that each subset of $B$ satisfying the 
finite intersection property is contained in an ultrafilter of $B$~\cite[Chapter 20, Exercise 12]{GH2009}. Proofs of this fact
use the Axiom of Choice, but the existence of ultrafilters is weaker than the Axiom of Choice.

The \emph{Stone space}\index{Stone space} $\CS(B)$ is the set of all ultrafilters in $B$, equipped with the subspace topology of $2^B$.
Although the topology of $\CS(B)$ is not of our concern in this book,
we note that this space is compact and totally disconnected; it is even perfect if $B$ is atomless.

The following lemma hints at the idea of identifying points in a space with ultrafilters of regular open sets.
The set of all  open neighborhoods of a point $x$ in a space $X$ will be denoted by $\Nbr_X(x)=\Nbr(x)$.

\begin{lem}\label{lem:haus-filter}
Let $X$ be a Hausdorff topological space. For each ultrafilter $F$ on the Boolean algebra $\Ro(X)$, we define
\[ 
A_F:=\bigcap_{U\in F}\cl U.\]
The following conclusions hold.
\be[(1)]
\item\label{p:single} The set $A_F$ contains at most one point.
\item If $F$ contains a relatively compact regular open set, then $A_F\ne\varnothing$. 
\item If $x\in A_F$ then $\operatorname{Nbr}(x)\cap\Ro(X)\sse F$.
\item\label{p:fx} For each $x\in X$ there exists some $F\in \Ro(X)$ such that $x\in A_{F}$. 
\ee
\end{lem}
\bp
(1) Suppose $x,y$ are distinct points in $A_F$. We can find disjoint regular open neighborhoods $U,V$ of them by part 
(\ref{p:uv-haus}) of Lemma~\ref{lem:ro-basic}. We may choose $V:=U^\perp$. 
By maximality of the ultrafilter, we have either $U\in F$ or $U^\perp\in F$. If $U\in F$ then we have that $\cl U\cap U^\perp=\varnothing$
and that $y\not\in A_F$, a contradiction. The case $U^\perp\in F$ is similar.

(2) Let $V\in F$ be relatively compact. 
Since $F$ satisfies the finite intersection property as a filter,
 the family of compact sets
\[
\{ \cl (U\cap V) \mid U\in F\}\]
satisfies the finite intersection property as a collection of subspaces in $X$.
It follows that
\[
\varnothing
\ne
\bigcap
\{ \cl (U\cap V) \mid U\in F\}\sse A_F.\]

(3) Let $V\in\Nbr(x)$ be regular open. For each finite subcollection \[\{U_1,\ldots,U_k\}\sse F,\]
we have that $U_0:=\bigcap_i U_i\in F$. Moreover, we have that
\[x\in A_F\cap V \sse(\cl U_0)\cap V.\]
It follows that $U_0\cap V\ne\varnothing$. We have that 
$F\cup\{V\}\sse \Ro(X)$ satisfies the finite intersection property.
Since $F$ is an ultrafilter, we see that $V\in F$.

(4) Using the finite intersection property 
one can find  an ultrafilter $F_x$ containing the set $\Nbr(x)\cap\Ro(X)$.
Assume for contradiction that the closure of some $U\cap F_x$
does not contain $x$. In other words, we have $x\in U^\perp$.
Then we have that $U^\perp\in \Nbr(x)\cap\Ro(X)\sse F_x$, which contradicts that $F_x$ is an ultrafilter.
\ep

The preceding lemma justifies the following definition.
\bd
Let $X$ be a Hausdorff topological space.
We say $F\in\CS(\Ro(X))$ is a \emph{topological filter}\index{filter!topological} of $X$ if 
\[\bigcap_{U\in F}\cl U\neq\varnothing.\] We denote $\TF(X)$ the set of all topological filters of $X$.
We write
\[
\rho_X\co \TF(\Ro(X))\longrightarrow X\]
for the surjection defined by the formula
\[\{\rho_X(F)\}=\bigcap_{U\in F}\cl U.\]
\ed

A topological filter is defined by the topology of the space being considered.
Rubin defined the following collection of ultrafilters, which is given purely in terms of a group action on a Boolean algebra. 
It will turn out that those two concepts coincide in the setting of Rubin's Theorem.

\bd
Let $B$ be a Boolean algebra, and let $G$ be a group of Boolean automorphisms of $B$. We say an ultrafilter $F$ of $B$ is a 
\emph{$G$--local filter}\index{local filter} if there exists some $b\in F$ such that \[B\restriction_b\sse \{0\}\cup G.F.\] 
We will write $\CS(B;G)$ for the set of all $G$--local filters in $B$.
\ed

It will be convenient for us to introduce the following construction.
\begin{notation}
Let $H$ be an automorphism group of a Boolean algebra $B$.
For an ultrafilter $F\in\CS(B)$, we define
\[
H\left\{F^\perp\right\}:=
\bigcup_{w\in F}
\{
g\in H\mid 
g(a)=a\text{ for all }a\in  B\restriction_w\}.\]\end{notation}
This is indeed a group, as can be seen by applying the finite intersection property of $F$.
In the case where $G\le\Homeo(X)$, we can also regard $G$ as a group of Boolean automorphisms of $\Ro(X)$ as in
Lemma~\ref{lem:top-bool}.
Then, for an ultrafilter $F\in\CS(\Ro(X))$, we obtain
\[
G\left\{F^\perp\right\}=
\bigcup_{W\in F} G\left[W^\perp\right].\]

The following theorem on a single space $X$ will provide us a connection from a Boolean isomorphism to a topological homeomorphism.

\begin{thm}\label{thm:bool-homeo}
Suppose $X$ is a locally compact, Hausdorff, perfect topological space,
and suppose that $G$ is a locally dense subgroup of $\Homeo(X)$.
Then the following conclusions hold.
\be[(1)]
\item\label{p:tf} We have that \[ \TF(X) = \CS(\Ro(X),G).\]
\item\label{p:rho-x}
For each $x\in X$, we have that
\[
\Nbr(x)\cap\Ro(X)
=\bigcap \rho_X^{-1}(x)
.\]
\item\label{p:rhoi-rho} 
For each $F\in \TF(X)$, we have that
 \[
 \TF(X)\setminus\rho_X^{-1}\circ \rho_X(F)
 = \CS\left(\Ro(X), G\left\{F^\perp\right\}\right).\]
\ee
\end{thm}

\begin{rem}When we use the notation $\CS(\Ro(X),H)$ for some group $H\le\Homeo(X)$, 
we are regarding $H$ as a group of Boolean automorphisms of $\Ro(X)$; this notation (as well as its ambiguity) is 
justified by Lemma~\ref{lem:top-bool}.\end{rem}

\bp[Proof of Theorem~\ref{thm:bool-homeo}, parts (\ref{p:tf}) and (\ref{p:rho-x})]
(\ref{p:tf}) 
Suppose $F\sse\Ro(X)$ is $G$--local filter of $\Ro(X)$.
By definition, we can find some $U\in F$ such that
\[\Ro(X)\restriction_U \sse \{\varnothing\}\cup G.F.\]
By local compactness, we can find a relatively compact nonempty regular open set $V\sse U$. 
Since $V\in G.F$, we see from Lemma~\ref{lem:haus-filter} that  $F\in\TF(X)$.

Conversely, suppose $F\in \TF(X)$ and let $x:=\rho_X(F)$.
By local density, we have a nonempty regular open set $U$ contained in the closure of $G.x\sse X$. 
Consider an arbitrary $V\in \Ro(X)\restriction_U\setminus\{\varnothing\}$.
Since $V$ is inside the closure of $G.x$, we can find a $g\in G$ such that
$g.x\in V$. 
By the same lemma as above, we have that
\[
x\in g^{-1}(V) \in \Nbr(x)\cap \Ro(X)\sse F.\]
This implies that $F$ is $G$--local.

(\ref{p:rho-x}) 
We see from Lemma~\ref{lem:haus-filter} that 
\[\Nbr(x)\cap\Ro(X)\sse F\]
for all $F\in \TF(X)$ satisfying  $\rho_X(F)=x$.
Conversely, suppose $U$ is an element of $\bigcap \rho_X^{-1}(x)$.
Assume for a contradiction that $x\not\in U$. 
For an arbitrary regular open neighborhood $V$ of $x$, we have that
$V\not\sse U$, which implies $V\cap U^\perp=\varnothing$.
The family
\[
F_0:=\left(\Nbr(x)\cap\Ro(X)\right)\bigcup \left\{ U^\perp\right\}\]
enjoys the finite intersection property, and hence extends to some $F\in \CS(\Ro(X))$.
Note that the intersection of $\cl V$ for all $V\in F_0$ is already a singleton, namely $\{x\}$. 
This implies that $\rho_X(F)=x$, and by hypothesis, $U\in F$.
This implies $\{U,U^\perp\}\sse F$, which is a contradiction.

\ep

We postpone the proof of the part (\ref{p:rhoi-rho}) above, first observing a general fact.
\begin{lem}\label{lem:pointed}
Let $X$ be a perfect Hausdorff topological space.
Fix a point $\infty\in X$, and let $Y:=X\setminus\{\infty\}$.
Then the following conclusions hold.
\be[(1)]
\item
There exists a Boolean isomorphism 
\[\psi\co \Ro(X)\longrightarrow\Ro(Y)\]
satisfying
\begin{align*}
\psi(U)&:=U\cap Y,&\text{ if }U\in\Ro(X),\\
\psi^{-1}(V)&:=\Int_X\cl_X V,&\text{ if }V\in\Ro(Y).\end{align*}
and inducing a bijection
\[
\TF(X)\setminus\rho_X^{-1}(\infty) \longrightarrow \TF(Y).\]
\item If $G\le\Homeo(X)$ is locally dense, then for each  $F\in\rho_X^{-1}(\infty)$
 the group
\[
G\left\{F^\perp\right\}\] fixes the point $\infty$ and acts locally densely on $Y$.
\ee
\end{lem}
\bp
Note first that every $U\in\Ro(X)$ satisfies
\begin{align*}
(\Int_X U)\cap Y&=\Int_Y (U\cap Y),\\
(\cl_X U)\cap Y&=\cl_Y (U\cap Y),\\
(\ext_X U)\cap Y&=\ext_Y (U\cap Y),\\
\Int_X \cl_X (U\cap Y)&=U.
\end{align*}
It is then easy to see then that $\psi$ as defined in the lemma is an isomorphism. 

Suppose $F\in\TF(X)$ satisfies \[y:=\rho_X(F)\ne \infty.\]
Then we have the equalities
\[
y\in\bigcap_{U\in F} \cl_X(U)\cap Y=\bigcap_{U\in F} \cl_Y(U\cap Y)
=\bigcap_{V\in \psi(F)} \cl_Y V,\]
which imply that $ \psi(F)\in \TF(Y)$.

Conversely if $F\in \TF(Y)$ then 
\[\infty\ne \rho_Y(F)\in \bigcap_{U\in F}\cl_Y(U)
\sse
\bigcap_{U\in F}\cl_X(U)
=\bigcap_{V\in\psi^{-1}(F)} \cl_X(V)
.\]
This implies that $ \psi^{-1}(F)\in\TF(X)\setminus\rho_X^{-1}(\infty)$.

To see part (2), let $g\in G\left\{F^\perp\right\}$. We have some $W\in F$
such that $\supp g\sse W^\perp$. 
Since $\infty\in \cl_X W$, we obtain $g\in \Fix\infty$.

Lastly, to verify the local density let us pick an arbitrary point $\infty\ne y\in Y$ and its open neighborhood $V\sse Y$. 
By shrinking $V$ if necessary (using the fact that $Y$ is locally compact and Hausdorff), we may assume that $V$ is regular
and that $\infty\in \ext_X V$.
 In particular, we see that $G[V].y\sse V$
 and that $G[V]$ fixes $\infty$.
By the local density of $G$, some element $V_0\in\Ro(X)$ satisfies
\[
\varnothing\ne V_0\sse \cl_X \left(G[V].y\right)\sse \cl_X V=\comp_X \ext_X V.\]
This implies that $\infty\not\in V_0$ and that
\[
V_0\sse Y\cap \cl_X \left(G[V].y\right)=\cl_Y \left(G[V].y\right).\]
We see that the action of $G\left\{F^\perp\right\}$ on $Y$ is  locally dense.
\ep

\bp[Proof of Theorem~\ref{thm:bool-homeo}, part (\ref{p:rhoi-rho})]
Set $\infty:=\rho_X(F)\in X$ and $Y:=X\setminus\{\infty\}$. 
By part (\ref{p:tf}) of Theorem~\ref{thm:bool-homeo} and by Lemma~\ref{lem:pointed}, we have a sequence of bijections
\[
\begin{tikzcd}
\TF(X)\setminus\rho_X^{-1}(\infty)
\arrow{r}{\psi} &
\TF(Y)\arrow{d}{=} \\
\CS(\Ro(X),G[F^\perp]) &
 \CS(\Ro(Y),G[F^\perp])\arrow{l}[swap]{\psi^{-1}}  
\end{tikzcd}\]
This completes the proof of the theorem.
\ep
The deduction of Rubin's Theorem from Theorem~\ref{thm:bool-homeo} will be mostly formal, just like  that of Theorem~\ref{thm:gp-bool2} from Theorem~\ref{thm:rubin-supp}.
\bp[Proof of Theorem~\ref{thm:rubin}]
We have constructed a Boolean isomorphism
\[
\Psi\co \Ro(X_1)\longrightarrow\Ro(X_2)\]
in Theorem~\ref{thm:gp-bool2} that is equivariant with respect to the group actions $G_i$.
We have the following equalities 
for  $F\in \TF(X_1)$:
\begin{align*}
&\Psi(\TF(X_1))
=\TF(X_2),\\
&\Psi(\rho_{X_1}^{-1}\circ\rho_{X_1}(F))
=\rho_{X_2}^{-1}\circ\rho_{X_2}(\Psi(F)),\\
&\Psi(\Nbr_{X_1}(\rho_{X_1}(F))\cap\Ro(X_1))
=\Nbr_{X_2}(\rho_{X_2}(\Psi(F))\cap\Ro(X_2).
\end{align*}
Indeed, the sides of each equality above can be expressed purely in terms of the Boolean actions $G_i$ on the Boolean algebras 
$\Ro(X_i)$ by Theorem~\ref{thm:bool-homeo}. Thus, they are preserved under the equivariant Boolean isomorphism $\Psi$.

We now define a map $\phi\co X_1\longrightarrow X_2$ by 
\[
\phi(\rho_{X_1}(F))=\rho_{X_2}(\Psi(F)),\]
for $F\in\TF(X_1)$. This is a well-defined bijection by the observation above.

For a regular open set $U\sse X_1$, we have
\begin{align*}
\phi(U)&=\{ \phi(y)\mid y\in U\}
=\{\phi\circ\rho_{X_1}(F)\mid  F\in\TF(X_1)\text{ and }\rho_{X_1}(F)\in U\}\\
&=\{\rho_{X_2}\circ\Psi(F)\mid  F\in\TF(X_1)\text{ and }U\in\bigcap\rho_{X_1}^{-1}\circ\rho_{X_1}(F)\}\\
&=\{\rho_{X_2}(F')\mid  
F'\in\TF(X_2)\text{ and }\rho_{X_2}(F')\in\Psi(U)\}=\Psi(U).
\end{align*}
This implies that $\phi$ is an open map, and by symmetry, a homeomorphism.

\ep
\begin{rem}
\be[(1)]
\item
In Rubin's Theorem, one cannot drop the perfectness hypothesis; for instance, one may let $G_1=G_2=\bZ/2\bZ$ nontrivially act on two finite discrete spaces of different cardinalities.
\item
Similarly, one cannot weaken the local density hypothesis to mere local movement.
To see an example, consider the action of $\Homeo_+[0,1]$ on the compact interval $X_1=[0,1]$ and also its one-point compactification $X_2=S^1$.
Both actions are locally moving, but the spaces are not homeomorphic.
\item One can weaken the hypothesis of local compactness 
to \emph{regional compactness}\index{regional compactness}, 
which means that every nonempty open set contains some nonempty compact neighborhood. We direct the reader to
\cite{Rubin1996} for details. 
\ee
\end{rem}

\subsection{Applications to manifold homeomorphism groups}
We can now
combine Rubin's Theorem (\ref{thm:rubin})
with Taken's Theorem (\ref{thm:takens})
and obtain a simple proof of (generalized) Filipkiewicz's Theorem (Theorem~\ref{thm:filip-gen}). 

\bp[Proof of Filipkiewicz's Theorem]
Note that for a smooth connected boundaryless manifold $X$
and for $p\in\bZ_{>0}\cup\{0,\infty\}$, the group 
$\Diff_c^p(X)_0$ acts on $X$ locally densely.
This is an immediate consequence of Lemma~\ref{lem:affine}.

Therefore, the groups $G$ and $H$ given in the hypothesis act locally densely on $M$ and $N$ respectively.
Rubin's Theorem implies that there exists a homeomorphism $w\co M\longrightarrow N$ intertwining these actions. 
We conclude from Takens' Theorem (Theorem~\ref{thm:takens}) that 
$p=q$ and that $w$ is a $C^p$ diffeomorphism.
\ep

Ben Ami and Rubin gave a different reconstruction theorem, as we now describe.
We will omit the proof, which relies on Theorem~\ref{thm:gp-bool2}
and which resembles that of Rubin's Theorem. 
We remark that in the paper of Ben Ami and Rubin,
the spaces are only assumed to be regular; 
however, the statement below is equivalent to the original, since the perfectness and local 
compactness are consequences of the group theoretic hypotheses,
even when only regularities are assumed.

\begin{thm}[\cite{BAR2010}]\label{thm:bar}
Let $X_1$ and $X_2$ be perfect, locally compact, Hausdorff topological space.
Suppose we have groups
\[H_i\le G_i\le \Homeo(X_i)\]
such that $H_i$ is fragmented and has no global fixed points for $i=1,2$. 
Assume further that $\cl(G_i.x)$ has nonempty interior for $i=1,2$ and for all $x\in X_i$.
If there exists a group isomorphism
\[
\Phi\co G_1\to G_2,\]
then there exists a homeomorphism
\[
\phi\co X_1\to X_2\]
such that 
\[
\phi\circ g=\Phi(g)\circ\phi\]
for all $g\in G_1$.\end{thm}

\bp[Deduction of Filipkiewicz's Theorem from Theorem~\ref{thm:bar}, for $p,q\ge1$]
Setting \[ H_1:=\Diff_c^p(M)_0, \quad H_2:=\Diff_c^q(N)_0,\]
we note that $H_i$ is fragmented (Lemma~\ref{lem:frag}). Furthermore, arbitrary orbits of these groups have nonempty interior. 
It follows from Theorem~\ref{thm:bar} that there exists a homeomorphism  $X_1\to X_2$ that is equivariant with the given group isomorphism 
$G_1\to G_2$. The rest of the proof is identical as above.
\ep

Another easy application of Rubin's Theorem is that the PL homeomorphism group determines the ambient PL manifold 
up to homeomorphisms~\cite{Rubin1989}. One also sees that the minimal action of the Thompson's group $F$
(more generally, a minimal chain group) on an interval is unique up to homeomorphism~\cite{KKL2019ASENS}; cf.~Section~\ref{sec:chain}.

Let us also note another consequence.
\begin{cor}\label{cor:normal}
Let $p\in\bZ_{>0}\cup\{\infty\}$.
\be[(1)]
\item
If two distinct normal subgroups of $\Diff^p(M)$ contain $\Diff^p_c(M)_0$,
then they are not isomorphic as groups.
\item
If $p\ne\dim M+1$,
then two distinct normal subgroups of $\Diff^p(M)$ are never isomorphic as groups.
\ee
\end{cor}
\bp
For part (1), assume for contradiction that for two such normal groups $N_1$ and $N_2$ 
there exists an isomorphism $\Phi\co N_1\longrightarrow N_2$.
By Theorem~\ref{thm:filip-gen}, this isomorphism is the restriction of an inner automorphism of $\Diff^p(M)$.
By normality, we obtain that $N_1=N_2$.

In part (2) note first that $\Diff^p_c(M)_0$ is simple (Theorem~\ref{thm:mather-thurston}); see~\cite{Mather1,Mather2,Thurston1974BAMS}. By Corollary~\ref{cor:minimal}, every nontrivial normal subgroup of $\Diff^p(M)$ contains 
the group
\[
\left[\Diff^p_c(M)_0,\Diff^p_c(M)_0\right]
=\Diff^p_c(M)_0.\]
Part 1 now implies that two distinct nontrivial normal subgroups can never be isomorphic.

\ep

In particular, the four groups in  Corollary~\ref{cor:minimal} are pairwise non-isomorphic.

\subsection{Locally moving groups obey no law}\label{ss:lawless}
The local movement hypothesis puts  strong restrictions on algebraic structure of the underlying group.
Let $G$ be a group.
Each word $w=w(a_1,\ldots,a_k)$ in the rank--$k$ free group $F(x_1,\ldots,x_k)$
can be regarded as a map
\[
w\co G^k\to G\]
sending $(g_1,\ldots,g_k)$ to $w(g_1,\ldots,g_k)$.
We say $G$ obeys a nontrivial \emph{law}\index{law} if $w(G^k)=\{1\}$ for some integer $k\ge1$ and 
nontrivial reduced word $w\in  F(a_1,\ldots,a_k)$.
For example, abelian
groups obey the law $a_1^{-1}a_2^{-1}a_1a_2$.

Brin and Squier proved that Thompson's group $F$ does not obey a law~\cite{BS1985}; see also Section~\ref{ss:sub-chain}.
Ab\'ert gave a short proof of this fact by establishing the result below.

\begin{thm}[{Ab\'ert, \cite{Abert2005}}]\label{thm:abert}
Let $G$ be a permutation group of a set $X$.
Assume that for all finite subsets $Y\sse X$ and for all $z\in X\setminus Y$,
there exists $g\in G$ such that  $g\restriction_Y=\Id$  and such that  $g(z)\ne z$.
Then $G$ does not satisfy a law.
\end{thm} 

Generalizing the proof of Lemma~\ref{lem:LM-basic} (\ref{p:move-fn}) given by Rubin~\cite{Rubin1996}, one can see that the local movement condition is closely related; this was observed in Nekrashevych's note~\cite{Nekrashevych-note}.

\begin{thm}\label{thm:lawless-top}
A locally moving group of homeomorphisms of a Hausdorff topological space does not obey a law.
\end{thm}

As an example, the group of $C^p$ diffeomorphisms of a smooth connected boundaryless manifold obeys no laws. 
One can also see that the commutator subgroup of a locally moving group acting on a Hausdorff topological space is locally moving~\cite{Nekrashevych-note}. To see this, suppose  $G\le\Homeo(X)$ is locally moving and let $\varnothing \ne U\in \Ro(X)$. Since $H:=G[U]$ is locally moving, 
this latter group is not abelian by Theorem~\ref{thm:lawless-top}. 
In particular, the group $[H,H]$ is nontrivial. It follows that \[[G,G][U]=G'[U]\] is nontrivial.
The same argument shows that an arbitrary term in the lower central or derived series of $G$ is locally moving, and in particular, nontrivial. 
One can also formulate a similar result on lawlessness in the context of a Boolean algebra.

We will deduce both of the preceding theorems from the following unifying fact.

\begin{thm}\label{thm:lawless}
Let $G$ be a permutation group of a set $X$.
Assume that for each finite subset $S\sse G$ containing $1$ and for each $g\in G$,
there exist elements $h\in G$ and $y\in X$ such that the following conditions hold.
\be[(i)]
\item $\#S(y)=\max_{x\in X}\#S(x)$;
\item $S(y)\setminus\{y\}\sse \Fix h$;
\item $h(y)\not\in gS(y)$.
\ee
Then $G$ does not obey a law.
\end{thm}
\bp[Proof of Ab\'ert's Theorem from Theorem~\ref{thm:lawless}]
We assume that $(G,X)$ satisfies the hypothesis of Ab\'ert's Theorem.
In order to verify conditions (i) through (iii) of Theorem~\ref{thm:lawless}, 
we let $S\sse G$ be a finite subset containing $1$, and let $g\in G$.
Pick $y\in X$ realizing the condition (i). 
We set
\[
Z:= S(y)\setminus\{y\}.\]
We claim that the orbit $C:=G[X\setminus Z].y$ is infinite.
If $C$ is finite, then consider the group 
\[H:= \{r\in G[X]\mid r(s)=s\text{ for all }s\in Z\cup C\setminus\{y\}\}.\]
Applying the hypothesis of Ab\'ert's Theorem to the finite set 
\[Z\cup C\setminus\{y\}\sse X\setminus\{y\},\] we see that $H$ cannot fix $y$.
On the other hand, the group $H\le G[X\setminus Z]$ stabilizes $C$ setwise. 
As $H$ fixes $C\setminus\{y\}$ pointwise,
we have a contradiction.

By the above claim, we can now pick an element $h\in G[X\setminus Z]$ such that
\[ h(y) \not\in gS(y).\]
This completes the verification of conditions (i) through (iii).\ep

\bp[Proof of Theorem~\ref{thm:lawless-top} from Theorem~\ref{thm:lawless}]
Let us again verify the conditions of Theorem~\ref{thm:lawless}.
We first pick $x$ such that $\#S(x)=\max_{z\in X}\#S(z)$. 
By considering a subset $S_0\sse S$ if necessary, we may assume that
$s(x)\ne s'(x)$ for all distinct $s,s'\in S$;
we can further assume for  some open neighborhood $U$ of $x$  that 
\[sU\cap s'U=\varnothing,\]
as we have seen in the proof of Lemma~\ref{lem:LM-basic}.

If $x\not\in gS(x)$, then the conditions are trivially satisfied by setting $y:=x$ and $h:=\Id$.
We assume $x=gq(x)$ for some $q\in S$; note that such a $q$ is unique.
We may shrink $U$ to another open set $V$ containing $x$
such that $V$ and $gs(V)$ are disjoint for all $s\in S\setminus\{q\}$,
and such that $gq(V)\sse U$.

We set $h\in G[V]$ and pick $y\in \supp h$. 
We first note that both $(h,y)$ and $(\Id,y)$ satisfy the conditions (i) and (ii).
Indeed, from the fact that $s(V)$ and $s'(V)$ are disjoint for all distinct $s,s'\in S$, 
the condition (i) is satisfied.
The condition (ii) is trivial for $(\Id,y)$. 
For each $s\in S\setminus\{1\}$ we have that $U\cap s(U)=\varnothing$, and that
$s(y)\not\in \supp h\sse V$. This implies that $hs(y)=s(y)$ and verifies (ii) for $(h,y)$.

We claim that either $y\not\in gS(y)$ or $h(y)\not\in gS(y)$.
Indeed, assume first that $y\in gS(y)$. Since $g(S\setminus\{q\})V\cap V=\varnothing$,
we have that 
\[
h(y)\ne y=gq(y).\]
We also have that $h(y)\not\in g(S\setminus\{q\})y$, implying $h(y)\not\in gS(y)$.
For the case $h(y)\in gS(y)$, we similarly obtain
\[
y\ne h(y)=gq(y),\]
and that $y\not\in gS(y)$. The claim is proved.
Summarizing, either $(\Id,y)$ or $(h,y)$ satisfies the condition (iii).\ep

Let us now verify the main result of this subsection.

\bp[Proof of Theorem~\ref{thm:lawless}]
To show that $G$ obeys no law, it suffices to show that it obeys no law in a free group $F_2$ on two generators 
$\{a,b\}$, since every finitely generated free
group embeds in a free group on two generators. 

Let $\ell$ be a positive integer, and let $w(a,b)=b_\ell\cdots b_1 \in F_2$ be a nontrivial reduced word.
In particular, we have $b_i\in\{a,b\}^{\pm1}$. 
We set
\[
w_i(a,b):=b_i\cdots b_1\]
and $w_0:=1$.

We use the induction on the length $\ell$  
to find elements $u,v\in G$ and a point $y\in X$ such that
the points in the collection
\[
\{w_i(u,v)(y)\}_{0\le i\le \ell}\]
are all distinct. 
If $\ell=1$, we have $w(a,b)$ is a single letter. 
Note that $G$ is nontrivial by condition (iii).
We can define $w(u,v)$ to be an arbitrary nontrivial element of $G$
and $y$ to be a point nontrivially moved by that element.

Assume now that the conclusion holds for $\ell-1$.
We have some $u,v\in G$ 
and a point $x_0\in X$ such that
the points
\[
x_i:= w_i(u,v)(x_0)\]
are all distinct for $i\in\{0,1,\ldots,\ell-1\}$.
We assume $x_\ell:=w(u,v)(x_0)=x_j$ for some $j\le \ell-1$, for otherwise the proof is done.
Without loss of generality, we may assume $b_\ell=a$.
Setting $x:=x_{\ell-1}$ we have $x_\ell=u(x)$.

Let us set
\[
S:=\{w_i(u,v)w_{\ell-1}(u,v)^{-1}\mid 0\le i\le \ell-1\}.\]
Applying the hypothesis to $(S,g=u^{-1})$ we obtain $h\in G$ and $y\in X$. In particular,
\[
\#S(y)=\ell=\#S.\]
We set
\[
y_i:=w_i(u,v)w_{\ell-1}(u,v)^{-1}(y)\]
for $0\le i\le \ell$. The points $y_0,\ldots,y_{\ell-1}$ are all distinct.

By the condition (ii), we have $h(y)_i=y_i$ for $0\le i\le \ell-2$.
Using the condition that $b_{\ell-1}\ne a^{-1}$, we easily see for each $0\le i\le \ell-1$ that
\[w_i(uh,v)(y_0)=w_i(u,v)(y_0)=y_i.\]
Moreover, we have that
\[
w_\ell(uh,v)(y_0)=uh(y)\not\in\{y_0,\ldots,y_{\ell-1}\}.\]
We thus obtain that \[\#\{w_i(uh,v)y_0\mid 0\le i\le \ell\}=\ell+1,\] as required.
\ep

Small amount of additional work will imply the following result on the abundance of free subgroups.
We note that a similar idea can be found in Ghys' article~\cite{Ghys1999}, where he proves that a generic two-generator subgroup of $\Homeo_+(S^1)$ is nonabelian free.
Recall that if a space is complete metrizable or locally compact Hausdorff
then it is \emph{Baire}\index{Baire space}, i.e. satisfies the conclusion of the Baire Category Theorem.
\begin{cor}\label{cor:free-abundant}
Let $G$ be a topological group, let $X$ be a space,
and let
\[ G\longrightarrow \Homeo(X)\]
be a continuous injective homomorphism.
Assume that for each nonempty open set $U\sse X$
and for each identity neighborhood $\VV\sse G$
we have that
\[
G[U]\cap\VV\ne\{1\}.\]
If $G$ and $X$ are Baire spaces,
then a generic pair of elements in $G$ generate a nonabelian free group.
\end{cor}
Here, a \emph{generic}\index{generic pair} pair means it belongs to a comeager subset of $G\times G$.
For instance, generic pair of diffeomorphisms in $\Diff^\infty(M)$ will generate $F_2$ for a smooth connected manifold $M$,
since the diffeomorphism group
is Baire (because it is a Fr\'echet space~\cite{conway-funct}).
Note that the group $G$ above is locally moving.
\bp[Proof of Corollary~\ref{cor:free-abundant}]
Since $G$ is Baire, so is the space \[Y=G\times G\times X.\] 
For each nontrivial reduced word $w(a,b)\in F_2:=F(a,b)$ we 
consider the open subset $U_w\sse Y$ given by triples
\[U_w=\{(u,v,x)\mid w(u,v)(x)\ne x\}.\] 
If we can show that $U_w$ is dense for all $1\neq w\in F_2$, 
then the Baire Category Theorem implies that
the $G_\delta$ set
\[\Delta=\bigcap_{1\neq w\in F_2} U_w\] is comeagre. 
If $(u,v,x)\in \Delta$ then for all nontrivial $w\in F_2$, we have that $w(u,v)(x)\neq x$,
and so the homeomorphisms $u$ and $v$ generate a free subgroup of $G$. This implies the conclusion.

It now remains for us to show that $U_w$ is dense, or equivalently,  $Z_w:=Y\setminus U_w$ has empty interior.
The idea is already contained in the proofs of Theorems~\ref{thm:lawless-top} and~\ref{thm:lawless}, so we will reuse notations and observations from those proofs.
We pick $(u,v,x_0)\in Z_w$, and consider arbitrary open neighborhoods  $\Id\in \VV\sse G$ and $x_0\in U\sse X$.
We wish to prove that $(u\VV,v\VV,U)$ contains an element of $U_w$, implying that $(u,v,x_0)$ is not an interior point of $Z_w$.
The case that $w$ is a single letter is obvious, as we can find $y\in U$ and $h\in G[U]\cap \VV$ such that $h(y)\ne y$,
and hence, either $w(y)\ne y$ or $wh(y)\ne w(y)=y$.

Inductively,  let $w$ have length $\ell>1$. 
We may assume that \[(u,v,x_0)\in Z_w\cap\left(\bigcap_g U_g\right),\] where $g$ ranges over all the reduced words shorter than $w$.
We define \[w_0=1,w_1,\ldots,w_\ell=w\] to be exactly the same as in Theorem~\ref{thm:lawless}. 
We also assume $w=a w_{\ell-1}$, as before.

The points $x_i:=w_i(u,v)(x_0)$ are all distinct for $0\le i<\ell$. We have that $x_\ell=x_j$ for some $j<\ell$ since $(u,v,x_0)\in Z_w$.
As noted in the proof of Theorem~\ref{thm:lawless-top}, we can find an open neighborhood $V_0$ of $x_0$ such that  $\{V_i:=w_i(u,v)(V_0)\}_{0\le i\le \ell-1}$ are all disjoint.
We can pick $h\in G[V_{\ell-1}]\cap \VV$ such that $h(y)\ne y$  for some $y\in V_{\ell-1}$. 
As was seen in the proof of Theorem~\ref{thm:lawless}, we deduce that 
either $(u,v,w_{\ell-1}(u,v)^{-1}(y))$ or
$(uh,v,w_{\ell-1}(u,v)^{-1}(y))$ does not belong to $Z_w$, completing the proof.\ep


\section{The original proof of Whittaker--Filipkiewicz}\label{sec:filip}
In this section we give another proof of Theorem~\ref{thm:filip-gen}, based on the original argument of 
Whittaker~\cite{Whittaker1963} and its generalization by Filipkiewicz~\cite{Filip82}. This method is much more 
involved than Rubin's Theorem and also significantly more limited in scope, but we present the original proof for the showcase of the 
idea of detecting stabilizer groups by group theoretic properties. Moreover, it has a slight advantage of dealing with possibly 
non-locally compact spaces, where Rubin's Theorem does not apply.
We will restrict ourselves to the differentiable world, namely the cases $p,q\ge1$, so that we may 
employ fragmentation techniques.

Throughout, we will let $M$ and $N$ be smooth connected boundaryless manifolds.
As we saw in the introduction of this chapter, this theorem implies the result of R.~P.~Filipkiewicz (Theorem~\ref{thm:filip}),
which shows that an isomorphism between the diffeomorphism groups of $M$ and of $N$ is induced by a diffeomorphism
between these manifolds themselves.

The idea of the proof of Theorem~\ref{thm:filip-gen} is to first find a bijection between $M$ and $N$, conjugation by which realizes the given
isomorphism $\Phi$ between the diffeomorphism groups. The bijection itself arises from a bijection between stabilizers of
points in $M$ and $N$. Specifically, one considers $\Stab_G(x)$ for a point $x\in M$, and concludes that there
is a point $y\in N$ such that $\Phi(\Stab_G(x))=\Stab_H(y)$. The assignment $x\mapsto y$ turns out to be a bijection
realizing the isomorphism $\Phi$.
As in Lemma~\ref{lem:bij-homeo}, the resulting bijection will necessarily be
continuous. The remaining work is in showing that the bijection is in fact a $C^p$ diffeomorphism,
and this last bit utilizes Theorem~\ref{thm:BM-simple} (the Bochner--Montgomery Theorem). 

This process can be more conceptually conveyed using abstract properties of group actions. 
A \emph{circular order}\index{circular order} on a set $X$ roughly means that
 arbitrary three distinct points $x,y,z$ in $X$
 can be listed as
\[
x<y<z<x\] after a suitable permutation,  and moreover this ordering is transitive;
see Section~\ref{ss:circular} for a precise definition.
The reader may only consider the case $X=S^1$ for the purpose of this section.

\bd
A \emph{dense order (without bounds)}\index{dense order} on a set $X$ is 
a total or circular order $\le$ such that
 every interval in $X$ contains a point in $X$
and such that $X$ does not have a maximum or a minimum.
\ed
We will often drop the phrase ``without bounds'' when discussing orders. 
By a classical result of Cantor, a countable set with a dense total order without bounds is order isomorphic to $\bQ$.
The same condition for a circular order yields $\bQ/\bZ$.

For a set $X$, we let
$\Sym(X)$ denote the group of all permutations of $X$, 
This group naturally acts on the configuration space of $n$ distinct points
\[
\operatorname{Conf}_n(X):=\{(x_1,\ldots,x_n)\mid x_i\in X\text{ and }x_i\ne x_j\}.\]

If $X$ is given with a linear order $\le$, then we define 
\[
\operatorname{Conf}^+_n(X):=\{(x_1,\ldots,x_n)\mid x_i\in X\text{ and }
x_1<x_2<\ldots<x_n\}.\]
If $X$ is equipped with a circular order $\le$, then we let
\[
\operatorname{Conf}^+_n(X):=\{(x_1,\ldots,x_n)\mid x_i\in X\text{ and }
x_1<x_2<\ldots<x_n<x_1\}.\]
When $X$ is equipped with either of these orders, then a linear or circular order preserving permutation of $X$
naturally acts on $\operatorname{Conf}^+_n(X)$.

\bd\label{defn:n-trans}
Let $X$ be a set.
\be[(1)]
\item A group $G\le\Sym(X)$ is \emph{$n$--transitive}\index{$n$--transitive}, or \emph{$T(n)$}, if
it acts transitively on $\operatorname{Conf}_n(X)$.
\item Suppose $X$ is equipped with a total or circular order.
A group $G\le\Sym(X)$ is \emph{positively $n$--transitive}\index{positively $n$--transitive},
or \emph{$T^+(n)$}, if it acts transitively on $\operatorname{Conf}^+_n(X)$.
\ee\ed
For brevity, we adopt the convention that when we discuss homeomorphisms of one--manifolds,
``$n$--transitivity'' actually means positive $n$--transitivity unless stated to the contrary.

\begin{rem}
\be[(1)]
\item A $2$--transitive action on a circularly ordered set is positively 2--transitive by definition.
\item A positively 3--transitive action on a set with a dense, total or circular order is automatically order preserving. 
\ee
\end{rem}

Roughly based on Banyaga's exposition~\cite{Banyaga1997}, we introduce the following axiomatic definitions for topological actions.
\bd\label{defn:axioms}
Let $X$ be a topological space, and let $G\le\Homeo(X)$.
\be[(1)]
\item The group $G$ is \emph{locally dilative (LD)}\index{locally dilative}
if every point in $X$ has an open local basis consisting of open sets $U$ such that some element  in $G_U$ has a connected closed support.
\item
The group $G$ is \emph{weakly fragmented (WF)}\index{weak fragmentation} 
if 
 for every open cover $\VV$
 of $X$ the group generated by the collection
\[
\big\{
\left[
G_V,G_V
\right]\;
\big\vert\;
V\in\VV
\big\}
\]
contains a nontrivial normal subgroup of $G$.
\ee
\ed

\begin{rem}\label{rem:banyaga}
\be[(1)]
\item
The property (LD) is much weaker than the \emph{property (B)}\index{property (B)} defined by Banyaga~\cite{Banyaga1997},
which states that some element $g\in G$ satisfies $\suppo g = U\setminus\{x_0\}$ for each element $U$ of the local basis at $x_0$.
\item
The property (WF) is trivially implied by the \emph{property (L)}\index{property (L)} of Banyaga, which states that 
\[
\form*{ \left[ (G_V)_0,(G_V)_0\right]\middle \vert V\in\VV}=[G_{c0},G_{c0}],
\]
assuming $G$ is a $C^k$ diffeomorphism group of a manifold
and $\VV$ is a subcover of $C^k$ open balls.
\item Banyaga defined $G\le\Homeo(M)$ to be \emph{path--transitive}\index{path transitive} if every neighborhood of a given path
$c\co[0,1]\longrightarrow M$ admits some $g\in G$ with $\suppc g\sse U$ and such that $g(c(0))=c(1)$. 
The path--transitivity of $G$ easily implies $(LT)$ and $T(n)$ or $T^+(n)$ (depending on the dimension of $M$).
\item\label{p:ptrans} Conversely, if $X$ is path--connected then every locally transitive action is path--transitive. This is obvious from a typical Lebesgue number argument.
\ee
\end{rem}

We can now state a generalization of Filipkiewicz Theorem.
Banyaga~\cite{Banyaga1997} proved the same conclusion with the stronger assumption that the properties (B), (L), and path--transitivity hold for the groups $G_1$ and $G_2$. 

\begin{thm}\label{thm:isom-homeo}
Let $X_1$ and $X_2$ be  connected Hausdorff topological spaces.
For each $i\in\{1,2\}$, we assume one of the following:
\be[(A)]
\item
A group $G_i$ acts faithfully and 3--transitively on $X$.
\item
A group $G_i$ acts faithfully and positively 3--transitively on $X_i$
with respect to some dense total or circular order on $X_i$.
\ee
Assume further that each $G_i$ is locally transitive, locally dilative, and weakly fragmented.
If there exists a group isomorphism
\[
\Phi\co G_1\longrightarrow G_2,\]
then there exists a homeomorphism
\[w\co X_1\longrightarrow X_2\]
such that each $g\in G_1$ satisfies
\[
\Phi(g)=wgw^{-1}.\]\end{thm}
We emphasize that it is not required a priori that the spaces $X_1$ and $X_2$ are simultaneously ordered,
even when one of them is assumed to be ordered.
We also do not need to assume that $X_i$ is locally compact, in comparison with Rubin's Theorem (Theorem~\ref{thm:rubin}).

\begin{rem}The 3--transitivity assumption is redundant when $X_1$ and $X_2$ are connected manifolds of dimension at least two. This is because for all distinct $p,q\in \Conf_3(X_i)$
one can find disjoint paths $P_1,P_2,P_3\sse X_i$ such that the path $(P_1,P_2,P_3)\sse \Conf_3(X_i)$ joins $p$ and $q$. One can then apply part~(\ref{p:ptrans}) of Remark~\ref{rem:banyaga} to see that the local transitivity hypothesis readily implies the 3--transitivity. The positive 3--transitivity assumption can also be dropped for a similar reason when $X_i$ is a connected one--manifold.\end{rem}

\begin{prop}\label{prop:wf}
Let $H$ be a group, and let $p\in\bZ_{>0}\cup\{\infty\}$.
\be[(1)]
\item If we have
\[\Diff^\infty_c(M)_0\le H\le\Homeo(M),\]
then $H$ is locally transitive, locally dilative and 3--transitive (positively, if $M$ is a one--manifold).
\item  If we have
\[\Diff^p_c(M)_0\le H\le\Diff^p(M),\]
then $H$ is weakly fragmented.
\ee
\end{prop}

The previous proposition implies that all of the four groups
\[
\Diff^p(M),\Diff^p_c(M),\Diff^p_{c}(M)_0,\Diff^p(M)_0\]
are locally transitive, locally dilative, 3--transitive and weakly fragmented.
 
The Generalized Filipkiewicz's theorem as stated in Theorem~\ref{thm:filip-gen} 
is an immediate consequence of Theorem~\ref{thm:isom-homeo}
and of Takens' Theorem (Theorem~\ref{thm:takens}).

\begin{rem}
Similar conclusions to Theorem~\ref{thm:isom-homeo} hold for the groups of contactomorphisms, symplectomorhisms
(when the manifolds are compact), volume-form preserving smooth diffeomoprhisms, and ``good-measure'' preserving homeomorphisms,
by establishing the properties described in the theorem above. The most difficult parts in the processes are usually the property (WF).
We will not delve into these groups as they require many definitions beyond the scope of this book. See~\cite{Banyaga1997} for details.
\end{rem}

\subsection{Transitivity, dilativity and weak fragmentation}\label{ss:ltld}
In this subsection, we establish the properties (LT), (LD), and $T(n)$ in  Proposition~\ref{prop:wf}. 
As these properties transfer from a subgroup $H\le G$ to a bigger group $G$, 
it suffices to verify them for the smallest group considered in the proposition, that is
\[
\Diff_c^\infty(M)_0.\]

For a subset $U\sse M$, we abbreviate that 
\[ \Diff^k_0(U):=\Diff^k_c(U)_0=(\Diff^k_c(M)_U)_0.\]

Returning to the proof of Proposition~\ref{prop:wf}, 
we observe from the above lemma that for each given $C^\infty$ open ball $U\sse M$
there exists $g\in\Diff_c^\infty(M)_0$ whose open support is $U$ itself
and which moves between two points given a priori.
It follows that the group
\[
\Diff_c^\infty(M)_0\]
has the properties (LT) and (LD). The property $T(n)$ of this group is given below.
\begin{lem}[\cite{Banyaga1997}]\label{lem:t3}
For each path $c\co I\longrightarrow M$
and for each open neighborhood $U$ of $c(I)$ there exists 
$g\in \Diff_0^\infty(U)$
such that $g(c(0))=c(1)$.
In particular, $\Diff_c^\infty(M)_0$ is $n$--transitive for all $n$.
\end{lem}
\bp
If $c(I)$ is contained in a $C^\infty$ open ball $U_0$ in $M$, 
and if $U_0\sse U$, then  Lemma~\ref{lem:affine} obviously implies 
the existence of $g\in\Diff_0^\infty(U_0)$ such that $g(c(0))=c(1)$.
In general, we may divide $c(I)$ into finitely small pieces $c_1,\ldots,c_m$, each piece of which is contained in a $C^\infty$ open ball $U_i$ inside $U$.
Composing maps $g_i$ sending the initial point of $\gamma_i$ to the terminal point as above, we have a map $g$ satisfying the first conclusion.

The second conclusion $T(n)$ is an easy inductive consequence of the first, along with the hypothesis that $M$ is path--connected.
\ep

The weak fragmentation property in Proposition~\ref{prop:wf}
is a consequence of a stronger result, that is Proposition~\ref{prop:wf2}.
This gives yet another perspective on fragmentation of homeomorphisms, which will later be crucial in the proof of Filipkiewicz's Theorem.

\begin{prop}[{cf. \cite[Theorem 2.2]{Filip82}}]\label{prop:wf2}
If $p\in\bZ_{>0}\cup\{\infty\}$ and if $\UU$ is an open cover of $M$, 
then the collection of groups
\[
\left\{
\left[
\Diff^p(U)_0,\Diff^p(U)_0
\right]\;
\middle\vert\;
U\in\UU
\right\}\]
contains the group  
\[\left[\Diff^p_c(M)_0,\Diff^p_c(M)_0\right].\]
\end{prop}

The weak fragmentation property in Proposition~\ref{prop:wf} would be an easy consequence.
\bp[Proof of Part 2 of Proposition~\ref{prop:wf}, assuming Proposition~\ref{prop:wf2}]
Let $H$ be as in the proposition,
and let $\VV$ be an open cover of $M$.
Denote by $G_1$ the group generated by the collection 
\[
\left[
H_V,H_V
\right]\]
for $V\in\VV$.
It suffices for us to prove that $G_1$ contains 
\[K_0:=\left[\Diff^p_c(M)_0,\Diff^p_c(M)_0\right]\]
since $K_0$ is nontrivial and normal in $H$.

Since the collection  $\BB^p(M)$
of $C^p$ open balls generates the topology of $M$, 
the collection 
\[
\UU:=\{U\in \BB^p(M)\mid U\sse V\text{ for some }V\in \VV\}\]
 still covers $M$. 
The group $G_2$ generated by the collection
\[
\left[
H_U,H_U
\right]\]
for $U\in\UU$ is contained in $G_1$. 
 
Since each $U\in\UU$ is relatively compact
and since $\Diff^p_c(M)_0\le H$, 
the group $G_2$ contains the group $K_L$
generated by the collection of the groups
\[
\left[
\Diff^p(U)_0,\Diff^p(U)_0
\right]\]
for $U\in\UU$.
This completes the proof that $H$ is weakly fragmented, as $K_0\le K_L$ by Proposition~\ref{prop:wf2}.\ep

The proof of Proposition~\ref{prop:wf2} requires two ingredient. The first is the  Fragmentation Lemma (Lemma~\ref{lem:frag}), and the second is the ``shrinking lemma'' of Filipkiewicz.
This latter lemma asserts that the closed unit ball in a Euclidean space
can be contracted onto a smaller ball
by the multiplication of commutators of smooth diffeomorphisms 
supported in smaller open balls.
Recall our notation that $B^d(a;r)$ denote the radius--$r$ ball centered at $a\in\bR^d$.

\begin{lem}[Filipkiewicz's Shrinking Lemma]\label{lem:filip-comm}
If $\UU$ is an open cover of the closed unit ball $B$ in $\bR^d$,
then for each $a\in (0,1]$
there exist elements \[\{f_i,g_i\}_{1\leq i\leq m}\sse\Diff_c^{\infty}(\bR^d)_0\] satisfying the following:
\be[(1)]
\item
For each $i$, there is a $U_i\in\UU$ such that $f_i,g_i\in\Diff^k_0(U_i)$.
\item
We have that \[[f_m,g_m]\cdots [f_1,g_1](B)\sse {B^d(0;a)}.\]
\end{enumerate}
\end{lem}

\bp
We let $A\sse[0,1]$ be the set of $a\in[0,1]$ satisfying the conclusion. This set is nonempty as it contains $1$.
Taking the infimum $a$ of this set, it suffices for us to show $a=0$. 

For a contradiction, assume $a>0$. 
Let us denote the sphere of radius $r>0$ as $S_r:=\partial B^d(0;r)$.
We claim that there exist elements $f_i, g_i$ satisfying the first condition of the lemma such that
\[\prod_i[f_i,g_i] \left(S_{a+\delta/2}\right)\sse S_{a-\delta/2}\] for some  $\delta>0$.
The desired contradiction would then follow from that $a+\delta/2\in A$.

We let $N_\delta$ denote the open $\delta$--neighborhood of $S_a$.
Pick a Lebesgue number $\epsilon\in(0,a/2)$ of the cover $\UU$ for $B$.
There exist points $x_1,\ldots,x_m\in S_a$
such that the  balls 
\[ V_i:=B^d(x_i;\epsilon/2)\] cover this sphere.
For some sufficiently small $\delta\in(0,\epsilon/4)$,
we also have that
\[
S_a\sse \overline{N_\delta}\sse \bigcup_i V_i\sse N_{\epsilon/2}.\]
We can pick an element
$f\in\Diff^\infty(\Int N_\delta)_0$ satisfying
\[ f\left(S_{a+\delta/2}\right)\sse S_{a-\delta/2}.\]

Each $U_i:=B(x_i;\epsilon)$ is contained in some element of $\UU$.
Moreover, 
we can find an element $g_i\in\Diff_0^\infty(U_i)$ such that 
\[g_i(V_i)\sse U_i\cap B(0;a-\epsilon/2)\sse\bR^d\setminus N_{\epsilon/2}.\]

Applying Lemma~\ref{lem:frag} to the manifold $\Int N_\delta$, we can find 
a fragmentation
\[
f=f_r\cdots f_1\]
such that each $f_i$ is supported in
 some set $W_i$ that coincides with one of $\{V_1,\ldots,V_m\}$.
Rearranging, we may assume $r=m$ and $W_i=V_i$ for each $i=1,\ldots,m$.
The map
\[
[f_i,g_i]=f_i\cdot (g_i f_i g_i^{-1})\]
is supported in $V_i\cup g_i(V_i)$.

If $x\in N_{\epsilon/2}$, then we have that $x\not\in g_i(V_i)$ and that
\[
[f_i,g_i](x)=f_i(x)\in N_{\epsilon/2}.\]
Inductively, we see that 
\[
[f_i,g_i]\cdots[f_1,g_1](x)=\prod_i f_i(x)\in N_{\epsilon/2}.\]
In particular, we have that
\[[f_m,g_m]\cdots[f_1,g_1](x)=f(x).\]
This proves the aforementioned claim.
\ep

Using Lemma~\ref{lem:filip-comm}, we can now establish Proposition~\ref{prop:wf2}.
This also gives yet another perspective on fragmentation of homeomorphisms (cf.~Lemma~\ref{lem:frag}).

\begin{proof}[Proposition~\ref{prop:wf2}]
As in our proof for the weak fragmentation part of Proposition~\ref{prop:wf},
we may assume that $\UU\sse\BB^p(M)$. 
We also define $K_L$ and $K_0$ as in that proof. 
Let \[K:=\form{ [\Diff_0^p(U),\Diff_0^p(U)]\mid U\in\BB^p(M)}.\] 
Since $\UU\sse\BB^p(M)$ we have that
\[K_L\le K\le K_0.\]
Moreover, $K$ is nontrivial and normal in $\Diff^p(M)$
since 
 $\BB^p(M)$ is $\Diff^p(M)$--invariant.
 By  Theorem~\ref{thm:epstein} and Corollary~\ref{cor:minimal}, we have that $K=K_0$.

It only remains to show that $K\le K_L$.
Let $d$ be the dimension of $M$
and let $U\in\BB^p(M)$.
There exists a $C^p$ embedding 
\[ h\co \bR^d\longrightarrow M\]
such that $U=h(B^d(0;1))$.
For each $x\in \overline U$, there is some neighborhood $U_x\in\UU$ of $x$ such that
\[[\Diff^p_0(U_x),\Diff^p_0(U_x)]\le K_L\] by the definition of $K_L$.
We may assume $U_x\sse h(B^d(0;2))$, without loss of generality.
Since $\overline{U}$ is compact, we can select finitely many of these neighborhoods, say $\{U_1,\ldots,U_m\}$, that cover $\overline{U}$.

We let $V_i=h^{-1}(U_i)$ for $i\in\{1,\ldots,m\}$. We have \[B^d(0;1)\sse\bigcup_{i=1}^m V_i.\] Relabeling if necessary, we assume that
$B^d(0;a)\sse  V_1\sse\bR^d$,
for some small $a>0$.
Lemma~\ref{lem:filip-comm} implies the existence of commutators
$\{[f_j,g_j]\}_{1\leq j\leq r}$, where $f_j,g_j\in \Diff^{\infty}_0(V_{i(j)})$ for some suitable labeling function $i(j)$, and such that
\[[f_r,g_r]\cdots[f_1,g_1](B^d(0;1))\sse B^d(0;a).\] 
Since all the diffeomorphisms $f_j$ and $g_j$ are supported in the interior of
$B^d(0;2)$, we can conjugate $f_j$ and $g_j$ by $h$ to define $C^p$ diffeomorphisms on the image of $h$, and extending to
all of $M$ by the identity. Writing $\tilde{f}_j$ and $\tilde{g}_j$ for the respective resulting diffeomorphisms of $M$, we obtain that
$\tilde{f}_j$ and $\tilde{g}_j $ are supported on $U_{i(j)}$, and \[[\tilde{f}_r,\tilde{g}_r]\cdots [\tilde{f}_1,\tilde{g}_1]\in K_L\]
sends $U$ into $U_1$. Thus,
$[\Diff_0^p(U),\Diff_0^p(U)]$ is conjugate by an element of $K_L$ into $[\Diff_0^p(U_1),\Diff_0^p(U_1)]\le K_L$.
This implies  the conclusion of the proposition.

\end{proof}

\begin{rem}
Strictly speaking, it suffices to show $K_L=K$ in the preceding proof for the purpose of establishing the weak fragmentation
properties in Proposition~\ref{prop:wf};
in other words, the fact that $K=K_0$ is not necessary for this purpose.
This is because the group $K$ in the proof is readily seen to be normal and nontrivial in $\Diff^p(M)$.
In that sense, Higman's Theorem and the Epstein--Ling Theorem are not necessary ingredients for the proof of Theorem~\ref{thm:filip-gen}, 
although those two theorems make the exposition more concrete by pinpointing exactly what the elements of $K$ are.
\end{rem}

\subsection{The pre-stabilizer subgroup}
The key to proving Theorem~\ref{thm:isom-homeo} is to show that $\Phi$ sends point stabilizers to point stabilizers,
from which one can build a
well--defined bijection from $X_1$ to $X_2$. The way one shows that $\Phi$ sends point stabilizers to point stabilizers is to characterize point
stabilizers as maximal proper subgroups in the corresponding homeomorphism groups.

Let us first consider group actions on abstract sets, without any reference to any topology.
The next definition records some useful double coset decompositions and their consequences for sufficiently transitive groups actions on sets.

\bd\label{defn:prestab}
Let $G$ be a group. A nontrivial proper subgroup $K$ of $G$ is called a \emph{pre-stabilizer subgroup}\index{pre-stabilizer subgroup}
if all of the following hold.
\be[(i)]
\item 
For all $f\in G\setminus K$, we have that 
\[
G\setminus K  =KfK\cup Kf^{-1}K.\]
\item For all $f\in G\setminus K$ and $g\in KfK$ satisfying $fg,gf\not\in G\setminus K$, we have that
\[ fg, gf\in KfK.\]
\item 
For all $f\in G\setminus K$ and $g_0,g_1\in KfK$ 
there exist some \[s_i,t_i\in f^{-1}Kg_i\cap K\]
for $i=0,1$ such that 
$
s_0s_1=t_1t_0$.
\ee
\ed

Being a pre-stablizer subgroup is a purely group theoretic property, which is obviously preserved under group isomorphisms. 
We first note a simple consequence.
\begin{lem}\label{lem:prestab-maximal}
A pre-stabilizer subgroup of a group $G$ is a maximal proper subgroup of $G$.
\end{lem}
\bp
Suppose we have
\[
K\lneq H\le G\]
for some pre-stabilizer subgroup $K$ of $G$.
Fix $h\in H\setminus K$. Then for all $g\in G\setminus K$ we have
\[
h\in KgK\cup Kg^{-1}K\sse HgH\cup Hg^{-1}H.\]
This implies that $g\in H$. We have shown that
\[
G\setminus K\sse H,\]
which implies that $H$ is not a proper subgroup.
\ep

The properties of Definition~\ref{defn:prestab} are possessed by point-stablizer groups of sufficiently transitive actions.

\begin{lem}[{cf.~\cite[Lemma 6]{Whittaker1963}}]\label{lem:prestab}
Assume one of the following.
\be[(A)]
\item
$X$ is an infinite set and 
a group $G\le \Sym(X)$ acts 3--transitively on $X$.
\item
$X$ is a set equipped with a dense order that is total or circular,
and a group $G\le \Sym(X)$ acts positively 3--transitively on $X$.
\ee
Then for each $x\in X$ the stabilizer group
\[
\Stab_G(x)\]
is a pre-stabilizer subgroup of $G$. 
\end{lem}
\bp
Set $K=\Stab_G(x)$. It is an easy consequence of  (positive) 3--transitivity that $1\ne K\ne G$. For instance,
if $X$ is circularly ordered,
then we pick some points $y,z,w$ such that
\[ x<y<z<w<x.\]
We have some $g\in G$ such that \[g(x)=x,\quad g(y)=z,\quad g(w)=w.\]
This implies that $g\in K\ne 1$. The other cases are similar.

Let us now consider the three cases separately. 

{\bf Case 1:  $G\le\Sym(X)$.}

We verify the three properties of Definition~\ref{defn:prestab}.
For condition (i) of the definition, we make a stronger claim that 
\[
G\setminus K = KfK.\]
Indeed, if $g\in G\setminus K$ then there exists an element $\sigma\in K$ such that $\sigma(g(x))=f(x)$, by 2--transitivity.
It follows that $\tau=f^{-1}\sigma g\in K$ and that $g=\sigma^{-1}f\tau\in KfK$.
This claim implies condition (ii) as well. 

For condition (iii), let us set
 \[u=f^{-1}(x),\quad v_0=g_0^{-1}(x),\quad v_1=g_1^{-1}(x).\] 
By the 2--transitivity of the action of $G$, the set below is nonempty:
\[f^{-1}Kg_i\cap K=
\{\sigma\in G\mid \sigma(x)=x\text{ and }\sigma(v_i)=u\}.\] 

Choosing a point $p\in X$, we will define maps based on the diagram below.
\[\begin{tikzcd}
p \arrow{r}{\sigma_0}\arrow[two heads]{d}[swap]{\sigma_1}  & v_1\arrow[two heads]{d}{\sigma_1}
& x\arrow[loop right]{}{\sigma_0,\sigma_1}  \\
v_0 \arrow{r}{\sigma_0} & u
\end{tikzcd}
\]
Namely, we use the 3--transitivity to find
$\sigma_0,\sigma_1\in K$ such that
\[\sigma_i(p)=v_{1-i},\quad \sigma_i(v_i)=u.\]
For such a $\sigma_i$ to exist, 
it is sufficient and necessary to require the following:
\begin{itemize}
\item $p=v_0$ if and only if $v_1=u$;
\item $p=v_1$ if and only if $v_0=u$.
\end{itemize}
After setting $\tau_0:=\sigma_0$ and
\[
\tau_1:=\sigma_0\sigma_1\sigma_0^{-1}\]
we obtain the condition (iii).

{\bf Case 2-1:  $X$ is equipped with a total order $\le$.}

The proof is very similar to Case 1. 
Given three distinct points $x,y,z$ of $X$, we say $y,z\in X$ are on the \emph{same side} of $x\in X$ 
if
either \[ x<y\text{ and }x<z\] or
either \[ y<x\text{ and }z<x.\]
To see condition (i) of Definition~\ref{defn:prestab}, let $g\in G\setminus K$. 
If $f(x)$ and $g(x)$ are on the same side of $x$, then we can find $\sigma$ as in  Case 1 and conclude that $g\in KfK$;
otherwise the points $f^{-1}(x)$ and $g(x)$ are on the same side and hence, $g\in Kf^{-1}K$.

For condition (ii), note that $f(x)$ and $g(x)$ are on the same side of $x$ by the hypothesis. For instance, assume
$f(x)<x$ and $g(x)<x$. Then we have
\[fg(x)<f(x)<x\quad \textrm{and}\quad gf(x)<g(x)<x.\] Condition (ii) follows from the preceding paragraph.

When applying the same proof of condition (iii) in Case 1 to the present case, 
note first that $u,v_0,v_1$ are on the same side of $x$. 
The only thing that remains to check is the existence of $\sigma_i$. For this, it suffices to impose the following additional conditions below:
\begin{itemize}
\item $p>v_1 \text{ if and only if }v_0>u$;
\item  $p>v_0 \text{ if and only if }v_1>u$;
\item $p,u,v_0,v_1$ are on the same side of $x$.
\end{itemize}
One may succinctly rewrite the first two conditions as 
\[
(p-v_0)(v_1-u)>0\text{ and }(p-v_1)(v_0-u)>0\]
in the special case when $X=\bR$.
Such a point $p$ exists by the density of $X$.
The rest of the proof proceeds in the exactly same manner as Case 1.
 
 {\bf Case 2-2:  $X$ is equipped with a circular order $\le$.}

Since $G$ is $2$--transitive in this case, 
the verification of the two conditions (i) and (ii) is exactly the same as in Case 1.
For the condition (iii), one may follow the proof for Case 2-1 and simply ignore the requirement that
\begin{center}
$p,u,v_0,v_1$ are on the same side of $x$.
 \end{center}
This verifies that $\Stab_G(x)$ is a pre-stabilizer subgroup in all cases.
\ep

We now consider locally transitivity and dilativity in relation to pre-stabilizer subgroups.
\begin{lem}\label{lem:four-ball}
Let $X$ be an infinite connected Hausdorff topological space,
on which a group $G$ acts faithfully by homeomorphisms.
Assume that $K\le G$ is a pre-stabilizer subgroup.
\be[(1)]
\item\label{p:three-orbit} If a normal subgroup $N$ of $G$ admits an orbit $N.x\sse X$ of cardinality at least two,
then there exists an open neighborhood $U$ of $x$ such that 
\[
\left[G_U,G_U\right]\le N.\]
\item\label{p:four-ball} If $G$ is locally dilative, then
some nontrivial element in $K$ fixes a nonempty open set. 
\item\label{p:singleton}
If $G$ is locally transitive, and if $A\sse X$ is a proper, nonempty, closed, $K$--invariant subset, then  $\partial A$ is a singleton
and $K=\Stab_G(\partial A)$.
\ee
\end{lem}
\bp
To prove part~\ref{p:three-orbit}, pick an element $g$ of $N$
and a nonempty open neighborhood $U\sse X$ of $x$ such that $U\cap g(U)=\varnothing$.
Let $h_1,h_2\in G_U$ be arbitrary. 
Since
\[\suppo gh_1^{-1}g^{-1}
= g\suppo h_1 \]
is disjoint from $U$, we see that
\[[[h_1,g],h_2]=[h_1gh_1^{-1}g^{-1},h_2]=[h_1,h_2].\]
By normality, the above commutator belongs to $N$, which proves the first part.

For part~\ref{p:four-ball}, let us fix four disjoint nonempty open sets $U_0,\ldots,U_3$ in $X$.
By local dilativity, we have some $g_i\in G$ such that 
\[\suppc g_i=\suppc g_i^{-1}\] is a connected subset of $U_i$. 
We see that the element $f:=g_2g_3$ is conjugate to none of $g_i^{\pm1}$, by comparing the number of
connected components of respective closed supports. 
Set $V:=U_2\cup U_3$.

Assume for contradiction that the fixed point set of each nontrivial element in $K$ has empty interior. 
Since $\Fix g_i$ contains an open set $U_j$ for $j\ne i$, we have that $f,g_i\not \in K$.
Replacing $g_i$ by $g_i^{-1}$ if necessary, we may 
apply the coset decompositions from Definition~\ref{defn:prestab}
and further require that
\[
g_0,g_1\in KfK.\]
We can also find some 
\[ s_i,t_i\in f^{-1}Kg_i\cap K\]
for $i=0,1$ such that
\[s_0s_1=t_1t_0.\]

We claim that if $\tau\in f^{-1}Kg_i\cap K$ for some $i=0,1$, 
then \[\tau^{-1}(X\setminus \overline V) \sse \suppc g_i\sse U_i.\]
Indeed, pick $\sigma\in K$ such that $\tau=f^{-1}\sigma g_i\in K$. 
If there exists a point 
\[
x\in \tau^{-1}(X\setminus \overline V)\cap (X\setminus  \suppc g_i)\]
then we would have
\[
\tau(x) = f\tau(x)= \sigma g_i(x)
= \sigma(x).\]
This implies that $x$ is an interior point of $\Fix \sigma^{-1}\tau$.
Since $f$ and $g_i$ are not conjugate we have that
 $\sigma^{-1}\tau$ a nontrivial element of $K$, a contradiction.

Applying the above claim to $s_i$ and $t_i$, we compute
\[
s_1^{-1}s_0^{-1}\left(X\setminus \overline V\right)
\sse s_1^{-1}\left({U_0}\right)
\sse s_1^{-1}\left(X\setminus \overline{V}\right)\sse {U_1}.\]
Similarly, we have $t_0^{-1}t_1^{-1}$ maps
$X\setminus \overline V$ into ${U_0}$. 
This contradicts $s_1^{-1}s_0^{-1}=t_0^{-1}t_1^{-1}$.
Part~\ref{p:four-ball} of the lemma follows.

Let us prove part~\ref{p:singleton}.
Since $X$ is connected, the set $A$ is not open. In particular, $\partial A$ is nonempty. If we show that $\partial A$ is a singleton, 
then the maximality of a pre-stabilizer group (Lemma~\ref{lem:prestab-maximal}) would imply that $K=\Stab_{G}(\partial A)$.
Here, note the trivial fact that $G(\partial A)\ne\partial A$.

Assume for contradiction that $\partial A$ contains distinct points $x_0$ and $x_1$. Pick disjoint open neighborhoods
$U_0$ and $U_1$ of $x_0$ and $x_1$. By local transitivity, we have some $h_i\in G_{U_i}$ such that \[h_i(x_i)\in U_i\setminus A.\] 
Since $A$ is $K$--invariant, we have  $h_i\not\in K$.
By the definition of a pre-stabilizer group, we have
\[
h_0\in K h_1^{\pm1}K.\]
Setting $g_0:=h_0$ and $g_1:=h_1^{\pm1}$, we have $g_0\in K g_1 K$.

Since $g_0g_1(x_0)=g_0(x_0)\not\in A$, we have that 
\[f:= g_1g_0=g_0g_1\not\in K.\]
It follows that $g_i\in K f K$ for $i=0,1$, again applying the definition of a pre-stabilizer group.
For $i=0,1$ we have
\[
x_i\in A_i:= g_i^{-1}(X\setminus A)\cap A.\]
It follows that $A_i\sse \suppo g_i\sse U_i$. Set 
\[B:=A_0\cup A_1=f^{-1}(X\setminus A)\cap A.\]

We claim that every $\sigma$ in the set
\[
f^{-1} K g_i\cap K,\]
satisfies $\sigma(A_i)=B$. To see this, let us write
\[
\sigma=f^{-1}\tau g_i\]
for some $\tau\in K$. Then
\[
f\sigma(A_i)=\tau g_i(A_i)\sse \tau(X\setminus A)=X\setminus A.\]
We have that
\[\sigma(A_i)\sse f^{-1}(X\setminus A).\]
Since $\sigma(A_i)\sse A$, we see from the definition of $B$ that 
$\sigma(A_i)\sse B$.
We also see that
\[
g_i\sigma^{-1}(B)=\tau^{-1}f(B)\sse \tau(X\setminus A)=X\setminus A.\]
This implies that $\sigma^{-1}(B)\sse A_i$. This proves the claim.

Since $A_0$ and $A_1$ are nonempty, so is $B$.
Using the condition for a pre-stabilizer group, we can find 
$s_i,t_i\in f^{-1}Kg_i$ such that $s_0s_1=t_1t_0$.
As in the previous part, we have
\[
s_1^{-1}s_0^{-1}(B)
\sse s_1^{-1}(A_0)
\sse s_1^{-1}(B)\sse A_1.\]
Similarly, we have $t_0^{-1}t_1^{-1}$ maps
$B$ into $A_0$.
This is a contradiction, completing the proof of the lemma.\ep

\subsection{Reconstructing a homeomorphism from pre-stabilizer groups}
We can now establish the aforementioned generalization of Filipkiewicz's theorem, namely Theorem~\ref{thm:isom-homeo}.

\bp[Proof of Theorem~\ref{thm:isom-homeo}]
Fix a $G_i$--invariant open basis $\UU_i$ of $X_i$.
For each $y\in X_2$, we set
\[
F^1_y:=\Phi^{-1}(\Stab_{G_2}(y))\le G_1.\]
Similarly for $x\in X_1$ put
\[
F^2_x:=\Phi(\Stab_{G_1}(x))\le G_2.\]
These are pre-stabilizer subgroups
in respective groups, by Lemma~\ref{lem:prestab}.

Define families of open sets
\[
\AAA_x^i:=\left\{ U\in \UU_i \middle\vert \left[ (G_i)_U,(G_i)_U\right]\le F_x^i\right\}\]
for each $i=1,2$ and $x\in X_{3-i}$.
For each $x\in X_1$, $f\in G_1$ and $U\in \AAA_x^2$, we have
\begin{align*}
& \left[ (G_2)_{\Phi(f)U}, (G_2)_{\Phi(f)U}\right]=
\Phi(f) \left[ (G_2)_U,(G_2)_U\right]\Phi(f)^{-1}\\
&\le
\Phi(f)F^2_{x}\Phi(f)^{-1}
=\Phi\left(\Stab_{G_1}(f(x))\right)
=F^2_{f(x)}.\end{align*}
For each $x\in X_{3-i}$, we set
\[
E_x^i=X_i\setminus \bigcup \AAA_x^i.\]
Then $E_x^i$ is a closed $F_x$--invariant  subset of $X_i$.

\begin{claim}
For each $i=1,2$ and $x\in X_i$ the set $E_x^i$ is nonempty.
\end{claim}

We may consider the case $i=2$ only, by symmetry. 
Assume for contradiction that $E_x^2$ is nonempty for some $x\in X_1$. 
The collection $\AAA_x^2$ is an open cover of $X_2$. 
Since $G_2$ is weakly fragmented, we have that 
$F_x^2$ contains a nontrivial normal subgroup $K_2$ of $G_2$.
Let $y\in X_1$ be arbitrary. By the simple transitivity (as is implied by the 3--transitivity) of $G_1$,
we have that some $f\in G_1$ sends $x$ to $y$. 
Then we see that
\[
K_2
=
\Phi(f)K_2\Phi(f^{-1})
= \Phi(f) F^2_x \Phi(f^{-1})=F^2_y.\]
It follows that 
\[
K_2\le \bigcap_{y\in X_1}F^2_y
=\Phi\left(\bigcap_{y\in X_1}\Stab_{G_1}(y)\right).\]
This is a contradiction, since the right hand side is trivial.
The claim is thus proved.

\begin{claim}
There exists some $i\in\{1,2\}$ and $x\in X_{3-i}$ such that
$\AAA_x^i$ is  nonempty.
\end{claim}

Let us pick $y\in X_2$. We have noted above that $F_y^1$ is a pre-stabilizer subgroup of $G_1$.
By part~\ref{p:four-ball} of Lemma~\ref{lem:four-ball}, 
there exists some nontrivial $g_y\in F_y^1$ fixing some open set, say $A_y\in \UU_1$.
Since $g_y$ is nontrivial, we have a proper nonempty closed set
\[
 B:=\Fix \Phi(g_y)\ni y.\] 
 From this,
we define two natural subgroups of $G_1$. Namely, we let
\[H=\Phi^{-1}(\Stab_{G_2}(B)),\quad K=\Phi^{-1}(G_2[X_2\setminus B]).\]
Here, $\Stab_{G_2}(B)$ is the
set-wise stabilizer of $B$.
 Observe that $y\in B$ and that
 \[ g_y\in K\unlhd H\le G_1.\]
Since $g_y$ is centralized by $G_1[A_y]$,
we have that $\Phi(G_1[A_y])$ fixes
$B$ set-wise, so that $G_1[A_y]\le H$. 

We now finish the proof of the claim by analyzing $K$.
Assume first that for some $a\in A_y$, the orbit $K.a$ has cardinality at least two. By part~\ref{p:three-orbit} of Lemma~\ref{lem:four-ball}, 
we can find an open set $A_y'\in \UU$ such that 
\[ a\in A_y'\sse A_y\] and such that 
\[
\left[
(G_1)_{A_y'},(G_1)_{A_y'}
\right]
\le
\left[
H_{A_y'},H_{A_y'}
\right]\le K\le F_y^1.\]
This fits the definition of $\AAA_y^1$, which is now shown to be nonempty.

Lastly, assume that $K$ fixes all points in $A_y$. 
Pick an arbitrary $a\in A_y$ so that $K\le\Stab_{G_1}(a)$.
Since $B$ is a proper closed subset of $X_2$, we have some $V\in \UU_2$ contained in $X_2\setminus B$.
We see that
\[
[(G_2)_V,(G_2)_V]\le (G_2)_V\le (G_2)_{X\setminus B}= \Phi(K)\le F_a^2.\]
This implies that $V\in \AAA_a^2$, completing the proof of the claim.

By symmetry, we may assume that $\AAA_a^2$ is nonempty for some $a\in X_1$. It follows that
$E_a^2$ is a proper nonempty closed $F_a$--invariant subset of $X_2$. 
Applying Lemma~\ref{lem:four-ball} to the group $F_a$ and to the set $E_a^2$, 
we obtain a point $b\in E_a^2$ such that
\[\Phi(\Stab_{G_1}(a))=\Stab_{G_2}(b).\]

We see that there exists a unique bijection
\[
w\co X_1\longrightarrow X_2\]
satisfying 
\[ w(g.a)=\Phi(g).b\]
for all $g\in G_1$. Indeed, $w$ is well-defined since whenever  $g.a=h.a$, the preceding paragraph implies that 
$\Phi(g).b=\Phi(h).b$. Similarly, $w$ is invertible and is therefore a bijection. 
It is immediate to see the relation \[\Phi(g)=wgw^{-1}\] for $g\in G_1$.
Note also that
\[
\Phi(\Stab_{G_1}(x))=\Stab_{G_2}(w(x))\]
for all $x\in X_1$. 
By Lemma~\ref{lem:bij-homeo}, we conclude that $w$ is a homeomorphism as claimed.\ep
This concludes the second proof of the Takens--Filpkiewicz Theorem, and thus the chapter.


%
%
%

\chapter{The $C^1$ and $C^2$ theory of diffeomorphism groups}\label{ch:c2-thry}
\begin{abstract}This chapter concentrates on the interplay of analytic and combinatorial methods for investigating finitely generated groups
of homeomorphisms of one--manifolds, with an emphasis on obstructions to $C^1$ and $C^2$ actions. Much (though not all) of this
chapter is self-contained. Whereas all the results presented herein are useful in various contexts, the $abt$--Lemma is the most important
theorem insofar as future applications are concerned, and its proof occupies the majority of the content. The reader will note that an
overarching theme in this chapter is that regularity greatly constrains the orbit structure of group actions in the presence of partial
commutation.\end{abstract}

\section{Kopell's Lemma}
Kopell's Lemma, formulated by N.~Kopell in~\cite{Kopell1970}, is one of the earliest and most powerful tools for controlling the relationship
between the algebraic structure of groups and regularity of group actions on one--manifolds. It says, roughly, that two commuting 
diffeomorphisms have to have equal or disjoint supports.
Originally formulated for $C^2$ actions, a generalization to $C^{1+\mathrm{bv}}$ actions was given by A.~Navas in his thesis.
We give a proof  here which is a mild generalization using an estimate due to Polterovich and Sodin~\cite{PS2004}.

Recall our convention from the introduction that the regularities we consider are 
always thought of as local properties of functions. This convention also applies to 
diffeomorphisms whose derivatives have bounded variation. The distinction between global and local estimates did not play a role in 
Chapter~\ref{sec:denjoy} where we only considered compact one--manifolds. In general, a 
$C^k$ diffeomorphism $f$ of a possibly non-compact one--manifold $M$ is said to be $C^{k,\mathrm{bv}}$ 
if each point in $M$ has a neighborhood $U$ where $f^{(k)}$ has a bounded total variation, denoted by
$\operatorname{Var}\left(f^{(k)};U\right)$. This way, we have that a $C^{k+1}$ 
diffeomorphism of a one manifold is necessarily $C^{k,\mathrm{bv}}$.

In the theorem below as well, each element $f$ in $\Diffb[0,1)$ satisfies \[ \Var(f';[0,\delta])<\infty\] for all $\delta\in(0,1)$.
We emphasize again that  the supremum
\[
\sup_{\delta>0}\Var(f';[0,\delta])\] is allowed to be infinite.

\begin{thm}[Kopell's Lemma]\label{thm:kopell}
Let $f\in\Diffb[0,1)$ and $g\in\Diff_+^1[0,1)$  be commuting diffeomorphisms such that \[\Fix(f)\cap (0,1)=\varnothing\quad \textrm{and}
\quad\Fix(g)\cap  (0,1)\neq\varnothing.\] Then $g$ is the identity.
\end{thm}

To establish the above, let us note two elementary facts about $C^1$ diffeomorphisms.

\begin{lem}\label{lem:diff-fix}
If $g$ is a $C^1$ diffeomorphism of a one--manifold $M$,
and if $x$ is an accumulation point of $\Fix g$ in $M$,
then we have that $g'(x)=1$.
\end{lem}
\bp
Let $\{x_n\}_{n\ge1}$ be a sequence of fixed points of $g$ converging to $x$. The Mean Value Theorem implies that for some $y_n$
lying between $x$ and $x_n$, we have \[g'(y_n)=\frac{g(x_n)-g(x)}{x_n-x}=1.\]
By the continuity of $g'$ we have that
 \[
 g'(x) = \lim_{n\to\infty} g'(y_n)=1,\]
 which establishes the lemma.
\ep

The following is essentially observed by Polterovich and Sodin~\cite{PS2004}.
\begin{lem}\label{lem:inf-finite}
If $g$ is a nontrivial orientation preserving $C^1$ diffeomorphism of a bounded interval $J$
then we have that
\[
\sum_{n\in\bZ}\inf\left\{\left(g^n\right)'(x)\mid x\in J\right\} <\infty.\]
\end{lem}
\bp
Switching $g$ and $g^{-1}$ if necessary
we may assume to have some $c\in J$ such that $g(c)>c$. 
Put $K:=[c,g(c)]$. By the Mean Value Theorem, we have some $x_n\in K$
such that
\[
|J|\ge \sum_{n\in\bZ} \abs*{g^n(K)}= \sum_{n\in\bZ}|K|\cdot \left(g^n\right)'(x_n)
\ge 
|K|\cdot \sum_{n\in\bZ}\inf\left\{\left(g^n\right)'(x)\mid x\in J\right\}.\]
This implies that the given infinite sum is at most $|J|/|K|$.
\ep

As was implicitly used in the proof of Lemma~\ref{lem:inf-finite},
whenever $f$  is a fixed point free diffeomorphism of $(0,1)$, then $f$ admits a ~\emph{fundamental domain}\index{fundamental domain},
which is a half-open interval
$J\sse (0,1)$ such that $f^n(J)\cap J=\varnothing$ if $n\in\Z\setminus\{0\}$, and such that \[(0,1)=\bigcup_{n\in\Z} f^n(J).\] In  fact,
one may set $J=[c,f(c))$ or $J=[f(c),c)$ for an arbitrary $c\in (0,1)$.

The following lemma considers a more general situation than Kopell's Lemma, as it 
only assumes that $g$ is $C^1$ in the open interval $(0,1)$.
Roughly, it shows that 
\[ \liminf_{y\to+0} \left(g^n\right)'(y)\]
is less than $1/|n|$ for almost all $n\in\bZ$, even when $g$ is not differentiable at $0$; the interested may wish to compare
with Lemma~\ref{l:sum-density}.

\begin{lem}\label{lem:kopell0}
If $f\in\Diffb[0,1)$ and $g\in\Diff_+^1(0,1)$ are nontrivial, commuting diffeomorphisms
such that \[\Fix(f)\cap (0,1)=\varnothing\quad \textrm{and}
\quad\Fix(g)\cap  (0,1)\neq\varnothing.\] 
Then we have the estimate
\[
\sum_{n\in\bZ}  \liminf_{y\to+0} \left(g^n\right)'(y)<\infty.\]
\end{lem}
\bp
Replacing $f$ by its inverse if necessary, we assume that $f(x)<x$ for all $x\in (0,1)$.
Pick a fundamental domain $J=[f(c),c)$ of $f$, for some $c\in (0,1)$ to be chosen later.
It is easy to see that $\log f'$ has bounded variation on $[0,1)$ whenever $f'$ does as well. 
So, for each $n\ge0 $ and for each $a,b\in J$ we have that
\[
\log\frac{(f^n)'(a)}{(f^n)'(b)}
=\sum_{i=0}^{n-1} \abs*{\log|f'(f^i(a))|-\log |f'(f^i(b))|}
\le v:=\Var\left(f;[0,c]\right)<\infty.\]

Let us now require \emph{a priori} that $c\in\Fix g$. 
We have that $g(J)=J$ since \[gf(c)=fg(c)=f(c).\]
Moreover, for arbitrary $m,n\in\bZ$ and $x\in J$,
we obtain the estimate

\[
\left(g^n\right)'(x) = \frac{\left(f^m\right)'(x)}{ \left(f^m\right)'\left( g^n(x)\right)}\cdot \left(g^n\right)'\left( f^m(x)\right)\ge e^{-v} \cdot \left(g^n\right)'\left( f^m(x)\right).\]
Here, the last inequality comes from the fact that $g^n(x)\in J$.
By fixing $n\in\bZ$ and sending $m\to\infty$, we see that
\[
\left(g^n\right)'(x) \ge 
e^{-v}  \liminf_{y\to+0} \left(g^n\right)'(y).\]
Assuming $g$ is nontrivial, we now have that
\[ \sum_n\liminf_{y\to+0} \left(g^n\right)'(y)
\le \sum_n e^v \inf\left\{\left(g^n\right)'(x) \mid x\in J\right\} <\infty,\]
by Lemma~\ref{lem:inf-finite}.
\ep

\begin{proof}[Theorem~\ref{thm:kopell}]
Let $c\in(0,1)$ be a fixed point of $g$.
Since $f$ and $g$ commute, the points $\{f^n(c)\}_{n\in\bZ}$ are also fixed by $g$. In particular, $0$ is an accumulation of $\Fix g =\Fix g^n$.
Lemma~\ref{lem:diff-fix} implies that $\left(g^n\right)'(0)=1$ for all $n$.
By Lemma~\ref{lem:kopell0}, this implies that $g$ must be trivial.\ep

One of the most immediate applications of Kopell's Lemma is the Plante--Thurston theorem~\cite{PT1976}, 
which asserts that $C^{1+\mathrm{bv}}$ actions of nilpotent groups on the interval or on the circle factor 
through abelian quotients. We closely follow an argument given by Navas~\cite{Navas2011} (cf.~\cite{KoberdaSurv20}).

\begin{thm}[Plante--Thurston Theorem]\label{thm:plante-thurston}
If $M\in\{I,S^1\}$ then every nilpotent subgroup of  $\Diff_+^{1+\mathrm{bv}}(M)$ is abelian.
\end{thm}
\begin{proof}
Let $N\le \Diff_+^{1+\mathrm{bv}}(M)$ be nilpotent.

{\bf Case 1: $M=[0,1]$.} It suffices to prove that $N$ acts freely on $I$, by H\"older's Theorem (Theorem~\ref{thm:holder}).
Clearly we may assume that $N$ has
no global fixed points in $(0,1)$, and we will use this to derive a contradiction. Since $N$ is nilpotent there is a nontrivial element
$g$ contained in the center $Z(N)\le N$.
Suppose now that $f\in N$ is nontrivial and has at least one fixed point in $(0,1)$. Since $f$ is assumed to
be nontrivial, we may choose a fixed point $x\in \partial\Fix(f)\cap (0,1)$.

We have that $g(x)=x$. Indeed otherwise we may iterate $g$ on $x$ and extract limit points $a$ and $b$ of $\{g^n(x)\}_{n\in\Z}$,
and $[a,b]$ is a nondegenerate $g$--invariant interval. Note that by construction, $g$ acts on $(a,b)$ without fixed points.
Since $f$ commutes with $g$, we have that $[a,b]$ is also $f$--invariant.
Kopell's Lemma now implies that $f$ is the identity, whence we obtain $g(x)=x$.

The previous argument shows that if the action  of $N$ is not free then every nontrivial element of $Z(N)$ has a fixed point in $(0,1)$.
So, let $y\in\partial\Fix(g)\cap (0,1)$. If $h\in N$ is nontrivial, then the argument of the previous paragraph (switching the roles of
$f$ and $g$ with $g$ and $h$ respectively) proves that $h(y)=y$. Since $h$ is arbitrary, $y$ is a global fixed point of the  action of $N$,
a contradiction.

{\bf Case 2: $M=S^1$.} First, we have that $N$ preserves a probability measure $\mu$ on $S^1$, as follows from the amenability of $N$.
Indeed, $N$ is nilpotent and hence solvable, and so Theorem~\ref{thm:kaku-mar} implies that $N$ preserves a probability measure on $S^1$.
Theorem~\ref{thm:invt-homo} implies that the
rotation number thus  gives a homomorphism from $N$ to the circle group. If there is an element $\nu\in N$ with irrational
rotation number, then $\nu$ is topologically conjugate to an irrational rotation, by Denjoy's Theorem (Theorem~\ref{thm:denjoy}.
Since $\nu$ is uniquely ergodic by Theorem~\ref{thm:irr-ue},
it follows that the invariant measure $\mu$
 is the pushforward of Lebesgue measure by a homeomorphism (see Lemma~\ref{lem:leb-push}). In this case, $N$ is conjugate to a group
of rotations of the circle and is therefore abelian. We leave the proof that $N$ is conjugate to a group of rotations as an exercise for the
reader, or refer the reader to Proposition 1.1.1 of~\cite{Navas2011}.

We may therefore assume that all elements of $N$ have rational rotation number.
If $\nu\in N$ has a fixed point in $S^1$ and if $J\sse S^1$ is an interval such that $J\cap \nu J=\varnothing$,
then $J$ cannot be given positive
measure under the $N$--invariant measure on $S^1$. It follows easily that the support of $\mu$ is contained in the intersection of all
periodic orbits of $N$.
If $N$ is nonabelian and nilpotent, then by taking a commutator of elements in the penultimate term of the lower central series of $N$,
we may find a nontrivial central element $h\in Z(N)$ which is a commutator in $N$.
We write $h=[f,g]$. Since the rotation number is a homomorphism on $N$, we have that the rotation number of $h$ is zero, whence
$h$ has a fixed point (see Proposition~\ref{prop:rot-easy}).
An easy computation shows that \[h^{mn}=f^{-n}g^{-m}f^ng^m.\] We suppose that $x_0$ is a point in the support
of $\mu$, so that $x_0$ is a periodic point of both $f$ and $g$. For suitably chosen
 values of $n$ and $m$, we see that $x_0$ is a fixed point of the subgroup $N_0=\form{ f^n, g^m}$. We have that $N_0$ is
 nilpotent and nonabelian. Indeed, if $N_0$ were abelian then $h^{mn}$ would be the identity. Since $h$ is nontrivial and fixes a point
 in $S^1$, it must have infinite order by Proposition~\ref{prop:l-order-homeo}, a contradiction.
 
 Cutting the circle open at $x_0$, we see that the subgroup of $N_0$
 acts on the resulting compact interval. The case where $M=[0,1]$ now implies that $N_0$ is abelian, and this is the desired contradiction.
 It follows that $N$ must be abelian.
\end{proof}

\section{(Residually) nilpotent groups acting by $C^1$ diffeomorphisms}\label{sec:nilpotent}

In this section, we will showcase some constructions of faithful $C^1$ actions of groups on $I$ and $S^1$. This section will therefore be
more expository than the others in this book, and the reader will find less detailed (or entirely absent) proofs the results discussed here.

The Plante--Thurston Theorem gives the first nontrivial algebraic restriction of groups that admit faithful representations into $\Diffb(M)$.
This means that the critical regularity of a nonabelian nilpotent group is at most two. It is therefore an interesting question to determine where
exactly the critical regularity of nonabelian nilpotent groups lies.

We have that $\Homeo_+[0,1]$ is torsion--free since a finite cyclic group cannot admit a left invariant ordering
(see Proposition~\ref{prop:l-order-homeo}), and combining this with H\"older's Theorem for the circle (Theorem~\ref{thm:holder-circle})
and the characterization of fixed points via rotation number (Propositions~\ref{prop:rot-easy} and~\ref{prop:rot-easy}), we have
that an arbitrary finite subgroup of $\Homeo_+(S^1)$ is cyclic.

If $N$ is a finitely generated nilpotent group, then certainly $N$ might have torsion. However, it is not difficult to prove that in $N$, the
product of two torsion elements is again torsion, and so the torsion elements of $N$ form a normal subgroup which is itself finite, as is
easily checked. It
follows that $N$ admits a torsion--free quotient with finite kernel. Since finitely generated nilpotent groups can be shown to be
residually finite (since they admit faithful homomorphisms into matrix groups, finitely generated subgroups of which are always
residually finite),
$N$ admits a finite index subgroup that is torsion--free. So, in studying actions of finitely generated nilpotent groups on
manifolds, one loses relatively little of the dynamical and algebraic richness by assuming that $N$ is torsion--free.

Now, if $N$ is an arbitrary torsion--free countable nilpotent group, then $N$ admits a left invariant ordering and hence embeds
into both $\Homeo_+[0,1]$ and $\Homeo_+(S^1)$ (Proposition~\ref{prop:l-order-homeo} and Proposition~\ref{prop:c-order-homeo}).
To see that $N$ is orderable,
it suffices to prove the claim for finitely generated subgroups of $N$, since
a group admits a left invariant ordering if and only if every finitely generated subgroup admits such an ordering. Then, we have that
$N$ admits quotient by a finitely generated
central subgroup $Z$ such that $N/Z$ is again torsion--free, and such that the lower central series of $N/Z$ is
strictly shorter than that of $N$. Since finitely generated torsion--free abelian groups admit left invariant orderings, we may assume that
$Z$ and $N/Z$ admit such orderings by induction. We then order $N$ lexicographically; that is $n_1<n_2$ if the image of $n_1$ is less
than the image of $n_2$ in $N/Z$. If the images are equal then $n_1<n_2$ if $n_1^{-1}n_2$ is positive in $Z$. It is easy to check that
the resulting ordering is left invariant on $N$.

A natural question then follows: if $N$ is finitely generated and torsion--free nilpotent, can $N$ be realized as a subgroup of $\Diff_+^1(M)$
for $M\in\{I,S^1\}$? The answer turns out to be yes, by a result of Farb--Franks~\cite{FF2003}, though this
fact is far from obvious. One can find a further generalization of this result in~\cite{Jorquera,CJN2014,JNR2018}.

\subsection{The Farb--Franks Theorem and the universal nilpotent group}
We will give an outline of the argument due to Farb and Franks that ansers the previous question.
Let $N_m\le \SL_{m+1}(\bZ)$ denote the group of lower triangular matrices
with integer entries, and where all diagonal entries are $1$. It is a standard computation that $N_m$ is nilpotent, and that the lower central
series of $N_m$ has length $m$. 
We have a finite generating 
 \[\{u_{1,m},\ldots,u_{m,m}\}\]
 of the group $N_m$,  
 where the matrix
$u_{i,m}\in\SL_{m+1}(\bZ)$ consists of ones down the diagonal, and whose only nonzero entry is a one in position $(i,i-1)$.

If $N$ is a finitely generated, torsion--free, nilpotent group, then there is an $m$ for which $N$ embeds into $N_m$, as follows from
a classical result of Mal'cev (see~\cite{Raghunathan1972}, for instance). One way to show this is to build a certain nilpotent real Lie group
$N\otimes\bR$, called the \emph{Mal'cev completion}\index{Mal'cev completion} 
of $N$. Then one shows that $N\otimes\bR$ embeds in the group $N_m(\bR)$
of lower triangular real $(m+1)\times (m+1)$ matrices with $1$ on the diagonal, and that $N$ lies in the integer points of this image.

Another perspective on this fact is as follows.
Recall, the \emph{Hirsch length}\index{Hirsch length} (or, the \emph{polycyclic length}\index{polycyclic length}) of a group $G_0$ is the 
 the length of a subnormal
series \[\{G_k\leq G_0\}_{k\ge1}\] such that $G_i/G_{i+1}$ is infinite cyclic. 
For a group $N$ as in the preceding paragraph, one can find an exact sequence of the form
\[1\longrightarrow K\longrightarrow N\longrightarrow \bZ\longrightarrow 1\]
such that $K$ is a nilpotent group whose Hirsch length is strictly shorter than that of $N$. 
By induction we have that $K$ embeds in $N_m$ for some $m$.
Since $\bZ$ is cyclic and torsion--free, the surjection $N\longrightarrow \bZ$ splits, and so $N$ has the structure of a semidirect product of
$K$ with $\bZ$. Choosing a generator $t$ for $\bZ$, we have that conjugation by $t$ is an automorphism of $K$. This automorphism
must act on the group $H_1(K,\bQ)$ by a unipotent matrix. This will guarantee that $N$ is nilpotent, as 
is easily verified by computing commutators
in the semidirect product $\yt{K}$ of $H_1(K,\bZ)$ with $\bZ$, the group $\yt{K}$ also necessarily being 
nilpotent since it is quotient of $N$. It is then possible
to embed $K$ in a (possibly larger) $N_m$, and to realize $t$ as global conjugation by a lower triangular matrix in $N_m$.

Let $P$ be a group theoretic property. We say a group $G$ is \emph{residually $P$}\emph{residual property} 
if for each nontrivial element $g\in G$ 
there exists some quotient $Q$ of $G$ having the property $P$ such that the image of $g$ in $Q$ is nontrivial.
Note that $G$ is residually $P$ if and only if $G$ embeds into the direct product of groups that have the property $P$. 

For instance, a residually torsion--free nilpotent group is a group $G$ such that every nontrivial element of $G$ survives in a 
torsion--free nilpotent quotient of $G$. The class of residually torsion--free nilpotent groups is quite extensive, and  includes 
free groups, fundamental groups of closed surfaces, right-angled Artin groups, pure braid groups, and Torelli groups of 
surfaces; see~\cite{BP2009JKTR} for instance. 

In~\cite{FF2003}, Farb and Franks proved the following remarkable result and initiated the study of 
smooth one--manifold actions of nilpotent groups, which resulted in many interesting 
discoveries~\cite{Navas2008GAFA,Jorquera,ParkheETDS,CJN2014,JNR2018}.

\begin{thm}[\cite{FF2003}]\label{thm:ff2003}
Every finitely generated residually torsion--free nilpotent group embeds into  $\Diff_+^1[0,1]$.
\end{thm}

Recall $\Diff^1_K(\bR)$ denotes the group of $C^1$--diffeomorphisms of $\bR$ supported in $K\sse\bR$.
The key step of the proof is the following.
\begin{prop}[\cite{FF2003}]\label{prop:ff2003}
For each $m\in\bN$ and $\epsilon>0$, there exists an embedding
\[
\rho=\rho_{m,\epsilon}\co N_m\longrightarrow \Diff^1_{[0,1]}(\bR)\]
such that for all $i=1,\ldots,m$ we have
\[
|\rho(u_{i,m})'-1|\le \epsilon.\]
\end{prop}

\bp[Proof of Theorem~\ref{thm:ff2003}, assuming Proposition~\ref{prop:ff2003}]
Let $G$ be a residually torsion-free nilpotent group with a finite generating set $\{a_1,\ldots,a_k\}$. We have a sequence of homomorphisms
\[
\phi_i\co G\longrightarrow N_{n(i)}\]
for each $i=1,2,\ldots$
and for some positive integers $n(1),n(2),\ldots$ 
 such that \[\bigcap_i\ker\phi_i=\{1\}.\]
By Proposition~\ref{prop:ff2003}, for each interval $K_i:=[1/(i+1),1/i]$
we can pick an embedding
\[
\psi_i\co N_{n(i)}\longrightarrow \Diff^1_{K_i}(\bR)\]
such that 
\[
|\left(\psi_i\circ\phi_i(a_j)\right)'-1|\le 1/i\]
for all $i\ge1$ and $1\le j\le k$. 

Let $g\in G$.
We define 
\[\phi\co G\longrightarrow \Homeo_{[0,1]}(\bR)\]
by the infinite product
\[
\phi(g):=\prod_i \psi_i\circ\phi_i(g).\]
This map $\phi$ is injective since each $g\in G$ survives under some $\phi_i$.

Since $\supp\phi(g)$ is contained in the disjoint union
\[K_0:=\bigsqcup_i \Int(K_i),\] it is obvious that $\phi(g)$ is $C^1$ on $K_0$. 
Since 
\[\psi_i(h)'(\partial K_i)=\{1\}\] for all $h\in N_{n(i)}$, we also have that 
$\phi(g)$ is $C^1$ at $x\in 1/i$, for each $i\ge1$. 

The only place remaining to check the $C^1$ regularity of $\phi(g)$ is $x=0$. 
Pick a point $t\in K_i$. We first have
\[
1\le\frac{\phi(g)(t)}{t} \le \frac{1/i}{1/(i+1)},\]
which implies that $\phi(g)'(0)=1$. 
Moreover, if we let $\ell$ be the word length of $g$ in the given generating set
then 
\[
(1-1/i)^{\ell}\le |\phi(g)'(t)|\le (1+1/i)^\ell.\]
We see that $\lim_{t\to+0} \phi(g)'(t)=1$, which establishes that $\phi(g)$ is $C^1$ at $x=0$.

\ep

Let us now sketch some of the ideas behind the proof of the Proposition~\ref{prop:ff2003}. 
The starting observation in~\cite{FF2003} is that the group $N_m$ acts faithfully on the real line. An even stronger conclusion
can be formulated as follows:
there exists a natural embedding $\iota\co N_m\longrightarrow N_{m+1}$ defined by
\[
\iota(A)_{ij}:=\begin{cases} 
A_{ij},&\text{ if }1\le i,j\le m+1\\
\delta_{ij},&\text{ if }i=m+2\text{ or }j=m+2,\end{cases}
\]
where $\delta_{ij}$ denotes the Kronecker delta.
There also exists a left inverse $\pi$ of $\iota$ which forgets the $(m+2)$--th row and column.
From these maps $\iota$ and $\pi$, we can define direct and inverse limits 
\begin{align*}
\Lambda&:=\varinjlim N_m,\\
\Gamma&:=\varprojlim N_m.\end{align*}
Each element of the group $\Gamma$ can be regarded as an $\bZ_{>0}\times\bZ_{>0}$ lower triangular matrix with ones on the diagonal.
We have natural inclusions
\[
N_m\hookrightarrow \Lambda\hookrightarrow \Gamma\hookrightarrow \prod_m N_m.\]

It is simple to check that the group $\Gamma$ acts faithfully on 
\[
X_\infty:=\{(1,x_0,x_1,\ldots)\mid x_i\in\bZ\}\sse \bZ^{\bN}.\]
It is important for us that the lexicographical order on $X_\infty$ is preserved by the action of $\Gamma$.

Let us introduce a notation that for two intervals $I,J\sse\bR$ we write $I\le J$ if $\sup I\le\inf J$. We partition $I^*=[0,1]$ into compact intervals $I_i$ as an interior--disjoint union
\[
I^*=\bigcup_{i\in\bZ} I_i\bigcup\{0,1\}\]
such that $I_i\le I_j$ for $i\le j$. 
For each vector $v\in\bZ^m$ with $m\ge1$ we partition again
\[
I_v=
 \partial I_v\cup
\bigcup_{i\in\bZ} I_{v,i}\]
as an interior--disjoint union.
We will also choose the lengths of the intervals such that
\[
\lim_{m\to\infty} \sup_{v\in\bZ^m} |I_v|=0.\]

If $v=(x_0,x_1,\ldots)\in\bZ^{\bN}$ then we have nested compact intervals
\[
I_{x_0}\supseteq I_{x_0,x_1}\supseteq I_{x_0,x_1,x_2}\supseteq\cdots.\]
So, there uniquely exists a point $z_v$ in the union of the above compact intervals. 
We have a faithful action of $\Gamma$ on the set $\{z_v\mid v\in\bZ^{\bN}\}$.
Namely,  for each $A\in\Gamma$ and $v\in\bZ^{\bN}$, we define
\[
A(z_v)=z_{A.v},\]
where $A.v$ is the image of $v\in X_\infty$ under $A$;
in other words,  we have an infinite matrix multiplication
\[
A\begin{pmatrix} 1\\ v\end{pmatrix}=\begin{pmatrix} 1\\ A.v\end{pmatrix}.\]
This action preserves the lexicographical order on $\bZ^{\bN}$,
and hence, naturally extends to the whole interval $I$, so that 
\[\Gamma\hookrightarrow\Homeo_+[0,1].\]

From this we see that the uncountable group $\Gamma$ admits a left invariant ordering; see Proposition~\ref{prop:l-order-homeo}.
 Because $\Gamma$ accommodates all finitely generated torsion--free nilpotent groups, we have that $\Gamma$ itself is not a nilpotent group.
We remark that every finitely generated residually torsion--free nilpotent group embeds into $\Gamma$ using similar ideas.

We now return to the proof of Proposition~\ref{prop:ff2003}.
For this we fix an integer $m\ge1$. 
We will consider the collection of intervals
\[
\II_m:=\bigcup_{v\in\bZ^m}\{I_v\}.\]
The goal is to define an action of $N_m$ on this collection of intervals
in a consistent way with the action of $\Gamma$,
and  extend this action to an action on $I$. For this, it will be useful to have a uniform way of 
choosing a diffeomorphism for an arbitrary pair of intervals. 

We use a certain equivariant family defined by Yoccoz (cf.~\cite{Navas2011}, Definition 4.1.23).
An \emph{equivariant family}\index{equivariant family} of homeomorphisms
\[\{\phi_{a,b}\colon [0,a]\longrightarrow [0,b]\mid a,b>0\},\] is one which satisfies
\[\phi_{b,c}\circ\phi_{a,b}=\phi_{a,c}\] for all $a,b,c>0$. We have already used one equivariant family when we constructed the first
examples of continuous Denjoy counterexamples in Chapter~\ref{sec:denjoy}; the family we used there was particularly simple, as they
were simply linear scalings. As we remarked then, such equivariant families are perfectly good for constructing continuous actions,
though they do not have good smoothness properties since the slope of $\phi_{a,b}$ near zero is $b/a$, and so if two different such
families are operating on intervals that share a boundary point, the derivatives on both sides will not agree unless the two families
scale by the same factor. This lattermost condition, as one might imagine, is usually too much to require.

Yoccoz's more sophisticated family is defined as follows. We define \[\phi_a\colon (0,a)\longrightarrow\bR\] by setting
\[\phi_a(x)=-\frac{1}{a}\cot\left(\frac{\pi\cdot x}{a}\right).\] The reader can check that this is a homeomorphism that preserves the usual
order on both $(0,a)$ and $\bR$. 

The map $\phi_{a,b}$ is defined by $\phi_b^{-1}\circ \phi_a$, and thus furnishes a $C^{\infty}$
diffeomorphism $(0,a)\longrightarrow (0,b)$, and is easily seen to be equivariant. In fact, the equivariance of this family has very little
to do with the specific nature of the map $\phi_a$; the definition of $\phi_{a,b}$ is used to analyze the asymptotic regularity properties
of the family. For $0<x<a$ we can compute
\[
\phi_{a,b}(x) =\frac{b}2-\frac{b\arctan(b/a\cdot \cot(\pi x/a))}{\pi}
=
 \int_0^x \left( 1 -
\frac{1-b^2/a^2}{1+b^2/a^2\cdot \cot^2(\pi x/a)}\right)dx.\]

An easy calculation using the chain rule shows that \[\phi'_{a,b}(x)=\frac{(\phi_a(x))^2+1/a^2}{(\phi_a(x))^2+1/b^2}.\] In particular,
\[\lim_{x\to 0}\phi'_{a,b}(x)=\lim_{x\to a}\phi'_{a,b}(x)=1,\] and so one can extend $\phi_{a,b}$ to a $C^1$ diffeomorphism $[0,a]
\longrightarrow [0,b]$. One can in fact prove that $\phi_{a,b}$ is a $C^2$ diffeomorphism, with second derivatives vanishing at $0$ and $a$.

The relevant properties of the family $\{\phi_{a,b}\}_{a,b>0}$ can be summarized as follows:

\begin{prop}\label{prop:yoccoz}
The equivariant family $\{\phi_{a,b}\}_{a,b>0}$ has the following properties.
\begin{enumerate}[(1)]
\item
For $a,b>0$, we have $\phi'_{a,b}(0)=\phi'_{a,b}(a)=1$.
\item
We have that \[\sup_{x\in[0,a]}|\phi'_{a,b}(x)-1|=\left|\frac{b^2}{a^2}-1\right|.\]
\item
For all $x\in [0,a]$, we have that \[|\phi''_{a,b}(x)|\leq \pi\frac{(\phi_a(x))^2+1/a^2}{(\phi_a(x))^2+1/b^2}\left|\frac{1}{b^2}-\frac{1}{a^2}\right|b.\]
\end{enumerate}
\end{prop}

We omit the details of Proposition~\ref{prop:yoccoz}, instead directing the reader to the discussion in~\cite{Navas2011,FF2003,KK2020-DCDS}. 

Now, to construct group actions on $I$, we use the Yoccoz family. If $I=[y,y+a]$ and $J=[z,z+b]$ are compact intervals with $a,b>0$,
then we obtain an identification of these intervals via the family
in the obvious way, namely by \[\phi^I_J\colon I\longrightarrow J,\quad \phi^{I}_{J}(x)=
\phi_{a,b}(x-y)+z.\]

Let us now assume that $A\in N_m\le\SL_{m+1}(\bZ)$, 
and $x\in I$. 
Suppose first that $x\in I_v$ for some $v\in\bZ^m$.
Write $w:=A.v\in\bZ^m$, for the action of $A$ on $\bZ^m$ described above.
Then we define
\[
A(x):=\phi^{I_v}_{I_w}(x).\]
We can write $I^*=\bigcup_{v\in\bZ^m}I_v\cup J$ for some countable set $J$,
and hence, this action $A$ extends to the whole interval $I^*$.

By Proposition~\ref{prop:yoccoz}, the moment
we piece elements in the Yoccoz equivariant family in a way that is continuous (which will be immediate from the construction),
the result will automatically be a diffeomorphism in the interior of the interval, at least away from the accumulation points of intervals
of the form $I_v$.
There is an issue with continuity of the derivative
at the accumulation points,
which is why the choices of lengths of the intervals (which we will not specify here)
and the particulars of the Yoccoz equivariant family are especially relevant.

By the faithfulness of the action of $\Gamma$ on $I$, we see that 
the action $x\mapsto A(x)$ defines an embedding
\[
N_m\longrightarrow\Homeo_+[0,1].\]

The rest of the proof of Proposition~\ref{prop:ff2003} consists of
technical assignment of lengths of $I_v$
and verifying (one by one!) the regularity of each generator of $N_m$.
We will not reproduce proofs of the more time--consuming calculations, though the ambitious reader may wish to try
them. See~\cite{FF2003,Jorquera,CJN2014,JNR2018} for related computations.
We will conclude the discussion by simply pointing out a key step for the proof.

\begin{lem}\label{lem:ff2003}
Let $m\ge2$ and $\epsilon>0$. 
Then for some assignment of lengths for $\{I_v\}_{v\in\bZ^m}$
and for the action 
of 
\[
N_m=\{ u_{1,m},\ldots,u_{m,m}\}\]
defined as above,
we have that
\[
\sup \left\{
|u_{i,m}'(x)-1|
\middle\vert
{i=1,\ldots,m}\text{ and }{x\in I^*} 
\right\}<\epsilon.\]
\end{lem}

\subsection{Topological conjugacy of virtually nilpotent group actions}\label{ss:parkhe}

Here we will briefly summarize some stronger results about nilpotent groups acting on one--manifolds, due to Parkhe, Castro, Jorquera,
Navas, and Rivas,
which appear in~\cite{ParkheETDS,CJN2014,JNR2018}.

In~\cite{ParkheETDS}, Parkhe proves the following remarkable
structure theorem for nilpotent groups actions on $I$ and $S^1$, which says that the construction of the action in Theorem~\ref{thm:ff2003}
is a common feature of all nilpotent group actions on one--manifolds:

\begin{thm}[See Theorem 1.1 in~\cite{ParkheETDS}]\label{thm:parkhe-main}
Let $M$ be a one--manifold and let $N\le \Homeo_+(M)$ be a finitely generated subgroup that contains a nilpotent group with finite index.
Then there exists countable collections of open sets $\{I_k\}_{k\ge1}$ and $\{J_k\}_{k\ge1}$ which satisfy the following conditions:
\begin{enumerate}[(1)]
\item
We have \[I=\bigcup_{k\ge1} I_k\quad \textrm{and}\quad J=\bigcup_{k\ge1} J_k\] are disjoint from each other.
\item
For $k\neq\ell$, we have \[I_k\cap I_{\ell}=J_k\cap J_{\ell}=\varnothing.\]
\item
For each $k$, the sets $I_k$ and $J_k$ are $N$--invariant.
\item
We can write \[I_k=\bigcup_{\ell} I_{k,\ell},\] where $I_{k,\ell}$ is a nonempty open interval and where the union defining $I_k$ is disjoint, such that
for all indices $\ell$ and $\ell'$, there
exists an element $g\in N$ satisfying that \[g(I_{k,\ell})=I_{k,\ell'},\] and that the stabilizer of $I_{k,\ell}$ in $N$ is trivial.
\item
For all $k$, the action of $N$ on $J_k$ is minimal, and the stabilizer in $N$ of each component of $J_k$ is abelian.
\item
We have $I\cup J$ is dense in $M$.
\end{enumerate}
\end{thm}

The components of $I$ are the ones occurring in the construction in the proof of Theorem~\ref{thm:ff2003}, and the components of $J$
account for the possibility that the group $N$ can surject to an abelian group of rank at least two, which can act minimally on the interval
and on the circle.

The usefulness of Theorem~\ref{thm:parkhe-main} to the problem of smoothability may not be immediately apparent, though it is evident
that it highly constrains the structure of orbits of actions by \emph{virtually}\index{virtual property}
nilpotent groups (i.e.~ones that contain nilpotent groups with finite
index).

To state the salient corollary to Theorem~\ref{thm:parkhe-main}, we recall a basic notion from geometric group theory. Let $G$ be a group
generated by a finite set $S$ satisfying $S=S^{-1}$. There is a natural metric on $G$, given by defining \[|g|_S=\min\{n\mid g=s_1\cdots s_n,\,
s_i\in S\},\] and then by defining $d_S(g,h)=|g^{-1}h|_S$. This metric depends on $S$, but its bi--Lipschitz equivalence class does not.

For $n\ge1$, we write \[B_{S}(n)=\{g\mid\,\, |g|_S\leq n\},\] and $b_{S}(n)=|B_{S}(n)|$.
The function $b_{S}(n)$ counts how many distinct group
elements in $G$ can be written as a product of at most $n$ generators. Of course $b_S(n)$ depends on $S$, but its coarse behavior does not.
We say that a finitely generated group $G$ has \emph{polynomial growth}\index{polynomial growth}
if $b_S(n)$ is dominated by a polynomial function $p(n)$ for
some (equivalently all) generating sets $S$. Since the coefficients of $p$ are sensitive to the choice of $S$ and since only the highest
degree term of $p$ is salient, we say that $G$ has polynomial growth of degree at most $d$ if $b_S(n)=O(n^d)$.

For a reader familiar with nilpotent groups, it is not difficult to prove that a finitely generated nilpotent group has polynomial growth, and
more generally that a finitely generated virtually nilpotent group has polynomial growth. It is a foundational result of Gromov that
the converse is also true.

\begin{thm}[See~\cite{gromov-poly,Kleiner-poly,TaoShalom-poly}]\label{thm:gromov-poly}
Let $N$ be a finitely generated group that has polynomial growth. Then $N$ is virtually nilpotent.
\end{thm}

Parkhe deduces the following consequence of Theorem~\ref{thm:parkhe-main}:

\begin{thm}[Theorem 1.4 in~\cite{ParkheETDS}]\label{thm:parkhe-conj}
Let $N\le \Homeo_+(M)$ be a finitely generated group of polynomial growth $O(n^d)$. Then for all $\tau<1/d$, the group $N$ is topologically
conjugate into $\Diff_+^{1,\tau}(M)$.
\end{thm}

Theorem~\ref{thm:parkhe-conj} shows that there are essentially no non--$C^1$ actions of a finitely generated nilpotent group $N$ on $M$;
any failure of differentiability results from a poor identification of the (abstract) support of $N$ with a subset of $M$, and that this can
be fixed by changing coordinates suitably.

Parkhe proves Theorem~\ref{thm:parkhe-conj} from Theorem~\ref{thm:parkhe-main} by choosing group elements whose $C^{1,\tau}$--norms
are bounded uniformly on every component of their supports, not unlike Jorquera. The diffeomorphisms are built from an equivariant
family which is more or less the same as Yoccoz's equivariant family.

Parkhe's result also has the remarkable consequence that the critical regularity of a finitely generated torsion--free nilpotent group
is always strictly greater than one. While this is a uniform result that holds for all nilpotent groups, the lower bound $1+1/d$ on the critical
regularity of a torsion--free nilpotent group of polynomial growth of degree $d$ is not sharp; this can be seen immediately from the
fact that torsion--free abelian groups act faithfully by $C^{\infty}$ diffeomorphisms, but $\bZ^d$ has growth $\sim n^d$.

Theorem~\ref{thm:parkhe-conj} does not furnish an upper bound on the critical regularity of a torsion--free nilpotent group, though the
Plante--Thurston Theorem (Theorem~\ref{thm:plante-thurston}) gives $2$ as an \emph{a priori} upper bound.
Recall that we write $N_n\le \SL_n(\bZ)$ for the group of unipotent upper triangular integer matrices.
Theorem~\ref{thm:parkhe-conj} is complemented by Theorem A of~\cite{CJN2014}, which shows
that the Farb--Franks construction of a faithful $C^1$ action of $N_n$ on $M$ is not topologically conjugate into $\Diff_+^{1,\tau}(M)$
for $\tau>2/((n-1)(n-2))$, provided $n\geq 4$.

Some better upper bounds
can be given for specific nilpotent groups.

\begin{thm}[~\cite{CJN2014}, Theorem C]\label{thm:cjn}
For all $\tau<1$ and $d\geq 1$, the group $\Diff_+^{1,\tau}(I)$ contains a metabelian group of nilpotence degree $d$.
\end{thm}

Here, a group is \emph{metabelian}\index{metabelian group}
if its commutator subgroup is abelian. Theorem~\ref{thm:cjn} implies that the integral Heisenberg
group $N_3$, which has polynomial growth of degree $4$, has critical regularity exactly two. Theorem~\ref{thm:cjn} does not say
anything about the groups $N_n$ for $n\geq 4$, since these groups are no longer metabelian.

\begin{thm}[Theorems A and B of~\cite{JNR2018}]\label{thm:jnr}
The group $N_4$ has critical regularity exactly $1+1/2$. That is, there is an injective homomorphism \[N_4\longrightarrow \Diff_+^{1,\tau}(I)\]
for all $\tau<1/2$ and no such homomorphism for $\tau>1/2$.
\end{thm}

It is unknown if there is an injective homomorphism of $N_4$ into $\Diff_+^{1,1/2}(I)$, and this is probably a very difficult question to resolve.
The critical regularity of $N_n$ for $n\geq 5$ is also unknown at the time of this book's writing.

\section{The Two-jumps lemma and the $abt$--lemma}\label{sec:abt}

In this section we turn away from constructions of $C^1$ actions of groups, and we
formulate and prove two fundamental dynamical obstructions for certain homeomorphisms (and certain groups
of homeomorphisms) to
be smoothable.

\subsection{The Two-jumps lemma}
The two-jumps lemma was originally formulated by H.~Baik and the authors. It can be viewed as a generalization of an unpublished result
of Bonatti--Crovisier--Wilkinson (\cite[Proposition 4.2.25]{Navas2011}), which concerns crossed $C^1$ diffeomorphisms. The statement is somewhat complicated, but
the result is very useful and the proof is
a fairly straightforward application of the Mean Value Theorem.

\begin{thm}[Two-jumps lemma; see~\cite{BKK2019JEMS}, and~\cite{KK2018JT}, Lemma 2.11]\label{thm:two-jumps}
Let $M$ denote either the compact interval or the circle, and let $f,g\co M\longrightarrow M$ be continuous maps.
Suppose  $(s_i), (t_i)$ and $(y_i)$ are infinite sequences of points in $M$
such that for each $i\ge1$, one of the following two conditions hold:
\be[(i)]
\item
$f(y_i)\le s_i = g(s_i) < y_i < t_i = f(t_i) \le g(y_i)$;
\item
$g(y_i)\le t_i = f(t_i) < y_i < s_i = g(s_i) \le f(y_i)$.
\ee
If $|g(y_i)-f(y_i)|$ converges to $0$ as $i$ goes to infinity,
then $f$ or $g$ fails to be $C^1$.
\end{thm}

Figure~\ref{f:fg} illustrates the case (i) of the Two-jumps Lemma.

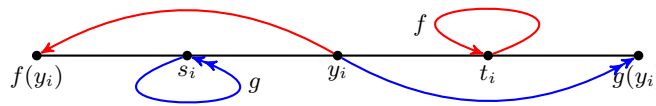
\begin{figure}[h!]
  \tikzstyle {bv}=[black,draw,shape=circle,fill=black,inner sep=1pt]
\centering
\begin{tikzpicture}[>=stealth',auto,node distance=3cm, thick]
\draw  (-4,0) node (1) [bv] {} node [below]  {\small $f(y_i)$} 
-- (-2,0)  node  (2) [bv] {} node [below] {\small $s_i$}
-- (0,0) node (3) [bv] {} node [below] {\small $y_i$} 
-- (2,0)  node (4) [bv] {} node [below] {\small $t_i$}
-- (4,0) node (5) [bv] {} node [below]  {\small $g(y_i)$};
\draw (-1.1,-.4) node {\small $g$};
\draw (1.1,.4) node {\small $f$};
\path (3) edge [->,bend right,red] node  {} (1);
\path (3) edge [->>,bend right,blue] node  {} (5);
\draw [->>, blue] (2)  edge [out = 200,in=-20,looseness=50] (2);
\draw [->, red] (4)  edge [out = 20,in=160,looseness=50] (4);
\end{tikzpicture}%
\caption{Two--jumps Lemma.}
\label{f:fg}
\end{figure}

\begin{proof}[Proof of Theorem~\ref{thm:two-jumps}]
Suppose the contrary, so that $f'$ and $g'$ are both continuous on $M$.
For each index $i$, we let $I_i$ be the closed interval whose endpoints are the points $f(y_i)$ and $g(y_i)$. We set $A_i$ and $B_i$
to be closed intervals characterized by the conditions \[I_i=A_i\cup B_i,\quad A_i\cap B_i= y_i,\quad f(y_i)\in A_i,\, g(y_i)\in B_i.\]
In Figure~\ref{f:fg}, the interval $A_i$ makes up the left half of the full interval and $B_i$ makes up the right half.
The assumptions of the theorem imply that $s_i\in A_i$ and $t_i\in B_i$, though not necessarily in the interiors of these intervals.

Since $M$ is compact, by passing to a subsequence if necessary, we may assume that  $\{y_i\}_{i\ge1}$ converges to a point $y\in M$.
Since by assumption the lengths of the intervals $I_i$ tend to zero as $i$ tends to infinity, we have that $y$ is fixed by both $f$ and $g$.
By Lemma~\ref{lem:diff-fix}, we have that $f'(y)=g'(y)=1$.

The Mean Value Theorem implies the existence of points $u_i\in (s_i,y_i)$ and $v_i\in (y_i, t_i)$ such that
\[g'(u_i)=\frac{g(y_i)-s_i}{y_i-s_i}= 1+\frac{|B_i|}{y_i-s_i}\geq 1+\frac{|B_i|}{|A_i|},\] and such that
\[f'(v_i)=\frac{t_i-f(y_i)}{t_i-y_i}= 1+\frac{|A_i|}{t_i-y_i}\geq 1+\frac{|A_i|}{|B_i|}.\] Note that \[\lim_{i\to\infty} u_i=\lim_{i\to\infty} v_i=y.\]

Multiplying these two expressions together, we see that
\[f'(v_i)g'(u_i)\geq 2+\frac{|B_i|}{|A_i|}+\frac{|A_i|}{|B_i|}\geq 4,\] since an easy calculation shows that
\[\frac{|B_i|}{|A_i|}+\frac{|A_i|}{|B_i|}=\frac{|A_i|^2+|B_i|^2}{|A_i|\cdot |B_i|}\geq 2.\] Lemma~\ref{lem:diff-fix} again shows that
\[\lim_{i\to\infty} f'(v_i)g'(u_i)\longrightarrow 1,\] a contradiction.
\end{proof}

\subsection{The $abt$--lemma}\label{ss:abt}

Suppose we are given three diffeomorphisms of the interval or of the circle, say $\{a,b,t\}$, subject to the proviso that $a$ and $b$ commute.
If these three diffeomorphisms are chosen otherwise in a sufficiently ``generic" manner, it would seem likely that the abstract
 group $\form{ A,B,T}$ they generate
should be one which $A$ and $B$ commute and in which there are no other relations, which is to say $\bZ^2*\bZ$. The main feature of the
$abt$--lemma is that this intuition is completely mistaken, at least in the case where the supports of $a$ and $b$ are disjoint. More
precisely:

\begin{thm}[The $abt$--lemma]\label{thm:abt}
Let $M$ denote the interval $I$ or the circle $S^1$, and suppose $a,b,t\in\Diff^1_+(M)$. Suppose that \[(\supp a)\cap(\supp b)=\varnothing.\]
Then the abstract group $\form{  a,b,t}$ is not isomorphic to $\bZ^2*\bZ$.
\end{thm}

In fact, one can often say significantly more beyond the mere fact that the group $\form{  a,b,t}$ is not isomorphic to $\bZ^2*\bZ$.
Typically, the group $\form{  a,b,t}$ will contain a copy of the lamplighter group $\bZ\wr\bZ$, whence it is clear that no isomorphism
with $\bZ^2*\bZ$ can exist.

Before embarking on a proof of Theorem~\ref{thm:abt}, we first show how it can fail for homeomorphisms.

\begin{prop}\label{prop:z2z-real}
There exist homeomorphisms $a,b,t\in\Homeo_+(M)$ such that \[(\supp a)\cap (\supp b)=\varnothing,\] and such that 
$\form{  a,b,t}\cong \bZ^2*\bZ$.
\end{prop}
\begin{proof}[Sketch]
Since the interval is a two--point compactification of $\bR$ and since $I$ is a submanifold of $S^1$, it suffices to give a proof of the proposition
for $M=\bR$. 
Moreover, for each sequence $\{\ell_n\}_{n\ge1}$ of positive real numbers, we may choose
a sequence $\{I_n\}_{n\ge1}$ of disjoint
intervals with $|I_j|=\ell_j$. Thus, it suffices to find an action of $\bZ^2*\bZ$ on $\bR$ satisfying the hypotheses of the proposition such
that each nontrivial $g\in\bZ^2*\bZ$ acts nontrivially on some interval $I_g$ whose endpoints are global fixed points of $\bZ^2*\bZ$.

Let $1\neq g=T^{n_k}W_k\cdots T^{n_1} W_1$ be a reduced expression for $g\in\bZ^2*\bZ=\form{  A,B,T}$, where each
\[W_i=A^{p_i}B^{q_i}\in \form{  A,B}=\bZ^2.\] For each $i$, choose a generator $X_i\in\{ A,B\}$ such that the exponent $r_i$ of
$X_i$ in $W_i$ is nonzero.

Now, let $J=[0,2k+1]$, let $J_i=[2i-2,2i]$, and let $K_i=[2i-1,2i+1]$
for $1\leq i\leq k$. We now define the action of $\bZ^2*\bZ$. We set $T$ to fix $\partial K_i$
for each $i$, and we set $T$ to act on each $K_i$ by an arbitrary homeomorphism subject to the requirement that \[T^{n_i}(2i-3/4)\geq 2i+3/4.\]
Similarly, we set $X_i$ to act on each $J_i$ by an arbitrary homeomorphism subject to the requirement that \[X_i^{r_i}(2i-7/4)\geq 2i-1/4.\]
We set $Y_i=\{A,B\}\setminus X_i$ to act by the identity on $J_i$. It is straightforward to check now that $g(1/4)\geq 2k+3/4$, so that
$g$ does not act by the identity on $J$.
\end{proof}

It is an informative exercise to generalize Proposition~\ref{prop:z2z-real} to a general right-angled Artin group
$A(\gam)$ (see Section~\ref{sec:raag} below). We sketch here
a few hints to guide the reader. First, consider a word $w$ written in the vertices of the defining graph $\Gamma$
and their inverses. The word $w$ represents the
identity in $A(\gam)$ if and only if it can be reduced to the identity by applying the following moves (see~\cite{cartier-foata,cgw2009,hm1995}):
\begin{itemize}
\item
Free reductions;
\item
Swapping of commuting generators, i.e.~ $v_1v_2\mapsto v_2v_1$ whenever $\{v_1,v_2\}$ is an edge of $\Gamma$.
\end{itemize}
Next, one can express $w$ as a product $w=w_kw_{k-1}\cdots w_1$ of subwords which satisfy the following two conditions:
\begin{itemize}
\item
For $1\leq i\leq k$, any two generators occurring in $w_i$ commute with each other;
\item
For $1\leq i\leq k-1$, given a generator $v_i$ occurring in $w_i$, there is a generator $v_{i+1}$ occurring in $w_{i+1}$ which does not
commute with $v_i$.
\end{itemize}
This last claim can be proved by implementing a sorting algorithm on the generators occurring in $w$. Given this setup, it is fairly
straightforward to prove that $A(\gam)$ acts faithfully by orientation preserving homeomorphisms of $\bR$.

Another variation of Proposition~\ref{prop:z2z-real} is the observation that for all subgroups $G$ and $H$ 
of $\Homeo_+(\bR)$ the free product $G\ast H$ embeds into $\Homeo_+(\bR)$.
This property is not true for $\Diff_+^1(I)$ or $\Diff_+^1(S^1)$, as we will see in Corollary~\ref{c:bs12}.
However, a similar ping-pong argument can show that 
for all subgroups $G$ and $H$ of $\Diff_+^\infty(I)$, the free product $G\ast H$ admits an 
\emph{eventually injective}\index{eventually injective} sequence of homomorphisms into $\Diff_+^\infty(I)$ in the following sense. 
One can also describe this property as that $G\ast H$ is \emph{residually $\Diff_+^\infty(I)$.}\index{residual property}

\begin{lem}[cf. \cite{BKK2014,KKFreeProd2017}]\label{l:cinfty-free-product}
If $G_0$ and $G_1$ are subgroups of $\Diff_+^\infty(I)$,
then for each nontrivial $g\in G_0\ast G_1$ 
 there exists a representation
\[
\phi=\phi_g\co G_0\ast G_1\longrightarrow \Diff_+^\infty(I)\]
having a connected support such that $\phi(g)\ne1$. Furthermore, one can require that
$\phi$ is injective on $G_0$ and on $G_1$.
\end{lem}
The proof will follow the classical  ping-pong argument of Klein~\cite{Klein1883MA}, Maskit~\cite{Maskit1988Springer} and
Tits~\cite{tits-jalg}. Indeed, the nontriviality of $g$ for some new action will be certified by exhibiting a point that is  moved
consecutively by factors of $g$. We will pay extra care so that the support of this new action is connected. 

The following notation will be used for the remainder book.
\begin{notation}\label{notation:diff-J}
For an interval $J$ in $\bR$, we define
\[
\Homeo_J(\bR):=\{f\in \Homeo_+(\bR)\mid \supp f\sse J\}.\]
We similarly define 
$\Diff_J^{k}(\bR)$ and $\Diff_J^{k,\alpha}(\bR)$
for an integer $k\ge1$ and for a concave modulus $\alpha$.
\end{notation}
If $J$ is bounded, one may identify each element of $\Diff_J^{k}(\bR)$ with a $C^{k}$ 
diffeomorphism of $J$ that is $C^k$--tangent to the identity at $\partial J$.

\bp[Proof of Lemma~\ref{l:cinfty-free-product}]
We let $I=[0,1]$. By Muller--Tsuboi trick (Theorem~\ref{t:muller-tsuboi}), we may assume that 
\[ G_0, G_1\le\Diff_I^\infty(\bR).\] 
It will be convenient for us to have a notation for the modulo--two remainder:
\[p(i)=i-2\floor{i/2}\]
 The conclusion for the case $g\in {G_0}\cup {G_1}$ is easy to prove, so we will assume to have
\[
g =v_{2k-1}\cdots v_1v_0\]
for some integer $k\ge1$ and for some nontrivial elements \[v_i\in G_{p(i)}\le\Diff_I^\infty(\bR).\]
Let $J_i$ be the closure of a connected component of $\supp G_{p(i)}$ such that $v_i\restriction_{J_i}$ is nontrivial.
Using  Muller--Tsuboi trick again, we can find a $C^0$--conjugate $H_{p(i)}\le\Diff_{ J_i}^\infty(\bR)$ of $G_{p(i)}\restriction_{J_i}$.

Note that the inversion $\sigma(x)=1-x$ of $I$ conjugates $\Diff_I^\infty(\bR)$ to itself.
We can construct an action
\[
\phi\co G_0\ast G_1\longrightarrow \Diff_{[0,2k+1]}^\infty(\bR)\]
satisfying the following.
\begin{itemize}
\item $\supp\phi(G_0)=\bigcup_{i=0}^{k-1} (2i,2i+2) \cup (2k+3/5,2k+4/5)$;
\item $\supp\phi(G_1)=\bigcup_{i=1}^{k} (2i-1,2i+1) \cup (1/5,2/5)$;
\item $\phi(G_0)\!\!\!\restriction_{[2k+3/5,2k+4/5]}$ is $C^\infty$--conjugate to $G_0\!\!\!\restriction_I$;
\item $\phi(G_1)\!\!\!\restriction_{[1/5,2/5]}$ is $C^\infty$--conjugate to $G_1\!\!\!\restriction_I$; 
\item $\phi(G_i)\!\!\!\restriction_{[i,i+2]}$ is $C^\infty$--conjugate to $H_{p(i)}\!\!\!\restriction_{J_i}$
or to $\left(\sigma H_{p(i)} \sigma\right)\!\!\!\restriction_{\sigma ( J_{i})}$
for $i=0,\ldots,2k-1$;
\item $\phi(v_i)(i+1/2)=i+3/2$ for $i=0,\ldots,2k-1$.
\end{itemize}
In particular, we see that $\phi(g)(1/2)=2k+1/2$, and the conclusion follows.
\ep

\subsubsection{Compact supports}\label{ss:compact supports}

Let $H\le \Homeo_+[0,1]$. Recall that $\supp H$ consists of $x\in M$ for which there exists a $h\in H$ with $h.x\neq x$, which is an open
subset of $(0,1)$. 
If $1\neq g\in H$, we say that $g$ is \emph{compactly supported (in $H$)}\index{compactly supported homeomorphism},
or has \emph{compact support}, if
$\overline{\supp g}$ is contained in $\supp H$. That is, $\supp g$ is a precompact subset of $\supp H$.

Elements with compact support are useful for producing ``unexpected" commutations among diffeomorphisms. This idea goes back
to the Zassenhaus Lemma (see~\cite{Raghunathan1972}),
and was crucial in Brin--Squier's proof that the group of piecewise linear homeomorphisms of the interval
does not contain a nonabelian free group (see Theorem~\ref{thm:f-subgp} below, where we will give a proof of the Brin--Squier Theorem).
 This idea will also play an important role in the construction of finitely generated groups of diffeomorphisms in
 Chapter~\ref{ch:optimal}, especially in Theorem~\ref{t:optimal-group}.

As a warmup, we consider the following.

\begin{prop}[See~\cite{KK2018JT}, Lemma 3.3]\label{prop:compact-conn}
Suppose $H\le \Homeo_+[0,1]$ with $\supp H$ connected, and suppose that $g\in H$ has compact support. Then there exists a nontrivial element
$t\in H$ such that $[g,tgt^{-1}]=1$.
\end{prop}
\begin{proof}
Since $\supp H$ is connected, there is a collection $\{J_1,\ldots,J_k\}$ of intervals such that:
\begin{enumerate}[(i)]
\item
We have \[\overline{\supp g}\sse\bigcup_{i=1}^k J_i.\]
\item
We have $J_i\cap J_{i+2}=\varnothing$ for $i<k-1$.
\item
We have $\inf J_i<\inf J_{i+1}<\sup J_i<\sup J_{i+1}$ for $i\le k-1$.
\item
For all $i$, the interval $J_i$ is a component of the support of an element $h_i\in H$.
\end{enumerate}
A collection satisfying the conditions (ii) and (iii) above will be called as a $k$--chain; see Definition~\ref{defn:chain}.

Write $K=\overline{\supp g}$. We may clearly assume that $z_1=\inf K\in J_1$ and $\sup K\in J_k$. There are
integral exponents $\{n_i\}_{1\leq i\leq k}$
such that:
\begin{itemize}
\item
We have $z_2=h_1^{n_1}(z_1)\in J_2$.
\item
Having inductively defined $z_i$, we have $z_{i+1}=h_i^{n_i}(z_i)\in J_{i+1}$ for all $i<k$.
\item
We have $h_k^{n_k}(z_k)\geq\sup K$.
\end{itemize}

We then set $t=h_k^{n_k}\cdots h_1^{n_1}$. It is easy to see then that $\supp g$ and $\supp tgt^{-1}$ are disjoint, whence the conclusion
of the proposition.
\end{proof}

In fact, Proposition~\ref{prop:compact-conn} implies that $H$ contains the \emph{lamplighter group}\index{lamplighter group}
\[
\bZ\wr\bZ= \left(\bigoplus_{\bZ}\bZ\right)\rtimes\bZ.\]
Recall that this is the
\emph{wreath product}\index{wreath product} of $\bZ$ with $\bZ$. 
This is built by taking \[N:=\bigoplus_{\bZ}\bZ,\] and letting $\bZ$ act on this direct sum
by shifting the index (which is to say via the natural action of $\bZ$ on itself), thus forming a semidirect product of $\bZ$ and $N$.
Letting $\lambda$ be a generator of the copy of $\bZ$ acting on $N$ and letting $\tau$ be a generator of the summand of $N$ corresponding to
$0\in\bZ$, it is easy to see the following group presentation:
\[\bZ\wr\bZ=\form{\tau,\lambda\mid \left[\tau,\lambda^i\tau \lambda^{-i}\right]=1\text{ for all }i\in\bZ}.\]
This wreath product is a finitely generated, nonabelian, solvable group
(in fact, it is a metabelian group), that is not finitely presentable. All but the last of these claims is an easy consequence of the definitions
save the last, which we will not discuss further here. 

Lamplighter groups occur 
quite commonly in the group of homeomorphisms of the interval.
Indeed, let $g\in\Homeo_+(\bR)$ be a homeomorphism such that $\supp g\sse (0,1)$, and let $f$ be translation by one. Then $f$ and $g$
are easily seen to generate a subgroup of $\Homeo_+(\bR)\cong \Homeo_+[0,1]$ that is isomorphic to $\bZ\wr\bZ$. Under the hypotheses
of Proposition~\ref{prop:compact-conn}, it is clear that the elements $g,t\in H$ furnished by the proposition generate a lamplighter group.

Proposition~\ref{prop:compact-conn} has the following corollary for homeomorphisms of the circle.

\begin{prop}\label{prop:comm-circle}
Let \[\{a,b,c,d\}\sse\Homeo_+(S^1)\] be nontrivial elements such that \[(\supp a)\cap(\supp b)=\varnothing\quad \textrm{and}\quad
(\supp c)\cap(\supp d)=\varnothing.\] Writing $G=\form{  a,b,c,d}$, if $\supp G= S^1$ then $G$ contains a copy of $\bZ\wr\bZ$.
\end{prop}

We leave the proof of Proposition~\ref{prop:comm-circle} as an exercise for the reader. To proceed with the proof of Theorem~\ref{thm:abt},
we need to state and prove some rather opaque facts about supports of commutators of homeomorphisms. It is difficult to give a
satisfactory intuitive explanation of the estimates we give; the reader should perhaps just keep in mind that we are merely aiming to force
a configuration of intervals under the hypotheses of Theorem~\ref{thm:abt} wherein the Two-jumps Lemma (Theorem~\ref{thm:two-jumps})
becomes applicable.

The reader may keep the following picture in their mind: suppose $\{a,b,t\}$ are as in the statement of Theorem~\ref{thm:abt}, with the
supports of $a$ and $b$ disjoint. It is not difficult to see that for $\form{  a,b,t}$ to be isomorphic to $\bZ*\bZ^2$, then both
$\supp a$ and $\supp b$ must have infinitely many components. Now, let $c=t^{-1}at$ and $d=t^{-1}bt$. If
$\form{  a,b,t}\cong\bZ*\bZ^2$ then $\form{  a,b,c,d}\cong\bZ^2*\bZ^2$. For the action of $\bZ^2*\bZ^2$ on $I$ via $\{a,b,c,d\}$
to have any hope of being faithful, then on infinitely many components of support, the diffeomorphisms $\{a,b,c,d\}$ have to translate points
a definite fraction of the length of their supports. This is a violation of the Mean Value Theorem (which is implicit in the Two-jumps Lemma).

The first technical result is the following.

\begin{lem}[\cite{KK2018JT}, Lemma 3.5]\label{lem:supp-comm}
Let $X$ be a Hausdorff topological space, and let $f,g\in\Homeo(X)$. Then, we have \[\overline{\supp [f,g]}\sse\supp f\cup\supp g\cup
\overline{\supp f\cap\supp g}.\]
\end{lem}
\begin{proof}
Suppose first that \[x\notin \supp f\cup\supp g\cup
\overline{\supp f\cap\supp g}.\] Then \[f(x)=g(x)=x,\] and so $x\notin \supp [f,g]$. It suffices to show that there is a neighborhood of $x$
that is in the fixed point set of $[f,g]$.

There is an open neighborhood $U$ containing $x$ such that \[U\cap (\supp f\cap\supp g)=\varnothing,\] since $x$ is not in the
closure of $\supp f\cap\supp g$. There is a further open sub-neighborhood $V$ of $x$ for which we have \[f^{\pm1}(V)\cup g^{\pm 1}(V)
\sse U.\] We claim that $V$ is fixed by $[f,g]$. If $y\in V$, then it is easy to check that $[f,g](y)=y$ by checking the cases
\[y\in V\cap\supp f,\quad y\in V\cap\supp g,\quad y\in V\cap \Fix f\cap\Fix g.\] This implies that \[x\notin\overline{\supp[f,g]},\] which
proves the lemma.
\end{proof}

The next lemma is a crucial estimate, which will allow us to finally apply the Two-jumps Lemma. It is established by a brute force
calculation.

\begin{lem}[\cite{KK2018JT}, Lemma 3.6]\label{lem:supp-phi}
Let $b,c,d$ be bijections of a set $X$ such that
\[\supp c\cap\supp d=\varnothing,\]
and write
\[\phi= [c,bdb^{-1}].\]
Then, we have
\[\supp\phi\sse\supp b\cup cb\left(\supp b\cap\supp d\right)
\cup db^{-1}\left(\supp b\cap\supp c\right).\] 
\end{lem}
\begin{proof}
We will follow the notation in~\cite{KK2018JT} and write $\yt g=\supp g$ for $g\in\{b,c,d\}$. Note that
\[[c,bdb^{-1}]=cbd(cb)^{-1}\cdot bd^{-1}b^{-1}=c\cdot b\cdot db^{-1}c^{-1}(db^{-1})^{-1}\cdot b^{-1}.\]

We have that \[\supp\phi\sse  \left(\yt{c}\cup b(\yt d)\right)\cap
\left(cb(\yt d)\cup b(\yt d)\right)\cap\left(\yt c\cup\yt b\cup db^{-1}(\yt c)\right).\]
Some tedious calculations, using the facts that \[b(\yt d)\sse  \yt b\cup\yt d\quad
\textrm{and}\quad \yt c\cap\yt d=\varnothing,\] all taken together imply that
\[ \left(\yt{c}\cup b(\yt d)\right)\cap \left(cb(\yt d)\cup b(\yt d)\right)\cap\left(\yt c\cup\yt b\cup db^{-1}(\yt c)\right)\sse 
\left(\yt c\cap cb(\yt d)\right)\cup\yt b\cup\left(\yt d\cap db^{-1}(\yt c)\right).\] We leave the details of these calculations
as an exercise for the reader.

Observe that it now suffices to prove the inclusions \[\yt c\cap \left(cb(\yt d)\right)\sse  cb\left(\yt b\cap\yt d\right),\] and
\[\yt d\cap \left(db^{-1}(\yt c)\right)\sse  db^{-1}\left(\yt b\cap\yt c\right).\] We will prove the first inclusion, as the second
will follow by the symmetry of $c$ and $d$.

So, let $x\in X$, and suppose that \[cb(x)\in\yt c\cap cb(\yt d).\] It follows then the $x\in\yt d$ and $cb(x)\in\yt c$. We have
that $x\notin\Fix b$, since $\yt c\cap\yt d=\varnothing$ and since $b(x)\in\yt c$. It follows that $x\in\yt b$, and so \[cb(c)\in cb(\yt b\cap\yt d).\]
This proves the desired inclusion, and with it, the lemma.
\end{proof}

The foregoing lemmata were purely combinatorial, and made no reference to diffeomorphisms. We will now use the
differentiability hypothesis to obtain a diffeomorphism whose algebraic provenance forces it to be compactly supported. We first see
that $\phi$ as in Lemma~\ref{lem:supp-phi}, if built on the interval $I$ using $C^1$
diffeomorphisms, is compactly supported ``modulo $\supp b$".

\begin{lem}[\cite{KK2018JT}, Lemma 3.7]\label{lem:diff-phi-b}
If we have orientation--preserving $C^1$ diffeomorphisms $b,c,d$ of a compact interval $I$ such that
\[\supp c\cap\supp d=\varnothing,\]
and if we write
\[\phi= [c,bdb^{-1}],\]
then we have that \[\overline{\supp\phi\setminus\supp b}\sse \supp c\cup\supp d.\]
\end{lem}
\begin{proof}
For the purposes of this proof, we will write $\BB$ for the set of components of $\supp b$ and $\CC$ for the set of components of $\supp c$.
For $B\in\BB$, we will write \[J_B=B\cup cb(B\cap\supp d)\cup db^{-1}(B\cap\supp c).\] Note that
\[\supp\phi\sse \bigcup \{J_B\mid B\in\BB\}=\bigcup\{J_B\setminus B\mid B\in\BB\}\cup\supp b,\] as follows from
Lemma~\ref{lem:supp-phi}. Observe also that \[\overline{J_B\setminus B}\sse  \overline{c(B)\setminus B}\cup\overline{d(B)\setminus B}
\sse \supp c\cup\supp d.\]

{\bf Claim: the collection \[\BB_0=\{B\in\BB\mid J_B\neq B\}\] consists of finitely many intervals.} This claim is where we will use the
differentiability hypothesis.

The conclusion of the Lemma follows from this claim. Indeed, we have \[\overline{\supp\phi\setminus\supp b}\sse \overline{\bigcup\{
J_B\setminus B\mid B\in\BB_0\}},\] as follows formally. If $\BB_0$ is a finite union of intervals then closure commutes with unions, and
so we get \[\overline{\supp\phi\setminus\supp b}\sse \bigcup\{\overline{J_B\setminus B}\mid B\in\BB_0\}\sse \supp c\cup \supp d.\]

To prove the claim,
we set \[\BB_1=\{B\in\BB\mid cb(B\cap\supp d)\setminus B\neq\varnothing\},\] and we write \[\BB_2=\{B\in\BB\mid db^{-1}(B\cap\supp c)
\setminus B\neq\varnothing\}.\] Observe that $\BB_0=\BB_1\cup\BB_2$, and suppose that this collection is infinite. We will suppose
that $\BB_1$ is infinite and derive a contradiction, the argument for $\BB_2$ being essentially the same.

Thus, we assume that \[\{B_i\}_{i\ge1}\sse \BB_1\] are distinct, and for all $i\ge1$ there are points $x_i\in B_i\cap \supp d$ such
that $cb(x_i)\notin B_i$. 
In particular, we have that $b(x_i)\in\supp c$ for $i\ge1$. We can find elements \[\{C_i\}_{i\ge1}\sse \CC\] such that for all $i$, we have \[b(x_i), cb(x_i)\in C_i.\]

Let $J_i=[x_i, cb(x_i)]$, possibly with the endpoints switched. Note that $b(x_i)$ lies in the interior of $J_i$. It follows that the interval
$(b(x_i),cb(x_i)]$ meets $\partial B_i$, and that the interval $[x_i,b(x_i))$ meets $\partial C_i$.

The claim now follows from the Two-jumps Lemma (Theorem~\ref{thm:two-jumps}), setting \[f=b^{-1},\quad g=c,\quad s_i=\partial C_i\cap B_i,
\quad t_i=\partial B_i\cap C_i,\quad y_i=b(x_i).\] This completes the proof of the lemma.
\end{proof}

As we have suggested before, the following lemma will conclude the proof of Theorem~\ref{thm:abt}.

\begin{lem}[\cite{KK2018JT}, Lemma 3.8]\label{lem:z2-z2}
Let $M\in\{I,S^1\}$, and let \[\{a,b,c,d\}\sse \Diff_+^1(M).\] Suppose that \[\supp a\cap\supp b=\supp c\cap\supp d=\varnothing.\]
Then the group $\form{  a,b,c,d}$ is not isomorphic to the group $\bZ^2*\bZ^2$.
\end{lem}

Indeed, let $\{a,b,t\}\sse \Diff_+^1(M)$ satisfy the hypotheses of Theorem~\ref{thm:abt}, and suppose that \[\form{  a,b,t}\cong\bZ*
\bZ^2.\] Then we have \[\form{  a,b,t^{-1}at,t^{-1}bt}\cong\bZ^2*\bZ^2.\] Setting $c=t^{-1}at$ and $d=t^{-1}bt$, we have that
\[\supp c\cap\supp d=\varnothing.\] It follows that the hypotheses of Lemma~\ref{lem:z2-z2} are satisfied, but the conclusion
is contradicted. It remains therefore to establish Lemma~\ref{lem:z2-z2}.

Before establishing Lemma~\ref{lem:z2-z2}, we record two basic facts about $\bZ^2*\bZ$, which we will prove fully here for the
purposes of self-containment.

\begin{prop}\label{prop:2g-hopf}
The following are properties of $\bZ^2*\bZ^2$.
\begin{enumerate}[(1)]
\item
If $g,h\in\bZ^2*\bZ^2$ then $\form{  g,h}$ is either free or abelian. In particular, $\bZ^2*\bZ$ contains no subgroup isomorphic to
$\bZ\wr\bZ$.
\item
Every surjective endomorphism of $\bZ^2*\bZ^2$ is an isomorphism; that is, $\bZ^2*\bZ^2$ is \emph{Hopfian}\index{Hopfian group}.
\end{enumerate}
\end{prop}

Both of the conclusions of Proposition~\ref{prop:2g-hopf} hold for general right-angled Artin groups, though we will not need such
generality here. See~\cite{Baudisch1981,KK2015GT,LM2010GD} for a discussion of two--generated subgroups,
and~\cite{dlHarpe2000,KS2020rfrp} for a
discussion of Hopficity and residual finiteness. The proof of Proposition~\ref{prop:2g-hopf} draws on ideas in~\cite{dlHarpe2000,KS2020rfrp}
especially.

\begin{proof}[Proof of Proposition~\ref{prop:2g-hopf}]
For the first item, we use a special case of
the Kurosh Subgroup Theorem, which says that a finitely generated subgroup of a free product $G*H$ of two groups
$G$ and $H$ is of the form \[K_1*K_2\cdots*K_n* F,\] where each $K_i$ is isomorphic to a finitely generated subgroup of $G$ or of $H$, and
where $F$ is a finitely generated free group. The Kurosh Subgroup Theorem itself follows easily from the covering space theory of
a wedge (i.e.~one--point union) of an arbitrary collection of connected cell complexes.
Thus, if \[K=\form{  g,f}\le \bZ^2*\bZ^2\] then if $K$ is itself not cyclic or trivial then either
$K\cong\bZ^2$, or $K$ is a free group of rank two.

For the second item, it suffices to show that $\bZ^2*\bZ^2$ is \emph{residually finite}\index{residually finite},
which is to say that every nontrivial element
$1\neq g\in\bZ^2*\bZ^2$ does not lie in the kernel of a homomorphism \[\phi_g\colon \bZ^2*\bZ^2\longrightarrow F,\] where
$F$ is a finite group. Indeed, suppose that \[\phi\colon \bZ^2*\bZ^2\longrightarrow \bZ^2*\bZ^2\] is a surjection with kernel $K$.
Since $\bZ^2*\bZ^2$ is finitely generated, there are only finitely many subgroups of $\bZ^2*\bZ^2$ of finite index $n$, say
$\{G_1,\ldots,G_m\}$. Then the subgroups \[\{\phi^{-1}(G_i)\mid 1\leq i\leq m\}\] coincide with the groups $\{G_1,\ldots,G_m\}$.
It follows that $K$ lies in the intersection of all finite index subgroups of $\bZ^2*\bZ^2$. So, if $\bZ^2*\bZ^2$ is residually finite then $K$
is trivial.

To see that $\bZ^2*\bZ^2$ is residually finite, we endow the usual $2$--dimensional torus $T$ with the standard flat metric coming from the
square. We set $X=X_0$ to be the one--point union of two copies of $T$, so that $\pi_1(X)\cong \bZ^2*\bZ^2$. We write
$p$ for the wedge point in $X$. We extend the metric
on $T$ to $X$ in the obvious (i.e.~the $\ell^1$) way. We define a tower
of covers $\{X_i\}_{i\ge1}$, where $X_{i+1}$ is the cover of $X_i$ classified by the surjection \[\pi_1(X_i)\longrightarrow
H_1(X_i,\bZ/2\bZ).\] 

Consider a shortest--length homotopically nontrivial loop $\gamma_i$ on $X_i$. If $\gamma_i$ lies entirely in a component of the
preimage of one of the two tori $T$ making up $X$, then clearly the homotopy class $[\gamma_i]$ is nontrivial in $H_1(T,\bZ)$
and is hence nontrivial in $H_1(X_i,\bZ)$. Since $\gamma_i$ is shortest, this homotopy class will obviously be primitive, and hence
will be nontrivial after reducing modulo two. It follows then that $\gamma_i$ does not lift to $X_{i+1}$.

If $\gamma_i$ does not lie in a component of the preimage of one of the two tori, then $\gamma_i$ traverses at least one component
of the preimage of $p$ in $X_i$. Since $\gamma_i$ is a shortest curve on $X_i$, there must be one such component of the preimage of $p$,
say $q$,
which is traversed exactly once. Traversing $q$
(with orientation) represents a nontrivial integral cohomology class, and counting the number of times a loop traverses $q$ modulo two
represents a nontrivial $\bZ/2\bZ$--valued cohomology class. It thus follows that $\gamma_i$
represents a nontrivial homology class in $H_1(X_i,\bZ/2\bZ)$. Again it follows that $\gamma_i$ does not lift to $X_{i+1}$.
Thus, an easy induction shows that the shortest loop on $X_i$ has length at least $i+1$. This proves that the intersection of all
finite index subgroups of $\pi_1(X)$ is trivial, whence $\pi_1(X)$ is residually finite.
\end{proof}

We now complete the proof of Lemma~\ref{lem:z2-z2} and hence of Theorem~\ref{thm:abt}.

\begin{proof}[Proof of Lemma~\ref{lem:z2-z2}]
We write \[\bZ^2*\bZ^2\cong \form{  w,x,y,z\mid [w,x]=[y,z]=1}.\] We define a surjection
\[\tau\colon\bZ^2*\bZ^2\longrightarrow\Diff_+^1(M)\]
by the rules \[w\mapsto a,\quad x\mapsto b,\quad y\mapsto c,\quad z\mapsto d.\] By Proposition~\ref{prop:2g-hopf}, it
suffices to show that $\form{  a,b,c,d}$ is not abstractly isomorphic to $\bZ^2*\bZ^2$. We assume the contrary,
namely that $\tau$ is an isomorphism, and so derive
a contradiction.

Since $\bZ^2*\bZ^2$ contains no copy of $\bZ\wr\bZ$, we may assume that $M\neq S^1$ by Proposition~\ref{prop:comm-circle}.
So, we will assume that $M=I$. Write \[\phi=[c,bdb^{-1}],\quad \psi=[\phi,a].\]
Applying Lemma~\ref{lem:supp-comm}, we see that \[\overline{\supp\psi}\sse \supp\phi\cup \supp a\cup\overline{\supp \phi\cap
\supp a}.\]

Lemma~\ref{lem:diff-phi-b} now implies that \[\overline{\supp\phi\cap\supp a}\sse \overline{\supp\phi\setminus\supp b}\sse 
\supp c\cup\supp d.\] We conclude that \[\overline{\supp\psi}\sse \supp G,\] so that $\psi$ is compact supported in
$\supp G$. Finally, we have that $\psi$ is the image
under $\tau$ of the element \[[[y,xzx^{-1}],w]\in\bZ^2*\bZ^2.\] This element is clearly seen to
be nontrivial in $\bZ^2*\bZ^2$, and so we have that if $\tau$ is an isomorphism then $\psi$ is nontrivial in $\Diff_+^1(I)$.
Proposition~\ref{prop:compact-conn} implies that $\bZ^2*\bZ^2$ contains a copy of $\bZ\wr\bZ$, which is a contradiction.
\end{proof}

The following is an immediate consequence from the proof of Lemma~\ref{lem:z2-z2}, after setting
 $c=a^t$ and $d=b^t$.
 We will employ this result crucially in the construction of groups of \emph{optimally expanding diffeomorphisms}\index{optimally expading
 diffeomorphism} in Chapter~\ref{ch:optimal}.
\begin{lem}\label{l:ucsd}
If $a,b,t$ are orientation--preserving $C^1$ diffeomorphisms of a compact interval $I$
and if $a$ and $b$ have disjoint supports,
then the closure of the support of the element
\[
\left[ [a^t,b b^t b^{-1}],a\right]\]
is contained in the support of the group $\form{a,b,t}$.
\end{lem}

\section{Crossed homeomorphisms and commutation}\label{sec:crossed}

One of the main difficulties in using Theorem~\ref{thm:abt} is in ensuring that the hypotheses are met. Typically, one does not
expect two homeomorphisms of a on--manifold to have disjoint supports, even if they commute with each other. 
We briefly recall
a standard construction: let $\tau\in\Homeo_+(\bR)$ be given by $\tau(x)=x+1$, and let $\sigma_0$ be an arbitrary homeomorphism
supported on $(0,1)$. The homeomorphism \[\sigma=\prod_{i\in\bZ} \tau^{-i}\sigma_0\tau^i\] commutes with $\tau$ by construction,
whereas the (closures of the) supports of both $\tau$ and $\sigma$ are the whole real line. To give applications of Theorem~\ref{thm:abt}.
we will give specific conditions under which the hypotheses will be satisfied; see Section~\ref{sec:misc} below, and especially
Section~\ref{sec:raag} below, where we will give a complete description of right-angled Artin groups that admit
faithful actions on $I$ and $S^1$ of regularity $C^{1+\mathrm{bv}}$ and better.

We first note an easy and almost immediate consequence of Theorem~\ref{thm:abt}. If $G\le \Homeo_+[0,1]$, we say that
the action of $G$ is \emph{overlapping}\index{overlapping action}
if for arbitrary nontrivial elements $g_1,g_2\in G$, we have \[\supp g_1\cap\supp g_2
\neq\varnothing.\]

\begin{prop}\label{prop:overlapping}
Let $G\le \Diff_+^1(I)$ be a subgroup such that $G*\bZ\le \Diff_+^1(I)$. Then the action of $G$ is overlapping.
\end{prop}

Overlapping actions of groups are closely related to crossed elements, Conradian actions, and Conradian orderability of groups
(cf.~Appendix~\ref{ch:append2} below).

Following Navas and Rivas~\cite{navasrivas09}, let $(\Omega,<)$ be a totally ordered space, let $G$ be a group of order preserving
permutations of $\Omega$, and let $f,g\in G$. We say that $f$ and $g$ are \emph{crossed}\index{crossed homeomorphisms}
if there are points
$u<w<v$ in $\Omega$ such that:
\begin{enumerate}[(1)]
\item
We have $g^n(u)<w<f^n(v)$ for all $n\in\bZ$.
\item
There is an $N\in\bZ$ such that $g^N(v)<w<f^N(u)$.
\end{enumerate}

If a group $G$ acts on $\Omega$ with no crossed pairs of elements, then the action of $G$ is said to be 
\emph{Conradian}\index{Conradian ordering}.

A pair of homeomorphisms of $I$ are crossed if they are crossed as order preserving permutations of the interval.
The reader will find many examples of groups generated by crossed homeomorphisms in Chapter~\ref{ch:chain-groups} below. It is
a useful exercise for the reader to draw several examples of crossed homeomorphisms of the interval, in order to develop an intuition
for the concept.

Conradian subgroups $G\le \Homeo_+[0,1]$ are generally easier to investigate than arbitrary subgroups. For one, it is easy to show that
if $G$ is Conradian and if $f,g\in G$ then if $J_x$ is a component of $\supp x$ for $x\in \{f,g\}$, then either $J_f\cap J_g=\varnothing$,
or there is an inclusion relation between $J_f$ and $J_g$.

The existence of crossed diffeomorphisms in a subgroup $G\le \Diff_+^1(I)$ also has some remarkable consequences. If $h\in\Diff_+^1(I)$
and if $p\in\Fix h$, we say that $p$ is a \emph{hyperbolic}\index{hyperbolic fixed point}
fixed point for $h$ if $|h'(p)|\neq 1$. The following fact can be found
in~\cite{DKN2007,BF2015}.

\begin{prop}\label{prop:hyp-fp}
Let $f,g\in\Diff_+^1(I)$ be crossed diffeomorphisms. Then there exists an element $h$ in the positive semigroup generated by $f$ and $g$
that has a hyperbolic fixed point.
\end{prop}

It is not difficult to imagine how a proof of Proposition~\ref{prop:hyp-fp} might go; one would use the definition of crossed homeomorphisms
and the contraction mapping principle to find a fixed point for an element in the positive semigroup generated by $f$ and $g$. One then
argues from some sort of uniform contraction and the $C^1$ hypothesis that the derivative at the fixed point must be either
strictly greater than one or strictly less than one. We omit further details.

To illustrate more precisely how crossed diffeomorphisms figure into the algebraic theory of diffeomorphisms, we have the following
fact:

\begin{prop}[See~\cite{Navas2011}, Proposition 4.2.2.25, and~\cite{KKR2020}, Lemma 3.10]\label{prop:crossed-comm}
Let $c\in\Diff_+^1(I)$ be a diffeomorphism such that $\supp c=(0,1)$. Then the centralizer of $c$ in $\Diff_+^1(I)$ has no pairs of
crossed diffeomorphisms.
\end{prop}
\begin{proof}
Suppose the contrary, so that $f$ and $g$ centralize $c$ and are crossed. Then the supports and dynamics of $f$ and $g$ are $c$--invariant,
and it is easy to check that then the hypotheses of Theorem~\ref{thm:two-jumps} are satisfied. It follows that either $f$ or $g$ fails to be
$C^1$.
\end{proof}

Often, quite a bit of manipulation is needed to show that a particular faithful group action on $I$ is forced to be Conradian, and in the interest
of space and focus we will not spell out many more details, directing the interested reader to~\cite{KKR2020} and the further
references therein, for instance. As we have already mentioned, once an action is known to be Conradian, one can often analyze the
combinatorics of orbits in order to prove that a purported faithful action cannot be faithful after all.

We close this section by stating the main results of~\cite{KKR2020} and~\cite{KKR2021},
a special case of which we will revisit as Theorem~\ref{thm:kkr-2020}
below. We recall quickly for the convenience of the reader that the \emph{derived series}\index{derived series}
$\left\{G^{(n)}\right\}_{n\in\bN}$ of $G$ is defined by $G^{(0)}=G$
and \[G^{(n+1)}=[G^{(n)},G^{(n)}].\] The group $G$ is \emph{not solvable of degree at most $n$}\index{non-solvable group}
if $G^{(n)}\neq \{1\}$.

\begin{thm}[See~\cite{KKR2020}, Theorem 1.1 and~\cite{KKR2021}, Theorem 3.1]\label{thm:kkr2020-gen}
Let $G$ and $H$ be groups.
\begin{enumerate}[(1)]
\item
Suppose $G$ and $H$ are not solvable of degree at most $n\geq 3$, let $M\in\{I,S^1\}$,
and suppose that $\tau\in [0,1)$ satisfied \[\tau(1+\tau)^{n-2}\geq 1.\]
Then there is no faithful homomorphism \[(G\times H)*\bZ\longrightarrow\Diff_+^{1,\tau}(M).\]
\item
If $G$ and $H$ are not solvable, then for all $\tau>0$ there is no faithful homomorphism
\[(G\times H)*\bZ\longrightarrow\Diff_+^{1,\tau}(M).\]
\end{enumerate}
\end{thm}

In Theorem~\ref{thm:kkr2020-gen}, the symmetry in the non--solvability of $G$ and $H$ is necessary. We do not know how to make
such nontrivial statements about $C^{1,\tau}$ actions of groups of the form $G\times\bZ$, for example, when $G$ has a relatively
long derived series. It is this technical difficulty which makes a full computation of the critical regularity of all right-angled Artin groups
out of reach with the present state of the available technology.

\section{Groups of $C^{1+\mathrm{bv}}$ diffeomorphisms}\label{sec:misc}

Here we record some further useful facts about elements of $\Diffb(S^1)$ which synthesize the foregoing discussion, and which
will be useful for us in the sequel. The following is an easy consequence of Kopell's Lemma (Theorem~\ref{thm:kopell}).

\begin{prop}[See Lemma 2.6 in~\cite{KK2018JT}, for instance]\label{prop:disj-ab}
The following hold.
\begin{enumerate}[(1)]
\item
Let $M\in\{I,S^1\}$ and let $a,b\in\Diffb(M)$ be commuting elements such that \[\Fix a,\Fix b\neq\varnothing.\]
Write $J_a$ and $J_b$ for components of $\supp a$ and $\supp b$
respectively. Then either $J_a=J_b$ or $J_a\cap J_b=\varnothing$.
\item
Let $a,b,c\in\Diffb(I)$ with $\Fix(a)\cap (0,1)=\varnothing$. Then commutation is transitive in the following sense:
if $a$ commutes with $b$ and $a$ commutes with $c$ then $b$ commutes with $c$.
\end{enumerate}
\end{prop}
\begin{proof}
(1) Suppose first that $M=I$. Suppose that $J_a\cap J_b=\varnothing$. It is easy to see that there must be an inclusion
relation among $J_a$ and $J_b$, say $J_a\sse  J_b$. Since
that $\partial J_b$ does not meet $\supp a$, we have that both $a$ and $b$ act on $J_b$, and $a$ acts with a fixed point in the interior
of $J_b$. Kopell's Lemma (Theorem~\ref{thm:kopell}) implies that $a$ must restrict to the identity on $J_b$, a contradiction.
If $M=S^1$, let $J_b$ be a component of $\supp b$, meeting $J_a$, a component of $\supp a$. Without loss of generality, $\partial J_a$
meets $J_b$. Since $a$ and $b$ commute, we have that if $x\in (\partial J_a)\cap J_b$, then $b^nx\to\partial J_b$ as $n\to\infty$, and each
of the points in this orbit is fixed by $a$ since $a$ and $b$ commute. We are then immediately reduced to the case $M=I$.

(2) Suppose that $b$ and $c$ do not commute. Then the group generated by $\form{  b,c}$ is not abelian, and therefore
cannot act freely on $(0,1)$, by H\"older's Theorem (Theorem~\ref{thm:holder}). It follows that there is a nontrivial $s\in \form{  b,c}$
with a fixed point in $(0,1)$. Since $[s,a]=1$, Theorem~\ref{thm:kopell} implies that $s$ is the identity, a contradiction.
\end{proof}

The two parts of Proposition~\ref{prop:disj-ab} are sometimes called the \emph{disjointness criterion}\index{disjointness criterion}
and the \emph{abelian criterion}\index{abelian criterion},
respectively.

The most substantial result in this section is the following, which can be found as Lemma 2.8 in~\cite{KK2018JT}.

\begin{prop}\label{prop:tame}
Let $f\in\Diffb(S^1)$ be an infinite order element, and let $Z\le\Diffb(S^1)$ denote its centralizer. We have the following statements.
\begin{enumerate}[(1)]
\item
If $\Fix f\neq\varnothing$ and if $H\le Z$ is generated by elements with nonempty fixed point sets, then every element $h\in H$ satisfies
$\Fix h\neq\varnothing$. Furthermore, we have \[\supp f\cap\supp H'=\varnothing.\]
\item
If the rotation number of $f$ is irrational then $Z$ is abelian.
\item
If the rotation number of $f$ is rational then the rotation number of each element of $Z$ is rational.
\item
If $f$ has rational rotation number then the rotation number restricted to $Z$ is a homomorphism.
In particular, if $h\in Z'$ then $\Fix h\neq\varnothing$.
\end{enumerate}
\end{prop}
\begin{proof}
(1)
Let $J$ be a component of $\supp f$. We have that $H$ acts on $J$, by Proposition~\ref{prop:disj-ab}. Since every element of $H$ preserves
$\partial J$, we have that every element of $H$ has a nonempty set of fixed points. Proposition~\ref{prop:disj-ab} shows that the restriction of
the commutator subgroup of $H$ to $J$ is the identity.

(2)
Denjoy's Theorem (Theorem~\ref{thm:denjoy}) implies that $f$ is topologically conjugate to an irrational rotation. If $1\neq g\in\Homeo_+(S^1)$
commutes with an irrational rotation then $\Fix g=\varnothing$. Indeed, we have that $\Fix g$ is closed and $f$--invariant,
and $f$ is minimal, whence
it follows that $\Fix g= S^1$. It follows that the centralizer of $f$ acts freely on the circle and is
therefore abelian by H\"older's Theorem for the circle (Theorem~\ref{thm:holder-circle}).

(3)
Suppose for a contradiction that $g\in Z$ has irrational rotation number. Then by the previous item, we have that no nontrivial element in
the centralizer of $g$ can have a fixed point. Since $f$ has infinite order,
it follows that $f$ has no periodic points and hence has irrational rotation number, a contradiction.

(4)
Passing to a positive power if necessary, we may assume that $\Fix f\neq\varnothing$. Indeed, passing to a power of $f$ cannot
decrease the size of the centralizer, so assuming $\Fix f\neq\varnothing$ is not a loss of generality. Let $G$ denote the centralizer of $f$, and
let $K\sse  G$ denote the set of elements of $G$ with zero rotation number. Item 1 of the proposition implies that $K$ is actually
a group. Since the rotation number depends only on the conjugacy class of a homeomorphism in $\Homeo_+(S^1)$
(Proposition~\ref{prop:rot-easy}), we have that
$K$ is normal subgroup of $G$.

We have that $\Fix f$ is $G$--invariant, and \[X=\partial \Fix f\] is a closed, $G$--invariant subset of $\Fix f$. We claim that $X$ is pointwise
fixed by $K$. Indeed, this follows from the fact that if $J$ is a component of $\supp f$ then an arbitrary element $k\in K$ fixes $\partial J$,
by Item 1. The set $X$ consists of the closure of the union of the sets $\partial J$, as $J$ ranges over $\supp f$. We thus obtain a map
\[\phi\colon G/K\longrightarrow \Homeo_+(X),\] just defined by restricting the action of $G$ on $S^1$. Since $g\in G\setminus K$ if and only if
$g$ has nonzero rotation number, Corollary~\ref{cor:holder-extend} implies that the action of $G/K$ on $X$ extends to
a map \[\Phi\colon G/K\longrightarrow\Homeo_+(S^1),\]
that this action is free, and that the rotation number is an injective homomorphism \[G/K\longrightarrow S^1.\] If $gK\in G/K$ then it is
straightforward to check that the rotation number of $\Phi(gK)$ coincides with the rotation number of $g$, completing the proof.
\end{proof}

In part (2) of Proposition~\ref{prop:tame}, one can in fact assert that $Z$ is conjugate into the group of rotations of the circle, though we will
not spell this out further. Similarly, in part (4) we may eschew the hypothesis that $f$ have rational rotation number.

\section{Classification of $C^1$--actions of solvable Baumslag--Solitar groups}\label{sec:bs}
For nonzero integers $p$ and $q$, we define the corresponding Baumslag--Solitar group\index{Baumslag--Solitar group} by
\[
\operatorname{BS}(p,q):=\form{a,e\mid ae^pa^{-1}=e^q}.\]
Having a Baumslag--Solitar subgroup is often an obstruction for a group from admitting certain geometric structures. 
For instance, word-hyperbolic groups, $\operatorname{CAT}(0)$ groups, and 3--manifold groups do not contain 
$\operatorname{BS}(p,q)$ for $0<|p| < |q|$ (see~\cite{Shalen-TAIA, Gromov1987, BH1999}). It is also known that 
$\operatorname{BS}(p,q)$ does not admit a faithful $C^2$ action on a compact interval for $1<p<q$ (see~\cite{FF2001}).

In this section, we will classify $C^1$--actions of \emph{non-abelian solvable}\index{Baumslag--Solitar groups}
Baumslag--Solitar groups; that is, we will consider the group
\[B:= \operatorname{BS}(1,m)\]
for a fixed $m>1$. As a result, we will deduce that every $C^1$--action of $\bZ\times\operatorname{BS}(1,2)$ on a compact 
connected one--manifold admits two nontrivial elements of the group with disjoint supports (Corollary~\ref{c:bs12}). 
This result will be crucially used in Chapter~\ref{ch:optimal}, namely for the construction of optimally expanding diffeomorphism groups.

We remark that the group $\operatorname{BS}(1,-m)$ is also solvable, but it does not yield an interesting $C^1$ 
action for us (see Proposition~\ref{p:bs-negative}). We also note that for $p$ and $q$ different from $\pm1$, the group 
$\operatorname{BS}(p,q)$ are non-solvable and contains $F_2$~\cite{KarrassSolitar1971CJM}.

\subsection{Topological smoothing of the standard affine action}\label{ss:bs-smooth}
The group of orientation preserving affine maps on the real line can be described as
\[\Aff_+(\bR)=e^\bR\ltimes \bR=
\left\{\begin{pmatrix} e^s&t\\ 0 & 1\end{pmatrix}\middle\vert s,t\in\bR\right\}.\]
This is naturally a subgroup of $\Diff_+^\infty(\bR)$, and hence of $\Diff_+^\infty(0,1)$. 
After a suitable topological smoothing, we can make this group act smoothly on a compact interval:

\begin{thm}\label{thm:bs-cpt}
There exists an injective group homomorphism
\[
\Aff_+(\bR)\longrightarrow \Diff_{[0,1]}^\infty(\bR)\]
defined as a topological conjugation by a homeomorphism $(0,1)\longrightarrow\bR $.
\end{thm}
More precisely, we will build a homeomorphism
 $\psi\co (0,1)\longrightarrow\bR $ such that  for each $g\in \Aff_+(\bR)$ 
the map 
\[\psi^{-1} g\psi\co (0,1)\longrightarrow(0,1)\]
extends to a map in $\Diff_{[0,1]}^\infty(\bR)$, the group of smooth real line diffeomorphisms fixing all the points outside $[0,1]$.

\bp[Proof of Theorem~\ref{thm:bs-cpt}]
Let us first consider a smooth diffeomorphism $\varphi_0\co (0,1)\longrightarrow\bR$ satisfying $\varphi_0(x) =-1/x$ near $x=0$
and that $\varphi_0(1-x)=\varphi_0(x)$ for all $x\in(0,1)$.
Then the affine transformation $f(x) = sx+t$ is conjugated to the map
\[
F(x):=\varphi_0^{-1}\circ f\circ\varphi_0(x) = \frac{x}{s-tx}\]
Since $s>0$ the map $F$ is smooth near $x=0$.
 It follows by symmetry at $x=1$ that $F\in \Diff_+^\infty[0,1]$.

By the Muller--Tsuboi trick (Theorem~\ref{t:muller-tsuboi}),
there exists a smooth diffeomorphism $\varphi$ of $(0,1)$ (which is not differentiable at $\{0,1\}$) such that 
\[\varphi^{-2}\Diff_+^\infty[0,1]\varphi^{2}\le\Diff_{[0,1]}^\infty(\bR).\]
Then the map
$\psi:=\varphi_0\varphi^2$
satisfies the condition above.
\ep

\begin{rem}Alternatively, one can view $\Aff_+(\bR)$ as a group of matrices acting on $\bR^2$
and preserving the $x$--axis. Projectivizing, this group acts on a compact interval contained in the unit circle. This action can be seen to be smooth~\cite{Tsuboi1987}.\end{rem}

The following definition will be handy for us.
\bd\label{defn:tangent-id}
Let $r\in\bZ_{>0}\cup\{\infty\}$.
A $C^r$--diffeomorphism $f$ of a manifold $M$
is said to be \emph{$C^r$--tangent to the identity}\index{tangent to the identity} at $p\in M$
if the following hold.
\begin{itemize}
\item $f(p)=p$;
\item $D^if(p)=D^i(\Id)(p)$ for $i=1,\ldots,r$.\end{itemize}\ed
Here, $D^i$ denotes the $i$--th derivative. 
In the context of one--manifold actions, this simply means that
$f(p)=p$, $f'(p)=1$ and $\der{f}{i}(p)=0$ for $i=2,\ldots,r$.

The group $\operatorname{BS}(1,m)=\form{a,e\mid aea^{-1}=e^m}$ admits the \emph{standard affine action}\index{affine group action} generated by
\[\bar a(x) =mx,\quad \bar e(x)=x+1.\]
By the theorem above, we see that this action is topologically conjugate to a smooth action on $[0,1]$  that is $C^\infty$--tangent to the identity at the endpoints.
We will continue to use the symbols $\bar a$ and $\bar e$ to represent the standard affine generators of $\operatorname{BS}(1,m)$.

\subsection{$C^1$--actions on intervals}
Let us now describe a topological classification of nonabelian $C^1$--actions of solvable Baumslag groups, which 
is a special case of results in~\cite{BMNR2017MZ}; see~\cite{TriestinoNote} for a nice exposition particularly for
 Baumslag--Solitar groups. A similar result can be found in~\cite{CC-unpublished}.

\begin{thm}\label{t:bmnr}
Let $m\ge2$ be an integer. Suppose we have a representation
\[
\rho\co \operatorname{BS}(1,m)\longrightarrow \Diff_+^1[0,1].\]
Then each $J\in\pi_0\supp \rho(e)$ is preserved by the whole image of $\rho$, 
and the restriction of the action $\rho$ on $J$ is topologically conjugate to the standard affine action of $\operatorname{BS}(1,m)$.
Moreover, the derivative of $\rho(a)$ at its unique fixed point in $J$ is $m$.
\end{thm}

We prove the preceding theorem in this subsection.
By the following simple observation, we see that the conclusion is vacuous when $\rho$ is not faithful, i.e. $\rho(e)=1$:
\begin{lem}\label{l:faithful}
Let $m\ge2$ be an integer.
If \[\phi\co \operatorname{BS}(1,m)\longrightarrow Q\]
is a quotient such that the image of $e$ has an infinite order, then $\phi$ is an isomorphism. 
\end{lem}
\bp
Pick an arbitrary element $g\in \ker\phi$. Up to conjugation, we can write $g=e^pa^q$ for some $p,q\in\bZ$. We have
\[
1=\phi(aga^{-1})=\phi(e^{mp}a^q)=\phi(e)^{(m-1)p}.\]
This implies that $p=0$ and $\rho(a)^q=1$.
We have that
\[
\rho(e) = \rho(a^qea^{-q})=\rho(e^{m^q}).\]
It follows that $m^q=1$ and $q=0$. We have shown $\ker\phi=\{1\}$.
\ep

Set
\[
H:=\fform{e}=\form{a^{-k}ea^k\mid k\ge0}\unlhd \BS(1,m).\]
We have an isomorphism from the additive group $\bZ[1/m]$ to $H$  that maps $1/m^k$ to $a^{-k}ea^k$.
For each $r\in \bZ[1/m]$, we let $e[r]$ to denote its image under this isomorphism.
Note that for $n\in\bZ$ we have $e[n]=e^n$.

Recall that an action $\rho\co K\longrightarrow \Homeo_+(\bR)$ of a group $K$ is \emph{semi-conjugate}\index{semi-conjugacy}
to another action
$\rho'\co K\longrightarrow\Homeo_+(\bR)$ if for some monotone increasing surjective (hence, continuous) map $q\co \bR\longrightarrow \bR$ 
we have a commutative diagram

\[\begin{tikzcd}
\bR\arrow{r}{\rho(g)} \arrow{d}{q} &
\bR\arrow{d}{q}\\
\bR\arrow{r}{\rho'(g)}& \bR
\end{tikzcd}
\]
 for all $g\in K$. The map $q$ is called a {semi-conjugacy}. 
One may call $\rho$ above as a \emph{blow-up}\index{blow-up} of $\rho'$ and define two representations are 
\emph{semi-conjugately equivalent} (or, \emph{semi-conjugate} for short) if they have a common blow-up (see~\cite{KKM2019}, cf.~Theorem~\ref{thm:irr-rot-semi}).

Let us first consider the setting of topological actions by $\operatorname{BS}(1,m)=\form{a,e}$.
\begin{lem}\label{l:bmnr-c0}
If $\BS(1,m)$ acts topologically on $\bR$ such that $\Fix e=\varnothing$, then this action is semi-conjugate to the standard affine action.
\end{lem}

\bp
For all $r\in \bZ[1/m]$, the element $e[r]$ has no fixed point since $e[r]$ and $e$ has a nontrivial common power.
After a topological conjugacy, we may assume $e(x)=x+1$ (up to inverting $e$ to $e^{-1}$ if necessary). 
Define a map $q\co H.0\longrightarrow \bZ[1/m]$ by
\[
q\left(e[r]\right)=r.\]
It is obvious that the map $q$ is well-defined, strictly increasing and surjective.
We consider its continuous extension
\[
\bar q\co \bR\longrightarrow \bR\]
satisfying
\[
\bar q(x) :=\sup q\left(H.0\cap(-\infty,x]\right)=\inf q\left(H.0\cap[x,\infty)\right).\]
It is not hard to see that $\bar q$ is well-defined, monotone increasing and surjective.
One can also see from the density of $\bZ[1/m]$ in $\bR$ that 
\[
\bar q(e[r].x)=\bar q(x)+r\]
for all $x\in\bR$ and $r\in \bZ[1/m]$.

Let us set $\sigma:=\bar q(a.0)$. Then we have that
\[
\bar q(a e[r].0)=\bar q (a e[r] a^{-1}\cdot a.0)=\bar q(e[{mr}]\cdot a.0)=mr+\sigma=m\bar q(e[r].0)+\sigma.\]
By the density of $\bZ[1/m]$ again, we have that 
\[\bar q(a.x)=m\bar q(x)+\sigma\]
for all $x$. Hence, the map $\bar q$ semi-conjugates the given action to the affine action
\[
\bar a(x)=mx+\sigma,\quad \bar e(x)=x+1.\]
We can then further conjugate $\bar a(x)$ to $x\mapsto m x$, as desired.
\ep

The last ingredient of the proof of Theorem~\ref{t:bmnr} is the following.
\begin{prop}\label{p:bmnr}
If $\rho\co \BS(1,m)\longrightarrow\Diff_+^1[0,1]$ is a faithful representation
and if $\rho$ does not have a global fixed point in $(0,1)$,
then the following hold.
\be[(1)]
\item\label{p:p:bmnr-a} $\Fix\rho(a)\cap(0,1)\ne \varnothing$ and $\Fix\rho(e)\cap(0,1)=\varnothing$;
\item $\rho(\BS(1,m))\restriction_{(0,1)}$ is conjugate to the standard affine action;
\item\label{p:bmnr-m}  At the unique fixed point $x_0$ of $\rho(a)$, we have that $\rho(a)'(x_0)=m$.
\ee
\end{prop}

The essence of the proof for the above proposition is the following counting argument.
\begin{lem}\label{l:counting}
Let $C>0$ be a constant, let $\mathcal{N}\sse\bZ$ be an infinite set,
 let $K\le\Diff_{[0,1]}^1(\bR)$ be a group with a finite generating set $S$,
and let $v\in K$ be a $C^1$ diffeomorphism satisfying $v(x)>x$ for all $x\in (0,1)$.
Assume we have subsets 
\[ S_1, S_2, \ldots\sse K\]
 and  compact intervals 
 \[ J_1, J_2, \ldots\sse(0,1)\]
satisfying the following three conditions for all $N\in\mathcal{N}$ and for all $n\ge1$:
\be[(i)]
\item\label{l:p:disj} The intervals in the set $\{g v^N(J_n)\mid g\in S_n\}$ have pairwise disjoint interiors;
\item\label{l:p:suffix} Each $g\in S_n$ can be written as 
\[g=s_\ell\cdots s_2s_1\]
for some $\ell\le Cn$ and $s_i\in S\cup S^{-1}$ such that 
\[ \left(\bigcup_{j=1,\ldots,\ell} s_j\cdots s_2s_1 v^N J_n \right)\sse [0,1]\setminus\left( v^{-|N|}\sup J_n, v^{|N|} \inf J_n\right)\]

\item\label{l:p:definite} The set $\bigcap_n \bigcup_{g\in S_n}gv^NJ_n$ has a nonempty interior.
\ee
Then we have that
\[ \lim_{n\to\infty} \left(\#S_n \cdot |J_n|\right)^{1/n}=1.\]
\end{lem}

The element $s_j\cdots s_2s_1$ above is called a \emph{suffix}\index{suffix} of $g$.
\begin{rem}\label{r:counting}
Note that the hypothesis of the lemma is invariant under a topological conjugacy.
The conclusion is also topologically invariant in the special case when all $J_n$'s coincide with a fixed interval $J_0$, since we   have
 \[\lim_{n\to\infty} |J_n|^{1/n}=\lim_{n\to\infty} |J_0|^{1/n}=1.\]
 Therefore in this special case, we can weaken the hypothesis so that $K\le\Homeo_+[0,1]$ is only topologically conjugate into
 $\Diff_+^1[0,1]$, and conclude that 
\[\lim_{n\to\infty} \left(\# S_n\right)^{1/n}=1.\]
\end{rem}

Intuitively speaking, the lemma addresses the situation where some set $S_n$ in the radius $O(n)$ ball of the Cayley graph 
disjointly translates the interval $v^N(J_n)$ in such a way that
the set $S_n\cdot v^N\cdot J_n$ has a definite length depending only on $N$,
and such that every suffix $h$ of $g\in S_n$ places $v^N(J_n)$ no farther from the boundary $\{0,1\}$ of $I$ than $v^N(J_n)$ itself.
If this happens arbitrarily near from the boundary (that is, for $|N|\gg0$) then
the growth rate of $S_n$ is equal to that of $1/|J_n|$.
It will be enlightening for the reader to give a proof of Proposition~\ref{p:bmnr} first, in order to illustrate concrete instances of the lemma.
\bp[Proof of Proposition~\ref{p:bmnr}, assuming Lemma~\ref{l:counting}]
As before, we simply write $g$ to denote the image of $g\in B$ under $\rho$ when the meaning is clear.
We will apply Lemma~\ref{l:counting} by setting \[K=\BS(1,m),\quad S=\{a,e\},\quad C=m+1,\] and varying the other parameters.

The conclusions of parts (1) and (2) are invariant under topological conjugacy. By Muller--Tsuboi trick, 
we may assume that the image of $\rho$ is in $\Diff_{[0,1]}^1(\bR)$ for those parts and apply Lemma~\ref{l:counting}.
To prove part (1), 
let us first note that the former half of the statement $\Fix a\cap(0,1)\ne\varnothing$ implies the latter half $\supp e =(0,1)$.
Indeed, if $\supp e\ne(0,1)$ then one can find some $x_0\in \Fix e\cap(0,1)$,
together with either the forward or the backward limit of $\{a^k(x_0)\}_k$ that is a global fixed point in the set $\Fix a\cap(0,1)$. 
This contradicts the absence of global fixed points. 

It therefore suffices for us to deduce a contradiction after assuming $\Fix a\cap(0,1)=\varnothing$.
We may assume $a(x)>x$; otherwise we can conjugate $\rho$ by the inversion $\sigma(x)=1-x$ 
and apply the same argument. Let $J\in\pi_0\supp e$; we may further assume $e(x)>x$ on $J$ by 
replacing $e$ by $e^{-1}$ if necessary.
Pick $J_0$ be a compact interval in $J$ such that $J_0$ and $e(J_0)$ are disjoint.
For each $n\ge1$ we set 
\[ S_n:=\left\{ e[r] \mid r\in(-m^n,0]\cap\bZ\right\}\sse H\]
and $J_n:=J_0$. 

We claim that the hypothesis of the lemma is satisfied with $v:=a$, for all  $N<0$.
Indeed, the condition (i) follows from that
\[ e^{-1} a^N(J_0)\cap a^N(J_0)
=a^N e^{-m^{-N}}(J_0)\cap a^N(J_0)=\varnothing.\]
For each $g\in S_n$ there exists some $k_0,\ldots,k_{n-1}\in\{0,\ldots,m-1\}$ such that
\[g= e\left[-\sum_{0\le i\le n-1} k_i m^i\right]= \prod_{i=0}^{n-1} a^i e^{-k_i} a^{-i}
=
a^{n-1} e^{-k_{n-1}}a^{-1} e^{-k_{n-2}}\cdots a^{-1}e^{-k_0}.
\]
The word length of $g$ is at most 
\[2(n-1)+\sum k_i \le (m+1)n.\]
Then for every suffix $h$ of the word $g$ and for every $x\in a^N(J)$, we have that \[h(x)\in a^{N-i}(J)\] for some  $i\in\{0,1,\ldots,n-1\}$. 
This verifies the condition (ii). The condition (iii) is trivial since $J_n=J_0$ is independent of $n$
and since each $S_n$ contains the identity.
As mentioned in Remark~\ref{r:counting}, we can compute
\[\lim_{n\to\infty} \left(\#S_n\cdot |J_n|\right)^{1/n}=\lim_{n\to\infty} \left(\#S_n\right)^{1/n}=m.\]
Since this contradicts Lemma~\ref{l:counting}, we conclude that $\Fix a\cap(0,1)\ne\varnothing$.

From now on, we will assume $e(x)>x$ without loss of generality.
To prove part (2), it suffices for us to show that the action $\rho(B)\restriction_{(0,1)}$ is minimal since we already know that 
 $\rho(B)$ is semi-conjugate to the standard affine action.
 Assume the contrary, and pick a wandering interval $J_0$ of this action. Now, apply Lemma~\ref{l:counting} for  $v=e$, $J_n:=J_0$
 and
\[ S_n:=\left\{ e[r] \middle\vert r\in[0,1)\cap \frac1{m^{n}}\bZ\right\}\sse H.\]
 
For all sufficiently large $N>0$ we have that $a(x)>x$ whenever $x\in a^N J_0$, by semi-conjugacy.
The conditions (i) and (iii) are obvious since $\rho(H)$ is a blow-up of a free affine action 
and since $J_n=J_0$ is independent of $n$.
For part (ii), note that
each $g\in S_n$ can be written as
\[
g =  e\left[\sum_{1\le i\le n} k_i m^{-i}\right]= \prod_{i=1}^{n} a^{-i} e^{k_i} a^{i}
=
a^{-n} e^{k_{n}}a e^{k_{n-1}}\cdots ae^{k_1}a.
\]
The word length of $g$ is again at most $(m+1)n$. 
We also see that for each  suffix $h$ of $g\in S_n$, we have that $ha^N(\inf J_0)\ge \inf a^N J_0$. 
This contradicts the lemma again, since we have
\[\lim_{n\to\infty} \left(\#S_n\cdot |J_n|\right)^{1/n}=\lim_{n\to\infty} \left(\#S_n\right)^{1/n}=m.\]

Finally, we establish part (3). Let us first assume that $\rho(B)\le\Diff_{[0,1]}^1(\bR)$.
We saw in part (2) that there exists a conjugacy $q\co (0,1)\longrightarrow\bR$ from $\rho$ to the standard affine action.
We set  \[x_0:=q^{-1}(0),\quad J_0=q^{-1}[0,1]\quad \textrm{and}\quad 
J_n:=q^{-1}[0,1/m^n]=a^{-n}(J_0).\]
We let $v=e$ and $S_n$ be the same as part (2) above,
and claim that the conditions (i) through (iii) hold for these choices and for all positive $N$.
Note that all these conditions are invariant under topological conjugacy and easy to check for the standard affine actions.
For instance, the condition (iii) is verified by the computation
\[
q\left( \bigcap_{n\ge1} \bigcup_{g\in S_n} ge^N J_n\right)
= \bigcap_{n\ge1}[N,N+1]=[N,N+1].\]

Using the Mean Value Theorem, it is elementary to see that
\[a'(x_0)=\lim_{n\to\infty} \left(\frac{|J_0|}{|a^{-n}(J_0)|}\right)^{1/n}=\lim_{n\to\infty} \frac1{|J_n|^{1/n}}.\]
We obtain from the lemma that
\[1=\lim_{n\to\infty} \left(\#S_n\cdot |J_n|\right)^{1/n}=\lim_{n\to\infty}m/ a'(x_0),\]
and that $a'(x_0)=m$.

In order to complete the proof, assume that $\rho(\BS(1,m))\le\Diff_+^1[0,1]$.
By the Muller--Tsuboi trick, we obtain a $C^\infty$ homeomorphism $h$ such that
\[ h\rho(\BS(1,m))h^{-1}\le\Diff_{[0,1]}^1(\bR).\]
We have seen so far  that at the unique interior fixed point $x_0$ of $h\rho(a)h^{-1}$, 
one can compute
\[h'\circ\rho(a)\circ h^{-1}(x_0)\cdot \rho(a)'\circ h^{-1}(x_0)/h'\circ h^{-1}(x_0)=m.\]
The point $h^{-1}(x_0)=y_0$ is the unique interior fixed point of $\rho(a)$
and satisfies
\[h'(y_0)\cdot \rho(a)'(y_0)/h'(y_0)=m.\]
This proves that $\rho(a)'(y_0)=m$.

\ep

To complete the proof of Theorem~\ref{t:bmnr}, it remains to show Lemma~\ref{l:counting}.
\bp[Proof of Lemma~\ref{l:counting}]
Assume the hypotheses of the lemma.
We will fix an arbitrarily small $\epsilon>0$, and let $n$ vary over $\bZ_{>0}$.
Since $K$ is $C^1$--tangent to the identity at the boundary, 
there exists some $N\in\mathcal{N}$ such that for each element $s\in S\cup S^{-1}$
and for each point
\[
x\in  [0,1]\setminus\left( v^{-|N|}\sup J_n, v^{|N|} \inf J_n\right)\]
we have that
\[
\abs*{s'(x)-1}\le\epsilon.\]

By  conditions~(\ref{l:p:disj}) and (\ref{l:p:definite}), we have some $c_N>0$ such that
\[
c_N \le \sum_{g\in S_n} \abs*{gv^NJ_n}\le 1\]
for all $n$.
For each $g=s_\ell\cdots s_1\in S_n$, there exists some $z=z_{N,n}\in v^NJ_n$ such that
\[
|gv^NJ_n|/|v^NJ_n| = g'(z) = \prod_{j=\ell}^1 s_j'\left(s_{j-1}\cdots s_1(z)\right).\]
Using  condition~(\ref{l:p:suffix}), we have
\[
 \# S_n\cdot |v^NJ_n|\cdot (1-\epsilon)^{Cn}\le
\sum_{g\in S_n} |gv^NJ_n|
\le \# S_n\cdot |v^NJ_n|\cdot (1+\epsilon)^{Cn}.\]
We also have
\[
\left(\inf v'\right)^N|J_n|\le |v^NJ_n| \le\left (\sup v'\right)^N|J_n|.\]
Combining the above inequalities, we obtain
\[
1=\liminf_{n\to\infty} c_N^{1/n}
\le 
\liminf_{n\to\infty} \left(\# S_n\cdot \left(\sup v'\right)^N\cdot |J_n|\right)^{1/n}\cdot (1+\epsilon)^C,\]
together with
\[
\limsup_{n\to\infty} \left(\# S_n\cdot \left(\inf v'\right)^N\cdot|J_n|\right)^{1/n}(1-\epsilon)^C
\le 
1.\]
Since the choice of $\epsilon>0$ is arbitrary, we conclude that 
\[\liminf_{n\to\infty} \left(\# S_n\cdot |J_n|\right)^{1/n}=\limsup_{n\to\infty} \left(\# S_n\cdot |J_n|\right)^{1/n}=1.\]
This completes the proof of the lemma.
\ep

In~\cite{BMNR2017MZ}, a different proof (in a more general setting) was given to part (1) of Proposition~\ref{p:bmnr}. 
The proof employs an idea of Thurston (see~\cite{Thurston1974Top}, cf. Appendix~\ref{ch:append3}) that
diffeomorphisms of a compact interval are ``asymptotically translations'' at the endpoints. 
Similar arguments are repeatedly used over 
the literature including~\cite{McCarthy2010,KKT2020}. To give the reader a feeling for this new idea, we
describe a self-contained proof along these lines.

\bp[Alternative proof of Proposition~\ref{p:bmnr} part (1)]
As we have seen in the previous proof of this proposition, 
it suffices for us to show that $a$ must fix some point in $(0,1)$.
Assume $a(x)>x$ for all $x\in (0,1)$ for contradiction.
For all $g\in B$ and $x\in(0,1)$, we introduce the notation
\[
\Delta_xg=g(x)-x.\]

Pick an arbitrary $\epsilon>0$. 
There exists a $\sigma\in(0,1)$ such that 
\[ s'(y)\in(1-\epsilon,1+\epsilon)\]
for all $s\in\{a^{\pm1},e^{\pm1}\}$ and 
for all $y\in(0,1)\setminus(\sigma,1-\sigma)$.
We choose $x$ to be sufficiently smaller than $\sigma$ so that 
$e^m(x)$ and $a^{\pm1}(x)$ are also less than $\sigma$. 
Then we have
\[
| \Delta_xe^m -m\Delta_xe |\le \sum_{i=1}^{m-1} \abs*{(e^i(e(x))-e^i(x))-(e(x)-x)}
\le
\epsilon'\Delta_xe \]
for $\epsilon'=(m-1) ((1+\epsilon)^{m-1}-1)$.
This inequality may be regarded as an asymptotic linearization \[\Delta_xe^m \approx m \Delta_xe .\]

We also have some $z$ between  $x$ and   $e^m(x)$ so that 
\begin{align*}
\abs*{\Delta_{a^{-1}(x)}e -\Delta_xe^m }
&=\abs*{ (a^{-1}e^m(x) - a^{-1}(x)) - (e^m(x)-x)}
=|\Delta_{e^m(x)}a^{-1}-\Delta_x a^{-1}|\\
&=|(a^{-1})'(z)-1|\cdot\Delta_xe^m \le \epsilon\cdot (m+\epsilon')\Delta_xe .\end{align*}
By choosing $\epsilon,\epsilon'>0$ to be sufficiently small from the beginning,
we obtain that
\[
\abs*{\Delta_{a^{-1}(x)}e }
\ge m|\Delta_xe |-
\abs*{\Delta_{a^{-1}(x)}e -m\Delta_xe }\ge
\left(m-
1/2\right)\Delta_xe .\]
Since $a^{-1}(x)$ is even closer to $0$ than $x$, we can iterate the estimate and see that
\[
\abs*{\Delta_{a^{-k}(x)}e }\ge (m-1/2)^k\Delta_xe \]
for all $k\ge1$. 
Since $\supp e$ is $a$--invariant, there exists an $x$ arbitrarily close to $0$ such that $\Delta_xe \ne0$.
We then have a contradiction, since the right-hand side diverges to the infinity as $k\to\infty$,
while $a^{-k}(x)$ converges to $0$.
\ep

We conclude this subsection by noting that the solvable Baumslag--Solitar group with a negative sign
\[
\operatorname{BS}(1,-m)=\form{a,e\mid aea^{-1}=e^{-m}}\]
does not admit any interesting $C^1$--faithful actions on a compact interval for $m>1$. 
We thank Crist\'obal Rivas for teaching us this proof.

\begin{prop}\label{p:bs-negative}
For $m>1$, every orientation--preserving $C^1$--action of $\operatorname{BS}(1,-m)$ on $[0,1]$ is abelian.
\end{prop}
\bp
We may assume that the action does not have a global fixed point in the interior $(0,1)$. 
We claim that the image of $e$ is trivial.
If $e(x)>x$ for all $x\in (0,1)$, we would obtain a contradiction since
\[
a a^{-1}(x) < aea^{-1}(x)=e^{-m}(x)<x.\]
By symmetry, it follows that $\Fix e\cap(0,1)\ne\varnothing$.

The rest of the proof goes almost identically to that of Proposition~\ref{p:bmnr} (\ref{p:p:bmnr-a}).
We consider the same set of words $S_n$, and 
note a slightly different expansion 
\[g= e\left[-\sum_{0\le i\le n-1} k_i m^i\right]= \prod_{i=0}^{n-1} a^i e^{-(-1)^ik_i} a^{-i}
=
a^{n-1} e^{-(-1)^{n-1}k_{n-1}}a^{-1} e^{-k_{n-2}}\cdots a^{-1}e^{-k_0}.
\]
Then the same estimates goes through and Lemma~\ref{l:counting} again applies to yield a contradiction.
We conclude that the image of $e$ is trivial.
\ep

\subsection{$C^1$--actions on circles}
We now consider $C^1$--actions of 
\[
\operatorname{BS}(1,2)=\form{a,e\mid aea^{-1}=e^2}\]
on the circle. 
First, we note the existence of a finite orbit.

\begin{thm}[\cite{GL2011}]\label{t:gl2011}
Every faithful $C^1$--action of $\operatorname{BS}(1,2)$ on a circle
admits a finite orbit.
\end{thm}
\begin{rem}\label{r:gl2011}
In the preceding theorem, note first that $e$ must fix a point since
\[\rot e = \rot \left(aea^{-1}\right)=\rot e^2=2\rot e.\]
In particular, if the action of $e$ has finite order then it must be trivial.
It follows that the faithfulness of the action is equivalent to the nontriviality of the action of $e$.
Moreover, if $x_0$ belongs to a finite orbit of $\BS(1,2)$ as in the conclusion of the theorem, then $x_0\in \Fix e$.
It follows that for some $m>0$, the point $x_0$ is fixed by the finite index subgroup \[\form{a^m,e}\cong\operatorname{BS}(1,2^m).\]
\end{rem}

\bp[Proof of Theorem~\ref{t:gl2011}]
We will use a counting argument as in Lemma~\ref{l:counting}. 
We have noted that $\Fix e\ne\varnothing$. 
Let us pick  $I_0\in\pi_0\supp e$.
Since
\[ a \Fix e = \Fix  \left(aea^{-1}\right)=\Fix e^2=\Fix e,\]
we have that $a^k(I_0)\in\pi_0\supp e$ for all $k$. 

\begin{claim}
We have that $a^{-m}(I_0)=I_0$ for some $m>0$.\end{claim}

Once the claim is proved, 
we can deduce the conclusion of the theorem by choosing a point $x_0\in\partial I_0$ and noting that
\[ \#\left(\form{a,e}x_0\right)=\#\left(\form{a}x_0\right)\le m.\]

In order to prove the claim,
assume for contradiction that
$\{a^{-k}(I_0)\}_{k\ge 0}$ is a disjoint collection of open intervals in $S^1$. 
Let $\epsilon>0$ be arbitrary.
Note that
\[ \lim_{k\to\infty} \abs*{a^{-k}(I_0)}=0.\]
By the uniform continuity of $e'$ and $\log a'$,
whenever $k$ is sufficiently large and $x,y\in a^{-k}(I_0)$, 
we have that
\[
{|e'(x)-1|}<\epsilon\]
and
\[
\abs*{a'(x)/a'(y)-1}<\epsilon.\]
After replacing $I_0$ by $a^{-k}(I_0)$ for some $k\gg0$ if necessary, 
we may assume that the above two inequalities hold for all $k\ge 0$
and for all $x,y\in a^{-k}(I_0)$.

We then proceed in a manner similar to that of Proposition~\ref{p:bmnr} (1).
Pick an arbitrary $n\ge1$.
We fix an open interval $J_0\sse I_0$ such that $e(J_0)\cap J_0=\varnothing$,
and set 
\[ S_n:=\{e^i\mid 0\le i<2^n\}.\]
Then $S_n(J_0)$ is a disjoint collection of open intervals in $I_0$. As in the proof of Proposition~\ref{p:bmnr}, we can write
\[
e^i=a^{n-1} e^{k_{n-1}}a^{-1}\cdots e^{k_1}a^{-1} e^{k_0}\]
for some $k_i\in\{0,1\}$. 
For each suffix $h$ of the above word, the interval $h(I_0)$ is contained in $a^{-t}(I_0)$ for some $t\ge0$. 
Hence,
we can find some $x_j,y_j,z_j\in a^{-j}(I_0)$ for $j=0,1,\ldots,n-1$ such that
\[
\abs*{e^i(J_0)}=\prod_{j=0}^{n-1} \left(e^{k_j}\right)'(x_j) \cdot \prod_{j=0}^{n-1}a'(y_j)/a'(z_j) \cdot|J_0| \]

It follows that
\[
1\ge \sum_{i=0}^{2^n-1} \abs*{e^i(J_0)}\ge 2^n (1-\epsilon)^{2n} |J_0| .\]
By choosing $2(1-\epsilon)^2>1$ and letting $n\to\infty$, we have the desired contradiction.
The claim is now proved.
\ep

The following corollary asserts that $\bZ\times\operatorname{BS}(1,2)$ admits pair of elements with 
``universally disjoint supports'', meaning they have disjoint supports for arbitrary $C^1$--actions on a compact connected one--manifold.
This will be a key ingredient for our construction of an optimally expanding diffeomorphism group.

\begin{cor}\label{c:bs12}
Let $M^1$ be a compact connected one--manifold. Then for every  representation
\[
\bZ\times \operatorname{BS}(1,2)=\form{c}\times \form{a,e\mid aea^{-1}=e^2}\longrightarrow\Diff_+^1(M^1),\]
there exists some $k\ge 1$ depending only on the image of $\form{a,e}$ such that 
\[
\supp c^k\cap \supp e=\varnothing.\]
In particular, for every representation
\[ \rho\co (\bZ\times \operatorname{BS}(1,2))\ast \bZ=(\form{c}\times\form{a,e})\ast\form{d}\longrightarrow\Diff_+^1(M^1),\]
 the support of the image of the group element 
\[
u_0:=\left[ [c^d,e e^d e^{-1}],c\right]\]
is contained in some compact subset of  $\supp\rho$.
\end{cor}
Note that if $M^1=I$ then we can choose $k=1$ since \[\supp c=\supp c^2=\cdots\] in this case.
Note also that the second part of the conclusion obviously follows from the first part
and from a consequence of the $abt$--Lemma, namely Lemma~\ref{l:ucsd}.

We will spend the remainder of this section in establishing the above corollary, beginning with the case $M^1=I$.

\begin{lem}\label{l:bs1m}
If $m\ge2$ is an integer, then for an arbitrary action
\[
\bZ\times\operatorname{BS}(1,m)=\form{c}\times\form{a,e}\longrightarrow\Diff_+^1[0,1],\]
the supports of $c$ and $e$ are disjoint.\end{lem}
\bp
We may assume the action does not have a global fixed point in $(0,1)$. 
Assume for contradiction that $J_0\in\pi_0\supp c$ and $J_1\in\pi_0\supp e$ nontrivially intersect. 
We then have either $J_0\sse J_1$ or $J_1\sse J_0$.

We note from Theorem~\ref{t:bmnr} that the restriction of the action of $\operatorname{BS}(1,m)$ to
$J_1$ is  topologically conjugate to the standard affine action.
Moreover, there exists $p\in J_0\cap \Fix a$ with $a'(p)=m$. 

Since $c^k(p)$ is fixed by $a$ for all $k\in\bZ$
and since $\Fix a\cap J_1$ is a singleton, we have that $J_0\ne J_1$.
If $J_1\subsetneq J_0$, then $\{c^k(p)\}_{k\in\bZ}$ is an infinite set of hyperbolic fixed points of $a$ with derivative $m$. 
This is a contradiction, since at accumulation points of $\Fix a$ the derivative of $a$ should be one. Hence, we have $J_0\subsetneq J_1$. 

Since the action of $\operatorname{BS}(1,m)$ on $J_1$ is topologically conjugate to the standard affine action, we should have $p\in \Fix c$. 
Moreover, the action of $H=\fform{e}$ on $J_1$ is minimal and preserves $\Fix c\cap J_1$. This implies that 
$\Fix c\cap J_1$ is dense in $J_1$, which contradicts the assumption $J_0\sse J_1$.
\ep

By plugging $m=2$, we obtain the conclusion of the corollary for the case $M^1=I$.
For the case of a circle, we prove a stronger result below.

\begin{lem}[\cite{KK2020crit}]\label{l:bs-z1}
Suppose we have a representation
\[ \operatorname{BS}(1,2)=\form{a,e}\longrightarrow\Diff_+^1(S^1).\]
We denote the centralizer by
\[ Z^1:=\{c\in \Diff_+^1(S^1)\mid [c,g]=1\text{ for all }g\in \operatorname{BS}(1,2)\}.\]
Then there exists a normal subgroup $Z_0\unlhd Z^1$ such that $Z^1/Z_0$ is a finite cyclic group and such that 
\[
\supp Z_0\cap \supp e=\varnothing.\]
\end{lem}
\bp 
We may assume that $\BS(1,2)$ acts faithfully, for otherwise $e$ acts trivially and there is nothing to show (cf.~Remark~\ref{r:gl2011}).
So, we may regard $\operatorname{BS}(1,2)$ as a subgroup of  $\Diff_+^1(S^1)$. 
Moreover, there exists an $m\ge1$ such that
\[
B_0:=\form{a^m,e}\cong\operatorname{BS}(1,2^m)\]
fixes a point. Applying Theorem~\ref{t:bmnr} to $B_0$, we see that each component $J\sse\supp B_0$ contains a point $p_J$ satisfying
\[
\left(a^m\right)'(p_J)=2^m.\]
We conclude that $\supp B_0$ has only finitely many components. 

Let $X$ be the union of  the boundary points of all the components of $\supp B_0$.
Since $Z^1$ preserves $\Fix B_0$, it permutes the finite set $X$ preserving the circular order on $X$. Therefore, there exists a representation
\[
\phi\co Z^1\longrightarrow \bZ/k\bZ\]
for $k:=\#X$; this is also an easy instance of Corollary~\ref{cor:holder-extend}.

Write $Z_0:=\ker\phi$. Then for every $z\in Z_0$ and for every $x_0\in X$, we obtain an action
\[
\form{z}\times\operatorname{BS}(1,2^m)\longrightarrow \Diff_+^1(S^1\setminus\{x_0\}).\]
By Lemma~\ref{l:bs1m}, the support of $z$ is disjoint from that of $e$. 
\ep

We thus obtain a proof of Corollary~\ref{c:bs12} by noting that 
\[
c^{K}\in Z_0,\]
where here $K=\# (Z^1/Z_0)$ in the above lemma.


%
%
%
\chapter{Chain groups}\label{ch:chain-groups}

\begin{abstract}In this chapter, we develop some ideas about a particularly natural class of diffeomorphism groups, which the authors and Lodha
named \emph{chain groups}\index{chain groups} in~\cite{KKL2019ASENS}.
 Roughly, a chain of intervals is a sequence of open intervals $\{J_1,\ldots,J_n\}$,
with $J_i\sse  I$ for all $i$, such that $J_i\cap J_k\neq\varnothing$ if and only if $|i-k|\leq 1$.
 One then considers the subgroup $G$
of $\Homeo_+[0,1]$ generated by homeomorphisms $\{f_1,\ldots,f_k\}$ such that $\supp f_i\subseteq J_i$. The group $G$ is called a chain
group if a certain mild dynamical condition is met.
It turns out that chain groups exhibit a combination of uniform properties, together with a remarkable diversity of behaviors. The flexibility
afforded by the diverse phenomena one can observe among chain groups, together with a unified theory that makes them tractable, makes
chain groups a powerful tool that plays a critical role in the construction of groups of homeomorphisms of a given critical regularity.\end{abstract}

\section{Chains and covering distances}
For a group $G$ acting on a topological space $X$, we recall the definition of the \emph{support}\index{support of a group element}
(which in this chapter will
always mean \emph{open support}\index{open support} unless otherwise noted)
of $g\in G$, that is to say the set $\supp g:=X\setminus\Fix g$.
We let $\suppc g$ denote its closure. We denote
\[\supp G:=\bigcup_{g\in G}\supp g.\]
If $X$ is an interval then each component of $\supp g$ is called a \emph{supporting interval}\index{supporting interval} of $g$.

In the case where $X=I$ is a (usually compact) interval and where $G$ is finitely generated, 
we will construct a combinatorial function that behaves much like a metric,
called the \emph{covering distance}\index{covering distance}.
Roughly speaking, this metric measures the ``topological complexity'' of coverings of a closed interval by supporting intervals of the generators.

In general, a  collection of subspaces in a topological space is said to have \emph{finite intersection multiplicity}\index{finite intersection
multiplicity}
if each point in the space belongs to at most finitely many elements of the collection (cf.~\cite{dMvS1993}).

\bd\label{defn:chain}
We say that a finite sequence of intervals
\[I_1,I_2,\ldots,I_m\]
in the real line is an \emph{$m$--chain}\index{$m$--chain} if
\[ \inf I_j <\inf I_{j+1} <\sup I_j<\sup I_{j+1}\]
for each $j=1,2,\ldots,m-1$. 
\ed
An infinite chain is also naturally defined by an infinite sequence of intervals such that every consecutive subsequence of 
$m$ intervals forms an $m$--chain. The parameter $m$ will sometimes be called the \emph{length}\index{chain length} of the chain.

\begin{lem}\label{l:m-chain-finite}
Fix $m\in\bZ_{>0}$.
If a collection $\VV$ of open intervals in $\bR$ has finite intersection multiplicity, then so does the collection
\[
\VV_m:=\{ J \mid J=J_1\cup\cdots \cup J_m\text{ for some }m\text{--chain }\{J_1,\ldots,J_m\}\text{ in }\VV\}.\]
\end{lem}
\bp
The proof is a simple induction, the base case $m=1$ being trivial. Suppose the conclusion holds for chains of length $m-1$ or less,
and assume for a contradiction that a point $x_0$ belongs to infinitely many distinct open intervals each of which is the union of an $m$--chain.
By the hypothesis, there exists an interval $J\in \VV$ containing $x_0$ and a positive integer $j\le m$ such that $J$ appears as the
$j^{th}$ term of infinitely many $m$--chains whose unions are all distinct. 
Since there exist only finitely many intervals in $\VV_{j-1}$ and in $\VV_{m-j}$ containing $\inf J$ and $\sup J$ respectively,
we obtain the desired contradiction. 
\ep

For a topological space $X$ and its covering $\UU$, we define 
 the \emph{covering length}\index{covering length} of a subset $A\sse X$ as the integer
\[
\CL_\UU(A):=\min\{ m \mid A\sse U_1\cup\cdots \cup U_m\text{ for some }U_i\in \UU\}.\]
Below we prove that under certain conditions, the covering length function for an interval is
topologically conjugate to the floor or the ceiling function.

\begin{lem}\label{l:cover}
Let $\mathcal{V}$ be a collection of bounded open intervals in $\bR$ with finite intersection multiplicity.
If $\VV$ covers $[0,\infty)$, then there exists a homeomorphism $h:[0,\infty)\longrightarrow[0,\infty)$ such that the following hold for all $x>0$:
\begin{itemize}
\item $\CL_{\mathcal{V}}[0,x)=\lceil h(x)\rceil$;
\item $\CL_\VV[0,x]=\lfloor h(x)\rfloor+1$.
\end{itemize}
\end{lem}

\bp
Consider the monotone function $f(z):=\CL_\VV[0,z)$, defined for $z>0$.
By Lemma~\ref{l:m-chain-finite}, for each $m>0$ we can define 
\[ z_m := \max \{ \sup J \mid J\text{ is the union of an }m\text{--chain from }\VV\text{ such that }0\in J\}<\infty.\]
Note that $f(z_1)=1$ and $f(z_m)\le m$. 

Assume inductively that $f(z_{m-1})=m-1$. By considering an interval from $\VV$ that contains $z_{m-1}$,  we see that
$f(z_{m-1}+\delta)=m$ for all sufficiently small $\delta>0$. It follows that \[f(z_m)\ge f(z_{m-1}+\delta)=m.\] That is, for all $m>0$ we have
\[f(z_{m-1},z_m]=\{m\}.\]
Choosing a homeomorphism $h\co [0,\infty)\longrightarrow[0,\infty)$ satisfying  $h(z_m)=m$, we obtain the first conclusion.
After noting that for all $z\in[z_{m-1},z_m)$ we have
\[
\CL_\VV[0,z]=m,\] the proof is complete.
\ep
\begin{rem}
One cannot drop the finite intersection multiplicity hypothesis. Indeed, if $\mathcal{V}$ consists of all intervals of the form $(-1,n)$ for $n\ge1$
then $\CL_{\mathcal{V}}[0,x)=1$ for all $x>0$.
\end{rem}

Let us now assume that $G\le\Homeo_+[0,1]$ is a group generated by a finite set $V$. The collection of open intervals
\[
\VV:=\bigcup_{v\in V} \pi_0\supp v\]
is a covering of $\supp G$ with finite intersection multiplicity.
For $x,y\in \Int I$, we define the \emph{covering distance}\index{covering distance}
\[
d_V(x,y) := \CL[x,y]\in\bN\cup\{\infty\}.\]

As mentioned above, the covering length is a measure of topological complexity. If
$G$ acts without global fixed points in $\Int I$, and if $f\in G$ is a compactly supported homeomorphism
(that is, the closure of $\supp f$ is a compact set in $\Int I$),
then $\CL({\suppc f})<\infty$. 
For two compactly supported homeomorphisms, the covering length of their supports allows one to 
quantitatively assert that one of them has a ``larger" support
than the other.

For $g\in G$, we define its \emph{syllable length (with respect to $V$)}\index{syllable length} as
\[
\|g\|_{\mathrm{syl}}:=\inf\{ \ell \mid g =  v_\ell^{n_\ell} \cdots v_2^{n_2}v_1^{n_1} \text{ for some }v_i\in V\text{ and }n_i\in\bZ\}.\]

The lemma below asserts that the covering length of the closure of an interval is the minimum syllable length of a 
group element that can move the interval off itself. Note that the former quantity depends on the action (dynamics), 
while the latter one is purely group theoretically defined. 
This observation may be considered as a manifestation of a recurring theme of the book, which is relating dynamical 
features of finitely generated group actions to group theoretic properties.

\begin{lem}\label{l:covering-translation-distance}
Let $G\le\Homeo_+[0,1]$ be a group with a finite generating set $V$,
and let $J\sse \supp G$ be a (compact or non-compact) nondegenerate interval.
Then we have that
\[\CL(\bar J)=\min\{ \|g\|_{\mathrm{syl}}\mid g\in G\text{ and }g(J)\cap J=\varnothing\}.\]
\end{lem}

\bp
In this proof, the covering distance and the syllable length are both defined by $V$.
We let $x<y$ be the endpoints of $J$.
Put $m:=\CL(\bar J)$.

There exists a minimal length $m$--chain 
\[ U_1, U_2, \ldots,U_m\]
(with the indices in the given order)
such that $U_i\in\pi_0\supp v_i$ for some $v_i\in V$
and such that \[\bigcup_i U_i\supset [x,y]=\bar J.\]
Depending on whether $v_i$  moves points in $U_i$ to the right or to the left, one can find integers $n_i\gg0$ or $n_i\ll0$ such that 
\[
v_i^{n_i}\cdots v_1^{n_1}(x)\in U_i\cap U_{i+1}\]
for each $i<m$. Since $y\in U_m$, we can pick $n_m$ such that 
$g_0:=v_m^{n_m}\cdots v_1^{n_1}$ moves $x$ further to the right of $y$.
Since $\|g_0\|_{\mathrm{syl}}\le m$, we have that
\[m\ge\min\{ \|g\|_{\mathrm{syl}}\mid g\in G\text{ and }g(J)\cap J=\varnothing\}.\]

For the opposite inequality, 
pick an arbitrary $g\in G$ satisfying $ g(J)\cap J=\varnothing$.
Without loss of generality, we may assume that 
\[ x<y \le g(x).\]
Writing $g=v^{n_\ell}_\ell \cdots v_1^{n_1}$ for some $v_i\in V$ and $n_i\in \bZ$,
we see that 
\[
[x,y]\sse [x,g(x)]\sse\bigcup_i \supp v_i.\]
This implies that $m=\CL(\bar J)\le \ell\le\|g\|_{\mathrm{syl}}$. This completes the proof.
\ep

We will make use of covering lengths crucially in Chapters~\ref{ch:slp} and~\ref{ch:optimal}.

\section{Generalities on chain groups}\label{sec:chain}

Let $\{J_1,\ldots,J_m\}$ be an $m$--chain, let \[J=\bigcup_{i=1}^m J_i\sse  I,\] and let $\{f_1,\ldots,f_m\}\sse \Homeo_+(J)$ such that
$\supp f_i\sse  J_i$. Following~\cite{KKL2019ASENS}, the group \[\form{  f_1,\ldots,f_m}=G\le \Homeo_+(J)\] is called a
\emph{pre-chain group}\index{pre-chain group}.
The group $G$ is called an \emph{$m$--chain group}\index{chain group} (or simply a \emph{chain group}) if for $i<m$ we have that 
\[f_{i+1}^{\pm1}(f_i^{\pm 1}(\inf J_{i+1}))\geq\sup J_i,\] where by the exponent $\pm 1$ we mean that this inequality holds for some choice
of exponents.
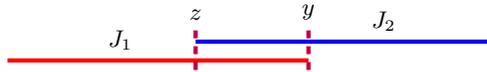
\begin{figure}[h!]
\centering
\begin{tikzpicture}[ultra thick,scale=.5]
\draw [red] (4,0) -- (12,0);
\draw [dashed,purple] (12,.8) -- (12,-.3);
\draw (9,.8) node [above] {\small $z$}; 
\draw [dashed,purple] (9,.8) -- (9,-.3);
\draw (12,.8) node [above] {\small $y$}; 
\draw (7,0) node [above] {\small $J_1$}; 
\draw [blue] (9,.5) -- (17,.5);
\draw (14,.5) node [above] {\small $J_2$}; 
\end{tikzpicture}
\caption{A chain of two intervals. The condition defining a chain group says that if $z=\inf J_2$ is the left endpoint of $J_2$
and $y=\sup J_1$ is the right endpoint of $J_1$, then
$f_2f_1(z)$ lies at least as far to the right as $y$.}
\label{f:coint}
\end{figure}

At first glance, this condition may look rather bizarre, though if $\supp f_i=J_i$ for all $i$ then the condition is dynamically
robust, in the sense that it holds after replacing each $f_i$ by $f_i^N$ for some $N\gg 0$. Indeed, this is obvious since (up to replacing
$f_i$ by its inverse) we have that $f_i^n(x)\to\sup J_i$ as $n\to\infty$, for all $x\in J_i$.
It turns out that in a chain group, the subgroup $\form{  f_i,f_{i+1}}$ is isomorphic to Thompson's group $F$. We will demystify
some of these definitions and apparent coincidences in this section.

We remark that if $f_1$ and $f_2$ are homeomorphisms generating a $2$--chain group, then $f_1$ and $f_2$ are 
\emph{crossed}\index{crossed interval}
in the sense of Section~\ref{sec:crossed} above.

\subsection{Two--chain groups and Thompson's group $F$}\label{ss:2-chain}

\emph{Thompson's group $F$}\index{Thompson's group $F$}
is a central object in the study of the group $\mathrm{PL}_+(I)$, the group of orientation preserving piecewise
linear homeomorphisms of the interval. The name \emph{piecewise linear}\index{piecewise linear}
is a bit of a misnomer, since elements of $\mathrm{PL}_+(I)$ are
locally described by affine functions. An element $f\in\PL_+(I)$ is required to be continuous,
and the derivative $f'$ is locally constant on the complement
of finitely many points in $I$, which we call the \emph{breakpoints}\index{breakpoint} of $f$.

The group $F$ is defined to be the full group of
orientation preserving piecewise linear homeomorphisms of $I=[0,1]$ such that all breakpoints of each element of $F$
are dyadic rationals, so that
all derivatives are powers of $2$. In this subsection we will discuss the following basic fact, which illustrates the ubiquity
with which Thompson's group is found within the study of groups of homeomorphisms of one-manifolds:

\begin{prop}\label{prop:two-chain}
Let $G$ be a two--chain group. Then $G$ is isomorphic to Thompson's group $F$.
\end{prop}

It is worth noting that there are no assumptions on $G$ other than it being a two--chain group. In particular, there are no restrictions
on the regularity of the generators of $G$; once a certain dynamical hypothesis is verified, the algebraic conclusion follows.

Since $(0,1)$ is homeomorphic to $\bR$, it is not surprising that the group $F$ is isomorphic to an easily describable subgroup of
$\Homeo_+(\bR)$. Indeed, it turns out that if we set
\[
a(x)=x+1\quad\text{and}\quad
b(x)=
\begin{cases}
 x&\text{ if }x\leq 0,\\
 2x&\text{ if }0<x< 1,\\
  x+1&\text{ if }1\leq x\\
\end{cases}
\]
then the subgroup of $\Homeo_+(\bR)$ generated by $a$ and $b$ is isomorphic to $F$. The subgroup $\form{  a,b}\le \Homeo_+(\bR)$
can be conjugated into $\mathrm{PL}_+(I)$ by an explicit homeomorphism $h$. See Figure~\ref{fig:gens}.

\begin{figure}
\centering
  \tikzstyle {a}=[postaction=decorate,decoration={%
    markings,%
    mark=at position .65 with {\arrow{stealth};}}]
{\begin{tikzpicture}[scale=.75]
\draw (0,0) -- (0,.5) -- (4,.5) -- (4,0) --cycle;
\draw [a] (2,.5) -- (1,0);
\draw [a] (3,.5) -- (2,0);
\draw [right] (0,.5) node [above] {\tiny $0$};
\draw [right] (4,.5) node [above] {\tiny $1$};
\draw  (2,0.5) node  [above] {\tiny $\frac12$} (3,.5) node [above] {\tiny $\frac34$};
\draw [below]  (1,0) node   {\tiny $\frac14$} (2,0) node {\tiny $\frac12$};
\draw [below]  (4,0) node   {\tiny $1$};
\draw [below]  (0,0) node   {\tiny $0$};
\draw [below]  (2,-1) node {\small $A$};
\end{tikzpicture}}
\quad\quad
{\begin{tikzpicture}[scale=.75]
\draw (0,0) -- (0,.5) -- (4,.5) -- (4,0) --cycle;
\draw [a] (2,.5) -- (2,0);
\draw [a]  (3,.5)--(2.5,0) ;
\draw [a]  (3.5,.5)--(3,0) ;
\draw [right] (0,.5) node [above] {\tiny $0$};
\draw [right] (4,.5) node [above] {\tiny $1$};
\draw [below]  (0,0) node   {\tiny $0$};
\draw [below]  (4,0) node   {\tiny $1$};
\draw  [below] (2,0) node {\tiny $\frac12$} (3,0) node  {\tiny $\frac34$} (2.5,0) node {\tiny $\frac58$};
\draw [above]  (2,.5) node {\tiny $\frac12$} (3,.5) node   {\tiny $\frac34$} (3.5,.5) node   {\tiny $\frac78$};
\draw [below]  (2,-1) node {\small $B$};
\end{tikzpicture}}
\quad\quad
{\begin{tikzpicture}[scale=.75]
\draw (0,0) -- (0,.5) -- (4,.5) -- (4,0) --cycle;
\foreach \i in {1,3,3.5} \draw [a] (\i,.5)--(\i,0);
\draw [right] (0,.5) node [above] {\tiny $0$};
\draw [right] (4,.5) node [above] {\tiny $1$};
\draw [below]  (0,0) node   {\tiny $-\infty$};
\draw [below]  (4,0) node   {\tiny $\infty$};
\draw  [below] (1,0) node {\tiny $-1$};
\draw  [below] (3,0) node  {\tiny $1$} (3.5,0) node  {\tiny $2$} ;
\draw [above]  
(1,.5) node {\tiny $\frac14$} (3,.5) node   {\tiny $\frac34$} (3.5,.5) node   {\tiny $\frac78$};
\draw [below]  (2,-1) node {\small $h$};
\end{tikzpicture}}
\caption{Elements $A,B\in\mathrm{PL}_+(I)$ generating Thompson's group $F$ and a conjugating homeomorphism $h$.}
\label{fig:gens}
\end{figure}
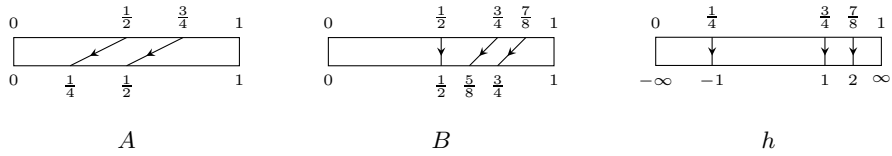

We have that \[a=hA^{-1}h^{-1}\quad \textrm{and}\quad  b=hB^{-1}h^{-1},\] and $F\cong\form{  A,B}$. In fact, one can say much more; 
it turns out that \[F\cong\form{  A,B\mid [AB^{-1}, A^{-1}BA],\, [AB^{-1},A^{-2}BA^2]}\] is a presentation for $F$. Since these claims about
presenting Thompson's group are not the primary concern of this book, we will omit proofs, directing the reader instead to some of the
standard references about Thompson's group, such as~\cite{CFP1996,BurilloBook}. The presentation for $F$
we have given here is essentially the only fact about
$F$ that we will not justify in this book.

Observe that there is a natural map \[D\colon \mathrm{PL}_+(I)\longrightarrow \bZ^2\] that takes a homeomorphism
and computes the right derivative at $0$ and the left derivative at $1$. Since $0$ and $1$ are global fixed points of the action of
$\mathrm{PL}_+(I)$, the chain rule implies that $D$ is in fact a homomorphism to the multiplicative group $\bR_+^2$.
Since elements of $\mathrm{PL}_+(I)$ only have finitely many
breakpoints, the derivative is continuous at all but finitely many points; it follows that breakpoints cannot accumulate,
and hence the homomorphism $D$ is well--defined since the derivative is constant in an open neighborhood of $0$ and $1$. For this reason,
the homomorphism $D$ is said to be computing derivatives of \emph{germs}\index{germ} of $\mathrm{PL}$ homeomorphisms; we will discuss
this concept in more generality shortly.

By considering
homeomorphisms that are the identity near one endpoint of $I$ but not the other, it is easy to see that $D$ is surjective. Moreover,
it is easy to see that the restriction of $D$ to $F$ is surjective. Since $F$ is two--generated, it follows that the commutator subgroup
$[F,F]$ consists of exactly the elements of $F$ that are the identity near $\{0,1\}$, which is to say homeomorphisms whose support
is compactly contained in $(0,1)$.

Note that since dyadic rational numbers are dense, there is an infinite sequence of elements
$\{f_n\}_{n\ge1}\sse  F$ such that \[{\suppc f_n}\sse  \supp f_{n+1}\] for all $n$, and such that if $J\sse  (0,1)$ is a compact
interval then $J\sse \supp f_n$ for $n\gg 0$. It follows that $[F,F]$ is not finitely generated.

The following fact has been encountered in another guise, and is an important feature of $F$.

\begin{lem}\label{lem:f-simple}
The group $F$ has trivial center, the commutator subgroup $F'=[F,F]$ is simple, and every proper quotient of $F$ is abelian.
\end{lem}
\begin{proof}
It is a straightforward exercise to show that $F'$ acts CO--transitively on $(0,1)$. It now follows from
Lemma~\ref{lem:co-trans} that $F''=[F',F']$ is nonabelian and simple, that the center of $F'$ is trivial, and that every proper quotient of $F'$ is
abelian.

Now, if $f\in F'$ is arbitrary, the ${\suppc f}$ is compactly contained in $(0,1)$. Thus, there is a compact interval $J\sse  (0,1)$ with
dyadic endpoints such that ${\suppc f}$ is contained in the interior of $J$. The subgroup $F_J\le \PL_+(J)$ consisting of homeomorphisms
with dyadic breakpoints and derivatives being powers of $2$ is isomorphic to $F$, since $J$ is homeomorphic to $[0,1]$ via an affine
homeomorphism that preserves dyadic rationals. Since $f\in F_J$ and since $f$ is compactly supported in $J$, we have that
$f\in F_J'=[F_J,F_J]$. Since $F_J$ consists of homeomorphisms compactly supported in $(0,1)$, we have that $F_J\le F'$. It follows that
if $f\in F'$ then $f\in F''$, so that $F'\le F''$ and hence $F'=F''$.

It is straightforward now to conclude that the center of $F$ is trivial. To see that every proper quotient of $F$ is abelian, let $K\le F$ be a
nontrivial normal subgroup. It follows that $[K,F]$ is a normal subgroup of $F'$, so that either $K$ commutes with every element of $F$
or $F'\le [K,F]$. The first possibility is ruled out since the center of $F$ is trivial, and the second possibility implies that $F/K$ is abelian.
\end{proof}

In the proof of Lemma~\ref{lem:f-simple}, we implicitly argued that $F'$ is perfect. This is part of a completely general phenomenon, which
we record now for future reference.

\begin{lem}\label{lem:chain-comm}
Let $G$ be a chain group. Then $G'$ is perfect and $G$ has trivial center.
\end{lem}

We recall that if $G$ acts by homeomorphisms on a Hausdorff topological space $X$ and if $x\in X$ is a global fixed point of $G$, we say
that the \emph{germ} of the $G$--action at $x$ is the equivalence relation on $G$ that says $g_1\sim g_2$ if $g_1$ and $g_2$ agree in
an open neighborhood of $x$.

\begin{proof}[Proof of Lemma~\ref{lem:chain-comm}]
Let $G$ be generated $\{f_1,\ldots,f_m\}$, where $\supp f_i=J_i$ and $\{J_1,\ldots,J_m\}$ forms an $m$--chain. Clearly we may assume
$m\geq 3$. Without loss of generality,
$\supp G=(0,1)$.

Let $u\in G$ be a nontrivial central element. We have that $\supp u$ is invariant under $G$ and accumulates at $0$ and $1$. Observe that
the germ of the action of $G$ at $0$ coincides with the cyclic group generated by $f_1$, and similarly the germ at $1$ coincides with the
cyclic group generated by $f_m$. Since the set $\supp u$ is invariant under $G$, we have that $\supp u=(0,1)$. Conjugating $f_1$ by a power
of $u$ results in a homeomorphism that does not commute with $f_1$, 
since $0=\inf J_1$ is fixed by $u$ and since $\sup J_1$ is not. This is a
contradiction, so the center of $G$ is trivial.

Now, we have that for $1\leq i\leq m-1$, the group $\form{  f_i,f_{i+1}}$ is isomorphic to $F$. Lemma~\ref{lem:f-simple} implies
that $F''=F'$, so \[[f_i,f_{i+1}]\in \form{  f_i,f_{i+1}}'=\form{  f_i,f_{i+1}}''\le G''.\] Since $G'$ is normally generated by elements of the
form $[f_i,f_{i+1}]$, we have that $G'=G''$ and so $G'$ is perfect.
\end{proof}

We can now prove the main result of this subsection.

\begin{proof}[Proof of Proposition~\ref{prop:two-chain}]
Let $\{J_1,J_2\}$ be a two--chain, as in Figure~\ref{f:coint}. Write $y=\sup J_1$ and $z=\inf J_2$. We write $G=\form{  f,g}$ for the
two--chain group in question, so that \[\supp f=J_1,\quad \supp g= J_2,\] and so that $g(f(z))\geq y$.

Observe that $k=g\circ f$ moves $J_2$ off of $J_1$, so that if $\supp h\sse J_2$ then \[(\supp khk^{-1})\cap J_1=\varnothing.\] In
particular, we have that \[[f^{-1},kg^{-1}k^{-1}]=[f^{-1},k^2g^{-1}k^{-2}]=1\] in the group $\form{  f,g}$.
Setting $f^{-1}g^{-1}=A$ and $g^{-1}=B$, we see that $AB^{-1}=f^{-1}$ and $A^{-1}=k$. It follows that $A$ and $B$ 
generate $G$, and we have
\[[f^{-1},kg^{-1}k^{-1}]=[AB^{-1},A^{-1}BA],\quad [f^{-1},k^2g^{-1}k^{-2}]=[AB^{-1},A^{-2}BA^2].\] 
It follows from the finite presentation of $F$ above that the map
\[\phi\colon F\longrightarrow G\] given by \[A\mapsto f^{-1}g^{-1}\quad \textrm{and}\quad B\mapsto g^{-1}\]
is a well-defined homomorphism of groups
that is clearly surjective.

It is immediate that $G$ is not abelian, since
on the one hand we have \[f(g(\inf J_2))=f(\inf J_2),\] but on the other hand since $f(\inf J_2)>\inf J_2$, we have that
\[g(f(\inf J_2))\neq f(\inf J_2).\] It follows immediately from Lemma~\ref{lem:f-simple} that $\phi$ is an isomorphism from $F$ to $G$.
\end{proof}

We remark that for $m$--chains with $m\geq 3$,
there is a similar result, except one replaces the group $F$ with the \emph{Higman--Thompson group}\index{Higman--Thompson group}
 $F_m$.
These groups form a natural generalization of Thompson's group $F$, and have the following uniform infinite presentations (though
they are all finitely presented):
\[ F_m\cong \form{  \{g_i\}_{i\ge1}\mid g_j^{g_i}=g_{j+m-1}\,\textrm{ for all } 0\leq i\leq j}.\] It turns out that the group
$F_2$ is just the group $F$, and the ambitious reader may try and prove this for themself.

Let $\{J_1,\ldots,J_m\}$ be an $m$--chain, and let $\{f_1,\ldots,f_m\}$ generate an $m$--chain group with $\supp f_i=J_i$.
If
\[f_m\cdots f_1(\inf J_2)\geq\sup J_{m-1},\] then we have that $\form{  f_1,\ldots,f_m}$ generates the $m^{th}$ Higman--Thompson
group $F_m$. Note that if $\{f_1,\ldots,f_m\}$ generate an $m$--chain group then for all sufficiently large $N\gg 0$, we have that 
$\{f_1^N,\ldots,f_m^N\}$ generate a copy of $F_m$. Thus, the group $F_m$ is the \emph{dynamical stabilization}\index{dynamical
stabilization} of all $m$--chain groups.

It is easy to see that the standard action of the Thompson's group $F\le\PL_+(I)$ on $(0,1)$ minimal, i.e. every orbit is dense 
under the action. 
We note a stronger fact regarding the \emph{positive $n$--transitivity}\index{positively $n$--transitive},
as defined in Definition~\ref{defn:n-trans}.
This fact is well-known; 
see~\cite[Theorem 3.1.2]{BurilloBook} for instance.
\begin{prop}\label{p:n-transitive}
For all $n\ge1$, the standard action of $F$ on $(0,1)$ is positively $n$--transitive on dyadic rationals.
In particular, every interval $U\sse (0,1)$ can be expanded arbitrarily close to $(0,1)$ by an element of $F$.
\end{prop}
The first statement means that given a pair of $n$--tuples
\[ x_1<x_2<\cdots<x_n\]
and \[y_1<y_2<\cdots<y_n\]
of dyadic rationals in $(0,1)$
there exists $g\in F$ such that $g(x_i)=y_i$ for all $i$. 
The second statement is an obvious consequence of the first by the density of dyadic rationals. 
\bp[Proof of Proposition~\ref{p:n-transitive}]
We follow an argument in~\cite{BurilloBook}. Set $x_0=y_0=0$ and $x_{n+1}=y_{n+1}=1$.
We claim that there exists a piecewise linear homeomorphism $g_i\co [x_i,x_{i+1}]\to[y_i,y_{i+1}]$ 
with breakpoints in $\bZ[1/2]$ and with slopes in $2^\bZ$. 
Then the map $g$ in the conclusion can be obtained by concatenating $g_0,\ldots,g_n$.

To see the claim, we fix $i$ and write \[x_{i+1}-x_i = p/2^m\quad\textrm{and}\quad y_{i+1}-y_i=q/2^n\quad \textrm{for some}\quad p,q,m,n\in\bZ_{>0}.\]
We may assume $p\le q$, for otherwise one can consider the inverse of $g_i$.
Write
\[
p/2^m =(p-1)/2^m + \sum_{i=1}^{q-p} 1/2^{m+i}+ 1/2^{m+q-p}.\]
So, $x_{i+1}-x_i$ can be partitioned into $q$ intervals whose lengths lie in $2^\bZ$. By mapping each of these intervals
to an interval of the form $[y_i+j/2^n,y_i+(j+1)/2^n]$, we obtain $g_i$ as claimed.
\ep

\subsection{Smooth realizations of Thompson's group}\label{ss:smooth}
This subsection follows Section 1.5.2 in~\cite{Navas2011}. The account here differs very little from Navas' exposition, and we include
it here because the smooth, locally dense realization of Thompson's group $F$ is critical in subsequent chapters.

Proposition~\ref{prop:two-chain} shows that an arbitrary two--chain group is automatically isomorphic to Thompson's group $F$. In particular,
it realizes many $C^{\infty}$ copies of $F$ inside of $\Homeo_+[0,1]$. While the smoothness of a copy of $F$ will be a critical tool in our
discussion of critical regularity later in this book, we will need to exercise some further control over the dynamical behavior of the
action of $F$ in a way that is not immediate as a consequence of smoothness.

Let $G\le \Diff_+^{\infty}(I)$ be a two--chain group such that $\supp G=(0,1)$.
If $G$ acts minimally on $(0,1)$ then Rubin's Theorem (Theorem~\ref{thm:rubin}) implies that this action is 
topologically conjugate to the standard piecewise linear action of $F$. However, there is little 
reason to believe that the action is minimal. If a given $C^{\infty}$ action of $F$ as a two--chain group
fails to be minimal, then one can apply a semi-conjugacy to obtain a minimal action (see Theorem~\ref{thm:min-set-r} below), though
in the process one may lose smoothness.

It turns out, however, that there is a copy of $F$
inside of $\Diff_+^{\infty}(I)$ that does act minimally on $(0,1)$, and therefore this action is topologically conjugate to the standard realization of
$F$ inside of $\PL_+(I)$. This fact was originally established by Ghys--Sergiescu~\cite{GS1987}. 
We remark that W.~P.~Thurston proved that $F$ is conjugate into the group $\Diff^{1+\mathrm{Lip}}_+(I)$, 
using Farey sequences and maps that are piecewise defined to be in $\PSL(2,\bZ)$ (see Section 1.5.1 of~\cite{Navas2011}).

Since we will use a $C^{\infty}$
action of $F$ on $I$ that is minimal on $(0,1)$ in an essential way, we will explain the construction of such an action, following
Navas' exposition of Ghys--Sergiescu's work in Section 1.5.2 of~\cite{Navas2011} (cf.~Subsection~\ref{sss:tsuboi}).
The idea is to associate representations of $F$ to homeomorphisms of $\bR$ that satisfy some
relatively mild hypotheses.

For the rest of this subsection, we will write $h\in\Homeo_+(\bR)$ for a homeomorphism that satisfies the following conditions.
\begin{enumerate}[(A)]
\item
We have $h(0)=0$.
\item
$h(x+1)=h(x)+2$ for all $x\in \bR$.
\item $|h(x)-h(y)|>|x-y|$ for all $x,y\in\bR$.
\end{enumerate}
An example of a homeomorphism of $\bR$ that satisfies the preceding three conditions is the map $x\mapsto 2x$.
For each such homeomorphism $h$, we will build an action \[\varphi_h\colon F\longrightarrow\Homeo_+(\bR)\] of $F$ on the real line.

We will write $\Aff_+^2(\bR)$ for the group or orientation preserving affine homeomorphisms that preserve the set of dyadic rational numbers.
It will be helpful to write $\bQ_2\le \bR$ for the additive subgroup of $\bR$ consisting of all dyadic rational numbers. 
For $y\in\bR$, 
we define $\lambda,\tau_y\in \Homeo_+(\bR)$ by
\begin{align*}
\lambda(x)&:=2x,\\
\tau_y(x)&:=x+y.\end{align*}

As in Section~\ref{ss:bs-smooth}, 
we have several equivalent expressions for $\Aff_+^2(\bR)$
\[
\Aff_+^2(\bR)=2^\bZ\times \bQ_2=
\left\{\begin{pmatrix} 2^n&p/2^q\\ 0 & 1\end{pmatrix}\middle\vert n,p,q\in\bZ\right\}=\form{a,e\mid aea^{-1}=e^2}.\]
The verification of the equivalence is trivial using the correspondence
\[
\lambda=\begin{pmatrix} 2&0\\ 0 & 1\end{pmatrix}=a,\qquad
\tau_1=\begin{pmatrix} 1&1\\ 0 & 1\end{pmatrix}=e.\]
The presentation for the group $\form{a,e}$ above is called the \emph{Baumslag--Solitar group} of type $(1,2)$, and denoted as $\BS(1,2)$; see Section~\ref{sec:bs} for an extensive discussion on this group. 
The following is transparent almost only by navigating the definitions.

\begin{lem}\label{lem:aff-homo}
There uniquely  exists a homeomorphism $\eta_h\in\Homeo_+(\bR)$ 
satisfying the following.
\be[(i)]
\item $\eta_h(0)=0$;
\item $\eta_h\tau_1=\tau_1\eta_h$;
\item $\eta_h\lambda=h\eta_h$;
\item $\eta_h(\bQ_2)$ is dense.
\ee
\end{lem}
\bp
Let us define a map
\[\Psi_h\colon \Aff_+^2(\bR)\longrightarrow\Homeo_+(\bR)\] 
by the formula
$\Psi_h(\lambda)=h$ and that $\Psi_h(\tau_1)=\tau_1$.
To see that $\Psi_h$ is well-defined and a homomorphism, it suffices to see
that the relation of $\BS(1,2)\cong \Aff_+^2(\bR)$ are preserved. Indeed, we note
for each $x\in \bR$ that
\[h\tau_1h^{-1}(x)=h(h^{-1}(x)+1)=x+2=\tau_1^2(x).\]

Let us now define
$\eta_h (g(0))=\Psi_h(g)(0)$ for all $g\in \Aff_+^2(\bR)$, which is easily seen to be well-defined and makes the diagram below commutes for all $g\in \Aff_+^2(\bR)$:
\[\begin{tikzcd}
\bQ_2\arrow{r}{g} \arrow{d}{\eta_h } &
\bQ_2\arrow{d}{\eta_h }\\
\bR\arrow{r}{\Psi_h(g)}& \bR
\end{tikzcd}
\]

The map $\eta_h$ preserves the order. Indeed, let us pick dyadic rationals $p/2^q$ and $p'/2^q$ for some $p<p'$.
We have that
\[
\eta_h (p/2^q)=\eta_h (\lambda^{-q}\tau_1^p\lambda^q(0))
=h^{-q}\tau_1^p h^q(0)=h^{-q}(p)
< h^{-q}(p')=\eta_h (p'/2^q), \]
implying that the order is preserved. We also see from the above computation that $\eta_h \restriction_\bZ$ is the identity.

We claim that the set $\eta_h(\bQ_2)$ is dense in $\bR$. Once this claim is established, we can simply extend $\eta_h$ to the closure of $\eta_h(\bQ_2)$ continuously. The uniqueness of $\eta_h$ will also follow since the condition
$\eta_h (g(0))=\Psi_h(g)(0)$ is necessary from that $\eta_h$ conjugates $\lambda$ to $h$.

To prove the claim, we assume that the set
\[
C:=\overline{\eta_h(\bQ_2)}=\overline{\{ h^{-q}(p)\mid p,q\in\bZ\}}\]
is proper in $\bR$.  
The set $C$ is also $1$--periodic since $\eta_h$ commutes with $\tau_1$.
In particular, we can choose a component $J$ of $\bR\setminus A$ with a maximal length.
On the other hand, we have that $h(C)=C$ since the set $\{h^{-q}(p)\mid p,q\in\bZ\}$ is $h$--invariant. It follows that $h(J)$ is also a component of $\bR\setminus C$.
By the condition (C) above we have
\[|h(J)|>|J|,\]
which contradicts the maximality.
\ep

Extending the above definition, we continue to denote by
\[\Psi_h\co \Homeo_+(\bR)\to\Homeo_+(\bR)\]
the conjugation by $\eta_h$. 
Namely,
\[
\Psi_h(g)=\eta_h g\eta_h^{-1}\]
for $g\in \Homeo_+(\bR)$.
We consider the restriction
\[
\varphi_h:=\Psi_h\restriction_F.\]
In other words, we let $\varphi_h$ be the topological conjugacy of the standard action of $F$ on $[0,1]$ by $\eta_h$, extending the map by the identity outside $[0,1]$.

Observe that if $h(x)=2x$ then $\varphi_h(F)$ is just the usual action of $F$ by piecewise linear homeomorphisms on $F$. 
 Indeed, let
$F=\form{  a,b}$ be the usual generators of $F$ (see Subsection~\ref{ss:2-chain}, especially Proposition~\ref{prop:two-chain}).
We have that $\varphi_h(a)=\Psi_h(\tau_1)=\tau_1=a$. Similarly we can check that 
$\varphi_h(b)$ agrees with the identity on $(-\infty,0)$ and with $\tau_1$ on $(1,\infty)$, and agrees with $h$ on $(0,1)$. Thus, $\varphi_h(b)=b$.

We recall the meaning of $C^r$--tangency to the identity from Definition~\ref{defn:tangent-id}.

\begin{thm}[\cite{GS1987}]\label{thm:phi-smooth}
If $h\in\Diff^r_+(\bR)$ satisfies the conditions (A), (B) and (C) above,
and if 
\addtocounter{enumi}{3}
\be[(D)]
\item
$h$ is $C^r$--tangent to the identity at $0$
for some $r\in{\bZ_{>0}}\cup\{\infty\}$,
\ee
then we have
\[
\varphi_h(F)\le\Diff_+^r(\bR).\]
\end{thm}
\begin{proof}
Let $\{a_n\}_{n\in\bZ}$ be an increasing sequence containing breakpoints of $g$, and let $\{b_n:=\eta_h(a_n)\}_{n\in\bZ}$. We write $f_n\in\Aff_+^2(\bR)$ describing $g$ on the interval $[a_n,a_{n+1}]$. 

Since $h\in\Diff^r_+(\bR)$, it is immediate from the definitions that 
\[\Psi_h(f_n)\in\form{\tau_1,h}\le \Diff^r_+(\bR)\]
 for all $n$. 
To prove that $\varphi_h(g)\in\Diff^r_+(\bR)$, we need
to check smoothness at each of the points $\{b_n\}_{n\in\bZ}$. Evidently, this is equivalent to the statement of equality of derivatives
\[\Psi_h(f_{n-1})^{(i)}(b_n)=\Psi_h(f_{n})^{(i)}(b_n),\] for all $i\leq r$.

Since $f_{n-1}(a_n)=f_n(a_n)$ by the continuity of $g$, we have that the
affine map $\tau_{a_n}^{-1}f_{n-1}^{-1}f_n\tau_{a_n}$ fixes the origin, and is hence given by multiplication by a power of two, say $2^p$.
Computing, we have that \[\Psi_h(\tau_{a_n}^{-1}f_{n-1}^{-1}f_n\tau_{a_n})=h^p.\] Taylor's Theorem says that $h$ is given
near the origin by the function $x+\omega(x)$, where all derivatives of $\omega$ at $0$ vanish up
to order $r$.  It follows that \[\Psi_h(\tau_{a_n}^{-1}f_{n-1}^{-1}f_n\tau_{a_n})^{(i)}(0)=
\begin{cases}
 1&\text{ if }i= 1,\\
 0&\text{ if }2\leq i\leq r.\\
\end{cases}\] 
It follows that the Taylor expansions of $\Psi_h(f_{n-1})$ and $\Psi_h(f_{n})$ agree up to order $r$ at $b_n$.
\end{proof}

\begin{cor}\label{cor:ghys-serg}
There is an embedding $F\longrightarrow\Diff_+^{\infty}(I)$ such that the action of $F$
on $(0,1)$ is minimal.
\end{cor}
\begin{proof}
We embed $F$ into $\PL_+(I)$ in the usual way, and extend by the identity to all of $\bR$. This realizes $F$ as a group of homeomorphisms
of $\bR$ and are locally described by elements of $\Aff_+^2(\bR)$, which are supported on $[0,1]$. Applying $\varphi_h$
for $h\in\Diff_+^{\infty}(\bR)$ 
satisfying the conditions (A) through (D) above,
we obtain a group of diffeomorphisms
of $\bR$ that are the identity outside of a compact interval. Restricting to this interval, we obtain an embedding \[F\longrightarrow\Diff_+^{\infty}
(I)\] such that each element of $F$ has derivatives of all orders that coincide with the identity at $0$ and $1$. The fact that this action of
$F$ is minimal on $(0,1)$ follows from the fact that the usual action of $F$ by piecewise linear homeomorphisms if minimal, and the
fact that $\varphi_h(F)$ is topologically conjugate to the usual action.
\end{proof}

\begin{rem}
In Lemma~\ref{lem:aff-homo}, if the condition
\be[(A)]
\setcounter{enumi}{2}
\item $|h(x)-h(y)|>|x-y|$ for all $x,y\in\bR$,
\ee
is dropped then the set $\eta_h(\bQ_2)$ may not be dense
and so, $\eta_h\co\bQ_2\to\bR$ may not extend to a continuous map of the real line.
We can still extend $\eta_h^{-1}$ to a periodic monotone increasing continuous map $\eta_h^{-1}\co\bR\to\bR$ in this case.
A smooth action 
\[
\Psi_h\co F\to\Diff_{[0,1]}^\infty(\bR)\]
can also still be defined in a piecewise manner, and moreover,
$\eta_h^{-1}$ semi-conjugates this action to the standard piecewise linear action.
This way, one can construct a faithful smooth exceptional action of $F$ on $[0,1]$.
The  construction of this subsection works equally well for the Thompson's group $G$, which piecewise linearly acts on the circle 
with dyadic breakpoints and with slopes powers of two~\cite{GS1987}.
\end{rem}

\subsection{Subgroups of chain groups}\label{ss:sub-chain}

In light of Proposition~\ref{prop:two-chain}, one might be tempted to believe that the structure of chain groups is similar to that of
$F$, with perhaps a na\"ive guess being that a general $m$--chain group is isomorphic to the Higman--Thompson group $F_m$. Though this
first guess is evidenced by the dynamical stabilization for chain groups we have already observed,
it is not right. In fact, the diversity of subgroups
of chain groups is so wild as to furnish continuum many isomorphism types of $m$--chain groups for all $m\geq 3$.

\begin{thm}[cf.~\cite{KKL2019ASENS}, Theorem 1.4]\label{thm:chain-subgp}
Let \[\gam=\form{  f_1,\ldots,f_n}\le \Homeo_+(\bR)\] be an arbitrary finitely generated subgroup.
\begin{enumerate}[(1)]
\item
The group $\Gamma$ is a subgroup of an $(n+2)$--chain group.
\item
If $\supp f_1$ has finitely many components then $\Gamma$ can be realized as a subgroup of an $(n+1)$--chain group.
\end{enumerate}
\end{thm}

Theorem~\ref{thm:chain-subgp} can be contrasted with the following fact about $F$, one that was originally due to Brin and Squier~\cite{BS1985},
and which severely limits the possible subgroups of $F$.

\begin{thm}\label{thm:f-subgp}
The group $F$ contains no nonabelian free subgroups and obeys no nontrivial law.
\end{thm}


The proof of Theorem~\ref{thm:f-subgp} is not difficult once one has made the correct observations, and the idea behind this theorem can be traced all the way back to at least the Zassenhaus Lemma (see~\cite{Raghunathan1972}). The reader may also compare with Proposition~\ref{prop:compact-conn}.
The same ideas will recur later in this book in a stronger form. 
We will give a proof here based on~\cite{BS1985}, more in this line of ideas.
In the local context of Theorem~\ref{thm:f-subgp} only, we will use the notation $\bF_2$ for the free group
of rank two, in order to distinguish it from Thompson's group.

\begin{proof}[Proof of Theorem~\ref{thm:f-subgp} (Sketch)]
Observing that $F\le\Homeo_+[0,1]$ is locally moving (Definition~\ref{d:locally-moving}), we immediately deduce from 
Theorem~\ref{thm:lawless} that $F$ obeys no law. We even note from the paragraph following the corollary that $F'$ obeys no law either. 

Alternatively, we can repeat the ideas in the proof of
Proposition~\ref{prop:z2z-real} to show that $F$ obeys no law. 
Briefly, using that argument for every nontrivial element 
of the free group $w\in\bF_2$, we can construct a finite sequence of compactly supported homeomorphisms
(i.e.~bumps) which together witness the fact that $w$ is not a law obeyed in $\Homeo_+(\bR)$.
It is not difficult to realize these bumps as elements of the commutator subgroup $F'$, which shows that $F'$ (and therefore $F$) obeys no law.

To show that $F$ contains no copy of $\bF_2$, let $f,g\in F$ be arbitrary elements that generate a subgroup $G\le F$.
Write $U=\supp f\cup\supp g$. Because $f$ and $g$ are piecewise linear, we have that $U$ consists of finitely many components,
and every element of the derived subgroup $G'\le G$ is the identity near $\partial U$.

Now, let $h\in G'$, and suppose that $V=\supp h$ meets $k>0$ components of $U$, say $\{U_1,\ldots,U_k\}$. Writing $V_1=V\cap U_1$,
there is a nontrivial elements $g_1\in G$ such that $g_1(V_1)\cap V_1=\varnothing$. It follows that $\supp [h,h^{g_1}]$ meets
$U$ in at most $k-1$ components. It follows that we may recursively construct a sequence of nested commutators which eventually
meet $U$ in no components, which is to say that the commutator becomes the identity. Since this commutator can easily be arranged to
be a reduced word in the generators $f$ and $g$ of $G$, this proves that $\form{  f,g}\ncong\bF_2$.
\end{proof}

We will postpone the proof of Theorem~\ref{thm:chain-subgp} until later in this chapter, as it will be useful to develop some further
dynamical and algebraic tools first.

\section{Simplicity}

Observe that since chain groups are natural generalizations of Thompson's group $F$, it is reasonable to suspect that they might share
some dynamical properties, even if their algebraic properties are not exactly analogous, as we have seen from 
Theorems~\ref{thm:chain-subgp}
and Theorem~\ref{thm:f-subgp}. Since the usual realization of $F\le \PL_+(I)$ is CO--transitive when restricted to $F'$, we were able to prove
that $F'$ is simple (see Lemma~\ref{lem:f-simple}).

It is not immediate, however, that arbitrary chain groups have simple commutator subgroups. Indeed, the simplicity of $F'$ follows from
the fact that $F$ has a particularly nice realization. There are many other realizations of $F$ that are not CO--transitive. Indeed,
consider the usual action of $F$ on $[0,1]$, and blow up the orbit of an arbitrary point $x\in (0,1)$, as we did to construct continuous
Denjoy counterexamples. That is, we glue in intervals of finite total length along points in the orbit of $x$ in order to get a new action
of $F$ on $I$ that has a wandering set $J$.

Let $y\in (0,1)\setminus J$, and let $K$ be a small compact neighborhood of $y$. It is not the
case that for all open sets $\varnothing\neq U\sse  (0,1)$ there is an element of $F$ that sends $K$ into $U$. Indeed, if $U$ is a connected
component of $J$ and $f\in F$, then $f(K)\cap U\neq\varnothing$ implies that $f(K)\cap\partial U\neq\varnothing$. Thus, such actions of
$F$ are clearly not CO--transitive. If the original orbit that was blown up avoided the endpoints of the intervals defining a two--chain group,
then this blown up copy of $F$ can be realized as a two--chain group. Through these observations
we see that the simplicity of $F$ really depended on
a minimal action of $F$ on $(0,1)$.

It is not surprising then that whether or not one can predict the simplicity of the commutator subgroup $G'$ of 
a chain group $G$ depends on the
particular realization of the chain group that is used.

In order to proceed, we need to analyze the general behavior of minimal sets for group actions
on the interval $I$, which unfortunately does not have as nice an answer as for
group actions on the circle (see Theorem~\ref{thm:minimal-set}).
The reason that the argument for the circle does not generalize verbatim to the interval is that $(0,1)$ is not compact. Since the circle
case relies in an essential way on compactness (via the finite intersection property), another argument is required.
In the end, one does find that for every finitely generated subgroup $G\le \Homeo_+[0,1]$,
there is a closed subset $C\sse  (0,1)\cong\bR$ on which
$G$ acts minimally, though unfortunately this set may not be unique. Precisely, we have the following.

\begin{thm}[\cite{Navas2011}, Proposition 2.1.12]\label{thm:min-set-r}
Let $G\le \Homeo_+[0,1]$ be a finitely generated group. Then there is a nonempty, closed, invariant subset $C\sse  (0,1)$ 
on which $G$ acts minimally.
Moreover, let $C'$ denote the derived subset of $C$ and let $\partial C$ be the boundary of $C$. Then we have
exactly one of the following conclusions:
\begin{enumerate}[(1)]
\item
We have $C'=\varnothing$. Then, $C$ is a discrete subset of $(0,1)$. If $C$ is finite then it consists of global fixed points of $G$. Otherwise,
$C$ is a bi-infinite sequence of points that accumulates exactly at $\{0,1\}$.
\item
We have $\partial C=\varnothing$. Then $C=(0,1)$ and $G$ acts minimally on $(0,1)$.
\item
We have $C'=\partial C=C$. Then, if $J\sse  (0,1)$ is an open interval that is compactly contained in $(0,1)$,
we have that $\overline{J\cap C}$ is a Cantor set. Moreover, there is a continuous, monotone function $h\colon (0,1)\to (0,1)$
and a quotient \[q\colon G\longrightarrow Q\] such that for all $g\in G$, we have \[q(g)\circ h=h\circ g.\] The map $h$ is surjective when
restricted to $C$ and constant on components of $(0,1)\setminus C$, and so $Q$ acts minimally on $(0,1)$.
\end{enumerate}
\end{thm}

As in the discussion surrounding Theorem~\ref{thm:irr-rot-semi}, the map $h$ is called a \emph{semi-conjugacy}\index{semi-conjugacy}.
As in the case of the circle, we say that group actions as in the last conclusion are \emph{exceptional}\index{exceptional minimal set},
and that $C$ is an \emph{exceptional
minimal set}.

\begin{proof}[Proof of Theorem~\ref{thm:min-set-r}]
Let $G$ be generated by a finite set $S$, which is assumed to be symmetric (i.e.~$S$ is closed under taking inverses). If there exists a
global fixed point for $G$ then clearly the conclusion of the theorem is true. Note that the global fixed set of $G$ can have nonempty
interior, and every point in the fixed point set is a closed, minimal, invariant subset. This illustrates the fact that the set $C$ may not be unique.

We assume for the rest of the proof that $G$ has no global fixed points. Let $x\in (0,1)$ is arbitrary, and let $\OO$ be the $G$--orbit of $x$.
We have that $\sup\OO$ and $\inf\OO$ are both global fixed points of $G$ and must therefore
coincide with the set $\{0,1\}$. Indeed, let $y=\sup\OO$,
and let $g\in G$. If $g(y)\neq y$ the $g(y)<y$, which implies that $g^{-1}(y)>y$, a contradiction. This also shows that if $C$ is discrete then
it consists of a bi-infinite sequence of points accumulating at $0$ and $1$.

To apply Zorn's Lemma as in Theorem~\ref{thm:minimal-set}, we build a compact subset $J$ of $(0,1)$ such that the closure of every orbit
meets $J$.
To this end,
if $x\in (0,1)$, let \[y=\max \{s(x)\mid s\in S\}.\] We claim that an arbitrary $G$--orbit $\OO$ must meet the interval $[x,y]$. Indeed, we have
that \[\{\sup\OO,\inf\OO\}=\{0,1\},\] and so there are points $w<x<y<z$ such that $w,z\in\OO$. Let $g(w)=z$. Writing $g$ as a product of
elements of $G$, we have \[g=s_k\cdots s_1.\] We let $\ell$ be the last index so that \[v=s_{\ell}\cdots s_1(w)<x.\] Then we have that 
\[x<s_{\ell+1}(v)<y,\] since $s_{\ell+1}(x)\leq y$ and $s_{\ell}$ is order preserving. Note that this is where we use the fact that $S$ is a
symmetric generating set.

Finally, let $J=[x,y]$. We write $\PP$ for the set of nonempty, $G$--invariant subsets of $(0,1)$, and for $C_1,C_2\in\PP$, write $C_1\geq C_2$
if \[C_1\cap J\sse C_2\cap J.\] If $C\in\PP$ is arbitrary, we have that $C\cap J\neq\varnothing$. The finite intersection property again
implies that chains have upper bounds in $\PP$, and hence Zorn's Lemma furnishes a maximal element, on which $G$ must act minimally.
Indeed, otherwise the closure of a $G$--orbit would be a proper closed invariant subset and would violate maximality.

We can now analyze the structure of $C$. Suppose first that $C$ contains a nonempty open interval $K$. The minimality of the $G$--action on
$C$ implies that if $x\in C$ is arbitrary, then there is an element $g\in G$ such that $g(x)\in K$, and so that $G$--invariance of $C$ implies that
a neighborhood of $x$ is contained in $C$. In particular, $\partial C=\varnothing$ and $C= (0,1)$, and so $G$ acts minimally on $(0,1)$.

Thus if $C\neq (0,1)$ then $C$ is totally disconnected. The minimality of the $G$--action on $C$ again implies that either every point of $C$
is isolated, or every point of $C$ is an accumulation point. If every point of $C$ is isolated then $C$ is discrete and accumulates only
at $0$ and $1$. If one point of $C$ is an accumulation point then $C$ has no isolated points and hence is perfect. It follows then that
if $J$ is a nonempty open interval that is compactly contained in $(0,1)$ then $\overline{J\cap C}$ is a Cantor set. 

In this last case, let $\mu$ be a $\sigma$--finite, regular, nonatomic Borel measure that is supported on $C$. Then the map
\[h\colon (0,1)\longrightarrow \bR\] given by \[h(x)=\frac{1}{2}+\int_{1/2}^x\,d\mu\] is a continuous map that is a surjection when
restricted to $C$, and which is constant on components of $(0,1)\setminus C$. For $y=h(x)$, we set $\overline{g}(y)=h(g(x))$.
This is well-defined, since $h$ is at most two--to--one, and is one--to--one outside of the set of boundary points of components of
$(0,1)\setminus C$. If $J$ is such a component then so if $g(J)$, and so $G$ preserves fibers of $h$. This establishes the well-definedness
of the action of $\overline{g}$.
Finally, we set $q(g)=\overline{g}$, which defines a quotient $Q$ of $G$.
\end{proof}

It is not difficult to show that $Q$ might be a proper quotient of $G$. Indeed, choose a
finitely generated group $G$ and a $G$--action on $(0,1)$ with an exceptional
minimal set $C$, and let $J$ be a component of $(0,1)\setminus C$. Choose an arbitrary finitely generated minimal $H$--action on $J$,
and propagate it to the $G$--orbit of $J$ be conjugation by $G$. Then $\form{  G,H}$ act on $(0,1)$, and $C$ is still an exceptional
minimal set. Building $h$ as in the proof of Theorem~\ref{thm:min-set-r} furnishes a quotient $Q\cong G$ of $\form{  G,H}$ that
is a proper quotient.

We now return to the subject of chain groups. First, a preliminary definition. Let $X$ be a Hausdorff topological space, and let
$G\le \Homeo(X)$ be a subgroup.
Following~\cite{KKL2019ASENS}, we say that $G$ is \emph{locally CO--transitive}\index{locally CO--transitive}
if for each proper compact $A\sse  X$, there exists
a point $x\in X$ such that for an arbitrary open neighborhood $U$ of $x$, we have $g(A)\sse  U$ for some $g\in G$.

\begin{lem}[\cite{KKL2019ASENS}, Lemma 3.6]\label{lem:chain-arb-dyn}
Let $G$ be an $m$--chain group, with $m\geq 2$.
\begin{enumerate}[(1)]
\item
There is a point $x\in\supp G$ such that every $G$--orbit accumulates at $x$.
\item
For $g\in G$ and $A\sse \supp G$ nonempty and compact, there is an element $h\in G'$ such that $g$ and $h$ agree as functions on $A$.
\item
Every $G$--orbit is also a $G'$--orbit.
\item
The action of $G'$ is locally CO--transitive.
\end{enumerate}
\end{lem}
\begin{proof}
Write $\{f_1,\ldots,f_m\}$ for generators of $G$, with $\supp f_i=J_i$. Without loss of generality, we may assume that
\[\bigcup_{i=1}^m J_i=(0,1).\]

Clearly, every $G$--orbit accumulates at $x=\sup J_1\in (0,1)$. This establishes the first part of the lemma.

For the second part of the lemma, write \[g=s_k\cdots s_1,\] where for each $i$ we have \[s_i\in\{f_1^{\pm1},\ldots,f_m^{\pm1}\},\]
and let $A\sse  (0,1)$ be nonempty and compact. Observe that there exist nonempty open intervals \[U_1=(0,a)\sse  J_1\quad
\textrm{and}\quad
U_2=(b,1)\sse  J_m\] such that for $1\leq i\leq k$, we have \[(s_i\cdots s_1(A))\cap (U_1\cup U_2)=\varnothing.\] It is a straightforward
exercise, using the connectedness of $\supp G$, to find elements $g_i\in G$ such that $g_i(J_i)\sse  U_1\cup U_2$.
Let \[h=\left(\prod_{i=1}^k g_is_i^{-1}g_i^{-1}\right)g.\] It is immediate that $h\in G'$, and an easy calculation shows that $h$ agrees
with $g$ on $A$.

The third part of the lemma follows from the second, letting $A$ be a single point.
For the fourth part, since local transitivity of $G'$ concerns the behavior of $G'$ on compact subsets of $(0,1)$, the second part of the
lemma implies that we may as well prove that $G$ is locally CO--transitive. We let $x=\sup J_1$ as in the first part, and we fix an open
neighborhood $U$ of $x$. Again, if $A\sse  (0,1)$ is nonempty and compact then we may find an element of
$g\in G$ such that $\sup A\in J_1$, so that we have $g(A)\sse  J_1$. Applying a sufficiently large power of $f_1$ sends $g(A)$ into
$U$, thus completing the proof.
\end{proof}

As a consequence of the preceding discussion, we have the following basic structural result about chain groups.

\begin{thm}[\cite{KKL2019ASENS}, Theorem 3.7]\label{thm:chain-dichotomy}
Let $G$ be an $m$--chain group, for $m\geq 2$. Then exactly one of the following conclusions holds.
\begin{enumerate}[(1)]
\item
The action of $G$ is minimal.
\item
The action of $G$ admits a unique exceptional minimal invariant set.
\end{enumerate}
In the first of these conclusions, we have that the commutator subgroup $G'$ is simple. In the second of these conclusions,
we have that $G$ surjects onto an $m$--chain group that acts minimally.
\end{thm}
\begin{proof}
As before, we assume that $\supp G=(0,1)$, and we write
$\{f_1,\ldots,f_m\}$ for generators of $G$, with $\supp f_i=J_i$. Suppose first that the action of $G$ is not minimal.
Then since $G$ obviously has no global fixed points,
Theorem~\ref{thm:min-set-r} implies that $G$ admits an (\emph{a priori} non--unique) exceptional minimal set $C$.

In the case of an exceptional minimal set, we have from Lemma~\ref{lem:chain-arb-dyn}, there is a point $x\in (0,1)$ that meets the
closure of an arbitrary orbit $\OO$ of $G$. It follows that the closure of the $G$--orbit of $x$ is contained in every
nonempty, closed, $G$--invariant
subset of $(0,1)$. It follows that the exceptional minimal set $C$ is unique. Theorem~\ref{thm:min-set-r} furnishes a quotient $Q$ of
$G$ acting on $(0,1)$ minimally. To complete the analysis of this case, it suffices to show that $Q$ is again an $m$--chain group.
Let $J_i$ and $J_{i+1}$ be two consecutive intervals in the $m$--chain defining $G$. We have that 
\[\inf J_i<\inf J_{i+1}<\sup J_i<\sup J_{i+1},\]
by definition. Write \[K_1=(\inf J_i,\inf J_{i+1}),\quad K_2=(\inf J_{i+1},\sup J_i),\quad K_3=(\sup J_i,\sup J_{i+1}).\]
It is a straightforward exercise to show that $C\cap K_{\ell}$ is infinite for $\ell\in \{1,2,3\}$. Since the map $h$ semi-conjugating the action
of $G$ to the action of $Q$ is at most two--to--one, we have that $h(J_i)$ and $h(J_{i+1})$ form a two--chain of nondegenerate intervals.
It follows that $Q$ is indeed an $m$--chain group.

If $G$ is a minimal chain group, then Lemma~\ref{lem:chain-arb-dyn} implies that the action of $G'$ on $(0,1)$ is also minimal and
CO--transitive. Indeed, the CO--transitivity can be shown as follows. If $A\sse  (0,1)$ is compact and $x=\sup J_1$, then for each $a<x$,
there is a $g_a\in G$ such that $g_a(A)\sse  (a,x)$. Since $G$ act minimally, for a 
nonempty open $U$, there is an $h$ such that $h(x)\in U$.
Continuity implies that there is an $a_0$ such that $h((a_0,x])\sse  U$. Then, $h\cdot g_{a_0}$ sends $A$ into $U$.

Now, elements of $G'$ have the property that their supports are compactly contained in $(0,1)$, since the germs of the $G$--action at
$0$ and $1$ are abelian groups. By Lemma~\ref{lem:chain-comm}, we have that $G'$ is perfect. Lemma~\ref{lem:co-trans} now implies
that $G'$ is simple.
\end{proof}

Theorem~\ref{thm:chain-dichotomy} is indeed a true dichotomy, at least for $m$--chain groups with $m\geq 3$. Precisely, we have the
following result, which we will not prove here but which the reader may try and prove for themself.

\begin{prop}[See~\cite{KKL2019ASENS}, Proposition 4.8]\label{prop:non-simple}
For $m\geq 3$, there exists an $m$--chain group with a non--simple commutator group.
\end{prop}

Proposition~\ref{prop:non-simple} is obviously false for $m=2$. The construction in the proposition is carried out 
by blowing up an orbit, much like in the
construction of a continuous Denjoy counterexample, and ``inserting" a group action in the wandering set. The construction fails when $m=2$
because one needs to build a nontrivial element in the chain group that fixes the whole wandering set, and this is not possible to achieve
with $F$.

\section{The chain group trick and the rank trick}

We now illustrate some algebraic and combinatorial tricks that can be performed in order to embed groups of homeomorphisms into
chain groups. These methods are quite flexible and broadly applicable, and will be useful in our proofs of the existence of finitely generated
groups with exotic regularity properties and which only admit abelian proper quotients.

For this section, it will be convenient to pass from $\Homeo_+[0,1]$ to $\Homeo_+(\bR)$, where of course the latter is isomorphic to the latter
by compactifying the real line with the points $\pm\infty$.

\subsection{Maximal rank and the chain group trick}
If $G$ is a group and $H$ is a subgroup with some interesting properties, we might hope that in homomorphic images of $G$, some relevant
features of $H$ are preserved. Chain groups provide a setting in which to arrange for pairs $(G,H)$ like this. 
For instance, if $G$ is a
chain group that acts minimally (either on $(0,1)$ or on $\bR$) and if $H$ lies in the commutator subgroup $G'$, then every nonabelian
quotient of $G$ will contain an isomorphic copy of $H$ by Theorem~\ref{thm:chain-dichotomy}.

We will retain notation from Subsection~\ref{ss:2-chain} and remind the reader that
\[
a(x)=x+1\quad\text{and}\quad
b(x)=
\begin{cases}
 x&\text{ if }x\leq 0,\\
 2x&\text{ if }0<x< 1,\\
  x+1&\text{ if }1\leq x\\
\end{cases}
\]

\begin{lem}[\cite{KKL2019ASENS}, Lemma 4.1]\label{lem:comm-embed}
Let $\gam\le \Homeo_+(0,1)$ be an $n$--generated group for $n\geq 1$, and suppose that \[\gam/\gam'\cong\bZ^n.\] Then the group $\Gamma$
embeds as a subgroup of \[\form{  \gam,a}'\le \Homeo_+(\bR).\]
\end{lem}
\begin{proof}
Choose generators $\{g_1,\ldots,g_n\}$ for $\Gamma$ such that $\gam/\gam'$ is freely generated as an abelian group
by the images of these generators.
Write \[h_i=a^ig_ia^{-1}\in\Homeo_+(i,i+1).\] The homeomorphism $h_i$ is just $g_i$, translated over by $i$. We now set $k_i=g_ih_i^{-1}$,
and we let \[K=\form{  k_1,\ldots,k_n}\le \form{  \gam,a}'.\] We have that $K\cong \gam$.
Indeed, we view $\bZ^n$ as the free abelian group
on $\{g_1,\ldots,g_n\}$, and we have a map
\[\phi\colon \gam\longrightarrow \gam\times\bZ^n,\] given by $g\mapsto (g,g^{-1})$.
Writing $p_i$ for the projections onto the two factors for $i\in\{1,2\}$,
we have $p_i\circ\phi$ is surjective and $p_1\circ\phi$ is an isomorphism. The homeomorphism $h_i$
acts on the interval $(i,i+1)$ by a copy of $\bZ$, and this copy of $\bZ$ is a well-defined quotient of $\Gamma$, given by sending
$g_i$ to $h_i^{-1}$ and sending $g_j$ to the identity for $j\neq i$. It follows that the map \[\psi\colon \gam\longrightarrow K\] is simply identifying
$\phi(g_i)$ with $k_i$, and this map is now clearly an isomorphism.
\end{proof}

Adjoining $a$ to a rather flexible class of groups of homeomorphisms of $\bR$ yields a profusion of chain groups. To this end, write $\DD$
for the set of homeomorphisms $g$ of $\bR$ that have the following property.
\[
g(x)=
\begin{cases}
 x&\text{ if }x\leq 0,\\
 g(x)\in (x,x+1)&\text{ if }0<x< 1,\\
  x+1&\text{ if }1\leq x\\
\end{cases}
\]

The homeomorphism $b$ from Subsection~\ref{ss:2-chain} is an element of $\DD$.

\begin{lem}[\cite{KKL2019ASENS}, Lemma 4.2]\label{lem:dd-chain}
Let $\{g_1,\ldots,g_n\}\sse \DD$. Then the group \[\form{  g_1,\ldots,g_n,a}\le \Homeo_+(\bR)\] 
is isomorphic to an $(n+1)$--chain group.
\end{lem}

Before giving a proof of Lemma~\ref{lem:dd-chain}, note that we did not assume that $\{g_1,\ldots,g_n\}$ are distinct homeomorphisms.
In particular, there is no harm in assuming, say, \[g_1=\cdots=g_{n-1}=b.\] Lemma~\ref{lem:dd-chain} immediately implies the following
curious fact, which illustrates the protean nature of chain groups.

\begin{cor}
For all $n\geq 2$, Thompson's group $F$ is isomorphic to an $n$--chain group.
\end{cor}

In fact, much more is true:
\begin{thm}[\cite{KKL2019ASENS}, Theorem 4.7]\label{thm:chain-protean}
Let $n\geq m\geq 2$, and let $G$ be an $m$--chain group. Then $G$ is isomorphic to an $n$--chain group. 
\end{thm}

We will not prove Theorem~\ref{thm:chain-protean} as it will not be relevant in the sequel, but an ambitious reader should be able to prove it
for themself without significant difficulty. Another bizarre phenomenon in the theory of chain groups which follows from
Theorem~\ref{thm:chain-protean} concerns the abelianizations of chain groups. Recall that for Thompson's group $F$, we have
$F/F'\cong\bZ^2$, where this isomorphism comes from computing the germs at the endpoints of the interval (or at $\pm\infty$ for an
action on $\bR$). For a general chain group $G$, it is much more difficult to compute the rank of the abelianization of $G$, and even harder
to
compute the full abelianization of $G$. Of course, we have that $G/G'$ surjects to $\bZ^2$, since we can compute the germs of the action
at $\{\inf\supp G,\sup\supp G\}$. However, the abelianization might be somewhat larger. Recall the
Higman--Thompson groups $\{F_n\}_{n\geq 2}$
from Subsection~\ref{ss:2-chain}, which are realizable as $n$--chain groups. From the presentation we gave, one can check easily that
$F_n/F_n'\cong\bZ^n$. Moreover, Theorem~\ref{thm:chain-protean} shows that $F_n$ is isomorphic to an $m$--chain group for all $m\geq n$.
This implies that the abelianization of an $m$--chain group can have rank anywhere from $2$ to $m$. We do not know if the
abelianization of a chain group can have torsion, but it seems not unlikely that they can.

\begin{proof}[Proof of Lemma~\ref{lem:dd-chain}]
Let \[S=\{g_1,\ldots,g_n\},\] and write \[f_0=g_1^{-1}a,\quad f_n=a^{n-1}g_na^{-n+1}.\] By construction, $f_0$ agrees with $a$ near $-\infty$
and with the identity to the right of $1$, and $f_n$ agrees with $a$ at $\infty$ and with the identity to the left of $n-1$.
Write \[f_i=(a^ig_{i+1}^{-1}a^{-i})(a^{i-1}g_ia^{-i+1}),\quad\textrm{for}\,\, 1\leq i\leq n-1.\] Notice that \[\form{ \{f_0,\ldots,f_n\}}=\form{  S,a
}.\]
Observe that for $1\leq i\leq n-1$, we have \[\supp f_i=(i-1,i+1).\] It follows that setting $J_i=\supp f_i$ for $0\leq i\leq n$, we have
that $\{J_0,\ldots,J_n\}$ forms an $(n+1)$--chain. To verify that $\{f_1,\ldots,f_n\}$ generate a chain group, we note that
\[f_{i+1}f_i(i)=i+1.\] Since \[i=\inf J_{i+1}\quad \textrm{and}\quad i+1=\sup J_i,\] this completes the proof.
\end{proof}

The proof of the following lemma is an easy computation, and we leave it as an exercise for the reader.
Recall from Notation~\ref{notation:diff-J} that for an interval $J\sse \bR$ we write
\[\Homeo_J(\bR):=\{ f\in \Homeo_+(\bR)\mid \supp f\sse J\}.\]
\begin{lem}[\cite{KKL2019ASENS}, Lemma 4.4]\label{lem:expansion}
Let $J:=(1/4,1/2)$.
If $f\in\DD$ satisfies $f(J)=(1/2,1)$, then $fg\in\DD$ for all $g\in\Homeo_J(\bR)$.
\end{lem}

We are finally able to give the proof of Theorem~\ref{thm:chain-subgp} we claimed above. We will actually prove a stronger statement:

\begin{thm}\label{thm:chain-subgp-full}
Let $\Gamma$ be a subgroup of $\Homeo_+(\bR)$ generated by $\{f_1,\ldots,f_n\}$, for $n\geq 1$.
\begin{enumerate}[(1)]
\item
The group $\Gamma$ embeds into an $(n+2)$--chain group $G$ such that $G'$ is simple.
\item
If $\supp f_1$ has finitely many components, then we may arrange in addition for $G$ to be an $(n+1)$--chain group.
\end{enumerate}
In either case, if we have that $\gam/\gam'\cong\bZ^n$ then we may arrange for $\Gamma$ to embed in the derived subgroup $G'$.
\end{thm}
\begin{proof}
We still let $J:=(1/4,1/2)$.
Write $S=\{f_1,\ldots,f_n\}$. Without loss of generality, we may assume that $G\le \Homeo_J(\bR)$, as this latter group is isomorphic
to $\Homeo_+(\bR)$.

Note that Lemma~\ref{lem:expansion} shows that $b\cdot S\sse \DD$, whence that group
\[G=\form{  S,b,a}=\form{  b\cdot S,b,a}\] is abstractly isomorphic to a chain group by Lemma~\ref{lem:dd-chain}.
Since the group $\form{  a,b}\cong F$ is topologically conjugate to the usual realization of $F$ as a subgroup of $\PL_+(I)$ and since
the latter group acts minimally on $(0,1)$, we have that $\form{  a,b}$ acts minimally on $\bR$. Therefore, $G$ acts minimally on
$\bR$, and Theorem~\ref{thm:chain-dichotomy} implies that $G'$ is simple. This proves the first assertion of the theorem.

If $f_1$ has finitely many components, then we may conjugate $S$ within the group $\Homeo_J(\bR)$ in order to arrange for
$f_1$ to be piecewise linear with dyadic breakpoints. In particular, in this way we may arrange for $f_1$ to be an element of
$F=\form{  a,b}$. Setting $S_0=S\setminus\{f_1\}$, we set \[G=\form{  S_0,b,a},\] which again is
isomorphic to a $(n+1)$--chain group by Lemma~\ref{lem:dd-chain}, and the action on $\bR$ remains minimal. This proves the
second assertion of the theorem.

For the third assertion, if $\gam/\gam'\cong\bZ^n$ then we may realize \[\gam\le \form{ \gam,a}'\le \form{ \gam,a,b}',\]
by Lemma~\ref{lem:comm-embed}. The conclusion follows immediately.
\end{proof}

We state the following fact, which in~\cite{KK2020crit} is called the chain group trick, and which will be useful in the sequel. We leave
the details of the proof to the reader.

\begin{cor}[The chain group trick]\label{cor:chain-gp-trick}
Let $\gam\le \Homeo_+(\bR)$ be an $m$--generated group such that $\supp\gam$ is compactly contained in $(0,1)$, and let $F_{GS}$
be a smooth realization of Thompson's group $F$ on $[0,1]$ that acts minimally on $(0,1)$, as furnished by 
Corollary~\ref{cor:ghys-serg}. Write
\[G=\form{ \gam,F_{GS}}.\] Then the following conclusions hold.
\begin{enumerate}[(1)]
\item
The group $G$ is an $(m+2)$--chain group acting minimally on $(0,1)$, so that $G'$ is simple and every proper quotient of $G$ is abelian.
\item
If $\gam/\gam'$ is free abelian, then there is an embedding of $\Gamma$ into $G'$.
\end{enumerate}
\end{cor}

We remark that in the second part of Corollary~\ref{cor:chain-gp-trick}, we do not assume that $\gam/\gam'$ has rank $m$. The reader will
observe that since $F_{GS}\le \Diff_+^{\infty}([0,1])$, we have that the regularity of $G$ coincides with the regularity of $\Gamma$. 
Since $\Gamma$
is compactly supported and since 
each element in $F_{GS}$ is tangent to the identity at $\partial I$,
 we have that the action of $G$ extends by the identity to $\bR$, without any loss of regularity. By identifying $0$ and $1$, we also obtain an action
of $G$ on $S^1$ without any loss of regularity.

\subsection{The rank trick}\label{ss:rank-trick}

In results such as Theorem~\ref{thm:chain-subgp-full}, and more fundamentally in Lemma~\ref{lem:comm-embed}, the maximal rank
hypothesis on the abelianization can be somewhat annoying. In this subsection, we illustrate a technical tool which will be of use later,
and which helps satisfy maximal rank hypotheses.

\begin{lem}[The rank trick; cf.~\cite{KK2020crit}, Lemma 6.1]\label{lem:rank-trick}
Let $G$ be a group such that $G/G'$ is a finitely generated torsion-free abelian group. 
If we have a representation
\[\rho\colon G\longrightarrow \Homeo_c(\bR),\] then there exists another representation
\[\tau\colon G\longrightarrow\form{ \rho(G),\Diff_c^{\infty}(\bR)}\] that has the following properties:
\begin{enumerate}[(1)]
\item
We have $\tau$ agrees with $\rho$ on $G'$;
\item
We have $\tau(G)/(\tau(G))'\cong G/G'$.
\end{enumerate}
\end{lem}
\begin{proof}
Let \[G/G'\cong\bZ^m,\] where $m\geq 1$; note that the case where $G=G'$ is trivial by setting $\rho=\tau$. We choose arbitrary
compactly supported
$C^{\infty}$ diffeomorphisms $\{f_1,\ldots,f_m\}$ such that \[\supp f_i\cap \supp f_j=\varnothing\] for $i\neq j$, and such that
\[\supp f_i\cap\supp\rho(G)=\varnothing\] for all $i$. Observe that there is a natural surjective map \[\phi\colon G\longrightarrow
\form{  f_1,\ldots,f_m}\cong \bZ^m.\] We simply set \[\tau\colon G\longrightarrow \form{  \rho(G),\Diff_+^{\infty}(\bR)},\quad
\tau(g)=\rho(g)\phi(g).\]

It is immediate that $\supp\tau(G)$ has compact closure in $\bR$, and that $\tau$ agrees with $\rho$ on $G'$.
The abelianization of $\tau(G)$ agrees with the abelianization of $G$, since we have a surjective composition of maps
\[G\longrightarrow \tau(G)=\rho(G)\times\form{  f_1,\ldots,f_m}\longrightarrow \form{  f_1,\ldots,f_m}\cong\bZ^m,\]
where the last arrow is just the projection onto the second coordinate. The reader may compare with the proof of
Lemma~\ref{lem:comm-embed} above.
\end{proof}

The reader will note that the hypothesis that ${\suppc \rho(G)}$ is compact can be relaxed; one need only assume that
$\supp \rho(G)$ be bounded from one side, and in the conclusion of the lemma, $\supp\tau(G)$ will also be bounded from the same side.


%
%
%

\chapter{The Slow Progress Lemma}\label{ch:slp}
\begin{abstract}In this chapter, we give a proof of the Slow Progress Lemma. This result, stated as Theorem~\ref{t:slp},
asserts that a certain iteration of smoother group elements makes slower progress of the orbit in the sense of covering distances than
corresponding iterations of less smooth group elements.
This lemma is a key dynamical ingredient of our main result (Theorem~\ref{t:optimal-all}), and also serves as a bridge between
analytic information (i.e. regularity) and dynamical data (i.e. the asymptotics of covering distances) associated to a group action.
The version of the Slow Progress Lemma we present here is significantly stronger and quantitatively more precise than the one 
in~\cite{KK2020crit}.\end{abstract}

\section{Statement of the result}
Recall that a \emph{concave modulus of continuity}\index{concave modulus of continuity} is simply a concave homeomorphism of $[0,\infty)$.
Often, a concave modulus is specified only near zero, though
such a local definition can always be replaced by a smooth concave modulus defined on the whole interval $[0,\infty)$;
see Lemma~\ref{lem:medvedev} for more detail.

We say a concave modulus $\beta$ is \emph{sup-tame}\index{sup-tame modulus} if
\[\lim_{t\to 0+} \sup_{x>0} t\beta(x)/\beta(tx)=0.\]
On the other hand, $\beta$ is said to be \emph{sub-tame}\index{sub-tame modulus} if 
\[\lim_{t\to 0+} \sup_{x>0} \beta(tx)/\beta(x)=0.\]
Intuitively, a sup-tame modulus is the one that is ``not too small'' near zero;
on the other hand, one may say that a sub-tame modulus is ``sufficiently small''.
For instance, if $\tau\in(0,1)$ then $\beta(x)=x^\tau$ is both sup-and sub-tame.
The Lipschitz modulus $\beta(x)=x$ is the smallest possible modulus; 
a Lipschitz continuous function is $\alpha$--continuous for all concave moduli $\alpha$. 
One may easily check that the Lipschitz modulus is sub-tame but not sup-tame.
A ``big'' concave modulus such as the one named $\beta$ below is sup-tame but not sub-tame:
\[
\beta(x) =1/\log(1/x).\]
Note that the above definition of $\beta$ makes sense only for small $x>0$.

We say that $(k,\beta)$ is a \emph{tame pair}\index{tame pair},
and write $\beta\succ_k0$,
if $k\in{\bZ_{>0}}\cup\{0\}$
and if $\beta$ is a concave modulus such that one of the following holds:
\begin{itemize}
\item $k\ge2$;
\item $k=1$ and $\beta$ is sub-tame (``sufficiently small'');
\item $k=0$ and $\beta$ is sup-tame (``not too small'').
\end{itemize}

For $k\ne1$, Mather~\cite{Mather1,Mather2} proved  the simplicity of the group $\Diff_c^{n+k}(M^n)_0$ of  compactly supported $C^{n+k}$ 
diffeomorphisms of an $n$--manifold $M$ that are isotopic to the identity by compactly supported isotopies. One can actually
extend Mather's arguments to prove that whenever $\beta\succ_k0$, the following group is simple (cf.~\cite{CKK2019}):
\[\Diff_c^{n+k,\beta}(M^n)_0:=\left\{ f\in \Diff_c^{n+k}(M^n)_0\mid f^{(n+k)}\text{ is locally }\beta\text{--continuous}\right\}.\]

Though Mather did not actually used this term, it is interesting to note that the precisely same notion of tame pairs
 appear again in quite a different setting as in Theorem~\ref{t:slp} below.

As usual, we let $I$ denote a compact, nondegenerate interval.
For a set $V\sse\Homeo_+(I)$, we will write $\CD_V(x,y)$ for the corresponding covering distance on $I$. That is,
\[
\CD_V(x,y):=\min\left\{m \mid [x,y]\sse J_1\cup \cdots\cup J_m\text{ for some }J_i\in\bigcup_{ v\in V} \pi_0\supp v\right\}.\]

Let $\{a_p\}_p$ and $\{b_p\}_p$ be real sequences indexed by a subset $P$ of $\bN$.
If there exists some $C>0$ such that $a_p\le Cb_p$ for all $p\in P$,
then we write 
\[a_p\preccurlyeq b_p.\]
We also add a phrase \emph{for almost all $i\in P$} if 
\[
\lim_{N\to\infty}\frac{\#\{ p\in P \mid a_p\le C b_p\}\cap[1,N]}{\#(P\cap[1,N])}=1.\]
We will often consider the case $P=\bN$.

\begin{thm}[Slow Progress Lemma]\label{t:slp}
Let $k\in{\bZ_{>0}}$ and let $\beta$ be a concave modulus such that $\beta\succ_k0$.
Suppose we have a sequence $\{N_i\}_{i\ge1}\sse  {\bZ_{>0}}$ satisfying  \[N_i\preccurlyeq i^{k-1}/\beta(1/i)\] for almost all $i\ge1$.
Suppose furthermore that $V$ is a finite set such that either \[V\sse \Diff_+^{k,\beta}(I)\quad\textrm{or}\quad 
V\sse \Diff_+^{k,\mathrm{bv}}(I),\]
and let $\{v_i\}_{i\ge1}$ be a sequence of elements of $V$.
Then there exists a positive constant $\epsilon$ depending only on $k$ such that
 for every point $x_0$ in $I$, we have 
\[
\limsup_{i\to\infty} \frac1i \CD_V\left(x_0,v_i^{N_i} \cdots v_2^{N_2} v_1^{N_1}(x_0)\right)\le1-\epsilon.\]
\end{thm}

We will  prove that one can set $\epsilon:=1/(4k+2)$.
As a special case of the Slow Progess Lemma, one obtains the same conclusion if $N_i\preccurlyeq i^k$ and if 
$V\sse\Diff_+^{k+1}(I)$ or $V\sse\Diff_+^{k+\mathrm{bv}}(I)$, since the hypothesis of the theorem holds for $\beta(x)=x$.

\begin{rem}
In~\cite{KK2020crit}, a more restrictive, qualitative version of the Slow Progress Lemma was established.
Namely, it was previously proved that
\[
\lim_{i\to\infty} \left(i -  \CD_V\left(x_0,v_i^{N_i} \cdots v_2^{N_2} v_1^{N_1}(x_0)\right)\right)=\infty,\]
under an additional hypothesis that the natural density of 
\[
\{i\in\bZ_{>0}\mid v_i=v\}\]
is well-defined for each $v\in V$. In Theorem~\ref{t:slp}, this additional hypothesis is dropped and the conclusion is strengthened.
\end{rem}

Observe that the hypotheses (and hence also the conclusion) of the Slow Progress Lemma are invariant under topological conjugacy.
Thus, one obtains the following topological non-smoothability criterion for group actions.

\begin{cor}\label{c:slp-0}
Let $k\in{\bZ_{>0}}$ and $\beta\succ_k0$.
Suppose $\{N_i\}_{i\ge1}\sse \bN$ is a sequence such that 
\[N_i\preccurlyeq i^{k-1}/\beta(1/i)\]
for almost all $i\ge1$. Assume we have an action 
\[
\psi\co G\to\Homeo_+(I)\]
of a group $G$ generated by a finite set $V$,
and a point $x_0\in I$ 
such that for some sequence
 $\{v_i\}_{i\ge1}$ of elements of $V$,
we have
\[\limsup_{i\to\infty} \frac1i {\CD_V\left(x_0,\psi\left(v_i^{N_i} \cdots v_2^{N_2} v_1^{N_1}\right)(x_0)\right)}=1.\]
Then $\psi$ is not topologically conjugate  into $\Diff_+^{k,\beta}(I)$ nor into $\Diff_+^{k,\mathrm{bv}}(I)$.
\end{cor}

\section{Natural densities}
For a set $A$, we with write $\#A$ for its cardinality.
\bd\label{d:density}
Let $P\sse \bN$ be a set. 
The \emph{upper density}\index{upper density} of $P$ is defined as
\[\overline{d}_\bN(P):=\limsup_{n\to\infty}  \#(P\cap [1,n])/n.\]
Similarly, the \emph{lower density}\index{lower density} of $P$ is
\[\underline{d}_\bN(P):=\limsup_{n\to\infty}  \#(P\cap [1,n])/n.\]
If $\overline{d}_\bN(P)$ coincides with $\underline{d}_\bN(P)$, this number is called the \emph{natural density}\index{natural density} of $P$.
\ed

The following is well-known~\cite{Salat1964MZ}; a particularly simple proof was given by Moser~\cite{Moser1958}, which we will reproduce
here for the convenience of the reader.
\begin{lem}\label{l:sum-density}
Let $Q\sse\bN$ be a set.  If \[\sum_{q\in Q} 1/q<\infty,\] then $d_\bN(Q)=0$.\end{lem}
\begin{proof}
For $n\in\bN$, we write $Q_n=Q\cap [1,n]$. Write \[\chi(n)=\sum_{q\in Q_n} 1,\quad H(n)=\sum_{q\in Q_n}\frac{1}{q},\] with $H(0)=0$ by
convention.
Observe that \[\chi(n)=\sum_{i=1}^n i\cdot (H(i)-H(i-1)).\]
Thus, we have that \[\frac{\chi(n)}{n}=H(n)-\frac{1}{n}\sum_{i=1}^n H(i-1).\] Note that the second term on the right hand sum is the $n^{th}$
Ces\`aro sum of the sequence $\{1/q\}_{q\in Q}$, and so converges to $\lim_{n\to\infty}H(n)$ if this latter
limit exists. If \[\sum_{q\in Q} 1/q<\infty,\]
we obtain $\lim_{n\to\infty}\chi(n)/n=0$, so that $Q$ has density zero.
\end{proof}

A set with a high upper density contains a long sequence of consecutive numbers. More precisely, we have the following.

\begin{lem}\label{l:upper-density}
For each subset $P\sse\bN$ and for each $C\in\bZ_{>0}$, we have that 
\[
\overline{d}_\bN\left( P\cap (P-1)\cap\cdots\cap (P-C+1)\right)\ge 1-C\left(1-\overline{d}_\bN(P)\right).\]\
\end{lem}

\bp
Set $\epsilon:=1-\overline{d}_\bN(P)$. 
For all $\epsilon'>\epsilon$, there  exist infinitely many $N$ such that 
\[\#(P\cap[1,N])>(1-\epsilon')N.\]
For such an $N$ we have
\[\#((P-1)\cap [1,N])=\#(P\cap [2,N+1])>   (1-\epsilon')N-1.\]
Using the equation $\#(A\cap B)=\#A+ \#B-\#(A\cup B)$, we obtain
\[\#(P\cap (P-1)\cap [1,N])> (1-2\epsilon')N-1.\]
We inductively deduce that 
\[\#(P\cap (P-1)\cap \cdots\cap (P-C+1)\cap [1,N])  > (1-C\epsilon')N-(1+2+\cdots+(C-1)).\]
As there are infinitely many such $N$, we have that
\[\overline{d}_\bN\left( P\cap (P-1)\cap\cdots\cap (P-C+1)\right)\ge 1-C\epsilon'.\]
Since $\epsilon'>\epsilon$ was arbitrary, the proof is complete.
\ep

If $\{d_i\}_i$ is a discrete walk in $\bR$ that can either stay put, move forward, or move backward at each step,
then we show below that the limit superior of the average speed $d_n/n$ can be realized by the upper density of the times at which
new ``progress'' is made. 

\begin{lem}\label{l:p-almost}
If $\{d_i\}_i$ is a sequence of nonnegative integers  such that $|d_{i+1}-d_i|\le 1$ for each $i\ge0$,
then we have that
\[
\overline{d}_\bN\left\{n\in\bN \middle\vert d_n > \max_{0\le m<n} d_m\right\}
= \limsup_{i\to\infty} {d_i}/{i}.\]
\end{lem}
\bp
We may normalize the sequence and assume $d_0=0$. We are only interested in the case when
\[P:=\left\{n\in\bN \mid d_n > \max_{0\le m<n} d_m\right\}\]
is an infinite set, for otherwise the claim is trivial. 

Let us enumerate $P$ in increasing order, so that
\[P=\{p_1<p_2<\cdots\}.\] For all sufficiently large $n\gg0$, we can find $j\in\bN$ such that $p_j\le n<p_{j+1}$. Then
\[
\#(P\cap[1,n]) = j = d_{p_j}\ge d_n.\] 
Dividing by $n$ and taking limit superiors over $n$, we obtain $\overline{d}_\bN(P)\ge \limsup_n d_n/n$.

To prove the opposite inequality, we may assume $\overline{d}_\bN(P)>0$
and pick an arbitrary $\epsilon'\in(0,\overline{d}_\bN(P))$. For each $p_j\le m<p_{j+1}$,
note that
\[ \#(P\cap[1,m])/m = j/m\le j/p_j.\]
By the definition of an upper density, we can find infinitely many $j$ such that 
\[
\epsilon'<j/p_j = d_{p_j}/p_j.\]
It follows that \[ \epsilon'\le \limsup_i \frac{d_i}{i},\] and the proof is complete.
\ep
\begin{rem}
If the sequence $\{d_i\}_i$ above is chosen as a random walk on $\bR$, then by the Law of Iterated Logarithm~\cite{Kolmogorov29},
one has that
$\limsup_i d_i/i=0$ \emph{almost surely}.
\end{rem}

\section{Probabilistic dynamical behavior}
The main idea of the Slow Progress Lemma can be roughly summarized by the phrase that ``smoother diffeomorphisms are slower''. Here, 
smoother diffeomorphisms are those of higher regularity. The ``slowness'' of such diffeomorphisms refers to the fact that the
displacement $|f(x)-x|$ is small in a certain
probabilistic sense described in Definition~\ref{d:expansive} and Lemma~\ref{l:prob-dyn} below.
\bd\label{d:k-fixed}
Let $k\in{\bZ_{>0}}$, and let $f\in\Homeo_+(I)$.
We say $f$ is \emph{$k$--fixed}\index{$k$--fixed homeomorphism} on an interval $J\sse I$
if at least one of the following conditions holds: 
\begin{itemize}
\item $f$ fixes more than $k$ points in $ J$.
\item $ J$ contains an accumulation point of $\Fix f$.
\end{itemize}
\ed

We will denote the identity function as $\Id(x)=x$. Note
 that \[\Id^{(0)}(x)=x,\quad \Id^{(1)}(x)=1,\quad \Id^{(i)}(x)=0\quad \textrm{for $i>1$}.\]
\begin{lem}\label{l:k-fixed}
Let $k\in{\bZ_{>0}}$. 
If $f\in\Diff_+^k(I)$ is $k$--fixed on a compact interval $J\sse I$,
then for each $i\in\{0,1,\ldots,k\}$ there exists a point $s_i\in J$ such that 
\[
f^{(i)}(s_i)=\Id^{(i)}(s_i).\]
\end{lem}

\bp
Let us first assume that $J$ contains more than $k$ fixed points of $f$. By the Mean Value Theorem (or Rolle's Theorem), one inductively
sees for each $i\in\{0,1,2,\ldots,k\}$  that
\[
\#\left\{ s\in J\mid f^{(i)}(s)=\Id^{(i)}(s)\right\}>k-i.\]
Thus, the proof is complete in this case.

Assume instead that $J$ contains an accumulation point $x_0$ of $\Fix f$.
We may further assume that $x_0\in \partial J$, for otherwise the proof is trivial by the previous paragraph.  
As above, for each $i\in\{0,1,\ldots,k\}$, we see that
$x_0$ is an accumulation point of 
\[\{ s\in J\mid f^{(i)}(s)=\Id^{(i)}(s)\}.\]
So, we may set $s_i=x_0$ for all $i\le k$.
\ep

For two subsets $A$ and $B$ of $\bR$, we let \[d(A,B)=\inf_{a\in A,b\in B} d(a,b)\] denote the  distance between $A$ and $B$.

\bd\label{d:expansive}
Let $f\co J\longrightarrow J$ be a homeomorphism of an interval $J$.
We say that $f$ is \emph{$\delta$--expansive on $J$}\index{$\delta$--expansive homeomorphism} for some $\delta>0$ if 
\[\sup_{x\in \Int J} \frac{|f(x)-x|}{d\left(\{x,f(x)\},\partial J\right)}\ge \delta.\]
We say that $f\co J\longrightarrow J$ is \emph{at most $\delta$--expansive on $J$} if
\[\sup_{x\in \Int J} \frac{|f(x)-x|}{d\left(\{x,f(x)\},\partial J\right)}\le \delta.\]
\ed

For instance, if $f$ acts on $J=[a,b]$ such that some $x\in (a,b)$ satisfies
\[\frac{f(x)-x}{x-a}\ge\delta,\]
then $f$ is $\delta$--expansive on $J$.
Note that the inverse of a $\delta$--expansive homeomorphism is also $\delta$--expansive.

For a concave modulus $\beta$ and  a $C^\beta$--continuous map $f$,
we  define its $\beta$--norm as
\[
[f]_\beta:=\sup_{x\ne y}\frac{|f(x)-f(y)|}{\beta\left(|x-y|\right)}.\]

\begin{lem}\label{l:expansive}
Let $k\in{\bZ_{>0}}$, and let $f\in \Diff_+^k(I)$.
If $f$ is $k$--fixed 
on some compact interval $J\sse I$, 
then $f$ is at most $\delta_1$--expansive on $J$ for
\[\delta_1:=  |J|^{k-1} \left(\sup_{x\in J} \left|f^{(k)}(x) - \Id^{(k)}(x)\right|+\sup_{x\in J} \left|\left(f^{-1}\right)^{(k)}(x) - \Id^{(k)}(x)\right|\right).\]
If we further assume that $f\in \Diff_+^{k,\beta}(I)$ for some concave modulus $\beta$,
then $f$ is at most $\delta_2$--expansive on $J$ for
\[\delta_2:= |J|^{k-1}\beta(|J|)\left(\left[f^{(k)}\right]_\beta+\left[\left(f^{-1}\right)^{(k)}\right]_\beta\right).\]
\end{lem}
\bp
By Lemma~\ref{l:k-fixed}, we have some $s_i\in J$ for each $i\in\{0,1,\ldots,k\}$ such that $f^{(i)}(s_i)=\Id^{(i)}(s_i)$.
Let $x\in \Int  J=(a,b)$. We compute that
\begin{align*}
f(x) - x &=\int_{s_0}^x (f'(t_1)-1)dt_1= \int_{s_0}^x \int_{s_1}^{t_1} f''(t_2)dt_2\;dt_1=\cdots\\
&=\int_{s_0}^{x}\int_{s_1}^{t_1}\cdots\int_{s_{k-1}}^{t_{k-1}} \left(f^{(k)}(t_k)-\Id^{(k)}(t_k)\right)\; dt_k\cdots dt_1.
\end{align*}
It follows that 
\[ |f(x)-x| \le \sup_{t\in J}\left|f^{(k)}(t)-\Id^{(k)}(t)\right| \cdot   |x-s_0|\cdot |J|^{k-1}  \le \delta_1  |x-s_0|.\]
After choosing $s_0$ to be either $a$ or $b$, and applying the same argument to $f^{-1}$ instead of 
$f$, we conclude that $f$ is at most $\delta_1$--expansive.

A similar argument shows that
\begin{align*}
|f(x) - x| &\le  \int_{s_0}^{x}\int_{s_1}^{t_1}\cdots\int_{s_{k-1}}^{t_{k-1}} \left|f^{(k)}(t_k)-f^{(k)}(s_k)\right|\; dt_k\cdots dt_1\\
&\le \left[f^{(k)}\right]_\beta\cdot\beta(|J|) \cdot  |x-s_0|\cdot |J|^{k-1}\le \delta_2 |x-s_0|,
\end{align*}
which implies the second required inequality.
\ep

Let us also note a simple consequence of the above proof for a later use (Section~\ref{sss:tsuboi}).

\begin{lem}\label{l:expansive2}
Let $k\in{\bZ_{>0}}$, and let $f\in \Diff_+^k(I)$.
If $f$ is $k$--fixed on some compact interval $J\sse I$, 
then we have
\[
\sup_{x\in J} |f'(x)-1|
\le
 \sup_{t\in J}\left|f^{(k)}(t)-\Id^{(k)}(t)\right|
 \cdot |J|^{k-1}.\]
\end{lem}
\bp
Continuing to use the notation as above, we have
\[ |f'(x)-1| \le \sup_{t\in J}\left|f^{(k)}(t)-\Id^{(k)}(t)\right| \cdot   |x-s_1|\cdot |J|^{k-2},\]
which implies the conclusion.
\ep

\begin{lem}\label{l:N-delta-exp}
Suppose that $f\in\Homeo_+(I)$ preserves a compact interval $J\sse I$.
If $f^N$ is $\delta$--expansive on $J$ for some $N\in\bN$ and $\delta>0$, then $f$ is $\frac1N\log(1+\delta)$--expansive on $J$.
\end{lem}
\bp
Let us write $J=[a,b]$. We can find some $x_0\in(a,b)$ such that either
\[
\delta\le \frac{f^N(x_0)-x_0}{x_0-a}
\] 
or
\[
\delta\le \frac{f^N(x_0)-x_0}{b-f^N(x_0)},
\] 
possibly after switching $f$ with $f^{-1}$. We may further assume to have the former of the above two cases,
as the latter case can be treated similarly.
Then we have
\[
1+\delta \le \frac{f^N(x_0)-a}{x_0-a}=\prod_{i=0}^{N-1} \frac{f^{i+1}(x_0)-a}{f^{i}(x_0)-a}.\]
Hence, for some $x_1:=f^{i}(x_0)$ we have that
\[
(1+\delta)^{1/N}\le \frac{f(x_1)-a}{x_1-a}.\]
We have an estiamte
\[
\frac{f(x_1)-x_1}{x_1-a}= \frac{f(x_1)-a}{x_1-a}-1 \ge(1+\delta)^{1/N}-1 \ge  \frac{\log(1+\delta)}{N},\]
which implies the conclusion.
\ep

\begin{lem}\label{l:decomposition}
Every collection of intervals in the real line with intersection multiplicity $K>0$ admits a partition into $K$ subcollections,
each of which consists of disjoint intervals.
\end{lem}
Here, we recall that the \emph{intersection multiplicity}\index{intersection multiplicity} of a collection $X$ of intervals
refers to the maximum number of intervals in $X$ that can contain a given point.
We thank Andreas Holmsen for teaching us the following argument.
\bp[Proof of Lemma~\ref{l:decomposition}]
Let $X$ be the given collection of intervals, and let $\Gamma$ be its \emph{intersection graph}\index{intersection graph};
that is, the vertices of $\Gamma$ 
correspond to intervals in $X$,
and two vertices are adjacent if the two corresponding intervals have nonempty intersection. 

Helly's Theorem (for $\bR$) says that
if $Y$ is a finite collection of intervals in $\bR$ such that every pair of elements of $Y$ intersect, then \[\bigcap_{J\in Y} J\neq\varnothing.\] We leave
the justification of Helly's Theorem as an exercise for the reader.
By Helly's theorem for $\bR$, we see that $K$ is the maximal size of a clique in $\Gamma$.

It now suffices for us to show that the chromatic number $K'$ of $\Gamma$ coincides with $K$.
It is obvious that $K'\ge K$.
It is also well-known that the intersection graph of finitely many intervals is \emph{chordal}\index{chordal graph} (i.e.~every cycle of
length $4$ or more admits a chord), and hence
\emph{perfect}~\cite{golumbic2004}\index{perfect graph}; this means that $K=K'$ when $X$ is finite. 
So, we may assume $X$ is infinite. By the De Brujin--Erd\"os Theorem (\cite{dBE1951}, cf.~Proposition~\ref{prop:db-e} below),
we can find a finite subcollection
$X'$ of $X$ such that the intersection graph $\Gamma'$ of $X'$ has the same chromatic number $K'$.
Also, the maximal size of a clique in $\Gamma'$ is at most $K$. By the result for the case when $X$ is finite,
we conclude that $K'\le K$.
\ep

For the convenience of the reader, we give a quick proof of the De Brujin--Erd\"os Theorem,
originally given by Gottschalk~\cite{Gottschalk51}.

\begin{prop}\label{prop:db-e}
Let $\Gamma$ be a graph, all of whose finite subgraphs are $k$--colorable. Then $\Gamma$ is itself $k$--colorable.
\end{prop}
\begin{proof}
Let $X=k^{V(\gam)}$ be the space of all assignments of $k$ colors to the vertices of $\Gamma$, where $k=\{0,1,\ldots,k-1\}$
is given the discrete topology.
Then $X$ is compact  by Tychonoff's Theorem. Valid colorings of finite subgraphs of $\Gamma$ are easily seen to be closed
subsets $X$.

If $F_1$ and $F_2$ are finite subgraphs of $\Gamma$ and if $X_{F_1}$ and $X_{F_2}$ are closed subsets of $X$ corresponding
to valid $k$--colorings of $F_1$ and $F_2$ respectively, then $X_{F_1}\cap X_{F_2}$ corresponds to valid $k$--colorings of $F_1\cup F_2$.
Moreover, this intersection is nonempty
by hypothesis. The finite intersection property characterization of compactness shows that if every finite subgraph of $\Gamma$ admits a valid
coloring, then the intersection of the corresponding closed subsets of $X$ is nonempty and consists precisely of valid $k$--colorings of $\Gamma$.
\end{proof}

\begin{lem}\label{l:mod-compare}
Let $k\in{\bZ_{>0}}$ and let $\beta$ be a concave modulus that satisfies $\beta\succ_k0$.
Suppose $P$ is a subset of $\bN$.
If two positive sequences $\{a_p\}_{p\in P}$ and $\{b_p\}_{p\in P}$ satisfy 
\[a_p^{k-1}\beta(a_p)\preccurlyeq b_p^{k-1}\beta(b_p),\]
then we have that $a_p\preccurlyeq b_p  $.
\end{lem}
\bp
Suppose not. We can find an increasing sequence $\{p_i\}$ in $P$ such that the sequence
$t_{p_i}:=b_{p_i}/a_{p_i}$ converges to zero.
Since $(k,\beta)$ is a tame pair, we obtain that 
\[
\frac{b_{p_i}^{k-1}\beta\left(b_{p_i}\right)}
{a_{p_i}^{k-1}\beta\left(a_{p_i}\right)}
=
t_{p_i}^{k-1}\cdot\frac{\beta\left(t_{p_i}a_{p_i}\right)}
{\beta\left(a_{p_i}\right)}\to 0.
\]
This contradicts the given hypothesis.
\ep

We can now prove the main result of this section.

\begin{lem}\label{l:prob-dyn}
Let $k\in{\bZ_{>0}}$ and let $\beta$ be a concave modulus that satisfies $\beta\succ_k0$.
Let $P\sse \bN$.
Suppose $V\sse\Homeo_+(I)$ is a finite set.
Assume that for each $p\in P$ we are given a natural number $N_p\in \bN$, an element $v_p\in V$, and a compact interval $J_p$ 
such that the following hold:
\be[(i)]
\item 
$N_p\preccurlyeq p^{k-1}/\beta(1/p)$ for $p\in P$;
\item $v_p$ is $k$--fixed on $J_p$;
\item\label{p:fim} The collection $\{J_p\mid p\in P\}$ has finite intersection multiplicity.
\ee
If  either   $V\sse \Diff_+^{k,\beta}(I)$ or  $V\sse\Diff_+^{k,\mathrm{bv}}(I)$,
then for all $\delta>0$ we have that
\[
d_\bN\left\{ p\in P\mid v_p^{N_p}\text{ is  }\delta\text{--expansive on }J_p\right\}=0.\]
\end{lem}

\bp
By the condition (\ref{p:fim}) above and by Lemma~\ref{l:decomposition}, the collection $\{J_p\mid p\in P\}$ can be partitioned into $K$ 
subcollections of disjoint intervals for some $K>0$. In particular, \[\sum_p |J_p|\le K|I|.\]

Let $Q\sse P$ be the set the density of which is claimed to be zero. 
For each $q\in Q$, Lemma~\ref{l:N-delta-exp} implies that $v_q$ is $\frac{\log(1+\delta)}{N_q}$--expansive on $J_q$.

{\bf Case 1: $V\sse\Diff_+^{k,\beta}(I)$.}
We obviously have that
\[
 \max_{v\in V} \left(\left[ v^{(k)}\right]_\beta +\left[ \left(v^{-1}\right)^{(k)}\right]_\beta 
\right)<\infty.\]
By Lemma~\ref{l:expansive}, we have for $q\in Q$  that
\[
\left(\frac1q\right)^{k-1}\beta\left(\frac1q\right)
\preccurlyeq \frac{\log(1+\delta)}{N_q}
\preccurlyeq|J_q|^{k-1}\beta(|J_q|).\]
From Lemma~\ref{l:mod-compare} we see that 
\[\frac{1}{q} \preccurlyeq |J_q|,\]
and that
\[
\sum_{q\in Q}1/q<C\sum_{q\in Q} |J_q|
\le CK|I|<\infty\]
for some $C>0$. 
 Lemma~\ref{l:sum-density} implies that  $d_\bN(Q)=0$.

{\bf Case 2: $V\sse\Diff_+^{k,\mathrm{bv}}(I)$.}
Note from the concavity that \[\frac{x}{\beta(x)}\le \frac{1}{\beta(1)}\] for  $x\le 1$. It follows for $q\in Q$ that
\[
\frac1{q^k}
\preccurlyeq \frac{\beta(1/q)}{q^{k-1}} \preccurlyeq  \frac1{N_q}.\]
Now, let $q\in Q$.
By Lemmas~\ref{l:k-fixed} and~\ref{l:expansive}, we have that
\begin{align*}
\frac{\log(1+\delta)}{N_q}&
\le |J_q|^{k-1} 
 \left(\sup_{x\in J_q} \left|v_q^{(k)}(x) - \Id^{(k)}(x)\right|+\sup_{x\in J_q} \left|\left(v_q^{-1}\right)^{(k)}(x) - \Id^{(k)}(x)\right|\right)\\
&\le
 |J_q|^{k-1} \left(
\operatorname{Var}\left(v_q^{(k)};J_q\right)+\operatorname{Var}\left(\left(v_q^{-1}\right)^{(k)};J_q\right)\right).
\end{align*}
By  H\"older's inequality, we have some $C,C'>0$ such that
\begin{align*}
\sum_{q\in Q}\frac1q
&\le
\sum_{q\in Q}\frac{C}{N_q^{1/k}}
\le C'
\sum_{q\in Q} |J_q|^{1-1/k} \left(\operatorname{Var}\left(v_q^{(k)};J_q\right)+
\operatorname{Var}\left(\left(v_q^{-1}\right)^{(k)};J_q\right)\right)^{1/k}\\
&\le C'
\left(\sum_{q\in Q} |J_q|\right)^{1-1/k}
\left(\sum_{q\in Q} \operatorname{Var}\left(v_q^{(k)};J_q\right)+
\sum_{q\in Q} \operatorname{Var}\left(\left(v_q^{-1}\right)^{(k)};J_q\right)\right)^{1/k}\\
&\le C'
\left(K\cdot |I|\right)^{1-1/k}
\left(
K\cdot
\sum_{v\in V\cup V^{-1}}\operatorname{Var}\left(v^{(k)};I\right)
\right)^{1/k}<\infty.
\end{align*}
It follows again from Lemma~\ref{l:sum-density} that $d_\bN(Q)=0$.
\ep
\section{Proof of the Slow Progress Lemma}
We let $k,\beta,V,\{v_i\}$ and $\{N_i\}$ be as in the hypotheses of the Slow Progress Lemma, 
such that one of the following holds.
\be[(i)]
\item $V\sse\Diff_+^{k,\beta}(I)$;
\item $V\sse \Diff_+^{k,\mathrm{bv}}(I)$.
\ee

Fix $x_0\in I$, and define
\[x_n:=v_n^{N_n}\cdots v_1^{N_1}(x_0).\]
We also define
\[P:=\left\{ n\in\bN\mid \CD_V(x_0,x_{n})>\max_{0\le m< n}\CD_V(x_0,x_m)\right\}.\]
By Lemma~\ref{l:p-almost}, we have that
\[\limsup_{n\to\infty} \frac{\CD_V(x_0,x_n)}n=\overline{d}_\bN(P)=:1-\epsilon\]
for some $\epsilon\in[0,1]$. It suffices for us to show that $\epsilon\ge 1/(4k+2)$.

We may assume $\epsilon\in[0,1)$ and $P$ is infinite. By symmetry that we can also assume that the map
$p\mapsto x_p$ is order--preserving (instead of order--reversing) for $p\in P$. 
Pick $p_0\gg0$ so that $\CD_V(x_0,x_p)>2k+1$ for all $p\in P\cap [p_0,\infty)$. 
For  such a  $p$, we let $J_p$ be the unique connected component of
$\supp v_p$
 containing $[x_{p-1},x_p]$.
We have at least $k$ fixed points of $v_p$ in $(x_0,x_p)$. Indeed, if 
\[ r:=\#(\Fix v_p\cap (x_0,x_p))\]
 for some $r<k$, then $(x_0,x_p)$ minus those $r$ fixed points can be covered by $r+1$ intervals in $\pi_0\supp v_p$.
 Combining these $r+1$ intervals with $r+2$ intervals from \[\bigcup_{v\in V}\pi_0\supp v\] that cover
 $\{x_0, x_p\}$ and the $r$ fixed points, we would have a contradiction that 
\[\CD_V(x_0,x_p)\le 2r+3\le 2k+1.\]

From the previous paragraph, we can define for each $p\ge p_0$ in $P$ the smallest compact intervals
$L_p$ and $R_p$ of the form $[x, \sup J_p]$ and $[\inf J_p,y]$ respectively, such that $v_p$ is $k$--fixed on $L_p$ and on $R_p$. 
For each $\delta>0$ and $C\in\bN$, we set
\begin{align*}
P_\delta&:=\left\{p\in P\cap(p_0,\infty)\middle\vert  v_p^{N_p}\text{ is at most }\delta\text{--expansive on }L_p \text{ and }R_p\right\},\\
P_{\delta,C}&:=\bigcap_{0\le i\le C-1} (P_\delta-i).
\end{align*}

By Lemmas~\ref{l:upper-density} and~\ref{l:prob-dyn}, we have that
\[
\overline{d}_\bN(P_{\delta,C}) \ge 1- C(1-\overline{d}_\bN(P_\delta))=1-C\epsilon.\]

Set $C:=4k+2$ and $\delta:=1/C$.
Assume for contradiction that \[\epsilon<\frac{1}{4k+2}=\frac{1}{C}.\] Then $P_{\delta,C}$ is nonempty and there exists some $p\in P_{\delta,C}$.
For each $i=-1,0,\ldots,C-1$, we have that
 \[ \CD_V(x_0,x_{p+i})=\CD_V(x_0,x_p)+i.\]   
Let $0\le i\le C/2-1=2k$. Since
 \[\CD_V(x_{p+i},\sup R_{p+i})\le 2k<\CD_V(x_{p+i},x_{p+i+2k+1}),\] we have that
\[
\sup R_{p+i}<x_{p+i+2k+1}\le x_{p+C-1}.\] 
Since $v_{p+i}^{N_{p+i}}$ is at most $\delta$--expansive on $R_{p+i}$ we see 
\[ x_{p+i}-x_{p+i-1}\le \delta (\sup R_{p+i}-x_{p+i})< \delta (x_{p+C-1} - x_p).\]
Similarly for $C/2\le i\le C-1$, we have that
\[
\inf L_{p+i}>x_{p+i-1-2k-1}\ge x_{p-1},\]
and that

\[ x_{p+i}-x_{p+i-1} \le \delta (x_{p+i-1}-\inf L_{p+i}) 
< \delta (x_{p+C-2} - x_{p-1}).\]
Summing up, we obtain
\[x_{p+C-1}-x_{p-1} =\sum_{i=0}^{C-1} (x_{p+i}-x_{p+i-1})< C\delta (x_{p+C-1}-x_{p-1}).\]
We thus obtain a contradiction, since $C\delta=1$. 
This completes the proof of the Slow Progress Lemma.

Let us now record a consequence of the slow progress lemma that will be useful for us later.
\begin{cor}\label{c:slp}
Let $k,\beta$ and $\{N_i\}$ be as in the hypothesis of the Slow Progress Lemma.
We let $G$ be a group with a finite generating set $V$,
and let
\[
\psi\co G\to \Homeo_+(I)\]
be a representation which is topologically conjugate either into 
$ \Diff_+^{k,\beta}(I)$ or $\Diff_+^{k,\mathrm{bv}}(I)$. 
Pick a sequence 
$\{v_i\}_{i\ge1}$ of elements
from $V$ and set
 \[w_i := v_i^{N_i}\cdots v_1^{N_1}.\]
If $K\sse\supp\psi$ is a compact interval,
then for all sufficiently large $i\ge1$
we have that
\[
\CL(\psi(w_i)K)<2i.\]
\end{cor}
\bp
The conclusion is only concerned with the action $\psi$ restricted 
to the component of $\supp\psi$ containing $K$.
So, for brevity we may assume that $I=[0,1]$ and $\supp\psi=(0,1)$.
We set $K=[x,y]$. 
We simply write $\CD$ and $\CL$ for the covering distance and the covering length defined by the set $\psi(V)$.

Put $T:=\CD(x,y)>0$.
By the Slow Progress Lemma, for all sufficiently large $i$ 
the covering distances 
$\CD(x,\psi(w_i)x)$
and $\CD(y,\psi(w_i)y)$ will be smaller than $i-T$. 
So, we have that
\[
\CD(\psi(w_i)x,\psi(w_i)y)< 2(i-T)+\CD(x,y)=2i-T<2i,\] as required.
\ep


%
%
%
\chapter{Algebraic obstructions for general regularities}\label{ch:optimal}
\begin{abstract}
The main goal of this chapter is to construct finitely generated groups of $C^{k,\alpha}$ diffeomorphisms of a compact 
one--manifold that cannot be embedded into the group of  $C^{k,\beta}$ diffeomorphisms for $\beta$ ``sufficiently smaller'' than $\alpha$. 
There are two crucial ingredients for such a construction. One is the Slow Progress Lemma from the previous chapter, 
which roughly asserts that for finitely generated group actions by diffeomorphisms, the smoother the action the slower it
expands covering lengths of an interval.
The other is an explicit construction of a finitely generated $C^{k,\alpha}$ diffeomorphism group that expands covering 
lengths faster than the rate allowed by the Slow Progress Lemma for $C^{k,\beta}$ diffeomorphism groups.
We will establish a result for this second ingredient in the case of a compact interval, and generalize the result to circles using 
the material from preceding chapters.
We conclude the chapter with consequences regarding H\"older regularities $C^r$ for $r\in[1,\infty)$. \end{abstract}

\section{Statement of the results}\label{s:statement}
For convenience of notation, we will write
$\CM$ for the set of all concave moduli.
Recall that we write $\beta\succ_k0$ for $\beta\in\CM$ if either $k\ge2$, or $k=1$ and $\beta$ is sub-tame.
For a manifold $M$, we will often use the notation 
\[
\Diff_+^{k,\beta'}(M)\]
with $\beta'\in\{\beta,\mathrm{bv}\}$. This group 
is either $\Diff_+^{k,\beta}(M)$ or $\Diff_+^{k,\mathrm{bv}}(M)$,
depending on the choice of $\beta'$.
The precise goal of this chapter is to establish the following theorem.

\begin{thm}\label{t:optimal-all}
If $k\in{\bZ_{>0}}$ and $\alpha,\beta\in\CM$ satisfy
that $\beta\succ_k0$, and that
\[
\int_0^1 \frac1x\left(\frac{\beta(x)}{\alpha(x)}\right)^{1/k}dx<\infty,\]
then there exists a finitely generated nonabelian group 
\[ R=R(k,\alpha,\beta)\le\Diff_+^{k,\alpha}(I)\]
such that
every homomorphism
\[
[R,R]\longrightarrow\Diff_+^{k,\beta'}(M^1)\]
is trivial for $\beta'\in\{\beta,\mathrm{bv}\}$ and for $M^1\in\{I,S^1\}$.
Moreover, $R$ can be chosen so that $[R,R]$ is simple and every proper quotient of $R$ is abelian.
\end{thm}
By the last condition, we see that every homomorphism from $R$
to $\Diff_+^{k,\beta'}(M^1)$ has abelian image.

\begin{rem}
We note that the conclusion regarding 
$\Diff_+^{k,\mathrm{bv}}(M^1)$
does not involve the modulus $\beta(x)$. 
Indeed, to obtain the desired conclusion, it suffices  to pick $\beta(x)$ as small as possible; that is,
we may assume $\beta(x)=x$ and verify the hypothesis of the theorem.
\end{rem}

\begin{exmp}\label{ex:hoelder-critical}
If we consider the H\"older moduli
\[
\alpha(x)=x^r,\quad\beta(x)=x^s\]
for some $0<r<s\le1$, then the hypothesis of Theorem~\ref{t:optimal-all} is satisfied for all $k\ge1$.
We thus obtain a finitely generated subgroup of $\Diff_+^{k,r}(I)$ that never embeds into $\Diff_+^{k,s}(I)$.\end{exmp}

We will actually deduce a much stronger consequence than Example~\ref{ex:hoelder-critical} regarding $\Diff_+^{k,r}(M)$. See Section~\ref{s:hoelder} for further discussion.

\begin{cor}\label{c:optimal-all}
For each real number $r\ge1$ we have the following.
\be[(1)]
\item
There exists a finitely generated nonabelian group $G_r\le\Diff_+^r(I)$ such that every homomorphism
\[ G_r\longrightarrow \bigcup_{s>r} \Diff_+^s(M^1)\] 
has an abelian image for $M^1\in\{I,S^1\}$.
\item
There exists a finitely generated nonabelian group \[H_r\le\bigcup_{s<r}\Diff_+^s(I)\] such that every homomorphism
\[
H_r\longrightarrow \Diff_+^r(M^1)\] has abelian image for $M^1\in\{I,S^1\}$.
\ee
Furthermore, we can require that $[G_r,G_r]$ and $[H_r,H_r]$ are simple.
\end{cor}

We will often write
\[\log^s(x):=\left(\log x\right)^s.\]
The groups $R(k,\alpha,\beta)$ in Theorem~\ref{t:optimal-all} will be chosen from a family of representations $\phi(k,\alpha,\{\ell_i\})$
 of a fixed group
\[
G^\dagger:=(\bZ\times\BS(1,2))\ast F_2,\] for a suitably chosen countable set of parameters $\{\ell_i\}$.
More precisely, we will consider a set of intervals whose lengths $\{\ell_i\}$ satisfy the following mild restrictions.
\bd\label{d:tame}
A positive real sequence $\{\ell_i\}_{i\ge1}$ is said to be (a sequence of) \emph{admissible lengths}\index{admissible lengths} 
if \[\sum_i \ell_i <\infty\]
and if
\[
0< \inf \ell_{i+1}/\ell_i \le \sup \ell_{i+1}/\ell_i <\infty.\]
\ed
The reader may find it useful to keep in mind the following concrete examples:
\begin{itemize}
\item $\ell_i = 1/\left({i\log^2i}\right)$;
\item $\ell_i = 1/\left(i\log i  (\log\log i)^{1+\epsilon}\right)$ for $\epsilon>0$.
\end{itemize}
We have  $\lim_i \ell_{i+1}/\ell_i=1$ in both of the cases.

As we have noted in the introduction of this chapter, the main content of the proof for Theorem~\ref{t:optimal-all}
is the following construction of an ``optimally expanding'' diffeomorphism group.

\begin{thm}\label{t:optimal-group}
Let $k\in{\bZ_{>0}}$, let $\alpha\in\CM$, and let $\{\ell_i\}$ be admissible lengths.
Then there exists a representation
\[\phi=\phi(k,\alpha,\{\ell_i\})\co G^\dagger\longrightarrow\Diff_+^{k,\alpha}(I)\]
satisfying the following:
\begin{itemize}
\item[]
If a concave modulus $\beta\succ_k0$ satisfies
\[\beta(1/i) \cdot(1/i)^k\preccurlyeq \alpha(1/i)\cdot \ell_i^k\]
for almost all $i\ge1$,
and if $\beta'\in\{\beta,\mathrm{bv}\}$,
then for every representation of the form
\[
\psi\co G^\dagger \longrightarrow \Diff_+^{k,\beta'}(I),\]

 we have that
\[
[G^\dagger,G^\dagger]\cap(\ker\psi\setminus\ker\phi)\ne\varnothing.\]
\end{itemize}
\end{thm}

\begin{rem}
The main result of~\cite{KK2020crit} asserts the conclusion of Theorem~\ref{t:optimal-group} under the hypothesis that 
\[ \lim_{x\to+0} \frac{\beta(x)\log^K (1/x)}{\alpha(x)}=0\]
for all $K>0$.
This hypothesis obviously implies\[\beta(1/i) \cdot (1/i)^k \preccurlyeq \alpha(1/i) \cdot \ell_i^k\]after choosing admissible lengths 
of the form $\ell_i := \frac1{i\log^2 i}$. That is, Theorem~\ref{t:optimal-group} recovers the main result of~\cite{KK2020crit}.
\end{rem}
\begin{rem}\label{r:non-Lipschitz}
If $\alpha$ is the Lipschitz modulus, i.e. if $\alpha(x)=x$
then the theorem is vacuous. 
Indeed, from $\sum_i \ell_i<\infty$, we have that $\frac1{i\ell_i}$ is unbounded.
If $\beta$ satisfies the hypothesis, then 
\[
\liminf_{x\to+0} \frac{\beta(x)}{x}
\le \liminf_{i\to\infty} \frac{\beta(1/i)}{1/i}
\le C\liminf_{i\to\infty} (i\ell_i)^k=0.\]
As $\beta(x)/x$ is monotone decreasing we have that $\beta(x)\equiv0$, which is a contradiction.
\end{rem}

The majority of this chapter is devoted to the proof of Theorem~\ref{t:optimal-group}.
We will then establish Theorem~\ref{t:optimal-all} in Section~\ref{s:circle}
by applying the Rank Trick and the Chain Group Trick from Chapter~\ref{ch:chain-groups}.

\section{Strategy for the proof of Theorem~\ref{t:optimal-group}}\label{s:strategy}
The proof of Theorem~\ref{t:optimal-group} involves results from several different parts of this book
and leverages somewhat complicated interactions between group theoretic, dynamical and analytic features.
Lemma~\ref{l:covering-translation-distance} relates syllable lengths (group theory) with covering distances (dynamics).
The Slow Progress Lemma is a bridge between the growth of covering distances (dynamics) and regularity (analysis).

In this section, we will collect essential consequences of the three main tools, stated for topological actions for ease of
understanding, and explain
how the proof of Theorem~\ref{t:optimal-group} will proceed.
We hope that this will allow the reader to grasp the main idea behind the proof before we embark on explicating the details.

To begin, suppose that $G$ is a group with a finite generating set $V$.
We will assume that $G$ contains a ``universally compactly supported element'' $u_0$ in the following sense.
The reader will recall that we have seen an instance of such a group in Corollary~\ref{c:bs12}.
\be[(A)]
\setcounter{enumi}{0}
\item\label{p:ucsd} \emph{There exists a nontrivial element $u_0\in [G,G]$
that is compactly supported under every action $\rho\co G\longrightarrow\Homeo_+(I)$. }
\ee

Assume now that we are given two actions 
\[\phi,\psi\co G\longrightarrow \Homeo_+(I).\] 
Let us denote the covering lengths corresponding to $\phi(V)$ and $\psi(V)$ as
$\CL_\phi$ and $\CL_\psi$, respectively.
We will require from the construction of $\phi$ that
\be[(A)]
\setcounter{enumi}{1}
\item\label{p:nontrivial}  \emph{$\supp\phi$ is connected, and $\phi(u_0)$ is nontrivial.}
\ee

We will also make the following assumption,
which is true for instance when $\phi(G)\le\Diff_+^\infty(\Int I)$ as illustrated in Proposition~\ref{prop:disj-ab}.
\be[(A)]
\setcounter{enumi}{2}
\item\label{p:noncommute}\emph{
If the images of $g_1,g_2\in G$ under $\phi$
 are compactly supported and commute with each other, 
and if $J_i$ is a component of $\supp\phi(g_i)$ for $i=1,2$,
then either $J_1=J_2$ or $J_1\cap J_2=\varnothing$.}
\ee

We assume furthermore that the following conditions hold for some choice of elements $\{w_i\}_{i\ge1}\sse G$.
\be[(A)]
\setcounter{enumi}{3}
\item\label{p:small} \emph{For each compact nondegenerate interval $K\sse\supp\psi$
and for all sufficiently large $i\ge1$, we have
\[\CL_\psi(\psi(w_i)K)<2i.\]}
\item\label{p:large} \emph{For each compact nondegenerate interval $K\sse\supp\phi$ there exists some $g\in G$ such that
\[\CL_\phi(\phi(w_ig)K)>2i\]}
whenever $i\ge1$.
\ee

In this abstract setting (without a concrete $G$ or $\phi$) we can deduce the conclusion of Theorem~\ref{t:optimal-group}.
\begin{lem}\label{lem:optimal-abstract}
Under the hypotheses (\ref{p:ucsd}) through (\ref{p:large}), the following set is nonempty:
\[ [G,G]\cap\ker\psi\setminus\ker\phi.\]
\end{lem}

A key idea of the proof is that if $G$ acts on an open interval
 the syllable length $\|g\|_{\mathrm{syl}}$ of $g\in G$ 
can be understood as the ``displacement energy'' of $g$,
which gives an upper bound for the covering length of an interval that can be moved off itself by $\rho(g)$. 
A precise statement of this observation is given in Lemma~\ref{l:covering-translation-distance}.
\bp[Proof of Lemma~\ref{lem:optimal-abstract}]
We may assume $\psi(u_0)\ne1$ for, otherwise there is nothing to show
by condition~(\ref{p:nontrivial}).
Set $u_1:=u_0$.
By condition~(\ref{p:ucsd}), we can find a minimal, finite collection
\[ U_1, U_2, \ldots, U_m\]
of components of $\supp\psi$ so that
\[{\suppc \psi(u_1)}\sse U_1\cup\cdots\cup U_m.\]

Let $J_1$ be the closure of a component of  $\supp\phi(u_1)$.
Using condition~(\ref{p:large}), we can find 
some $g_1\in G$ such that
\[\CL_\phi(\phi(w_ig)J_1)>2i\]
for all $i\ge1$. Setting $u_1':=g_1u_1g_1^{-1}$ we see that 
$J_1':=\phi(g)J_1$ is the closure of a component of $\phi(u_1')$.

Since each $U_i$ is $\psi(G)$--invariant, the homeomorphism $\psi(u_1')$ is still
compactly supported in $\bigcup_i U_i$. 
Pick a compact interval $K_1$ such that
\[\supp\psi(u_1')\cap U_1\sse K_1\sse U_1.\]
Applying condition~(\ref{p:small}) we can fix some $i\gg0$ such that
\[\CL_\psi(\psi(w_i)K_1)<2i.\]
By Lemma~\ref{l:covering-translation-distance}, there exists some $h_1\in G$ with syllable length less than $2i$
such that $\psi(h_1)$ moves $\psi(w_i)K_1$ off itself.

If $\phi(h_1w_i)(J_1')$ is not equal to $\phi(w_i)(J_1')$,
then we set $v=1$; otherwise,
pick some $v\in V$ such that $\partial \phi(h_1w_i)(J_1')\cap \supp\phi(v)\ne\varnothing$ .
For this choice of  $v\in\{1\}\cup V$, we have that
\[\phi(v^{\pm1}h_1)\phi(w_i)(J_1')\ne \phi(w_i)(J_1').\]

We claim that some $s\in\{-1,1\}$ satisfies
\[\psi(v^sh_1w_i)K_1\cap \psi(w_i)K_1=\varnothing.\]
Indeed, consider the case when
\[p:=\sup \psi(h_1w_i)K_1<\inf \psi(w_i)K_1.\]
For some $s\in\{-1,1\}$, the map $\psi(v^{s})$ will move $p$ to the left or stabilize it.
In this case, the interval
$\psi(v^{s} h_1w_i)K_1$ will still be disjoint from $ \psi(w_i)K_1$.
The remaining cases can be treated in the same manner, and thus the claim is proved.

It follows from the claim that 
\[U_1\cap \supp\psi\left[w_iu_1'w_i^{-1},h_1'w_iu_1'(h_1'w_i)^{-1}\right]=\varnothing\]
for $h_1':=v^s h_1$.

Since $\|h_1'\|_{\mathrm{syl}}\le 2i$, 
Lemma~\ref{l:covering-translation-distance} implies that $\phi(h_1')$ cannot move $\phi(w_i)(\Int J_1')$ off itself.
Note that  $\phi(w_i)(\Int J_1')$  is a component of $\supp\phi(w_iu_1'w_i^{-1})$.
So, we see from condition~(\ref{p:noncommute}) that the image of
\[u_2:=
\left[w_iu_1'w_i^{-1}, 
\left(h_1'w_i \right) u_1' \left(h_1'w_i \right)^{-1}\right]\in [G,G]\]
under the action $\phi$ is nontrivial.
Note also that $\phi(u_2)$ and $\psi(u_2)$ are still compactly supported
since $u_2\in\fform{u_1}$.

Since $\psi(u_2)$ acts trivially on $U_1$, 
its support is now contained in $U_2\cup\cdots \cup U_m$.
After rearranging if necessary, we assume that $\psi(u_2)$ acts nontrivally on $U_2$.
Proceeding as above, we continue producing
\[J_2, g_2,u_2', J_2',K_2,h_2,h_2',\]
in this order, using the given hypotheses. 
We then find an element $u_3\in [G,G]\setminus\ker\phi$
such that the support of $\psi(u_3)$ is precompact in $U_3\cup\cdots\cup U_m$.
This process terminates in $i\le m$ steps, and we eventually find 
some element \[u_{i+1}\in ( [G,G]\cap\ker\psi)\setminus\ker\phi.\]
Note that the final word $u_{i+1}$ in the above proof belongs $\fform{u_0}$.

\ep

The preceding proof is yet another instance of finding a kernel element of a representation by taking successive
commutators; see the beginning of Subsection~\ref{ss:compact supports} for  relevant remarks.

In this chapter, we will prove Theorem~\ref{t:optimal-group} by
verifying conditions (A) through (E) of the above lemma
for the group
\[
G:=G^\dagger=(\bZ\times \BS(1,2))\ast F_2=
(\form{c}\times\form{a,e\mid aea^{-1}=e^2})\ast\form{b,d}.\]
Condition~(\ref{p:ucsd}) holds for
 \[u_0:=\left[
[c^d,ee^de^{-1}],c\right]\in G^\dagger\]
by  Corollary~\ref{c:bs12}.
Under the hypotheses of the theorem, we let
\[ N_i :=\ceil*{\frac1{\ell_i^{k-1}\alpha(\ell_i)}},\]
and set
\[ w_i:= v_i^{N_i}\cdots v_2^{N_2}v_1^{N_1}\in G^\dagger\]
for \[a=v_1=v_3=\cdots\quad \textrm{and}\quad b=v_2=v_4=\cdots.\]
Let $\beta$ and $\psi$ be as given in the theorem.
Note from  Lemma~\ref{l:sum-density}  that the set \[P:=\{i\in\bN\mid i\ell_i\le 1\}\]
has natural density one.
For each $i\in P$,  we have 
\[\frac{\ell_i}{\alpha(\ell_i)}\le\frac{1/i}{\alpha(1/i)}.\]

It follows for $i\in P$ that
\[
 \frac{N_i\beta(1/i)}{i^{k-1}}\preccurlyeq \left( \frac{\beta(1/i)}{(i\ell_i)^k\alpha(1/i)}\cdot 
 \frac{i\ell_i\alpha(1/i)}{\alpha(\ell_i)}\right)\le \left( \frac{\beta(1/i)}{(i\ell_i)^k\alpha(1/i)}\right).
\]
Therefore, we have
\[
N_i\preccurlyeq i^{k-1}/\beta(1/i)\quad\text{ for almost all }i.\]
So, we can apply Corollary~\ref{c:slp} and see that condition~(\ref{p:small}) holds.

From the next section on, we will construct a representation
\[\phi\co G^\dagger\longrightarrow\Diff_+^{k,\alpha}(I)\cap\Diff_+^\infty(\Int I).\]
In particular, condition (\ref{p:noncommute}) will be automatic by Proposition~\ref{prop:disj-ab}.
It will then remain to engineer $\phi$ so that
the two conditions (\ref{p:nontrivial})  and (\ref{p:large}) also hold, which will finish the proof of the theorem. 
We informally call $\phi(G^\dagger)$ as an \emph{optimally expanding group}\index{optimally expanding group}
because of condition (\ref{p:large}), in contrast with (\ref{p:small}).

\section{A single diffeomorphism of optimal expansion}\label{s:optimal-single} 
In order to find a representation $\phi$ satisfying  condition~(\ref{p:large}), 
we will describe a method of constructing a single diffeomorphism that moves ``sufficiently fast'' on each supporting intervals.
Let us fix a constant
 \begin{equation}\delta_0\in[9/10,1),\end{equation}
and set $D_0:=(1-\delta_0)/2\in(0,1/20]$.

The following condition is a key dynamical feature that will provide us a necessary condition for the regularity of a diffeomoprhism.
\bd\label{d:fast}
Let $f\co J\longrightarrow J$ be a homeomorphism of an interval $J$.
We say $f$ is \emph{$\delta$--fast}\index{$\delta$--fast} for some $\delta$ if there exists some $x_0\in I$ such that $|f(x_0)-x_0|\ge \delta |J|$.
\ed
This notion is stronger than the $\delta'$--expansiveness (Definition~\ref{d:expansive}) for some $\delta'>0$.

\begin{lem}\label{l:fast-expansive}
If $\delta\in(0,1)$,
then a $\delta$--fast homeomorphism of an interval is $\delta'$--expansive
for 
\[
\delta':=2\delta/(1-\delta).\]
\end{lem}
\bp
Let $f$ denote a $\delta$--fast homeomorphism of an interval $J$.
We may assume that $f(x_0)-x_0 \ge \delta |J|$ for some $x_0\in J$. 
For brevity, let us write
\[
A := x_0-\inf J, \quad
B :=f(x_0)-x_0,\quad
C := \sup J - f(x_0).\]
Then we have that $B/(A+B+C)\ge \delta$. By elementary computation, 
we deduce that 
\[ \max \left(B/A, B/C\right)\ge 2/\left(1/\delta-1\right).\]
This implies that $f$ is $\delta'$--expansive.
\ep

The starting point of our construction is to explicitly find a single $C^{k,\alpha}$--diffeomorphism that is $\delta_0$--fast, with
suitably bounded $C^k$ and $C^{k,\alpha}$ norms. Then we will suitably concatenate such diffeomorphisms in order to obtain
sequence of diffeomorphisms that can expand an arbitrary open interval intersecting the union of supports. 
Recall our notation that for an interval $J\sse \bR$, we write
\[
\Diff_J^{k,\alpha}(\bR):=\{f\in \Diff_+^{k,\alpha}(\bR)\mid \supp f\sse J\}.\]
If $J$ is bounded, one may identify each element of $\Diff_J^{k,\alpha}(\bR)$ with a $C^{k,\alpha}$ 
diffeomorphism of $J$ that is $C^k$--tangent to the identity at $\partial J$.

Pick a smooth function $\Phi\co \bR\longrightarrow\bR$ satisfying the following conditions:
\begin{itemize}
\item $\Phi(x)=-1$ for $x\le 0$.
\item $\Phi(x)+\Phi(1-x)=0$; in particular,  $\Phi(x)=1$ for $x\ge1$ and $\Phi(1/2)=0$.
\item $\Phi(x)$ is strictly increasing on $(0,1)$.
\end{itemize}
The precise choice of $\Phi$ will not matter for us. For instance, one may construct $\Phi$ from a 
smooth bump function supported in $(0,1)$ as follows:
\begin{equation}
\Phi(x) := 2\frac{\int_{-\infty}^x e^{-1/(t(1-t))} \chi_{(0,1)}(t)dt}{\int_{\bR} e^{-1/(t(1-t))} \chi_{(0,1)}(t)dt} -1.
\end{equation}

For a function $f$, we let $\|f\|$ denote its uniform norm.
Let $r\ge0$ be an integer.
If $f$ is a real-valued $C^r$--map defined on some set $U$, then  its $C^r$--norm (or $C^r$--metric) is 
\[
\| f\|_{C^r}:=\sup_{0\le i\le r} \abss*{f^{(i)}}.\]
We define the $C^{r,\alpha}$--norm of $f$ by
\[[f]_{r,\alpha}:=\left[f^{(r)}\right]_\alpha=\sup_{x\ne y\text{ in }U}\frac{ |f(x)-f(y)|}{\alpha(|x-y|)}.\]

We define a constant
\begin{equation}
K_0=K_0(k):=\frac1{D_0^{k+1}}\abss*{\Phi}_{C^{k+1}}\end{equation}
and pick $\ell_0=\ell_0(k,\alpha,\delta_0)>0$ such that 
\begin{equation}\ell_0+K_0\alpha(\ell_0)\le 1.\end{equation}

Let us consider an arbitrary $\ell\in(0,\ell_0]$. 
We put $\Delta:=\delta_0\ell^k\alpha(\ell)$, and 
define 
\[
g_\ell(x) := \frac\Delta2 \left( \Phi\left(\frac{t}{D_0\ell}\right)+\Phi\left(\frac{\ell-t}{D_0\ell}\right)\right).\]
Then we have the following.
\be[(i)]
\item $g_\ell(x)=0$ outside $(0,\ell)$;
\item $g_\ell(x)=\Delta$ on $[D_0\ell,(1-D_0)\ell]$;
\item $g_\ell$ is strictly increasing on $(0,D_0\ell)$ and strictly decreasing on $((1-D_0)\ell,\ell)$;
\item\label{p:gkal} $\left[{g_\ell^{(k)}}\right]_\alpha\le K_0$.
\ee
The observations (i) through (iii) are immediate. 
To establish (\ref{p:gkal}), we first note that
\begin{align*}
\left[g_\ell^{(k)}\right]_\alpha
&\le \Delta \left[ \left(\Phi\left(\frac{t}{D_0\ell}\right)\right)^{(k)}\right]_\alpha
= \frac{\Delta}{(D_0\ell)^k}
\sup_{0\le x<y\le D_0\ell} \frac{\abs*{\Phi^{(k)}\left(\frac{y}{D_0\ell}\right)-\Phi^{(k)}\left(\frac{x}{D_0\ell}\right)}}{\alpha(y-x)}\\
&\le
\frac{\delta_0 \alpha(\ell)}{D_0^k} \cdot
\abss*{\Phi^{(k+1)}}\cdot
\sup_{0\le x<y\le D_0\ell} \ \frac{y-x}{D_0\ell \alpha(y-x)}.\end{align*}
Using the fact that $x/\alpha(x)$ is monotone increasing, we have that
\[
\left[g_\ell^{(k)}\right]_\alpha
\le 
\frac{\alpha(\ell)}{D_0^k\alpha(D_0\ell)} 
\abss*{\Phi^{(k+1)}}
\le
\frac{1}{D_0^{k+1}}
\abss*{\Phi^{(k+1)}}=K_0.\]

Using $g_\ell$, let us construct a $C^{k,\alpha}$ diffeomorphism supported on a single interval that is $\delta_0$--fast.

\begin{lem}\label{l:optimal-diffeo}
For each $\ell\in(0,\ell_0]$ there exists some $f=f_\ell\in\Diff_+^\infty(\bR)$ such that the following hold.
\be[(i)]
\item $f(x)=\Id$ outside $(0,\ell)$;
\item $f(x)>x$ on $(0,\ell)$;
\item\label{p:fkal} $\left[f\right]_{k,\alpha}\le K_0$;
\item\label{p:ck} $\|f-\Id\|_{C^k}\le K_0\alpha(\ell)$;
\item\label{p:fnfast} $f^N$ is $\delta_0$--fast for all $N\ge 1 / (\ell^{k-1}\alpha(\ell))$.
\ee
\end{lem}

\bp
Setting $f:=g_\ell+\Id$,
we immediately see the conditions (i) through (iii).
The condition~(\ref{p:ck}) is actually a consequence of (iii). Indeed, for each $x\in[0,\ell]$ we have  
\[\left|f^{(k)}(x)-\Id^{(k)}(x)\right|=\left|f^{(k)}(x)-f^{(k)}(0)\right|\le K_0\alpha(\ell).\]
For $1\le i<k$, we have
\begin{align*}
|f^{(i)}(x)-\Id^{(i)}(x)|&\le
 \int_{s_1=0}^x\int_{s_2=0}^{s_1} \cdots\int_{s_{k-i}=0}^{s_{k-i-1}} |f^{(k)}(s_{k-i})-\Id^{(k)}(s_{k-i})|
 ds_{k-i} \cdots ds_1\\
&\le
\ell^{k-i} K_0\alpha(\ell)\le K_0\alpha(\ell).
 \end{align*}
We see that  $f$ is a diffeomorphism since
\[
f'(x)\ge 1- \|f'-1\|\ge 1 - K_0\alpha(\ell)>0.\]

It only remains to show part (\ref{p:fnfast}). 
We first set
\[
N_0:=\left\lceil \frac1{\ell^{k-1}\alpha(\ell)}\right\rceil=\left\lceil \frac{(1-2D_0)\ell}{\Delta}\right\rceil.\]
We then have that 
\[
D_0\ell+(N_0-1)\Delta<(1-D_0)\ell\le D_0\ell+N_0\Delta.\]

Using an induction on $i=1,\ldots, N_0$, we see that 
\[
f^i(D_0\ell) = g\circ f^{i-1}(D_0\ell) + f^{i-1}(D_0\ell)=\Delta+ f^{i-1}(D_0\ell)=D_0\ell+ i\Delta.\]
For all $N\ge N_0$, it follows that
\[
f^N(D_0\ell)-D_0\ell
\ge 
f^{N_0}(D_0\ell)-D_0\ell= N_0\Delta\ge (1-2D_0)\ell =\delta_0\ell.\]
This proves part (\ref{p:fnfast}).
\ep

The following general observation is an easy consequence of the concavity of $\alpha$.
Here, we emphasize that $\prod_{i=1}^m f_i$ denotes a composition of diffeomorphisms.
\begin{lem}\label{l:disj-kal}
Let $\{J_i\}_{1\le i\le m}$ be a disjoint collection of compact intervals in $\bR$.
If $f_1,\ldots,f_m$ are orientation-preserving $C^{k,\alpha}$ diffeomorphisms of the real line 
satisfying $\supp f_i\sse J_i$ for each $i$,  then
\[
\left[\prod_{i=1}^m f_i\right]_{k,\alpha}\le 2\sup_i [f_i]_{k,\alpha}.\]
\end{lem}
\bp
Let us set $F:=\prod_{i=1}^m f_i$ and $K:=\sup_i[f_i]_{k,\alpha}$. 
Pick real numbers $x<y$. It suffices for us to show that
\[ \abs*{ F^{(k)}(y) -  F^{(k)}(x)} \le 2K\alpha(y-x).\]
This is immediate when $x,y\in J_i$ for some $i$. 
Suppose $x\in J_i$ and $y\in J_j$ for some $i\ne j$.
We can find $x_0\in\partial J_i$ and $y_0\in\partial J_j$ such that
\[x\le x_0<y_0\le y.\]
Since $f_t=\Id$ at $\partial J_t$ for each $t$, we have
\begin{align*}
\abs*{ F^{(k)}(y) -  F^{(k)}(x)} 
&\le\abs*{ f_j^{(k)}(y) -  f_j^{(k)}(y_0)} +\abs*{ \Id^{(k)}(y_0) -  \Id^{(k)}(x_0)} +\abs*{ f_i^{(k)}(x_0) -  f_i^{(k)}(x)}\\
&\le K\alpha(y-y_0)+K\alpha(x_0-x)\le 2K\alpha(y-x).
\end{align*}
If $x\in J_i$ and $y\not\in \bigcup_j J_j$ then we can find $x_0\in \partial J_i$ such that
$x\le x_0<y$. Then we have 
\[ \abs*{ F^{(k)}(y) -  F^{(k)}(x)} 
=\abs*{ \Id^{(k)}(y) -  F^{(k)}(x)} 
=\abs*{F^{(k)}(x_0) -  F^{(k)}(x)} 
\le K\alpha(y-x).\]
\ep


The following is a standard consequence of Arzel\`a--Ascoli theorem. We briefly sketch the proof 
for readers' convenience; a detailed proof can be found in~\cite{CKK2019}. 
\begin{lem}[cf. {\cite[Lemma 4.19]{CKK2019}}]\label{l:ck-compact}
Let $I\sse\bR$ be a compact interval, let $r\in(0,\infty)$ and let $\epsilon\in(0,1)$.
Then the set
\[
\{f\in \Diff_I^{k,\alpha}(\bR)\mid  [f]_{k,\alpha}\le r\text{ and }\|f'-1\|\le \epsilon\}\]
is isometric to a compact convex subset of the Banach space
\[C_I^{k,\alpha}(\bR):=\{g\in C^{k,\alpha}(\bR) \mid g(x)=0\text{ for }x\not\in I\},\]
both equipped with the $C^k$--metrics.
\end{lem}
\bp
Let us set
\[
B:=\{g \in C_I^{k,\alpha}(\bR)\mid [g]_{k,\alpha}\le r\text{ and }\|g'\|\le\epsilon\}.\]
Then Arzel\`a--Ascoli theorem implies that $B$ is compact with respect to the $C^k$--topology in the Banach space of 
$C^k$--maps supported on $I$.
In particular, the map $f\mapsto f-\Id$ isometrically maps
the first set given in the hypothesis
 onto the $C^k$--compact convex set $B$. We note that the condition $\|f'-1\|\le\epsilon$ was used to guarantee that this map is surjective.
\ep

We can now proceed to the main construction of this section.

\begin{thm}\label{t:optimal-diffeo}
Let $k\in{\bZ_{>0}}$, let $\alpha\in\CM$,
and let $\delta_0\in[9/10,1)$.
If $\{J_i\}_{i\ge1}$ is a collection of disjoint compact intervals contained in some compact subset of $\bR$,
then there exists $f\in \Diff_+^{k,\alpha}(\bR)$ supported in $\bigcup_i J_i$
such that the following hold for all $i$:
\be[(i)]
\item $f(x)>x$ for $x\in J_i\setminus\partial J_i$;
\item $f^N$ is $\delta_0$--fast on $J_i$ if $N\ge 1/\left(|J_i|^{k-1}\alpha(|J_i|)\right)$;
\item If an open interval $U\sse \bR$ intersects finitely many $J_i$'s, then $f$ is $C^\infty$ on $U$.
\ee
\end{thm}

\bp
We first reduce the proof to the case when $|J_i|\le \ell_0$ for all $i\ge1$. 
Pick a sufficiently large $i_0$ such that $|J_i|\le \ell_0$ for $i\ge i_0$.
It is trivial that we can find $F_0\in \Diff_+^{\infty}(\bR)$ supported in $\bigcup_{i< i_0} J_i$ 
such that the conditions (i) through (iii) hold for $i< i_0$.
If we can find $F$ supported in $\bigcup_{i\ge i_0} J_i$ such that 
the conditions (i) through (iii) hold for $i\ge i_0$, then $f:=F\circ F_0$ will satisfy the desired conclusion.
So, we may now assume that $\{J_i\}$ is an infinite collection and that $i_0=1$.

Let us consider a compact interval $I$ the interior of which contains the closure of $\bigcup_i J_i$.
Using Lemma~\ref{l:optimal-diffeo}, we can find $f_i\in\Diff_+^\infty(\bR)$ supported in $J_i$ such that 
\begin{itemize}
\item $f_i(x)>x$ on $J_i\setminus\partial J_i$;
\item $[f_i]_{k,\alpha}\le K_0$;
\item $f_i^N$ is $\delta_0$--fast for $N\ge 1/ \left( |J_i|^{k-1}\alpha(|J_i|)\right)$.
\end{itemize}
Since $\{J_i\}$ is a disjoint collection, we have a well-defined homeomorphism
\[
f:=\prod_{i=1}^\infty f_i. \]
The conditions (i) and (ii) are obvious for $f$. 

Let us now define a sequence $\{F_m\}$ of $C^\infty$ diffeomorphisms supported in $I$ as follows.
\[
F_m:= \prod_{1\le i\le m} f_i.\]
The condition (iii) follows from that $f=F_m$ on $U$ for some $m\gg0$.

By Lemma~\ref{l:disj-kal} we have
\[[F_m]_{k,\alpha}\le 2K_0.\]
We also have that
\[
\|F_m'-1\|=\sup_{i\le m} \|f_i'-1\|\le K_0\alpha(\ell_0)<1.\]

By Lemma~\ref{l:ck-compact}, it follows that $F_m$ converges to $f$ in the $C^k$--metric, and that $f\in \Diff_+^{k,\alpha}(I)$. 
Since the closure of $\bigcup_i J_i$ is contained in the interior of $I$, we can extend $f$ outside $I$ by the identity. 
In particular, we have $f\in \Diff_+^{k,\alpha}(\bR)$. 
\ep

\section{Construction of optimally expanding diffeomorphism groups}\label{s:optimal}
As in the hypotheses of Theorem~\ref{t:optimal-group}, let us fix 
 $k\in{\bZ_{>0}}$, $\alpha\in\CM$,
and a sequence of admissible lengths $\{\ell_i\}_{i\ge1}$. 
We will now construct a representation $\phi$ as promised in Section~\ref{s:strategy}.

There are positive constants $c_1, c_2$ satisfying
\[
c_1\le \ell_{i+1}/\ell_i \le c_2\]
for all $i\ge1$.
We let $\kappa>0$ be small enough that $\kappa\ell_1<1$
and that
\[\frac{\ell_i}{\ell_i+\ell_{i+1}}\ge \frac1{1+c_2}>\kappa\]
for all $i\ge1$.
We then pick a constant $\delta_0\in[9/10,1)$ sufficiently near from $1$ so that
\[ \frac{\ell_{i+1}}{\ell_i+\ell_{i+1}}\kappa \ge \frac{c_1\kappa}{1+c_1}>1-\delta_0,\]
and we set \[ N_i:=\ceil*{\frac1{\ell_i^{k-1}\alpha(\ell_i)}}.\]

\subsection{Specifying the images of generators}\label{ss:specifying}
The symbols $a,b,c,d,e$ will mean the corresponding generators of 
\[G^\dagger=(\bZ\times \BS(1,2))\ast F_2=(\form{c}\times\form{a,e\mid aea^{-1}=e^2})\ast\form{b,d}.\]
The images of the generators under $\phi$ will be supported in some intervals from a collection $\FF$ defined below. See Figure~\ref{fig:config} also.

\begin{figure}[h!]
{
\tikzstyle {Av}=[blue,draw,shape=circle,fill=red,inner sep=1pt]
\tikzstyle {Bv}=[red,draw,shape=circle,fill=blue,inner sep=1pt]
\tikzstyle {Cv}=[brown,draw,shape=circle,fill=Maroon,inner sep=1pt]
\tikzstyle {Dv}=[teal,draw,shape=circle,fill=PineGreen,inner sep=1pt]
\tikzstyle {Ev}=[violet,draw,shape=circle,fill=Plum,inner sep=1pt]
\tikzstyle {A}=[blue,postaction=decorate,decoration={%
    markings,%
    mark=at position .7 with {\arrow[blue]{stealth};}}]
\tikzstyle {B}=[red,postaction=decorate,decoration={%
    markings,%
    mark=at position .7 with {\arrow[red]{stealth};}}]
\tikzstyle {C}=[brown,postaction=decorate,decoration={%
    markings,%
    mark=at position .7 with {\arrow[brown]{stealth};}}]
\tikzstyle {D}=[teal,postaction=decorate,decoration={%
    markings,%
    mark=at position .7 with {\arrow[teal]{stealth};}}]
\tikzstyle {E}=[violet,postaction=decorate,decoration={%
    markings,%
    mark=at position .7 with {\arrow[violet]{stealth};}}]
\begin{tikzpicture}[ultra thick,scale=.38]
\path (1,0) edge [A] node  {}  (4,0);
\draw (2.5,0) node [below] {\small $b$}; 

\path (3,1) edge [C] node  {}   (6,1);
\draw (4.5,1) node [above] {\small $c$};

\path (5,0)  edge [D] node  {}   (8,0);
\draw (6.5,0) node [below] {\small $d$}; 
\path    (7,1)  edge [B] node  {}   (10,1);
\draw (8.5,1) node [above] {\small $a$};

\path  (9,0)  edge [A] node  {}   (12,0);
\draw (10.5,0) node [below] {\small $b$};

\path  (11,1)  edge [B] node  {}   (14,1);
\draw (12.5,1) node [above] {\small $a$};
\path  (13,0)  edge [A] node  {}   (16,0);
\draw (14.5,0) node [below] {\small $b$};

\path   (-1.7,1)  edge node  {}   (1.7,1);
\draw (0,1) node [above] {\small $a,c,d,e$};

\path   (-1,0)  edge [A] node  {}   (-4,0);
\draw (-2.5,0) node [below] {\small $b$}; 
\path   (-3,1)  edge [C] node  {}   (-6,1);
\draw (-4.5,1) node [above] {\small $c$};
\draw (16,.2) node [above] {$\cdots$};
\draw (-16,.2) node [above] {$\cdots$};

\path (-5,0)  edge [D] node  {}   (-8,0);
\draw (-6.5,0) node [below] {\small $d$}; 

\path  (-7,1)  edge [B] node  {}   (-10,1);
\draw (-8.5,1) node [above] {\small $a$};

\path  (-9,0)  edge [A] node  {}   (-12,0);
\draw (-10.5,0) node [below] {\small $b$};

\path  (-11,1)  edge [B] node  {}   (-14,1);
\draw (-12.5,1) node [above] {\small $a$};
\path  (-13,0)  edge [A] node  {}   (-16,0);
\draw (-14.5,0) node [below] {\small $b$};

\draw  [thin,black] (1.7,1) -- (1.7,2) node [above,black] {\tiny $3/4$};
\draw  [thin,black] (-1.7,1) -- (-1.7,2) node [above,black] {\tiny $-3/4$};

\end{tikzpicture}}
\caption{The infinite chain
 $\mathcal{F}$
 and the diffeomorphism supported in each interval of  $\mathcal{F}$. }
\label{fig:config}
\end{figure}
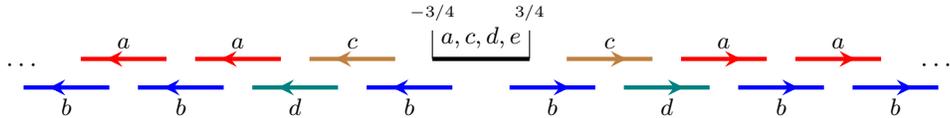

\begin{claim}
There exists an infinite chain of open intervals
\[ \FF= \{ \ldots, L_{-2},L_{-1},J_{-3},J_{-2},J_{-1},J_0,J_1,J_2,J_3,L_1,L_2,\ldots\},\]
in the given order, specified by the following conditions:
\begin{itemize}
\item $J_i=(i-3/4,i+3/4)$ for $-3\le i\le 3$;
\item $|J_3\cap L_1|=\kappa\ell_1$;
\item $|L_i\cap L_{i+1}|=\kappa\ell_{i+1}$ for $i\ge1$;
\item $L_{-i}=L_i$ for $i\ge1$.
\end{itemize}
\end{claim}
\bp
It is clear that $\{J_i\}_{-3\le i\le 3}$ forms a 7-chain. 
To see that $\FF$ is a chain, it remains to verify the following conditions.
\be[(i)]
\item $|J_2\cap J_3|+|J_3\cap L_1|<|J_3|$;
\item $|J_3\cap L_1|+|L_1\cap L_2|<|L_1|$;
\item $|L_i\cap L_{i+1}|+|L_{i+1}\cap L_{i+2}|<|L_{i+1}|$ for $i\ge1$.
\ee
The inequality (i) follows from \[ 1/2 + \kappa\ell_1 < 1/2+1=|J_3|.\]
The inequalities (ii) and (iii) are consequences of the bound
\[
\kappa\ell_i + \kappa \ell_{i+1}<\ell_i.\]\ep

We now define $I$ to be the closure $\bigcup\FF$, which is a compact interval in $\bR$.
Recall we have chosen in Section~\ref{s:strategy} that
\[u_0:=\left[
[c^d,ee^de^{-1}],c\right]\in \form{a,c,d,e}.\]
Using Theorem~\ref{thm:bs-cpt} and Lemma~\ref{l:cinfty-free-product},
we can find a representation
\[
\phi_0\co G^\dagger\longrightarrow \Diff_{J_0}^\infty(\bR)\]
such that $\phi_0(b)=1$
and such that 
$\phi_0(u_0)\ne1$. We can further require that $\supp\phi_0=J_0$, verifying condition~(\ref{p:nontrivial}) in Section~\ref{s:strategy}.

Let  $b^+\in\Diff_{J_1}^\infty(\bR)$  satisfy that $b^+(x)>x$ for $x\in J_1$.
Using Corollary~\ref{cor:ghys-serg}, we can find $c^+\in \Diff_{J_2}^\infty(\bR)$ and 
$d^+\in\Diff_{J_3}^\infty(\bR)$ such that the following hold:
\begin{itemize}
\item $c^+(x)>x$ for $x\in J_2$;
\item $d^+(x)>x$ for $x\in J_3$;
\item the action of $\form{c^+,d^+}$ on $J_2\cup J_3$ is topologically conjugate to the standard action of Thompson's group $F$ on $(0,1)$.
\end{itemize}

We define $b^-,c^-,d^-$ by symmetry
\[ 
v^-(x)=-v^+(-x)\]
for $v\in\{b,c,d\}$ and $x\in \bR$. We have another representation
\[
\phi_1\co G^\dagger\longrightarrow \Diff_I^\infty(\bR)\]
defined by $\phi_1(a)=\phi_1(e)=1$
and by 
\[
\phi_1(v)=v^+v^-\]
for $v\in\{b,c,d\}$.

We apply Theorem~\ref{t:optimal-diffeo}
to the intervals
\[ L_1,L_3,L_5,\ldots\]
to obtain \[a^+\in \Diff_I^{k,\alpha}(\bR)\cap \Diff_+^\infty( I\setminus\partial I)\]
such that 
$a^+(x)>x$ for each $x\in \bigcup_{i\ge1} L_{2i-1}$
and $a^+(x)=x$ otherwise; moreover, we require that
$\left(a^+\right)^{N_{2i-1}}$ is $\delta_0$--fast on $L_{2i-1}$ for all $i\ge1$. 
We define $a^-(x)=-a^+(-x)$.
Similarly, we define \[b^+\in \Diff_I^{k,\alpha}(\bR)\cap \Diff_+^\infty( I\setminus\partial I)\] by applying the same theorem to the intervals 
$\{L_{2i}\}_{i\ge1}$ and set $b^-(x)=-b^+(-x)$. Then we define
\[
\phi_2\co G^\dagger\longrightarrow  \Diff_I^\infty(\bR)\]
by $\phi_2\form{c,d,e}=1$
and by
\[
\phi_2(v)=v^+v^-\]
for $v\in\{a,b\}$.

Summing up, we have a representation
\[
\phi\co G^\dagger\longrightarrow  \Diff_I^\infty(\bR)\]
defined by 
\[\phi(v)=\phi_0(v)\phi_1(v)\phi_2(v)\]
for all $v\in\{a,b,c,d,e\}$. For all such $v$, the supports of $\phi_0(v), \phi_1(v), \phi_2(v)$ are all disjoint and so, the representation $\phi$ is well-defined. 
After rescaling if necessary, we may assume $I=[0,1]$.


\subsection{Verifying the optimal expansion of $\phi$}\label{ss:top-smooth}
As above, we let $v_{2i-1}:=a$ and  $v_{2i}:=b$ and for $i\ge1$. 
Since the action of $\phi(v_i)^{N_i}$ is $\delta_0$--fast on $L_i$, it is 
$2\delta_0/(1-\delta_0)$--expansive by Lemma~\ref{l:fast-expansive}.
We have already seen that
\[
N_i\preccurlyeq i^{k-1}/\beta(1/i)\quad\text{ for almost all }i.\]
By applying Lemma~\ref{l:prob-dyn} with $V=\{\phi(a)\}$ or $V=\{\phi(b)\}$,
we deduce that
the elements  \[\phi(a),\phi(b)\not\in \Diff_+^{k,\beta}(I)\cup\Diff_+^{k,\mathrm{bv}}(I).\]

We will now use covering distance estimates as dynamical obstructions that prohibits $\phi(G^\dagger)$ from being topologically conjugate into $\Diff_+^{k,\beta}(I)$
or into $\Diff_+^{k,\mathrm{bv}}(I)$.
Let us first make an easy general observation.
\begin{lem}\label{l:1-delta}
Let $\ell>0$ and $\delta\in(0,1)$.
Suppose $f\in\Homeo_+(0,\ell)$ is $\delta$--fast
and satisfies $f(x)\ge x$ for all $x\in(0,\ell)$.
Then for each \[ s\in((1-\delta)\ell,\ell),\] we have that
$f(s)>\delta\ell$.
\end{lem}
\bp
Pick $z\in(0,\ell)$ so that $f(z)-z \ge\delta\ell$. We have
\[
z\le f(z)-\delta\ell < (1-\delta)\ell < s.\]
It follows that $f(s)>f(z)\ge z+\delta\ell\ge \delta\ell$.
\ep

Returning to the proof of Theorem~\ref{t:optimal-group}, 
we will from now on define \[
w_i := v_i^{N_i}\cdots v_2^{N_2} v_1^{N_1}\in G^\dagger.\]

\begin{lem}\label{l:jumps}
If $s_1\in I$ satisfies
\[\inf L_1+(1-\delta_0)\ell_1 < s_1 < \inf L_1+\kappa\ell_1,\]
then for all $i\ge1$ we have that $\phi(w_i)(s_1)\in L_i\cap L_{i+1}$.
\end{lem}
\bp
Set $s_{i+1}:=\phi(w_i)s_1$ for $i\ge1$. We will prove a stronger claim that
\[ s_i - \inf L_i\in ( (1-\delta_0)\ell_i, \kappa\ell_i)\]
for all $i\ge1$. The case $i=1$ is obvious from the hypothesis. 
Suppose the claim is proved for $i\ge1$. 
Since $\phi(v_i)^{N_i}$ is $\delta_0$--fast on $L_i$,
we can apply Lemma~\ref{l:1-delta} to see that 
\begin{align*}
s_{i+1}&=\phi(v_i)^{N_i} s_i > \inf L_i +\delta\ell_i =\sup L_i-(1-\delta_0)\ell_i\\
&=\inf L_{i+1}+\kappa\ell_{i+1}-(1-\delta_0)\ell_i>\inf L_{i+1}+(1-\delta_0)\ell_{i+1}.
\end{align*}
Since $s_{i+1}\in L_i\cap L_{i+1}$, the claim is proved.
\ep

Let $\CD_\FF$ denote the covering distance on $I$ defined by the intervals in $\FF$.
The above lemma implies that 
\[
\CD_\FF(s_1,\phi(w_i)(s_1))=i\text{ or }i+1.\]
By the Slow Progress Lemma (Theorem~\ref{t:slp}), we see that the action $\phi$ is 
topologically conjugate neither into $\Diff_+^{k,\beta}(I)$ nor $\Diff_+^{k,\mathrm{bv}}(I)$.
What will be important for us is that every nondegenerate interval $U$ in $I$ can be arbitrarily expanded under the action of $\phi$:
\begin{lem}\label{l:expand}
For every nondegenerate interval $U\sse I$ there exists some $g\in G^\dagger$ such that 
\[\phi(g)(U)\cap L_1\cap J_3\ne\varnothing\ne \phi(g)(U)\cap L_{-1} \cap J_{-3}.\]
\end{lem}
From now on,
we will write $g(x)$ to mean $\phi(g)(x)$ whenever the meaning is clear from the context.
\bp[Proof of Lemma~\ref{l:expand}]
First note that the action $\phi$ on the interior of $I$ is sufficiently transitive; namely, 
for $x\in \Int I$ and for all $J$ in the collection $\FF$ there exists some $g\in G^\dagger$ that moves $x$ into $J$.
This immediately follows from that $\FF$ is an infinite chain and from that the restriction of the action 
$\form{a,c,d,e}$ onto $J_0$ has no global fixed points. So, we may assume that some 
$x_1\in U$ belongs to $L_1\cap J_3$ possibly after moving $U$ by an element of $G^\dagger$ 
if necessary. Let us consider four (overlaping) cases.

\emph{Case 1. Some $x_{-1}\in U$ belongs to $\bigcup_{i\le-1} L_i$.}\\
The conclusion is trivial in this case since $\inf U\le \sup L_{-1}$ and $\sup U\ge \inf L_1$.

\emph{Case 2. Some $x_{-1}\in U$ belongs to $\bigcup_{-3\le i\le -1} J_i$.}\\
For some $P,Q,R\ge0$, the element \[g=d^Rc^Qb^P\] moves $x_{-1}$ into $L_{-1}\cap J_{-3}$.
Since $g(x_1)=d^R(x_1)\ge x_1$, we have that 
\[
g(U)\cap (L_{\pm1}\cap J_{\pm3})\ne\varnothing.\]

\emph{Case 3. Some $x_{-1}\in U$ belongs to $J_0\cup J_1$.}\\
For some $P,Q\ge0$ and $w\in \form{a,c,d,e}$, the element $g_1=b^Q w b^{-P}$ moves $x_{-1}$ into $J_{-1}\cap J_{-2}$.
Then we have \[g_1(x_1)=b^Qw(x_1)\in J_2\cup J_3\cup L_1\cup L_2.\]
Using that $\form{c,d}\cong F$ acts on $J_2\cup J_3$ minimally, we can find $g_2\in \form{c,d}$ such that $g_2g_1(x_1)\in L_1\cup L_2$.
Since \[g_2g_1(x_{-1})\in \bigcup_{-3\le i\le -1} J_i,\] the proof reduces to Case 2.

\emph{Case 4. Some $x_{-1}\in U$ belongs to $J_2\cup J_3$.}\\
Using the 2--transitivity of $F$ as in Proposition~\ref{p:n-transitive},
we can find $g\in \form{c,d}$ such that $g(x_{-1})\in J_1\cap J_2$ and $g(x_1)\in J_3\cap L_1$.
The proof reduces to Case 3.

Since $x_1\in J_3\cap L_1$, the above argument exhaust all possible cases.
\ep

The above lemma verifies condition~(\ref{p:large}) in Section~\ref{s:strategy}, hence completing the proof of Theorem~\ref{t:optimal-group}.

\begin{lem}\label{l:non-commute}
For every nondegenerate interval $K\sse\supp\phi$, 
there exists some $g\in G^\dagger$
such that
\[
\CL(\phi(w_i g) K)>2i\]
whenever $i\ge1$.
\end{lem}

\bp
By Lemma~\ref{l:expand}, we can find a $g\in G^\dagger$
such that some $s_1\in L_1\cap J_3$ and $s_{-1}:=-s_1$ both belong to $\phi(g)(K)$.
Lemma~\ref{l:jumps} then implies that 
\[\CL(\phi(w_i g) K)\ge
\CL\left(\phi(w_i)[s_{-1},s_1]\right)>2i\] for all $i$.\ep

\section{Promotion to circles}\label{s:circle}
In this section, we will promote optimally expanding groups for compact intervals to those for circles.
The key idea of such a process is the following lemma.

\begin{lem}\label{l:circle-bs12}
Let $k\ge1$ be an integer and let $\beta$ be a concave modulus.
If a torsion-free simple group $H$ does not admit an embedding into $\Diff_+^{k,\beta}(I)$,
then $H\times\BS(1,2)$ does not admit an embedding into $\Diff_+^{k,\beta}(S^1)$.
\end{lem}

\bp
As usual, let us write \[\BS(1,2)=\form{a,e\mid aea^{-1}=e^2}.\]
Suppose we have a faithful representation
\[
\psi\co H\times\BS(1,2)\longrightarrow \Diff_+^{k,\beta}(S^1).\]
By Corollary~\ref{c:bs12}, there exists an integer $k\ge1$ depending on $\psi$ such that
\[
H_0:=\form{h^k\mid h\in H}\unlhd H\]
fixes all the points in $\supp e$. 
In particular, we have a restriction map
\[\psi\restriction_{H_0}\co H_0\longrightarrow \Diff_+^{k,\beta}(I).\]
On the other hand,  the group $H_0$ is nontrivial since $H$ is torsion-free.
By simplicity, we have that $H_0=H$, contradicting that $H$ does not faithfully act on $I$ by $C^{k,\beta}$ diffeomorphisms.
\ep

The following theorem generalizes the discussion so far for Theorem~\ref{t:optimal-group} to all compact connected one--manifolds.
Theorem~\ref{t:optimal-all} will then be an immediate consequence.

\begin{thm}\label{t:optimal-all2}
Let $k\in{\bZ_{>0}}$, let $\alpha\in\CM$, and let $\{\ell_i\}$ be admissible lengths.
There exists a finitely generated nonabelian group 
\[ Q=Q(k,\alpha,\{\ell_i\})\le\Diff_+^{k,\alpha}(I)\]
such that
if a concave modulus $\beta\succ_k0$ satisfies
\[\beta(1/i) \cdot(1/i)^k\preccurlyeq \alpha(1/i)\cdot \ell_i^k\]
for almost all $i\ge1$,
and then every homomorphism
\[
[Q,Q]\longrightarrow\Diff_+^{k,\beta'}(M^1)\]
is trivial for $\beta'\in\{\beta,\mathrm{bv}\}$ and for $M^1\in\{I,S^1\}$.
Furthermore, $Q$ can be chosen so that $[Q,Q]$ is simple and every proper quotient of $Q$ is abelian.
\end{thm}

\bp
Given the parameter $(k,\alpha,\{\ell_i\})$, we apply Theorem~\ref{t:optimal-group}
to define
\[\phi_1:=\phi(k,\alpha,\{\ell_i\})\co G^\dagger\longrightarrow\Diff_c^{k,\alpha}(\bR).\]
We will fix a concave modulus $\beta\succ_k0$ satisfies
\[\beta(1/i) \cdot(1/i)^k\preccurlyeq \alpha(1/i)\cdot \ell_i^k\]
for almost all $i\ge1$. We let $\beta'\in\{\beta,\mathrm{bv}\}$.

Let $T_1$ be the image of $\phi_1$, which is a five--generated group.
We will successively upgrade $T_1$ through a sequence of finitely generated groups
such that each group in the sequence enjoys stronger algebraic properties than the previous.
The Chain Group Trick~\ref{cor:chain-gp-trick} will be used twice for this purpose,
but it requires that the group being promoted to have torsion-free abelianization.
So, we will apply the Rank Trick (Lemma~\ref{lem:rank-trick})  twice 
to guarantee that the group's abelianization is torsion-free.

We now spell out this sketched argument a bit more. 
Using the Rank Trick to $\phi_1$,
we can promote $\phi_1$ to another representation
\[
\phi_2\co G^\dagger \longrightarrow \form{T_1,\Diff_c^\infty(\bR)}\le\Diff_c^{k,\alpha}(\bR)\]
 such that $\phi_1(g)=\phi_2(g)$ for $g\in [G^\dagger,G^\dagger]$,
and such that the abelianization of $T_2:=\phi_2(G^\dagger)$ is torsion-free.
We may assume ${\suppc T_2}\sse(0,1)$.
By the Chain Group Trick,
we see that 
\[ 
T_3:=\form{T_2,\rho_{\mathrm{GS}}(F)}\]
is a seven--chain group acting minimally on $(0,1)$.
Furthermore, there exist an embedding
\[u\co T_2\into [T_3,T_3].\]

\begin{claim}
Every homomorphism $[T_3,T_3]\longrightarrow\Diff_+^{k,\beta'}(I)$ is trivial.
\end{claim}

Suppose we have a homomorphism
\[\rho\co[T_3,T_3]\longrightarrow \Diff_+^{k,\beta'}(I).\] 
Since $[T_3,T_3]$ is simple, it suffices to show that $\rho$ is not injective.
Comparing the map $\phi_1$ with the composition
\[\begin{tikzcd}
G^\dagger\arrow{r}{\phi_2} &
T_2\arrow[hookrightarrow]{r}{u} &
{[T_3,T_3]}\arrow{r}{\rho}& \Diff_+^{k,\beta'}(I),
\end{tikzcd}
\]
we see from Theorem~\ref{t:optimal-group} that some $g\in[G^\dagger,G^\dagger]$ satisfies $\rho\circ u\circ  \phi_2(g)=1$ but $\phi_1(g)\ne1$.
Since $u\circ\phi_2(g)=u\circ\phi_1(g)\ne1$, the claim is proved.

The group $T_3$ may not have a torsion-free abelianization.
So, we will apply the Rank Trick to the quotient $F_7\onto T_3$
in order to obtain a representation
\[
\phi_4\co F_7\longrightarrow \Diff_c^{k,\alpha}(\bR)\]
such that the abelianization of $\phi_4(F_7)$ is torsion-free.
Setting $T_4:=\phi_4(F_7)$, we can further require that
\[
[T_3,T_3]=[T_4,T_4].\]
We define the nine--generated group
\[
T_5:=T_4\times \BS(1,2).\]

\begin{claim}
There does not exist a faithful representation
\[
T_5\longrightarrow \Diff_+^{k,\beta'}(S^1).\]
\end{claim}

Suppose we have a homomorphism
\[
\rho\co T_5\longrightarrow \Diff_+^{k,\beta'}(S^1).\]
Applying the preceding claim
and Lemma~\ref{l:circle-bs12} to 
\[
H:=[T_3,T_3]=[T_4,T_4],\]
we see that $H\times \BS(1,2)$ cannot embed into $\Diff_+^{k,\beta'}(I)$.
This is a contradiction, and the claim is proved.

Since $T_5$ has a torsion-free abelianization, we can apply the Chain Group Trick 
and embed $T_5$ into the commutator subgroup of some eleven--chain group
\[
R=R(k,\alpha,\beta)\le\Diff_c^{k,\alpha}(\bR)\]
acting minimally on its support.
By Theorem~\ref{thm:chain-dichotomy} and by the last claim above, we complete the proof of Theorem~\ref{t:optimal-all2}.
\ep

\bp[Proof of Theorem~\ref{t:optimal-all}]
Assuming the hypotheses of Theorem~\ref{t:optimal-all}, we set
\[
\ell_i:=  \frac1i \cdot \left(\frac{\beta(1/i) }{\alpha(1/i)}\right)^{1/k}.\]

\begin{claim}
The lengths $\{\ell_i\}$ are admissible.\end{claim}
By the concavity of $\alpha$, we see that
\[
1\le \alpha(y)\le \frac{y}{x}\cdot \alpha(x)\]
for all $0<x\le y$. 
In particular, we have for $i\ge1$ that
\[
1\le
\frac{\alpha(1/i)}{\alpha(1/(i+1))}\le\frac{i+1}{i}\le 2.\]
Since $\beta(x)$ satisfies the same bounds, we see that
$\{ \ell_i/\ell_{i+1}\}$ is bounded, and bounded away from zero.
Moreover, whenever $i\le y\le i+1$ we have that
\[
\frac1y\left(\frac{\beta(1/y)}{\alpha(1/y)}\right)^{1/k}
\ge
\frac1{i+1}\left(\frac{\beta(1/(i+1))}{\alpha(1/i)}\right)^{1/k}
\ge
\frac1{2^{1+1/k}}\cdot \frac1i \left(\frac{\beta(1/i)}{\alpha(1/i)}\right)^{1/k}.
\]
Since the leftmost term is integrable on $[1,\infty)$, we have that \[\sum_i\ell_i<\infty.\]
The claim is thus proved.

By the above claim, we can apply Theorem~\ref{t:optimal-all2}
and see that the group
\[R(k,\alpha,\beta):=Q(k,\alpha,\{\ell_i\})\]
satisfies the conclusion of the theorem.\ep

\begin{rem}\label{rem:g1}
Theorem~\ref{t:optimal-all} appears weaker than Theorem~\ref{t:optimal-all2}.
There is a subtle point that highlights this difference when $k=1$.
Consider the concave modulus
\[
\alpha(x)=\frac1{\log(1/x)}.\]
Applying Theorem~\ref{t:optimal-all2}, we can find some group
\[
G_1:=Q(1,\alpha,\{1/(i\log^2 i)\})\le\Diff_+^{1,\alpha}(I)\le\Diff_+^1(I),\]
such that
\[
[G_1,G_1]\not\into \Diff_+^{1,s}(M^1)\]
for all real number $0<s<1$. This is because $\beta_s(x)=x^s$ is sub-tame,
and because
\[
x^s \log^2(1/x) \preccurlyeq 1/\log(1/x)\]
for all $x$ sufficiently near from 0. 
The existence of such a group $G_1$ does not immediately follow from Theorem~\ref{t:optimal-all}. Indeed,
an obvious choice of a concave modulus $\beta(x)$
satisfying
\[\int_0^1 \frac1{x}\left(\frac{\beta(x)}{\alpha(x)}\right)<\infty\]
and $\beta(x)\succcurlyeq x^s$ for all $s\in(0,1)$
would be 
\[\beta(x):= \alpha(x)/\log^2(1/x).\]
However, this modulus $\beta$ is not sub-tame and one cannot apply this theorem for the pair $(\alpha,\beta)$.
\end{rem}

\section{Continua of groups of prescribed critical H\"older regularity}\label{s:hoelder}
Recall from Section~\ref{ss:defn-rem} that the critical regularity of a group $G$ with respect to a manifold $M$ is defined as
\[\CR_M(G):=\sup\{r\ge1\mid G\into\Diff_c^r(M)_0\}\in\{-\infty\}\cup[1,\infty].\]
In this section,
we  prove that for each real number $r\ge1$ and for every compact connected one--manifold $M$
there exist continuum--many distinct isomorphism classes of finitely generated groups $G$ such that 
\[
\CR_M(G)=r.\] Here, \emph{continuum}\index{continuum} refers to the cardinality of $\bR$.
We will actually prove a much stronger result below, which obviously implies Corollary~\ref{c:optimal-all} as well.

\begin{thm}\label{t:continuum}
For each real number $r\ge1$, there exist continua $X_r, Y_r$ of finitely generated groups such that the following hold.
\be[(i)]
\item
Each group in $X_r\cup Y_r$  is a subgroup of $\Homeo_+[0,1]$,
which has a simple commutator group
and  every proper quotient of which is abelian.
\item No two groups in $X_r\cup Y_r$ have isomorphic commutator subgroups.
\item 
For each group $A\in X_r$ and for each $M^1\in\{I,S^1\}$, we have that $A\le\Diff_+^r(I)$ and that
\[
[A,A]\not\into \bigcup_{s>r} \Diff_+^s(M^1).\]
\item
For each group $B\in Y_r$ and for $M^1\in\{I,S^1\}$, we have that \[B\le\bigcap_{r<s} \Diff_+^r(I)\] and that
\[
[B,B]\not\into \Diff_+^r(M^1).\]
\ee
\end{thm}

In order to construct such continua, we will apply Theorem~\ref{t:optimal-all} for suitable pairs of concave moduli 
in the following form:
\[\omega_{s,t,u}(x):=\frac{x^s}{\log^u(1/x)} \exp\left(- t \sqrt{\log(1/x)}\right)
.\]
Roughly speaking, 
$\omega_{s,t,u}$ is a small perturbation of the H\"older modulus $x^s$, as we can write
\[\omega_{s,t,u}(x)=
{\exp\left(-s\log(1/x) - t \sqrt{\log(1/x)}-u\log\log(1/x)\right)}.\]

\begin{lem}\label{l:omega-st}
The following hold for $u\ge 0$.
\be[(1)]
\item
The map $\omega_{s,t,u}$ is a concave modulus (defined near $x=0$) in the following cases:
\be[(i)]
\item $s=0$ and $t>0$.
\item $0<s<1$ and $t\in\bR$; in this case, $\omega_{s,t,u}$ is sub-tame.
\item $s=1$ and $t\le0$; in this case, $\omega_{s,t,u}$ is sub-tame.
\ee
\item
Let $(s,t)$ and $(S,T)$ be in the range described part (1).
If $(s,t)<(S,T)$ in the lexicographical order, then we have that
\[\lim_{x\to+0}\omega_{S,T,0}(x)/\omega_{s,t,u}(x)=0.\]
\ee
\end{lem}
 \bp
(1)
The precise verification of the concavity is a bit tedious, but one can first have an intuition by noting that the behavior of $\omega_{s,t,u}$ for $s>0$ is governed by $x^s$ near $x=0$.
Rigorously, we compute (by Wolfram Mathematica 12) that
\[
\omega'_{s,t,u}(x)=
\frac{1}{2} x^{s-1} e^{-t \sqrt{\log \left(\frac{1}{x}\right)}} \log ^{-u-1}\left(\frac{1}{x}\right) \left(2 s \log \left(\frac{1}{x}\right)+t \sqrt{\log \left(\frac{1}{x}\right)}+2 u\right).\]
If $s>0$ or $s=0$ and $t>0$, then we see that $\omega_{s,t,u}'>0$ near $x=0$ from that
\[ \sqrt{\log(1/x)}=o(\log (1/x)).\]
Similarly, to determine the sign of $\omega_{s,t,u}''$, we first write
\[
\omega''_{s,t,u}(x)=
\frac{1}{4} x^{s-2} e^{-t \sqrt{\log \left(\frac{1}{x}\right)}} \log ^{-u-2}\left(\frac{1}{x}\right)F(x),\]
where
\begin{align*}
F(x)=&4 (s-1) s \log ^2\left(\frac{1}{x}\right)
+2 (2 s-1) t \log ^{\frac{3}{2}}\left(\frac{1}{x}\right)+
\left(8 s u+t^2-4 u\right)\log \left(\frac{1}{x}\right) 
\\
&+(4 t u+t) \sqrt{\log \left(\frac{1}{x}\right)}+4 u (u+1).\end{align*}
Hence, if $0<s<1$ then $F(x)<0$ and $\omega''$ is concave.
The same is true when $s=0$ and $t>0$, or $s=1$ and $t<0$.

To verify that $\omega_{s,t,u}$ is sub-tame in the cases (ii) and (iii), 
let us write
\[
\frac{\omega_{s,t,u}(vx)}{\omega_{s,t,u}(x)}=v^s \exp\left(-t\left(\sqrt{\log(1/(vx))}-
\sqrt{\log(1/x)}\right)\right)\cdot \left(\frac{\log(1/x)}{\log(1/v)+\log(1/x)}\right)^u.\]
Here, we assume $v>0$ is small enough. 
If $t\ge 0$ it is obvious that 
\[
\frac{\omega_{s,t,u}(vx)}{\omega_{s,t,u}(x)}\le v^s.\]
If $t<0$ we see
\[\frac{\omega_{s,t,u}(vx)}{\omega_{s,t,u}(x)}\le v^s\exp(-t \sqrt{\log(1/v)}).\]
In either case, we see that $\omega_{s,t,u}(x)$ is sub-tame.

Part (2) is straightforward after noticing that 
\[
\log(1/x)\gg \sqrt{\log(1/x)}\gg \log\log(1/x).\]\ep

\bp[Proof of Theorem~\ref{t:continuum}]
Let us write $r = k+s$ for $k=\floor{r}\ge1$.

\emph{Case $0<s<1$.}

By Lemma~\ref{l:omega-st} and Theorem~\ref{t:optimal-all}, we can define collections of finitely generated groups as below.
\begin{align*}
X_r&:=\{ R(k,\omega_{s,t,0},\omega_{s,t,k+1}) \mid t>0)\},\\
Y_r&:=\{ R(k,\omega_{s,t,0},\omega_{s,t,k+1}) \mid t<0\}.
\end{align*}
For $t>0$, we have
\[
\bigcup_{ S\in(s,1]}\Diff_+^{k,S}(M^1)\le
\Diff_+^{k,\omega_{s,t,k+1}}(M^1)\le
\Diff_+^{k,\omega_{s,t,0}}(M^1)\le \Diff_+^{k,s}(M^1).\]
If $t<0$, then
\[
 \Diff_+^{k,s}(M^1)
 \le
 \Diff_+^{k,\omega_{s,t,k+1}}(M^1)\le
\Diff_+^{k,\omega_{s,t,0}}(M^1)\le
\bigcap_{ S\in(0,s)}\Diff_+^{k,S}(M^1).\]
Thus, parts (i), (iii), and (iv) of the theorem are readily implied by Theorem~\ref{t:optimal-all}.

For real numbers $t<T$, we saw in Lemma~\ref{l:omega-st} that
\[\omega_{s,T,0}(x)\le\omega_{s,t,k+1}(x)\]
for all $x>0$ near $0$. In particular, we have
\[
\Diff_+^{k,\omega_{s,T,0}}(M^1)\le \Diff_+^{k,\omega_{s,t,k+1}}(M^1).\]
This implies that
\[
[R(k,\omega_{s,t,0},\omega_{s,t,k+1}),R(k,\omega_{s,t,0},\omega_{s,t,k+1})]\not\into
R(k,\omega_{s,T,0},\omega_{s,T,k+1}).\]
This proves part (ii).

\emph{Case $s=0$ and $k\ge2$.}

We have $r=k\ge2$. Similarly to the previous case, we let
\begin{align*}
X_k&:=\{ R(k,\omega_{0,t,0},\omega_{0,t,k+1}) \mid t>0)\},\\
Y_k&:=\{ R(k-1,\omega_{1,t,0},\omega_{1,t,k+1}) \mid t<0\}.
\end{align*}
These are well-defined since $\omega_{0,t,k+1}\succ_k 0$
and
$\omega_{1,t,k+1}\succ_{k-1}0$.
The conclusion is again obvious.
Since $\omega_{1,t,k+1}(x)\ge x$ for $t<0$ and for small $x>0$, we have an even stronger fact:
\[
[R(k-1,\omega_{1,t,0},\omega_{1,t,k+1}),R(k-1,\omega_{1,t,0},\omega_{1,t,k+1})]\not\into \Diff_+^{k-1,\mathrm{Lip}}(M^1).\]

\emph{Case $s=0$ and $k=1$.}

Recall in Remark~\ref{rem:g1} we have found
\[ 
G_1\le\Diff_+^1(I)\]
such that $[G_1,G_1]\not\into\Diff_+^{1,V}(M^1)$ for all $V>0$.
For each $V>1$, let us fix $G_V\in X_V$. 

Let us consider $T\ge1$. We  apply the Rank Trick for the natural surjection $F_m\onto G_T$ with some $m\ge2$,
and find a finitely generated group $\tilde G_T\le\Diff_+^1(I)$ such that 
\[
[G_T,G_T]=[\tilde G_T,\tilde G_T],\]
and such that the abelianization of $\tilde G_T$ is free abelian.
Since $\tilde G_T$ is finitely generated, we have that
 \[\tilde G_T\not\into \bigcup_{S>T}\Diff_+^{S}(M^1).\]

Let $V>1$.
We apply the Chain Group Trick to the product $\tilde G_1\times\tilde G_V$ 
and obtain a minimal chain group $\Gamma(V)\le\Diff_+^1(I)$ satisfying that
\[
\tilde G_1\times\tilde G_V\into [\Gamma(V),\Gamma(V)].\]
The group $[\Gamma(V),\Gamma(V)]$ does not embed into \[\bigcup_{S>1}\Diff_+^S(M^1),\]
since neither does $\tilde G_1$.

If there were only countably many isomorphism classes in the collection
\[\{[\Gamma(S),\Gamma(S)]\mid S\ge1\},\]
then the same is true for
the collection of finitely generated subgroups in that collection.
This would contradict that
 $\CR_I(\tilde G_S)=S$ for all $S\ge1$. 
 It follows that there exists a continuum $X^*\sse(1,\infty)$ such that 
whenever $S\ne T$ belong to $X^*$ we have
\[
[\Gamma(S),\Gamma(S)]\not\cong[\Gamma(T),\Gamma(T)].\]
This completes the construction of $X_1$.

We now construct the collection $Y_1$.
Note that every countable subgroups of $\Homeo_+(M^1)$ is topologically conjugate into 
$\Diff_+^{\mathrm{Lip}}(M^1)$ by~\cite{DKN2007}; see Theorem~\ref{thm:lip-conj} above.
So, the condition (iv) of the current theorem is now equivalent to that $B\le\Homeo_+(M^1)$ and 
\[
[B,B]\not\into\Diff^1_+(M^1)\]
for each $B$ in the collection $Y_1$ that we are about to build.
We start with an observation of Kropholler and Thurston  
that the following is a nontrivial finitely generated perfect orderable group:
\[\tilde\Delta=\form{a,b,c,t\mid a^2=b^3=c^7=abc=t}\le\widetilde{\PSL}(2,\bR)\le\Homeo_+(\bR).\]
We direct the reader to Appendix~\ref{ch:append3} and to~\cite{Thurston1974Top,Bergman1991PJM} for a further discussion of this
example and its perfectness.
Note that $\tilde\Delta$ is the fundamental group of a Seifert fibered 3--manifold, which fibers over the Fuchsian $(2,3,7)$--triangle group 
\[
\Delta=\form{a,b,c\mid a^2=b^3=c^7=abc=1}\le\PSL(2,\bR)\le\Homeo_+(S^1)\]
By Thurston Stability (Theorem~\ref{thm:thurston-stab}), the group $\tilde\Delta$ does not embed into $\Diff_+^1(I)$.
We feed $T_2:=\tilde\Delta$ to our proof of Theorem~\ref{t:optimal-all}.
We proceed to obtain $T_3,T_4,T_5$ and $R$ as in Section~\ref{s:circle},
and see that $R$ is a minimal chain group on $(0,1)$
and $[R,R]\not\into\Diff_+^1(S^1)$.
Seting $H_1:=R$, follow the construction of $X_1$.
Namely, we pick some $H_S\in Y_S$ for each $S>1$, and build finitely generated groups $\tilde H_1,\tilde H_S$ and $\Lambda(S)$ such that
\[
\tilde H_1\times\tilde H_S\into[\Lambda(S),\Lambda(S)]\not\into\Diff_+^1(M^1).\]
In particular, $\Lambda(S)$ will be chosen to be a minimal chain group.
By the same reasoning, we can find a continuum $Y^*\sse(0,1)$ such that 
for every distinct pair $S,T$ in $Y^*$ we have
\[
[\Lambda(S),\Lambda(S)]\not\cong [\Lambda(T),\Lambda(T)].\]
This completes the construction of $Y_1$.\ep
\begin{rem}
Calegari~\cite{Calegari2006TAMS} provided an example a group $Q\le\Homeo_+(S^1)$
that admits no embedding into $\Diff_+^1(S^1)$.
He extended the embedding
\[
\tilde\Delta\into \Homeo_+(\bR)\into\Homeo_+(S^1)\]
to a homomorphism
\[
\rho\co \form{\tilde\Delta,x,y,z\mid xax^{-1}=a^2,yby^{-1}=b^2,zcz^{-1}=c^2}\longrightarrow\Homeo_+(S^1),\]
and proved that the image of $\rho$ is such an example.
Theorem~\ref{t:continuum} generalizes his result to provide continuum--many such examples.
\end{rem}


%
%
%

\chapter{Applications}\label{ch:app}
\begin{abstract}In this chapter, we will outline some of the main consequences of the machinery developed in this book, following mostly
applications that are discussed in~\cite{KK2020crit,BKK2019JEMS,KK2018JT}. We will begin with foliation theory, which
was the driving force behind the study of regularity of group actions on manifolds for much of the latter's history, and then move on
to the problem of computing the critical regularity of various classes of groups.
This chapter will be of a different flavor than most of the rest of the book, serving more as a survey in a more informal tone.
As such, less detailed proofs and a more
informal exposition will be given.\end{abstract}

One of the major consequences of
the construction of groups with specified critical regularity is that they yield strong unsmoothability
results for certain foliations on $3$--manifolds, coming from foliations on trivial $I$--bundles over surfaces whose monodromy groups
are groups with specified critical regularity (see Corollary~\ref{cor:tsuboi-ref} and Corollary~\ref{cor:3-mfld}).

We will also describe some results that illustrate the interplay between algebraic structure of groups and possible regularities of group
actions one one--manifolds, concentrating on right-angled Artin groups and mapping class groups of surfaces. The main results in that
discussion will by Theorem~\ref{thm:kharlamov} and Theorem~\ref{thm:mcg-c2}.

Throughout this chapter, we will adopt the following standing convention. If $ r\geq 1$ is a real number of the form $ r=k+\tau$, where
$k$ is an integer and $\tau\in[0,1)$, then we will write $C^{ r}$ and $\Diff^{ r}$ for $C^{k,\tau}$
and $\Diff^{k,\tau}$, respectively. This will be notationally convenient for the
discussion of constructions that are uniform in $ r$, so that we do not need to decompose each real number as a sum of its floor and
fractional part.

\section{Foliation theory}\label{sec:foliation}

One of the original sources of interest in regularity of group actions comes from the theory of foliations on manifolds.
In fact, the Plante--Thurston Theorem (Theorem~\ref{thm:plante-thurston}) and the Thurston Stability Theorem
(Theorem~\ref{thm:thurston-stab}) were both originally formulated as results about foliations, and their content as facts about regularity
is mostly abstracted after the fact. In this section, we
shall go in the other direction, by discussing the applicability of critical regularity to the non--smoothability of
certain codimension one foliations on manifolds. 

\subsection{Foliations and suspensions of group actions}

A foliation $\FF$ on an $n$--manifold $M$ is a local decomposition of $M$ modeled on the direct product decomposition
$\bR^n\cong\bR^p\times\bR^q$, where $p,q\geq 1$ and $p+q=n$. Precisely, let $\{x_1,\ldots,x_n\}$ be coordinates for $\bR^n$;
the manifold $M$ is equipped with a
$q$--dimensional \emph{foliation
atlas}\index{foliation atlas},
which is a collection of charts \[\varphi_i\colon U_i\longrightarrow \bR^n\] for some suitable cover $\UU=\{U_i\}_{i\in I}$ of $M$,
so that the transition functions $\varphi_{ij}=\varphi_i\circ\varphi_j^{-1}$ decompose in coordinates as
\[\varphi_{ij}(x_1,\ldots,x_n)=(\varphi_{ij}^1(x_1,\ldots,x_n), \varphi_{ij}^2(x_{q+1},\ldots,x_n)),\] where 
\[\varphi^1_{ij}\colon\bR^n\longrightarrow
\bR^q,\quad \varphi^2_{ij}\colon\bR^p\longrightarrow\bR^p.\]

The intuitive meaning behind this decomposition is that a foliation is locally modeled by a submersion $\bR^n\longrightarrow\bR^q$. More
generally, a submersion from an $n$--dimensional manifold $M$ to a $q$--dimensional manifold induces a foliated structure on $M$ via
the Implicit Function Theorem. Transition maps between foliation charts preserve leaves of the foliation, which is why $\varphi^1_{ij}$
is allowed to depend on all variables but $\varphi^2_{ij}$ may not.

Thus locally, $M$ is an $\bR^p$--worth of $q$--dimensional
immersed manifolds. The natural number $p$ is the \emph{codimension}\index{codimension of a foliation} of the foliation, 
and $q$ is the \emph{dimension}\index{dimension of a foliation} of the foliation.
The regularity of the functions \[\{\varphi_{ij}^k\}_{i,j\in I, k\in\{1,2\}},\] which is always assumed to be at least $C^0$,
is the \emph{regularity}\index{regularity of a foliation} of the foliation. The regularities of $\varphi_{ij}^1$ and $\varphi_{ij}^2$
need not coincide, and so the regularity of the foliation is generally an ordered pair. When one regularity is finite and the other is infinite,
then the regularity is simply the finite member of the pair, and this will usually be the case for us; specifically, unless otherwise noted, we
will always assume that $\varphi_{ij}^1$ is $C^{\infty}$, and so the regularity of the foliation is simply the regularity of the functions
$\varphi_{ij}^2$ occurring in the atlas.

The regularity of the foliation
may be much lower than the regularity of $M$. Indeed, $(M,\FF)$ may be such that $M$ itself may admit a $C^{\infty}$ atlas, but the
foliation atlas describing $\FF$ is merely required to be $C^0$. If $ r\geq 1$, a \emph{$C^{ r}$--smoothing}\index{smoothing of a foliation} 
of $(M,\FF)$ is a homeomorphism
\[f\colon (M,\FF)\longrightarrow (M',\FF'),\] where $\FF'$ is a $C^{ r}$ foliation.

Beyond the basic definitions, we will not survey much of the vast theory of foliations, and instead we 
will content ourselves with directing the reader
to some standard references such as~\cite{CandelConlonI,CandelConlonII}.

We will use one fundamental construction in the theory of foliations, namely that of the \emph{suspension of a group action}\index{suspension
of a group action}. Let $B$ and $M$
be fixed connected, smooth manifolds, and let $\yt B\longrightarrow B$ be the universal covering map. Fix a representation
\[\psi\colon\pi_1(B)\longrightarrow\Diff_+^{ r}(M).\] Notice that $\yt B\times M$ is naturally a manifold of dimension \[n=\dim B+\dim M,\]
and admits a foliation of dimension $\dim B$ and codimension $\dim M$.
 We also have a natural action of $\pi_1(B)$ on
$\yt B\times M$ via the diagonal action, where $g\in\pi_1(B)$ acts on $\yt B$ by a deck transformation and on $M$ via $\psi(g)$.
Observe that since the action of $\pi_1(B)$ on $\yt B$ is free and properly discontinuous, the induced action on $\yt B\times M$ is as
well, and so \[\yt B\times M\longrightarrow E(\psi)=(\yt B\times M)/\pi_1(B)\] is a covering map. The space $E(\psi)$ is called the
\emph{suspension} of $\psi$. It has the structure of a $C^{ r}$ fiber bundle, which is in turn a foliation 
$\FF$; this is called a \emph{$C^r$ foliated bundle}\index{foliated bundle}.
If $g\in\pi_1(B)$ then $g$ can be represented by a loop in
$B$.
 If $b_0\in \yt B$ is a basepoint, then the fiber over $b_0$ is identified with the fiber over $g.b_0$ by the map $\psi(g)$. The map $\psi$ is
called the \emph{monodromy representation}\index{monodromy representation} of the bundle. 
Suspensions of group actions give rise to rich families of examples.
Consider, for instance, an action of the fundamental group $\pi_1(S)$ of a surface $S$ on the interval $I$, and suppose that there
is a point $x\in I$ whose orbit is free; that is, the map $\pi_1(S)\longrightarrow I$ given by $g\mapsto g(x)$ is injective. Then
the suspension of the action of $\pi_1(S)$ will be homeomorphic to $S\times I$, though the induced foliation will not be the product foliation:
instead, it will admit a contractible leaf that is identified with the universal cover $\yt{S}$ of $S$. The reader is directed to~\cite{SoutoMarq18}
for some interesting issues surrounding smooth foliations of $S\times I$.

Suspensions of group actions can be understood systematically, as they
are classified up to homeomorphism by the conjugacy class of the monodromy representation:

\begin{prop}[See Theorem 3.1.5 of~\cite{CandelConlonI}]\label{prop:suspension-unique}
Let $B$ be a smooth manifold, let $M$ be a manifold of regularity at least $C^r$,
and let \[\psi,\psi'\colon\pi_1(B)\longrightarrow\Diff_+^{ r}(M).\] Then $E(\psi)$ is homeomorphic
as a bundle to $E(\psi')$ if and only if $\psi$ and
$\psi'$ are conjugate by a homeomorphism of $M$. The bundles $E(\psi)$ and $E(\psi')$ are $C^{r}$--diffeomorphic for
$r$ if and only if $\psi$ and $\psi'$ are conjugate by a $C^{r}$ diffeomorphism.
\end{prop}

In Proposition~\ref{prop:suspension-unique}, a bundle homeomorphism $E(\psi)\longrightarrow E(\psi')$ is a homeomorphism
of the total spaces that takes fibers to fibers. We leave the proof of the proposition as an exercise for the reader.


%
%
%

\subsection{Non-smoothability of codimension one foliations}\label{sss:tsuboi}

Proposition~\ref{prop:suspension-unique} gives a tool for building foliations that are $C^0$ but not $C^1$, or in general
$C^r$ for some $ r\geq 1$ but not $C^{s }$ whenever $s > r$. The first examples of group actions
on $M=I$ that are $C^k$ for some $k\in{\bZ_{>0}}$ but not $C^{k+1}$ were constructed by  Tsuboi~\cite{Tsuboi1987},
and
 independently
by Cantwell and Conlon~\cite{CC1988}.

We sketch Tsuboi's construction here, since it is simple and easy to understand, and since it makes good use of the Muller--Tsuboi
trick (Lemma~\ref{l:muller-tsuboi}). 

Consider the maps 
\[
\tau_t(x)=x+t,\quad \lambda_t(x)=e^tx\]
for $t\in\bR$. 
These generate the orientation preserving affine transformation group $\Aff_+(\bR)$. 
Recall from Theorem~\ref{thm:bs-cpt} that we have an injection
\[
\Aff_+(\bR)\longrightarrow\Diff_{[0,1]}^\infty(\bR),\]
defined using a certain topological conjugacy $\psi\co (0,1)\longrightarrow\bR$.
More precisely, the affine map $x\mapsto ax+b$ gets sent first to \[x\mapsto \frac{x}{a-bx},\] and then
this map is made infinitely tangent to the identity at the endpoints by the Muller--Tsuboi trick.
The map $\psi$ is chosen to be smooth in the open interval $(0,1)$.
For the rest of this subsection,
we will denote by $A_0\le \Diff_{[0,1]}^\infty(\bR)$ the image of such an embedding.

The vector fields  $\partial/\partial x$ and $x\partial/\partial x$ generate the flows $\tau_t(x)$ and 
$\lambda_t(x)$. We write $\tau^*$ and $\lambda^*$ for the pull-backs of these two vector fields by $\psi$.
Theorem~\ref{thm:bs-cpt} then amounts to showing that 
\begin{align*}
T(t,x)&:=\Phi_{\tau^*}(t,x)=\psi^{-1}\circ\tau_t\circ\psi(x),\\
L(t,x)&:=\Phi_{\lambda^*}(t,x)=\psi^{-1}\circ\lambda_t\circ\psi(x)\end{align*}
belong to $A_0\le \Diff_{[0,1]}^\infty(\bR)$ for each $t$.
Here, $\Phi_\rho(t,x)$ denotes the flow of an arbitrary vector field $\rho\co\bR\longrightarrow\bR$.

More rigorously, the maps $T$ and $L$ are first defined on $\bR\times (0,1)$. 
Theorem~\ref{thm:bs-cpt} then implies that for each $t\in \bR $ the maps $x\mapsto T(t,x)$ and 
$x\mapsto L(t,x)$ extend to smooth maps on $x\in \bR$, by the identity outside of $(0,1)$. It is also 
clear from the definition that the extended maps are simultaneously continuous on $(t,x)$.
Hence, we can apply the Bochner--Montgomery Theorem (Theorem~\ref{thm:BM-simple})
to deduce that $T$ and $L$ are smooth on $\bR\times\bR$.

In order to describe Tsuboi's example, let us now fix an integer number $r\geq 1$, and $\eps>0$; we further require that \[\eps(r-1)<1\] in the case when $r\ge2$.

We select a strictly increasing sequence of points $\{y_n\}_{n\in\bN}\sse [0,1]$ such that $y_0=0$, which converges to $1$, and which satisfies \[y_n-y_{n-1}=n^{-1-\eps}\] for $n\gg 0$. Fix an affine homeomorphism
\[\omega_n\co [0,1]\longrightarrow J_n:=[y_{n-1},y_n].\]
We also fix a positive sequence $\{u_n\}$ such that $u_n\in(0,1/n^r)$.

We define $\bar T, \bar L\in\Homeo_+[0,1]$ by requiring that $\bar T(y_n)=\bar L(y_n)=y_n$
and that
\begin{align*}
\bar T\restriction_{J_n}&:=
\omega_n\circ  T(u_n,  \underbar{\phantom{x}})\circ \omega_n^{-1},\\
\bar L\restriction_{J_n}&:=
\omega_n\circ L\left(\frac{\log 2}{n^r}, \underbar{\phantom{x}}\right)\circ  \omega_n^{-1},
\end{align*}
for each $n$.
It is clear that by setting $\bar T(1)= \bar L(1)=1$, we obtain a homeomorphism of $[0,1]$.

Since $A_0\le \Diff_{[0,1]}^{\infty}(I)$, the homeomorphisms $\bar T$ and $\bar L$ are $C^{\infty}$ at $y_n$ for all $n$, and
so lie in $\Diff_+^{\infty}(-\infty,1)$. The particular choices involving the number $r$ will guarantee that the resulting homeomorphisms $\bar T$ and $\bar L$  are in fact $C^r$ at the point $1$; this follows from a computation that we now spell out.

For each smooth vector field  $\rho\co \bR\longrightarrow\bR$  supported in $[0,1]$
and for each $i\ge1$, we set
\[
P_{i,\rho}:=\sup\left\{ \abs*{
\frac\partial{\partial t}
\left(
\frac{\partial^i\Phi_\rho(t,x)}{\partial x^i}
\right)
}\co
{-1\le t\le 1, x\in\bR}\right\}.
\]
We let $P_i:=\max(P_{i,\tau^*},P_{i,\lambda^*})$.
For each $x\in [0,1]$, we have
\[
\abs*{\frac{\partial T}{\partial x}(t,x)-1}
=\abs*{\frac{\partial T}{\partial x}(t,x)-\frac{\partial T}{\partial x}(0,x)}
\le P_1\cdot  t.\]
We similarly obtain the following estimates for all $i\ge1$:
\begin{align*}
\abs*{\frac{\partial^i T}{\partial x^i}(t,x)-\delta_{i1}}
&\le P_i\cdot t,\\
\abs*{\frac{\partial^i L}{\partial x^i}(t,x)-\delta_{i1}}
&\le P_i\cdot t,\end{align*} where $\delta_{i1}$ denotes the Kronecker delta.

We now have the following estimates for each $n,i\ge1$:
\begin{align*}
\sup_{J_n} \abs*{
\bar T^{(i)}-\delta_{i1}
}&\le P_i\cdot  \frac{u_n}{|J_n|^{i-1}}\le P_i\cdot n^{(i-1)(1+\epsilon)-r},\\
\sup_{J_n} \abs*{
\bar L^{(i)}-\delta_{i1}
}&\le P_i\cdot  \frac{\log 2}{n^r} \cdot \frac{1}{|J_n|^{i-1}}\le \log 2 \cdot P_i\cdot n^{(i-1)(1+\epsilon)-r}.
\end{align*}
It follows that for each $1\le i\le r$, the $C^i$--distances from the identity to $\bar T$ and  to $\bar L$ on $J_n$ tend to zero as $n\to\infty$. 
Thus, we  conclude that \[\bar T, \bar L\in \Diff_{[0,1]}^r(\bR).\]

The restriction of the group $\form{\bar T,\bar L}$  to each interval $J_n$ is simply an action of an affine group. 
As \[\lambda_s^i \tau_t\lambda_s^{-i}=\tau_{e^{si} t}\] for all $s,t\in\bR$ and $i\in\bZ$,
we have that
\[
[\bar T,\bar L^i\bar T\bar L^{-i}]=1.\]
Recall from Section~\ref{sec:abt} that the \emph{lamplighter group}\index{lamplighter group} is a metabelian group defined by the
presentation
\[\bZ\wr\bZ:=\form{ \tau,\lambda\mid \left[\tau,\lambda^i\tau\lambda^{-i}\right]=1\text{ for all }i\in\bZ}.\]
We have just constructed a representation
 \[\psi_r\colon\bZ\wr\bZ\longrightarrow\Diff_+^r(I),\]  defined by $\tau\mapsto \bar T$ and
$\lambda\mapsto \bar L$. 

We note the following non-$C^k$ criterion, which is due to Tsuboi. 
\begin{lem}\label{lem:tsuboi-2}
Let $k\ge2$ be an integer.
Let $\{x_n\}_{n\in\bN}\sse I$ be a strictly increasing sequence such that $x_0=0$ and such that $\lim_{n\to\infty} x_n=1$. If
$f\in\Diff_+^k[0,1]$ and if $f(x_n)=x_n$ for all $n$, then we have \[\sum_{n=0}^{\infty}
\sup_{x\in [x_n,x_{n+1}]}
|f'(x)-1|^{1/(k-1)}<k\cdot \sup_{x\in[0,1]} |f^{(k)}(x)|^{1/(k-1)}<\infty.\]
\end{lem}
\bp 
Set $P:=\sup_{x\in[0,1]} |f^{(k)}(x)|$. We note that on the interval
\[
J_{n,k}:=[x_n,x_{n+k}]\]
the map $f$ is $k$--fixed (Definition~\ref{d:k-fixed}).
It follows  from Lemma~\ref{l:expansive2} that
\begin{align*}
&
\sum_{n=0}^{\infty}
\sup_{[x_n,x_{n+1}]}
|\bar f'(x)-1|^{1/(k-1)}
\le
\sum_{n=0}^{\infty}
\sup_{J_{n,k}}
|\bar f'(x)-1|^{1/(k-1)}
\\
&\le
\sum_n \sup_{J_{n,k}} |\bar f^{(k)}(x)|^{1/(k-1)} \cdot |J_{n,k}| \le k\cdot P^{1/(k-1)}<\infty.
\end{align*}
This completes the proof.\ep

We now argue that $\psi_r$ is not conjugate to a representation into $\Diff_+^{r+1}[0,1]$.
Indeed, suppose the contrary. 
We call the corresponding $C^{r+1}$ diffeomorphisms $\yt T$ and $\yt L$, and their common fixed points
$\{z_n\}_{n\in\bN}$.
For convenience of notation, write $T_n$ and $L_n$ for the restrictions of $\yt T$ and $\yt L$ to $K_n=[z_{n-1},z_n]$.

We claim that there exists some $v_n\in K_n$ such that
\[
\left(L_n^{n^r}\right)'(v_n)=2.\]
Indeed, 
from the equation
\[
\lambda_s \tau_t\lambda_s^{-1}=\tau_{te^s}\]
we have that
\[
L_n^{n^r} T_n L_n^{-n^r} = T_n^2.\]
In other words, we have a representation
\[
\BS(1,2)=\form{a,e\mid aea^{-1}=e^2}\longrightarrow\Diff_{K_n}^\infty(\bR)\]
defined as $a\mapsto L_n^{n_r}$ and $e\mapsto T_n$.
We see from Proposition~\ref{p:bmnr} 
 \ref{p:bmnr-m} that the existence of the required point $v_n$, as claimed.

From the claim and from the chain rule, we can find some point
$w_n\in K_n$ such that \[L_n'(w_n)\geq 2^{1/n^r}.\]
On the other hand, since we are assuming $\form{\yt T,\yt L}\le\Diff^{r+1}_+[0,1]$ we can apply Lemma~\ref{lem:tsuboi-2} and have that
\[
\sum_{n\ge1}
|2^{1/n^r}-1|^{1/r}\le
\sum_{n\ge1} |L_n'(w_n)-1|^{1/r}<\infty.\]
Note the estimate \[(1+t)^u-1\ge ut/2\] for $t\in[0,1]$ and $u\in[0,1]$.
Therefore, we have
\[
\sum_{n\ge1}
|2^{1/n^r}-1|^{1/r}
= \sum_{n\ge1} |(1+1)^{1/n^r}-1|^{1/r}\ge 2^{-1/r}\cdot \sum_{n\ge1} 1/n.\]
This is a contradiction, proving that $\form{\bar T,\bar L}\le\Diff_+^r[0,1]$ is not topologically conjugate into $\Diff_+^{r+1}[0,1]$.
We now have the following.

\begin{thm}
For each integer $r\ge1$ there exists a faithful representation
\[
\psi_r\co  \bZ\wr\bZ\longrightarrow \Diff_+^r[0,1]\]
that is not topologically conjugate into $\Diff_+^{r+1}[0,1]$.\end{thm}
\bp
It only remains for us to show the faithfulness of the representation. 
Consider an arbitrary nontrivial element $g\in \bZ\wr\bZ$. We can write $g$ in a normal form
\[
g = \prod_{i=-N}^N \lambda^i \tau^{k_i}\lambda^{-i} \cdot \lambda^m\]
for $N\ge0$ and suitable integers $k_i$ and $m$. 
Let us define for each $n\ge1$ that
\[
c_n:=
\left( \sum_{i=-N}^N \left(2^{1/n^r}\right)^ik_i
\right)u_n.\]
On each interval $J_n$, the map $g$ is topologically conjugate to the affine map
\[
x\mapsto 2^{m/n^r}x + c_n.\]

If $m\ne0$ then $\psi_r(g)$ is topologically conjugate to an
affine map with a nontrivial scaling part on each interval $J_n$.
This implies that $\psi_r(g)\ne1$. 

Let us now assume that  $m=0$ and that $k_i\ne0$ for some $i$. Clearly, we may assume $k_N\ne0$.
The map $\psi_r(g)$ is topologically conjugate to the affine translation by $c_n$
on each interval $J_n$.  Since the polynomial equation
\[
\sum_{i=-N}^N k_i x^i=0\]
has at most finitely many solutions, the value $c_n$ cannot be zero for all $n$. 
Together, these observations imply that $\psi_r(g)\ne1$, and that $\psi_r$ is faithful.
\ep

We remark that in Tsuboi's original paper~\cite{Tsuboi1987}, the exact isomorphism type (namely, $\bZ\wr\bZ$) of the group $\form{\bar T,\bar L}$ is not specified; rather, he only noted that it is a quotient of the group
\[
\form{a,b\mid [a,[a,b]]=1}.\]
As an application, let us consider a closed orientable surface $S_g$ of genus $g$. The preceding theorem can be combined with
Proposition~\ref{prop:suspension-unique} and allows us to conclude:

\begin{cor}\label{cor:tsuboi}
Let $g\geq 2$. For all integer $r\geq 1$, there exists a codimension one
$C^r$ foliation on $S_g\times I$ that is not smoothable to a $C^{r+1}$ foliation. 
\end{cor}

Our results on the existence of groups with prescribed critical regularity implies that for $g\geq 5$ and $ r\in\bR_{\geq 1}$, there
is a representation \[\psi_{ r}\in\Hom(\pi_1(S_g),\Diff_{[0,1]}^{ r}(\bR))\] that is not topologically conjugate to a representation into
\[\bigcup_{s > r}\Diff_{[0,1]}^{ r}(\bR).\]
We thus obtain a different perspective Corollary~\ref{cor:tsuboi}. In fact,
more is true: we get monodromy representation for $C^{ r}$ foliated structures on $S_g\times I$ that occur in no $C^{s }$ foliated
structure, for $s > r$.

\begin{cor}\label{cor:tsuboi-ref}
Let $g\geq 2$. For all $ r\geq 1$, there exists a
codimension one $C^{ r}$ foliation on $S_g\times I$ that is not smoothable to a $C^{s }$ foliation
whenever $s > r$.
\end{cor}

More precisely, the foliations furnished in Corollary~\ref{cor:tsuboi-ref} come from $I$--bundles over $S$ whose monodromy groups
do not arise from a $C^{s }$ group action on the interval, for $s > r$.

We can push the non--smoothability of foliated structures on $I$--bundles of the form $S\times I$ (where here $S$ is a surface)
even farther, producing non--smoothable foliations of closed $3$--manifolds subject to mild topological hypotheses.

Let $M$ be a closed, orientable $3$--manifold, and suppose that $H_2(M,\bZ)\neq 0$. By Poincar\'e duality, we have that
$H^1(M,\bZ)\neq 0$, and a nontrivial element $\phi\in H^1(M,\bZ)$ is classified by a smooth map \[\Phi\colon M\longrightarrow S^1.\]
Standard transversality arguments imply that $M$ admits an embedded, two--sided surface $S=\Phi^{-1}(p)$
for some suitable $p\in S^1$, possibly after modifying $\Phi$ by a homotopy, and that $S$ is Poincar\'e dual
to $\phi$. Fairly standard methods, which we shall not spell out in further detail here, imply that we may increase the genus of $S$
arbitrarily, so that for all $g\gg 0$, the manifold $M$ admits an embedded, two--sided copy of $S_g=S$. Since $S$ is two--sided, this
means that a tubular neighborhood of $S$ is homeomorphic to a trivial $I$--bundle $S\times I$.

Removing a copy of $S\times (0,1)$ from $M$ results in a manifold with boundary $N$. The following is a result of
Goodman~\cite{Goodman1975} (cf.~\cite{Thurston1976AM,CC1982}), whose proof we will not provide here.

\begin{prop}\label{prop:goodman}
Let $N$ be as above. The manifold $N$ admits a smooth foliation $\FF$ with both components of $\partial N$ as leaves.
\end{prop}

By gluing in a copy of $S\times I$ with an unsmoothable $C^{ r}$ foliation as furnished by Corollary~\ref{cor:tsuboi-ref}, we obtain
the following:

\begin{cor}\label{cor:3-mfld}
Let $M$ be a closed, orientable $3$--manifold with $H_2(M,\bZ)\neq 0$. Then there exists a codimension one $C^{ r}$ foliation
on $M$ that is not homeomorphic to a $\bigcup_{s > r} C^{s }$ foliation.
\end{cor}

\section{Right-angled Artin groups}\label{sec:raag}
Let $\Gamma$ be a finite, simplicial graph. That is, $\Gamma$ has a finite set of vertices $V(\gam)$, and a finite set of
undirected edges $E(\gam)$, so that
the $1$--complex given by the data of $V(\gam)$ and $E(\gam)$ forms a simplicial complex. One defines the 
\emph{right-angled Artin group}\index{right-angled Artin group}
\[A(\gam)=\form{V(\gam)\mid [v_i,v_j]=1\,\textrm{ whenever }\, \{v_i,v_j\}\in E(\gam)}.\]

The reader will easily check that the class of right-angled Artin groups contains free abelian groups and free groups, and can be
seen as interpolating between these two classes of groups.

\subsection{Abelian groups, free groups, and smooth actions}
Part of our interest in right-angled Artin groups comes from their actions on one--manifolds.
For some right-angled Artin groups, it is very easy to find actions on one--manifolds of high regularity.

For instance,
consider the free abelian group $\bZ^n$. Suppose that $0\neq X$ is a smooth vector field on $M\in\{I,S^1\}$. In the case where $M=I$,
assume that $X$ vanishes at $\partial I$. 
By integrating this vector field, we get a flow on $M$, which is to say a continuous injective homomorphism
\[\phi\colon \bR\longrightarrow\Diff_+^{\infty}(M).\] It follows that for all times $t,s\in\bR$, the corresponding diffeomorphisms $\phi(t)$ and
$\phi(s)$ commute with each other. By choosing $n$ times that are linearly independent over $\bQ$, say $\{t_1,\ldots,t_n\}$, we have that
the diffeomorphisms $\{\phi(t_1),\ldots,\phi(t_n)\}$ form a subgroup of $\Diff_+^{\infty}(M)$ that is isomorphic to $\bZ^n$.

At the other extreme, consider the free group $F_n$ of rank $n$. It is a bit less obvious how to find copies of 
$F_n$ inside of $\Diff_+^{\infty}(M)$,
though as it happens these subgroups are very common. We already observed this fact in Corollary~\ref{cor:free-abundant}.

For right-angled Artin groups that are not as structurally uncomplicated as free groups or abelian groups, it is significantly less
obvious why (or indeed if) there should be smooth actions on $I$ or on $S^1$. The reader will recall from Subsection~\ref{ss:abt}
that one can build faithful actions of $A(\gam)$ by homeomorphisms on $\bR$, by building a homomorphism $\phi_w$ for every element
$1\neq w\in A(\gam)$ that witnesses that $\phi_w(w)$ is nontrivial. Once such a sequence of homomorphisms is constructed, it is not
difficult to improve them to homomorphisms to $\Diff_+^{\infty}(\bR)$. Thus, we obtain:

\begin{thm}[See~\cite{BKK2014}]\label{thm:bkk-israel}
Let $\Gamma$ be a finite simplicial graph. Then there is an injective homomorphism \[A(\gam)\longrightarrow\Diff_+^{\infty}(\bR).\]
\end{thm}

If one tries to compactify the construction in Theorem~\ref{thm:bkk-israel} in order to get a faithful homomorphisms into $\Diff_+^{\infty}(I)$,
one encounters an analytic difficulty. The reader can check without great difficulty that, even though this construction furnishes a homomorphism
from $A(\gam)$ into $\Homeo_+(I)$, the Mean Value Theorem precludes this action from being $C^1$. By pushing the methods of constructing
$C^{\infty}$ actions outlined here to their extreme, we can prove the following.

\begin{prop}[cf.~\cite{KK2018JT}]\label{prop:c-infty-raag}
Let $A(\gam)$ be a right-angled Artin group that decomposes as a direct product of free products of free abelian groups. Then
there is an injective homomorphism  \[A(\gam)\longrightarrow\Diff_+^{\infty}(M),\] where $M\in\{I,S^1\}$.
\end{prop}

The graphs under the purview of Proposition~\ref{prop:c-infty-raag} are easy to describe combinatorially; these ideas will be developed
in Subsection~\ref{ss:p4} below. We will not justify some of the claims we make about passing between graphs and groups, though
the reader will not find them to be controversial. More details can be found in~\cite{Charney2007,KK2018JT,koberda21survey}.

A right-angled Artin group $A(\gam)$ is abelian if and only if the defining graph $\Gamma$ is complete. The group $A(\gam)$ is a free product
of free abelian groups if and only if every connected component of $\Gamma$ is complete. We have that $A(\gam)$ decomposes as a
(nontrivial) direct
product if and only if $\Gamma$ decomposes as a
(nontrivial) \emph{join}\index{join of graphs}, which is to say $V(\gam)=V(J_1)\cup V(J_2)$, where the graphs $J_1$ and $J_2$
spanned by $V(J_1)$ and $V(J_2)$ share no vertices (and are both nonempty),
and where every vertex $v_1\in V(J_1)$ spans an edge with every vertex
$v_2\in V(J_2)$. Therefore, we have that $A(\gam)$ decomposes as a direct product of free products of free abelian groups if and only
if $\Gamma$ decomposes as a (possibly trivial) join, where each component of the join is a disjoint union of complete graphs.

\begin{proof}[Proof of Proposition~\ref{prop:c-infty-raag}]
In the preceding paragraphs, we have shown how to obtain nonabelian free groups and free abelian groups in $\Diff_+^{\infty}(I)$, and
it is easy to see that one can arrange such an action to be supported on an arbitrary prescribed nondegenerate subinterval $J\sse I$.
Now, let \[A(\gam)\cong A(\gam_1)\times\cdots\times A(\gam_k),\] where $A(\gam_i)$ is a free product of free abelian groups for
$1\leq i\leq k$. We will choose $k$ disjoint intervals $\{J_1,\ldots,J_k\}$ of $I$ such that $(\partial J_i)\cap(\partial J_{\ell})=\varnothing$ 
for $i\neq\ell$.

Let $N_i$ be the maximal rank of a free abelian subgroup of $A(\gam_i)$. An easy applications of the
Kurosh Subgroup Theorem~\cite{Serre1977} shows that the group $\bZ*\bZ^{N_i}$ contains all groups of the form
\[\bZ^{n_1}*\cdots*\bZ^{n_{\ell}},\] provided that \[N_i\geq \max \{n_1,\ldots,n_{\ell}\}.\] Clearly, to prove the result for $M=I$, it suffices
to establish that $\bZ*\bZ^N\le \Diff_+^{\infty}(I)$ for $N\in\bN$,
in such a way that this group extends by the identity past the endpoints of $I$. 
Indeed,
then we can build an action of $\bZ*\bZ^{N_i}$ on $J_i$ for $1\leq i\leq k$, furnishing a faithful action of $A(\gam)$ by $C^{\infty}$
diffeomorphisms on $I$. Clearly we may then identify the endpoints of $I$ to get an action on $S^1$.

Choosing a smooth flow on $I$ whose derivatives of all orders agree with those of the identity at $0$ and $1$, we may find a
group of diffeomorphisms isomorphic to $\bZ^N$. Clearly we may assume that the vector field giving rise to this flow has no zeros
in the interior of the interval, so that this copy of $\bZ^N$ acts freely on $I$. We then consider all diffeomorphisms of $I$
whose derivatives of all orders agree with those of the identity at $0$ and $1$. 
By a Baire Category argument very similar to that used
to find a copy of $F_2\le \Diff_+^{\infty}(I)$ above, a generic choice of diffeomorphism of $I$ will
combine with the abelian group we have produced in order to generate a copy of $\bZ*\bZ^N$. We leave the remaining details to
the reader.\end{proof}

In the proof of Proposition~\ref{prop:c-infty-raag}, the fact that $\bZ^N$ is acting freely is essential; it can be weakened to an assumption
of \emph{full support}\index{fully supported homeomorphism},
i.e.~that a nontrivial element of $\bZ^N$ cannot have an interval's worth of fixed points. If one does not have 
such an assumption, then the induction fails. Indeed, let $t$ generate the free $\bZ$--factor of $\bZ*\bZ^N$. If $x\in I$ and $t(x)$ is fixed
by $g\in\bZ^N$ then $(t^{-1}gt)(x)=x$. If $t(x)$ lies in the interior of $\Fix(g)$ then a small perturbation of $t$ will not prevent $x$ from being
a fixed point of $t^{-1}gt$. Thus, the set of ``good" choices of $t$ and $x$ is no longer dense in the set $\Diff_+^{\infty}(I)\times I$, and
the Baire Category Theorem is not applicable. This issue is not an artifact of the proof; indeed, we will see in the next section that
the right-angled Artin groups covered by Proposition~\ref{prop:c-infty-raag} are exactly the ones which embed in $\Diff_+^{\infty}(I)$.

\subsection{$P_4$, the cograph hierarchy, and $C^{1+\mathrm{bv}}$ actions}\label{ss:p4}
Given that this natural construction cannot have good regularity properties, we are motivated to consider the following question, which
M.~Kapovich attributes to Kharlamov~\cite{Kapovich2012}.

\begin{que}\label{que:kharlamov}
For which graphs $\Gamma$ is there an injective homomorphism \[A(\gam)\longrightarrow \Diff^{\infty}(M),\] where $M\in\{I,S^1\}$?
\end{que}

Proposition~\ref{prop:c-infty-raag} gives us a partial answer to Question~\ref{que:kharlamov}.
The machinery developed in this book allows us to give a complete answer.

Historically, Question~\ref{que:kharlamov} was resolved in two steps. Strictly speaking, the second step subsumed the first, and so it
will suffice to only carry out the second. It will be useful for us to retain the historical progression of the resolution of
Question~\ref{que:kharlamov}, because of connections to mapping class group we will explore in Section~\ref{sec:mcg} below.

For $n\geq 1$, we will write $P_n$ for the path of length $n-1$, which is a graph with exactly $n$ vertices which can be thought of
as ordered linearly, with the linear ordering inducing the edge relation. We have illustrated the graph $P_4$ in Figure~\ref{f:p4}.

\begin{figure}[h!]
  \tikzstyle {bv}=[black,draw,shape=circle,fill=black,inner sep=1pt]
  \begin{center}
\begin{tikzpicture}[main/.style = {draw, circle}]
\node[main] (1) {$a$};
\node[main] (2) [right of=1] {$b$};
\node[main] (3) [right of=2] {$c$}; 
\node[main] (4) [right of=3] {$d$};

\draw (1)--(2)--(3)--(4);
\end{tikzpicture}%
\caption{The graph $P_4$.}
\label{f:p4}
\end{center}
\end{figure}
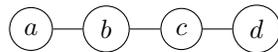

The graphs $\{P_n\}_{n\geq 1}$ play an important role in the theory of right-angled Artin groups. Trivially, we have $A(P_1)\cong\bZ$ and
$A(P_2)\cong\bZ^2$. Slightly less trivially, we have $A(P_3)\cong F_2\times \bZ$. The right-angled Artin group $A(P_4)$ is the first
interesting right-angled Artin group on a path. It was observed by Droms~\cite{Droms1987} that if $A(\gam)$ is a right-angled Artin group
then every finitely generated subgroup of $A(\gam)$ is a right-angled Artin group if and only if $\Gamma$ has no full subgraph isomorphic
to the path $P_4$ and no full subgraph isomorphic to the cycle of length four $C_4$. It is true though nonobvious that $\Gamma$ contains
one of these graphs as a full subgraph if and only if $A(\gam)$ contains one of the corresponding right-angled Artin groups as a
subgroup (see~\cite{Kambites2009,KK2013}, cf.~Proposition~\ref{prop:gam-p4}).

The right-angled Artin group $A(P_4)$ turns out to be universal for right-angled Artin groups on finite forests, in the following sense:

\begin{prop}[See~\cite{KK2013}]\label{prop:raag-forest}
Let $F$ be a finite forest. Then there is an injective homomorphism $A(F)\longrightarrow A(P_4)$.
\end{prop}

The first progress on resolving Question~\ref{que:kharlamov} was obtained by the authors with Baik.

\begin{thm}[See~\cite{BKK2019JEMS}]\label{thm:a4-c2}
There is no injective homomorphism \[A(P_4)\longrightarrow \Diffb(M),\] where $M\in\{I,S^1\}$.
\end{thm}

It follows that if $A(\gam)$ admits a faithful action on a compact one--manifold $M$ with regularity $C^{1+\mathrm{bv}}$ or higher
then not only is $P_4$ not a full subgraph of $\Gamma$, but in fact $A(\gam)$ cannot contain a copy of $A(P_4)$ as a subgroup. This latter
condition may seem weaker than the former, but it turns out that they are equivalent to each other.

\begin{prop}[See~\cite{KK2013}]\label{prop:gam-p4}
There is an injective homomorphism $A(P_4)\longrightarrow A(\gam)$ if and only if $P_4$ occurs as a full subgraph of $\Gamma$.
\end{prop}

Proposition~\ref{prop:gam-p4} says that if $A(\gam)$ contains a copy of $A(P_4)$ then it contains an ``obvious" copy of $A(P_4)$ (though
of course not every copy of $A(P_4)$ need be ``visible" in the graph). Thus, to fully answer Question~\ref{que:kharlamov}, we need only
consider graphs containing no full paths of length three. It turns out that these are well--understood by graph theorists.

A finite, simplicial graph $\Gamma$ is called a \emph{cograph}\index{cograph} 
if it contains no full copy of $P_4$. Such graphs fall within an inductive hierarchy,
which makes them much easier to understand. The key observations are as follows: if $\gam_1$ and $\gam_2$ are cographs, then so is
their disjoint union. Moreover, it is easy to check by hand that the join $\gam_1*\gam_2$ is also a cograph. We recall briefly that
$\gam_1*\gam_2$ can be thought of as the disjoint union of $\gam_1$ and $\gam_2$, augmented by adding an edge between each vertex
of $\gam_1$ and each vertex of $\gam_2$. It turns out that these two operations completely characterize cographs. We build a
hierarchy of graphs as follows.

\begin{enumerate}
\item
Let $\KK_0$ be a singleton vertex.
\item
For $n=2i+1$ for $i\geq 0$, we set $\KK_n$ to consist of finite joins of elements in $\KK_{j}$ for $j\leq n-1$.
\item
For $n=2i+2$ for $i\geq 0$, we set $\KK_n$ to consist of finite disjoint unions of elements in $\KK_{j}$ for $j\leq n-1$.
\item
We set \[\KK=\bigcup_{n\geq 0} \KK_n.\]
\end{enumerate}

We have the following fact, which is well--known from graph theory~\cite{cograph1,cograph2,cograph3}.
The reader may consult~\cite{KK2013} for an algebraic proof that
uses right-angled Artin groups and relies fundamentally on Proposition~\ref{prop:gam-p4}.

\begin{prop}\label{prop:cograph}
If $\Gamma$ is a cograph then $\gam\in\KK_n$ for some $n\geq 0$.
\end{prop}

Proposition~\ref{prop:cograph} gives an algebraic description of right-angled Artin groups on cographs.

\begin{prop}\label{prop:cograph-alg}
The class of groups \[\{A(\gam)\mid \gam \textrm{ is a cograph}\}\] coincides with the
smallest class of finitely generated groups that contains
$\bZ$, is closed under taking direct products, and is closed under taking free products.
\end{prop}

Clearly if $\gam\in\KK_0$ then $A(\gam)\cong \bZ$. If $\gam\in\KK_1$ then $A(\gam)$ is a free abelian group. If $\gam\in\KK_2$ then
$A(\gam)$ is a free product of free abelian groups. If $\gam\in\KK_3$ then $A(\gam)$ is a direct product of free products of free abelian
groups. The groups under the purview of Proposition~\ref{prop:c-infty-raag} therefore coincide exactly with $A(\gam)$ for $\gam\in\KK_3$.

It turns out that the group $(F_2\times\bZ)*\Z$, which has a defining graph $\Lambda$ given by a path of length two together with
an isolated vertex (see Figure~\ref{f:k4}), characterizes right-angled Artin groups whose defining graphs lie in $\KK_4$.

\begin{figure}[h!]
  \tikzstyle {bv}=[black,draw,shape=circle,fill=black,inner sep=1pt]
  \begin{center}
\begin{tikzpicture}[main/.style = {draw, circle}]
\node[main] (1) {$a$};
\node[main] (2) [right of=1] {$b$};
\node[main] (3) [right of=2] {$c$}; 
\node[main] (4) [right of=3] {$d$};

\draw (1)--(2)--(3);
\end{tikzpicture}%
\caption{The defining graph $\Lambda$ for $(F_2\times\bZ)*\bZ$.}
\label{f:k4}
\end{center}
\end{figure}
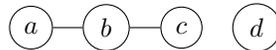

\begin{prop}[See~\cite{KK2018JT}]\label{prop:k4}
We have that $\gam\in\KK_n$ for $n\geq 4$ if and only if $A(\gam)$ contains a copy of $(F_2\times\bZ)*\bZ$.
\end{prop}

The reader might note that $P_4$ does not contain a full subgraph isomorphic to $\Lambda$. However, the reader may check the following
facts.
Consider the labeling of the vertices of $P_4$ as in Figure~\ref{f:p4}. Let \[\phi_d\colon A(P_4)\longrightarrow\bZ/2\bZ\] be the
homomorphism sending $d$ to the generator of $\bZ/2\bZ$, and where $\{a,b,c,d^2\}\sse\ker\phi_d$. Then $\ker\phi_d$ is itself
isomorphic to a right-angled Artin group whose defining graph is a path of length four, with an extra degree vertex sprouting from the
middle vertex. Thus, we get an explicit realization \[A(P_5)\cong \form{ a,b,c,b^d,a^d}.\] Since $\Lambda$ is a full subgraph
of $P_5$, we get that $(F_2\times\bZ)*\bZ$ is a subgroup of $A(P_4)$.

\begin{proof}[Proof of Proposition~\ref{prop:k4}]
Clearly if $\gam\in\KK_n$ for $n\leq 1$ then $(F_2\times\bZ)*\bZ$ does not occur as a subgroup of $A(\gam)$.
It is an easy consequence of the Kurosh Subgroup Theorem that if $\gam\in\KK_2$ then $(F_2\times\bZ)*\bZ$
also does not occur as a subgroup
of $A(\gam)$. To see that if $\gam\in\KK_3$ then $A(\gam)$ does not contain a copy of $(F_2\times\bZ)*\bZ$, let $G$ and $H$ be
arbitrary groups and let \[\psi\colon (F_2\times\bZ)*\bZ\longrightarrow G\times H\] be a homomorphism. Let $p_G$ and $p_H$
be the projections of $G\times H$ onto the two factors. It is an easy exercise in combinatorial group theory that if
$\psi$ is injective then either $p_g\circ\psi$ or $p_H\circ\psi$ is injective. Thus, if $\gam\in\KK_n$ for $n\leq 3$ then $(F_2\times\bZ)*\bZ$
is not a subgroup of $A(\gam)$.

To establish the other direction, suppose that $\Gamma$ in $\KK_n\setminus\KK_{n-1}$ for $n\geq 4$.
Decomposing $\Gamma$ maximally into its join factors,
we have that some join factor of $\Gamma$ also lies in $\KK_{n-1}\setminus\KK_{n-2}$ for $n-1\geq 4$.
Indeed, if each join factor of $\Gamma$ lies in $\KK_3$
then $\Gamma$ already lies in $\KK_3$. Let $\gam_0$ be a join factor of $\Gamma$ that does not lie in $\KK_3$.
Then $\gam_0$ is a nontrivial disjoint union of graphs,
and at least one of the components lies in  \[\bigcup_{n\geq 3}\KK_n\setminus\KK_2.\]
Let $\gam_0^1$ be such a component, and let $\gam_0^2$
be another component of $\gam_0$. An arbitrary vertex of $\gam_0^2$ generates a copy of $A(\gam_0^1)*\bZ$ with $\gam_0^1$.
Since $\gam_0^1$ is connected and lies in  \[\bigcup_{n\geq 3}\KK_n\setminus\KK_2,\]
we have that $A(\gam_0^1)$ is nonabelian and hence
$\gam_0^1$ admits two nonadjacent vertices. 
If $v$ and $w$
are two vertices of $\gam_0^1$ that are not adjacent, then  an arbitrary path between them (which exists since $\gam_0^1$ is connected)
gives a copy of $A(\Lambda)\cong (F_2\times\bZ)*\bZ$ in $A(\gam)$.
\end{proof}

The group $(F_2\times\bZ)*\bZ$ is a poison subgroup for group actions on compact one--manifolds in regularity at least $C^{1+\mathrm{bv}}$.
We begin with the case of the interval.

\begin{thm}\label{thm:f2-int}
There is no injective homomorphism $(F_2\times\bZ)*\bZ\longrightarrow\Diffb(I)$.
\end{thm}
\begin{proof}
Suppose the contrary, and let such an injective homomorphism be given. We will write \[(F_2\times\bZ)*\bZ=\form{ a,b,c,t},\]
where $\form{ a,b} \cong F_2$, where $c$ commutes with $a$ and $b$, and where the cyclic subgroup generated by $t$
splits off as a free factor. We will abuse notation and use the names for the generators of $(F_2\times\bZ)*\bZ$ to denote diffeomorphisms
of $I$.

Let $J_c$ be a component of the support of $c$, and let $J_s$ be a component of the support of $s\in\form{ a,b}$.
If $J_c\cap J_s\neq\varnothing$
then we must have $J_c=J_s$. This follows from an easy application of Kopell's Lemma
(Theorem~\ref{thm:kopell}), and explicitly from Proposition~\ref{prop:disj-ab}.

If $J_c=J_s$ for some nontrivial $s\in\form{ a,b}$ then $J_c=J_a$ or $J_c=J_b$ for
suitable components of the support of $a$ or of $b$.
If $s$ is a nontrivial product of commutators, then we must in fact have $J_c=J_a=J_b$. If $s$ has a fixed point in the interior of $J_c$
then Theorem~\ref{thm:kopell} again implies that $s$ is the identity, and so $\form{ a,b}$ acts freely on $J_c$. It follows from
H\"older's Theorem (see~\ref{thm:holder}) that the action of $\form{ a,b}$ on $J_c$ factors through the abelianization.
It follows that if $s$ lies in the commutator subgroup of $\form{ a,b}$ then each component $J_s$ of $\supp s$ is disjoint from $\supp c$.

Setting $s=[a,b]$, we have that $\form{ s,c}\cong\bZ^2$ and \[(\supp s)\cap(\supp c)=\varnothing.\] It follows from the $abt$--Lemma
(Theorem~\ref{thm:abt}) that \[\form{ s,c,t} \ncong\bZ^2*\bZ.\] Since the elements $\{s,c,t\}$ generate a copy of
$\bZ^2*\bZ\le (F_2\times\bZ)*\bZ$, it follows that the initial action of $(F_2\times\bZ)*\bZ$ was not
faithful.
\end{proof}

The case $M=S^1$ is slightly more complicated, since elements of $\Diffb(S^1)$ may have no fixed points, and so the support of a
diffeomorphism may be the whole circle and not just a disjoint union of intervals.

\begin{thm}\label{thm:f2-circ}
There is no injective homomorphism $(F_2\times\bZ)*\bZ\longrightarrow\Diffb(S^1)$.
\end{thm}

Before proving Theorem~\ref{thm:f2-circ}, we reduce to the case where the $abt$--Lemma can be applied.

\begin{lem}\label{lem:tame}
Let \[\phi\colon F_2\times\bZ\longrightarrow\Diffb(S^1)\] be an injective homomorphism. Then there exists a copy $\bZ^2\le F_2\times\bZ$
whose image under $\phi$ is generated by diffeomorphisms $a$ and $b$ such that \[(\supp a)\cap (\supp b)=\varnothing.\]
\end{lem}
\begin{proof}
Write \[F_2\times\bZ\cong\form{ a,b}\times\form{c},\] and as before we identify the names of the generators with elements in
$\Diffb(S^1)$. Proposition~\ref{prop:tame} shows that the rotation number of $c$ must be rational. In particular,
after replacing $c$ by a positive power if necessary, we may assume that $\Fix c\neq\varnothing$ (Proposition~\ref{prop:rot-easy}).
By Item 4 of Proposition~\ref{prop:tame}, we have that the rotation number restricts to a homomorphism
\[F_2=\form{ a,b}\longrightarrow S^1,\] and so if $h\in F_2'$ then $h$ has a nonempty fixed point set.
Applying Proposition~\ref{prop:tame} again, we have that \[(\supp c)\cap(\supp F_2'')=\varnothing.\] Since $F_2$ is a nonabelian free group,
we have that $F_2''$ is nontrivial. The conclusion of the lemma follows.
\end{proof}

\begin{proof}[Proof of Theorem~\ref{thm:f2-circ}]
Suppose the contrary. Lemma~\ref{lem:tame} furnishes a copy of $\bZ^2*\bZ$ in the image of the homomorphism such that $\bZ^2$ is
generated by diffeomorphisms with disjoint support. This violates Theorem~\ref{thm:abt}.
\end{proof}

Combining all the preceding discussion, we have the following complete answer to Question~\ref{que:kharlamov}:

\begin{thm}\label{thm:kharlamov}
Let $\Gamma$ be a finite simplicial graph. Then the following are equivalent.
\begin{enumerate}[(1)]
\item
We have $A(\gam)\le \Diff_+^{\infty}(I)$.
\item
We have $A(\gam)\le \Diff_+^{\infty}(S^1)$.
\item
We have $A(\gam)\le \Diffb(I)$.
\item
We have $A(\gam)\le \Diffb(S^1)$.
\item
We have $\gam\in\KK_3$.
\end{enumerate}
\end{thm}

It is a trivial consequence of the fact that right-angled Artin groups are residually torsion--free nilpotent~\cite{DK1992a}
and Theorem~\ref{thm:ff2003}
that for a finite simplicial graph $\Gamma$, we have $A(\gam)\le \Diff_+^1(M)$ for $M\in\{I,S^1\}$. Determining the exact
critical regularity of right-angled Artin groups is an open question in general. The following is the state of the art at the time of
this book's writing.

\begin{thm}[\cite{KKR2020,KKR2021}]\label{thm:kkr-2020}
Suppose $(F_2\times F_2)*\bZ\le A(\gam)$. Then for $\tau>0$, we have that $A(\gam)$ admits no faithful $C^{1+\tau}$ action on $I$
or $S^1$.
\end{thm}

In particular, Theorem~\ref{thm:kkr-2020} implies that if $\Gamma$ has a square and is not a nontrivial join then the critical regularity of
$A(\gam)$ on a compact one--manifold is exactly one. The proof of Theorem~\ref{thm:kkr-2020} is beyond the scope of this book.

\section{Mapping class groups}\label{sec:mcg}

Part of the authors' interest in groups acting on one manifolds arose naturally from their investigations of the relationship between mapping 
class groups of surfaces and right-angled Artin groups. Let $S$ be an orientable surface of finite type. That is, $S$ is an orientable,
two--dimensional real manifold such that the Euler characteristic of $S$ is finite. It is a standard fact from combinatorial topology that
$S$ is entirely determined by its genus $g$, its number of punctures $n$, and its number of boundary components $b$.

The mapping class group $\Mod(S)$ is defined to be the group of isotopy classes of orientation preserving homeomorphisms of $S$,
i.e. \[\Mod(S)=\pi_0(\Homeo_+(S)).\] Observe that, by definition, the group $\Mod(S)$ is identified with a group of outer automorphisms
of $\pi_1(S)$. We only obtain a group of outer automorphisms, since usually there is no canonical basepoint that is preserved by
homeomorphisms. If one chooses a preferred
marked point $p\in S$ and requires all homeomorphisms and isotopies to preserve $p$, then
one obtains a subgroup of $\Aut(\pi_1(S,p))$. It is a deep result of Dehn--Nielsen--Baer that if $S$ is closed, then the mapping class group
of $S$ is canonically identified with a subgroup of $\Out(\pi_1(S))$ of index two. If one chooses a marked point, then the mapping class group
of $S$ preserving $p$ is canonically identified with a subgroup of $\Aut(\pi_1(S))$ of index two. These results
admit suitable (but not necessarily obvious) generalizations to surfaces
that are not closed~\cite{FM2012}.

\subsection{Continuous actions}
The relationship between mapping class groups and homeomorphisms arises from the following
result originally due to Nielsen~\cite{Nielsen1927,HT1985}.

\begin{thm}\label{thm:nielsen}
Let $S$ be an orientable surface of finite type such that $g\geq 2$, with $n=1$ and $b=0$. Then there is an injective homomorphism
\[\phi\colon\Mod(S)\longrightarrow\Homeo_+(S^1).\]
\end{thm}
\begin{proof}[Sketch of proof of Theorem~\ref{thm:nielsen}]
We view $S$ as a closed surface of genus $g$ with a single marked point $p$. Let $\psi\in \Mod(S)$. We choose $\Psi\in\Homeo_+(S)$
lifting $\psi$, such that $\Psi(p)=p$. We lift $\Psi$ to the universal cover $\bH^2$ of $S$, so that we get a homeomorphism
$\yt\Psi\in\Homeo_+(\bH^2)$ commuting the action of $\pi_1(S)$ on $\bH^2$. The standard Morse Lemma from hyperbolic geometry
implies that $\yt\Psi$ induces a homeomorphism $\partial\yt\Psi$ on $\partial\bH^2\cong S^1$, which is unique. Since $S$ is compact,
modifying $\Psi$ by a homotopy does not affect $\partial\yt\Psi$, and so this latter
map depends only on the mapping class $\psi$, and a choice
of lift to $\bH^2$. By choosing a preferred preimage $q\in\bH^2$ lying over $p$, we obtain a preferred lift of $\Psi$ fixing $q$.
This furnishes a homomorphism from $\Mod(S)$ to $\Homeo_+(S^1)$, where $\phi(\psi)$ is $\partial\yt\Psi$ for the preferred lift.

If $\psi\in\Mod(S)$ is nontrivial then there is a nontrivial free homotopy class on $S$ that is not fixed by $\psi$. Taking a geodesic representative
$\gamma$ of such a free homotopy class, we have that the endpoints in $\partial\bH^2$ of $\gamma$ and $\psi(\gamma)$ do not coincide.
It follows that $\psi$ acts nontrivially on $\partial\bH^2$ and hence induces a nontrivial homeomorphism of $S^1$.
\end{proof}

One can say significantly more about $\phi$. For one, the action of $\Mod(S)$ on $S^1$ under $\phi$ is minimal. The regularity properties
of $\phi$ are somewhat more difficult to understand. One fact which is known is that the image of $\phi$ consists of 
\emph{quasi--symmetric}\index{quasi--symmetric homeomorphism}
homeomorphisms. 
Here, a homeomorphism $f\in\Homeo_+(S^1)$ is called quasi--symmetric if there exists an increasing
function \[\omega\colon [0,\infty)\longrightarrow [0,\infty)\] such
that for all triples
of distinct points $\{x,y,z\}\sse S^1$, we have \[\frac{|f(x)-f(y)|}{|f(x)-f(z)|}\leq\omega\left(\frac{|x-y|}{|x-z|}\right).\] The fact that
the map $\phi$ in Theorem~\ref{thm:nielsen} is a consequence of the fact if $f$ is a quasi--isometry of hyperbolic space $\bH^2$,
then the map $\partial f$ induced on $S^1$ is quasi--symmetric (see~\cite{kap-lect}, for instance). 
Quasi--symmetry is a generalized notion of bi--Lipschitz
homeomorphism, and we have from Theorem~\ref{thm:lip-conj} that one can find an action of $\Mod(S)$ on $S^1$ by bi--Lipschitz
homeomorphisms.

Theorem~\ref{thm:nielsen} has the consequence that mapping class groups of surfaces preserving a marked point are \emph{circularly
orderable}\index{circular ordering},
which for countable groups is equivalent to being a subgroup of $\Homeo_+(S^1)$. See Appendix~\ref{ch:append2}. Circular
orderability of groups of homeomorphisms of the circle implies that mapping class groups of closed surfaces (without marked points)
cannot act faithfully on the circle, since they contain non-cyclic torsion. Since this torsion disappears in a finite index subgroup, it
is unknown whether or not there exists a finite index subgroups of the closed mapping class group acting faithfully on the circle.

If $S$ has a boundary component $B$ then we may consider the mapping class group of $S$ preserving $B$ pointwise, which we will
write $\Mod(S,B)$. It turns out that the mapping class group $\Mod(S,B)$ can be made to act faithfully on the real line, by a result
that is analogous to Theorem~\ref{thm:nielsen}. When $S$ has a boundary component, then one can pick a preferred lift of the boundary
component to $\bH^2$, which then becomes an invariant copy of the real line connecting two points in $\partial \bH^2$, and the mapping
class group acts on it faithfully.

\begin{thm}[See~\cite{HT1985}]\label{thm:thurston-handel}
If $S$ has a boundary component $B$ then there is an injective homomorphism \[\Mod(S,B)\longrightarrow\Homeo_+(\bR).\]
\end{thm}

Thus, mapping class groups of surfaces with boundary are \emph{orderable}\index{linear ordering},
which like in the case of the circle, characterizes countable
subgroups of $\Homeo_+[0,1]$ (see Appendix~\ref{ch:append2}). The mapping class group of a surface with a marked point has torsion
(though only of a cyclic kind), and hence is not orderable. This torsion disappears in a finite index subgroup of the mapping class group,
and it is an open problem to determine whether this mapping class group admits a finite index subgroup acting faithfully on the interval.

\subsection{Smoothing to $C^2$ and beyond and connections to right-angled Artin groups}
It was a conjecture for over a decade whether or not the action of
a finite index subgroup of $\Mod(S)$ could be smoothed, which is to say topologically conjugated
to a differentiable action. A more general question allows one to abandon the constraints of Nielsen's action,
and so, one can ask: is there a finite index subgroup $G\le \Mod(S)$ and an injective homomorphism $G\longrightarrow\Diff_+^k(M)$, 
for $M\in\{I,S^1\}$?

This question can be fit into the general setup of the Zimmer Conjecture, which roughly asserts that ``large" groups cannot
act in interesting ways on ``small" compact manifolds. The meaning of ``large" and ``small" are usually context dependent; for instance,
one can consider a group to be large if it is a lattice in a real simple Lie group of rank $n\geq 2$, and a compact manifold to be small if
is has dimension smaller than $n$. 
A lattice in a Lie group, as a mathematical object, is only well--defined up to finite index; so, one usually considers a group
up to commensurability, for the purposes of Zimmer's Conjecture.
Interesting actions are generally continuous actions (though one often considers actions of higher
regularity, or actions that in addition preserve some extra data like a measure or a form) that do not factor through finite quotients.

For mapping class groups, deciding whether they are large or small is a matter of philosophical debate. Nielsen's action does furnish
a faithful action of a mapping class groups on the circle. On the other hand, results of Ghys, Farb--Franks, and 
Parwani~\cite{FF2001,Ghys1999,Parwani2008} have shown
that the full mapping class group does not admit faithful $C^2$ actions on $I$ or $S^1$ when the genus of $S$ is at least three, and
does not admit faithful $C^1$ actions when the genus of $S$ is at least $6$. These latter results rely on the full
power of the mapping class group, and the methods are not robust under passing to finite index subgroups.

Usually, finite index subgroups of mapping class groups are quite difficult to understand. For our purposes, right-angled Artin groups
come to the rescue. Let $G$ be an arbitrary group, and let $A(\gam)\le G$ be a right-angled Artin subgroup. It turns out that every finite index
subgroup of $G$ also contains a copy of $A(\gam)$. Indeed, if $H\le G$ has index $n$, then there is an $N$ such
that for all $g\in G$, we have $g^N\in H$. If $\{v_1,\ldots,v_k\}$ are the vertices of a graph $\Gamma$ viewed as generators of $A(\gam)$, the
the elements $\{v_1^N,\ldots,v_k^N\}$ also generate a group isomorphic to $A(\gam)$. In short,
the property of containing a right-angled Artin group
of a particular isomorphism type is stable under passing to finite index subgroups.

Since by Theorem~\ref{thm:a4-c2} we know that $A(P_4)$ is not a subgroup of $\Diffb(M)$ in order to rule out a faithful action of
a finite index subgroup of a mapping class group by $C^{1+\mathrm{bv}}$ diffeomorphisms on $M$, it would suffice to find copies of
$A(P_4)$ inside of these mapping class groups. These are furnished in abundance by a result of the second author.

Let $f\in\Mod(S)$ be a mapping class. We say that $f$ is \emph{reducible}\index{reducible mapping class} if there is a
nonempty collection $R$ of (isotopy classes of)
essential, nonperipheral (i.e.~not parallel to a boundary component of puncture), simple
closed curves on $S$ such that $f$ fixes $R$. The collection $R$ is called a \emph{reduction system}\index{reduction system}
for $f$. The minimal such $R$
(with respect to inclusion) always exists~\cite{BLM1983} and is called the \emph{canonical reduction system}\index{canonical reduction
system} for $f$.

We say that $f$ is \emph{pure}\index{pure mapping class} if whenever $S_0\sse S\setminus
R$ is a component, then
$f$ preserves $S_0$, preserves $R$ component-wise,
and the restriction of $f$ to $S_0$ is either the identity on $S_0$ (away from the boundary) or acts on $S_0$
irreducibly. It is a standard fact that every mapping class in $\Mod(S)$ admits a power that is pure.

If $f$ is pure and $\gamma$ is a component of $R$, then $f$ may act nontrivially in an annular neighborhood of $\gamma$. Namely,
$f$ may cut open $S$ along $\gamma$ and reglue with a full twist. Such a mapping class is called a \emph{Dehn twist}\index{Dehn twist}.

If $f$ is pure, we say that $f$ has \emph{connected support}\index{connected support}
if $f$ satisfies one of the following conditions:
\begin{enumerate}[(1)]
\item
The mapping class $f$ is a Dehn twist about a simple closed curve.
\item
The mapping class $f$ is not a Dehn twist about a simple closed curve, and the restriction of $f$ to $S\setminus R$ is irreducible
for exactly one component $S_0\sse S\setminus R$.
Moreover, we require that, essentially, $f$ does not perform Dehn twists about curves in its
canonical reduction system that do not border on $S_0$. In other words, $R=\partial S_0$.
\end{enumerate}

We direct the reader to~\cite{FM2012} for a more detailed discussion of the foregoing facts about mapping class groups.

We say that a collection of pure mapping classes with connected support is \emph{irredundant}\index{irredundant mapping classes}
if no two mapping classes in it
generate a virtually cyclic subgroup of $\Mod(S)$.

Let $X=\{f_1,\ldots,f_k\}$ be
irredundant pure mapping classes with connected support. We build a graph $\gam_X$ by taking one vertex for each element of $X$,
and drawing an edge between two vertices if the two mapping classes commute in $\Mod(S)$. The graph $\gam_X$
is called the \emph{co--intersection}\index{co--intersection graph}
graph; it is so named since generally two vertices will be adjacent in $\gam_X$ if and only if their supports are disjoint, up to isotopy
in $S$.
The following result, whose proof is beyond the scope of this book,
gives a systematic way of producing right-angled Artin groups in mapping class groups.

\begin{thm}[See~\cite{Koberda2012}]\label{thm:mcg-raag}
Let $X=\{f_1,\ldots,f_k\}$ be
irredundant pure mapping classes with connected support with co--intersection graph $\gam_X$. Then there exists an $N>0$ such
that for all $n\gg N$, we have \[\form{ f_1^n,\ldots,f_k^n}\cong A(\gam_X)\le \Mod(S).\]
\end{thm}

Theorem~\ref{thm:mcg-raag} was the first comprehensive algebraic results about subgroups of mapping class groups generated by
sufficiently high powers of mapping classes. There are many results that are closely related to Theorem~\ref{thm:mcg-raag}, some which
preceded it~\cite{CP2001,CW2007}, some which were contemporary and focussed more on the geometric aspects of the embedding
$A(\gam_X)\longrightarrow\Mod(S)$~\cite{CLM2012},
and some which focussed on effectivization~\cite{Seo2021,Runnels2021}. For the sake of space,
we will not survey the literature on this topic in greater depth.

In order to find copies of $A(P_4)$ in the largest variety of mapping class groups, we consider copies of $A(P_4)$ generated by the
mapping classes that take up the ``least amount of space", which is to say they each are supported on a single annulus. If
we consider four distinct (i.e.~pairwise non--isotopic), essential, nonperipheral, simple closed
curves $\{\gamma_1,\ldots,\gamma_4\}$ which have the property that $\gamma_i\cap\gamma_j\neq\varnothing$
if and only if $|i-j|=1$ (where this intersection is minimized over the isotopy classes of these curves), then the curves
$\{\gamma_1,\ldots,\gamma_4\}$ form a \emph{chain}\index{chain of curves}
of four curves. Considering the Dehn twists $\{T_1,\ldots,T_4\}$ about these
curves, we may apply Theorem~\ref{thm:mcg-raag} to find powers of these twists which generate a right-angled Artin group.
The resulting group is in fact $A(P_4)$; this is a consequence of the fact that the co--intersection graph of $\{T_1,\ldots,T_4\}$ is
isomorphic to the graph whose vertices are $\{\gamma_1,\ldots,\gamma_4\}$ and whose adjacency relation is given by
nontrivial intersection. In other words, the graph $P_4$ is self--complementary.

If $S_{g,n}$ denotes the surface of genus $g$ with $n$ punctures and boundary components (as for the purpose of this discussion,
the difference between punctures and boundary components is immaterial), then there is a natural measure of complexity of
$S_{g,n}$ given by \[c(S_{g,n})=3g-3+n.\] The reader may note that $c(S_{g,n})$ resembles the Euler characteristic, though it is not
exactly the same thing. An elementary exercise in surface topology shows that $c(S_{g,n})$ is exactly the size of a maximal
collection of (isotopy classes of) distinct, pairwise disjoint, essential, nonperipheral, simple closed curves on $S$.

The following is also an easy exercise in combinatorial topology:

\begin{prop}\label{prop:surface-complex}
Let $S=S_{g,n}$. Then $S$ admits a chain of four curves if and only if $c(S)\geq 2$.
\end{prop}

Thus, the moment that $S$ admits two disjoint simple closed curves, it admits a chain of four of them. We can characterize such surfaces
as ones where $g\geq 2$, or $g= 1$ and $n\geq 2$ or $g=0$ and $n\geq 5$. 
It turns out that, if $c(S)\leq 1$ then every right-angled Artin
subgroup of $\Mod(S)$ is (virtually) a product of a free group and an abelian group, 
and is therefore relatively uninteresting from the point of view
of regularity of group actions on one-manifolds.
We have illustrated two different
chains of four curves in Figure~\ref{f:chain-surface}.

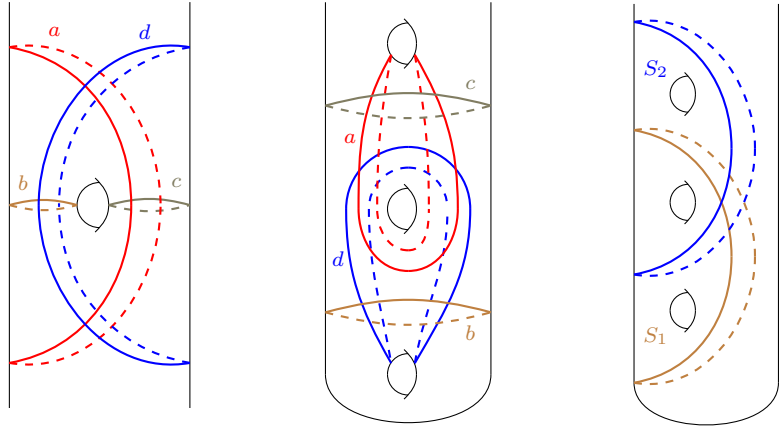
\begin{figure}[htb!]
\centering
{
\begin{tikzpicture}[scale=.6,rotate=90]
\draw (-.5,0) node {};
\draw (0,-2) -- (9,-2);
\draw (0,2) -- (9,2);

\draw (4.5,.5) edge [out=0,in=90]  (5,0) edge [out=180,in=90] (4,0) ;
\draw (4.5,-.2) edge [out=0,in=-135]  (5.1,.1) edge [out=180,in=-45] (3.9,.1);

\draw  (4.5,2) edge [brown,thick,out=-75,in=75] (4.5,.5) edge [brown,thick,dashed,out=-105,in=105](4.5,.5);

\draw  (4.5,-.2) edge [ucl2dkgreen,thick,out=-75,in=75] (4.5,-2) edge [ucl2dkgreen,thick,dashed,out=-105,in=105]  (4.5,-2);

\draw (1,2) edge [red,thick,out=-80, in=180] (4.5,-.7);
\draw [dashed] (1,2) edge [red,thick,out=-100, in=180] (4.5,-1.35);
\draw (8,2) edge [red,thick,out=-100, in=0] (4.5,-.7);
\draw [dashed] (8,2) edge [red,thick,out=-80, in=0] (4.5,-1.35);

\draw [dashed] (1,-2) edge [blue,thick,out=80, in=-180] (4.5,.9);
\draw  (1,-2) edge [blue,thick,out=100, in=-180] (4.5,1.35);
\draw [dashed] (8,-2) edge [blue,thick,out=100, in=0] (4.5,.9);
\draw (8,-2) edge [blue,thick,out=80, in=0] (4.5,1.35);

\node [ucl2dkgreen] at (5,-1.7) {\footnotesize $c$};
\node [brown] at (5,1.7) {\footnotesize $b$};
\node [blue] at (8.35,-1) {\footnotesize $d$};
\node [red] at (8.35,1) {\footnotesize $a$};
\end{tikzpicture}
}
$\qquad\qquad$
{\begin{tikzpicture}[scale=.55,rotate=90]
\draw (2,-2) -- (11,-2);
\draw (2,2) edge [out=180,in=180] (2,-2) edge (11,2);

\draw (2,.5) edge [out=0,in=90]  (2.5,0) edge [out=180,in=90] (1.5,0) ;
\draw (2,-.2) edge [out=0,in=-135]  (2.6,.1) edge [out=180,in=-45] (1.4,.1);
\draw (6,.5) edge [out=0,in=90] (6.5,0) edge [out=180,in=90] (5.5,0);
\draw (6,-.2) edge [out=0,in=-135] (6.6,.1) edge [out=180,in=-45] (5.4,.1);
\draw (10,.5) edge [out=0,in=90] (10.5,0) edge [out=180,in=90] (9.5,0);
\draw (10,-.2) edge [out=0,in=-135] (10.6,.1) edge [out=180,in=-45] (9.4,.1);

\draw (6,1.5)  edge  [blue,thick,out=180, in=30] (2.27,.41)  edge [thick,blue,out=0, in=90]  (7.5,0);
\draw (6,-1.5) edge [thick,blue,out=180, in=-30] (2.27,-.16) edge [thick,blue,out=0, in=-90] (7.5,0);
\draw (6,.95) edge [thick,blue,dashed, out=180, in=10] (2.27,.41);
\draw [thick,blue,dashed, out=-10, in=180] (2.27,-.16) edge (6,-.95);
\draw (7,0) edge [thick,blue,dashed,out=90, in=0] (6,.95) edge [thick,blue,dashed,out=-90, in=0] (6,-.95);

\draw (9.73,.41) edge [thick,red,out=150, in=0] (6,1.2) edge  [thick,red,dashed, out=170, in=0] (6,.75);
\draw  (9.73,-.16) edge  [thick,red,out=-150, in=0] (6,-1.2) edge [thick,red,dashed, out=-170, in=0] (6,-.5);
\draw (4.5,0) edge  [thick,red,out=-90, in=180] (6,-1.2) edge  [thick,red,out=90, in=180] (6,1.2);
\draw (5,0) edge [thick,red,dashed,out=90, in=180] (6,.75) edge [thick,red,dashed,out=-90, in=180] (6,-.5);

\draw  (3.5,2) edge [brown,thick,out=-75,in=75] (3.5,-2) edge [brown,thick,dashed,out=-105,in=105](3.5,-2);
\draw  (8.5,2) edge [ucl2dkgreen,thick,out=-75,in=75] (8.5,-2) edge [ucl2dkgreen,thick,dashed,out=-105,in=105]  (8.5,-2);

\node [brown] at (3,-1.5) {\footnotesize $b$};
\node [blue] at (4.8, 1.7) {\footnotesize $d$};
\node [red] at (7.7,1.4) {\footnotesize $a$};
\node [ucl2dkgreen] at (9,-1.5) {\footnotesize $c$};
\end{tikzpicture}
}
$\qquad\qquad$
{
\begin{tikzpicture}[scale=.48,rotate=90]
\draw (1,-2) -- (11.5,-2);
\draw (1,2) edge [out=180,in=180] (1,-2) edge (11.5,2);

\draw (3,1) edge [out=0,in=90]  (3.5,.5) edge [out=180,in=90] (2.5,.5) ;
\draw (3,.3) edge [out=0,in=-135]  (3.6,.6) edge [out=180,in=-45] (2.4,.6);
\draw (6,1) edge [out=0,in=90] (6.5,.5) edge [out=180,in=90] (5.5,.5);
\draw (6,.3) edge [out=0,in=-135] (6.6,.6) edge [out=180,in=-45] (5.4,.6);
\draw (9,1) edge [out=0,in=90] (9.5,.5) edge [out=180,in=90] (8.5,.5);
\draw (9,.3) edge [out=0,in=-135] (9.6,.6) edge [out=180,in=-45] (8.4,.6);

\draw (1,2) edge [brown,thick,out=-80, in=180] (4.5,-.7);
\draw [dashed] (1,2) edge [brown,thick,out=-100, in=180] (4.5,-1.35);
\draw (8,2) edge [brown,thick,out=-100, in=0] (4.5,-.7);
\draw [dashed] (8,2) edge [brown,thick,out=-80, in=0] (4.5,-1.35);
\node [brown] at (2.3, 1.4) {\footnotesize $S_1$};

\draw (4,2) edge [blue,thick,out=-80, in=180] (7.5,-.7);
\draw [dashed] (4,2) edge [blue,thick,out=-100, in=180] (7.5,-1.35);
\draw (11,2) edge [blue,thick,out=-100, in=0] (7.5,-.7);
\draw [dashed] (11,2) edge [blue,thick,out=-80, in=0] (7.5,-1.35);
\node [blue] at (9.7, 1.4) {\footnotesize $S_2$};
\end{tikzpicture}
}
\caption{Various ways to realize $A(P_4)$ as a subgroup of $\Mod(S)$. The first two pictures illustrate chains of four curves, where in
the first picture they can be nonseparating in $S$, and in the second picture they are all separating in $S$. The third picture illustrates
the first two surfaces in a chain of subsurfaces of $S$.}
\label{f:chain-surface}
\end{figure}

We thus obtain a corollary to Theorem~\ref{thm:mcg-raag}, Proposition~\ref{prop:surface-complex}, and Theorem~\ref{thm:a4-c2},
which completely answers the $C^2$--smoothability question for finite index subgroups of mapping class groups of surfaces:

\begin{thm}\label{thm:mcg-c2}
Let $M\in\{I,S^1\}$. There exists a finite index subgroup $G\le \Mod(S)$ and an embedding $G\longrightarrow \Diffb(M)$ if and only if
$c(S)\leq 1$.
\end{thm}

\subsection{Smooth actions of related groups}
Theorem~\ref{thm:mcg-raag} gives us analogues of Theorem~\ref{thm:mcg-c2} for many other classes of groups other than $\Mod(S)$.
Mapping class groups of surfaces with punctures or boundary components can be identified with groups of automorphisms (or outer
automorphisms) of free groups. So, Theorem~\ref{thm:mcg-c2} carries over verbatim with mapping class groups replaced by $\Aut(F_n)$
and $\Out(F_n)$ for $n\geq 3$.

The mapping class group of a surface $S$ admits a natural linear representation, arising from the action of $\Mod(S)$ on the
first homology
of $S$. Since the homology of $S$ is an abelian invariant, the action of the automorphism group of $\pi_1(S)$ on $H_1(S,\bZ)$
factors through $\Out(\pi_1(S))$. We therefore have a map \[\rho\colon \Mod(S)\longrightarrow \GL_n(\bZ),\] where $n$ is the
rank of $H_1(S,\bZ)$. It is not difficult to see that the image of $\rho$ is infinite; indeed, if $\gamma$ is a simple closed curve whose
homology class $[\gamma]$ is nonzero, then the Dehn twist $T$ about $\gamma$ is easily seen to have an infinite order action on
$H_1(S,\bZ)$. If $S$ is a closed surface of genus $g\geq 1$, then in fact the image of $\rho$ is the full group $\mathrm{Sp}_{2g}(\bZ)$,
the $2g\times 2g$ symplectic group over $\bZ$ (see Chapter 6 of~\cite{FM2012}).
The symplectic form preserved by the action of $\Mod(S)$ is the algebraic intersection
pairing on $H_1(S,\bZ)$.

The kernel of the map $\rho$ is called the \emph{Torelli group}\index{Torelli group} of $S$, and is written $\II(S)$.
Unless one has had some experience dealing with mapping class groups
of surfaces, it may not be obvious that $\rho$ is not injective, and that
consequently $\II(S)$ is nontrivial. Indeed, it is true and not completely trivial that
if $S$ has genus one and at most one puncture, then the map $\rho$ is injective.
The moment that $S$ admits an essential simple closed curve
$\gamma$ that is nonperipheral and separating (i.e.~$S\setminus\gamma$ is disconnected), then $\II(S)$ is nontrivial. Indeed, it is
easily checked that if $\gamma$ is such a curve then the Dehn twist about $\gamma$ acts trivially on $H_1(S,\bZ)$. Some argument
is needed to show that the Dehn twist about $\gamma$ represents a nontrivial mapping class, though this can be established
without too much difficulty by computing
a lift of the Dehn twist to a finite cover of $S$ and showing that the lift acts nontrivially on the homology of the
cover (see~\cite{KoberdaGD,HadariGT20,KobMan15,LiuJAMS20} for instance).

If $S$ admits a chain of separating curves, then $\II(S)$ contains a copy of $A(P_4)$. If $S$ is closed and of genus at least three, or
has genus two and two punctures or boundary components, then such a chain exists. See the middle picture in Figure~\ref{f:chain-surface}.
We thus have the following consequence:

\begin{cor}\label{cor:torelli}
Let $S=S_{g,n}$, and let $G\le \II(S)$ be a finite index subgroup. Then $G$ admits no faithful action by $C^{1+\mathrm{bv}}$ diffeomorphisms
on $M$, for $M\in\{I,S^1\}$, provided that one of the following conditions holds:
\begin{enumerate}[(1)]
\item
We have $g\geq 3$.
\item
We have $g=2$ and $n\geq 2$.
\end{enumerate}
\end{cor}

One can dig even deeper in the mapping class group than the Torelli group and obtain conclusions analogous to Corollary~\ref{cor:torelli}.
For simplicity, assume that $S$ has a fixed marked point $p$ that is preserved by $\Mod(S)$, so that $\Mod(S)$ is naturally a subgroup
of $\Aut(\pi_1(S))$. Write $\gamma_i(\pi_1(S))$ for the $i^{th}$ term of the lower central series, where $i\geq 1$. That is,
$\gamma_1(\pi_1(S))=\pi_1(S)$ and \[\gamma_{i+1}(\pi_1(S))=[\pi_1(S),\gamma_i(\pi_1(S))].\] We write
\[N_i=\pi_1(S)/\gamma_{i+1}(\pi_1(S)),\]
so that $N_i$ is the largest $i$--step nilpotent quotient of $\pi_1(S)$. It is a fact that we have alluded to before that
\[\bigcap_i \gamma_i(\pi_1(S))=\{1\}.\] This last equality is implied by the assertion that $\pi_1(S)$ is residually torsion--free nilpotent.
We will not prove this fact here, and instead direct the reader to~\cite{MKS04} or a combination of~\cite{DK1992a}
and~\cite{CW2004}, or leave it as a good exercise.

Since the subgroup $\gamma_i(\pi_1(S))$ is invariant under all automorphisms of $\pi_1(S)$, we obtain a map
\[\rho_i\colon \Aut(\pi_1(S))\longrightarrow\Aut(N_i).\] Restricting the map $\rho_i$ to $\Mod(S)$, we write $\JJ_i(S)=\ker(\rho_i)\cap\Mod(S)$.
The sequence of nested subgroups $\{\JJ_i(S)\}_{i\ge1}$ is called the \emph{Johnson filtration}\index{Johnson filtration}
of $S$. Note that $\JJ_1(S)=\II(S)$.
Since the intersection of the terms
of the lower central series of $\pi_1(S)$ intersect in the identity, we obtain \[\bigcap_i\JJ_i(S)=\{1\}\] as a formal consequence.

Nontriviality of the terns in the Johnson filtration for $i\geq 2$,
however, is not a formal consequence of the definitions. It is true that, if $S$ has genus
at least two, then each term of $\JJ_i(S)$ is nontrivial, and $\JJ_{i+1}(S)$ has infinite index in $\JJ_i(S)$. If $S_0\longrightarrow S$
is an inclusion of connected surfaces, then (ignoring some basepoint issues) we obtain a map \[\JJ_i(S_0)\longrightarrow\JJ_i(S)\] for
all $i$.

These considerations allow us to find many copies of $A(P_4)$ inside of each term of the Johnson filtration of $S$, provided that
$S$ is large enough. A connected subsurface $S_0\sse S$ is \emph{essential}\index{essential subsurface}
if $S_0$ is not contractible and if the inclusion map is $\pi_1$--injective.
Like in the case of a collection of simple closed curves, a collection
$\{S_1,\ldots,S_k\}$ of pairwise non--isotopic essential subsurfaces of $S$ is a \emph{chain}\index{chain of subsurfaces}
if $S_i\cap S_j\neq\varnothing$ if and only
if $|i-j|\leq 1$, where this intersection is minimized within the respective isotopy classes. See the third picture in Figure~\ref{f:chain-surface}.
The preceding discussion implies easily that if $S$ admits a chain of four essential subsurfaces,
each of which has genus at least two, then $\JJ_i(S)$ contains a copy of $A(P_4)$ for all $i\geq 1$. Such a chain exists, provided that the
genus of $S$ is at least five.

\begin{cor}\label{cor:johnson}
Let $i\geq 1$, let $S=S_{g,n}$, and let $G\le \JJ_i(S)$ be a finite index subgroup.
Then $G$ admits no faithful action by $C^{1+\mathrm{bv}}$ diffeomorphisms
on $M$, for $M\in\{I,S^1\}$, provided that $g\geq 5$.
\end{cor}

\subsection{$C^{1+\tau}$ actions}
To analyze regularities that are weaker than $C^{1+\mathrm{bv}}$, we need to consider more complicated subgroups of $\Mod(S)$.
We have that $\Mod(S)$ contains a copy of $(F_2\times F_2)*\bZ$ provided that:
\begin{enumerate}
\item
The genus of $S$ is at least two.
\item
The genus of $S$ is one and $S$ has at least three punctures or boundary components.
\item
The genus of $S$ is zero and $S$ has at least six punctures or boundary components.
\end{enumerate}

This is easily seen from the fact that these surfaces' mapping class groups contain copies of $F_2\times F_2$ generated by
the squares of two pairs of noncommuting
Dehn twists supported on disjoint surfaces,
say $\{T_1,\ldots,T_4\}$, and by adding another Dehn twist that commutes with none of $\{T_1,\ldots,T_4\}$. The fact that
sufficiently high powers of these Dehn twists then generate a copy of $(F_2\times F_2)*\bZ$ follows from Theorem~\ref{thm:mcg-raag}.

Theorem~\ref{thm:kkr-2020} implies that for all $\tau>0$, the group $(F_2\times F_2)*\bZ$ cannot act on the interval $I$
or $S^1$ by $C^{1+\tau}$
diffeomorphisms. Thus, we obtain:

\begin{cor}\label{cor:mcg-tau}
Let $S=S_{g,n}$, with $g\geq 2$, or $g= 1$ and $n\geq 3$, or $g=0$ and $n\geq 6$, and let $G\le \Mod(S)$ be a finite index subgroup.
Then for $\tau>0$, there is no faithful homomorphism $G\longrightarrow\Diff^{1+\tau}(M)$ for $M\in\{I,S^1\}$.
\end{cor}

In particular, if $S$ is sufficiently complicated and
if $\Mod(S)$ admits a finite index subgroup with a faithful $C^0$ action on $M$, then the critical regularity of that finite index
subgroup of $\Mod(S)$ is exactly one.
When $S$ is closed, or when it has a puncture or marked point and $M=I$, then the content of
Corollary~\ref{cor:mcg-tau} is unclear since we do not know if there
are such subgroups $G$ acting faithfully by homeomorphisms. However, if $S$ has a boundary component then
Theorem~\ref{thm:thurston-handel} does furnish a faithful action of $\Mod(S)$ on $I$ by homeomorphisms, which lends content to
Corollary~\ref{cor:mcg-tau}.

\subsection{$C^1$ actions and Ivanov's Conjecture}
For $C^1$--diffeomorphisms, there is little one can say for finite index subgroups of $\Mod(S)$, given the current state of knowledge.
For the whole mapping class group,
Mann--Wolff proved that for $S=S_{g,1}$, every action of $\Mod(S)$ on $S^1$ 
is semi-conjugate to Nielsen's action~\cite{MWGT20}. Previously, Parwani~\cite{Parwani2008} showed that every $C^1$ action of $\Mod(S)$ on $S^1$ is trivial, 
provided that $g\geq 6$. The basic
idea behind Parwani's proof is the following fact:

\begin{thm}[See~\cite{Parwani2008}, Theorem 1.4]\label{thm:parwani-prod}
Let $G$ and $H$ be finitely generated groups with trivial abelianizations. If \[\phi\colon G\times H\longrightarrow\Diff^1_+(S^1)\]
is a homomorphism, then the restriction of $\phi$ to $\{1\}\times H$ or $G\times\{1\}$ is trivial.
\end{thm}

The basic idea behind the proof of Theorem~\ref{thm:parwani-prod} is to find an invariant measure for $G\times H$, and in the absence
of one, to find an invariant measure for one the factors $\{1\}\times H$ or $G\times\{1\}$. In either case, one finds a global fixed point
for the action of $\{1\}\times H$ or $G\times\{1\}$, which by Thurston's Stability Theorem~\ref{thm:thurston-stab} implies that the
action of one of these two factors is trivial.

If $g\geq 6$ then $\Mod(S)$ contains two disjoint essential subsurfaces $S_1$ and $S_2$, each of which has genus at least $3$.
We immediately have that the inclusion of subsurfaces $S_1$ and $S_2$ into $S$ induces an injection
\[\Mod(S_1)\times \Mod(S_2)\longrightarrow \Mod(S).\] A basic result about mapping class groups shows that if $S=S_{g,n}$ has
genus at least three then $\Mod(S)$ has trivial abelianization (see for example Chapter 5 of~\cite{FM2012}, also~\cite{Harer83}).
Combining with Theorem~\ref{thm:parwani-prod}, we see that
if $S$ has genus at least six then $\Mod(S)$ does not admit a faithful action on $S^1$ by $C^1$ diffeomorphisms. 
To get triviality
of an arbitrary $C^1$ action of $\Mod(S)$, one simply argues that a generating set for $\Mod(S)$ must lie in the kernel of the action.
For closed surfaces, this follows from the fact that a Dehn twist about a nonseparating simple closed curve on $S$ normally generates
the whole mapping class group of $S$. For
surfaces with punctures, some extra argument is needed.

As we have suggested already,
Parwani's result relies on the fact that most mapping class groups have trivial abelianizations. This fact is a consequence of several
facts which fail for finite index subgroups of mapping class groups. The first is that the whole mapping class group is generated by
the conjugacy class of a Dehn twist about a nonseparating simple closed curve, provided that the surface is closed and has positive
genus. The other is that there is a particular relation among seven Dehn twists, known as the \emph{lantern relation}\index{lantern
relation}. In the abelianization
of the whole mapping class group in genus three or more, this relation can be made to read that the image of the cube of a Dehn twist
is equal to its fourth power.

Both of these fact fail for finite index subgroups of $\Mod(S)$,
and it is not known whether if $G\le \Mod(S)$ has finite index then $G$ has finite abelianization. The positive version of this last statement
is known as \emph{Ivanov's Conjecture}\index{Ivanov's Conjecture},
and it follows from Parwani's work that if Ivanov's Conjecture is true then there is no faithful
$C^1$ action by a finite index subgroup of $\Mod(S)$ for $g\geq 6$.


%
%
%
\appendix
\chapter{Concave moduli of continuity}\label{ch:append1}
\begin{abstract}In this appendix, we gather some relevant facts relating to concave moduli of continuity. These topics are included for the convenience of the
reader, and because some of them are difficult to locate in the literature.\end{abstract}

\section{Smooth optimal concave moduli}
In Section~\ref{ss:integrability}, we needed the fact that an arbitrary concave modulus of continuity can be replaced by a smooth (i.e. $C^\infty$) modulus of continuity without changing the functions that are continuous with respect to that modulus. 
We will prove a stronger fact below, which was previously given as Lemma~\ref{lem:modulus-existence}.

\begin{lem}[cf.~\cite{CKK2019}, Proposition 2.7]\label{lem:modulus-existence2}
Let $X$ and $Y$ be metric spaces. If $X$ is a geodesic space and if $f\co X\longrightarrow Y$ is a uniformly continuous function, then there exists a concave modulus\[\beta\co[0,\infty)\longrightarrow[0,\infty)\] which is smooth on $(0,\infty)$ such that the following hold:
\be[(i)]
\item $f$ is $\beta$--continuous;
\item if $\alpha$ is a concave modulus such that $f$ is $\alpha$--continuous,
then  there exist constants $\delta>0$ and $C\geq 1$ such that
 \[ \beta(x)\le C \alpha(x)\] for all $x\in[0,\delta]$.
\ee
\end{lem}
\begin{rem}
It follows from part (ii) that
such an optimal concave modulus is unique up to bi-Lipschitz equivalence.
That is, if $\beta_1$ and $\beta_2$ are such moduli
then we have
\[\frac{1}{C}\leq\frac{\beta_1(x)}{\beta_2(x)}\leq C\] for some $C\ge1$ and for all 
sufficiently small $x>0$.
\end{rem}

A key technical step in the proof is to replace an arbitrary concave modulus by a smooth concave modulus:
\begin{lem}[cf. \cite{Medvedev2001}]\label{lem:medvedev}
If $\alpha$ is an arbitrary concave modulus of continuity, then there exists a concave modulus $\beta$ which is $C^{\infty}$ on $(0,\infty)$ such that 
\[\frac12\alpha(x)\leq \beta(x)\le 2\alpha(x)\] for all $x\in[0,1]$.
Furthermore, we can require that $\alpha(x)\le \beta(x)$ for $x\ge1$.
\end{lem}

This lemma implies that a \emph{locally}\index{locally $\alpha$--continuous} $\alpha$--continuous function
is locally $\beta$--continuous, and vice versa; see Subsection~\ref{ss:smooth-denjoy}.
The choice of the interval $[0,1]$ above is arbitrary, and can be replaced by $[0,c]$  for $c>0$.

Let us establish this lemma first, based on the idea of~\cite{Medvedev2001}.
We remark that in~\cite{Medvedev2001}, Medvedev allows for more general moduli of continuity than concave ones. 
Namely, moduli need only be continuous, subadditive, and non--decreasing.
There is no harm in assuming that all moduli of continuity are in fact concave, as we will see in the proof of 
Lemma~\ref{lem:modulus-existence2}.

Approximation of arbitrary real valued functions by smooth functions can be done via usual techniques from analysis. 
We will use a smooth approximation $\{\phi_\delta\}_{\delta>0}$ to a Dirac mass at the origin.

For this, let us first fix a smooth real valued function
\[
\phi\co \bR\longrightarrow\bR\]
satisfying the following:
\be
\item
the function $\phi$ is positive on an interval $(-1,1)$ and vanishes outside of this interval;
\item
the function $\phi$ satisfies $\phi(-x)=\phi(x)$ for all $x\in\bR$;
\item
we have $\int_{\bR}\phi=1$.
\end{enumerate}
For each $\delta>0$, we then define
\[
\phi_\delta(x):=\frac1\delta\phi\left(\frac{x}{\delta}\right).\]
With this definition, $\phi_\delta$ is an even, smooth bump function supported in $[-\delta,\delta]$ that integrates to one.

Let $f$ be a real--valued continuous function on $(a,b)\sse\bR$.
For each sufficiently small $\delta>0$ we define
\[f_{\delta}(x):=\int_{\bR}\phi_{\delta}(t-x)f(t)\,dt
=\int_{-\delta}^{\delta} \phi_\delta(s)f(s+x)ds.\] 

The Leibniz Rule (Theorem~\ref{thm:leibniz}) implies that $f_{\delta}(x)$ is a smooth function on 
$(a+\delta,b-\delta)$. If $f$ is strictly increasing, then so is $f_\delta$.

Let $[c,d]\sse(a,b)$ be a compact subinterval.
For all $\epsilon>0$ we can find $\delta>0$ such that $|f(x)-f(y)|\le\epsilon$ whenever \[x,y\in[c-\delta,d+\delta]\sse(a,b)\]
satisfies $|x-y|\le\delta$.
If $x\in[c,d]$, 
we have an estimate
\[
|f_\delta(x)-f(x)| \le 
\int_{-\delta}^{\delta} \phi_\delta(s)|f(s+x)-f(x)|ds\le\epsilon.\] 
In particular, $f_\delta$ converges uniformly to $f$ on a compact subset of $(a,b)$
as $\delta\to0$.

Assume that $f$ agrees with (possibly different) affine functions
in $\delta'$--neighborhoods of $c$ and $d$ respectively. 
Then for $\delta>0$ sufficiently small, the function $f_{\delta}$ agrees with $f$ in 
$(\delta'-\delta)$--neighborhoods of $c$ and $d$ as well. This follows easily from the fact that
the function $x\cdot\phi_{\delta}(x)$ is a
compactly supported odd function and therefore integrates to zero, and so that if $f$ is affine then $f_{\delta}=f$.

The concavity of $f$ is also inherited by $f_\delta$.
Namely, suppose $f$ is concave
and $\delta>0$ is sufficiently small.
Let \[a+\delta<x<y<b-\delta,\] and let $u\in[0,1]$.
Setting $z:=ux+(1-u)y$,   we have that
\begin{align*}
f_{\delta}(z)=\int_{-\delta}^\delta \phi_{\delta}(s)f(s+z)\,ds
\geq
\int_{-\delta}^\delta \phi_{\delta}(s)
\left(
u\cdot f(x+s)+(1-u)\cdot f(y+s)\right)\, ds\\=uf_{\delta}(x)+(1-u)f_{\delta}(y),
\end{align*}
whence it follows that $f_{\delta}$ is concave on $(a+\delta,b-\delta)$.

Whereas a concave function may not be $C^1$, 
its left and right derivatives are defined at every point in the interior of the domain. 
Suppose $f(x)$ is a concave function that is strictly increasing on $(a,b)$.
For a point $c\in (a,b)$ let us consider the tangent line defined by the right derivative $f'(c+)$:
\[
g(x):=f'(c+)(x-c)+f(c).\]
This line is above $f(x)$ for all $x\in(a,b)$, and not flat.
We will exploit this observation in order to majorize a given convex function by a
smooth one.

\bp[Proof of Lemma~\ref{lem:medvedev}]
Let us pick a sequence 
\[ a_0=1>a_1>a_2>\cdots\] such that $\lim_{n\to\infty} a_n=0$.
We consider the affine functions
\[
\psi_i(x):=\alpha'(a_i+) (x-a_i)+\alpha(a_i)\]
for $i\ge0$. 
We will define concave moduli $\gamma$ and $\beta$ 
such that 
\[\gamma(x)=\beta(x)=\psi_i(x)\]
in small neighborhoods of $a_i$.

For this purpose, we 
 let
 \[\{\delta_i\},\{\delta_i'\},\{\epsilon_i\}>0\] be positive decreasing sequences which will be determined later. We set 
\[
\gamma(x)=\beta(x)=\psi_0(x)\]
for all $x\ge a_0$. We also set $\gamma(0)=\beta(0)=0$.

We will define $\gamma$ and $\beta$ on each interval $[a_i,a_{i-1}]$ piecewise.
For each $x$ in that interval, we set
\[
\gamma(x):=\min\{ \psi_i(x),\psi_{i-1}(x),\alpha(x)+\epsilon_i\}.\]
We will choose a sufficiently small $\delta_i'>0$ depending on $\epsilon_i$;
more precisely,
for all 
\[x\in [a_i-\delta_i',a_i+\delta_i']\sse(a_{i+1},a_{i-1})\] we  require
\[
\alpha(x)\le\psi_i(x)\le \alpha(x)+\epsilon_{i+1}<\alpha(x)+\epsilon_i.\]
As it is defined as the minimum of strictly increasing continuous concave functions on each interval $[a_i,a_{i-1}]$, we see that
the combined function \[\gamma\co[0,\infty)\longrightarrow[0,\infty)\] is a concave modulus satisfying the following.
\begin{itemize}
\item
$\gamma=\psi_i$ on each interval $[a_i-\delta_i',a_i+\delta_i']$;
\item $0\le \gamma(x)-\alpha(x)\le\epsilon_i$ for $x\in [a_i,a_{i-1}]$.
\end{itemize}

Let $i\ge1$.
For each $x\in[a_i,a_{i-1}]$ we define
\[
\beta(x):=\gamma_{\delta_i}(x)=\int_{-\delta_i}^{\delta_i} \phi_{\delta_i}(s) \gamma(s+x)\,ds.\]
By requiring that $0<\delta_i<\delta_i'$, we have that
$
\beta=\psi_i$
on each $(\delta_i'-\delta_i)$--neighborhood of $a_i$.
As a strictly increasing surjective concave function, 
the map \[\beta\co[0,\infty)\longrightarrow[0,\infty)\] is a concave modulus
that is smooth on $(0,\infty)$.

By further decreasing $\delta_i>0$, we can also require that
whenever 
\[x,y\in[a_i-\delta_i',a_{i-1}+\delta_{i-1}']\] satisfies $|x-y|\le \delta_i$,
we have 
\[|\gamma(x)-\gamma(y)|\le \epsilon_i.\]
Then each $x\in[a_i,a_{i-1}]$ satisfies that
\[
|\beta(x)-\gamma(x)|=|\gamma_{\delta_i}(x)-\gamma(x)|\le \epsilon_i\]
as we have seen above.

Setting $\epsilon_i:=\frac12\alpha(a_i)$, we have for each $x\in[a_i,a_{i-1}]$ that
\[
\frac12\alpha(x)\le\alpha(x)-\epsilon_i\le \gamma(x)-\epsilon_i\le \beta(x)\le\gamma(x)+\epsilon_i\le \alpha(x)+2\epsilon_i\le 2\alpha(x).\]
Since for $x\ge 1$ we have
\[
\beta(x)=\psi_0(x)\ge\alpha(x),\]
 the proof is complete.
\ep

Let us now establish the main result of this section.
\begin{proof}[Proof of Lemma~\ref{lem:modulus-existence2}]
For brevity, we will often write $|x-y|$ to mean $d_X(x,y)$ or $d_Y(x,y)$.
Consider the modulus of continuity function for $f$ defined as
\[\mu(t):=\sup_{|x-y|\leq t} |f(x)-f(y)|.\]  
We leave it to the reader to check that 
$
\mu\co [0,\infty)\longrightarrow[0,\infty)$ is nondecreasing and continuous.

Let $s,t>0$.
Since $X$ is a geodesic space, 
for each pair $x,y\in X$ satisfying \[|x-y|\le s+t\] we can find $z\in X$ such that
$|x-z|\le s$ and $|y-z|\le t$. It follows that $\mu$ is \emph{subadditive}\index{subadditive}:
\[
\mu(s+t)\le\mu(s)+\mu(t).\]
Writing $t=sm+r$ for $m\in\bZ_{\ge0}$ and $r\in[0,s)$, we also note that
\[
\mu(t)\le\mu((m+1)s)\le (m+1)\mu(s) \le \left(1+\frac{t}s\right)\mu(s)=:\mu_s(t).\]

To produce the modulus $\beta$, we consider the family $\MM$ of continuous, nondecreasing, concave functions 
\[\nu\co[0,\infty)\longrightarrow[0,\infty)\] such that
$\mu(t)\leq\nu(t)$ for all $t>0$. 
The family $\MM$ is nonempty since it contains $\{\mu_s\}_{s>0}$.

We define \[\beta_0(t)=\inf\{\nu(t)\mid {\nu\in\MM}\}.\]
It is immediate that $\beta_0(0)=0$
and that $\beta_0$ is concave. For each $t>0$ we have
\[\mu(t)\le \beta_0(t)\le \mu_t(t)=2\mu(t).\]
It follows that $\beta_0(t)$ is continuous at $t=0$, and by concavity, at arbitrary $t>0$. It follows that
\[\beta_1(t)=\beta_0(t)+t\]
is a concave modulus. Feeding $\beta_1$ to Lemma~\ref{lem:medvedev}, 
we obtain a concave modulus $\beta$ which is smooth on $(0,\infty)$ 
such that
\[
\frac12\mu(t)\le \frac12\beta_1(t)\le\beta(t)\le 2\beta_1(t)\le 4\mu(t)+2t\]
for all $t\in[0,1]$.
Moreover, for $t\ge1$ we can require that 
\[
\mu(t)\le\beta_1(t)\le\beta(t).\]

Part (i) of the conclusion follows from that
\[
[f]_\beta
=\sup_{|x-y|>0} \frac{|f(x)-f(y)|}{\beta(|x-y|)}
\le
\sup_{|x-y|>0} \frac{\mu(|x-y|)}{\beta(|x-y|)}\le2.\]

For part (ii), assume that $f\in C^\alpha(X,Y)$ for some concave modulus $\alpha$.
We have that $\mu(t)\le [f]_\alpha \alpha(t)$ for all $t\ge0$.
Recall also that $t/\alpha(t)$ is monotone increasing in $t$ by concavity.
If $t\in[0,1]$ we have that 
\[\beta(t)\le 4\mu(t)+2t \le \left( 4[f]_\alpha+\frac2{\alpha(1)}\right)\alpha(t).\]
This completes the proof.

\end{proof}

\section{The group structure on $\Diff^{k,\alpha}(M)$}
Let $M$ be a smooth connected manifold of dimension $n\ge1$.
For an integer $k\ge0$ and a concave modulus $\alpha$, we 
denote by $\Diff^{k,\alpha}(M)$ the set of $C^k$ diffeomorphisms $f\co M\longrightarrow M$ 
such that its $k$--th derivative $D^kf$ is locally $\alpha$--continuous.
Here, we locally regard $D^kf$ as a multilinear map from $(\bR^n)^k$ to $\bR^n$.
We say such a diffeomorphism $f$ is $C^{k,\alpha}$. The case $\alpha(x)=x$ is referred to as $C^{k,\mathrm{Lip}}$.

When $\dim M=1$, we also allow $\alpha=\mathrm{bv}$.
That is,  $\Diff^{k,\mathrm{bv}}(M)$ denotes the set 
of $C^k$ diffeomorphisms $f\co M\longrightarrow M$ such that every point in $M$ admits an open neighborhood 
$U$ in which the variation $\Var\left(\der{f}{k};U\right)$  is finite;
such diffeomorphisms are called $C^{k,\mathrm{bv}}$.

In this section, we record a proof that all of the sets defined above are closed under composition and thus have the structure of a group,
provided that $k\ge1$. 
As remarked in the introduction, the same conclusion does not hold for $k=0$;
for instance, the map $f(x)=\sqrt{x}$ lies in  
the set $\Diff^{1/2}[0,1]$ of $1/2$--H\"older continuous diffeormophisms, though $f^2\not\in \Diff^{1/2}[0,1]$.
The concavity of $\alpha$ is also a reasonable assumption, by Lemma~\ref{lem:medvedev}.

The discussion here is already contained in~\cite{Mather1} for the $C^{k,\alpha}$ case, and in the 
appendix of~\cite{KK2020crit} for the $C^{k,\mathrm{bv}}$ cases.
We summarize those proofs here in the interest of making the book more self--contained.

\begin{prop}
Let $M$ be a smooth connected manifold, and let $\alpha$ be a concave modulus of continuity. Then for integers 
$k\geq 1$, the sets $\Diff^{k}(M)$ and  $\Diff^{k,\alpha}(M)$ are closed under composition and inversion. 
The same conclusion holds for $\Diff^{k,\mathrm{bv}}(M)$ when $\dim M=1$. In particular, these sets are groups.
\end{prop}
\begin{proof}
Let us first consider the case of $C^{k,\alpha}$ diffeomorphisms.
We have the following observation~\cite{Mather1}.

\begin{claim}
Let $X,X',Y,Z,W$ are open subsets of Euclidean spaces. The following conclusions hold.
\be[(1)]
\item
If $F\co X\longrightarrow Y$ is $C^\alpha$,
and if $G\co X'\longrightarrow Z$ is $C^{\mathrm{Lip}}$, 
then $F\circ G$ and $G\circ F$ are $C^\alpha$, whenever these compositions are defined.
\item
If $F\co X\longrightarrow Y$
and  $G\co X\longrightarrow Z$ is $C^\alpha$,
and if $B\co Y\times Z\longrightarrow W$ is a bilinear map,
then the map
\[
B(F,G)\co X\longrightarrow W\]
is $C^\alpha$.
\ee
\end{claim}
For part (1), it suffices to use the following local estimates.
\begin{align*}
|F\circ G(x)-F\circ G(y)|
&\le [F]_\alpha \cdot \alpha(|G(x)-G(y)|)
\le [F]_\alpha \left(1+[G]_{\mathrm{Lip}}\right)\alpha( |x-y|),\\
|G\circ F(x)-G\circ F(y)|
&\le [G]_{\mathrm{Lip}}\cdot |F(x)-F(y)|
\le [G]_{\mathrm{Lip}}\cdot [F]_\alpha \cdot \alpha( |x-y|).
\end{align*}
We note that the concavity of $\alpha$ was crucially used here.
Part (2) is an obvious consequence of part (1),
combined with the fact that $B$ is locally Lipschitz. This proves the claim.

Let us now assume $f,g\in \Diff^{k,\alpha}(M)$.
We have
\begin{align*}
D(f\circ g) =& (Df\circ g) Dg,\\
D^2(f\circ g)=& (D^2f\circ g) (Dg\times Dg)+(Df\circ g) D^2g,\\
D^k(f\circ g) =& (D^k f\circ g)\prod_{t=1}^k Dg  + (Df\circ g)D^kg
 \\
&+\sum_{\gamma=(j_1,\ldots,j_i)} C_\gamma (D^i f\circ g) \prod_{t=1}^i D^{j_t}g.
\end{align*}
for suitable nonnegative integers $C_\gamma$.
In the last formula,  $k$ is assumed to be at least three, and the sum is taken over
\[
\gamma=(j_1,\ldots,j_i)\]
for some $i\in[2,k-1]$ and $j_1,\ldots,j_i\in[1,k-1]$ such that $\sum_t j_t = k$.
This sum is sometimes called the Fa\`a di Bruno formula~\cite{Fraenkel1978}.
We inductively deduce that $D^k(f\circ g)$ is $C^\alpha$ from the above claim, since it is a sum of multiplications of
matrices with $C^\alpha$ entries.

Let $\Inv\co \GL_n(\bR)\longrightarrow \GL_n(\bR)$ be the inversion map.
By the inverse function theorem, we have that

\[
D(f^{-1})= \Inv\circ Df\circ f^{-1}.\]
Since $\GL_n(\bR)$ is analytic, since $Df$ is $C^{k-1,\alpha}$,
and since $f^{-1}$ is $C^k$
we again see from the above claim that $D(f^{-1})$ is $C^{k-1,\alpha}$.
This completes the proof that $\Diff^{k,\alpha}(M)$ is a group.

The set $\Diff^k(M)$ is also a group since 
\[
\Diff^k(M) = \bigcap \Diff^{k,\alpha}(M),\]
where the intersection is taken over all smooth concave moduli $\alpha$; see Lemma~\ref{lem:medvedev}.

Let us finally consider the case $\dim M=1$ and $\alpha=\{\mathrm{bv}\}$.
We note the following.

\begin{claim}
Let $k\ge1$ be an integer.
\be[(1)]
\item 
If $F,G\co M\longrightarrow \bR$ are $C^{k-1,\mathrm{bv}}$  
then so is $F\cdot G\co M\longrightarrow\bR$.
\item
If $F\co M\longrightarrow\bR$ is $C^{k-1,\mathrm{bv}}$ and
if $G\co  M\longrightarrow M$ is $C^k$,
then $F\circ G\co M\longrightarrow \bR$ is $C^{k-1,\mathrm{bv}}$.
\item
Let $U\sse\bR$ be an interval.
If $F\co U\longrightarrow\bR$ is $C^{k}$ and
if $G\co  M\longrightarrow U$ is $C^{k-1,\mathrm{bv}}$,
then $F\circ G\co M\longrightarrow \bR$ is $C^{k-1,\mathrm{bv}}$.
\ee
\end{claim}
To see part (1), we first let $k=1$. 
For $x$ and $y$ sufficiently near in $M$, we have
\[
|F(x)G(x)-F(y)G(y)| \le \sup_U (|F| + |G|) ( |F(x)-F(y)|+|G(x)-G(y)|),\]
where the suprema are taken on a suitable relatively compact neighborhood $U$.
We deduce that
\[
\Var(F\cdot G; U) \le \sup_U (|F| + |G|) (\Var(F;U)+\Var(G;U)).\]
It follows that $F\cdot G$ is $C^{\mathrm{bv}}$.

In the case $k\ge2$, we note that
\[ (F\cdot G)'= F'\cdot G + F\cdot G'.\]
By induction, we see that the right hand side is $C^{k-2,\mathrm{bv}}$,
proving part (1);
here, we used the fact that a $C^{k-1}$ map is $C^{k-2,\mathrm{bv}}$,
which is true because we defined the class $C^{\mathrm{bv}}$ to consist
of functions locally having bounded variations. Parts (2) and (3) follow from an easy induction argument and the chain rule 
\[ (F\circ G)'=(F'\circ G)\cdot G',\]
along with an application of part (1).

Let us now assume $f,g\in \Diff^{k,\mathrm{bv}}(M)$.
Since
\[
(f\circ g)'=(f'\circ g)\cdot g',\]
we can apply the above claim to the right hand side and see that 
$(f\circ g)'$ is $C^{k-1,\mathrm{bv}}$. This implies that $f\circ g$ is $C^{k,\mathrm{bv}}$.
For the inverse, we again write
\[
\left(f^{-1}\right)' =\Inv\circ f'\circ f^{-1},\]
with the map $\Inv(x)=1/x$.
By the inverse function theorem, we have $f^{-1}$ is $C^k$.
Part (3) of the above claim implies that
\[
\Inv\circ f'\co M\longrightarrow\bR\]
is $C^{k-1,\mathrm{bv}}$. We conclude that 
$f^{-1}$ is $C^{k,\mathrm{bv}}$
by 
applying part (2) of the same claim.\ep

\section{The Muller--Tsuboi trick for non--integer regularities}
Recall our notation that $\Diff_J^{k,\alpha}(\bR)$ denote the set of $C^{k,\alpha}$ diffeomorphisms of $\bR$ which are the identity outside $J$.
We have a natural inclusion
\[\Diff_{[0,1]}^{\infty}(\bR)\hookrightarrow\Diff_+^{\infty}[0,1].\]
This inclusion is strict, as the elements of $\Diff_+^{\infty}[0,1]$ are not required to be $C^k$--tangent to the identity. 

Let us fix a $C^\infty$--homeomorphism \[\varphi\co [0,1]\longrightarrow[0,1]\] such that $\varphi(x)=e^{-1/x}$ near $x=0$ and such that 
\[\varphi(1-x)=1-\varphi(x)\] for all $x\in[0,1]$.  
Muller~\cite{Muller} and Tsuboi~\cite{Tsuboi1984Asterisque} observed that one can reverse the direction of the above embedding
after applying a topological conjugacy by a $C^\infty$ homeomorphism $\varphi$.
We will generalize their result to cover regularities of the form $(k,\alpha)$.

\begin{thm}[{Muller--Tsuboi trick}]\label{t:muller-tsuboi}
If $k\ge1$ be an integer, and let $\alpha$ be a concave modulus. 
Then we have the following inclusion relations.
\be[(1)]
\item $\varphi^{-2}\Diff_+^{k,\alpha}[0,1]\varphi^{2}\le\Diff_{[0,1]}^{k,\alpha}(\bR)$;
\item $\varphi^{-2}\Diff_+^{k}[0,1]\varphi^{2}\le\Diff_{[0,1]}^{k}(\bR)$;
\item $\varphi^{-2}\Diff_+^\infty[0,1]\varphi^{2}\le\Diff_{[0,1]}^\infty(\bR)$.
\ee
\end{thm}

We give a proof following \cite{Tsuboi1984} and \cite[Appendix A.2]{KK2020crit}. It will be convenient for us to say a function $f(x)$ is \emph{defined for $x\ge0$} when $f(x)$ is defined on an interval $[0,A]$ for some $A>0$.

\begin{lem}\label{l:muller-tsuboi}
Let $k\in{\bZ_{>0}}$, and let $\alpha$ be a concave modulus.
Suppose $f$ and $g$ are $C^{k,\alpha}$ maps defined for $x\ge0$.
\be[(1)]
\item
If $f$ satisfies
\[
f(0)=f'(0)=\cdots={f}^{(k)}(0)=0,\]
then $f(x)/x$ extends to a $C^{k-1,\alpha}$ map defined for $x\ge0$.
\item
If $g$ satisfies
\[g(0)=0,\quad g'(0)>0,\]
then $g_1:=\varphi^{-1}\circ g\circ\varphi$ extends to a $C^{k,\alpha}$ map defined for $x\ge0$
such that
\[ g_1(0)=0,\quad g_1'(0)=1.\]
\item
In part 2, the map
$
g_2:=\varphi^{-2}\circ g\circ\varphi^2$ extends to a $C^{k,\alpha}$ map defined for $x\ge0$
such that
\[ g_2^{(i)}(0)=\Id^{(i)}(0)\]
for $i=0,1,\ldots,k$.
\ee
\end{lem}

\bp
(1) By an iterated application of the Fundamental Theorem of Calculus, we have that
\[
f(x)=\int_{x_1=0}^{x}\int_{x_2=0}^{x_1}\cdots\int_{x_k=0}^{x_{k-1}} f^{(k)}(x_k)\,dx_k\cdots dx_1.\]

\begin{claim}
The map $f(x)/x^k$ extends to a $C^\alpha$--map defined for $x\ge0$.
\end{claim}

We prove the claim after declaring the value of $f(x)/x^k$ to be zero at $x=0$.
We first note that if $x>0$ then 
\[\abs*{f^{(k)}(x)}=\abs*{f^{(k)}(x)-f^{(k)}(0)}\le \left[ f^{(k)}\right]_\alpha \alpha(x).\]
To see that $f(x)/x^k$ is $C^\alpha$ at $x=0$, it suffices to check
 for $x>0$ that
\[\abs*{\frac{f(x)}{x^k}}
\le\frac1{x^k}\int_0^{x}\int_0^{x_1}\cdots\int_0^{x_{k-1}} \left[{f^{(k)}}\right]_\alpha \alpha(x_k)\, dx_k\cdots dx_1
\le\left[{f^{(k)}}\right]_\alpha\alpha(x).\]

To verify the $C^\alpha$ condition at $x>0$ we note for small $h>0$ that
\begin{align*}
\abs*{\frac{f(x+h)}{(x+h)^k}-\frac{f(x)}{x^k}}
&\le\left(
\frac1{x^k}-\frac1{(x+h)^k}\right)
\int_0^{x}\int_0^{x_1}\cdots\int_0^{x_{k-1}} \left[{f^{(k)}}\right]_\alpha \alpha(x_k) \,dx_k\cdots dx_1
\\&+
\frac1{(x+h)^k}
\int_x^{x+h}\int_0^{x_1}\cdots\int_0^{x_{k-1}} \left[{f^{(k)}}\right]_\alpha  \alpha(x_k) \,dx_k\cdots dx_1
\end{align*}
From the general inequality \[(1+t)^{-k}\ge 1-kt\] for $t\ge0$, 
we can bound first term from above by
\[
\left(1-
\frac{x^k}{(x+h)^k}\right)
\left[{f^{(k)}}\right]_\alpha  \alpha(x)
\le
 \frac{kh}{x}\left[{f^{(k)}}\right]_\alpha  \alpha(x).
\]
The second term is bounded from above by
\[
\frac{h(x+h)^{k-1}}{(x+h)^k}
\left[{f^{(k)}}\right]_\alpha  \alpha(x+h).
\]
By choosing $0<h\le x$ and using the fact that $x/\alpha(x)$ is monotone increasing, we have that
\[
\abs*{\frac{f(x+h)}{(x+h)^k}-\frac{f(x)}{x^k}}
\le
 (k+1)\left[{f^{(k)}}\right]_\alpha \alpha(h).\]
This proves the claim.

To complete the proof of part (1), we use the product rule to see that
\[
\left(\frac{f(x)}{x}\right)^{(k-1)}=\sum_{i=1}^{k} a_i \frac{f^{(k-i)}(x)}{x^i}\]
for some integers $a_i$. Using the above claim for the functions $f^{(k-i)}$, where $i\in\{1,\ldots,k\}$, 
we see that each term $f^{(k-i)}(x)/x^i$ is $C^\alpha$ for $x\ge0$.
It follows that $f(x)/x$ is $C^{k-1,\alpha}$, as asserted. 

(2)
It is clear that $g_1$ is $C^{k,\alpha}$ for $x>0$.
Put $h(x)=g(x)/x$. Applying part (1) to the function 
\[ g(x) - (\text{the }k\text{--th Taylor polynomial of }g)=g(x) - g'(0)x-g''(0)x^2/2-\cdots\]
we see that $h$ is $C^{k-1,\alpha}$ map for $x\ge0$.
Using $\varphi^{-1}(y)=-1/\log y$, we have
\[
g_1(x)=\frac{-1}{\log(g\circ\varphi)}=\frac{x}{1-x\log h\circ\varphi}.\]
As $x\to+0$, the denominator of the last expression is bounded away from zero
since $h(0)=g'(0)>0$.  Since $h$ is $C^{k-1,\alpha}$ for $x\ge0$, so is $g_1$.

It only remains to verify that $g_1'$ is $C^{k-1,\alpha}$ at $x=0$.
It is trivial that $g_1'(0)=1$.
For $x>0$, we have
\[
g_1'(x) = \frac{1+x^2 (h'\circ\varphi)\varphi'/(h\circ\varphi)}{\left(1-x\log(h\circ\varphi)\right)^2}
= \frac{1+ \varphi\cdot (h'\circ\varphi)}{\left(1-x\log(h\circ\varphi)\right)^2}
= \frac{g'\circ\varphi}{(h\circ\varphi)\left(1-x\log(h\circ\varphi)\right)^2}
.\]
Since 
the denominator and the numerator extend to $C^{k-1,\alpha}$ maps for $x\ge0$
maintaining that the former is positive,
so does $g_1'$.
Since \[\lim_{x\to+0} g_1'(x)=1=g_1'(0),\] such an extension coincides with $g_1'$ itself.
We conclude that $g_1'$ is also $C^{k-1,\alpha}$ for $x\ge0$.

(3) By part (2), we know that $g_2$ extends to a $C^{k,\alpha}$ map defined for $x\ge0$ such that
\[g_2(0)=0,\quad g_2'(0)=1.\]

We claim that 
\[
\lim_{x\to+0} \frac{g_2(x)-x}{\varphi(x)}=0.\]
Setting \[y:=\varphi^2(x),\quad u=\log(1/y),\quad v=\log(1/g(y)),\] we have that
\[
\frac{g_2(x)-x}{\varphi}= \frac{\varphi^{-2}\circ g(y)-\varphi^{-2}(y)}{\varphi^{-1}(y)}
=u\left(\frac1{\log v} - \frac1{\log u}\right)
=\frac{u}{(\log u)^2} \left(\log\frac{u}v\right)\left( \frac{\log u}{\log v}\right)
.\]
Note that as $x\to+0$, we have $y\to+0$ and $u,v\to+\infty$.
Using L'H\^opital's rule, we have 
\[
\lim_{y\to+0} \frac{v}{u}=
\lim_{y\to+0} \frac{\log g}{\log y}
=\lim_{y\to+0} \frac{yg'}{g}=1.\]
So, we have
\[\lim_{y\to+0} \frac{\log u}{\log v}
=\lim_{y\to+0} \frac{\log(-\log y)}{\log(-\log g )}
=\lim_{y\to+0} \frac{g\log g}{yg' \log y}=1.\]
We also see that
\[\lim_{y\to+0} (u-v) =\lim_{y\to+0} \log(g/y)=\log g'(0).\]
Using \[\lim_{t\to1} \frac{t-1}{\log t} = 1,\] we deduce that
\[
\lim_{x\to+0}
\frac{g_2(x)-x}{\varphi}
=
\lim_{y\to+0}
\frac{u}{(\log u)^2} \left(\frac{u}v-1\right)
=
\lim_{y\to+0}
\frac{u-v}{(\log u)^2} =0,
\]
and so the claim is proved.

Now, we have that
\[ |g_2(x)-x| = o(\varphi(x)) = o(x^i)\]
for all $i\ge0$. By L'H\^opital's rule again, we see that $g_2$ is $C^k$--tangent to the identity at $x=0$.

\ep

\bp[Proof of Theorem~\ref{t:muller-tsuboi}]
(1) Let $g\in \Diff_+^{k,\alpha}[0,1]$. By Lemma~\ref{l:muller-tsuboi}, 
we see that $g_2:=\varphi^{-2}\circ g\circ\varphi^2$ is a $C^{k,\alpha}$ diffeomorphism near $x=0$ and $C^k$--tangent to the identity at $x=0$.
By symmetry, the same is true near $x=1$ and we obtain $g_2\in \Diff_{[0,1]}^{k,\alpha}(\bR)$.

(2) Let $g\in \Diff_+^{k}[0,1]$. We can choose a concave modulus $\alpha$ such that $g^{(k)}$ is $C^\alpha$, as given in Lemma~\ref{lem:medvedev}.

Now, we see from part 1 that
\[
\varphi^{-2}\circ g\circ\varphi^2 \in \Diff_{[0,1]}^{k,\alpha}(\bR)\le\Diff_{[0,1]}^{k}(\bR).\]

Part (3) is obvious from the fact that \[\Diff_+^\infty[0,1]=\bigcap_{k\ge1} \Diff_+^k[0,1],\] whereby the proof is complete.\ep


%
%
%
\chapter{Orderability and H\"older's Theorem}\label{ch:append2}

In this appendix we include some classical facts about groups acting on one-manifolds which are implicit in or even essential to the discussion 
above, but which do not fit cleanly into the narrative of the book. We record these facts for the convenience of the reader, in the interest
of a self-contained discussion. The reader will find here just a basic account of orderability of groups;
much more detailed accounts can be found
in~\cite{BMR77,CR2016,DDRW08,DNR2014,Navas2011}, for instance.

\section{Orderability of groups and homeomorphism groups}

A basic fact about $\Homeo_+(M)$, where $M\in \{I,S^1\}$, is that countable subgroups of the full group of homeomorphisms can be
characterized in a purely algebraic manner. Here, we establish these characterizations.

\subsection{Linear orderability and the interval}
Let $G$ be a group. An ordering on $G$ is simply a total order on the elements of $G$. An ordering $<$ is called
\emph{left invariant}\index{left invariant ordering}
if it is compatible with the action of $G$ on itself by multiplication. That is, for $\{a,b,c\}\sse  G$, we have \[a<b \quad \textrm{if and only if}\quad
c\cdot a<c\cdot b.\] A \emph{right invariant}\index{right invariant ordering} 
ordering on $G$ is defined similarly, and an ordering is \emph{bi-invariant}\index{bi-invariant ordering}
if it is both left and right invariant. It is straightforward
to check that an ordering on $G$ is bi-invariant if and only if it is left invariant and conjugation invariant. A group is called
\emph{orderable}\index{linear ordering} if it admits a left invariant or right invariant ordering.

A left invariant ordering $<$ on $G$ has associated to it a \emph{positive cone}\index{positive cone},
which is defined to be \[\PP_<=\{g\in G\mid 1<g\}.\] The
subscript $<$ may be omitted if the ordering relation is clear from context. The reader may check that $\PP$ is a semigroup, and that
\[G\setminus\{1\}=\PP\cup\PP^{-1}\quad \textrm{and}\quad \PP\cap\PP^{-1}=\varnothing,\]
where $\PP^{-1}$ consists of inverses of elements of $\PP$. Conversely, given a semigroup $\PP\sse  G$
such that $\PP\cup\PP^{-1}$ coincides with the non-identity elements of $G$ and such that $\PP\cap\PP^{-1}=\varnothing$, then $\PP$ defines
a left invariant ordering on $G$ by \[a<b\quad \textrm{if and only if}\quad  a^{-1}b\in\PP.\]

The view of orderability via positive cones is useful from several perspectives. For one, it does not privilege left invariant orderings over
right invariant ones. Indeed, given a positive cone $\PP$, we may define a right invariant ordering on $G$ by setting \[a<b\quad
\textrm{if and only if}\quad ab^{-1}\in\PP^{-1}.\] We will concentrate on left invariant orderings as opposed to right invariant ones,
since this perspective
is compatible with group
actions on the left. A positive cone arises from a bi--invariant ordering if and only if it is also conjugation invariant, which is to say
for all $g\in G$, we have $g^{-1}\PP g=\PP$.

Positive cones also allow one to identify a left invariant ordering on a group with a subset of $G$. Since subsets of
$G$ are identified with $2^G$, and since the latter of which is topologized by giving $2$ and $G$ the discrete topology and the $2^G$ the
product topology, we obtain a compact topological space in which left invariant orderings of $G$ live. The space of all left invariant orderings
is of interest in its own right, being a closed subset of a Cantor set equipped with a natural action of $G$ by conjugation. We will not comment
on this last topic any further here, contenting ourselves with directing the interested reader to some 
references~\cite{CZ2017,DDRW08,KobNYJM11,Morris12NYJM,NavasAIF10,NWblms11,Sikora04}.

Unpacking the definition of a left invariant ordering on a group $G$, it is not surprising that the existence of such an ordering is closely related
to actions of $G$ on $I$ (or $\bR$) by homeomorphisms. Indeed, the very definition of the left invariant ordering says that one can
place the elements of $G$ on the real line (if $G$ is countable, at least), in such a way that multiplication on the left preserves the order
in which they were arranged.

\begin{prop}\label{prop:l-order-homeo}
Let $G$ be a group. 
\be[(1)]
\item
If there exists a faithful homomorphism \[G\longrightarrow \Homeo_+(\bR),\]
then $G$ admits a left
invariant ordering.
\item
Conversely, if $G$ is countable and admits a left invariant ordering,
then 
 there exists a faithful homomorphism \[G\longrightarrow \Homeo_+(\bR).\]
\ee
\end{prop}
\begin{proof}
For part (1)
suppose that $G\le \Homeo_+(\bR)$. Let $\{a_n\}_{n\in\bN}$ be an enumeration of $\bQ$ (though an arbitrary countable dense subset would suffice).
For $g,h\in G$ distinct elements, we set \[g<h\quad \textrm{if}\quad g(a_n)<h(a_n),\]
where here $n$ is the first index for which $g$ and $h$ disagree. Since $\{a_n\}_{n\in\bN}$ is dense in $\bR$, we obtain a total ordering
on $G$. If $f\in\Homeo_+(\bR)$ is arbitrary, we have \[g(a_n)\neq h(a_n)\quad \textrm{if and only if}\quad f(g(a_n))\neq f(h(a_n)),\] and
so this ordering is also left invariant.

For part 2,  suppose that $G$ admits a left invariant ordering $<$, and let $\{g_n\}_{n\ge0}$ be an enumeration of the elements of $G$.
We first build an order preserving function $\phi\colon G\longrightarrow \bR$;
this process is the standard  method of proof
 that every countable totally ordered set 
is order--isomorphic to a subset of $\bQ$.
We begin by sending $\phi\colon g_0\mapsto 0$. If
$\phi$ has been defined on $\{g_0,\ldots,g_n\}$, then we consider $g_{n+1}$. If $g_i<g_{n+1}$ for $i\leq n$ then we define
\[\phi(g_{n+1})=\sup_{0\leq i\leq n} \phi(g_i)+1.\] Similarly, if $g_{n+1}<g_i$ for all $i\leq n$ then we define
\[\phi(g_{n+1})=\inf_{0\leq i\leq n} \phi(g_i)-1.\] Otherwise, there are indices $0\leq i\neq j\leq n$ such that $g_i<g_{n+1}<g_j$ and such
that no $g_k$ lies between $g_{n+1}$ and $g_i$ nor $g_j$ for $0\leq k\leq n$. We then set \[\phi(g_{n+1})=\frac{\phi(g_i)+\phi(g_j)}{2}.\]
This recursively defines an order preserving function from $G$ to $\bR$.

Let $C$ be the closure of $\phi(G)$. For all $g\in G$, we have a map \[\mu_g\colon \phi(G)\longrightarrow\phi(G)\] given by 
$\mu_g(\phi(g_i))=\phi(gg_i)$. 
We note the following.
\be[(1)]
\item 
We have $\sup\phi(G)=\infty$ and $\inf\phi(G)=-\infty$;
\item
For each interval $J$ in the connected component of $\bR\setminus C$,
we have $\partial J\sse\phi(G)$.
\ee
From this 
it is not difficult to deduce that for all $g\in G$, the map $\mu_g$ extends to a map $C\longrightarrow C$, which we
will also call $\mu_g$. Furthermore, $\mu_g$ is in fact a homeomorphism of $C$ preserving the natural order on points inherited from $\bR$.
The association $g\mapsto\mu_g$ is easily seen to be a homomorphism of groups.

Since $C$ is closed, we have that $\bR\setminus C$ (if not empty) is a countable disjoint union of open intervals. We extend $\mu_g$ to
$\bR\setminus C$ by affine homeomorphisms, which thus furnishes a well-defined homomorphisms $G\longrightarrow\Homeo_+(\bR)$.
\end{proof}

We mention, though we will not explicate, that the constructions in Proposition~\ref{prop:l-order-homeo} are mutual inverses of each other,
suitably interpreted (cf.~Theorem 2.2.19 in~\cite{Navas2011}).
Moreover, the passage from a left invariant ordering on $G$ to an action of $G$ on $\bR$ is essentially unique,
up to semi-conjugacy.

It is easy to see that a group is left-orderable if and only if all of its finitely generated subgroups are left-orderable. 
One can prove this using the fact that the product space $2^G$ for a group $G$ is compact,
and the fact that the subspace $C(H)\sse 2^G$ of subsets inducing positive cones 
for a finitely generated subgroup $H$ is a closed subspace; see~\cite{Calegari2007,CR2016}, for instance. 
In fact, a stronger theorem due to Burns and Hale~\cite{BH1972}
holds: a group is orderable if and only if all of its nontrivial finitely generated subgroups surject onto some nontrivial left-orderable groups.

We can now observe the following from Proposition~\ref{prop:l-order-homeo}.
\begin{cor}[\cite{BH1972}]
A group is left-orderable if and only if all of its finitely generated (or equivalently, countable) 
subgroups admit faithful topological actions on the real line.
\end{cor}

\subsection{Circular orderability and the circle}\label{ss:circular}

For groups acting on the circle, there is an algebraic characterization that is 
analogous to that for groups acting on the real line. A \emph{circular
order}\index{circular ordering}
on $G$ is somewhat more complicated to define than a (linear) ordering on $G$, having to do with the fact that there is no canonical
way to order pairs of points on the circle. Instead, a circular ordering is a function $\omega$ that takes on values in the set $\{-1,0,1\}$,
and which is defined on triples of points in $G$. The function $\omega$ can be thought of as choosing an orientation on $S^1$. Let
$(a,b,c)$ be a triple of points on $S^1$. If traveling from $a$ to $b$ to $c$ coincides with the chosen orientation on $S^1$, then $\omega$
assigns the value $1$ to the triple, and $-1$ otherwise. Clearly, permuting the triple $(a,b,c)$ cyclically should not change the output
of $\omega$, though switching two points should change the sign by $-1$. Therefore, we require a circular ordering $\omega$
on an arbitrary set $X$ to satisfy the following
axioms:
\begin{enumerate}
\item
The function $\omega$ is defined on all triples of $X$, and is nonzero precisely on triples of distinct elements of $X$.
\item
For $\{a,b,c\}\sse  X$, we have \[\omega(a,b,c)=-\omega(a,c,b).\]
\item
For $\{a,b,c\}\sse  X$, we have \[\omega(a,b,c)=\omega(b,c,a).\]
\item
For $\{a,b,c,d\}\sse  X$, if \[\omega(a,b,c)=\omega(a,c,d)=1\] then $\omega(a,b,d)=1$.
\end{enumerate}

The second and third axioms say that for $\sigma\in S_3$, we have \[\omega(\sigma(a),\sigma(b),\sigma(c))=
\mathrm{sgn}(\sigma)\omega(a,b,c),\] where $\mathrm{sgn}$ denotes the sign representation of $S_3\longrightarrow\bZ/2\bZ$.

The last of these axioms, sometimes called the \emph{transitivity}\index{transitivity of a circular ordering} axiom,
is a compatibility condition for pairs of triples, and is intuitive given the interpretation of $\omega$ we have
given. Sometimes, this axiom is written in an invariant way as a \emph{cocycle}\index{cocycle}
axiom. That is, for $\{a,b,c,d\}\sse  X$, we have
\[\omega(b,c,d)-\omega(a,c,d)+\omega(a,b,d)-\omega(a,b,c)=0.\] To verify the equivalence, it is easy to see that if the quadruple is
degenerate (e.g.~if $c=d$) then the cocycle condition just says $0=0$. Otherwise, the cocycle axiom is just saying that
\[\omega(b,c,d)+\omega(a,b,d)=\omega(a,b,c)+\omega(a,c,d).\] Suppose first that \[\omega(a,b,c)+\omega(a,c,d)=2,\] with
the case of the opposite sign being analogous. The transitivity axiom 
combines with the other three to determine $\omega$ on all triples in $\{a,b,c,d\}$, and so
we have that \[\omega(d,c,b)=\omega(d,b,a)=-1,\] and so we get \[\omega(b,c,d)+\omega(a,b,d)=2\] by the change of sign under
transposition rule. Similarly, if \[\omega(a,b,c)+\omega(a,c,d)=0\] then \[\omega(a,b,d)+\omega(b,c,d)=0.\] Thus, the axioms for the
circular ordering imply the cocycle axiom. Conversely, suppose that the first three axioms together with the cocycle axiom hold.
The transitivity axiom follows easily, since if $\omega(a,b,c)$ and $\omega(a,c,d)$ are both equal to one then \[\omega(a,b,c)+\omega(a,c,d)
=2=\omega(a,b,d)+\omega(b,c,d),\] and so neither of the terms on the right hand side can be equal to $-1$.

The reader may wonder whether the suggestively named cocycle axiom is related to group cohomology, 
and indeed it is. We refer the reader to
Section 2.7 of~\cite{Calegari2007}, or to~\cite{FrigerioBook}.

Circular orderings are trickier than linear orderings, and so there is not a straightforward characterization of circular orderings as there
is for linear orderings via positive cones. See~\cite{BS2015,CGjlms19,CMR18}, for example.
We say that a circular ordering on a group $G$ is \emph{left invariant} if
for all $\{a,b,c,d\}\sse  G$, we have \[\omega(a,b,c)=\omega(da,db,dc).\] As with linear orderings, we have analogous notions of right
invariant and bi--invariant circular orderings on a group.
The following result characterizes countable circularly orderable
groups as subgroups of $\Homeo_+(S^1)$. The argument is essentially identical, and we leave the details to the reader.

\begin{prop}\label{prop:c-order-homeo}
Let $G$ be a countable group. Then $G$ admits a left invariant circular ordering if and only if $G$ admits a faithful homomorphism
into $\Homeo_+(S^1)$.
\end{prop}

We note that in light of Proposition~\ref{prop:l-order-homeo} and Proposition~\ref{prop:c-order-homeo}, and since $\Homeo_+[0,1]$ is identified
with a subgroup of $\Homeo_+(S^1)$, we have that linearly orderable groups are circularly orderable. To see this without
reference to homeomorphisms, if $G$ admits a left invariant linear ordering $<$, then for $\{a,b,c\}\sse  G$, we define $\omega(a,b,c)=1$
if $a<b<c$. The transformation of $\omega$ under permutations via the sign determines $\omega$ on all other triples of elements
of $G$.

\section{H\"older's Theorem}

H\"older's Theorem is one of the most fundamental results that relates dynamics of group actions on $I$ and $S^1$ with their algebraic
structure. We have relegated the discussion of H\"older's Theorem to this appendix since its proof relies in an essential way on orderability.

Let $M$ be a connected boundaryless one--manifold.
We say a group $G\le \Homeo_+(M)$ acts
\emph{freely}\index{free action} on $M$ if for all $g\in G$,
the existence of an $x\in M$ such that $g(x)=x$ implies that $g$ is the identity element. H\"older's
Theorem algebraically characterizes subgroups of $\Homeo_+(M)$ that act freely on $M$.

Let $<$ be a left invariant ordering on a group $G$. 
We say a subgroup $H$ of $G$ is \emph{bounded}\index{bounded subgroup} if there exists $g_0,g_1\in G$
such that every element in $H$ satisfies $g_0\le h\le g_1$.
It is easy to build orderings on groups that have nontrivial  bounded subgroups.
For instance, let $\tau\in\Homeo_+(\bR)$ denote
translation by one, and let $\sigma$ be an arbitrary nontrivial homeomorphism of $\bR$ that is supported in the interval $(0,1)$. We then get
a homeomorphism $\yt\sigma\in\Homeo_+(\bR)$ given by \[\yt\sigma=\prod_{i\in\bZ}\tau^{-i}\sigma\tau^i,\] and $\yt\sigma$ is easily
seen to commute with $\tau$. The subgroup of $\Homeo_+(\bR)$ generated by $\tau$ and $\yt\sigma$ is isomorphic to $\bZ^2$,
and $\bZ^2$ acquires an ordering as in Proposition~\ref{prop:l-order-homeo}. If $n\in\bZ$ and $x\in\bR$ then $|\yt\sigma^n(x)-x|<1$,
whereas $\tau(x)-x=1$. 
It follows easily now that the resulting ordering on $\bZ^2$ is has a bounded subgroup, namely $\form{\yt\sigma}$.

We say that $<$ is \emph{Archimedean}\index{Archimedean ordering} if the group does not contain a nontrivial bounded subgroup. Archimedean orderings are natural in the sense that they mirror usual intuitions about orders. 
If $G$ is realized as an additive subgroup of $\bR$, then the natural ordering on $G$ coming from Proposition~\ref{prop:l-order-homeo} is Archimedean.

Suppose $G\le \Homeo_+(\bR)$ acts freely on $\bR$. We can run the construction of a left invariant ordering from Proposition~\ref{prop:l-order-homeo} on $G$, though the argument is simplified. Indeed, if $a,b\in\bR$ and $g,h\in G$
are distinct elements, then $g(a)\neq h(a)$ and \[g(a)<h(a)\quad\textrm{if and only if}\quad g(b)<h(b).\] The first of these claims follows immediately from the freeness
of the action of $G$. The second holds since otherwise the Intermediate Value Theorem implies that for some $c\in (a,b)$, we have $g(c)=h(c)$ and so $g^{-1}h(c)=c$. Thus, one can order $G$ by comparing actions of group elements on an arbitrary point of $\bR$. It is easy to see then that the resulting ordering is Archimedean. By precomposing with an arbitrary homeomorphism of $\bR$, it is
immediate that the ordering is also right invariant. Thus, we see:

\begin{prop}
If a group acts freely and topologically on the real line, then it admits
an Archimedean bi--invariant ordering.
\end{prop}

We now have:

\begin{lem}\label{lem:bounded}
If a group $G$ admits a bi--invariant ordering then it admits a group homomorphism $\bar e$ into the additive group of the reals such that
\[
\ker\bar e=\{g\in G\mid \form{g}\text{ is bounded}\}.\]
\end{lem}

This result trivially implies that if a group admits a bi--invariant Archimedean ordering then it is isomorphic to an additive subgroup of the reals.
In particular, we deduce the following.

\begin{thm}[H\"older's Theorem]\label{thm:holder}
A group acting freely on the real line is abelian.\end{thm}

The basic idea of Lemma~\ref{lem:bounded} is to show that each element $g$ in the given group $G$ have a translation of $\bR$ 
associated to it in a canonical way; a similar idea appears again the proof of the Thurston Stability 
Lemma (Theorem~\ref{thm:thurston-stab}). Once such an association is made, the fact that one gets an embedding of $G$ 
into the additive group of reals is almost a formality. To find the suitable translation for $g$, one simply looks at the behavior 
of high powers of $g$ and renormalizes, with the hope that one
gets a nonzero translation for a nontrivial element of $G$.

\begin{proof}[Proof of Lemma~\ref{lem:bounded}]
Let $G$ be a group with a bi--invariant ordering $<$.
If every cyclic subgroup of $G$ is bounded, then we may set $\bar e$ to be trivial.
So, we may fix an element $h\in G$ such that $\form{h}$ is not bounded.
Let us first define
for all $g\in G$ that
\[
e(g):=\max\{n\in\bZ\mid h^n\le g\}.\]
For $f,g\in G$, we note from the bi--invariance of the ordering that
\[
h^{e(f)+e(g)}
= h^{e(f)}\cdot h^{e(g)}\le fg
< h^{e(f)+1}\cdot h^{e(g)+1}=h^{e(f)+e(g)+2}.\]
We have that
\[ 
e(f)+e(g)\le e(fg)\le e(f)+e(g)+1,\]
which implies 
 $e\co G\longrightarrow\bZ$ is a quasimorphism; see Lemma~\ref{lem:qm}.
It follows from the same lemma that we have a homogenization
\[
\bar e(g):=\lim_{n\to\infty} \frac{e(g^n)}n,\]
for each $g\in G$. Note that $\bar e(g^m)=m\bar e(g)$ for all $m\in\bZ$.

We note that $e\co G\longrightarrow \bR$ is monotone; that is, if $f<g$ then $e(f)\le e(g)$.
Let $f,g\in G$, assuming $fg\le gf$ without loss of generality.
We have
\[
e(f^n)+e(g^n)\le e(f^ng^n)\le e( (fg)^n)\le e(g^nf^n)\le  e(g^n) + e(f^n)+1.\]
Homogenizing, we have that $\bar e$ is a homomorphism. Obviously, $\bar e$ is monotone.

Let $K$ be the set of $f\in G$ such that $\form{f}$ is bounded.
Suppose $g\not\in K$.  There exists an integer $m$ such that $h\leq g^m$. This implies that \[
m\bar e(g)=\bar e(g^m)\ge 1.\]
We have shown $\ker\bar e\sse K$.The reverse direction is also easy to check.\ep

\begin{rem}\label{rem:conrad}
A bi-invariant ordering is a special case of a \emph{Conradian ordering}\index{Conradian ordering},
which is a left invariant ordering $<$ on a group such that all positive elements $f,g$ satisfy $fg^2>g$.
A group is  \emph{locally indicable}\index{locally indicable} if every nontrivial finitely generated surjects onto $\bZ$.
The idea behind Lemma~\ref{lem:bounded} can be strengthened to prove that a group is locally indicable 
if and only if it admits a Conradian ordering~\cite{Brodskii1984,RR2002,Navas2010Conrad}.\end{rem}

To prove an analogous result for the circle, suppose now that $G\le\Homeo_+(S^1)$ acts freely on $S^1$. Lifting $G$ to the real line, we obtain
a subgroup $\yt G\le \Homeo_+(\bR)$ commuting with $\tau=\tau_1$, i.e.~translation by one, and $\yt G$ acts freely on $\bR$. We may
repeat the argument for Theorem~\ref{thm:holder}
for the group $\form{\yt G,\tau}$, setting $\tau=h$ for the normalization. We thus obtain a homomorphism \[\phi\colon\yt G
\longrightarrow \bR\] which is injective, and which commutes with $\phi(\tau)=\tau$. It follows that the image of $\phi$ projects to
$\SO(2)$, so that $G$ is isomorphic to a subgroup of the circle group. The rotation number restricts to a homomorphism on $G$,
since $G$ preserves a (conjugate of) Lebesgue measure and so Theorem~\ref{thm:invt-homo} applies. Since $G$ acts freely on
$S^1$, Proposition~\ref{prop:rot-easy} implies that the rotation number is injective when restricted to $G$.

\begin{thm}\label{thm:holder-circle}
Let $G\le \Homeo_+(S^1)$ act freely on $S^1$. Then $G$ is a subgroup of the circle group $S^1$, and this isomorphism is realized
by the rotation number. \end{thm}
It obviously follows from the theorem that such a $G$ is abelian, and every
finite subgroup of $G$ is cyclic.

We have the following consequence of Theorem~\ref{thm:holder-circle}. The same statement can be found as Corollary 2.3
in~\cite{KK2018JT}, and a similar one is found in~\cite{FF2001} as Theorem 2.2.

\begin{cor}\label{cor:holder-extend}
Let $\varnothing\neq X\sse  S^1$ be closed, and let $G\le \Homeo_+(X)$, where the orientation on $X$ is just the circular order
inherited from $S^1$. If $G$ acts freely on $X$ then $G$ extends to a free action on $S^1$. Moreover, the rotation number restricts
to an injective homomorphism on $G$.
\end{cor}
\begin{proof}
If $J$ is a component of $S^1\setminus X$, then so is $g(J)$ for $g\in G$. We may therefore use the usual affine 
equivariant family to extend the action
of $G$ on $X$ to an action \[\phi\colon G\longrightarrow \Homeo_+(S^1).\] It is straightforward to check that this action is free, and so
Theorem~\ref{thm:holder-circle} implies that the rotation number furnishes an injective homomorphism $G\longrightarrow S^1$.
\end{proof}

\chapter{The Thurston Stability Theorem}\label{ch:append3}
In 1970s, W. P. Thurston announced the following striking theorem regarding classification of $C^1$--foliations.
\begin{thm}[Thurston Stability Theorem~\cite{Thurston1974Top}]
If $(M,F)$ is a compact, connected, transversely oriented codimension--one $C^1$--foliation with a compact leaf $L$ satisfying $H^1(L;\bR)=0$,
then  as a foliation it is isomorphic to an $L$--bundle over $S^1$ or to $L\times I$.
\end{thm}

Using a standard foliation theoretic technique, Thurston deduced the theorem from the following simpler result.

\begin{thm}
Let $G$ be a finitely generated group acting nontrivially by $C^1$ diffeomorphisms on $\bR^n$ fixing the origin.
If $g'=\Id$ for all $g\in G$, then $G$ surjects onto $\bZ$.\end{thm}

At the heart of Thurston's proof lies in the intuition that each generator of $G$ 
``becomes more and more nearly a translation'' (see~\cite{Thurston1974Top}) near the origin.
Using this phenomenon, he builds a homomorphism from $G$ to $\bR^n$ corresponding to such an approximate translation.
This idea is well represented by a one-dimensional result below, which is also sometimes called the Thurston Stability Theorem
 in the literature.
Recall from Remark~\ref{rem:conrad} that a group $G$ is \emph{locally indicable}\index{locally indicable} if every {nontrivial} 
finitely generated subgroup of $G$ admits a surjective homomorphism to
$\bZ$.

\begin{thm}[Thurston Stability Theorem in dimension one]\label{thm:thurston-stab}
The group $\Diff_+^1[0,1)$ is locally indicable.
\end{thm}

This result fits naturally with the $C^1$--theory of diffeomorphisms as discussed in Chapter~\ref{ch:c2-thry}, though since it is peripheral to the discussion of critical regularity, we have relegated it to this appendix. Thurston's Stability Theorem has several important immediate corollaries.

\begin{cor}
Let $G$ be a finitely generated nontrivial group acting on $S^1$ with a global fixed point. Then $G$ admits a surjective homomorphism to $\bZ$.
\end{cor}

The following consequence of Thurston's Stability Theorem shows why finding finitely generated simple subgroups of $\Homeo_+(I)$
is relatively difficult (see~\cite{HL19,HLNR21,MBT20}).
\begin{cor}
Let $G\le \Homeo_+(I)$ be finitely generated and perfect. Then $G$ admits no nontrivial homomorphism to $\Diff_+^1(I)$.
\end{cor}

Indeed, most functions that can be ``written down" in a straightforward way are going to be differentiable, but if one wants to write
down homeomorphisms of the interval that generate a simple group, they must exhibit non-differentiability in an essential way.
Thurston's Stability Theorem also illustrates why the commutator subgroup of Thompson's group $F$ must fail to be finitely generated.
We know it is simple by Lemma~\ref{lem:f-simple}, and we know that $F$ can be realized as a group of $C^{\infty}$ diffeomorphisms of
the interval by Corollary~\ref{cor:ghys-serg}. The failure of finite generation of the commutator subgroup $F'$ is the only way to prevent
a contradiction with Theorem~\ref{thm:thurston-stab}.


The starting point is the same as the dynamical interpretation of the abelianization of Thompson's group $F$ (see Subsection~\ref{ss:2-chain}).
For $F$, we realize $F$ in the usual way as a subgroup of $\PL_+(I)$, and we can evaluate the slope of an element $f\in F$ at the endpoints of the interval $I$. The chain rule says that composition of elements of $F$ results in multiplication of slopes, and taking the logarithm of the slope furnishes a homomorphism to the additive group $\bR$.
Because of the nature of the explicit description of $F$, it is immediate that this homomorphism is nontrivial at both endpoints of $I$
and hence furnishes a surjective homomorphism $F\longrightarrow\bZ^2$, which in turn coincides with the abelianization of $F$.

We also consider the map $\phi\co g\mapsto \log g'(0)$ in Theorem~\ref{thm:thurston-stab}. 
The key case is when this map is trivial, where we use the idea of ``approximate linearization''. This roughly means that if $g'(0)=1$ the $g$ behaves like a translation; see the remark following the proof below.

\begin{proof}[Proof of Theorem~\ref{thm:thurston-stab}]
Let $G$ be a subgroup of $\Diff_+^1[0,1)$ with a finite generating set $S$. We may assume $S$ is symmetric, that is $S=S^{-1}$.
If some $s\in S$ is not $C^1$--tangent to the identity (i.e. $s'(0)\ne1$), then the homomorphism $\phi$ above is already nontrivial.
So, we may assume that $s'(0)=1$ for all $s\in S$; alternatively, this can be guaranteed by the Muller--Tsuboi trick (Theorem~\ref{t:muller-tsuboi}).

We can also assume that at least one element of $S$ does not agree with the
identity in an any neighborhood of $0$ (i.e.~the germ of some $s\in S$ at $0$ is nontrivial).  
This is not a loss of generality, since otherwise
we replace $0$ with 
\[y:=\sup \{x\in[0,1)\mid  s\restriction_{[0,x]}\text{ is the identity map for all }s\in S\}.\]
It is clear that $y\in [0,1)$ provided that $G$ is nontrivial.

For each $\ell\ge0$ and $x\in[0,1)$, we define
\begin{align*}
A_\ell(x)&:=\max \{ t_1\cdots t_i(x) \mid 0\le i\le \ell\text{ and }t_1,\ldots,t_i\in S\}\in[0,1),\\
M(x)&:=\max_{s\in S}|s(x)-x|,\\
T(x)&:=\max_{s\in S}\sup_{0\le t\le x} |s'(t)-1|.
\end{align*}
For each fixed $\ell\in\bN$, we have $\lim_{x\to+0} A_\ell(x)=0$
since $0$ is the global fixed point of $G$.
We can  find some sequence $\{x_n\}\sse(0,1)$ converging to $0$ such that $M(x_n)\ne0$ for all $n$.
Note from our assumption that $\lim_{x\to+0}T(x)=0$.

For each $g\in G$ and $x\in [0,1)$ we denote the displacement of $g$ at $x$ as
\[
\Delta_x g:=g(x)-x.\]
For $f,g\in G$, we see from the Mean Value Theorem that
\[
\Delta_x fg=\Delta_{g(x)}f+\Delta_xg=\Delta_xf+\Delta_xg+(f'(z)-1)\Delta_xg\]
for some $z$ between $x$ and $g(x)$. Inductively, if $t_1,\ldots,t_\ell\in S$ we have
\[
\Delta_x t_1\cdots t_\ell
=\sum_{i=1}^\ell \Delta_x t_i + \sum_{i=1}^{\ell-1} \left( t_1'(z_1)\cdots t_i'(z_i)-1\right) \Delta_x t_{i+1}\]
for some $z_1,\ldots,z_{\ell-1}\in [0,A_\ell(x)]$.
Note the general inequality \[|(1+t)^n-1|\le \left(2+2|t|\right)^n\cdot |t|.\]
Setting
\[\eta(\ell,x):= \ell \left(2+2T(A_\ell(x))\right)^\ell T(A_\ell(x)),\]
we deduce that
\[
\left|\Delta_xt_1\cdots t_\ell -\sum_i \Delta_x t_i\right|
\le \eta(\ell,x) M(x).\]
Moreover, for each fixed $\ell$ we see
\[
\lim_{x\to+0}\eta(\ell,x)=0.\]

Note that the sequence $\{\Delta_{x_n}(s)/M(x_n)\}$ is bounded for each $s\in S$.
We can take a subsequence if necessary and assume that the following limit exists for each $s\in S$:
\[
\bar\Delta(s):=\lim_{n\to\infty} \frac{\Delta_{x_n}(s)}{M(x_n)}.\]

Let $g\in G$ be arbitrary.
Writing $g=t_1\cdots t_\ell$ for some $t_i\in S$, we see from the above estimate that
\[
\left| \frac{\Delta_{x_n}(g)}{M(x_n)} - \sum_i \frac{\Delta_{x_n}(t_i)}{M(x_n)}\right| \le\eta(\ell,x_n).\]
It follows that we have a limit
\[
\bar\Delta(g):=\lim_{n\to\infty} \frac{\Delta_{x_n}(g)}{M(x_n)}=\sum_i \bar\Delta(t_i).\]
This expression implies that the map $\bar\Delta\co G\longrightarrow\bR$ is a homomorphism,
which is nontrivial since $\max_{s\in S}\bar\Delta(s)=1$ by design.\ep
\begin{rem}
The estimate above comes from an interpretation of Thurston Stability Theorem due to Bonatti~\cite{BonattiProches,BMNR2017MZ}.
In fact, he proved that for an arbitrary compact manifold $M$, for each real number $\eta>0$, and for each integer $\ell\in\bN$,
 there exists a neighborhood of the identity $\VV\sse \Diff^1(M)$ such that for all points $x\in M$, for all diffeomorphisms
$f_1,\ldots,f_\ell \in \VV$
and for all $\epsilon_1,\ldots,\epsilon_\ell\in \{-1,1\}$,
we have that
\[
\left \|
f_k^{\epsilon_k}\circ \cdots\circ f_1^{\epsilon_1}(x) - x
- \sum_{i=1}^k \epsilon_i (f_i(x)-x)\right \|\le \eta \max_{i=1,\ldots,k} \left\|f_i(x)-x\right \|.\]
This estimate may be regarded as an approximate linearization of diffeomorphisms in small $C^1$--neighborhoods of the identity.
\end{rem}

As might be imagined, Theorem~\ref{thm:thurston-stab} furnishes finitely generated subgroups of $\Homeo_+(I)$ that are not
isomorphic to subgroups of $\Diff_+^1(I)$, though of course they can be conjugated to groups of Lipschitz homeomorphism by
Theorem~\ref{thm:lip-conj}. For an explicit such example, let $G_0\le \PSL_2(\bR)$ denote the $(2,3,7)$ triangle group.
Since $\PSL_2(\bR)\le \Homeo_+(S^1)$, we lift $G_0$ to a subgroup $\Homeo_+^{\bZ}(\bR)$, consisting of all lifts of elements of $G_0$.
This group has the presentation \[G=\form{ x,y,z,t\mid x^2=y^3=z^7=t,\, xyz=t}.\] We have that $G\le \Homeo_+(I)$ and
the commutator subgroup $G'$ coincides with $G$, so that $G$ admits no nontrivial homomorphisms of abelian groups.
It follows that $G$ cannot be realized as a group of diffeomorphisms of $I$. This example was first noticed by Kropholler and Thurston
 (cf.~Bergman~\cite{Bergman1991PJM} and Example 2.120
of~\cite{Calegari2007}).

\backmatter


\printindex

\bibliographystyle{amsplain}
\bibliography{ref}
\end{document}